\documentclass[12pt,a4paper]{report}
\RequirePackage[l2tabu, orthodox]{nag}
\usepackage{geometry}
\usepackage[english, strings]{babel}
\usepackage{nchairx}
\usepackage[utf8]{inputenc}
\usepackage[T1]{fontenc}       % also for the font encoding
\usepackage{longtable}         % tables longer than one page
\usepackage{exscale}           % large summation signs in 11pt
\usepackage[sort]{cite}        % nicer citations
\usepackage{xspace}            % better spacing after macros
\usepackage{tikz}              % for commutative diagrams and stuff
\usetikzlibrary{cd, decorations.pathmorphing}
\usepackage{ifdraft}           % to determine whether draft mode
\usepackage{appendix}
\usepackage[expansion=false    % no font expansion
]{microtype}        % only protrusion
\usepackage{imakeidx} % index
\usepackage[backref=page,      % backrefs in the bibliography
final=true,         % always treat as final
pdfpagelabels,       % use pdf page labels
]{hyperref}         % hyperrefs are cool!
\mathtoolsset{showonlyrefs,showmanualtags}
\usepackage{enumitem}
\usepackage{epigraph} %for quotes
\usepackage{nicematrix} %nicer matrices
\usepackage{graphicx} %figures
\counterwithout{figure}{chapter} %Continuous figure numbering
\usepackage{afterpage} %Empty pages
\usepackage{float}

\usepackage{parskip}
\usepackage[symbols,nogroupskip,sort=use,toc,section=section]{glossaries-extra}
\usepackage{glossary-longragged} % table style for list of symbols
\usepackage{glossary-longbooktabs} %more table styles for list of symbols
\glsdisablehyper % no hyperlinks to list of symbols

\usepackage{booktabs}% for Terminology table

\usepackage[en-US]{datetime2} % date
\DTMlangsetup{showdayofmonth=false} % only month and year

\newenvironment{dedication}
{%\clearpage           % we want a new page
    \thispagestyle{empty}% no header and footer
    \vspace*{\stretch{1}}% some space at the top
    \itshape             % the text is in italics
    \raggedleft          % flush to the right margin
}
{\par % end the paragraph
    \vspace{\stretch{3}} % space at bottom is three times that at the top
    \clearpage           % finish off the page
}

\makeatletter

\newcommand{\bibnote}[2]{\nocite{#1}\@namedef{#1chairxnote}{#2}}
\makeatother

\newcommand\blankpage{%
    \null
    \thispagestyle{empty}%
    \addtocounter{page}{-1}%
    \newpage}

\renewcommand{\tilde}{\widetilde}

\newcommand{\disk}{\mathbb{D}}

\newcommand{\C}{\field{C}}

\newcommand{\R}{\field{R}}

\newcommand{\N}{\field{N}}

\newcommand{\Z}{\field{Z}}

\renewcommand{\epsilon}{\varepsilon}

\newcommand{\quant}{\operatorname{quant}}

\DeclareFontFamily{U}{mathx}{}
\DeclareFontShape{U}{mathx}{m}{n}{ <-> mathx10 }{}
\DeclareSymbolFont{mathx}{U}{mathx}{m}{n}
\DeclareFontSubstitution{U}{mathx}{m}{n}
\DeclareMathAccent{\widecheck}{0}{mathx}{"71}

\renewcommand{\Holomorphic}{\mathscr{H}}

\newcommand{\pos}{{\operatorname{pos}}}

\newcommand{\mom}{{\operatorname{mom}}}

\newcommand{\Linear}{\operatorname{L}}

\newcommand{\GL}{\operatorname{GL}}

\newcommand{\cs}{\operatorname{cs}}

\newcommand{\polar}{\textrm{\tiny \fontencoding{U}\fontfamily{ding}\selectfont\symbol{'136}}}

\newcommand{\HS}{\operatorname{HS}}

\renewcommand{\Sec}[1][\infty]{\Gamma^{#1}}

\newcommand{\Bdd}{\mathfrak{B}}

\newcommand{\absconv}{\operatorname{absconv}}

\newcommand{\Gutt}{\mathrm{G}}

\newcommand{\liegroup}[1]{\operatorname{#1}}

\newcommand{\Taylor}{\operatorname{T}}

\newcommand{\Majorant}{\operatorname{F}}

\newcommand{\LieAlg}{\operatorname{Lie}}

\newcommand{\Entire}{\mathscr{E}} % a stylefile in the making
\glsxtrnewsymbol[description={Composition of
mappings}]{Composition}{\ensuremath{\, \circ \,}}

\glsxtrnewsymbol[description={Topological dual space of topological vector space
$V$, space of continuous linear mappings $V \longrightarrow
\C$}]{DualSpace}{\ensuremath{V'}}

\glsxtrnewsymbol[description={Space of all mappings with complex
values}]{Mappings}{\ensuremath{\Map}}

\glsxtrnewsymbol[description={Euler's number}]{Euler}{\ensuremath{e}}

\glsxtrnewsymbol[description={Limit taken within the positive
real numbers}]{LimitFromAbove}{\ensuremath{\; \downarrow \;}}

\glsxtrnewsymbol[description={Supremum norm of a sequence $\gamma$, $\sup_{n \in \N}
\abs{\gamma(n)}$}]{Supnorm}{\ensuremath{\norm{\gamma}_\infty}}

\glsxtrnewsymbol[description={Supremum norm of a vector or
sequence}]{SupnormVector}{\ensuremath{\norm{\argument}_\infty}}

\glsxtrnewsymbol[description={Supremum norm of a function $\phi$ on a compact set
$K$}]{SupnormCompact}{\ensuremath{\norm{\phi}_K}}

\glsxtrnewsymbol[description={Natural pairing, scalar product on Hilbert
space}]{Pairing}{\ensuremath{\langle \argument, \argument \rangle}}

\glsxtrnewsymbol[description={Completion of topological (injective, projective,
Hilbert) tensor product}]{TensorCompletion}{\ensuremath{\widehat{\tensor}}}

\glsxtrnewsymbol[description={Poisson
bracket}]{PoissonBracket}{\ensuremath{\{\argument, \argument\}}}

\glsxtrnewsymbol[description={Densely defined linear adjointable operators on some
Hilbert space $H$}]{Adjointable}{\ensuremath{\mathfrak{A}(H)}}

\glsxtrnewsymbol[description={Quantization map}]{QuantizationMap}{\ensuremath{Q}}

\glsxtrnewsymbol[description={Binomial coefficient, $n$ choose
$k$,$\tfrac{n!}{(n-k)! \cdot k!}$}]{Binomial}{\ensuremath{\binom{n}{k}}}

\glsxtrnewsymbol[description={Constant function with value $1 \in
\C$}]{ConstantFunction}{\ensuremath{\mathbb{1}}}

\glsxtrnewsymbol[description={The factor
$\xi_p$ is absent}]{AbsentFactor}{\ensuremath{\overset{\xi_p}{\cdots}}}

\glsxtrnewsymbol[description={Kostant-Kirillov-Souriau Poisson bracket on
$\Pol(\liealg{g}^*)$}]{KostantKirillovSouriau}{\ensuremath{\{\argument,\argument\}_\KKS}}

\glsxtrnewsymbol[description={Seminorms corresponding to the topology of uniform
convergence of bounded subsets, \eqref{eq:BoundedSeminorms}}]{SeminormsBounded}
{\ensuremath{\seminorm{p}_{B,\seminorm{q}}}}

\glsxtrnewsymbol[description={Alternative seminorms corresponding to the topology
of uniform convergence of bounded subsets on the space of bounded Fréchet
holorphic functions}]{SeminormsBoundedAlternative}
{\ensuremath{\seminorm{r}_{B, \seminorm{q}, v_0}}}

\glsxtrnewsymbol[description={Seminorms corresponding to the $R$-topology,
\eqref{eq:RTopologySeminorms}}]{SeminormsRTopology}
{\ensuremath{\seminorm{p}_{R,c}}}

\glsxtrnewsymbol[description={Seminorms corresponding to the topology of the entire
functions~$\Entire_0(G)$, \eqref{eq:EntireSeminorms}}]{SeminormsEntire}
{\ensuremath{\seminorm{q}_{0,c}}}

\glsxtrnewsymbol[description={Seminorms corresponding to smooth functions from a
Lie group $G$ to a vector space $V$,
\eqref{eq:SmoothVectorsSeminorms}}]{SeminormsSmooth}
{\ensuremath{\seminorm{q}_{\xi, K, \seminorm{p}}}}

\glsxtrnewsymbol[description={Seminorms corresponding to entire
vectors}]{SeminormsEntireVectors}
{\ensuremath{\seminorm{q}_{c,\seminorm{p}}}}

\glsxtrnewsymbol[description={Seminorms corresponding to holomorphic extensions of
entire vectors}]{SeminormsEntireVectorsExtensions}
{\ensuremath{\seminorm{q}_{K,\seminorm{p}}}}

\glsxtrnewsymbol[description={Complex plane}]{Complex}{\ensuremath{\C}}

\glsxtrnewsymbol[description={Real part of a complex
number}]{RealPart}{\ensuremath{\RE(\argument)}}

\glsxtrnewsymbol[description={Imaginary part of a complex
number}]{ImaginaryPart}{\ensuremath{\IM(\argument)}}

\glsxtrnewsymbol[description={Dimension of a manifold}]{Dimension}{\ensuremath{\dim}}

\glsxtrnewsymbol[description={Reduced Planck's constant, element of deformation domain
$\mathfrak{D} \subseteq \C$}]{Planck}{\ensuremath{\hbar}}

\glsxtrnewsymbol[description={Deformation
domain, subdomain of $\C$}]{DeformationDomain}{\ensuremath{\mathfrak{D}}}

\glsxtrnewsymbol[description={Landau Big O, Order of
approximation}]{Order}{\ensuremath{\mathcal{O}}}

\glsxtrnewsymbol[description={Set of continuous seminorms defined on a locally convex
space $V$}]{ContinuousSeminorms}{\ensuremath{\cs(V)}}

\glsxtrnewsymbol[description={Zero set of a linear map or seminorm,
kernel}]{Kernel}{\ensuremath{\ker \,}}

\glsxtrnewsymbol[description={Local Banach space $W / \ker \seminorm{q}$ at
$\seminorm{q} \in \cs(W)$}]{LocalBanach}{\ensuremath{W_{\seminorm{q}}}}

\glsxtrnewsymbol[description={Direct Limit}]{DirectLimit}{\ensuremath{\varinjlim
\,}}

\glsxtrnewsymbol[description={Permutation group in the letters $\{1,\ldots,n\}$}]{PermutationGroup}{\ensuremath{S_n}}

\glsxtrnewsymbol[description={Left action of a group on some
set}]{Action}{\ensuremath{\; \acts \;}}

\glsxtrnewsymbol[description={Polynomial, Polynomial function, diagonal of a
multilinear map, motive of
Section~\ref{sec:Polynomials}}]{Polynomial}{\ensuremath{P}}

\glsxtrnewsymbol[description={Polynomials in the
indeterminants $z_1, \ldots, z_n$ and complex
coefficients}]{PolynomialsAlgebraic}{\ensuremath{\C[z_1, \ldots, z_k]}}

\glsxtrnewsymbol[description={Symmetric part of a multilinear mapping,
\eqref{eq:SymmetricPart}}]{SymmetricPart}{\ensuremath{L_s}}

\glsxtrnewsymbol[description={Unique symmetric multilinear mapping inducing the
polynomial $P$}]{InducingPolynomial}{\ensuremath{\widecheck{P}}}

\glsxtrnewsymbol[description={Space of continuous
polynomials of homogeneity $n$ between locally convex spaces $V$ and
$W$}]{PolynomialsHomogeneous}{\ensuremath{\Pol^n(V,W)}}

\glsxtrnewsymbol[description={Graded vector space of continuous
    polynomials between locally convex spaces $V$ and
    $W$}]{Polynomials}{\ensuremath{\Pol^\bullet(V,W)}}

\glsxtrnewsymbol[description={Space of continuous
    polynomials of homogeneity $n$ between locally convex spaces $V$ and
    $W$, endowed with topology of uniform convergence on bounded
    sets}]{PolynomialsBounded}{\ensuremath{\Pol_\beta^n(V,W)}}

\glsxtrnewsymbol[description={Space of continuous functions with compact support,
test functions}]{TestFunctions}{\ensuremath{\Continuous_c}}

\glsxtrnewsymbol[description={Space of smooth functions with compact support,
    test functions}]{TestFunctionsSmooth}{\ensuremath{\Cinfty_c}}

\glsxtrnewsymbol[description={Taylor polynomial of degree $n \in \N_0$ with
expansion point $v_0$}]{TaylorPolynomial}{\ensuremath{P_{n,v_0}}}

\glsxtrnewsymbol[description={Gâteaux derivative of a function $f$ with expansion
point $v_0$ and in direction $v$}]{GateauxDerivative}{\ensuremath{df(v_0,v)}}

\glsxtrnewsymbol[description={Real analytic functions on some manifold
    $M$ with values in $\C$}]{RealAnalytic}{\ensuremath{\Comega(M)}}

\glsxtrnewsymbol[description={Real analytic functions on a Lie group $G$ with minimum
radius of convergence $r>0$ for the Lie-Taylor series around any $g \in
G$}]{RealAnalyticMinumumRadius}{\ensuremath{\Comega_r(G)}}

\glsxtrnewsymbol[description={Holomorphic functions on some complex manifold
$M_\C$ with values in $\C$}]{Holomorphic}{\ensuremath{\Holomorphic(M_\C)}}

\glsxtrnewsymbol[description={Frechet holomorphic functions defined on some open
    subset $U \subseteq V$ with values in locally convex space
    $W$}]{HolomorphicFrechet}{\ensuremath{\Holomorphic(U,W)}}

\glsxtrnewsymbol[description={Bounded Frechet holomorphic functions defined on
some open subset $U \subseteq V$ with values in locally convex space
    $W$, endowed with the topology of uniform convergence on bounded
    sets}]{HolomorphicBounded}{\ensuremath{\Holomorphic_\beta(U,W)}}

\glsxtrnewsymbol[description={Gâteaux holomorphic functions defined on some open
subset $U \subseteq V$ with values in locally convex space
$W$}]{HolomorphicGateaux}{\ensuremath{\Holomorphic_G(U,W)}}

\glsxtrnewsymbol[description={Symmetric tensor power of order $n \in \N$ of a vector
space $V$}]{SymmetricPower}{\ensuremath{\Sym^n(V)}}

\glsxtrnewsymbol[description={Projective symmetric tensor power of order $n \in \N$ of
a vector space $V$}]{SymmetricPowerProjective}{\ensuremath{\Sym^n_\pi(V)}}

\glsxtrnewsymbol[description={Tensor power of order $n \in
\N$}]{TensorPower}{\ensuremath{\Tensor^n}}

\glsxtrnewsymbol[description={Symmetric algebra over a vector space
$V$}]{SymmetricAlgebra}{\ensuremath{\Sym^\bullet(V)}}

\glsxtrnewsymbol[description={Symmetric algebra over a vector space
    $V$, endowed with the
    $R$-topology}]{RTopology}{\ensuremath{\Sym^\bullet_R(V)}}

\glsxtrnewsymbol[description={Symmetric algebra over a vector space
    $V$, endowed with the $0$-topology built from injective tensor
    products}]{RTopologyInjective}{\ensuremath{\Sym_\epsilon^\bullet(V)}}

\glsxtrnewsymbol[description={Symmetric tensor power of order $n \in \N$ of a
seminorm~$\seminorm{p}$}]{SymmetricSeminormPower}{\ensuremath{\seminorm{p}^n}}

\glsxtrnewsymbol[description={Completion of a Hausdorff locally convex space
$W$}]{Completion}{\ensuremath{\widehat{W}}}

\glsxtrnewsymbol[description={Symmetric tensors as polynomials on the dual,
canonical embedding, zero section}]{CanonicalEmbedding}{\ensuremath{\iota}}

\glsxtrnewsymbol[description={Continuous extension of the canonical embedding
$\iota$}]{CanonicalEmbeddingExtension}{\ensuremath{\hat{\iota}}}

\glsxtrnewsymbol[description={Semidirect product}]{Semidirect}{\ensuremath{\; \ltimes
\;}}

\glsxtrnewsymbol[description={Direct
sum over $\N_0$}]{DirectSum}{\ensuremath{\bigoplus_{n=0}^\infty} \,}

\glsxtrnewsymbol[description={Cartesian product over $\N_0$}]{CartesianProduct}
{\ensuremath{\prod_{n=0}^\infty} \,}

\glsxtrnewsymbol[description={Disjoint union over $g \in G$}]{DisjointUnion}
{\ensuremath{\bigsqcup_{g \in G}} \,}

\glsxtrnewsymbol[description={Open unit disk, $\{z \in \C \colon \abs{z} <
1\}$}]{Disk}{\ensuremath{\mathbb{D}}}

\glsxtrnewsymbol[description={Open disk with radius $r > 0$, $\{z \in \C \colon
\abs{z} <
    r\}$}]{DiskRadiusR}{\ensuremath{\mathbb{D}_r}}

\glsxtrnewsymbol[description={Set of continuous functions between topological
spaces~$X$ and $Y$}]{Continuous}{\ensuremath{\Continuous(X,Y)}}

\glsxtrnewsymbol[description={Topology of uniform convergence on bounded
sets, strong topology}]{BoundedUniformConvergence}{\ensuremath{\beta}}

\glsxtrnewsymbol[description={Topology of pointwise convergence, uniform
convergence on finite sets, weak
topology}]{PointwiseConvergence}{\ensuremath{\sigma}}

\glsxtrnewsymbol[description={Norm, operator
norm}]{Norm}{\ensuremath{\norm{\argument}}}

\glsxtrnewsymbol[description={Banach space of zero sequences indexed by the natural
numbers $\N$ and norm $\supnorm{\argument}$}]{ZeroSequences}{\ensuremath{c_0}}

\glsxtrnewsymbol[description={Banach space of absolutely summable sequences indexed by
the natural numbers $\N$ and norm
$\norm{\argument}_1$}]{AbsolutelySummable}{\ensuremath{\ell^1}}

\glsxtrnewsymbol[description={Banach space of bounded sequences indexed by the natural
numbers $\N$ and norm
$\supnorm{\argument}$}]{BoundedSequences}{\ensuremath{\ell^\infty}}

\glsxtrnewsymbol[description={Sequence with entries $e_n(k) \coloneqq
\delta_{n,k}$}]{UnitSequence}{\ensuremath{e_n}}

\glsxtrnewsymbol[description={Cartesian product over an index set $J$ with factors
$\C$, sequences indexed by $J$ with entries in
$\C$}]{CartesianPower}{\ensuremath{\C^J}}

\glsxtrnewsymbol[description={Space of all sequences with entries in $\C$,
Cartesian product of countably many copies of
$\C$}]{Sequences}{\ensuremath{\C^\N}}

\glsxtrnewsymbol[description={Hilbert space of absolutely square summable
sequences with norm $\norm{\argument}_2$}]{SquareSummable}{\ensuremath{\ell^2}}

\glsxtrnewsymbol[description={Sequences indexed by $J$ with only finitely many
nonzero entries}]{CartesianDirectPower}{\ensuremath{\C^{(J)}}}

\glsxtrnewsymbol[description={Absolutely convex hull of a subset of a vector
space}]{AbsolutelyConvexHull}{\ensuremath{\absconv}}

\glsxtrnewsymbol[description={Absolute value of a complex number $\hbar \in
\C$}]{AbsoluteValue}{\ensuremath{\abs{\hbar}}}

\glsxtrnewsymbol[description={Dirac delta functional, evaluation
character}]{Dirac}{\ensuremath{\delta}}

\glsxtrnewsymbol[description={Partial derivative with respect to a real or complex
variable $\hbar$}]{PartialDerivative}{\ensuremath{\frac{\partial}{\partial \hbar}}}

\glsxtrnewsymbol[description={Pullback
with a function $f$}]{Pullback}{\ensuremath{f^*}}

\glsxtrnewsymbol[description={Set of absolutely convex bounded sets $B$ such
that the translation $v_0 + B \subseteq
U$.}]{Polydisks}{\ensuremath{\mathfrak{B}(U,v_0)}}

\glsxtrnewsymbol[description={Line segment between $v_0,v \in V$ in a vector
space $V$}]{LineSegment}{\ensuremath{[v_0,v]}}

\glsxtrnewsymbol[description={Hilbert-Schmidt
operators}]{HilbertSchmidt}{\ensuremath{\HS}}

\glsxtrnewsymbol[description={Special linear group over
$\R^2$}]{SpecialLinear}{\ensuremath{\liegroup{SL}_2(\R)}}

\glsxtrnewsymbol[description={Heisenberg group in one
spcial dimension}]{Heisenberg}{\ensuremath{H_1}}

\glsxtrnewsymbol[description={Space of two-by-two matrices
}]{Matrices}{\ensuremath{\Mat_2}}

\glsxtrnewsymbol[description={Lie algebra of $\liegroup{SL}_2(\R)$, traceless
two-by-two matrices}]{SpecialLinearAlgebra}{\ensuremath{\liealg{sl}_2(\R)}}

\glsxtrnewsymbol[description={Symmetrizer}]{Symmetrizer}{\ensuremath{\Symmetrizer}}

\glsxtrnewsymbol[description={Space of bounded
mappings between an open subset $U \subseteq V$ of a locally convex space with values
in another locally convex space $W$}]{Bounded}{\ensuremath{\Bounded}(U,W)}

\glsxtrnewsymbol[description={Space of all mappings}]{Map}{\ensuremath{\Map}}

\glsxtrnewsymbol[description={$\ell^1$-norm of a sequence
$\gamma$}]{Ell1Norm}{\ensuremath{\norm{\gamma}_1}}

\glsxtrnewsymbol[description={Open metric ball with radius $r > 0$ and center
$x$}]{OpenBall}{\ensuremath{\Ball_r(x)}}

\glsxtrnewsymbol[description={Open interval from $-R$ to
$R$}]{Interval}{\ensuremath{(-R,R)}}

\glsxtrnewsymbol[description={Closure of a set $S$}]{Closure}{\ensuremath{S^\cl}}

\glsxtrnewsymbol[description={Open cylinder with radius $r > 0$ and center $v_0$
corresponding to a seminorm
$\seminorm{p}$}]{OpenCylinder}{\ensuremath{\Ball_{\seminorm{p},r}(v_0)}}

\glsxtrnewsymbol[description={Polar of a subset $A$ of a topological vector space
$V$}]{Polar}{\ensuremath{A^\polar}}

\glsxtrnewsymbol[description={Tensor product}]{Tensor}{\ensuremath{\, \tensor \,}}

\glsxtrnewsymbol[description={Projective tensor
product}]{ProjectiveTensor}{\ensuremath{\, \tensor_\pi \,}}

\glsxtrnewsymbol[description={Injective tensor
product}]{InjectiveTensor}{\ensuremath{\tensor_\epsilon}}

\glsxtrnewsymbol[description={Identity mapping $\id_S(x) = x$ of a set $S$}]{Identity}{\ensuremath{\id_S}}

\glsxtrnewsymbol[description={Real linear mappings between vector spaces $V$ and $W$}]{Linear}{\ensuremath{\Linear(V,W)}}

\glsxtrnewsymbol[description={Complex linear mappings between complex vector spaces $V$ and $W$}]{CLinear}{\ensuremath{\Linear_\C(V,W)}}

\glsxtrnewsymbol[description={Complex antilinear mappings between complex vector spaces $V$ and $W$}]{CALinear}{\ensuremath{\cc{\Linear}_\C(V,W)}}

\glsxtrnewsymbol[description={Formal power series with coefficients in a ring
$\ring{R}$ and indeterminate
$\lambda$}]{FormalPowerSeries}{\ensuremath{\ring{R}\formal{\lambda}}}

\glsxtrnewsymbol[description={Product of algebra $\algebra{A}$,
pointwise product}]{Product}{\ensuremath{\cdot}}

\glsxtrnewsymbol[description={Star product}]{StarProduct}{\ensuremath{\star}}

\glsxtrnewsymbol[description={Convolution product}]{Convolution}{\ensuremath{\ast}}

\glsxtrnewsymbol[description={Wick star
product}]{StarProductWick}{\ensuremath{\star_{\operatorname{Wick}}}}

\glsxtrnewsymbol[description={Complex conjugation of a universal
complexification~$G_\C$}]{ComplexConjugationGroup}{\ensuremath{\sigma}}

\glsxtrnewsymbol[description={Complex conjugated scalar
multiplication}]{ScalarConjugated}{\ensuremath{\odot}}

\glsxtrnewsymbol[description={Automorphism
group}]{AutomorphismGroup}{\ensuremath{\Aut}}

\glsxtrnewsymbol[description={One sphere, boundary of the unit circle, unitary group
in one dimension}]{OneSphere}{\ensuremath{\mathbb{S}^1}}

\glsxtrnewsymbol[description={Covariant derivative, connection for tangent
bundle}]{Connection}{\ensuremath{\nabla}}

\glsxtrnewsymbol[description={Associative
algebra}]{Algebra}{\ensuremath{\algebra{A}}}

\glsxtrnewsymbol[description={Quantization of the algebra
$\algebra{A}$}]{AlgebraQuantum}{\ensuremath{\algebra{A}_{\quant}}}

\glsxtrnewsymbol[description={Subalgebra of $\algebra{A}$ with convergent star
product}]{AlgebraConvergent}{\ensuremath{\algebra{A}_{\conv}}}

\glsxtrnewsymbol[description={Tangent functor}]{Tangent}{\ensuremath{T}}

\glsxtrnewsymbol[description={Cotangent bundle of a
manifold}]{Cotangent}{\ensuremath{T^*}}

\glsxtrnewsymbol[description={Complexified tangent
bundle}]{TangentComplexified}{\ensuremath{T^\C}}

\glsxtrnewsymbol[description={Smooth functions on a manifold
$M$}]{SmoothFunctions}{\ensuremath{\Fun[\infty](M)}}

\glsxtrnewsymbol[description={Smooth sections of a vector
bundle}]{Sections}{\ensuremath{\Sec[\infty]}}

\glsxtrnewsymbol[description={Multiplication operator with algebra element $a \in
\algebra{A}$}]{MultiplicationOperator}{\ensuremath{M_a}}

\glsxtrnewsymbol[description={Lie bracket,
commutator}]{LieBracket}{\ensuremath{[\argument,
\argument]}}

\glsxtrnewsymbol[description={Filtered algebra of differential
operators}]{DiffOps}{\ensuremath{\Diffop^\bullet}}

\glsxtrnewsymbol[description={Smooth functions on a vector bundle, which are
polynomial in every fiber}]{PolynomialsVectorBundle}{\ensuremath{\Pol}}

\glsxtrnewsymbol[description={Unit of a group
$G$}]{GroupUnit}{\ensuremath{\E}}

\glsxtrnewsymbol[description={Left multiplication with a group element $g \in G$}]{LeftMultiplication}{\ensuremath{\ell_g}}

\glsxtrnewsymbol[description={Right multiplication with a group element $g \in G$}]{RightMultiplication}{\ensuremath{r_g}}

\glsxtrnewsymbol[description={Conjugation with a group element $g \in G$}]{Conjugation}{\ensuremath{\Conj_g}}

\glsxtrnewsymbol[description={Adjoint action of a group element $g \in G$ on
$\liealg{g}$}]{AdjointGroup}{\ensuremath{\Ad_g}}

\glsxtrnewsymbol[description={Group inversion, $g \mapsto
g^{-1}$}]{Inversion}{\ensuremath{\operatorname{inv}}}

\glsxtrnewsymbol[description={Representative function corresponding
to a representation $\pi \colon G \longrightarrow \GL(V)$ and $v \in V$, $\varphi \in
V'$}]{RepresentativeFunction}{\ensuremath{\pi_{v,\varphi}}}

\glsxtrnewsymbol[description={Adjoint action of a Lie algebra element $\xi \in \liealg{g}$ on $\liealg{g}$}]{AdjointAlgebra}{\ensuremath{\ad_\xi}}

\glsxtrnewsymbol[description={Left invariant vector field corresponding to $\xi \in \liealg{g}$}]{LeftInvariantVectorField}{\ensuremath{X_\xi}}

\glsxtrnewsymbol[description={Left invariant one-form corresponding to $\alpha \in
\liealg{g}^*$}]{LeftInvariantOneForm}{\ensuremath{\theta^\alpha}}

\glsxtrnewsymbol[description={Basis of a Lie
algebra}]{Basis}{\ensuremath{(\basis{e}_1, \ldots, \basis{e}_n)}}

\glsxtrnewsymbol[description={Center of a group $G$}]{Center}{\ensuremath{\Center(G)}}

\glsxtrnewsymbol[description={Category of Lie groups}]{LieCat}{\ensuremath{\categoryname{Lie}}}

\glsxtrnewsymbol[description={Category of complex Lie groups}]{LieCatComplex}{\ensuremath{\categoryname{Lie_\C}}}

\glsxtrnewsymbol[description={Exponential series, Lie
exponential}]{ExponentialSeries}{\ensuremath{\exp}}

\glsxtrnewsymbol[description={Exponential mapping of a Lie group $G$}]{Exponential}{\ensuremath{\exp_G}}

\glsxtrnewsymbol[description={Lie algebra of a Lie group
$G$}]{LieAlgStandalone}{\ensuremath{\liealg{g}}}

\glsxtrnewsymbol[description={Complex Lie algebra of a complex Lie group
$G$}]{LieAlgStandaloneComplex}{\ensuremath{\hat{\liealg{g}}}}

\glsxtrnewsymbol[description={Lie algebra of a Lie group
$G$}]{LieAlg}{\ensuremath{\LieAlg}(G)}

\glsxtrnewsymbol[description={Complexication of a Lie algebra
$\liealg{g}$}]{ComplexifiedLieAlgebra}{\ensuremath{\liealg{g}_\C}}

\glsxtrnewsymbol[description={Lie group with Lie algebra $\liealg{g} =
    \LieAlg(G)$}]{LieGroup}{\ensuremath{G}}

\glsxtrnewsymbol[description={Real form of a complex Lie group
    $G$}]{RealForm}{\ensuremath{G_\R}}

\glsxtrnewsymbol[description={Connected component of the group unit within a Lie group
$G$}]{UnitComponent}{\ensuremath{G_0}}

\glsxtrnewsymbol[description={Lie derivative}]{LieDer}{\ensuremath{\Lie}}

\glsxtrnewsymbol[description={Standard ordered quantization
map}]{StandardOrdering}{\ensuremath{\varrho_\std}}

\glsxtrnewsymbol[description={Standard ordered star
product}]{StarProductStandardOrdered}{\ensuremath{\star_\std}}

\glsxtrnewsymbol[description={Gutt star
    product, Lie algebra star product}]{StarProductGutt}{\ensuremath{\star_\Gutt}}

\glsxtrnewsymbol[description={Taylor series of a function
$\phi$}]{Taylor}{\ensuremath{\Taylor_\phi}}

\glsxtrnewsymbol[description={Formal Lie-Taylor series at $g$ with formal
parameter~$\underline{z}$}]{LieTaylor}{\ensuremath{\Taylor(\underline{z};g)}}

\glsxtrnewsymbol[description={Lie-Taylor majorant of $\phi \in \Comega(G)$ at a group
element~$g \in G$}]{Majorant}{\ensuremath{\Majorant_\phi(\underline{z};g)}}

\glsxtrnewsymbol[description={Lie-Taylor majorant of $\phi \in \Comega(G)$ at the
group unit~$\E$}]{MajorantAtUnit}{\ensuremath{\Majorant_\phi(\underline{z})}}

\glsxtrnewsymbol[description={Coefficients of the Lie-Taylor majorant of $\phi \in \Comega(G)$ at a group element $g \in G$}]{MajorantCoefficients}{\ensuremath{c_k(\phi;g)}}

\glsxtrnewsymbol[description={Coefficients of the Lie-Taylor majorant of $\phi
\in \Comega(G)$ at the group unit
$\E$}]{MajorantAtUnitCoefficients}{\ensuremath{c_k(\phi)}}

\glsxtrnewsymbol[description={Leading symbol of a differential
operator}]{SymbolLeading}{\ensuremath{\sigma_1}}

\glsxtrnewsymbol[description={Entire functions}]{EntireGeneral}{\ensuremath{\Entire}}

\glsxtrnewsymbol[description={Entire functions on a Lie group
$G$}]{Entire}{\ensuremath{\Entire_0(G)}}

\glsxtrnewsymbol[description={Holomorphic polynomial algebra of a complex Lie group
$G$}]{HolomorphicPolynomials}{\ensuremath{\Holomorphic_\Pol(T^* G)}}

\glsxtrnewsymbol[description={Complexification of a vector space $V$, $V_\C = V
\tensor_\R \C$}]{ComplexifiedVectorSpace}{\ensuremath{V_\C}}

\glsxtrnewsymbol[description={Universal complexification $\eta_G \colon G
\longrightarrow G_\C$ of a Lie group $G$}]{Complexification}{\ensuremath{G_\C}}

\glsxtrnewsymbol[description={Universal complexification $\eta_G \colon G
\longrightarrow G_\C$ of a Lie group $G$}]{ComplexificationMorphism}{\ensuremath{\eta}}

\glsxtrnewsymbol[description={Universal complexification
functor}]{ComplexificationFunctor}{\ensuremath{\argument_\C}}

\glsxtrnewsymbol[description={Covering projection}]{Covering}{\ensuremath{p}}

\glsxtrnewsymbol[description={Universal covering group of a Lie group
$G$}]{CoveringGroup}{\ensuremath{\widetilde{G}}}

\glsxtrnewsymbol[description={Universal enveloping algebra of a Lie algebra
$\liealg{g}$}]{UniversalEnvelope}{\ensuremath{\Universal^\bullet(\liealg{g})}}

\glsxtrnewsymbol[description={Smallest closed complex subgroup generated by a subset $S$ of a complex Lie group $G$}]{GroupSpanComplex}{\ensuremath{\left<S\right>_\C}}

\glsxtrnewsymbol[description={Fundamental Group of topological group $G$}]{FundamentalGroup}{\ensuremath{\pi_1(G)}}

\glsxtrnewsymbol[description={Complex conjugated group of complex Lie group
$G$}]{ComplexConjugatedGroup}{\ensuremath{\cc{G}}}

\glsxtrnewsymbol[description={Complex conjugation $\operatorname{cc}(z) = \cc{z}$ on $\liealg{g}_\C$}]{ComplexConjugation}{\ensuremath{\operatorname{cc}}}

\glsxtrnewsymbol[description={Vector space generated by a set over the field $\R$}]{SpanR}{\ensuremath{\Span_\field{R}}}

\glsxtrnewsymbol[description={Vector space generated by a set over the field
$\C$}]{SpanC}{\ensuremath{\Span_\field{C}}}

\glsxtrnewsymbol[description={Gårding space}]{Garding}{\ensuremath{\mathfrak{G}}}

\glsxtrnewsymbol[description={General linear group of a vector space
$V$}]{GeneralLinear}{\ensuremath{\GL(V)}}

\glsxtrnewsymbol[description={General linear group of
$\C^d$}]{GeneralLinearComplex}{\ensuremath{\GL_d(\C)}}

\glsxtrnewsymbol[description={Smooth vectors of a represention $\pi \colon G \longrightarrow \GL(V)$}]{SmoothVectorGroup}{\ensuremath{\Cinfty(\pi)}}

\glsxtrnewsymbol[description={Orbit map of a vector $v \in V$ under a represention
$\pi \colon G \longrightarrow \GL(V)$}]{OrbitMap}{\ensuremath{\pi_v}}

\glsxtrnewsymbol[description={Orbit maps associated to smooth vectors of a represention $\pi \colon G \longrightarrow \GL(V)$}]{SmoothVectorGroupAsFunctions}{\ensuremath{\Cinfty_\pi(G,V)}}

\glsxtrnewsymbol[description={Mapping a vector to its orbit map under a representation
$\pi \colon G \longrightarrow \GL(V)$}]{VectorToOrbitMap}{\ensuremath{\Pi}}

\glsxtrnewsymbol[description={Smooth vectors of a represention $\varrho \colon
\liealg{g} \longrightarrow
\Linear(V)$}]{SmoothVectorAlgebra}{\ensuremath{\Cinfty(\varrho)}}

\glsxtrnewsymbol[description={Analytic vectors of a represention $\varrho \colon
\liealg{g} \longrightarrow
\Linear(V)$}]{AnalyticVectorAlgebra}{\ensuremath{\Comega(\varrho)}}

\glsxtrnewsymbol[description={Analytic vectors of a represention $\pi \colon G
\longrightarrow \GL(V)$}]{AnalyticVectorGroup}{\ensuremath{\Comega(\pi)}}

\glsxtrnewsymbol[description={Entire vectors of a represention $\varrho \colon
    \liealg{g} \longrightarrow
    \Linear(V)$}]{EntireVectorAlgebra}{\ensuremath{\underline{\Entire}(\varrho)}}

\glsxtrnewsymbol[description={Entire vectors of a represention $\pi \colon G
    \longrightarrow \GL(V)$}]{EntireVectorGroup}{\ensuremath{\underline{\Entire}(\pi)}}

\glsxtrnewsymbol[description={Strongly entire vectors of a
represention~$\varrho$}]{StronglyEntireVectorAlgebra}{\ensuremath{\Entire(\varrho)}}

\glsxtrnewsymbol[description={Strongly entire vectors of a
represention~$\pi$}]{StronglyEntireVectorGroup}{\ensuremath{\Entire(\pi)}}

\glsxtrnewsymbol[description={Schwartz space}]{Schwartz}{\ensuremath{\Schwartz}}

\glsxtrnewsymbol[description={Lebesgue measure of a
set}]{measure}{\ensuremath{\operatorname{meas}}} % glossary entries, mention with \gls

\author{Michael Heins}
\title{A Holomorphic perspective \\ of Strict Deformation Quantization}
\date{\today}

\makenoidxglossaries % initiate list of symbols
\makeindex[program=makeindex, options=-s Auxiliary/IndexStyle,
columns=2] %initiate index

\addto\captionsenglish{% Replace "english" with the language you use
    \renewcommand{\contentsname}%
    {Table of Contents}%
}

\begin{document}
\afterpage{\blankpage} %Empty page after titlepage
\begin{titlepage}
    \begin{center}
        \vspace*{-2cm}

        \Huge
        \textbf{A Holomorphic perspective \\ of Strict Deformation Quantization}

        \vspace{.5cm}

        \textbf{Michael Heins}

        \vfill

        \includegraphics[width=0.9\textwidth]{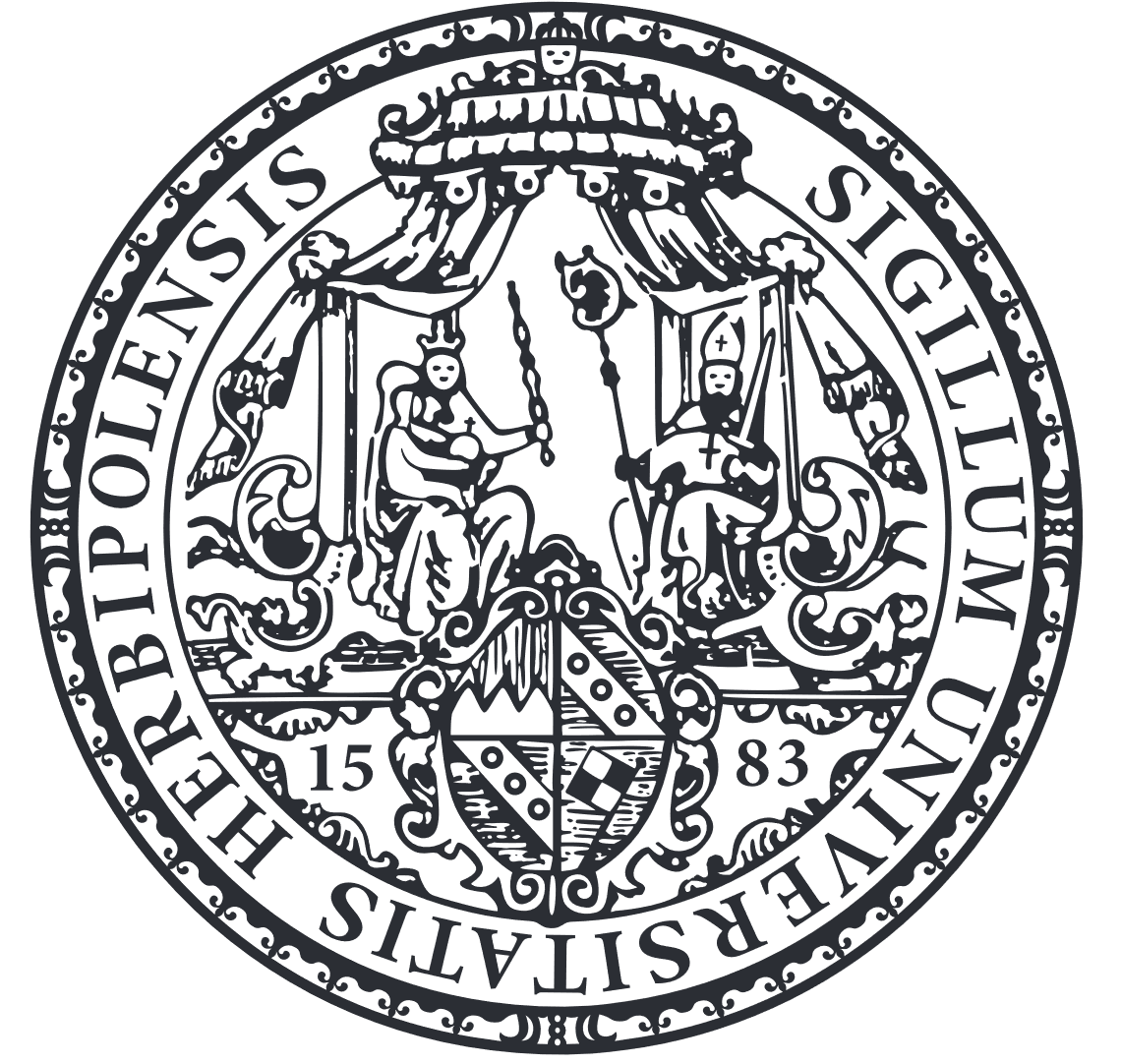}

        Dissertation in Mathematics

        \vspace{0.5cm}

        \large
        Julius-Maximilians-Universität Würzburg \\
        Faculty of Mathematics and Computer Science \\
        \begin{center}
            Doctoral Advisors: \\
             Prof. Dr. Oliver Roth \\
             Prof. Dr. Stefan Waldmann
        \end{center}
        October 2024
    \end{center}
\end{titlepage}

\begin{dedication}
    For Michael Werner Heins,
    \par   %% or a blank line
%    \vspace{2\baselineskip}
    my father, who gave up everything \\
    -- just to be always there for me.
\end{dedication}

\begin{dedication}
    \epigraph{In the Ramtop village where they dance the real Morris dance, for
        example, they believe that no one is finally dead until the ripples they cause
        in
        the world die away -- until the clock he wound up winds down, until the wine
        she
        made has finished its ferment, until the crop they planted is harvested. The
        span
        of someone’s life, they say, is only the core of their actual
        existence.}{\emph{Reaper Man} -- Terry Pratchett}
\end{dedication}

\tableofcontents
\addcontentsline{toc}{chapter}{Table of Contents}
\addtocounter{page}{-3} %Adjusting the counter accordingly

\chapter{Deformation Quantization}
\label{ch:StrictDeformationQuantization}
\epigraph{The trouble with having an open mind of course, is that people insist on
coming along and trying to put things in it.}{\emph{Diggers} -- Terry Pratchett}

% !TeX root = ../Dissertation.tex
This is a treatise concerning the field of \gls{Complex}omplex analysis. Unless
specified otherwise, this means that all vector spaces are complex\footnote{As they
should be.}, the division by $n!$ may be pursued with reckless abandon, linear maps are
complex linear and, perhaps most importantly, functions take complex values. As we
shall see in many places throughout the text, this should not at all be
regarded as a restriction, but rather as the natural setting of our considerations.
Moreover, the central theme of \emph{complexification} will allow us to enhance real
situations to complex ones, giving us access to the versatile toolbox of complex
analysis throughout. This, in turn, simplifies many of our arguments dramatically. In
the timeless\footnote{And often misquoted.} words of Jacques Hadamard\footnote{Jacques
Salomon Hadamard (1865-1963)
was a French
complex analyst, differential geometer and philosopher with a particular interest in
partial differential equations. He pioneered the notion of well-posedness and invented
the method of descent. His doctoral advisors were Èmile Picard and Jules Tannery, and
notably Hadamard supervised in turn Paul Lévy, Maurice René Fréchet and André Weil.}
\cite[p.123]{hadamard:1996a}:

\epigraph{It has been written that the shortest and best way between two truths of
    the real domain often passes through the imaginary one.}{\emph{The mathematician's
    mind}}

That being said, this is also a thesis residing in the wondrous, expansive and
sometimes eldritch realm of
mathematical physics, and ultimately aims for results therein. In the pursuit of the
corresponding questions, we forgo some of the elegance and brevity of the complex
analyst to admire and cherish the view of the mathematical landscape our problems have
lead us into. Indeed, the following is my credo:
\vspace{0.5cm}
\begin{center}
    \textbf{The question dictates the methods we learn, and not the other way
    around.}
\end{center}
\vspace{0.75cm}

Of course, one only really knows in retrospect what these methods precisely
constitute. As such, we let our minds wander within the mathematical maze, we get
lost, pursue both the interesting and uninteresting, converse with the necessary and
the pointless. And in due time, some of the fog residing over the mathematical
landscape starts to lift, a bit of clarity is achieved and some territory is charted.
Invariably, we then ask another question, and dive back into the allure of the
unknown and challenging. Writing this text, I have attempted to conserve the spirit of
this -- at times rather daring and arguably unwise -- journey,
including its more entertaining pitfalls and the spontaneous diversions~I invariably
followed in the pursuit of mastery and butterflies. It is my ambition that the
included material constitutes both a learning experience and provides at least some
mild amusement. Hopefully, this coherently incoherent\footnote{Or is it incoherently
coherent?} rambling explains the almost ungraceful length of the present disquisition.
Surely, it does.

\newpage

This chapter serves an introduction and sets the stage for our later considerations.
In Section~\ref{sec:QuantizationFormal}, we
briefly present a rather naïve incarnation of the quantization problem, leading to the
notion of a formal star
product. Replacing the formal parameter with a complex number $\gls{Planck} \in \C$,
which carries the physical interpretation of Planck's constant, then provides the
problem of convergence of the formal series.

Following the program outlined in \cite{waldmann:2019a}, we then strengthen mere
convergence to continuity in Section~\ref{sec:QuantizationStrict}. This ostensibly
stronger requirement has the advantage that
it may be established on some smaller subalgebra -- where, say, convergence is trivial
by termination of the power series or elements are factorizing in a convenient manner.
Then, in a second step, one may extend the product to the completion of that
subalgebra, and the hope is that this process obtains a wealth of interesting and
non-trivial observables.

As a guiding example, one should think of the passage from
polynomials to power series. The bridge between their global behaviour and their
coefficients consists in the classical \emph{Cauchy\footnote{Augustin-Louis Cauchy
(1789-1857) was an extraordinarily productive French mathematician, who wrote more
than eight hundred research articles. He pioneered the rigorous formulation of real
and complex analysis, and introduced a vast number of the concepts constituting the
basis of modern analysis.} estimates},
which
have become a common theme within strict deformation quantization.
This is the principal philosophy of strict deformation
quantization. In an attempt to include the growing wealth of known examples, we also
propose a general, analytic notion of star product in
Definition~\ref{def:StarProduct}, which captures more general dependencies on
Planck's constant.

Finally, Section~\ref{sec:QuantizationHolomorphic} attempts to take these ideas to
their logical conclusion: One should work with holomorphic objects from the start
instead of the real analytic shadows they cast onto the original geometry.
Indeed, in
all known examples such an approach turned out to be viable, at least a posteriori.
This has the advantage of separating the continuity of the product from the Cauchy
estimates, which are inherent to the geometry. Philosophically, the ideas we present
are close to the ones central to \cite{dito.schapira:2007a, schapira:2008a,
kashiwara.schapira:2012a}, albeit with vastly different practical implementation.

Throughout the text, we attempt to provide some biographical background information
on the scientists we encounter in the form of short footnotes -- lest we forget they were
people like us. As such, the sources used
for this particular endeavour are rather diffuse, and we do not even make the attempt at
the herculean task of listing all them in their entirety. As most of our heroes have since
passed away, obituaries by various sources such as the American Mathematical
Society have however been invaluable.

\newpage

\section{Formal Deformation Quantization}
\epigraph{Light thinks it travels faster than anything but it is wrong. No matter how
    fast light travels, it finds the darkness has always got there first, and is
    waiting for it.
}{\emph{Reaper Man} -- Terry Pratchett}
\label{sec:QuantizationFormal}
% !TeX root = ../Dissertation.tex

The problem of quantization consists in formalizing the passage from a classical
mechanical system, mathematically encoded by its observables as a
Poisson\footnote{Siméon Denis
Poisson (1781-1840) was a French mathematician and physicist. His memoirs on
electricity and magnetism were groundbreaking and many of his ideas still resonate
within modern textbooks. Moreover, he proved several fundamental theorems on Fourier
integrals, the calculus of variations and probability theory.} algebra, to its quantum
analogue. Recall that a Poisson algebra~$(\gls{Algebra}, \{\argument, \argument\})$ is
a
pair of an associative, commutative and unital algebra~$(\algebra{A},+,\cdot)$ over the
complex numbers $\C$ and a Lie\footnote{Marius Sophus Lie (1842-1899) was a Norwegian
mathematician and the father of modern Lie theory, the geometric study of continuous
symmetry groups. His fundamental theorems, two of which we will meet in
Theorem~\ref{thm:Lie2} and Theorem~\ref{thm:Lie3} provided a precise correspondence
between Lie groups -- and their rich global geometry -- with Lie algebras, which
encode the local geometry in an algebraic fashion.} bracket
    \begin{equation}
        \index{Poisson algebra}
        \index{Poisson bracket}
        \label{eq:PoissonBracket}
        \gls{PoissonBracket}
        \colon
        \algebra{A} \tensor \algebra{A}
        \longrightarrow
        \algebra{A}
    \end{equation}
    fulfilling the Leibniz rule
    \begin{equation}
        \{f, gh\}
        =
        g
        \cdot
        \{
        f,h
        \}
        +
        \{
        f,g
        \}
        \cdot
        h
        \qquad
        \textrm{for all }
        f,g,h \in \algebra{A}.
\end{equation}
The principal example of a Poisson algebra is the algebra of smooth functions
\gls{SmoothFunctions} on a symplectic manifold $M$, endowed with the canonical Poisson
bracket
induced by the symplectic form. A lovely overview over the connection between
classical mechanics symplectic geometry by means of Lagrangian submanifolds and
illustrated by examples from celestial mechanics is \cite{weinstein:1981a}. The
corresponding technical details can be found in the textbooks
\cite{abraham.marsden:1985a, arnold:1989a, honerkamp.roemer:1993a,
marsden.ratiu:1999a}. For a comprehensive discussion of the mathematical theory of
symplectic and Poisson geometry we moreover refer to
\cite{mcduff.salamon:1998a,fomenko:1995a,dasilva:2001a},
respectively~\cite{weinstein:1998a,
waldmann:2007a}.

The quantization problem may then be resolved by constructing a \emph{quantization
map}, that is, a linear bijection
\begin{equation}
    \index{Quantization map}
    \label{eq:QuantizationMap}
    \gls{QuantizationMap}
    \colon
    \algebra{A}
    \longrightarrow
    \gls{AlgebraQuantum}
    \subseteq
    \mathfrak{A}(H)
\end{equation}
onto a \emph{subalgebra} $\algebra{A}_{\quant} \subseteq \mathfrak{A}(H)$,
where $H$ is a Hilbert space and \gls{Adjointable} denotes the space of
densely defined linear adjointable maps on $H$. In particular, such maps need
\emph{not} be continuous, unless their domains are all of the ambient Hilbert space
$H$. Its presence is demanded to make the powerful tools of spectral
theory available, as they are wonderfully laid out and explored in
\cite{schmuedgen:2020a}. Indeed, physically, the spectrum of a quantum
observable~$O \in
\algebra{A}_{\quant}$
corresponds to the possible outcomes of its measurement, and as such is an
indispensable structure within any formulation of a reasonable quantum theory.

The bijectivity of \eqref{eq:QuantizationMap} is an ultimately necessary technical
assumption. Indeed, physical reality \emph{is inherently quantum} and consequently
classical mechanics is merely an -- extraordinarily useful and highly successful --
approximation with a limited range of applicability. This means that every classical
observable necessarily possesses a unique quantum analogue. Conversely,
there may however in principle be quantum observables that cannot be observed
classically. Indeed,
historically precisely this occurred with the discovery of quantum mechanical spin.
That being said, in this particular situation, by \cite{marsden:1968a, marsden:1968b},
there exists an enhanced version of Hamiltonian mechanics that allows for a classical
version of spin that quantizes correctly.

The philosophy underlying these ideas puts the precise correspondence between
classical and quantum observables at its core. There is no ambiguity concerning
the quantum versions of positions, momenta and the Hamiltonian.
Moreover, \eqref{eq:QuantizationMap} implements the, physically not at all harmless,
\emph{superposition principle} simply by its linearity. Finally, the
assumption~$\algebra{A}_{\quant} \subseteq \mathfrak{A}(H)$ as a \emph{subalgebra}
endows the quantum observable algebra $\algebra{A}_{\quant}$ with the associative
product $\circ$ of operator composition and the corresponding commutator
\begin{equation}
    [O_1,O_2]
    \coloneqq
    O_1 \gls{Composition} O_2
    -
    O_2 \circ O_1
    \qquad
    \textrm{for all }
    O_1, O_2 \in \mathfrak{A}(H).
\end{equation}
That being said, the resulting problem is still far too general and, conversely,
glosses over a number of crucial details. Indeed, the Poisson
bracket~\eqref{eq:PoissonBracket} has not yet entered at all into our discussion.
Taking a closer look, a natural demand would be that the quantization map $Q$
preserves brackets -- up to some scalar multiple involving Planck's constant $\hbar$
-- which mathematically translates to $Q$ being an isomorphism
\begin{equation}
    Q
    \colon
    \bigl(
        \algebra{A},
        \{\argument,\argument\}
    \bigr)
    \longrightarrow
    \bigl(
        \algebra{A}_{\quant},
        \gls{LieBracket}
    \bigr)
\end{equation}
of Lie algebras. This then would have the pleasant consequence of $Q$ also respecting
the canonical commutation relations. However, these demands turn out to be impossible
to realize already for the Schrödinger quantization of $\R^{2n}$ by the classical
Groenewold\footnote{Hilbrand Johannes Groenewold
    (1910-1996) was a Dutch theoretical quantum physicist. He pioneered the phase
    space formulation of quantum mechanics, which works with functions in $(q,p)$
    instead of the Hilbert space operators. This ingenious idea should sound vaguely
    familiar by now.}-Van Hove\footnote{Léon Charles Prudent Van Hove (1924-1990)
    was a Belgian mathematical physicist, who made contributions to solid-state and
    nuclear physics, the theory of elementary particles as well as cosmology.}
    Theorem \cite{groenewold:1946a, vanhove:1951b,
    vanhove:1951a}.
Another intrinsic problem of this formulation is that the rôle of the Hilbert space in
\eqref{eq:QuantizationMap} is unclear. Indeed, it ultimately corresponds to the choice
of a concrete representation of the quantum system. However, there may well be other,
inequivalent representations required to capture physical reality. A concrete example
of this problem arises within the description of the Aharonov-Bohm effect
\cite{aharonov.bohm:1959a}.

Consequently, one has proceed in a different manner. This has lead to the approaches
of geometric quantization\footnote{A comprehensive account of which is the monograph
\cite{woodhouse:1992a}.}, whose focus lies on the space of states associated to the
observable algebra, and formal deformation quantization, which keeps the observable
algebra as the central piece. The latter was conceived in the seminal paper
\cite{bayen.et.al:1978a}, which defines the abstract notion of a formal star product,
and proves that this approach is consistent with established quantization schemes.
Some of the underlying ideas already appear within Dirac's foundational work
\cite{dirac:1947a} and were subsequently refined by Berezin\footnote{Felix
Alexandrovich Berezin
(1931-1980) was a Russian mathematical physicist. He made various fundamental
contributions to supermathematics, which is based on anticommuting Grassman
variables, and is the mathematical language of supersymmetry.} \cite{berezin:1975a,
berezin:1975b, berezin:1975c}.

What follows is an adaptation of and an overview over
\cite[Sec.~6.1.2]{waldmann:2007a}. We will keep the discussion brief and as
non-technical as possible. We denote the ring of formal power series in the
indeterminate $\lambda$ with coefficients in another ring $\ring{R}$ by
\begin{equation}
    \index{Formal!power series}
    \gls{FormalPowerSeries}
    \coloneqq
    \biggl\{
        \sum_{n=0}^{\infty}
        a_n
%        \cdot
        \lambda^n
        \colon
        a_n \in \ring{R}
        \textrm{ for all }
        n \in \N_0
    \biggr\}.
\end{equation}
Note that, canonically, $\ring{R} \subseteq \ring{R}\formal{\lambda}$ as the
coefficient of $\lambda^0$. Its product is given by the convolution
\begin{equation}
    \index{Formal!Convolution}
    \label{eq:FormalConvolution}
    \biggl(
        \sum_{n=0}^{\infty}
        a_n
        \lambda^n
    \biggr)
    \cdot
    \biggl(
        \sum_{n=0}^{\infty}
        b_n
        \lambda^n
    \biggr)
    \coloneqq
    \sum_{n=0}^{\infty}
    \lambda^n
    \sum_{k=0}^n
    a_k b_{n-k},
\end{equation}
where $(a_n), (b_n) \subseteq \ring{R}$ and we use the ring structure of $\ring{R}$ in
the inner summation. If $\algebra{A}$ is an associative algebra,
then so is $\algebra{A}\formal{\lambda}$. Moreover, if $L \colon \algebra{A}
\longrightarrow \algebra{A}$ is a linear mapping, then setting
\begin{equation}
    \label{eq:FormalExtensionLinear}
    L
    \biggl(
        \sum_{n=0}^{\infty}
        a_n
        \lambda^n
    \biggr)
    \coloneqq
    \sum_{n=0}^{\infty}
    L(a_n)
    \lambda^n
\end{equation}
provides a $\C\formal{\lambda}$-linear extension of $L$ to a self-map of
$\algebra{A}\formal{\lambda}$. Analogously, one may extend any linear mapping $B \colon
\algebra{A} \tensor \algebra{A} \longrightarrow \algebra{A}$ to a
$\C\formal{\lambda}$-linear mapping
\begin{equation}
    B
    \colon
    \bigl(
        \algebra{A}
        \tensor
        \algebra{A}
    \bigr)
    \formal{\lambda}
    \longrightarrow
    \algebra{A}\formal{\lambda}
\end{equation}
by setting
\begin{equation}
    \label{eq:FormalExtensionBilinearGeneral}
    B
    \biggl(
        \sum_{n=0}^{\infty}
        \lambda^n
        \sum_k
        a_k^{(n)}
        \tensor
        b_k^{(n)}
    \biggr)
    \coloneqq
    \sum_{n=0}^{\infty}
    \lambda^n
    \sum_k
    B
    \bigl(
        a_k^{(n)}
        \tensor
        b_k^{(n)}
    \bigr).
\end{equation}
Note that the formal convolution product \eqref{eq:FormalConvolution} provides a
canonical mapping
\begin{equation}
    \algebra{A}\formal{\lambda}
    \tensor_{\C\formal{\lambda}}
    \algebra{A}\formal{\lambda}
    \longrightarrow
    \bigl(
    \algebra{A}
    \tensor
    \algebra{A}
    \bigr)
    \formal{\lambda},
\end{equation}
which fails to be surjective in general. Suppressing this mapping,
\eqref{eq:FormalExtensionBilinearGeneral} then takes the form
\begin{equation}
    \label{eq:FormalExtensionBilinear}
    B
    \biggl(
        \sum_{n=0}^{\infty}
        a_n
        \lambda^n
        \tensor
        \sum_{k=0}^{\infty}
        b_k
        \lambda^k
    \biggr)
    \coloneqq
    \sum_{n=0}^{\infty}
    \lambda^n
    \sum_{k=0}^n
    B
    \bigl(
        a_k \tensor b_{n-k}
    \bigr).
\end{equation}
Conceptually, one should work with $(\algebra{A} \tensor \algebra{A})
\formal{\lambda}$ instead of the more naive $\algebra{A}\formal{\lambda}
\tensor_{\C\formal{\lambda}} \algebra{A}\formal{\lambda}$, as the former is complete
with respect to the $\lambda$-adic topology, whereas the other is typically not. In
passing, we also note that all our extensions are continuous in this setting, and
refer to \cite[Sec.~6.2.1]{waldmann:2007a} for further topological discussion.

Moreover, it is a straightforward exercise to prove
that~\eqref{eq:FormalExtensionLinear} and \eqref{eq:FormalExtensionBilinearGeneral}
are already the general forms of $\C\formal{\lambda}$-linear and bilinear mappings, a
precise formulation of which can be found in \cite[Prop.~2.1]{dewilde.lecomte:1983a}.
Having established these preliminaries, the definition of a formal star product is the
following.
\begin{definition}[Formal star product, {\cite[Def.~21]{bayen.et.al:1978a}}]
    \label{def:FormalStarProduct}
    \index{Formal!star product}
    \index{Star product!Formal}
    \index{Deformation quantization!Formal}
    \index{Formal!Deformation quantization}
    Let $(\algebra{A}, \{\argument, \argument\})$ be a Poisson algebra. A formal star
    product is a $\C\formal{\lambda}$-linear associative product
    \begin{equation}
        \gls{StarProduct}
        \colon
        \algebra{A}\formal{\lambda}
        \tensor
        \algebra{A}\formal{\lambda}
        \longrightarrow
        \algebra{A}\formal{\lambda}
    \end{equation}
    fulfilling the following conditions:
    \begin{definitionlist}
        \item There are linear mappings
        \begin{equation}
            \label{eq:StarProductDiffops}
            C_n
            \colon
            \algebra{A} \tensor \algebra{A}
            \longrightarrow
            \algebra{A}
            \qquad
            \textrm{for all }
            n \in \N_0
        \end{equation}
        such that
        \begin{equation}
            \label{eq:StarProductFormal}
            a \star b
            =
            \sum_{n=0}^{\infty}
            C_n(a,b)
            \lambda^n
            \qquad
            \textrm{holds for all }
            a,b \in \algebra{A}\formal{\lambda},
        \end{equation}
        where we extend the mappings $C_n$ to $\algebra{A}\formal{\lambda} \tensor
        \algebra{A}\formal{\lambda}$ as in \eqref{eq:FormalExtensionBilinear}.
        \item Let $f,g \in \algebra{A}$. Then
        \begin{equation}
            \index{Classical limit!Formal}
            \label{eq:ClassicalLimit}
            C_0(f,g)
            =
            fg
        \end{equation}
        and
        \begin{equation}
            \index{Semiclassical limit!Formal}
            \label{eq:SemiclassicalLimit}
            C_1(f,g)
            -
            C_1(g,f)
            =
            \bigl\{
                f,
                g
            \bigr\}.
        \end{equation}
        \item The unit $1 \in \algebra{A} \subseteq \algebra{A}\formal{\lambda}$ acts
        as neutral element for $\star$, i.e.
        \begin{equation}
            \label{eq:StarProductUnits}
            1 \star a
            =
            a
            =
            a \star 1
            \qquad
            \textrm{for all }
            a \in \algebra{A}.
        \end{equation}
    \end{definitionlist}
    In this case, we call the pair $(\algebra{A}\formal{\lambda}, \star)$ a formal
    deformation quantization of $(\algebra{A}, \{\argument, \argument\})$.
\end{definition}

Notably, \eqref{eq:StarProductUnits} is both easy to achieve and to break by means of
equivalence transformations once one has obtained a formal star product at all.
Demanding the equality of classical and quantum unit is however fundamental from a
physical point of view, as they correspond to the empty measurement. That is to say,
not asking a question at all, which in turn has to correspond to receiving no answer
either.

Typically, the mappings
\eqref{eq:StarProductDiffops} come from bidifferential operators of order $(n,n)$, a
notion we shall make precise in Remark~\ref{rem:DiffOps}. The condition
\eqref{eq:ClassicalLimit} is
known as the \emph{classical limit} and \eqref{eq:SemiclassicalLimit} as the
\emph{semiclassical} limit. Algebraically, one thus also speaks of a deformation of
the product in direction of the Poisson bracket, an idea that may be
formalized by means of Gerstenhaber's\footnote{Murray Gerstenhaber (1927-2024) was an
American mathematician, physicist and, notably, lawyer. He earned his law degree
during his first sabbatical as a professor and later lectured on the statistics of law.
His seminal contributions to algebraic deformation theory have since found
applications and generalizations in the context of operads, deformation quantization
and also of quantum groups.} theory
\cite{gerstenhaber:1964a,gerstenhaber:1966a,gerstenhaber:1968a,gerstenhaber:1974a} of
algebraic deformations.

The tricky part for the construction of formal star products is somewhat hidden within
Definition~\ref{def:FormalStarProduct}:
It is the \emph{associativity} of the product. Indeed, all other conditions are linear
with respect to the operators $C_n$, whereas associativity has a quadratic dependence.
Concretely, this yields the system of equations
\begin{equation}
    \label{eq:Associativity}
    \sum_{k=0}^n
    C_k
    \bigl(
        C_{n-k}(a,b), c
    \bigr)
    =
    \sum_{k=0}^n
    C_k
    \bigl(
        a, C_{n-k}(b,c)
    \bigr)
    \qquad
    \textrm{for all }
    a,b,c \in \algebra{A}\formal{\lambda}, \,
    n \in \N_0.
\end{equation}
As \eqref{eq:ClassicalLimit} and \eqref{eq:SemiclassicalLimit} fix
$C_0$ and the antisymmetric part of $C_1$, one may then in principle attempt to solve
\eqref{eq:Associativity} inductively. Indeed, at a very informal level and
disregarding many substantial difficulties requiring additional technology, this idea
is at the heart of Kontsevich's\footnote{Maxim Lwowitsch Kontsevich (born 1964) is a
French-Russian differential geometer also interested in rigorous formulations of
quantum field theory. In the year 1998, he received the Fields Medal for his
contributions to ``algebraic geometry, topology, and mathematical physics''.}
and also Fedosov's construction.\footnote{Boris Vasil’evich Fedosov (1938-2011) was a
Russian
differential geometer. Already early on in his career, he established a theory of
formal symbols and algebraic indices for pseudodifferential operators as an
alternative approach to index theory in the sense of Atiyah-Singer. His ideas however
lacked the modern notion of cyclic cohomology. Building on his symbol calculus
later lead him to his fundamental construction within formal deformation
quantization, which is the content of Remark~\ref{rem:Fedosov}.}

Returning to our somewhat unrefined ideas of quantization maps, we may
provide a toy example. We shall return to it twice to illustrate and exemplify our
general philosophy.
\begin{example}[Standard ordered star product I: Formal,
{\cite[Sec.~5.2.2]{waldmann:2007a}}]
    \index{Standard ordered!Star product}
    \index{Formal!Standard ordered star product}
    \index{Star product!Standard ordered formal}
    \label{ex:StdOrdI}
    \; \\ Consider
    \begin{equation}
        \algebra{A}
        \coloneqq
        \gls{PolynomialsVectorBundle}(\gls{Cotangent}\R)
        \coloneqq
        \Cinfty(\R)
        \tensor
        \R[p]
        \cong
        \Cinfty(\R)[p],
    \end{equation}
    where we model the cotangent bundle $T^*\R$ as $\R^2$ by means of the global
    coordinates $(q,p)$ and endow the polynomial algebra $\Pol(T^* \R)$ with the
    standard Poisson bracket
    \begin{equation}
        \{
            f, g
        \}
        \coloneqq
        \frac{\partial f}{\partial q}
        \cdot
        \frac{\partial g}{\partial p}
        -
        \frac{\partial f}{\partial p}
        \cdot
        \frac{\partial g}{\partial q}
        \qquad
        \textrm{for all }
        f,g \in \Pol(T^* \R).
    \end{equation}
    The physical interpretation of $(q,p)$ are the position and the corresponding
    momentum of a particle confined to the real axis $\R$. Associating a polynomial
    with the -- typically non constant -- differential operator
    \begin{equation}
        \label{eq:QuantizationStd}
        Q(\phi \cdot p^n)
        \coloneqq
        \lambda^n
        \cdot
        \phi
        \cdot
        \frac{\D^n}{\D x^n}
    \end{equation}
    provides a quantization map, as $Q$ is clearly injective. Note that we have
    made the choice of an \emph{ordering} here, namely to put the derivatives all
    the way to the right-hand side to avoid having them act on their own
    coefficient functions. We refer to this as the \emph{standard ordering}, and
    to $Q$ as the standard ordered quantization map.

    \index{Weyl, Hermann}
    \index{Groenewold, Hilbrand Johannes}
    \index{Moyal, José Enrique}
    This choice results in comparatively simple formulas, but there are of course
    other options. For instance, one could have put all derivatives to the left
    instead, or symmetrized over all possible orders. The latter would have resulted
    in the so-called Weyl ordering,\footnote{Hermann Klaus Hugo Weyl (1885-1955) was a
    German mathematical physicist, widely regarded as one of the last polymaths of
    both disciplines. His work covered the realms of general relativity, matter
    itself, point set topology, group and representation theory, harmonic analysis as
    well as number theory, but also the very foundations of mathematics.} which is
    notably compatible with the $^*$-involution given by complex conjugation and
    already appeared in \cite{weyl:1927a, weyl:1931a} and expanded on within
    \cite{wigner:1932a},
    albeit without a precise specification of
    the involved function spaces. This, and various other technical subtleties, were
    later rectified independently within \cite{groenewold:1946a} and \cite{moyal:1949a}.
    Consequently, one should really speak of the
    Weyl-Wigner\footnote{Eugene Paul Wigner (1902-1995) was a Hungarian-American
    mathematical physicist. In 1930, he was recruited by Princeton and emigrated from
    Germany to the United States of America, where he remained until his death. During
    the second world war, he participated in the Manhattan project, and he remained
    politically active afterwards. His wonderful essay \cite{wigner:1960a} is still well worth
    a read.}-Groenewold-Moyal\footnote{José Enrique Moyal (1910-1998) was an
        Australian mathematical physicist, who connected quantum mechanics with
        classical statistical mechanics in \cite{moyal:1949a}. This facilitated the
        phase space formulation of quantum mechanics. Moreover, he was a pioneer in
        the field of stochastic processes.} ordering.\footnote{The inclined reader is
        invited to
    say this out aloud five times in a row. Quickly.}

    All such ordering prescriptions -- of which there is actually a continuum -- agree
    on the generators~$p, q$ and consequently respect the canonical commutation
    relations
    \begin{equation}
        \label{eq:CanonicalCommutation}
        \lambda
        \cdot
        Q
        \bigl(
            \{q,p\}
        \bigr)
        =
        \lambda \cdot Q(1)
        =
        \lambda
        =
        -
        \biggl[
            x, \lambda \frac{\D}{\D x}
        \biggr]
        =
        -
        \bigl[
            Q(q), Q(p)
        \bigr],
    \end{equation}
    and thus neither choice is canonical. We will return to this observation in
    Remark~\ref{rem:Ordering}.

    The operators \eqref{eq:QuantizationStd} indeed constitute
    adjointable linear operators on the space of smooth functions with compact
    support \gls{TestFunctionsSmooth}$(\R) \subseteq \Ltwo(\R)$, a situation
    which we will come back to in Example~\ref{ex:Schroedinger}.\footnote{Where
    we notably also provide the missing technical details as a by-product of our
    considerations.} Notably, the composition fulfils
    \begin{align}
        Q(\phi \cdot p^n)
        \circ
        Q(\psi \cdot p^m)
        &=
        \lambda^{n+m}
        \cdot
        \phi
        \cdot
        \frac{\D^n}{\D x^n}
        \biggl(
            \psi
            \cdot
            \frac{\D^m}{\D x^m}
        \biggr) \\
        &=
        \lambda^{n+m}
        \cdot
        \phi
        \cdot
        \sum_{k=0}^{n}
        \gls{Binomial}
        \frac{\D^k \psi}{\D x^k}
        \frac{\D^{n+m-k}}{\D x^{n+m-k}} \\
        &=
        \sum_{k=0}^{n}
        \lambda^k
        \binom{n}{k}
        \biggl(
            \phi
            \cdot
            \frac{\D^k \psi}{\D x^k}
        \biggr)
        \cdot
        \biggl(
            \lambda^{n+m-k}
            \frac{\D^{n+m-k}}{\D x^{n+m-k}}
        \biggr) \\
        &=
        \sum_{k=0}^{n}
        \lambda^k
        \cdot
        Q
        \biggl(
            \binom{n}{k}
            \phi
            \cdot
            \frac{\D^k \psi}{\D x^k}
            \cdot
            p^{n+m-k}
        \biggr)
    \end{align}
    for $\phi,\psi \in \Cinfty(\R)$ and $n,m \in \N_0$. Thus, the image of $Q$ is
    closed under composition,
    and we may define another \emph{associative} product on $\Pol(T^*\R)$ by pulling
    back the composition of differential operators. That is, we define
    \begin{equation}
        f \star g
        \coloneqq
        Q^{-1}
        \bigl(
            Q(f) \circ Q(g)
        \bigr)
        \qquad
        \textrm{for all }
        f,g \in \Pol(T^*\R).
    \end{equation}
    The associativity of $\star$ follows immediately from the associativity
    of function composition and we get to avoid dealing with \eqref{eq:Associativity}.
    Our computation provides the explicit formula
    \begin{equation}
        (\phi \cdot p^n)
        \star
        (\psi \cdot p^m)
        =
        \sum_{k=0}^{n}
        \lambda^k
        \cdot
        \binom{n}{k}
        \cdot
        \phi
        \cdot
        \frac{\D^k \psi}{\D q^k}
        \cdot
        p^{n+m-k}
    \end{equation}
    for $\phi,\psi \in \Cinfty(\R)$ and $n,m \in \N_0$. Setting $f(q,p) \coloneqq \phi(q)
    \cdot p^n$ and $g(q,p) \coloneqq \psi(q) \cdot p^m$ this results in the aesthetically
    pleasing expression
    \begin{equation}
        \label{eq:StdExplicit}
        f \star g
        =
        \sum_{k=0}^{\infty}
        \frac{\lambda^k}{k!}
        \cdot
        \frac{\partial^k f}{\partial p^k}
        \cdot
        \frac{\partial^k g}{\partial q^k},
    \end{equation}
    where the series terminates after finitely many terms and is thus only formally
    infinite. Consequently, one may treat the formal parameter $\lambda$ as a complex
    number for the moment. By
    linearity of $Q$, \eqref{eq:StdExplicit} holds for all $f,g \in
    \Pol(T^*\R)$. Taking a closer look, we see that~\eqref{eq:StdExplicit} fulfils
    all of \eqref{eq:ClassicalLimit}, \eqref{eq:SemiclassicalLimit} and
    \eqref{eq:StarProductUnits}. Moreover, the operators \eqref{eq:StarProductDiffops}
    are given by
    \begin{equation}
        C_n
        \coloneqq
        \frac{1}{n!}
        \frac{\partial^n}{\partial p^n}
        \gls{Tensor}
        \frac{\partial^n}{\partial q^n}
        \qquad
        \textrm{for all }
        n \in \N_0,
    \end{equation}
    which are indeed bidifferential of order $(n,n)$ and acting on $\Cinfty(\R)$.
    Finally, we notice that
    \begin{equation}
        \label{eq:Std}
        \phi
        \star
        g
        =
        \phi
        \cdot
        g
        \qquad
        \textrm{for all }
        \phi \in \Cinfty(\R)
        \textrm{ and }
        g \in \Pol(T^*\R),
    \end{equation}
    which means that $\star$ is left $\Cinfty(\R)$-linear. We call $\star$ the
    standard ordered star product. Plugging $f,g \in \Cinfty(T^*\R)$ into
    \eqref{eq:StdExplicit} results in a formal version of the standard ordered star
    product. Unlike for the polynomial algebra $\Pol(T^*\R)$, this is now really only
    a formal power series. Indeed, by the classical Borel Lemma as discussed in
    \cite[§1.5]{narasimhan:1985a},
    every sequence of complex numbers arises as the Taylor data of some smooth
    function, making the
    convergence of \eqref{eq:StdExplicit} on all of $\Cinfty(T^*\R)$ hopeless.
    This necessitates the passage to a proper subalgebra to achieve convergence, an
    endeavour which we shall pursue in Example~\ref{ex:StdOrdII}.
\end{example}

Having established the archetypical example of a formal star product, we turn towards
the other extreme of the story, namely the general existence of formal deformations in
the setting of Poisson manifolds.
\begin{theorem}[Kontsevich's Formality Theorem, {\cite{kontsevich:2003a}}]
    \index{Kontsevich, Maxim}
    \index{Formality theorem}
    \label{thm:Formality}
    Let $M$ be a Poisson manifold.\footnote{That is, a smooth manifold $M$ such that
    its algebra of smooth functions $\Cinfty(M)$ carries the structure of a Poisson
    algebra.} Then
    there exists a formal
    star product quantizing $(\Cinfty(M), \{\argument, \argument\})$.
\end{theorem}
\begin{proof}[Sketchy sketch]
    First, one tackles the local problem, that is to say the manifold $M = \R^n$.
    This is where the main difficulty lies. The principal
    idea is to make an ansatz regarding the~$C_k$ from \eqref{eq:StarProductDiffops}
    as a linear combination of \emph{all possible}
    bidifferential operators of order~$(k,k)$ acting on $\Cinfty(\R^n)$. Encoding these
    bidifferential operators
    as graphs, namely a particular type of quiver, Kontsevich then realized the coefficients
    in the linear combination as certain integrals over these graphs. Astonishingly,
    the coefficients split into a universal and a combinatorial part. The latter
    arises from the graphs and the Poisson tensor in a simple combinatorial manner akin
    to translating Feynman\footnote{Richard Feynman (1918-1988) was an American
    theoretical physicist. His path integral formulation of quantum field theory provides a
    graphical calculus, which organizes the -- ostensibly divergent -- perturbation
    expansions in a systematic fashion. This kind of perspective was characteristic of his
    general outlook: ``If you cannot explain something in simple terms, you don't
    understand it.''}
    diagrams to
    formal integrals.\footnote{Indeed, the graphs
    Kontsevich considers can be understood as the Feynman graphs associated to a
    Poisson sigma model \cite{schaller.strobl:1994a} on a disk by
    \cite{cattaneo.felder:2000a}, which is formalized
    as a quantum field theory by means of path integrals.} The former is then
    independent of the Poisson tensor and may be expressed as integer values of
    multiple zeta functions by \cite{banks.panzer.pym:2020a}. Having established the
    local result, one then globalizes by gluing between charts, a precise account of
    which is \cite{cattaneo.felder.tomassini:2002b}.
\end{proof}

Notably, Kontsevich's result also includes a classification of star products up to
their equivalence. As we shall not need this notion, we have refrained from
introducing it. A modern incarnation of these classification results can be found in
\cite{dolgushev:2005a,dolgushev:2006a,dolgushev:2011a}. This concludes our brief
overview
of formal deformation quantization. For comprehensive -- albeit in some aspects
somewhat dated -- discussions, we refer to the excellent surveys
\cite{weinstein:1995a, sternheimer:1998a, bordemann:2008a} and monographs
\cite{fedosov:1996a, waldmann:2007a}.

\section{Strict Deformation Quantization}
\epigraph{In the beginning there was nothing, which exploded.}{\emph{Lords and
Ladies} -- Terry Pratchett}
\label{sec:QuantizationStrict}
% !TeX root = ../Dissertation.tex

Our discussion, apart from Definition~\ref{def:StarProduct}, loosely follows the
survey \cite{waldmann:2019a}. The principal purpose
of strict deformation quantization is to overcome the formal character of
Definition~\ref{def:FormalStarProduct}. Indeed, to even have the faintest
hope of returning to physically meaningful applications, one has to replace the formal
parameter $\lambda$ with Planck's\footnote{Max Karl Ernst Ludwig Planck (1858-1947)
was a German theoretical physicist, who introduced his constant $h = 2\pi \hbar$ in
his fundamental considerations of energy quanta in the context of black-body radiation.
Not yet fully aware of its fundamental significance, he called it $h$ for the German
\emph{Hilfskonstante} (auxiliary constant), as it played the role of a proportionality
constant. In 1918, he received the Nobel Prize in Physics for his discoveries.}
constant \gls{Planck} in a specified
system of units and to ask for the convergence of the resulting power series.

Taking a closer look,~$\hbar$ carries the unit of an \emph{action}. Going back to
the power series \eqref{eq:StarProductFormal}, this means that the coefficients
$C_n(a,b)$ necessarily carry the unit $[\text{action}]^{-n}$ to make the summation
meaningful. This may be achieved by having $n$~derivatives with respect to position
and~$n$~derivatives with respect to momentum, as derivatives with respect to some
quantity may be thought of carrying inverse units. Indeed, if $f \colon \C
\longrightarrow \C$ is some differentiable function with respect to $z \in \C$, then
its derivative is given by a limit over difference quotients, which implies
\begin{equation}
    \biggl[
        \frac{\D f}{\D z}
    \biggr]
    =
    [f']
    =
    \frac{[f]}{[z]}.
\end{equation}
Formally solving for the differentiation operator gives
\begin{equation}
    \biggl[
        \frac{\D}{\D z}
    \biggr]
    =
    \frac{[f']}{[f]}
    =
    \frac{1}{[z]}.
\end{equation}

In Example~\ref{ex:StdOrdI} we have seen that the classical Borel Lemma implies that
there are always smooth functions $f,g \in \Cinfty(T^* \R)$ such that the power series
\eqref{eq:StdExplicit} has radius of convergence equal to zero. This turns out to be
a typical feature of concrete examples. In most constructions of formal star products
-- including Kontsevich's formality Theorem~\cite{kontsevich:2003a} -- the mappings
$C_n$ from
\eqref{eq:StarProductDiffops} can be chosen as bidifferential operators of order
$(n,n)$ for all $n \in \N_0$, which matches with the preceding dimensional
analysis. For a proper discussion of different types of formal star products and their
classification, we refer to \cite[Sec.~6.1.2 \& 6.1.3]{waldmann:2007a}.

There are at least two feasible ways to resolve the issue of convergence. One is to
view the formal power series \eqref{eq:StarProductFormal} as an \emph{asymptotic
expansion} of a $\hbar$-dependent product $\circ$ in the limit~$\hbar
\gls{LimitFromAbove} 0$. This
essentially means that the difference
\begin{equation}
    \index{Asymptotic expansion}
    a \circ b
    -
    \sum_{n=0}^{N}
    \hbar^n
    \cdot
    C_n(a,b)
\end{equation}
is of order $\gls{Order}(\hbar^{N+1})$ in the limit $\hbar \downarrow 0$ for all $N \in
\N_0$ with respect to some suitable system of seminorms on~$\algebra{A}$.
Historically, this is actually how many of the early
formal star products -- such as the Kähler\footnote{Erich Kähler (1906-2000) was a
German mathematical physicist, who worked on celestial mechanics and the
algebraic incarnation of complex differential geometry. This made him a
progenitor within algebraic geometry, even though his ideas on what would become
schemes were never widely adopted.} type
quantizations
\cite{berezin:1975a, berezin:1975b, berezin:1975c, rawnsley:1977a, rawnsley:1978a,
bordemann.meinrenken.schlichenmaier:1991a,
cahen.gutt.rawnsley:1990a, cahen.gutt.rawnsley:1993a, cahen.gutt.rawnsley:1994a,
cahen.gutt.rawnsley:1995a, cahen.gutt.rawnsley:1996a}, which are based on
Toeplitz\footnote{Otto Toeplitz (1881-1940) was a German functional analyst. In 1911,
he posed the inscribed square problem, which asks whether every Jordan curve
contains an inscribed square. Many particular cases have since been solved, but the
general case remains open, see also the survey \cite{matschke:2014a}.} operators
acting on Bergman\footnote{Stefan Bergman
(1895-1977) was a Polish functional and complex analyst, who after fleeing from the
Nazis ultimately became an American citizen. Using operator theoretic methods,
including his Bergman kernel, he solved problems in complex analysis by means of
operator theoretic techniques, and derived various fundamental integral formulas for
holomorphic functions.} spaces,
see also the more
recent monograph \cite{ma.marinescu:2007a}, -- were constructed: As asymptotic
expansions of integral formulas.

Refining these ideas has lead to Rieffel's\footnote{Marc
Rieffel (born 1937) is a Professor at the department of mathematics at the University of
California, Berkeley. He has made seminal contributions to the theory of
C$^*$-algebras, quantum
groups and pioneered early noncommutative geometry.} theory of
deformations associated to the Lie group $\R^n$ initiated in his seminal paper
\cite{rieffel:1993a}. A comprehensive discussion of this can be found in the
monograph \cite{landsman:1998a}, which also covers aspects of strict deformation
quantization. More recently, advances in replacing $\R^n$ with
general Lie groups have lead to \cite{bieliavsky.massar:2001b, bieliavsky:2002a,
bieliavsky.gayral:2015a}.

In \cite[Thm.~6.8]{heins.moucha.roth.sugawa:2024a}, we
proved that the formal star product explicitly constructed in
\cite{bordemann.brischle.emmrich.waldmann:1996a,
bordemann.brischle.emmrich.waldmann:1996b}, which relied on the foundations laid in
\cite{moreno.ortega-navarro:1983a, moreno.ortega-navarro:1983b, moreno:1986a,
moreno:1987a}, provides an asymptotic expansion for the
Wick star product \cite[(6.13)]{heins.moucha.roth.sugawa:2024a}. Its real analytic
versions
are the content of \cite{esposito.schmitt.waldmann:2019a,
kraus.roth.schoetz.waldmann:2019a, schmitt.schoetz:2022a}. The Wick star product is
given by a factorial series rather than a power series with respect to~$\hbar$ and
exhibits
singularities at $\hbar = -1/n$ for all $n \in \N$ and depends holomorphically on
$\hbar$ otherwise. In particular, the corresponding
power series expansion around the origin has radius of convergence equal to zero.

Another approach for dealing with convergence is due to \cite{waldmann:2019a} and at
the heart of \cite{beiser.roemer.waldmann:2007a, beiser.waldmann:2014a, waldmann:2014a,
esposito.stapor.waldmann:2017a, schoetz.waldmann:2018a,
esposito.schmitt.waldmann:2019a,
kraus.roth.schoetz.waldmann:2019a, schmitt:2021a,
barmeier.schmitt:2022a, heins.roth.waldmann:2023a}. The principal concession is to
put restrictions on the allowed observables. That is, one passes to a proper
subalgebra of the observable algebra. In Example~\ref{ex:StdOrdI}, we have seen that
the formal star product $f \star g$ from \eqref{eq:StdExplicit} converges for all
polynomial functions $f,g \in \Pol(T^*\R)$ irrespective of the topology one considers,
as the series simply terminates after finitely many terms. In concrete examples with
sufficiently explicit formulas, one often has good candidates $\algebra{A}_{\conv}
\subseteq \algebra{A}$ for such subalgebras -- usually some flavour of polynomials.
However, these subalgebras tend to be too small to be truly interesting or useful. In
the particular case of quantum dynamics on cotangent bundles such as $T^*\R$, one
would for instance like to allow also observables exhibiting Gaussian decay with
respect to the momenta.

To remedy this, one then goes one step further and endows \gls{AlgebraConvergent}
with a locally convex topology and asks for the \emph{continuity} of the star product as a
bilinear mapping
\begin{equation}
    \label{eq:StarProductConv}
    \star
    \colon
    \algebra{A}_{\conv}
    \times
    \algebra{A}_{\conv}
    \longrightarrow
    \algebra{A}_{\conv}.
\end{equation}
Taking a closer look, it would suffice to endow
$\algebra{A}_{\conv}$ with the structure of a topological algebra. However, local
convexity has the pleasant consequence that many questions may be resolved by
means of seminorm estimates, making it an extraordinarily useful technical assumption.
The principal problem here is, of course, that it is typically rather difficult to
determine
suitable topologies and there may well be several inequivalent choices available.

Nevertheless, having established the continuity of \eqref{eq:StarProductConv} by some
educated guess, we may pass to the completion $\widehat{\algebra{A}}_{\conv}$ of
$\algebra{A}_{\conv}$: As a continuous bilinear mapping, the product is locally
uniformly continuous and thus has a unique locally uniformly continuous extension
\begin{equation}
    \star
    \colon
    \widehat{\algebra{A}}_{\conv}
    \times
    \widehat{\algebra{A}}_{\conv}
    \longrightarrow
    \widehat{\algebra{A}}_{\conv}.
\end{equation}
Note that we have implicitly fixed a value of $\hbar$ in the preceding discussion.
This was, as it turns out, unreasonable. Taking another look at
\eqref{eq:StarProductFormal}, formal star products are given by vector-valued power
series. Convergence of $a \star b$ for some $\hbar \in \C \setminus \{0\}$ thus
implies the convergence of $a \star b$ on the entirety of the open disk with radius $r
= \gls{AbsoluteValue}$.

Consequently, it is natural to treat $\hbar$ as a \emph{complex} parameter, and to try
to find a common domain of convergence for all $a,b \in \algebra{A}_{\conv}$. From a
physical
vantage point, one should think of this version of $\hbar$ no longer as a fixed
constant, but rather the ratio of $\hbar$ with the \emph{typical action} of the system
one is considering. This fraction may then, in principle, take arbitrary positive
values and at this point, complex values for $\hbar$ come along for free. Having a
common domain of convergence along the algebra then simply translates to the demand
that our model can be applied on the chosen scale at all.

Moreover,
this interpretation gives the limiting procedure $\hbar \downarrow 0$ actual meaning.
It studies the behaviour of the system when $\hbar$ is negligible compared to the
typical actions occurring within the system. This is precisely the situation, in which
one expects behaviour accurately described by the classical observable algebra with
its product and Poisson bracket.

In order to subsume both of the above approaches and, in particular, the observations
from already established examples, we propose the following definition.
\begin{definition}[Star product]
    \label{def:StarProduct}
    \index{Star product}
    \index{Deformation quantization!Strict}
    \index{Deformation domain}
    Let $(\algebra{A}, \{\argument, \argument\})$ be a Poisson algebra endowed
    with a Fréchet topology\footnote{A complete Hausdorff locally convex topology
    induced by a countable system of seminorms.} such that the product \gls{Product}
    and the
    Poisson bracket $\{\argument, \argument\}$ are continuous. A star product is a
    mapping
    \begin{equation}
        \star
        \colon
        \gls{DeformationDomain}
        \times
        \algebra{A}
        \tensor
        \algebra{A}
        \longrightarrow
        \algebra{A},
    \end{equation}
    where $\mathfrak{D} \subseteq \C$ is a domain with accumulation point $0$,
    fulfilling the following conditions:
    \begin{definitionlist}
        \item For every $\hbar \in \mathfrak{D}$, the restriction
        \begin{equation}
            \label{eq:StarProductFixedHbar}
            \star_\hbar
            \coloneqq
            \star
            \at[\Big]{\{\hbar\} \times \algebra{A} \tensor \algebra{A}}
            \colon
            \algebra{A}
            \tensor
            \algebra{A}
            \longrightarrow
            \algebra{A}
        \end{equation}
        endows $\algebra{A}$ with the structure of a Fréchet algebra.\footnote{That
        is to say, \eqref{eq:StarProductFixedHbar} constitutes an associative and
        continuous
        product on $\algebra{A}$.}
        \item The units of $(\algebra{A},\cdot)$ and $(\algebra{A},\star_\hbar)$ coincide for
        all $\hbar \in \mathfrak{D}$.
        \item The mapping
        \begin{equation}
            \label{eq:StarProductHolomorphy}
            \mathfrak{D}
            \ni
            \hbar
            \; \mapsto \;
            f \star_\hbar g
            \in
            \algebra{A}
        \end{equation}
        is holomorphic for all
        $f,g \in \algebra{A}$.\footnote{One may demand either Gâteaux holomorphy, which
        is the content of Section~\ref{sec:HolomorphicGateaux}, or Fréchet holomorphy,
        which we shall study in Section~\ref{sec:HolomorphicFrechet}. In
        Proposition~\ref{prop:HartogRedux} we will prove that the resulting notion of
        holomorphy is the same for functions with finite dimensional domains such as
        \eqref{eq:StarProductHolomorphy}.}
        \item For all $f,g \in \algebra{A}$, we have
        \begin{equation}
            \index{Classical limit}
            \index{Classical limit!Strict}
            \label{eq:ClassicalLimitContinuous}
            f \star_\hbar g
            \overset{\hbar \rightarrow 0}{\longrightarrow}
            f \cdot g
        \end{equation}
        and
        \begin{equation}
            \index{Semiclassical limit}
            \index{Semiclassical limit!Strict}
            \label{eq:SemiclassicalLimitContinuous}
            \gls{PartialDerivative}
            \bigl(
                f \star_\hbar g
                -
                g \star_\hbar f
            \bigr)
            \overset{\hbar \rightarrow 0}{\longrightarrow}
            \I
            \bigl\{
                f,
                g
            \bigr\},
        \end{equation}
        where the limits are taken for $\hbar \in \mathfrak{D}$ and within the
        topology of $\algebra{A}$.
    \end{definitionlist}
    In this case, we call $(\algebra{A},\star)$ a strict deformation quantization of
    $(\algebra{A}, \{\argument,\argument\})$ with deformation domain
    $\mathfrak{D}$.
\end{definition}

It should be noted that, at the time of writing and unlike for formal star products,
there is no general theory concerning the existence of star products. Instead,
Definition~\ref{def:StarProduct} attempts to capture the observations from a growing
list of examples. Generalizing further, one may allow for an asymmetry between the two
factors. Roughly speaking, demanding better properties of one factor affords the other
more freedom. A concrete realization of this idea in the context of the hyperbolic disk
is \cite[Sec.~7]{heins.moucha.roth.sugawa:2024a}. The underlying abstract concept of
locally convex modules within deformation quantization is discussed comprehensively
within \cite{lechner.waldmann:2016a}.

Having the origin as an accumulation point is necessary to formulate
\eqref{eq:ClassicalLimitContinuous} and \eqref{eq:SemiclassicalLimitContinuous}, which
subsume the formal classical and semiclassical limits from \eqref{eq:ClassicalLimit}
and \eqref{eq:SemiclassicalLimit}. Consequently, we refer to them as classical
and semiclassical limits, as well. The additional multiplication with the
imaginary unit $\I$ in
\eqref{eq:SemiclassicalLimitContinuous} is a convention ensuring that positive values
for $\hbar$ correspond to physical reality in the situations, where other
quantizations schemes are available and one thus knows what to expect. That is to say,
we replace the formal parameter~$\lambda$ with~$\I \hbar$ instead of $\hbar$. For
this reason, some authors include an additional factor $\I$ already
in~\eqref{eq:SemiclassicalLimit}.

While in many examples such as \cite{waldmann:2014a, esposito.stapor.waldmann:2017a,
heins.roth.waldmann:2023a} the
deformation domain may be chosen as the full complex plane and one is thus dealing
with an \emph{entire} dependence on $\hbar$, there are situations with $\mathfrak{D}
\subsetneq \C$. Indeed, this appears to go hand in hand with quotient constructions
and in particular coadjoint orbits, as studied in
\cite{esposito.schmitt.waldmann:2019a, kraus.roth.schoetz.waldmann:2019a,
schmitt:2021a}. The Wick star product on the unit disk and sphere even
exhibits simple poles accumulating at the origin, making the classical and
semiclassical limits rather subtle, see
\cite[Sec.~4.1]{kraus.roth.schoetz.waldmann:2019a} and
\cite[Thm.~6.8]{heins.moucha.roth.sugawa:2024a} as well as the discussion thereafter.
\begin{example}[Standard ordered star product II: Convergence]
    \index{Standard ordered!Star product}
    \index{Star product!Standard ordered}
    \label{ex:StdOrdII}
    We return to \\ the standard ordered star product from Example~\ref{ex:StdOrdI},
    where
    we replace
    $\lambda$ by $\I \hbar$. Its continuity was established and the consequences
    thereof were studied in
    \cite{omori.maeda.miyazaki.yoshiaki:1999a, omori.maeda.miyazaki.yoshioka:2002a,
    omori.maeda.miyazaki.yoshioka:2007a} and extended to vastly greater
    generality and including non-trivial infinite dimensional situations in
    \cite{waldmann:2014a}. It is instructive to retrace their proof, for which we shall
    tacitly utilize some of the technology we have yet to establish. We begin with the
    approach taken in \cite{waldmann:2014a}.

    We have already seen that $\algebra{A}_{\conv} \coloneqq
    \Pol(T^*\R)$ provides a Poisson subalgebra of $\Cinfty(T^*\R)$, on which the
    standard ordered
    star product \eqref{eq:StdExplicit} converges. Motivated by the explicit
    formula, we define seminorms
    \begin{equation}
        \index{Seminorms!Standard ordered}
        \label{eq:SeminormsStd}
        \seminorm{p}_c
        (\phi)
        \coloneqq
        \sum_{n=0}^{\infty}
        \frac{c^n}{n!^{1/2}}
        \abs[\big]
        {\phi^{(n)}(0)}
    \end{equation}
    for any $c \ge 0$ and $\phi \in \Cinfty(\R)$. Here, $\phi^{(n)}$ denotes the
    $n$-th derivative of $\phi$. Invoking once again the Borel Lemma, we see that
    \eqref{eq:SeminormsStd} has no reason to converge in general. This leads us to
    consider
    \begin{equation}
        \label{eq:EntireIntro}
        \gls{EntireGeneral}
        \coloneqq
        \big\{
            \phi \in \Comega(\R)
            \colon
            \seminorm{p}_c(\phi)
            <
            \infty
            \textrm{ for all }
            c \ge 0
        \big\},
    \end{equation}
    which we endow with the locally convex topology induced by the seminorms
    \eqref{eq:SeminormsStd}. The symbol $\Comega(\R)$ denotes the space of all
    real analytic functions $\phi \colon \R \longrightarrow \C$, to which we restrict
    in
    \eqref{eq:EntireIntro}
    to ensure the Hausdorff property. Taking another look at
    \eqref{eq:SeminormsStd}, we get that $\Pol(\R) \subseteq \Entire$. The power
    $1/2$ is ultimately a consequence of treating positions and momenta equally as it
    was done in \cite{omori.maeda.miyazaki.yoshiaki:1999a, waldmann:2014a}, which in
    this situation is not at all necessary
    by \cite[Lem.~6.2]{heins.roth.waldmann:2023a}, but simplifies the bookkeeping and
    the estimates.

    We claim that the pair
    $(\Entire \, \widehat{\tensor} \, \Entire, \star)$ with the standard ordered star product
    from
    \eqref{eq:StdExplicit} is a strict deformation quantization. Here, we endow the
    tensor square $\Entire \tensor \Entire$ with the projective tensor product
    topology, which we will discuss in detail in
    Section~\ref{sec:RTopologiesPolynomial}, and $\widehat{\tensor}$ denotes passage
    to the completion. For now, all we need is that it suffices to provide continuity
    estimates on
    factorizing tensors within $\Entire \tensor \Entire$ to infer the continuity of the product
    by Proposition~\ref{prop:InfimumArgument}. Taking
    \begin{equation}
        f(q,p)
        =
        \phi_1(q)
        \psi_1(p)
        \quad \textrm{and} \quad
        g(q,p)
        =
        \phi_2(q) \psi_2(p),
    \end{equation}
    we note the factorization
    \begin{align}
        f \star g
        &=
        \sum_{k=0}^{\infty}
        \frac{(\I \hbar)^k}{k!}
        (\phi_1 \cdot \phi_2^{(k)})
        \tensor
        (\psi_1^{(k)} \cdot \psi_2) \\
        &=
        (\phi_1 \tensor 1)
        \star
        (1 \tensor \psi_1)
        \star
        (\phi_2 \tensor 1)
        \star
        (1 \tensor \psi_2) \\
        &=
        (\phi_1 \tensor 1)
        \cdot
        (1 \tensor \psi_1)
        \star
        (\phi_2 \tensor 1)
        \cdot
        (1 \tensor \psi_2)
        \label{eq:StdFactorization}
    \end{align}
    by associativity of $\star$ and \eqref{eq:Std}, where $1$ denotes the
    constant function with value $1 \in \C$. Consequently, it suffices to prove the
    continuity of the pointwise product for $\Entire$ and the continuity of the mixed
    term in the middle. Using the crude estimate $\binom{n}{k} \le 2^n$ for non negative
    integers $n,k \in \N_0$ with $k \le n$ and the Cauchy product formula provides the
    continuity estimate
    \begin{align}
        \seminorm{p}_c
        (\phi \cdot \psi)
        &=
        \sum_{n=0}^{\infty}
        \frac{c^n}{n!^{1/2}}
        \abs[\Big]
        {
            \bigl(
            \phi \cdot \psi
            \bigr)^{(n)}
            (0)
        } \\
        &\le
        \sum_{n=0}^{\infty}
        \frac{c^n}{n!^{1/2}}
        \sum_{k=0}^n
        \binom{n}{k}
        \abs[\big]
        {
            \phi^{(k)}
            (0)
        }
        \cdot
        \abs[\big]
        {
            \psi^{(n-k)}
            (0)
        } \\
        &=
        \sum_{n=0}^{\infty}
        c^n
        \sum_{k=0}^n
        \binom{n}{k}^{1/2}
        \frac{1}{(n-k)!^{1/2} k!^{1/2}}
        \abs[\big]
        {
            \phi^{(k)}
            (0)
        }
        \cdot
        \abs[\big]
        {
            \psi^{(n-k)}
            (0)
        } \\
        &\le
        \sum_{n=0}^{\infty}
        \sum_{k=0}^n
        \frac{(\sqrt{2}c)^{k}}{k!^{1/2}}
        \cdot
        \abs[\big]
        {
            \phi^{(k)}
            (0)
        }
        \cdot
        \frac{(\sqrt{2}c)^{n-k}}{(n-k)!^{1/2}}
        \cdot
        \abs[\big]
        {
            \psi^{(n-k)}
            (0)
        } \\
        &=
        \biggl(
            \sum_{n=0}^{\infty}
            \frac{(\sqrt{2}c)^n}{n!^{1/2}}
            \cdot
            \abs[\big]
            {
                \phi^{(n)}
                (0)
            }
        \biggr)
        \cdot
        \biggl(
        \sum_{n=0}^{\infty}
        \frac{(\sqrt{2}c)^n}{n!^{1/2}}
        \cdot
        \abs[\big]
        {
            \psi^{(n)}
            (0)
        }
        \biggr) \\
        &=
        \seminorm{p}_{\sqrt{2}c}(\phi)
        \cdot
        \seminorm{p}_{\sqrt{2}c}(\psi)
    \end{align}
    for any $c \ge 0$ and $\phi, \psi \in \Entire$. Note that, heuristically, it is reasonable to
    expect the Cauchy product formula to make an appearance. Indeed,
    \eqref{eq:SeminormsStd} is, up to the square root, the Taylor series of $\phi$ with
    absolute values pulled all the way inside, and the Taylor series of the pointwise
    product is precisely given by the Cauchy product of the Taylor series. Similarly, we
    may estimate the mixed term in \eqref{eq:StdFactorization} by
    \begin{align}
        &\bigl(
            \seminorm{p}_{c_1}
            \tensor
            \seminorm{p}_{c_2}
        \bigr)
        \bigl(
            (1 \tensor \psi)
            \star
            (\phi \tensor 1)
        \bigr) \\
        &\le
        \sum_{k=0}^{\infty}
        \frac{\abs{\hbar}^k}{k!}
        \bigl(
           \seminorm{p}_{c_1}
           \tensor
           \seminorm{p}_{c_2}
        \bigr)
        (
            \phi^{(k)}
            \tensor
            \psi^{(k)}
        ) \\
        &=
        \sum_{k=0}^{\infty}
        \frac{\abs{\hbar}^k}{k!}
        \sum_{n,m=0}^{\infty}
        \frac{c_1^n \cdot c_2^m}{n!^{1/2} \cdot m!^{1/2}}
        \abs[\big]
        {\phi^{(k+n)}(0)}
        \cdot
        \abs[\big]
        {\psi^{(k+m)}(0)} \\
        &=
        \sum_{k=0}^{\infty}
        \frac{\abs{\hbar}^k}{c_1^k \cdot c_2^k}
        \sum_{n,m=k}^{\infty}
        \frac{c_1^n \cdot c_2^m}{(n-k)!^{1/2} \cdot (m-k)!^{1/2}}
        \biggl(
            \frac{n!}{n! \cdot k!}
            \cdot
            \frac{m!}{m! \cdot k!}
        \biggr)^{1/2}
        \abs[\big]
        {\phi^{(n)}(0)}
        \cdot
        \abs[\big]
        {\psi^{(m)}(0)} \\
        &=
        \sum_{k=0}^{\infty}
        \frac{\abs{\hbar}^k}{c_1^k \cdot c_2^k}
        \sum_{n,m=k}^{\infty}
        \frac{c_1^n \cdot c_2^m}{n!^{1/2} \cdot m!^{1/2}}
        \binom{n}{k}^{1/2}
        \binom{m}{k}^{1/2}
        \abs[\big]
        {\phi^{(n)}(0)}
        \cdot
        \abs[\big]
        {\psi^{(m)}(0)} \\
        &\le
        \sum_{k=0}^{\infty}
        \frac{\abs{\hbar}^k}{c_1^k \cdot c_2^k}
        \sum_{n,m=0}^{\infty}
        \frac{(\sqrt{2}c_1)^n \cdot (\sqrt{2}c_2)^m}{n!^{1/2} \cdot m!^{1/2}}
        \abs[\big]
        {\phi^{(n)}(0)}
        \cdot
        \abs[\big]
        {\psi^{(m)}(0)} \\
        &=
        \seminorm{p}_{\sqrt{2}c_1}(\phi)
        \cdot
        \seminorm{p}_{\sqrt{2}c_1}(\psi)
        \cdot
        \sum_{k=0}^{\infty}
        \frac{\abs{\hbar}^k}{c_1^k \cdot c_2^k}
    \end{align}
    for any $c_1, c_2 \ge 0$ and $\phi,\psi \in \Entire$. Hence, noting that
    \begin{equation}
        \seminorm{p}_{c}
        \le
        \seminorm{p}_{C}
        \qquad
        \textrm{whenever}
        \qquad
        c \le C,
    \end{equation}
    we may choose $c_1 \ge \abs{2 \hbar}^{1/2}$ and $c_2 \ge \abs{2 \hbar}^{1/2}$ to
    arrive at the estimate
    \begin{equation}
        \bigl(
        \seminorm{p}_{c_1}
        \tensor
        \seminorm{p}_{c_2}
        \bigr)
        \bigl(
        (1 \tensor \psi)
        \star
        (\phi \tensor 1)
        \bigr)
        \le
        2
        \cdot
        \seminorm{p}_{\sqrt{2}c_1}(\phi)
        \cdot
        \seminorm{p}_{\sqrt{2}c_1}(\psi),
    \end{equation}
    which is locally uniform with respect to all involved quantities. In
    Corollary~\ref{cor:FrechetPowerSeries} we will see that this implies the holomorphy
    of \eqref{eq:StarProductHolomorphy}. As the standard ordered star product fulfils
    \eqref{eq:SemiclassicalLimitContinuous}, a closer look at our estimate also proves
    the continuity of the Poisson bracket. Overall, we have thus -- up to some technical
    subtleties we will discuss in Section~\ref{sec:RTopologiesPolynomial} and the proof
    of Theorem~\ref{thm:LieGroupStarProductContinuity} in the more general context of
    Lie groups -- shown that the standard ordered star product is a star product in the
    sense of Definition~\ref{def:StarProduct}. We shall see in
    Example~\ref{ex:StdOrdIII}
    that the estimation of the pointwise product simplifies drastically if one works with
    holomorphic extensions instead of the real analytic functions within $\Entire$. Indeed,
    taking another look at~\eqref{eq:SeminormsStd}, we see that the Taylor series
    \begin{equation}
        \label{eq:StdTaylor}
        \gls{Taylor}(z)
        \coloneqq
        \sum_{n=0}^{\infty}
        \frac{z^n}{n!}
        \phi^{(n)}(0)
    \end{equation}
    of $\phi \in \Entire$ converges absolutely for all $z \in \C$. This defines an
    entire function $\Taylor_\phi \colon \C \longrightarrow \C$.
\end{example}

In the forms of Theorem~\ref{thm:LieAlgebraStarProductContinuity} and
Theorem~\ref{thm:LieGroupStarProductContinuity} we shall provide two more examples of
star products, both constructed within the context of Lie theory.

\section{Holomorphic Deformation Quantization}
\label{sec:QuantizationHolomorphic}
\epigraph{If you ignore the rules people will, half the time, quietly rewrite them so that
    they don’t apply to you.}{\emph{Equal Rites} -- Terry Pratchett}
% !TeX root = ../Dissertation.tex

In this section, we bring the idea that star products are inherently holomorphic
objects and should thus be treated as such to its logical conclusion. In
Example~\ref{ex:StdOrdII}, we have seen that weighted Taylor series provide a system of
seminorms, for which the standard ordered star product is continuous. Employing
suitable versions
of the Taylor formula -- such as the Lie-Taylor formula we shall meet in
\eqref{eq:LieTaylor} -- the phenomenon extends to other examples.

This motivates the following strategy, which assumes that the abstract Poisson algebra
$\algebra{A}$ is actually an algebra of functions defined on some smooth manifold $M$.
Then, we embed
\begin{equation}
    \label{eq:Complexification}
    M \hookrightarrow M_\C
\end{equation}
into a suitable complex manifold, which we refer to as its
complexification~$M_\C$. Afterwards, we study
its algebra of holomorphic functions \gls{Holomorphic} instead of the algebra of
real analytic functions~\gls{RealAnalytic}, which is much more complicated and
less well behaved from a functional analytic point of view. Pulling back with
\eqref{eq:Complexification} then provides a canonical linear mapping
\begin{equation}
    \Holomorphic(M_\field{C})
    \longrightarrow
    \Comega(M).
\end{equation}

\begin{figure}
    \begin{center}
        \includegraphics[width = 10cm]{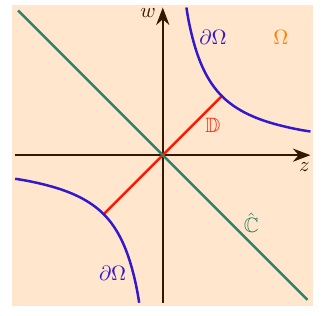}
        \caption{A schematic picture of the domain $\Omega \subseteq
        \widehat{\C}^2$}
        \label{fig:Omega}
    \end{center}
\end{figure}

Notably, there is -- not yet -- a general recipe to obtain $M_\C$ from $M$, but the existing
literature on holomorphic extensions and analytic continuation provides educated
guesses within concrete situations: For vector spaces $V$, it is natural to consider
$V_\C \coloneqq V \tensor \C$, as we shall at the start of
Section~\ref{sec:UniversalComplexification}, where we also discuss Hochschild's
\cite{hochschild:1966a} notion
of universal complexification \gls{Complexification} of Lie groups $G$. For
the star
products on the hyperbolic disk and the Riemann sphere $\widehat{\C}$
\cite{cahen.gutt.rawnsley:1994a, kraus.roth.schoetz.waldmann:2019a,
schmitt.schoetz:2022a} the natural complexification turned out to be
\begin{equation}
    \label{eq:Omega}
    \Omega
    \coloneqq
    \bigl\{
        (z,w)
        \in
        \widehat{\C}^2
        \colon
        zw
        \neq
        1
    \bigr\},
\end{equation}
where we extend the usual arithmetic of $\C$ to $\widehat{\C}$ by
\begin{equation}
    0 \cdot \infty
    \coloneqq
    1
    \eqqcolon
    \infty \cdot 0
    \quad \textrm{and} \quad
    \infty \cdot z
    \coloneqq
    \infty \eqqcolon
    \infty \cdot z
    \qquad
    \textrm{for all }
    z \in \C \setminus \{0\}.
\end{equation}
A schematic picture of $\Omega$
can be found in Figure~\ref{fig:Omega}. Its precise form may be recovered by
considering the hyperbolic disk and Riemann sphere as homogeneous spaces of real Lie
groups and then passing to their universal complexifications as in
Example~\ref{ex:ComplexificationOfComplex}. In the literature
\cite{kroetz.opdam:2008a, kroetz:2009a, kroetz.schlichtenkrull:2009a,
comporesi.kroetz:2012a}, the manifold $\Omega$ is known as the \emph{crown
domain}
associated to
the group~\gls{SpecialLinear} of real two by two matrices with determinant one. It
constitutes the natural domains of holomorphy for Laplace eigenfunctions, which was
independently proved in~\cite{kroetz:2009a} from the point of view of symmetric
spaces
and square integrable eigenfunctions and in \cite{heins.moucha.roth:2024b} in the
pursuit of Möbius invariant eigenspaces, generalizing results of \cite{rudin:1983a}.

Having obtained $M_\C$ in some manner, one may then investigate the given formal star
product \eqref{eq:StarProductFormal} as a \emph{single} holomorphic function
\begin{equation}
    \star
    \colon
    \mathfrak{D} \times \Holomorphic(M_\C) \tensor \Holomorphic(M_\C)
    \longrightarrow
    \Holomorphic(M_\C).
\end{equation}
Fixing $f,g \in \Holomorphic(M_\C)$ yields in particular holomorphic mappings
\begin{equation}
    \mathfrak{D} \times M_\C
    \ni
    (\hbar, z)
    \mapsto
    \bigl(
        f \star g
    \bigr)(z),
\end{equation}
now within the usual finite dimensional framework. This simple shift in perspective
allows for the utilization of complex analytic
techniques and turns out to be extraordinarily useful. As formal star products
are, by
definition, already given by a power series in $\hbar$, one should then
try to prove its \emph{absolute} convergence within the Fréchet space
$\Holomorphic(\mathfrak{D} \times M_\C)$ endowed with its natural topology of locally
uniform convergence, or at least a weighted version thereof. In practice, the required
inequalities have then always been based on a suitable version of the \emph{Cauchy
estimates}. Let us recall their classical incarnation for an entire function on~$\C^n$,
which can be readily proved by means of the Cauchy integral formula.
\begin{lemma}[Cauchy estimates]
    Let $\phi \in \Holomorphic(\C^n)$ and $z_0, v \in \C^n$. Then
    \begin{equation}
        \label{eq:CauchyEstimatesClassical}
        \index{Cauchy estimates!Classical}
        \abs[\big]
        {\phi^{(n)}(z_0)}
        \le
        \frac{n!}{r^n}
        \cdot
        \max_{\abs{z} = r}
        \abs[\big]
        {\phi(z_0 + z \cdot v)}
        \qquad
        \textrm{for all }
        r > 0
        \textrm{ and }
        n \in \N_0.
    \end{equation}
\end{lemma}

Conceptually, the Cauchy estimates intertwine the local data at $z_0$, namely the
Taylor coefficients of $\phi$, with its global behaviour in a quantitative manner. The
entirety of $\phi$ means that one gets to choose the expansion point $z_0$ and the
radius $r > 0$ freely. In practice, one makes use of this flexibility for each $n \in \N_0$
\emph{individually}, i.e. by choosing~$z_0$ and~$r$ as functions of $n \in \N_0$. This
sleight of
hand will be at the core of Example~\ref{ex:StdOrdIII} and a recurring theme within the
text. We are going to derive a locally convex version of the Cauchy estimates in
Proposition~\ref{prop:CauchyEstimatesLocallyConvex} and a Lie theoretic one in
Lemma~\ref{lem:CauchyEstimates}.

Moreover, one may take a step back and consider the function theory of $\algebra{A}$
as holomorphic functions. In the case of $\Omega$ this has been extraordinarily rich
and fruitful. The Peschl\footnote{Ernst Ferdinand Peschl (1906-1986) was a German
complex analyst of both one and several complex variables. His survey
\cite{peschl:1937a} on schlicht functions is still an interesting read. Academically,
Peschl is soon-to-be the authors grand-grandfather.}-Minda\footnote{David Minda is
an American complex analyst working as a professor at the University of Cincinnati. He
is interested in hyperbolic Riemann surfaces and geometric incarnations of classical
function-theoretic inequalities.} differential
operators, which take a role analogous to the operators~\eqref{eq:StarProductDiffops},
provide a bridge between strict
deformation quantization and complex analysis by \cite[Lem.~6.5 \&
Thm.~6.6]{heins.moucha.roth.sugawa:2024a}. Moreover, the geometric and
function-theoretic study of~$\Holomorphic(\Omega)$ has lead to the papers
\cite{schmitt.schoetz:2022a, heins.moucha.roth:2024b,
kraus.roth.schleissinger.waldmann:2023a,
heins.moucha.roth:2024a, moucha:2023a}.

The Wick\footnote{Gian Carlo Wick (1909-1992) was an Italian quantum field theorist.
The Wick rotation is a sleight of hand using a complex time variable to relate analytically
continued observables on Minkowski spacetime with functions defined on Euclidean
spacetime.} star product \gls{StarProductWick} of $f,g \in \Holomorphic(\Omega)$ is by
\cite[(6.13)]{heins.moucha.roth.sugawa:2024a} given by
\begin{equation}
    f \star_{\operatorname{Wick}} g
    =
    \sum_{n=0}^{\infty}
    \frac{(-1)^n}{n!}
    \frac{1}{(-1/\hbar) \cdot (-1/\hbar - 1) \cdots (-1/\hbar - n)}
    B_n(g,f)
\end{equation}
with the Peschl-Minda bidifferential operators $B_n$
\cite[(6.7)]{heins.moucha.roth.sugawa:2024a}, which already appeared in
\cite{bordemann.brischle.emmrich.waldmann:1996a,
bordemann.brischle.emmrich.waldmann:1996b}, albeit without their function-theoretic
interpretation.

Their precise form does not matter for the following observation: This
is not an entire function with respect to~$\hbar$! Instead, and as already mentioned
before, it exhibits simple poles at~$\hbar = -1/n$ for all~$n \in \N$. Similar effects
may be observed in the context of coadjoint orbits, see
\cite{esposito.schmitt.waldmann:2019a, schmitt:2021a}, and conjecturally arise in the
context of quotient constructions.

The standard ordered star product also allows for a
realization of the purely holomorphic strategy we have laid out.
\begin{example}[Standard ordered star product III: Holomorphy]
    \label{ex:StdOrdIII}
    \index{Standard ordered!Star product}
    \index{Star product!Standard ordered holomorphic}
    \index{Holomorphic!Standard ordered star product}
    We once again return to the scene of the crime, that is Example~\ref{ex:StdOrdII}.
    As before, we follow ideas from
    \cite[Prop.~6.1]{omori.maeda.miyazaki.yoshiaki:1999a}, which
    were generalized in \cite[Prop.~5.15]{waldmann:2014a}, but also adapt some of the
    discussion from \cite[Sec.~4.2]{heins.roth.waldmann:2023a}. It should be noted
    that, this time, we are quite close to the original treatment
    \cite{omori.maeda.miyazaki.yoshiaki:1999a} and essentially reverting the change of
    perspective \cite{waldmann:2014a} introduced.

    By convergence of the Taylor series \eqref{eq:StdTaylor}, we may extend every
    $\phi \in \Entire$ to a unique \emph{entire} function $\phi \in \Holomorphic(\C)$.
    We prove that these extensions are of finite order bounded above by two and of
    minimal type, which means that
    \begin{equation}
        \index{Seminorms!Standard ordered}
        \label{eq:SeminormsStdAlternative}
        \norm{\phi}_\epsilon
        \coloneqq
        \sup_{z \in \C}
        \abs[\big]
        {\phi(z)}
        \exp
        \bigl(
            -
            \epsilon
            \abs{z}^2
        \bigr)
        <
        \infty
    \end{equation}
    for all $\epsilon > 0$. Notice that $\norm{\argument}_\epsilon \le
    \norm{\argument}_\delta$ whenever $\delta \le \epsilon$. That is to say, the norms
    corresponding to \emph{small} values of $\epsilon$ already generate the topology. To
    be precise, $\epsilon \in (0,R)$ for any $R > 0$ suffices. This observation will be
    useful later.

    A comprehensive discussion of these concepts, including
    the techniques we will employ in a moment, can be found
    in the form of the monograph \cite{lelong.gruman:1986a}.
    More precisely, we shall show that the locally convex
    topologies induced by the systems of seminorms~\eqref{eq:SeminormsStd}
    and~\eqref{eq:SeminormsStdAlternative} are one and the same on $\Entire$. This
    is, essentially, a quantitative version of the identity principle: The finiteness of
    \eqref{eq:SeminormsStd} is a condition on the size of the Taylor data \emph{at
    one single point}, which turns out to be equivalent to a \emph{global growth
    condition} of the function itself. The classical Cauchy
    estimates \eqref{eq:CauchyEstimatesClassical} may be viewed as a statement of the
    very same flavour.

    Returning to our claim, there are mutual estimates to be derived. To this end, we
    fix~$\epsilon > 0$ and consider the auxiliary functions
    \begin{equation}
        f_n
        \colon
        (0,\infty) \longrightarrow (0,\infty), \quad
        f_n(r)
        \coloneqq
        r^n
        \cdot
        \exp
        \bigl(
            - \epsilon r^2
        \bigr)
        \qquad
        \textrm{for all }
        n \in \N_0.
    \end{equation}
    It is clear that each of the functions $f_n$ is bounded, vanishes for $r
    \rightarrow \infty$ and thus attains a global maximum at a unique $r_n \in
    (0,\infty)$. Computing the derivative reveals that
    \begin{equation}
        r_n
        =
        \sqrt
        {
            \frac{n}{2 \epsilon}
        }
        \qquad \textrm{and} \qquad
        f_n(r_n)
        =
        \biggl(
            \frac{n}{2 \epsilon}
        \biggr)^{n/2}
        \exp
        \bigl(
            -
            n/2
        \bigr)
        =
        \biggl(
            \frac{n}{2 \epsilon e}
        \biggr)^{n/2}.
    \end{equation}
    Applying the classical Cauchy estimates \eqref{eq:CauchyEstimatesClassical} to
    $f(z) \coloneqq \exp z$ with $r = n$ provides the elementary estimates
    \begin{equation}
        \label{eq:CauchyEstimateFactorial}
        \frac{1}{n!}
        \le
        \frac{e^n}{n^n}
        \qquad \iff \qquad
        \frac{n^n}{n! \cdot e^n}
        \le
        1
    \end{equation}
    for all $n \in \N_0$. This may also be proved directly by observing that the power
    series representation of $e^n$ consists of only non negative terms, the $n$-th
    being exactly $n^n/n!$. Either way, plugging in the Taylor series
    \eqref{eq:StdTaylor}, our preliminary considerations yield
    \begin{align}
        \norm{\phi}_\epsilon
        &\le
        \sum_{n=0}^{\infty}
        \frac{\abs{\phi^{(n)}(0)}}{n!}
        \sup_{z \in \C}
        \abs{z}^n
        \exp(-\epsilon \abs{z}^2) \\
        &=
        \sum_{n=0}^{\infty}
        \frac{\abs{\phi^{(n)}(0)}}{n!}
        \max_{r > 0}
        f_n(r) \\
        &=
        \sum_{n=0}^{\infty}
        \frac{\abs{\phi^{(n)}(0)}}{n!^{1/2}}
        \biggl(
            \frac{n^n}{n! \cdot \gls{Euler}^n}
        \biggr)^{1/2}
        \biggl(
            \frac{1}{2 \epsilon}
        \biggr)^{n/2} \\
        &\le
        \sum_{n=0}^{\infty}
        \frac{\abs{\phi^{(n)}(0)}}{n!^{1/2}}
        \biggl(
            \frac{1}{2 \epsilon}
        \biggr)^{n/2} \\
        &=
        \seminorm{p}_{1/\sqrt{2\epsilon}}(\phi).
    \end{align}
    We turn towards the converse inequality, for which we fix $c \ge 0$ and, by
    divine precognition, also $\epsilon \coloneqq 1/(8c^2)$. The strategy is to
    imitate the proof of the classical
    Cauchy estimates, using the additional information that $\norm{\phi}_\epsilon <
    \infty$. By Cauchy's integral formula, we have
    \begin{equation}
        \partial^n \phi(z_0)
        =
        \frac{n!}{2 \pi r^n}
        \int_{0}^{2\pi}
        \phi
        \bigl(
            z_0 + r e^{\I t}
        \bigr)
        e^{- \I n t}
        \D t
        \qquad
        \textrm{for all }
        n \in \N_0, \;
        r > 0,
        z_0 \in \C,
    \end{equation}
    which provides the estimate
    \begin{equation}
        \label{eq:CauchyEstimatesStandardOrderedPrototype}
        \abs[\big]
        {\partial^n \phi(z_0)}
        \le
        \frac{n!}{2\pi r^n}
        \int_0^{2\pi}
        \abs[\big]
        {\phi(z_0 + r e^{\I t})}
        \D t
        \le
        \frac{n!}{r^n}
        \max_{\abs{z} = r}
        \abs[\big]
        {\phi(z_0 + z)}
        \le
        n!
%        \cdot
        \frac{e^{\epsilon (\abs{z_0} + r)^2}}{r^n}
%        \cdot
        \norm{\phi}_\epsilon
    \end{equation}
    for all $n \in \N_0$, $r > 0$ and $z_0 \in \C$. Setting $z_0 \coloneqq 0$, the
    prefactor is minimal for $r = \sqrt{n/(2\epsilon)}$. Consequently,
    \begin{equation}
        \abs[\big]
        {\partial^n \phi(0)}
        \le
        \norm{\phi}_\epsilon
        \cdot
        n!
        \cdot
        \biggl(
            \frac{2 \epsilon e}{n}
        \biggr)^{n/2}
        \qquad
        \textrm{for all }
        n \in \N_0.
    \end{equation}
    Invoking once more \eqref{eq:CauchyEstimateFactorial} and our choice of
    $\epsilon$, this results in the desired estimate
    \begin{equation}
        \seminorm{p}_c
        (\phi)
        \le
        \norm{\phi}_\epsilon
        \sum_{n=0}^{\infty}
        \biggl(
            \frac{n!}{n^n \cdot e^n}
        \biggr)^{1/2}
%        \cdot
        \bigl(
            2 \epsilon c^2
        \bigr)^{n/2}
        \le
        2
        \cdot
        \norm{\phi}_\epsilon.
    \end{equation}
    Thus, the topologies induced by the systems of seminorms \eqref{eq:SeminormsStd}
    and \eqref{eq:SeminormsStdAlternative} are indeed one and the same. This has a
    number of pleasant consequences. The most immediate is that the continuity of the
    pointwise product, which was a non-trivial undertaking in
    Example~\ref{ex:StdOrdII}, is
    now obvious: Indeed, we simply have
    \begin{equation}
        \norm{\phi \cdot \psi}_\epsilon
        \le
        \norm{\phi}_{\epsilon/2}
        \cdot
        \norm{\psi}_{\epsilon/2}
        \qquad
        \textrm{for all }
        \phi, \psi \in \Entire, \;
        \epsilon > 0.
    \end{equation}
    Similarly, the continuity estimate for the inner product within the factorization
    \eqref{eq:StdFactorization} reduces to understanding
    \begin{equation}
        \label{eq:StdHolomorphicHalfway}
        \bigl(
        \norm{\argument}_{\epsilon_1}
        \tensor
        \norm{\argument}_{\epsilon_2}
        \bigr)
        \bigl(
        (1 \tensor \psi)
        \star
        (\phi \tensor 1)
        \bigr)
        \le
        \sum_{k=0}^{\infty}
        \frac{\abs{\hbar}^k}{k!}
        \norm[\big]
        {\phi^{(k)}}_{\epsilon_1}
        \cdot
        \norm[\big]
        {\psi^{(k)}}_{\epsilon_2}.
    \end{equation}
    To do this, we return to the penultimate step within
    \eqref{eq:CauchyEstimatesStandardOrderedPrototype} to estimate
    \begin{align}
        \sup_{z_0 \in \C}
        \abs[\big]
        {\partial^n \phi(z_0)}
        \cdot
        e^{- \epsilon \abs{z_0}^2}
        &\le
        \frac{n!}{r^n}
        \cdot
        \sup_{z_0 \in \C}
        e^{- \epsilon \abs{z_0}^2}
        \max_{\abs{z} = r}
        \frac{e^{- \epsilon \abs{z + z_0}^2/2}}
        {e^{- \epsilon \abs{z + z_0}^2/2}}
        \abs[\big]
        {\phi(z_0 + z)} \\
        &\le
        \frac{n!}{r^n}
        \biggl(
            \sup_{z_0 \in \C}
            \max_{\abs{z} = r}
            \frac{e^{- \epsilon \abs{z_0}^2}}{e^{- \epsilon \abs{z + z_0}^2/2}}
        \biggr)
        \biggl(
            \sup_{z_0 \in \C}
            \max_{\abs{z} = r}
            \abs[\big]
            {\phi(z_0 + z)}
            e^{- \epsilon \abs{z + z_0}^2/2}
        \biggr) \\
        &\le
        \norm[\big]
        {\phi}_{\epsilon/2}
        \frac{n!}{r^n}
            \sup_{z_0 \in \C}
            \exp
            \biggl(
                    -
                    \frac{\epsilon r^2}{2}
                    +
                    \epsilon r \abs{z_0}
                    +
                    \frac{\epsilon \abs{z_0}^2}{2}
            \biggr)  \\
        &=
        \norm[\big]
        {\phi}_{\epsilon/2}
        \frac{n!}{r^n}
        e^{\epsilon r^2},
    \end{align}
    as is readily verified by maximizing the quadratic polynomial within the exponent,
    leading to $r_0 = r$. This now has the same prefactor as
    \eqref{eq:CauchyEstimatesStandardOrderedPrototype}, so minimizing $r$ as before
    yields the
    Cauchy type estimate
    \begin{equation}
        \index{Cauchy estimates!Standard ordered}
        \label{eq:CauchyEstimatesStandardOrdered}
        \norm{\phi^{(k)}}_\epsilon
        \le
        \norm[\big]
        {\phi}_{\epsilon/2}
        \cdot
        n!
        \cdot
        \biggl(
            \frac{2 \epsilon e}{n}
        \biggr)^{n/2}
        \qquad
        \textrm{for all }
        n \in \N_0, \;
        \epsilon > 0.
    \end{equation}
    Continuing with \eqref{eq:StdHolomorphicHalfway} and assuming without loss of
    generality\footnote{But with a renewed and refreshing stroke of divine
    precognition. The discussion after \eqref{eq:SeminormsStdAlternative} explains
    why this is no loss of generality.} that the product
    \begin{equation}
        \epsilon_1 \epsilon_2
        \le
        \frac{1}{4e\abs{\hbar}},
    \end{equation}
    we get
    \begin{align}
        \sum_{k=0}^{\infty}
        \frac{\abs{\hbar}^k}{k!}
        \norm{\phi^{(k)}}_{\epsilon_1}
        \cdot
        \norm{\psi^{(k)}}_{\epsilon_2}
        &\le
        \norm[\big]
        {\phi}_{\epsilon_1/2}
        \cdot
        \norm[\big]
        {\psi}_{\epsilon_2/2}
        \sum_{k=0}^{\infty}
        \frac{\abs{\hbar}^k}{k!}
        k!^2
        \cdot
        \biggl(
            \frac{2 \epsilon_1 e}{k}
        \biggr)^{k/2}
        \biggl(
            \frac{2 \epsilon_2 e}{k}
        \biggr)^{k/2} \\
        &=
        \norm[\big]
        {\phi}_{\epsilon_1/2}
        \cdot
        \norm[\big]
        {\psi}_{\epsilon_2/2}
        \sum_{k=0}^{\infty}
        \bigl(
            2e \abs{\hbar} \epsilon_1 \epsilon_2
        \bigr)^k
        \cdot
        \frac{k!}{k^k} \\
        &\le
        2
        \cdot
        \norm[\big]
        {\phi}_{\epsilon_1/2}
        \cdot
        \norm[\big]
        {\psi}_{\epsilon_2/2}.
    \end{align}
    This provides a continuity estimate for $\star$ with respect to the alternative system of
    seminorms \eqref{eq:SeminormsStdAlternative}, which like the one we provided
    originally, is locally uniform regarding $\hbar$. While this derivation
    was arguably still technical, it should once again be stressed that the
    estimates we have derived are \emph{intrinsic} to entire functions we have
    considered. As such, they constitute very classical and well charted territory and
    can
    be found in various sources such as \cite[Ch.~1]{szego.polya:1976a},
    \cite[Ch.~1]{jank.volkmann:2014a} or the monograph \cite{levine:1997a}.
    For instance, the Cauchy type estimate \eqref{eq:CauchyEstimatesStandardOrdered}
    is a
    quantitative version of the fact that the order of an entire function is stable
    under differentiation, which is remarkable in its own right. We have thus
    separated matters of continuity from the study of function properties. This
    concludes our considerations of the standard ordered star product.
\end{example}

It is the principal purpose of Chapter~\ref{ch:InfiniteDimensions} to conduct
such an analysis in the context of the infinite dimensional generalization of the
standard ordered star product and $\Entire$ constructed in \cite{waldmann:2014a}.
Here, our focus in a thorough investigation of the type of holomorphy one
obtains. Similarly, Chapter~\ref{ch:Lie} concerns a holomorphic extension of the
real theory established in \cite{heins.roth.waldmann:2023a}. In this Lie
theoretic context, Hochschild's notion of universal complexification
\cite{hochschild:1966a} is going to provide the appropriate version of
complexification.

\chapter{Entire Functions on Vector Spaces}
\label{ch:InfiniteDimensions}
\epigraph{“It’s a lot more complicated than that—”

    “No. It ain’t. When people say things are a lot more complicated than that, they
    means they’re getting worried that they won’t like the truth.”}{\emph{Carpe
    Jugulum} -- Terry Pratchett}

% !TeX root = ../Dissertation.tex

The setting of this chapter is the one of locally convex spaces $V$ and holomorphic
mappings between them. Our principal interest lies in proving that the tensor
algebra of $V$ and its completions arising within strict deformation quantization may
be understood as \emph{bounded entire holomorphic functions} on the strong dual
$V'_\beta$. This is the content of
Theorem~\ref{thm:TensorsAsPolynomials}, Theorem~\ref{thm:IotaEmbedding} and
Theorem~\ref{thm:IotaEmbedding2}.

To facilitate these results, we give a streamlined
account of the theory of holomorphic functions between open subsets of locally convex
spaces with the goal of collecting various folklore results, while also supplying the
technical details or at least precise references. One major diversion we shall take
from most textbooks such as \cite{nachbin:1970a, dineen:1981a, dineen:1999a} is to
admit open subsets as the domains of our functions
throughout. The resulting theory looks much the same as for globally defined
functions, but comes at the cost of some additional technical
complications.\footnote{Which is really another way of saying \emph{interesting
phenomenona}.}

\index{Fréchet!René Maurice}
\index{Gâteaux!René}
After fixing some common notation at the end of this preliminary section, we proceed
with the study of polynomials between vector spaces and their continuity in
Section~\ref{sec:Polynomials}, where we choose an algebraic approach based on
multilinear mappings and focus on the resulting generalizations of the linear theory.
Afterwards,
we discuss the notion of Gâteaux holomorphy\footnote{Named
in honour of the French mathematician René Gâteaux (1889-1914), who was quite
tragically surprised by and subsequently lost to the first world war.} in
Section~\ref{sec:HolomorphicGateaux},
which once again results in a rather algebraic theory. Of particular importance will
be the validity of Taylor expansions for Gâteaux holomorphic functions, which
we
discuss in Theorem~\ref{thm:GateauxTaylor}. Conversely, it will be crucial that
pointwisely convergent power series, and thus in particular polynomials, are Gâteaux
holomorphic. We make this precise within Corollary~\ref{cor:GateauxPowerSeries},
which will be instrumental to interpret expressions such as
\eqref{eq:StarProductFormal} in an analytic manner.
Finally,
we discuss the well known connection between Gâteaux holomorphy and the existence of
directional -- or, indeed, Gâteaux -- derivatives in
Corollary~\ref{cor:GateauxVsDirectional}.

Adding the additional assumption of continuity then naturally leads to the
notion of Fréchet holomorphy\footnote{René Maurice Fréchet (1878-1973) is, among
many other important contributions to mathematics, the father of the modern notions of
metric space and compactness. }, which we study in
Section~\ref{sec:HolomorphicFrechet}. Here, the main
goals are the derivation of sufficient criteria for Fréchet holomorphy. Once again,
our particular interest lies in functions given by power series, which results in
Corollary~\ref{cor:FrechetPowerSeries}. Moreover, we shall utilize Baire's
\cite{baire:1899a, baire:1905a} classical insights into points of continuity of
pointwise limits of continuous functions. Indeed, as we shall see in
Proposition~\ref{prop:FrechetBaire}, the continuity of the Taylor polynomials of a
function defined on a Baire space at a \emph{singular} expansion point already implies
continuity of the full function. Using the explicit form of the Taylor
polynomials~\eqref{eq:GateauxTaylorPolynomials} as vector valued Cauchy integrals, we
moreover study how their continuity properties relate to the continuity of the full
function, and establish various continuity results with respect to the natural
function space topologies.

Along the way, it will become more and more apparent that the additional
hypothesis of \emph{boundedness} leads to a more fruitful theory, which we shall
explore in Section~\ref{sec:UniformConvergenceOnBoundedSets}. Here, we adopt the
notion customarily used in the context of linear operators, namely
that a bounded function should map bounded sets to bounded sets. For polynomials
defined on bornological spaces, boundedness turns out to be equivalent to
continuity, and thus it is natural to endow spaces of polynomials with the
topology of uniform convergence on bounded sets.

The situation for Fréchet holomorphic functions is then more complicated, as they need
not be bounded. We provide a concrete example of this behaviour within
Example~\ref{ex:FrechetLocallyUnbounded}. Taking a step back, we shall first assert in
Proposition~\ref{prop:BoundedCompleteness} that the space of all bounded functions is
complete with respect to the topology of uniform convergence on bounded sets,
provided
the codomain is complete itself.\footnote{Which, in view of the existence of constant
functions, means there are no non-trivial assumptions.} This insight provides
conceptually pleasing
and transparent proofs for the completeness of the spaces of polynomials with a
fixed degree of homogeneity, and the space of bounded Fréchet holomorphic functions.
The former is Theorem~\ref{thm:CompletenessHomogeneousPolynomials}, the latter is
the content of Theorem~\ref{thm:CompletenessBetaFrechet}. To facilitate the proof
of the latter, we make some significant additional assumptions on the involved
locally convex spaces. This then ties back to the criteria for Fréchet holomorphy
we have established within Section~\ref{sec:HolomorphicFrechet}.

Section~\ref{sec:RTopologiesPolynomial} then uses the theory we have established
to study the $\Sym_R$-topologies from strict deformation
quantization from the vantage point of infinite dimensional complex analysis. Most of
our results concern the case
$R=0$, which corresponds to the largest amount of functions that still behave in a
reasonable manner. In Example~\ref{ex:StdOrdII} we have studied the case $R = 1/2$
and
$V = \R$. In particular, it turns out that replacing projective
with injective tensor products within the construction -- or, alternatively
demanding nuclearity of the underlying space -- results in getting a canonical
\emph{embedding}\footnote{An
injective linear mapping that is a homeomorphism onto its image.} of the completed
symmetric algebra of $V$ into the space of bounded Fréchet holomorphic functions
defined on the strong dual space.

Having established the abstract theory, we study its implications in a
number of classical examples in Section~\ref{sec:EntireVectorSpacesExamples}. The
common theme here is that we use more detailed knowledge about the structure of the
strong duals $V'_\beta$ to make the abstract objects more concrete. This, in turn,
allows us to make and confirm educated guesses about how our constructions look in
detail. This also brings us back to the historical motivations behind the abstract theory,
namely the formalization of holomorphic functions depending on an infinite number of
complex variables.

\index{Riesz, Frigyes}
One such case is the one of $V$ being a Hilbert space, where $V'_\beta \cong V$
antilinearly and isometrically by the Riesz\footnote{The Hungarian mathematician
Frigyes Riesz (1880-1956) and his brother Marcel Riesz (1886-1969) pioneered
many central theorems of functional analysis. Their shared expertise makes it an
amusing game to guess which one of the two invented any particular Riesz concept
one comes across. Similar ideas apply to the Cartans in the context of
differential geometry as well as Lie theory.}
isomorphism. Working with Hilbert spaces, we also replace projective and injective
tensor products with the -- in this context much more natural -- Hilbert tensor
product. It turns out that most of our results from
Section~\ref{sec:RTopologiesPolynomial} pass into this
alternative setting unharmed, with the notable exception of conjecturally not getting
an
embedding, which we leave as Question~\ref{q:Embedding}.

\index{Locally convex space}
\index{Vector space!Locally convex}
As already indicated, we fix some notation before getting started with the study
of
polynomials. As customary, a locally convex space $V$ is a \emph{complex} vector
space, whose topology is induced by a system of seminorms $\seminorm{P}$ by taking
\begin{equation}
    \bigl\{
        \Ball_{\seminorm{p},r}(v_0)
        \colon
        \seminorm{p} \in \seminorm{P},
        v_0 \in V,
        r > 0
    \bigr\}
\end{equation}
as a subbasis of the topology, where
\begin{equation}
    \label{eq:OpenCylinder}
    \gls{OpenCylinder}
    \coloneqq
    \bigl\{
        v \in V
        \colon
        \seminorm{p}(v - v_0)
        <
        r
    \bigr\}
    =
    \Ball_{\seminorm{p},r}(0)
    +
    v_0
\end{equation}
denotes the open $\seminorm{p}$-cylinder around $v_0 \in V$ with radius $r > 0$.
The crucial feature of these cylinders is that, if $v_0 = 0$, they are
\emph{absolutely convex}, i.e. linear combinations of the form
\begin{equation}
    \lambda_1 v_1
    +
    \cdots
    +
    \lambda_n v_n
\end{equation}
satisfying $\lambda_1, \ldots, \lambda_n \in \C$ with $\abs{\lambda_1}, \ldots,
\abs{\lambda_n} \le 1$ and $v_1, \ldots, v_n \in \Ball_{\seminorm{p},r}(0)$ are
again contained in $\Ball_{\seminorm{p},r}(0)$. Taking another look at
\eqref{eq:OpenCylinder}, this in particular means that
$\Ball_{\seminorm{p},r}(v_0)$ is convex for any $v_0 \in V$. This also explains
why one speaks of \emph{local convexity} of such a space~$V$. The resulting
topology makes both the vector space addition and the multiplication by scalars
continuous as maps
\begin{equation}
    \label{eq:Plus}
    +
    \colon
    V \times V \longrightarrow V
\end{equation}
and
\begin{equation}
    \label{eq:ScalarMult}
    \cdot
    \colon
    \field{C} \times V
    \longrightarrow
    V.
    \index{Topological vector space}
    \index{Vector space!Topological}
\end{equation}
More generally, any vector space endowed with a topology making \eqref{eq:Plus} and
\eqref{eq:ScalarMult} continuous is called \emph{topological vector} space.
We denote the set of all continuous seminorms on a locally convex space $V$ by
\gls{ContinuousSeminorms},
which constitutes the maximal set of seminorms inducing the topology of $V$.
Finally, we note the following elementary Lemma on boundedness, which we shall
use in
many places throughout the text, sometimes without explicit mention.
\begin{lemma}
    \label{lem:BoundednessVsAbsoluteConvexity}
    Let $V$ be locally convex and $B \subseteq V$ be bounded. Then its absolutely
    convex hull $\gls{AbsolutelyConvexHull}(B)$, i.e. the smallest absolutely convex
    set containing
    $B$,
    is bounded as well.
\end{lemma}
\begin{proof}
    The crucial point is that the absolutely convex hull is explicitly given by
    \begin{equation}
        \absconv(B)
        =
        \bigl\{
            \lambda_1 v_1
            +
            \cdots
            +
            \lambda_n v_n
            \in
            V
            \colon
            n \in \N, \,
            \abs{\lambda_1}
            + \cdots +
            \abs{\lambda_n}
            \le
            1
        \bigr\}.
    \end{equation}
    Moreover, if $\seminorm{p} \in \cs(V)$, then $C \coloneqq \sup_{v \in B}
    \seminorm{p}(v) < \infty$ by boundedness of $B$. Thus, given vectors $v_1,
    \ldots, v_n \in B$ and scalars $\lambda_1,
    \ldots, \lambda_n \in \C$ with $\abs{\lambda_1} + \cdots +
    \abs{\lambda_n} \le 1$, then also
    \begin{equation}
        \seminorm{p}
        \biggl(
            \sum_{k=1}^{n}
            \lambda_k
            v_k
        \biggr)
        \le
        \abs{\lambda_1}
        \cdot
        \seminorm{p}(v_1)
        +
        \cdots
        +
        \abs{\lambda_n}
        \cdot
        \seminorm{p}(v_n)
        \le
        C
        <
        \infty,
    \end{equation}
    which proves the boundedness of $\absconv(B)$ by taking the supremum.
\end{proof}

As it should be becoming clear, we will require a variety of terminology from locally
convex analysis, a systematic discussion of which threatens to divert the reader's
attention from the crucial ideas. To declutter the main text, especially once we also add
Lie theory to our bubbling concoction in Chapter~\ref{ch:Lie}, we recall the precise
definitions of the supporting cast within numerous footnotes. Moreover, we have
gathered them centrally as a cheatsheet in the
Appendix~\hyperref[sec:LocallyConvexVocabulary]{``Terminology from locally convex
analysis''}. For a comprehensive discussion of the aspects of locally convex analysis
we have touched on so far and shall require throughout the text, we refer to the
textbooks \cite{osborne:2014a,
schaefer:1999a, koethe:1969a, koethe:1979a, meise.vogt:1992a, treves:2006a}
as well as the encyclopaedic \cite{jarchow:1981a}.

\newpage
\section{Polynomials}
\label{sec:Polynomials}
\epigraph{It's not worth doing something unless someone, somewhere, would much
rather you weren't doing it.}{\emph{Thief of Time} -- Terry Pratchett}
% !TeX root = ../Dissertation.tex
\index{Polynomials}
Recall that, if $L \colon V \times \cdots \times V \longrightarrow W$ is an $n$-linear
mapping between complex vector spaces $V$ and $W$, then
\begin{equation}
    \label{eq:PolynomialFromMultilinear}
    \gls{Polynomial}
    \colon
    V \longrightarrow W, \quad
    P(v)
    \coloneqq
    L(v,\ldots,v)
\end{equation}
defines an $n$-homogeneous mapping in the sense that
\begin{equation}
    \index{Homogeneity}
    \label{eq:NHomogeneity}
    P(\lambda v)
    =
    \lambda^n
    \cdot
    P(v)
    \qquad
    \textrm{for all $v \in V$ and $\lambda \in \C$.}
\end{equation}
We call $P$ a \emph{polynomial} of degree $n \in \field{N}_0$ and
say that the multilinear mapping $L$ induces~$P$.
\begin{example}[${\C[z_1, \ldots, z_k]}$]
    We quickly convince ourselves that this notion of polynomial is consistent with
    algebraically defined polynomials \gls{PolynomialsAlgebraic} in finitely many
    variables $z_1, \ldots, z_k$. Indeed, let
    \begin{equation}
        p
        \coloneqq
        z_{j_1}
        \cdots
        z_{j_n}
        \in
        \C[z_1, \ldots, z_k]
    \end{equation}
    for some $1 \le j_1, \ldots, j_n \le k$. Then the corresponding polynomial function $P \colon \C^k \longrightarrow \C$ is $n$-homogeneous and arises from the $n$-linear mapping
    \begin{equation}
        L
        \colon
        (\C^k)^n
        \longrightarrow
        \C, \quad
        L(v_1, \ldots, v_n)
        \coloneqq
        (v_1)^{j_1}
        \cdots
        (v_n)^{j_n},
    \end{equation}
    where the upper index denotes taking the corresponding coordinate. As we are
    dealing with complex vector spaces throughout the text and the ground field
    $\C$ is of
    characteristic zero, we shall not distinguish between polynomials and
    polynomial functions in the sequel. This also has the pleasant consequence
    that we may divide by $n!$ with reckless abandon.
\end{example}

Remarkably, one may always reconstruct the symmetric part
\begin{equation}
    \index{Symmetric!Part}
    \label{eq:SymmetricPart}
    \gls{SymmetricPart}
    \colon
    V^n \longrightarrow W, \quad
    L_s(v_1, \ldots, v_n)
    \coloneqq
    \frac{1}{n!}
    \sum_{\sigma \in S_n}
    L
    \bigl(
    v_{\sigma(1)}, \ldots, v_{\sigma(n)}
    \bigr)
\end{equation}
of $L$ from $P$ in an explicit manner by means of polarization. Here,
\gls{PermutationGroup} denotes the \emph{permutation or symmetric group} in the
letters $\{1,\ldots,n\}$. In the sequel, if $(s_1, \ldots, s_n)$ is some finite
ordered set, we also write
\begin{equation}
    \index{Symmetric!Group}
    \index{Permutation!Group}
    \index{Permutation!Action}
    \sigma
    \gls{Action}
    (s_1, \ldots, s_n)
    \coloneqq
    \bigl(
        s_{\sigma(1)}, \ldots, s_{\sigma(n)}
    \bigr)
\end{equation}
for the natural action of the permutation group.
\begin{proposition}[Polarization identity, {\cite[Cor.~1.6]{dineen:1999a}}]
    \index{Polarization!Identity}
    \label{prop:Polarization}
    Let $V, W$ be vector spaces and $P \colon V \longrightarrow W$ an
    $n$-homogeneous polynomial induced by $L \colon V^n \longrightarrow W$. Then
    \begin{equation}
        \label{eq:Polarization}
        L_s(v_1, \ldots, v_n)
        =
        \frac{1}{2^n \cdot n!}
        \sum_{\epsilon_1, \ldots, \epsilon_n = \pm 1}
        \epsilon_1 \cdots \epsilon_n
        \cdot
        P
        \biggl(
        \sum_{k=1}^{n}
        \epsilon_k v_k
        \biggr)
    \end{equation}
    for all $v_1, \ldots, v_n \in V$.
\end{proposition}
\begin{proof}
    In \cite[Prop.~1.5]{dineen:1999a} the identity \eqref{eq:Polarization} is given a
    probabilistic generalization, from which it may be derived by considering
    Bernoulli random variables with all probabilities given by $p = 1/2$. We go
    down a
    different
    path and verify \eqref{eq:Polarization} directly. Fixing vectors~$v_1,
    \ldots, v_n \in V$ and using \eqref{eq:PolynomialFromMultilinear} as well as
    the multilinearity of $L$, we first note
    \begin{align}
        P
        \biggl(
        \sum_{k=1}^{n}
        \epsilon_k v_k
        \biggr)
        &=
        L
        \biggl(
        \sum_{k_1=1}^{n}
        \epsilon_{k_1} v_{k_1},
        \ldots,
        \sum_{k_n=1}^{n}
        \epsilon_{k_n} v_{k_n}
        \biggr) \\
        &=
        \sum_{k_1,\ldots,k_n=1}^{n}
        \epsilon_{k_1} \cdots \epsilon_{k_n}
        \cdot
        L
        \bigl(
        v_{k_1},
        \ldots,
        v_{k_n}
        \bigr).
    \end{align}
    Fixing for the moment $k_1,\ldots,k_n \in \{1,\ldots,n\}$ and interchanging
    the summations within \eqref{eq:Polarization} leads thus
    \begin{equation}
        \sum_{\epsilon_1, \ldots, \epsilon_n = \pm 1}
        \epsilon_1 \cdots \epsilon_n
        \cdot
        \epsilon_{k_1} \cdots \epsilon_{k_n}
        \cdot
        L
        \bigl(
            v_{k_1},
            \ldots,
            v_{k_n}
        \bigr).
    \end{equation}
    If now $(k_1,\ldots,k_n)$ is not a permutation of $(1,\ldots,n)$, then there
    exists some index
    \begin{equation}
        j
        \in
        \{1,\ldots,n\} \setminus \{k_1,\ldots,k_n\}
    \end{equation}
    and thus
    \begin{equation}
        \sum_{\epsilon_j = \pm 1}
        \epsilon_j
        \cdot
        \epsilon_{k_1} \cdots \epsilon_{k_n}
        \cdot
        L
        \bigl(
        v_{k_1},
        \ldots,
        v_{k_n}
        \bigr)
        =
        \epsilon_{k_1} \cdots \epsilon_{k_n}
        \underbrace
        {
            \sum_{\epsilon_j = \pm 1}
            \epsilon_j
            \cdot
            L
            \bigl(
            v_{k_1},
            \ldots,
            v_{k_n}
            \bigr)
        }_{= \; 0}
        =
        0.
    \end{equation}
    Consequently, it suffices to sum over permutations $(k_1, \ldots, k_n)$ of
    $(1,\ldots, n)$. This yields
    \begin{equation}
        \epsilon_{k_1} \cdots \epsilon_{k_n}
        =
        \epsilon_1 \cdots \epsilon_n,
    \end{equation}
    and in view of \eqref{eq:SymmetricPart} we arrive at the desired result
    \begin{align}
        \frac{1}{2^n \cdot n!}
        \sum_{\epsilon_1, \ldots, \epsilon_n = \pm 1}
        \epsilon_1 \cdots \epsilon_n
        \cdot
        P
        \biggl(
        \sum_{k=1}^{n}
        \epsilon_k v_k
        \biggr)
        &=
        \sum_{\sigma \in S_n}
        \frac{1}{2^n \cdot n!}
        \sum_{\epsilon_1, \ldots, \epsilon_n = \pm 1}
        L
        \bigl(
        v_{\sigma(1)}, \ldots, v_{\sigma(n)}
        \bigr) \\
        &=
        \frac{1}{n!}
        \sum_{\sigma \in S_n}
        L
        \bigl(
        v_{\sigma(1)}, \ldots, v_{\sigma(n)}
        \bigr) \\
        &=
        L_s(v_1, \ldots, v_n).
        \tag*{$\qed$}
    \end{align}
\end{proof}

In the sequel, we call an $n$-linear mapping $L$ symmetric if $L_s = L$.
Conversely, the polynomial induced by
\eqref{eq:PolynomialFromMultilinear} only depends on the symmetric part $L_s$
of~$L$:
\begin{corollary}
    \label{cor:PolynomialInducing}
    Let $P \colon V \longrightarrow W$ be an homogeneous polynomial of degree $n \in
    \N_0$ between locally convex spaces $V$ and $W$. Then there exists a unique
    symmetric $n$-linear mapping inducing $P$.
\end{corollary}
\begin{proof}
    By definition of a polynomial, there exists some $n$-linear mapping $L \colon V^n
    \longrightarrow W$ inducing $P$; denote the polynomial induced by $L_s$ as $Q$.
    Invoking \eqref{eq:PolynomialFromMultilinear} and \eqref{eq:SymmetricPart}, we
    get
    \begin{equation}
        Q(v)
        =
        L_s(v,\ldots,v)
        =
        \frac{1}{n!}
        \sum_{\sigma \in S_n}
        L
        \bigl(
            \sigma
            \acts
            (v,\ldots,v)
        \bigr)
        =
        L(v,\ldots,v)
        =
        P(v)
    \end{equation}
    for $v \in V$. Thus the symmetric $n$-linear mapping $L_s$ induces $P$. If now $L' \colon V^n \longrightarrow W$ is another symmetric $n$-linear mapping inducing $P$, then \eqref{eq:Polarization} yields $L' = L_s$ at once.
\end{proof}

This warrants the shorthand notation \gls{InducingPolynomial} for the unique symmetric
multilinear mapping inducing a homogeneous polynomial $P$.\footnote{Reading
    $\widecheck{P}$ out loud yields $P$-check. This also reveals the corresponding
    \LaTeX-command.} We proceed with our investigation of the algebraic properties of
    polynomials. The following is essentially
\cite[Lemma~1.9]{dineen:1999a}.
\begin{proposition}
    \label{prop:PolynomialsAlgebraic}
    \index{Polynomials!Algebraic}
    Let $P \colon V \longrightarrow W$ be an homogeneous polynomial of degree $n \in
    \N_0$ between locally convex spaces $V$ and $W$ as well as $v,w \in V$.
    \begin{propositionlist}
        \item \label{item:PolynomialFixesZero}
        If $n>0$, then $P(0) = 0$.
        \item We have
        \begin{equation}
            \label{eq:PolynomialTaylorSeries}
            P(v + w)
            =
            \sum_{k=0}^{n}
            \binom{n}{k}
            \widecheck{P}
            \bigl(
            \underbrace{v,\ldots,v}_{k\text{-times}},
            \underbrace{w,\ldots,w}_{(n-k)\text{-times}}
            \bigr).
        \end{equation}
        \item It holds that
        \begin{equation}
            \label{eq:PolynomialsDifference}
            P(v) - P(w)
            =
            \sum_{k=0}^{n-1}
            \binom{n}{k}
            \widecheck{P}
            \bigl(
                \underbrace{w,\ldots,w}_{k\text{-times}},
                \underbrace{v-w,\ldots,v-w}_{(n-k)\text{-times}}
            \bigr).
        \end{equation}
    \end{propositionlist}
\end{proposition}
\begin{proof}
    The first statement is clear from \eqref{eq:NHomogeneity}. To see \eqref{eq:PolynomialTaylorSeries}, we note
    \begin{align}
        P(v_1 + v_2)
        &=
        \widecheck{P}(v_1+v_2, \ldots, v_1+v_2) \\
        &=
        \sum_{k_1, \ldots, k_n = 1}^{2}
        \widecheck{P}(v_{k_1}, \ldots, v_{k_n}) \\
        &=
        \sum_{k=0}^{n}
        \binom{n}{k}
        \widecheck{P}
        \bigl(
        \underbrace{v_1,\ldots,v_1}_{k\text{-times}},
        \underbrace{v_2,\ldots,v_2}_{(n-k)\text{-times}}
        \bigr)
    \end{align}
    for $v_1, v_2 \in V$ by $n$-linearity and symmetry of $\widecheck{P}$. Plugging in $v_1 = w$ and $v_2 = v - w$ yields in particular
    \begin{align}
        P(v)
        &=
        P(v_1 + v_2) \\
        &=
        \sum_{k=0}^{n}
        \binom{n}{k}
        \widecheck{P}
        \bigl(
        \underbrace{v_1,\ldots,v_1}_{k\text{-times}},
        \underbrace{v_2,\ldots,v_2}_{(n-k)\text{-times}}
        \bigr) \\
        &=
        \sum_{k=0}^{n}
        \binom{n}{k}
        \widecheck{P}
        \bigl(
        \underbrace{w,\ldots,w}_{k\text{-times}},
        \underbrace{v-w,\ldots,v-w}_{(n-k)\text{-times}}
        \bigr) \\
        &=
        P(w)
        +
        \sum_{k=0}^{n-1}
        \binom{n}{k}
        \widecheck{P}
        \bigl(
        \underbrace{w,\ldots,w}_{k\text{-times}},
        \underbrace{v-w,\ldots,v-w}_{(n-k)\text{-times}}
        \bigr),
    \end{align}
    which is \eqref{eq:PolynomialsDifference}.
\end{proof}

One may view \eqref{eq:PolynomialTaylorSeries} as the Taylor expansion of $P$
around $v$ evaluated at $w$; that is, the $k$-th Taylor polynomial of $P$ arises
from $P$ by plugging in $k$ copies of the evaluation point, leaving a homogeneous
polynomial of degree $n-k$. This of course matches with our intuition concerning
polynomials in finitely many variables. We will come back to this point of view in
Theorem~\ref{thm:GateauxTaylor} and Corollary~\ref{cor:GateauxVsDirectional}.
Identifying polynomials with the corresponding multilinear mappings allows us to
infer some elementary linear algebraic facts about polynomials.
\begin{proposition}[Linear Algebra of Polynomials]
    \label{prop:Polynomials}
    Let $V,W,U$ be locally convex spaces and $P \colon V \longrightarrow W$ a
    homogeneous polynomial of degree $n \in \N_0$.
    \begin{propositionlist}
        \item \label{item:PolynomialsHomogeneousVectorSpace}
        The set of all polynomials of degree $n$ mapping $V$ into $W$ forms a
        vector space.
        \item The set of all polynomials mapping $V$ into $W$ forms a vector space. If
        the codomain~$W$ is an algebra $\algebra{A}$, so is the set of all polynomials
        mapping $V$ into $\algebra{A}$.
        \item \index{Polynomials!Compositions}
        \label{item:PolynomialsComposition}
         If $Q \colon W \longrightarrow U$ is another polynomial of degree $m$,
        then the composition
        \begin{equation}
            Q \circ P \colon V \longrightarrow U
        \end{equation}
        is a polynomial of degree $n \cdot m$.
        \item \label{item:PolynomialsFiniteDimensions}
        If $V$ is finite dimensional, then the image $P(V)$ of $P$ is contained within
        a
        finite dimensional subspace of $W$.
    \end{propositionlist}
\end{proposition}
\begin{proof}
    The first two parts are obvious and so is the third, as the $(nm)$-linear
    mapping
    \begin{equation*}
       L(v_1, \ldots, v_{nm})
       \coloneqq
       \widecheck{Q}
       \bigl(
            \widecheck{P}(v_1, v_2, \ldots, v_m),
            \ldots,
            \widecheck{P}(v_{(n-1)m+1}, v_{(n-1)m+2}, \ldots, v_{nm})
       \bigr)
    \end{equation*}
    induces $Q \circ P$. Finally, part~\ref{item:PolynomialsFiniteDimensions} is a
    consequence of
    \begin{equation}
        P(V)
        \subseteq
        \widecheck{P}(V^n)
    \end{equation}
    and the fact that multilinear maps map finite dimensional spaces onto finite
    dimensional subspaces.
\end{proof}

\begin{example}[Multiplications as polynomials]
    Let $\algebra{A}$ be a $\C$-algebra. Then its multiplication
    \begin{equation}
        \mu
        \colon
        \algebra{A} \times \algebra{A} \longrightarrow \algebra{A}, \quad
        \mu(a,b)
        \coloneqq
        ab
    \end{equation}
    is bilinear and thus $P(a) \coloneqq \mu(a,a) = a^2$ constitutes a quadratic
    polynomial. Passing back to $\widecheck{P}$ yields, essentially, the
    \emph{anticommutator}
    \begin{equation}
        \widecheck{P}
        (a,b)
        =
        \frac
        {
            \mu(a,b)
            +
            \mu(b,a)
        }{2}
        =
        \frac{ab + ba}{2}.
    \end{equation}
    Note that the above works without assuming associativity, hence one may
    consider non-associative algebras, as well. However, if the corresponding
    multiplication $\mu$ is antisymmetric like in the case of Lie algebras, the
    corresponding polynomial simply vanishes and is thus not very interesting. For a
    commutative algebra $\algebra{A}$, the polarization identity
    \eqref{eq:Polarization} implies that all products are already determined by
    squares, a fact that is sometimes useful within the study of positivity. The
    explicit formula is also straightforward to guess; indeed, using the
    commutativity, we get
    \begin{equation}
        ab
        =
        \frac{1}{2}
        \bigl(
            (a+b)^2
            -
            a^2
            -
            b^2
        \bigr)
        \qquad
        \textrm{for all }
        a,b \in \algebra{A}.
    \end{equation}
\end{example}

Having noted these algebraic preliminaries and some first examples, we turn our
attention towards continuity, where we write
\begin{equation}
    \gls{PolynomialsHomogeneous}
    \coloneqq
    \bigl\{
        V
        \ni
        v
        \mapsto
        L(v,\ldots,v)
        \in
        W
        \;\big|\;
        L
        \colon
        V^n \longrightarrow W
        \; \textrm{is $n$-linear and continuous}
    \bigr\}
\end{equation}
for the set of all $n$-homogeneous \emph{continuous polynomials} between two topological vector spaces~$V$ and $W$. For $n=0$ this is to be understood as $\Pol^0(V,W) \coloneqq W$. Moreover, we define
\begin{equation}
    \label{eq:PolynomialAlgebra}
    \gls{Polynomials}
    \coloneqq
    \gls{DirectSum}
    \Pol^n(V,W).
\end{equation}
By Proposition~\ref{prop:Polynomials},
\ref{item:PolynomialsHomogeneousVectorSpace}, the set $\Pol^n(V,W)$ is a vector
space and the sum in \eqref{eq:PolynomialAlgebra} is indeed direct by virtue of
\eqref{eq:NHomogeneity}. If $W = \C$, we use the shorthand notations $\Pol^n(V)$
and~$\Pol^\bullet(V)$ instead.

Remarkably, not every mapping fulfilling
\eqref{eq:NHomogeneity} arises via \eqref{eq:PolynomialFromMultilinear} and
we exclude such pathologies from the start. In
Corollary~\ref{cor:GateauxNHomogeneousDuck} we will however establish that this
phenomenon is purely real and does not occur in the complex setting.

Note that the set of linear
polynomials $\Pol^1(V,W)$ is exactly the space of continuous linear mappings
from~$V$ to~$W$ and $\Pol^1(V)$ is the continuous dual space \gls{DualSpace} of
$V$. Using symmetric
tensor products~$\vee$, one may replace the multilinear mappings by linear ones, which
then identifies $\Pol^\bullet(V)$ as the continuous dual space of the $n$-th symmetric
power of $V^{\vee n}$. For a comprehensive discussion of this duality we refer to
\cite[Prop.~1.17]{dineen:1999a}.

As our notation suggests, the space $\Pol^\bullet(V)$ is a graded algebra with
respect to composition\footnote{Note that the grading is then done over
    the monoid $(\N_0,\times)$ rather than the usual $(\N_0,+)$. We shall provide a
    precise definition in \eqref{eq:Grading}.}, see again
    Proposition~\ref{prop:Polynomials}, \ref{item:PolynomialsComposition}, and also the
pointwise product, i.e.
\begin{equation}
    \cdot
    \colon
    \Pol^n(V)
    \times
    \Pol^m(V)
    \longrightarrow
    \Pol^{m+n}(V)
\end{equation}
is well defined for all $m,n \in \N_0$. This has the useful consequence that
homogeneous polynomials generate all of $\Pol(V)$, a fact we have already been
tacitly employing. The following inequality is an immediate consequence of
Proposition~\ref{prop:Polarization}.
\begin{corollary}[Polarization estimate, {\cite[(1.14)]{dineen:1999a}}]
    \label{cor:PolynomialContinuityVsMultilinear}
    \index{Polarization!Estimate}
    Let $P \colon V \longrightarrow W$ be an homogeneous polynomial of degree $n \in
    \N_0$ between locally convex spaces $V$ and $W$. Then $P$ is continuous if and
    only if $\widecheck{P}$ is. Moreover,
    \begin{equation}
        \label{eq:PolarizationEstimate}
        \sup_{v \in B}
        \seminorm{q}
        \bigl(
            P(v)
        \bigr)
        \le
        \sup_{v_1, \ldots, v_n \in B}
        \seminorm{q}
        \Bigl(
            \widecheck{P}(v_1, \ldots, v_n)
        \Bigr)
        \le
        \frac{n^n}{n!}
        \cdot
        \sup_{v \in B}
        \seminorm{q}
        \bigl(
        P(v)
        \bigr)
    \end{equation}
    holds for absolutely convex subsets $B \subseteq V$ and $\seminorm{q} \in \cs(W)$.
\end{corollary}
\begin{proof}
    The equivalence of the continuity of $P$ and $\widecheck{P}$ is a consequence of
    the polarization identity~\eqref{eq:Polarization}. As $\widecheck{P}(v,\ldots,v) =
    P(v)$, the first estimate in \eqref{eq:PolarizationEstimate} is clear. For the
    second, we invoke
    \eqref{eq:Polarization} and the absolute convexity of $B$, which yields
    \begin{align}
        \sup_{v_1, \ldots, v_n \in B}
        \seminorm{q}
        \Bigl(
            \widecheck{P}(v_1, \ldots, v_n)
        \Bigr)
        &\le
        \frac{1}{n! \cdot 2^n}
        \sum_{\epsilon_1, \ldots, \epsilon_n = \pm 1}
        \sup_{v_1, \ldots, v_n \in B}
        \seminorm{q}
        \biggl(
            P
            \Bigl(
                \sum_{k=1}^{n}
                \epsilon_k v_k
            \Bigr)
        \biggr) \\
        &=
        \frac{n^n}{n! \cdot 2^n}
        \sum_{\epsilon_1, \ldots, \epsilon_n = \pm 1}
        \sup_{v_1, \ldots, v_n \in B}
        \seminorm{q}
        \biggl(
        P
        \Bigl(
        \underbrace
        {
            \frac{1}{n}
            \sum_{k=1}^{n}
            \epsilon_k v_k
        }_{\in B}
        \Bigr)
        \biggr) \\
        &\le
        \frac{n^n}{n!}
        \cdot
        \sup_{v \in B}
        \seminorm{q}
        \bigl(
            P(v)
        \bigr),
    \end{align}
    which is the desired inequality.
\end{proof}

The constant $n!/n^n$ we picked up during the passage from the multilinear mapping to
the polynomial is best possible, see \cite[Example~1.39]{dineen:1999a} for a
comprehensive analysis. As the notation suggests, one is typically interested in
bounded subsets $B$ in \eqref{eq:PolarizationEstimate}. We will prove in
Lemma~\ref{lem:PolynomialContinuityVsBoundedness} that, for bornological domains
$V$, the finiteness of all expressions in~\eqref{eq:PolarizationEstimate} for bounded
$B \subseteq V$ is equivalent to the continuity of both~$P$ and~$\widecheck{P}$. In
analogy to linear mappings, many notions of
continuity turn out to coincide for polynomials. The following is a minor
generalization of \cite[Prop.~1.11]{dineen:1999a}.
\begin{proposition}[Continuity of Polynomials]
    \index{Polynomials!Continuity}
    \label{prop:PolynomialContinuity}
    Let $P \colon V \longrightarrow W$ be a homogeneous polynomial of degree $n \in
    \N_0$ between locally convex spaces $V$ and $W$. Then the following are
    equivalent:
    \begin{propositionlist}
        \item \label{item:PolynomialLocallyUniformlyContinuous}
        For every $v \in V$ there exists a neighbourhood $U$ of $v$ such that the restriction $P \at{U}$ is uniformly continuous.
        \item \label{item:PolynomialContinuous}
        The polynomial $P$ is continuous, i.e. $P \in \Pol^n(V,W)$.
        \item \label{item:PolynomialContinuousSomewhere}
        The polynomial $P$ is continuous at some $v \in V$.
        \item \label{item:PolynomialContinuousAtZero}
        The polynomial $P$ is continuous at zero.
        \item \label{item:PolynomialContinuityEstimate}
        For every continuous seminorm $\seminorm{q} \in \cs(W)$ there exist a constant $c > 0$ and continuous seminorm $\seminorm{p} \in \cs(V)$ such that
        \begin{equation}
            \label{eq:PolynomialContinuityEstimate}
            \seminorm{q}
            \bigl(
            P(v)
            \bigr)
            \le
            c
            \cdot
            \seminorm{p}(v)^n
            \qquad
            \textrm{for all $v \in V$.}
        \end{equation}
    \end{propositionlist}
\end{proposition}
\begin{proof}
    Clearly $\ref{item:PolynomialLocallyUniformlyContinuous}
    \implies
    \ref{item:PolynomialContinuous}
    \implies
    \ref{item:PolynomialContinuousSomewhere}
    \implies
    \ref{item:PolynomialContinuousAtZero}$.
    Assume \ref{item:PolynomialContinuousAtZero}. Invoking the polarization
    identity \eqref{eq:Polarization} once more, this means that $\widecheck{P}$ is
    continuous at zero, as well. Fix a continuous seminorm $\seminorm{q} \in \cs(W)$.
    Within the realm of locally
    convex spaces, the continuity at zero translates to the continuity estimate
    \begin{equation}
        \label{eq:PolynomialsContinuityProof}
        \seminorm{q}
        \bigl(
            \widecheck{P}(v_1, \ldots, v_n)
        \bigr)
        \le
        \seminorm{r}(v_1)
        \cdots
        \seminorm{r}(v_n)
    \end{equation}
    for some $\seminorm{r} \in \cs(V)$ and all $v_1, \ldots, v_n \in V$. Here we
    use the fact that powers of seminorms in $\cs(V)$ already generate the
    topology of the finite Cartesian product $V^n$. Consider now the --
    continuous! -- composition $\seminorm{s} \coloneqq (\seminorm{q} \circ
    P)^{1/n}$. Let $v,w \in V$ and $\lambda \in \C$. By construction, we then
    have $\seminorm{s}(\lambda
    v) = \abs{\lambda} \cdot \seminorm{s}(v)$ and by virtue of
    \eqref{eq:PolynomialTaylorSeries} also
    \begin{align}
        \seminorm{s}
        (v+w)^n
        &=
        \seminorm{q}
        \bigl(
            P
            (v+w)
        \bigr) \\
        &=
        \seminorm{q}
        \biggl(
            \sum_{k=0}^{n}
            \binom{n}{k}
            \widecheck{P}
            \bigl(
            \underbrace{v,\ldots,v}_{k\text{-times}},
            \underbrace{w,\ldots,w}_{(n-k)\text{-times}}
            \bigr)
        \biggr) \\
        &\le
        \sum_{k=0}^{n}
        \binom{n}{k}
        \seminorm{q}
        \bigl(
        \widecheck{P}
        \Bigl(
        \underbrace{v,\ldots,v}_{k\text{-times}},
        \underbrace{w,\ldots,w}_{(n-k)\text{-times}}
        \bigr)
        \Bigr) \\
        &\le
        \sum_{k=0}^{n}
        \binom{n}{k}
        \seminorm{r}(v)^k
        \cdot
        \seminorm{r}(w)^{n-k} \\
        &=
        \bigl(
            \seminorm{r}(v)
            +
            \seminorm{r}(w)
        \bigr)^n.
    \end{align}
    As this is not quite the triangle inequality, we are lead to define $\seminorm{p}
    \coloneqq \max\{\seminorm{r},\seminorm{q}\}$, which is now really a continuous
    seminorm on $V$. By
    construction, $\seminorm{p}$ fulfils \eqref{eq:PolynomialContinuityEstimate}
    with $c \coloneqq 1$. Thus we have shown
    \ref{item:PolynomialContinuityEstimate} and turn towards the converse
    implication \ref{item:PolynomialContinuityEstimate} $\implies$
    \ref{item:PolynomialContinuousAtZero}. As the collection of all unit cylinders
    $\Ball_{\seminorm{q},r}(0)$ for $\seminorm{q} \in \cs(W)$ and $r > 0$ generates the
    locally convex topology of $W$, it suffices to check that their preimages are
    zero neighbourhoods in $V$. Now, our assumption
    \eqref{eq:PolynomialContinuityEstimate} translates to
    \begin{equation*}
        \Ball_{\seminorm{p},r}(0)
        \subseteq
        P^{-1}
        \bigl(
        \Ball_{\seminorm{q},c r^n}(0)
        \bigr),
    \end{equation*}
    and thus $P^{-1}(\Ball_{\seminorm{q},c r^n}(0))$ is indeed a neighbourhood of $0$. We have established
    \begin{equation}
        \ref{item:PolynomialLocallyUniformlyContinuous}
        \implies
        \ref{item:PolynomialContinuous}
        \implies
        \ref{item:PolynomialContinuousSomewhere}
        \implies
        \ref{item:PolynomialContinuousAtZero}
        \iff
        \ref{item:PolynomialContinuityEstimate}
    \end{equation}
    and finally infer \ref{item:PolynomialLocallyUniformlyContinuous} from
    \ref{item:PolynomialContinuousAtZero}. To this end, fix again $\seminorm{q} \in
    \cs(W)$ and let $\seminorm{r} \in \cs(V)$ such
    that~\eqref{eq:PolynomialsContinuityProof}
    holds. Given $v_0 \in V$, we get for $v,w \in \Ball_{\seminorm{r},1/2}(v_0)$
    by using \eqref{eq:PolynomialsDifference} and
    \eqref{eq:PolynomialsContinuityProof}
    \begin{align}
        \seminorm{q}
        \bigl(
            P(v) - P(w)
        \bigr)
        &\le
        \sum_{k=0}^{n-1}
        \binom{n}{k}
        \seminorm{q}
        \Bigl(
            \widecheck{P}
            \bigl(
            \underbrace{w,\ldots,w}_{k\text{-times}},
            \underbrace{v-w,\ldots,v-w}_{(n-k)\text{-times}}
            \bigr)
        \Bigr) \\
        &\le
        \sum_{k=0}^{n-1}
        \binom{n}{k}
        \seminorm{r}(w)^k
        \cdot
        \seminorm{r}(v-w)^{n-k} \\
        &\le
        \seminorm{r}(v-w)
        \cdot
        \sum_{k=0}^{n-1}
        \binom{n}{k}
        \seminorm{r}(w)^k \\
        &\le
        \seminorm{r}(v-w)
        \cdot
        \bigl(
            1 + \seminorm{r}(w)
        \bigr)^{n} \\
        &\le
        \seminorm{r}(v-w)
        \cdot
        \bigl(
            3/2 + \seminorm{r}(v_0)
        \bigr)^n,
    \end{align}
    where we have also used that $\seminorm{r}(v-w) \le 1$ by the triangle inequality. This estimate is uniform in the open cylinder $U \coloneqq \Ball_{\seminorm{r},1/2}(v_0)$ and thus we have proved \ref{item:PolynomialLocallyUniformlyContinuous}.
\end{proof}

Having to restrict to locally uniform continuity in
\ref{item:PolynomialLocallyUniformlyContinuous} already appears for $V = \C = W$ and
the quadratic polynomial $z^2$. For constant and linear mappings,
one may of course upgrade back to uniform continuity on all of $V$.

\index{Bounded mapping}
Recall that a mapping $f \colon V \longrightarrow W$ between locally convex spaces $V$ and $W$ is called bounded if it maps bounded subsets of $V$ to bounded subsets of $W$. This is equivalent to
\index{Seminorms!$\beta$}
\begin{equation}
    \label{eq:BoundedSeminorms}
    \gls{SeminormsBounded}(f)
    \coloneqq
    \sup_{v \in B}
    \seminorm{q}
    \bigl(
    f(v)
    \bigr)
    <
    \infty
\end{equation}
for all bounded subsets $B \subseteq V$ and $\seminorm{q} \in \cs(W)$. If $W$ is
normed, we abbreviate $\seminorm{p}_{B,\norm{\argument}}$ to~$\seminorm{p}_{B}$.
Mimicking the linear theory, we get the following relationship between
boundedness and continuity for polynomials. \newpage
\begin{lemma}[Boundedness and Polynomials]
    \index{Bounded mapping!Polynomials}
    \index{Polynomials!Boundedness Vs. Continuity}
    \label{lem:PolynomialContinuityVsBoundedness}%
    Let $P \colon V \longrightarrow W$ be a polynomial between locally convex spaces $V$ and $W$.
    \begin{lemmalist}
        \item \label{item:PolynomialContinuityImpliesBoundedness}%
        If $P$ is continuous, then $P$ is bounded.
        \item \label{item:PolynomialBoundednessImpliesContinuity}%
        Let $V$ be bornological.\footnote{If $U \subseteq V$ is bornivorous, i.e. for every
        bounded set $B \subseteq V$ there is some $r>0$ such that $B \subseteq rU$,
        then it is a zero neighbourhood. See \cite[II.~7]{schaefer:1999a} or
        \cite[Sec.~4.1]{osborne:2014a}.} If $P$ is bounded, then it is continuous.
    \end{lemmalist}
\end{lemma}
\begin{proof}
    Without loss of generality, we may assume that $P$ is $n$-homogeneous. Let
    moreover $P$ be continuous, $B \subseteq V$ be bounded and $\seminorm{q} \in
    \cs(W)$. By Proposition~\ref{prop:PolynomialContinuity},
    \ref{item:PolynomialContinuityEstimate} there exist $c > 0$ and $\seminorm{p} \in
    \cs(V)$ such that \eqref{eq:PolynomialContinuityEstimate} holds for all $v \in V$. As
    $B$ is bounded, we have
    \begin{equation}
        \sup_{v \in B}
        \seminorm{p}(v)
        <
        \infty,
    \end{equation}
    and thus
    \begin{equation}
        \seminorm{p}_{B,\seminorm{q}}(P)
        =
        \sup_{v \in B}
        \seminorm{q}
        \bigl(
            P(v)
        \bigr)
        \le
        c
        \cdot
        \sup_{v \in B}
        \seminorm{p}(v)^n
        =
        c
        \cdot
        \bigl(
            \sup_{v \in B}
            \seminorm{p}(v)
        \bigr)^n
        <
        \infty.
    \end{equation}
    Consequently, $P$ is indeed bounded. Turning to
    \ref{item:PolynomialBoundednessImpliesContinuity}, we assume that $V$ is
    bornological and~$P$ is bounded. By
    Proposition~\ref{prop:PolynomialContinuity},
    \ref{item:PolynomialContinuousAtZero} it suffices to check continuity at the
    origin. To this end, let $U \subseteq W$ be an absolutely convex zero
    neighbourhood and $B \subseteq V$ be bounded. By assumption, the image $P(B)
    \subseteq W$ is bounded in $W$. Consequently, there is a radius $r > 0$ such that
    $P(B) \subseteq rU$, and thus
    \begin{equation*}
        B
        \subseteq
        P^{-1}
        \bigl(
        P(B)
        \bigr)
        \subseteq
        P^{-1}(rU)
        =
        r^{n}
        \cdot
        P^{-1}(U).
    \end{equation*}
    This means that the preimage $P^{-1}(U) \subseteq V$ is bornivorous. As $V$ is
    bornological, this implies that $P^{-1}(U)$ is a zero neighbourhood, proving the
    continuity of $P$ at the origin in view of
    Proposition~\ref{prop:PolynomialsAlgebraic},
    \ref{item:PolynomialFixesZero}.
\end{proof}

In the sequel, we endow $\Pol^n(V,W)$ and the polynomial algebra $\Pol^\bullet(V,W)$
with the locally convex topology induced by the seminorms
$\seminorm{p}_{B,\seminorm{q}}$ from \eqref{eq:BoundedSeminorms}. We indicate
this by
adding a subscript \gls{BoundedUniformConvergence} and call it the \emph{topology
of uniform convergence on
bounded subsets}, the strong topology or just the $\beta$-topology. We are going to
study the properties of the resulting locally convex spaces \gls{PolynomialsBounded}
and $\Pol_\beta^\bullet(V,W)$ in detail in
Section~\ref{sec:UniformConvergenceOnBoundedSets}, after establishing and thoroughly
discussing the notions of Gâteaux and Fréchet holomorphy.

\section{G\^ateaux Holomorphic Functions}
\label{sec:HolomorphicGateaux}
\epigraph{The young man is also an idealist. He has yet to find out that what’s in the
    public interest is not what the public is interested in.}{\emph{The Truth} --
    Terry Pratchett}
% !TeX root = ../Dissertation.tex
Our presentation mostly follows \cite[Sec.~3.1]{dineen:1999a}, but provides some
additional algebraic considerations and technical details in various places. We begin
by fixing the terminology we shall use.
\begin{definition}[Gâteaux Holomorphy]
    \index{Holomorphy!Gâteaux}
    \index{Gâteaux!Differentiability}
    \index{Gâteaux!Holomorphy}
    \label{def:Gateaux}
    Let $V,W$ be locally convex spaces, $U \subseteq V$ be open and $f \colon U
    \longrightarrow W$ a mapping.
    \begin{definitionlist}
        \item \label{item:GateauxDifferentiability}
        The map $f$ is called Gâteaux differentiable at $v_0 \in U$ if for
        every $v \in V$ and $\varphi \in W'$, the composition
        \begin{equation}
            \label{eq:Gateaux}
            z
            \quad \mapsto \quad
            \varphi
            \bigl(
                f(v_0+zv)
            \bigr)
            \in
            \C
        \end{equation}
        is holomorphic on some open neighbourhood of $z=0$ in $\C$.
        \item \label{item:GateauxHolomorphy}
        The map $f$ is called Gâteaux holomorphic if it is Gâteaux
        differentiable at all $v_0 \in U$.
        \item The set of all Gâteaux holomorphic mappings from $U$ to $W$ is denoted by
        \begin{equation}
            \gls{HolomorphicGateaux}
            \coloneqq
            \bigl\{
                f \colon U \longrightarrow W
                \colon
                f \text{ is Gâteaux holomorphic}
            \bigr\}.
        \end{equation}
    \end{definitionlist}
\end{definition}

We collect some first algebraic properties of the set of Gâteaux holomorphic functions.
\begin{proposition}
    \label{prop:Gateaux}
    \index{Gâteaux!Algebraic}
    Let $V,W$ be locally convex spaces and $U \subseteq V$ be open.
    \begin{propositionlist}
        \item \label{item:GateauxVectorSpace}
        The set $\Holomorphic_G(U,W)$ is a vector space with respect to the pointwise
        operations.
        \item \label{item:GateauxAlgebra}
        If $\algebra{A}$ an algebra, so is $\Holomorphic_G(U,\algebra{A})$ with respect to
        the pointwise operations.
        \item \label{item:GateauxSheaf}
        Sending open subsets $U' \subseteq V$ to $\Holomorphic_G(U',W)$ yields a sheaf
        $\Holomorphic_G(\argument,W)$ of vector spaces with restrictions given by the
        restriction of functions.
        \item \label{item:GateauxPolynomials}
        Every polynomial is Gâteaux holomorphic. In particular,
        \begin{equation}
            \label{eq:GateauxPolynomials}
            \index{Polynomials!Gâteaux holomorphy}
            \Pol^\bullet(V,W)
            \subseteq
            \Holomorphic_G(V,W).
        \end{equation}
    \end{propositionlist}
\end{proposition}
\begin{proof}
    The condition \eqref{eq:Gateaux} is linear and so \ref{item:GateauxVectorSpace} and
    \ref{item:GateauxAlgebra} are
    obvious. As we are dealing with functions, have defined Gâteaux holomorphy by
    Gâteaux
    differentiability at all points of the domain and the condition \eqref{eq:Gateaux}
    only depends on the values of $f$ in a neighbourhood of $v_0$, it is clear that
    $\Holomorphic_G(\argument,W)$ is a sheaf.
    %    It thus remains to study gluing. To this end, let $\{U_\alpha\}_{\alpha \in J}$ be a
    %collection of open subsets of $V$ with $f_\alpha \in \Holomorphic_G(U_\alpha,W)$
    %such that
    %    \begin{equation}
        %        f_\alpha
        %        \at{U_\alpha \cap U_\beta}
        %        =
        %        f_\beta
        %        \at{U_\alpha \cap U_\beta}
        %    \end{equation}
    %    whenever $\alpha,\beta \in J$ fulfil $U_\alpha \cap U_\beta \neq \emptyset$.
    To see \ref{item:GateauxPolynomials}, let $P \colon V \longrightarrow W$ be a
    -- not necessarily continuous -- $n$-homogeneous polynomial, $\varphi \in W'$
    and $v_0 \in
    V$. By~\eqref{eq:PolynomialTaylorSeries}, we then get
    \begin{align}
        \varphi
        \bigl(
        P(v_0 + zv)
        \bigr)
        &=
        \varphi
        \biggl(
        \sum_{k=0}^{n}
        \binom{n}{k}
        \widecheck{P}
        \bigl(
        \underbrace{v_0,\ldots,v_0}_{k\text{-times}},
        \underbrace{zv,\ldots,zv}_{(n-k)\text{-times}}
        \bigr)
        \biggr) \\
        &=
        \sum_{k=0}^{n}
        \binom{n}{k}
        z^{n-k}
        \cdot
        \varphi
        \Bigl(
        \widecheck{P}
        \bigl(
        \underbrace{v_0,\ldots,v_0}_{k\text{-times}},
        \underbrace{v,\ldots,v}_{(n-k)\text{-times}}
        \bigr)
        \Bigr)
    \end{align}
    for $z \in \C$ and $v \in V$. As this is simply a complex polynomial of degree $n$ with
    respect to~$z$, it is certainly holomorphic, proving $P \in \Holomorphic_G(V,W)$.
\end{proof}

\begin{remark}[Hartogs' Theorem on separate holomorphy]
    \index{Hartogs' Theorem}
    \index{Hartogs' Theorem!Finite dimensional}
    \index{Hartogs', Friedrich Moritz}
    \label{rem:Hartog} Recall
    Hartogs'\footnote{Friedrich Moritz Hartogs (1874-1943) was a German Jewish
    mathematician, who was ultimately discriminated against and killed by the
    Nazi regime. He
    proved several groundbreaking results concerning power series expansion and
    singularities of holomorphic functions in several complex variables.} \\ seminal
    Theorem on separate holomorphy: If $f \colon \C^n \longrightarrow \C$ is such that
    \begin{equation}
        \C
        \ni
        z_j
        \quad \mapsto \quad
        f(z_1, \ldots, z_n)
        \in
        \C,
    \end{equation}
    where the remaining entries of $z = (z_1, \ldots, z_n)$ are fixed, is holomorphic
    for $j=1, \ldots, n$, then the full function $f$ is holomorphic. This means that it is
    totally differentiable with complex linear differential. In particular, this then implies
    the continuity of $f$. The proof of this remarkable
    result is based on Baire's category theorem, which is used to infer the
    boundedness of $f$ on bounded sets. Having established this, one may then
    construct local Taylor expansions for such functions by means of Cauchy
    integrals, which implies the holomorphy. For a detailed discussion and the proof,
    we refer to \cite[Thm.~1.2.5]{krantz:2001a}. The first major goal of this
    section is to adapt this argument to construct Taylor expansions of Gâteaux
    holomorphic functions as series over homogeneous polynomials,
    culminating in Theorem~\ref{thm:GateauxTaylor}. One may think of these expansions
    as power series in infinitely many variables, and historically this is indeed the
    idea they make precise.

    To this end, we first note the following: By virtue of Hartogs' Theorem, we
    may rephrase Gâteaux holomorphy of some mapping $f \colon V \longrightarrow W$
    as all restrictions $f \at{F} \colon F \longrightarrow W$ of~$f$ to finite
    dimensional subspaces being holomorphic. This yields the following Gâteaux
    holomorphic version of Hartogs' Theorem.
\end{remark}

\begin{proposition}[Hartogs' Theorem for Gâteaux holomorphy]
    \label{prop:GateauxHartog}
    \index{Hartogs' Theorem!Gâteaux}
    Let $V_1, V_2, W$ be locally convex spaces with open subsets $U_1
    \subseteq V_1$ and $U_2 \subseteq V_2$. Let moreover
    \begin{equation}
        f
        \colon
        U_1 \times U_2
        \longrightarrow
        W
    \end{equation}
    be such that both $f(v_1, \argument)$ and $f(\argument,
    v_2)$ are Gâteaux holomorphic for all vectors $v_1 \in U_1$ and $v_2 \in
    U_2$. Then $f$ is Gâteaux holomorphic.
\end{proposition}
\begin{proof}
    If $F \subseteq V_1 \times V_2$ is finite dimensional, then so is
    \begin{equation}
        F'
        \coloneqq
        \pi_1(F)
        \times
        \pi_2(F)
    \end{equation}
    with the canonical linear projections $\pi_1 \colon V_1 \times V_2
    \longrightarrow V_1$ and $\pi_2 \colon V_1 \times V_2 \longrightarrow V_2$. By
    assumption, $f \at{F'}$ is separately holomorphic and thus holomorphic by the
    classical Hartogs' Theorem. Consequently, $f$ is Gâteaux holomorphic, as well.
\end{proof}

This exposes the natural topology of $\Holomorphic_G(V,W)$ as the one of
\emph{locally uniform convergence on finite dimensional subspaces}. However,
trying to recast this as a seminorm condition, one runs into the complication
that, in principle, the functions \eqref{eq:Gateaux} corresponding to a sequence
$(f_n)_{n \in \N_0} \subseteq \Holomorphic_G(U,W)$ need not be defined on a common
domain. We rule out this inconvenience with a preparatory lemma.
\begin{lemma}
    \label{lem:GateauxZeroNeighbourhoods}
    Let $V,W$ be locally convex spaces, $U \subseteq V$ be open and $v_0 \in V$.
    Writing
    \begin{equation}
        \label{eq:GateauxMaximalNeighbourhood}
        U_{v}
        \coloneqq
        \bigl\{
            z \in \C \colon v_0 + zv \in U
        \bigr\}
        \subseteq
        \C
    \end{equation}
    for $v \in U$, the mapping
    \begin{equation}
        \label{eq:GateauxMaximal}
        U_v
        \ni
        z
        \quad \mapsto \quad
        \varphi
        \bigl(
        f(v_0+zv)
        \bigr)
        \in
        \C
    \end{equation}
    is holomorphic for all $f \in \Holomorphic_G(U,W)$ and $\varphi \in W'$.
\end{lemma}
\begin{proof}
    The crucial point is that we have assumed Gâteaux holomorphy and not just Gâteaux
    differentiability at a single point. Fix $v_0 \in U$, $z_0 \in U_v$ and $v \in V$.
    We define
    \begin{equation}
        v_1 \coloneqq v_0 + z_0v
    \end{equation}
    and denote the mappings from \eqref{eq:GateauxMaximal} corresponding to $v_j$
    by $g_j \colon U_v \longrightarrow \C$ for $j=0,1$. Then
    \begin{equation}
        g_{0}(z)
        =
        \varphi
        \bigl(
            f(v_0+zv)
        \bigr)
        =
        \varphi
        \bigl(
            f
            (v_1 + (z-z_0)v)
        \bigr)
        =
        g_{1}(z-z_0)
    \end{equation}
    holds for all $z \in U_v$. But by assumption, $g_1$ is holomorphic in a zero
    neighbourhood and thus $g_0$ is holomorphic in a neighbourhood of $z_0$. As
    this holds for all $z_0 \in U_v$, we have completed the proof.
\end{proof}

The advantage of our, a priori less restrictive, definition is that in practice
one often checks the holomorphy of \eqref{eq:Gateaux} by producing a local power
series expansion, whose convergence one then does not need to check on all of
$U_v$. Conversely, Lemma~\ref{lem:GateauxZeroNeighbourhoods} ensures that the
domains of \eqref{eq:Gateaux} are intrinsic to the open set $U$ and independent of
the particular function one considers. Indeed, in Theorem~\ref{thm:GateauxTaylor},
we will ultimately see that a similar statement is true about Taylor expansions for
Gâteaux holomorphic functions on $U$.

Having dealt with this subtlety, we may establish the completeness of the set of
Gâteaux holomorphic functions, which in the sequel we will always endow with the
topology of locally uniform convergence on finite dimensional subspaces. A generating
system of seminorms is given by
\index{Seminorms!Gâteaux}
\begin{equation}
    \label{eq:GateauxSeminorms}
    \Holomorphic_G(U,W)
    \ni
    f
    \quad \mapsto \quad
    \max_{v \in K}
    \seminorm{q}
    \bigl(
        f(v)
    \bigr),
\end{equation}
where we vary over compact sets $K \subseteq F \cap U$ contained in finite dimensional
subspaces~$F$ of $V$ intersected with the open set $U$ and continuous
seminorms $\seminorm{q} \in \cs(W)$.
\begin{proposition}
    \index{Gâteaux!Completeness}
    \index{Completeness!Gâteaux}
    \label{prop:GateauxTopology}
    Let $V,W$ be locally convex spaces such that $W$ is complete, and let~$U \subseteq
    V$ be open. Then the space $\Holomorphic_G(U,W)$ is complete with respect to the
    topology of locally uniform convergence on finite dimensional subspaces. If $W$ is
    merely sequentially complete, the same is true for $\Holomorphic_G(U,W)$.
\end{proposition}
\begin{proof}
    Let $(f_\alpha)_{\alpha \in J} \subseteq \Holomorphic_G(U,W)$ be a Cauchy net. As
    $W$ is complete and the topology of locally uniform convergence on finite
    dimensional subspaces is finer than the topology of pointwise convergence, we may
    consider the pointwise limit
    \begin{equation}
        f
        \colon
        U \longrightarrow W, \quad
        f(v)
        \coloneqq
        \lim_{\alpha \in J}
        f_\alpha(v).
    \end{equation}
    Fix $v_0 \in U$, $v \in V$ and $\varphi \in W'$. By
    Lemma~\ref{lem:GateauxZeroNeighbourhoods} the net $(g_\alpha)_{\alpha \in J}$
    given by
    \begin{equation}
        g_\alpha
        \colon
        U_v \longrightarrow \C, \quad
        g_\alpha(v)
        \coloneqq
        \varphi
        \bigl(
            f_\alpha(v_0+zv)
        \bigr)
    \end{equation}
    is well defined. By assumption, $(g_\alpha)_{\alpha \in J}$ is Cauchy in
    $\Holomorphic(U_v)$ and thus locally uniformly convergent to some $g \in
    \Holomorphic(U_v)$. By continuity of $\varphi$, it is clear that
    \begin{equation}
        g(v)
        =
        \varphi
        \bigl(
            f(v_0+zv)
        \bigr)
        \qquad
        \textrm{holds for all }
        v \in U_v.
    \end{equation}
    We have thus established $f \in \Holomorphic_G(U,W)$ and that $f_\alpha \rightarrow
    f$ locally uniformly on finite dimensional subspaces.
\end{proof}

Using the Gâteaux holomorphy of polynomials it is easy to establish that power
series, i.e. pointwisely convergent series of homogeneous polynomials, are Gâteaux
holomorphic.
\begin{corollary}
    \label{cor:GateauxPowerSeries}
    Let $(P_n)_{n \in \N_0}$ be a sequence of polynomials
    between locally convex spaces~$V$ and $W$ such that $P_n$ is homogeneous of
    degree
    $n$ for all $n \in \N_0$. Assume moreover that the series
    \begin{equation}
        \label{eq:GateauxPowerSeries}
        f(v)
        \coloneqq
        \sum_{n=0}^{\infty}
        P_n(v)
    \end{equation}
    converges in $W$ for all $v$ in some open subset $U \subseteq V$. Then the
    resulting mapping
    \begin{equation}
        f
        \colon
        U \longrightarrow W
    \end{equation}
    is Gâteaux holomorphic.
\end{corollary}
\begin{proof}
    The idea is that the partial sums of \eqref{eq:GateauxPowerSeries} are Cauchy in
    $\Holomorphic_G(U,W)$. Indeed, let $F \subseteq V$ be a finite dimensional
    subspace. Then, by Proposition~\ref{prop:Gateaux}, \ref{item:GateauxPolynomials}
    and Remark~\ref{rem:Hartog}, the restrictions
    \begin{equation}
        P_n
        \at{F}
        \colon
        F
        \longrightarrow
        W
    \end{equation}
    are holomorphic for all $n \in \N_0$. Consequently, after choosing coordinates, they
    are \emph{globally} given by locally uniformly convergent power series on $F$. This
    shows that the
    partial sums of \eqref{eq:GateauxPowerSeries} are indeed Cauchy in
    $\Holomorphic_G(U,W)$. By Proposition~\ref{prop:GateauxTopology}, the series
    \eqref{eq:GateauxPowerSeries} thus converges to the pointwise limit $f$ in
    $\Holomorphic(U,\widehat{W})$, where~$\widehat{W}$ denotes the completion of $W$.
    As by assumption, all necessary limits are already contained in $W$, i.e.~$f$
    actually maps into $W$, this completes the proof. Note that by virtue of the
    continuous dual $W'$ and the continuous dual
    $\widehat{W}'$ of the completion being essentially the same space,\footnote{By
    means of restriction and the extension principle for locally uniformly continuous
    mappings.}
    the co-restriction back to
    $W$ does not interfere with Gâteaux differentiability.
\end{proof}

Another natural question is whether compositions of Gâteaux holomorphic functions are
again Gâteaux holomorphic. This turns out to be a bit more subtle: Indeed, the
condition~\eqref{eq:Gateaux} is linear and thus, at least a priori, incompatible with
precomposition with general maps. In Corollary~\ref{cor:GateauxPullbackLinear}, we
will conversely prove that linearity suffices for good behaviour of the pullback.

Returning to the general problem in the guise of restrictions to finite dimensional
subspaces, see again Remark~\ref{rem:Hartog}, we run into a similar obstacle: Gâteaux
holomorphic functions
\begin{equation}
    f
    \colon
    \C \supseteq U
    \longrightarrow
    W
\end{equation}
may already possess an infinite dimensional image, even when restricted to compact
subsets.\footnote{This is not a bug, but rather a feature: Otherwise, the theory of
vector-valued holomorphic functions would be rather dull.}
\begin{example}
    \index{Gâteaux!Infinite dimensional range}
    \label{ex:GateauxInfiniteDimensionalRange}
    Let $V \coloneqq \gls{ZeroSequences}$, the Banach space of zero sequences
    indexed by $\N_0$, and
    consider the standard basis vectors $e_n \in c_0$ with $e_n(k) \coloneqq
    \delta_{n,k}$ for $n,k \in \N_0$. Write
    \begin{equation}
         \gls{Disk}
         \coloneqq
         \bigl\{
            z \in \C \colon \abs{z} < 1
         \bigr\}
    \end{equation}
    for the open unit disk and consider the mapping $f \colon \disk
    \longrightarrow c_0$,
    \begin{equation}
        f(z)
        \coloneqq
        \sum_{n=0}^{\infty}
        e_n
        \cdot
        z^n
        =
        \bigl(
            1,
            z,
            z^2,
            \ldots
        \bigr),
    \end{equation}
    which converges within the Banach\footnote{Stefan Banach (1892-1945) was a
    Polish mathematician, and one of the founding fathers of modern functional
    analysis. He survived both world wars, but tragically succumbed young to lung
    cancer. Hugo Dionizy Steinhaus, who ultimately became Banach's doctoral
    advisor, remarked ``Banach was my greatest scientific
    discovery.''} space $c_0$, as
    \begin{equation}
        \sum_{n=0}^{\infty}
        \norm[\big]
        {
            e_n
            \cdot
            z^n
        }_\infty
        =
        \sum_{n=0}^{\infty}
        \abs{z}^n
        =
        \frac{1}{1 - \abs{z}}
        <
        \infty.
    \end{equation}
    Note that the $n$-th term is induced by the $n$-linear mapping
    \begin{equation}
        L_n
        \bigl(
            z_1,\ldots,z_n
        \bigr)
        \coloneqq
        z_1 \cdots z_n \cdot e_n
    \end{equation}
    and thus $f$ is
    Gâteaux holomorphic by Corollary~\ref{cor:GateauxPowerSeries}. Moreover, if $z_1,
    \ldots, z_n \in \disk$ are distinct, then the set $\{f(z_1),\ldots,f(z_n)\}
    \subseteq c_0$ is
    linearly independent. Hence, the image of any smaller open disk
    $\gls{DiskRadiusR} = \{z \in \C
    \colon \abs{z} < r\}$ with $0 < r < 1$ is not contained within any finite dimensional
    subspace of $c_0$.
\end{example}

Ruling out this behaviour facilitates a proof of the continuity of a composition. To avoid
complications regarding the openness of domains of definitions, we work with globally
defined functions.
\begin{lemma}
    \label{lem:GateauxComposition}
    Let $f \in \Holomorphic_G(V_1,V_2)$ and $g \in \Holomorphic_G(V_2,V_3)$ be Gâteaux
    holomorphic maps between locally convex spaces $V_1, V_2$ and $V_3$. Assume
    moreover that for every finite dimensional subspace $F_1 \subseteq V_1$, the
    restriction $f \at{F_1}$ maps into a finite dimensional subspace~$F_2$
    of~$V_2$. Then $g \circ f \in \Holomorphic_G(V_1,V_3)$. Moreover, the pullback
    \begin{equation}
        \label{eq:GateauxPullback}
        \gls{Pullback}
        \colon
        \Holomorphic_G(V_2,V_3)
        \longrightarrow
        \Holomorphic_G(V_1,V_3), \quad
        f^*g
        \coloneqq
        g \circ f
    \end{equation}
    with $f$ is a linear and continuous mapping.
\end{lemma}
\begin{proof}
    We use the reformulation of Gâteaux holomorphy from Remark~\ref{rem:Hartog}. To
    this end, let $F_1 \subseteq V_1$ be finite dimensional. By assumption, there
    exists a corresponding finite dimensional subspace $F_2 \subseteq V_2$ with $f
    \at{F_1} \subseteq F_2$. Hence,
    \begin{equation}
        g \circ f
        \at{F_1}
        =
        g
        \at{F_2}
        \circ
        f
        \at{F_1}
    \end{equation}
    is holomorphic by virtue of $f \in \Holomorphic_G(V_1,V_2)$ and $g \in
    \Holomorphic_G(V_2,V_3)$. This proves
    \begin{equation}
        g \circ f \in \Holomorphic_G(V_1,V_3).
    \end{equation}
    The continuity of \eqref{eq:GateauxPullback} follows from the fact that $f$
    maps compact subsets of finite dimensional subspaces of $V_1$ to compact
    subsets of finite dimensional subspaces of $V_2$. Here, we use that
    restrictions of $f$ to finite dimensional subspaces of $V_1$ are, in
    particular, continuous as holomorphic maps.
\end{proof}

We take a closer look two particular cases with additional structure.
\begin{corollary}
    \label{cor:GateauxPullbackLinear}
    Let $V_1, V_2, V_3$ be locally convex spaces and $L \colon V_1 \longrightarrow V_2$ be linear. Then the pullback
    \begin{equation}
        \label{eq:GateauxPullbackLinear}
        L^*
        \colon
        \Holomorphic_G(V_2,V_3)
        \longrightarrow
        \Holomorphic_G(V_1,V_3), \quad
        L^*f
        \coloneqq
        f \circ L
    \end{equation}
    with $L$ is a well defined and continuous mapping.
\end{corollary}
\begin{proof}
    As a linear mapping, $L \at{F_1}$ has finite dimensional image for any finite
    dimensional subspace
    $F_1 \subseteq V_1$. Hence, all claims follow from
    Lemma~\ref{lem:GateauxComposition}. The well definedness
    of~\eqref{eq:GateauxPullbackLinear} may also be readily verified directly.
    Indeed, we have
    \begin{equation}
        \varphi
        \bigl(
            L^*f
            (v_0 + zv)
        \bigr)
        =
        \varphi
        \bigl(
        f
        (Lv_0 + zLv)
        \bigr)
    \end{equation}
    for $f \in \Holomorphic_G(V_2,V_3)$, $v_0,v \in V_1$, $z \in \C$ and $\varphi \in
    V_3'$. Hence, Gâteaux differentiability of $f$ at $Lv_0$ implies Gâteaux
    differentiability of $L^*f$ at $v_0$.
\end{proof}

As a yet more particular case, we may consider $L \colon V \longrightarrow V$ to be
the translation by some fixed $v_0 \in V$, i.e. $Lv \coloneqq v + v_0$. This allows
one to restrict many considerations about Gâteaux holomorphic functions to
considerations around the origin. In particular, one gets to utilize
\eqref{eq:NHomogeneity} freely when dealing with homogeneous polynomials. This matches
with \eqref{eq:PolynomialTaylorSeries}, which states that the pullback $L^*P$ of an
$n$-homogeneous polynomial $P$ is a linear combination of homogeneous
polynomials of degrees $0, \ldots, n$. Refining our methods, we may moreover
establish the following variant of Corollary~\ref{cor:GateauxPullbackLinear} for
\emph{continuous
polynomials} and general domains of definitions.
\begin{proposition}
    \label{prop:PullbackPolynomials}
    Let $V_1, V_2$ and $V_3$ be locally convex spaces, $P \in \Pol(V_1,V_2)$ be a
    continuous polynomial and $U \subseteq V_2$ be open. Then the pullback
    \begin{equation}
        \label{eq:GateauxPullbackPolynomial}
        P^*
        \colon
        \Holomorphic_G(U,V_3)
        \longrightarrow
        \Holomorphic_G
        \bigl(
            P^{-1}(U),V_3
        \bigr), \quad
        P^*f
        \coloneqq
        f \circ P
    \end{equation}
    with $P$ is a  well defined and continuous mapping.
\end{proposition}
\begin{proof}
    By continuity of $P$, the preimage $P^{-1}(U) \subseteq V_1$ is open and so we may
    speak of the locally convex space $\Holomorphic_G(P^{-1}(U),V_3)$ in the first
    place. Let $F_1 \subseteq V_1$ be a finite dimensional subspace. Then
    \begin{equation}
        P
        \bigl(
            P^{-1}(U) \cap F_1
        \bigr)
        \subseteq
        U \cap P(F_1)
        \subseteq
        P(F_1),
    \end{equation}
    where $P(F_1)$ is contained in a finite-dimensional subspace $F_2$ of $V_2$,
    see again Proposition~\ref{prop:Polynomials},
    \ref{item:PolynomialsFiniteDimensions}. Thus, if $f \in \Holomorphic_G(U, V_3)$,
    then the composition
    \begin{equation*}
        P^* f = P \circ f
        \colon
        F_1 \cap P^{-1}(U)
        \longrightarrow
        V_3
    \end{equation*}
    is holomorphic. As $F_1$ was arbitrary, this proves $P^*f \in
    \Holomorphic_G(P^{-1}(U),V_3)$. Taking another look at the seminorms
    \eqref{eq:GateauxSeminorms}, the continuity of $P^*$ is now clear.
\end{proof}

We have gathered all ingredients to prove the existence of Taylor expansions for
Gâteaux holomorphic functions $f \in \Holomorphic_G(U,V)$. The idea is that
restrictions of $f$ to finite dimensional subspaces do possess such expansions by the
finite dimensional theory. Our derivation is notably different from the one presented
\cite[Sec.~3.1]{dineen:1999a} and avoids using some of the heavier functional analytic
machinery such as Grothendieck's\footnote{Alexander Grothendieck (1928-2014)
    was, among many other astonishing achievements, the father of modern algebraic
    geometry. He was awarded the Fields medal in 1966. After a dispute over military
    funding, he withdrew from the public, while still pursuing both mathematical and
    non-mathematical writing, and spent his final years secluded in a small village
    at the foot of the Pyrenees.} Theorem on equicontinuity.
\begin{theorem}[Taylor Expansions of Gâteaux holomorphic functions]
    \index{Taylor!Gâteaux}
    \index{Gâteaux!Taylor expansion}
    \index{Taylor!Polynomials}
    \label{thm:GateauxTaylor} \;\\
    Let $V,W$ be locally convex spaces, $U \subseteq V$ open, $v_0 \in U$ and $f \in
    \Holomorphic_G(U,W)$. Denote the completion
    of $W$ by \gls{Completion}.\footnote{If $W$ is not Hausdorff, we pass to its
    Hausdorffization before completing.}
    \begin{theoremlist}
        \item \label{item:GateauxTaylorExpansion}
        There exist unique $n$-homogeneous polynomials $P_n \colon V
        \longrightarrow \widehat{W}$ for all $n \in \N_0$ and a balanced\footnote{For all $v \in
        B$ and $t \in \R$, we have $e^{\I t} \cdot v \in B$. Also referred to as circled in the
        literature and some authors demand the condition for multiplications with any complex
        number $z \in \disk^\cl$.} zero neighbourhood $U_0 \subseteq V$ such that for all $v \in
        U_0$
        \begin{equation}
            \label{eq:GateauxTaylor}
            f(v_0 + v)
            =
            \sum_{n=0}^{\infty}
            P_n(v).
        \end{equation}
        \item \label{item:GateauxTaylorDomain}
        If $U = U_0 + v_0$ for some balanced zero neighbourhood $U_0
        \subseteq V$, then \eqref{eq:GateauxTaylor} holds for all $v \in U_0$.
        \item The polynomials $P_n$ are explicitly given by the vector-valued Cauchy integrals
        \begin{equation}
            \label{eq:GateauxTaylorPolynomials}
            \index{Gâteaux!Taylor polynomials}
            \index{Gâteaux!Cauchy integral formula}
            P_n(v)
            =
            \frac{1}{2\pi \I}
            \int_{\boundary \mathbb{D}}
            \frac{f(v_0 + zv)}{z^{n+1}}
            \D z
            \in
            \widehat{W}
        \end{equation}
        for $n \in \N_0$ and $v \in V$ with $v_0 + zv \in U$ for all $z \in \boundary
        \mathbb{D}$.
        \item \label{item:GateauxTaylorConvergence}
        The series \eqref{eq:GateauxTaylor} converges in
        $\Holomorphic_G(U_0,\widehat{W})$,
         i.e. locally uniformly on the intersections of finite dimensional subspaces
         with $U_0$.
    \end{theoremlist}
\end{theorem}
\begin{proof}
    By our considerations around Corollary~\ref{cor:GateauxPullbackLinear}, it
    suffices to work with $v_0 = 0$.\footnote{This will not simplify any of our
    arguments, but at least ease the pain of bookkeeping.} Moreover, we use local
    convexity of $V$
    to replace $U$ with a balanced neighbourhood of zero. We begin by establishing the
    uniqueness of the expansion. It suffices to prove that, assuming \eqref{eq:GateauxTaylor}
    holds and $P_n(v) \neq 0$ for some $n \in \N_0$ and $v \in V$, then $f$ is not the zero
    function. Indeed, as $P_n(v) \neq 0$, there exists a $\varphi \in W'$ such that
    $\varphi(P_n(v)) \neq 0$. Consider the auxiliary function $h \colon \C \longrightarrow \C$,
    \begin{equation}
        h(z)
        \coloneqq
        \varphi
        \bigl(
            f(zv)
        \bigr)
        =
        \sum_{k=0}^{\infty}
        \varphi
        \bigl(
            P_k(zv)
        \bigr)
        =
        \sum_{k=0}^{\infty}
        z^k
        \cdot
        \varphi
        \bigl(
            P_k(v)
        \bigr),
    \end{equation}
    where the series converges by virtue of Gâteaux differentiability of $f$ at
    $v_0$ and the
    continuity of~$\varphi$. Taking another look, we notice that $h$ is a power series
    whose $n$-th coefficient is ${\varphi(P_n(v)) \neq 0}$. Hence, $h$ is not the zero
    function, and thus the same is true for $f$.

    We turn towards proving that choosing $P_n$ as in
    \eqref{eq:GateauxTaylorPolynomials} yields $n$-homogeneous polynomials such that
    \eqref{eq:GateauxTaylor} holds. To this end, we fix $v \in V$ as well as
    $\varphi \in W'$ and check that the integrals
    \eqref{eq:GateauxTaylorPolynomials}
    actually converge in the completion $\widehat{W}$ of $W$. Indeed, if $\seminorm{q}
    \in \cs(W)$, then we may estimate
    \begin{equation*}
        \seminorm{q}
        \biggl(
            \int_{\boundary \mathbb{D}}
            \frac{f(zv)}{z^n}
            \D z
        \biggr)
        \le
        2 \pi
        \cdot
        \max_{z \in \boundary \mathbb{D}}
        \seminorm{q}
        \bigl(
            f(zv)
        \bigr)
        <
        \infty,
    \end{equation*}
    as $f$ is continuous on the one-dimensional subspace spanned by $v$. Sequential
    completeness of $\widehat{W}$ thus implies the convergence of
    \eqref{eq:GateauxTaylorPolynomials} as vector-valued Riemann integrals. By
    absolute convexity of $U$, this yields functions
    \begin{equation}
        P_n \colon U \longrightarrow \widehat{W}, \quad
        P_n(v)
        \coloneqq
        \frac{1}{2\pi \I}
        \int_{\boundary \mathbb{D}}
        \frac{f(zv)}{z^{n+1}}
        \D z
    \end{equation}
    The substitution $w \coloneqq \lambda \cdot z$ and path independence of the integral yields
    \begin{equation}
        P_n(\lambda v)
        =
        \frac{1}{2\pi \I}
        \int_{\boundary \mathbb{D}}
        \frac{f(\lambda z \cdot v)}{z^{n+1}}
        \D z
        =
        \frac{\lambda^{n}}{2 \pi \I}
        \int_{\lambda \cdot \boundary \mathbb{D}}
        \frac{f(w\cdot v)}{w^{n+1}}
        \D w
%        =
%        \frac{\lambda^{n}}{2 \pi \I}
%        \int_{\boundary \mathbb{D}}
%        \frac{f(w\cdot v)}{w^{n+1}}
%        \D w
        =
        \lambda^n
        \cdot
        P_n(v)
    \end{equation}
    for all $\lambda \in \C$ with $\lambda v \in U$ and we extend $P_n$ to $V$ by enforcing
    this homogeneity, i.e. by setting
    \begin{equation}
        P_n(\lambda v)
        \coloneqq
        \lambda^n
        \cdot
        P_n(v)
        \qquad
        \textrm{for all $\lambda \in \C$ and $v \in U$.}
    \end{equation}
    As $U$ is open and thus in particular absorbing, this indeed yields an
    $n$-homogeneous extension $P_n \colon V \longrightarrow \widehat{W}$. If $F
    \subseteq V$ is a finite dimensional subspace, then the restriction~$P_n \at{F}$
    is still $n$-homogeneous and, as $f$ is holomorphic on $F$, also smooth. Euler's
    Theorem on homogeneity thus implies that $P_n$ is a homogeneous polynomial of
    degree $n$. That is, by Corollary~\ref{cor:PolynomialInducing} there exists a
    unique corresponding symmetric $n$-linear mapping
    \begin{equation}
        L_F
        \coloneqq
        \widecheck{Q}_{n,F}
        \colon
        F^n \longrightarrow \widehat{W},
    \end{equation}
    which we extend by zero to $V^n$. This yields a net $(L_F)_F$ of multilinear mappings
    $V^n \longrightarrow W$ indexed by the finite dimensional subspaces of $V$ endowed
    with the direction of upwards set inclusion. If
    \begin{equation}
        F \cap F'
        \neq
        \emptyset,
        \qquad
        \textrm{then}
        \quad
        L_F
        \at[\Big]{F \cap F'}
        =
        L_{F'}
        \at[\Big]{F \cap F'}
    \end{equation}
    by construction. Thus, this net is Cauchy within
    the space
    of all $n$-linear mappings, endowed with the topology of pointwise convergence. It
    converges to some multilinear mapping $L \colon V \longrightarrow W$, which by
    construction induces $P_n$. Thus, $P_n$ is a polynomial itself. By Gâteaux holomorphy,
    we moreover get the locally uniformly convergent power series expansion
    \begin{equation}
        \varphi
        \bigl(
            f(zv)
        \bigr)
        =
        \sum_{n=0}^{\infty}
        \frac{z^n}{2\pi \I}
        \int_{\boundary \mathbb{D}}
        \frac{\varphi(f(wv))}{w^{n+1}}
        \D w
        =
        \varphi
        \biggl(
            \sum_{n=0}^{\infty}
            \frac{z^n}{2\pi \I}
            \int_{\boundary \mathbb{D}}
            \frac{f(wv)}{w^{n+1}}
            \D w
        \biggr)
    \end{equation}
    for any $\varphi' \in W'$ and $z \in U_v$, see again
    \eqref{eq:GateauxMaximalNeighbourhood} with $v_0 \coloneqq 0$. As $W'$
    separates the points of the completion $\widehat{W}$ and by the Hahn-Banach
    Theorem, this implies
    \begin{equation}
        f(zv)
        =
        \sum_{n=0}^{\infty}
        \frac{z^n}{2\pi \I}
        \int_{\boundary \mathbb{D}}
        \frac{f(wv)}{w^{n+1}}
        \D w
        =
        \sum_{n=0}^{\infty}
        z^n
        \cdot
        P_n(v)
    \end{equation}
    for $z \in U_v$, and the series converges locally uniformly. By construction,
    \eqref{eq:GateauxTaylor} holds on all of $U$. By our choice of $U$, we have thus proved
    \ref{item:GateauxTaylorDomain}. If finally $F
    \subseteq V$ is a finite dimensional subspace, then \eqref{eq:GateauxTaylor}
    reduces
    to the usual vector-valued Taylor expansion of ${f \at{F} \colon F \longrightarrow W}$
    by construction. In particular, it converges locally uniformly.
\end{proof}

In the process, we have proved the following ``duck type''\footnote{If it looks like a duck and quacks like a duck, then it is a duck.} lemma for polynomials.
\begin{corollary}
    \label{cor:GateauxNHomogeneousDuck}
    \index{Euler's Theorem}
    \index{Gâteaux!Euler's Theorem}
    Let $V,W$ be locally convex spaces and $U \subseteq V$ an open neighbourhood of zero.
    Assume that $f \in \Holomorphic_G(U,V)$ is $n$-homogeneous in the sense of
    \eqref{eq:NHomogeneity} for some integer $n \in \N_0$. Then $f$ is a polynomial of
    homogeneity $n$.
\end{corollary}

Plugging in $v = 0$ in \eqref{eq:GateauxTaylor} and
\eqref{eq:GateauxTaylorPolynomials} implies the validity of the mean value property for
Gâteaux holomorphic functions.
\begin{corollary}[Mean value property]
    \index{Mean value property}
    \index{Gâteaux!Mean value property}
    Let $V,W$ be locally convex spaces, $U \subseteq V$ be open and $f \in
    \Holomorphic_G(U,W)$. Then
    \begin{equation}
        \label{eq:MeanValue}
        f(v_0)
        =
        \frac{1}{2 \pi}
        \int_{0}^{2\pi}
        f
        \bigl(
            v_0 + r e^{\I t}v
        \bigr)
        \D t
        =
        \frac{1}{2 \pi \I}
        \int_{\boundary \disk}
        \frac{f(v_0 + zv)}{z}
        \D z
    \end{equation}
    for any $r > 0$ and $v \in V$ such that $v_0 + r e^{\I t} v \in U$ for all $t
    \in [0,2\pi]$.
\end{corollary}

We call the expansion \eqref{eq:GateauxTaylor} the \emph{Taylor expansion} of $f$
around $v_0$ and the $n$-homogeneous polynomials \eqref{eq:GateauxTaylorPolynomials}
the \emph{Taylor polynomials} of $f$ around $v_0$. Whenever necessary, we will use the
more precise notation \gls{TaylorPolynomial} instead, but stick to the less
cumbersome $P_n$
whenever we may. Taking another look at \eqref{eq:GateauxTaylorPolynomials}, we note that
$P_n$ takes values in the \emph{first sequential} completion of
$W$\footnote{Which needs not at all be sequentially complete!}, i.e. one does not
need to pass to the full completion. As for holomorphic functions in finite
dimensions, our construction shows
that the natural domains of convergence only depend on the geometry of the domain and not
the particular function.\footnote{Analytic continuation is, of course, a much more
complicated endeavour.} We have already encountered this phenomenon in
Lemma~\ref{lem:GateauxZeroNeighbourhoods}.

 The case of $U = V$ is particularly nice, as the Taylor expansion simply
 represents $f$ globally. Following esteemed tradition, we call functions $f \in
 \Holomorphic_G(U,V)$ \emph{entire} Gâteaux holomorphic functions. Taking into
 account Corollary~\ref{cor:GateauxPowerSeries} there are the following equivalent
 characterizations of entire Gâteaux holomorphic functions.
\begin{corollary}[Entire Gâteaux holomorphic functions]
    \index{Entire!Gâteaux holomorphic functions}
    \index{Gâteaux!Entire function}
    \label{cor:GateauxEntire}
    Let $V,W$ be locally convex spaces and $f \colon V \longrightarrow W$ be a
    mapping. Then the following are equivalent:
    \begin{corollarylist}
        \item The function $f$ is an entire Gâteaux holomorphic function, i.e. $f \in
        \Holomorphic_G(V,W)$.
        \item For every $v_0 \in V$ there exist unique $n$-homogeneous polynomials
        $P_n \colon V \longrightarrow W$ for all~$n \in \N_0$ such that
        \begin{equation}
            \label{eq:GateauxTaylorEntire}
            f(v_0 + v)
            =
            \sum_{n=0}^{\infty}
            P_n(v)
            \qquad
            \textrm{for all }
            v \in V.
        \end{equation}
        \item There exists some $v_0 \in V$ and unique $n$-homogeneous polynomials
        $P_n \colon V \longrightarrow W$ for all $n \in \N_0$ such that
        \eqref{eq:GateauxTaylorEntire} holds for all $v \in V$.
        \item There exist unique $n$-homogeneous polynomials $P_n \colon V \longrightarrow W$ for all $n \in \N_0$ such that for all $v \in V$.
        \begin{equation}
            f(v)
            =
            \sum_{n=0}^{\infty}
            P_n(v).
        \end{equation}
    \end{corollarylist}
\end{corollary}

Taking a closer look at \eqref{eq:GateauxTaylorPolynomials}, we get the following
locally convex version of the Cauchy estimates, which will be instrumental for our
considerations in Section~\ref{sec:RTopologiesPolynomial} and
Section~\ref{sec:EntireVectors}:
\begin{proposition}[Locally convex Cauchy estimates, {\cite[(3.12)]{dineen:1999a}}]
    \index{Cauchy estimates!Locally convex}
    \label{prop:CauchyEstimatesLocallyConvex}
    Let $V,W$ be locally convex spaces, $U \subseteq V$ be open and $f \in
    \Holomorphic_G(U,W)$ with Taylor expansion \eqref{eq:GateauxTaylor} around $v_0
    \in U$. Then, for every balanced subset $B \subseteq V$, we have
    \begin{equation}
        \label{eq:CauchyEstimatesLocallyConvex}
        \sup_{v \in B}
        \seminorm{q}
        \bigl(
            P_n(v)
        \bigr)
        \le
        r^{-n}
%        \cdot
        \sup_{v \in v_0 + rB}
        \seminorm{q}
        \bigl(
            f(v)
        \bigr)
    \end{equation}
    for all $\seminorm{q} \in \cs(W)$, $n \in \N_0$ and $r >0$ with $v_0 + rB \subseteq U$.
\end{proposition}
\begin{proof}
    We use the Cauchy integrals \eqref{eq:GateauxTaylorPolynomials} to estimate
    \begin{align}
        \sup_{v \in B}
        \seminorm{q}
        \bigl(
            P_n(v)
        \bigr)
        &\le
        \sup_{v \in B}
        \frac{1}{2 \pi r^n}
        \int_0^{2\pi}
        \seminorm{q}
        \bigl(
            f(v_0 + r e^{\I t} v)
        \bigr)
        \D t \\
        &\le
        r^{-n}
%        \cdot
        \sup_{v \in B}
        \sup_{t \in [0,2\pi]}
        \seminorm{q}
        \bigl(
            f(v_0 + r e^{\I t} v)
        \bigr) \\
        &\le
        r^{-n}
%        \cdot
        \sup_{v \in v_0 + rB}
        \seminorm{q}
        \bigl(
        f(v)
        \bigr)
    \end{align}
    for all $n \in \N_0$ and $r > 0$ with $v_0 + rB \subseteq U$.
\end{proof}

The locally convex Cauchy estimates facilitate several continuity statements regarding
the passage from holomorphic functions and expansion points towards Taylor
polynomials and their values. A first instance is the following:
\begin{corollary}
    \label{cor:GateauxTaylorContinuity}
    Let $V,W$ be locally convex spaces such that $V$ is Hausdorff, $U \subseteq V$ be
    absolutely convex, $v_0 \in U$ and $n \in \N_0$. Then the linear mapping
    \begin{equation}
        P_{n,v_0}
        \colon
        \Holomorphic_G(U,W)
        \longrightarrow
        \Holomorphic_G(V,W)
    \end{equation}
    is continuous with respect to the topology of locally uniform convergence on finite
    dimensional subspaces.
\end{corollary}
\begin{proof}
    Let $F \subseteq V$ be a finite dimensional subspace. Without loss of generality, we
    may assume $v_0 \in F$. Then $F \cap U \subseteq F$ is an open neighbourhood of
    $v_0$ and thus contains a compact ball $\Ball_r(v_0)^\cl$ with respect to some
    auxiliary norm. The locally convex Cauchy estimates
    \eqref{eq:CauchyEstimatesLocallyConvex} then yield
    \begin{equation}
        \sup_{v \in \Ball_r(0)}
        \seminorm{q}
        \bigl(
            P_{n,v_0}(v)
        \bigr)
        \le
        \sup_{v \in \Ball_r(v_0)}
        \seminorm{q}
        \bigl(
            f(v)
        \bigr)
    \end{equation}
    for all $f \in \Holomorphic_G(U,W)$. If now $K \subseteq F$ is an arbitrary compact
    set, then $K \subseteq \Ball_R(0)$ for some $R > 0$ and thus
    \eqref{eq:NHomogeneity} and what we have already proved yield
    \begin{equation}
        \sup_{v \in K}
        \seminorm{q}
        \bigl(
        P_{n,v_0}(v)
        \bigr)
        \le
        \frac{R^n}{r^n}
        \sup_{v \in \Ball_r(v_0)}
        \seminorm{q}
        \bigl(
        f(v)
        \bigr).
    \end{equation}
    Variation of $F$ completes the proof.
\end{proof}

Having established the existence of Taylor expansions, we investigate how one can
relate
Gâteaux holomorphy to the existence of suitable differential quotients. Indeed, geometrically,
Gâteaux holomorphy corresponds to the existence of all \emph{directional derivatives} in
every direction with complex homogeneous dependence on the direction. The following
adaptation of \cite[Lemma~3.3]{dineen:1999a} makes this precise in order one.
\begin{corollary}[Gâteaux Differentiability vs. Directional Derivatives]
    \index{Gâteaux!vs. directional derivatives}
    \index{Directional derivative}
    \index{Derivative!Directional}
    \label{cor:GateauxVsDirectional} \; \\
    Let $V,W$ be locally convex spaces, $U \subseteq V$ be open and $f \colon U
    \longrightarrow W$ a mapping. Then $f$ is Gâteaux holomorphic if and only if the limit
    \begin{equation}
        \label{eq:GateauxVsDirectional}
        \lim_{z \rightarrow 0}
        \frac{f(v_0 + zv) - f(v_0)}{z}
    \end{equation}
    exists for all $v_0 \in U$ and $v \in V$ in the completion \gls{Completion} of
    $W$. In this case, Taylor expanding
    \begin{equation}
        f(v_0 + v)
        =
        \sum_{n=0}^{\infty}
        P_n(v)
    \end{equation}
    around $v_0$, we get
    \begin{equation}
        \label{eq:DirectionalDerivativeVsTaylor}
        P_1(v)
        =
        \lim_{z \rightarrow 0}
        \frac{f(v_0 + zv) - f(v_0)}{z}
        \qquad
        \textrm{for all }
        v \in V.
    \end{equation}
\end{corollary}
\begin{proof}
    Assume that the limit \eqref{eq:GateauxVsDirectional} exists in the completion
    $\widehat{W}$ of $W$ for all $v_0 \in U$ and $v \in V$. We fix $v_0 \in U$. Then
    for every $\varphi \in W'$, the limit
    \begin{equation}
        \lim_{z \rightarrow 0}
        \varphi
        \biggl(
        \frac{f(v_0 + zv) - f(v_0)}{z}
        \biggr)
        =
        \varphi
        \biggl(
        \lim_{z \rightarrow 0}
        \frac{f(v_0 + zv) - f(v_0)}{z}
        \biggr)
    \end{equation}
    exists by continuity. That is, the function \eqref{eq:Gateaux} is complex
    differentiable at zero. If now $\lambda \in \C$ is such that $v_0 + \lambda v \in
    U$, then replacing $v_0$ with $v_0 + z \lambda$ yields the complex
    differentiability of
    \begin{equation}
        z
        \mapsto
        f\bigl(
        (v_0 + \lambda v)
        +
        zv
        \bigr)
        =
        f
        \bigl(
        v_0 + (z + \lambda)v
        \bigr)
    \end{equation}
    at the origin, which implies complex differentiability of \eqref{eq:Gateaux}
    at
    $\lambda$. Varying $\lambda$ as above yields holomorphy of \eqref{eq:Gateaux} in
    an open neighbourhood of zero by openness of $U$ and thus~$f$ is Gâteaux
    differentiable at $v_0$.

    Assume conversely that $f$ is Gâteaux holomorphic and Taylor expand $f$ around
    $v_0$ with Taylor polynomials $P_n \colon V \longrightarrow \widehat{W}$
    such that \eqref{eq:GateauxTaylor} holds. Plugging the expansion into the
    difference quotient yields then, if $v \in V$ and $z \in \C$ with $v_0 + zv
    \in U$,
    \begin{equation}
        \frac{f(v_0 + zv) - f(v_0)}{z}
        =
        \frac{1}{z}
        \sum_{n=1}^{\infty}
        P_n(zv)
        =
        P_1(v)
        +
        \sum_{n=2}^{\infty}
        z^{n-1}
        P_n(v)
        \overset{z \rightarrow 0}{\longrightarrow}
        P_1(v),
    \end{equation}
    where the convergence is locally uniform on intersections of $U$ with finite dimensional subspaces of $V$. In particular, the limit exists pointwisely, i.e. within $\widehat{W}$.
\end{proof}

If $W$ is complete and Hausdorff, the limit \eqref{eq:GateauxVsDirectional} is unique and it is
customary to write
\begin{equation}
    \gls{GateauxDerivative}
    \coloneqq
    \lim_{z \rightarrow 0}
    \frac{f(v_0 + zv) - f(v_0)}{z}
    \qquad
    \textrm{for all $v_0 \in U$ and $v \in V$.}
\end{equation}
One calls $df(v_0, \argument)$ the \emph{Gâteaux derivative} of $f$ at $v_0$. By
\eqref{eq:DirectionalDerivativeVsTaylor}, Gâteaux
derivatives of functions between complex locally convex spaces always have a linear
dependence on the direction. This was first established in \cite{zorn:1945a} in
the setting of Banach spaces.
\index{Zorn, Max}
Our ``duck lemma'' Corollary~\ref{cor:GateauxNHomogeneousDuck} applied to the
Taylor polynomials \eqref{eq:GateauxTaylorPolynomials} may be viewed as the higher
order generalization of this observation.
\begin{example}[A non-linear real Gâteaux derivative]
    Notably, if one were to consider real differentiability instead, Gâteaux
    derivatives need not be linear, already in finite dimensions. Indeed, consider
    \begin{equation}
        F \colon \R^2 \longrightarrow \R, \quad
        F(x,y)
        \coloneqq
        \begin{cases}
            \frac{x^3}{x^2 + y^2}
            \; &\textrm{for} \;
            (x,y) \neq 0, \\
            0
            \; &\textrm{for} \;
            (x,y) = 0.
        \end{cases}
    \end{equation}
    At the origin, we have
    \begin{equation}
        \frac{F(tx,ty) - F(0,0)}{t}
        =
        \frac{1}{t}
        \cdot
        \frac{(tx)^3}{(tx)^2 + (ty)^2} - 0
        =
        \frac{x^3}{x^2 + y^2}
    \end{equation}
    for $(x,y) \neq 0$ and $t \neq 0$, which clearly fails to be linear with respect to $(x,y)$.
\end{example}

\begin{remark}[Gâteaux holomorphy and Continuity]
    \index{Gâteaux!Vs. Continuity}
    We provide an alternative proof of Proposition~\ref{prop:Gateaux},
    \ref{item:GateauxPolynomials}. Indeed, combining \eqref{eq:PolynomialsDifference}
    with Corollary~\ref{cor:GateauxVsDirectional} yields
    \begin{equation}
        \label{eq:PolynomialDerivative}
        dP(v_0,v)
        =
        \frac{1}{z}
        \sum_{k=0}^{n-1}
        \binom{n}{k}
        \widecheck{P}
        \bigl(
        \underbrace{v_0,\ldots,v_0}_{k\text{-times}},
        \underbrace{zv,\ldots,zv}_{(n-k)\text{-times}}
        \bigr)
        \at[\bigg]{z=0}
        =
        n
        \widecheck{P}
        \bigl(
        v_0,
        \underbrace{v, \ldots, v}_{(n-1)\text{-times}}
        \bigr)
    \end{equation}
    as the Gâteaux derivative of $P$. As in real analysis of several variables,
    existence of directional derivatives does not necessarily imply continuity, a
    drastic instance of which follows from combining Proposition~\ref{prop:Gateaux},
    \ref{item:GateauxPolynomials} with
    Corollary~\ref{cor:PolynomialContinuityVsMultilinear}: Every discontinuous
    multilinear symmetric mapping induces discontinuous polynomials, which are
    nevertheless Gâteaux holomorphic.
\end{remark}

This completes our study of Gâteaux holomorphy and motivates the definition of Fréchet
holomorphy: We simply add the assumption of continuity.

\section{Fréchet Holomorphic Functions}
\label{sec:HolomorphicFrechet}
\epigraph{\begin{flushright}
        It’s still magic even if \\ you know how it’s done.
    \end{flushright} \, \vspace{-0.5cm}}{\emph{A Hat Full of
Sky} -- Terry Pratchett}
% !TeX root = ../Dissertation.tex

Our presentation is guided by the textbooks \cite[Ch.~2-3]{dineen:1981a},
\cite[Sec.~7-8]{kriegl.michor:1997a}, \cite[Sec.~3.1]{dineen:1999a}, again with
the added self-imposed enrichment of working with arbitrary open domains.
\begin{definition}[Fréchet Holomorphy]
    \index{Holomorphy!Fréchet}
    \index{Fréchet!Differentiability}
    \index{Fréchet!Holomorphy}
    \label{def:Frechet}
    Let $V,W$ be locally convex spaces, $U \subseteq V$ be open and $f \colon U
    \longrightarrow W$ a mapping.
    \begin{definitionlist}
        \item \label{item:FrechetDifferentiability}
        The map $f$ is called \emph{Fréchet differentiable} at $v_0 \in U$ if $f$ is Gâteaux
        differentiable and continuous at $v_0$.
        \item \label{item:FrechetHolomorphy}
        The map $f$ is called \emph{Fréchet holomorphic} if it is Fréchet differentiable at
        all
        $v_0 \in U$.
        \item The set of all Fréchet holomorphic mappings from $U$ to $W$ is denoted by
        \begin{equation}
            \label{eq:FrechetHolomorphicSet}
            \gls{HolomorphicFrechet}
            \coloneqq
            \bigl\{
            f \colon U \longrightarrow W
            \colon
            f \text{ is Fréchet holomorphic}
            \bigr\}
            =
            \Holomorphic_G(U,W)
            \cap
            \Continuous(U,W),
        \end{equation}
        where \gls{Continuous} denotes the set of continuous functions between
        topological spaces $X$ and $Y$.
    \end{definitionlist}
\end{definition}

As $\Holomorphic(U,W)$ is defined as an intersection, we immediately get the following analogue of Proposition~\ref{prop:Gateaux} for Fréchet holomorphic functions.
\begin{proposition}
    \label{prop:Frechet}
    \index{Fréchet!Algebraic}
    Let $V,W$ be locally convex spaces and $U \subseteq V$ be open.
    \begin{propositionlist}
        \item \label{item:FrechetVectorSpace}
        The set $\Holomorphic(U,W)$ is a vector space with respect to the pointwise operations.
        \item \label{item:FrechetAlgebra}
        If $\algebra{A}$ is a locally convex algebra, then the set
        $\Holomorphic(U,\algebra{A})$ is a locally convex algebra with respect to the
        pointwise operations and the subspace topology inherited from
        the continuous mappings $\Continuous(U,W)$.
        \item \label{item:FrechetSheaf}
        Sending open subsets $U' \subseteq V$ to $\Holomorphic(U',W)$ yields a sheaf $\Holomorphic(\argument,W)$ of vector spaces with restrictions given by the restriction of functions.
        \item \label{item:FrechetPolynomials}
        We have
        \begin{equation}
            \label{eq:FrechetPolynomials}
            \Pol^\bullet(V,W)
            \subseteq
            \Holomorphic(V,W).
        \end{equation}
    \end{propositionlist}
\end{proposition}

Taking another look at the Cauchy integrals
\eqref{eq:GateauxTaylorPolynomials}, we get the following first regularity result for Taylor
polynomials of Fréchet holomorphic functions:
\begin{corollary}[Continuity of Taylor polynomials]
    \label{cor:TaylorFrechetContinuity}
    \index{Fréchet!Taylor polynomials continuity}
    Let $V$ and $W$ be locally convex spaces, $U \subseteq V$ be open and $f
    \in \Holomorphic(U,W)$. Then the Taylor polynomials of $f$ are continuous as
    mappings
    \begin{equation}
        P_{n,\bullet}
        \colon
        U
        \times
        V
        \longrightarrow
        \widehat{W}
    \end{equation}
    with respect to the product topology and all $n \in \N_0$. In particular,
    $P_{n,v_0} \in \Pol(V,\widehat{W})$ for all $v_0 \in U$.
\end{corollary}
\begin{proof}
    Let $U_0$ be the absolutely convex zero neighbourhood corresponding to the Taylor
    expansion of $f$ around $v_0$ as in Theorem~\ref{thm:GateauxTaylor},
    \ref{item:GateauxTaylorExpansion} and fix $n \in \N$. Using the continuity of
    addition, we find
    another absolutely convex zero neighbourhood $U_0'$ such that $U_0' + U_0'
    \subseteq U_0$. We define
    \begin{equation}
        U_0''
        \coloneqq
        v_0
        +
        U_0'.
    \end{equation}
    Invoking \eqref{eq:GateauxTaylorPolynomials} for the expansion point $v_0' \in
    U_0''$, we may write
    \begin{equation}
        P_{n,v_0'}(v)
        =
        \frac{1}{2 \pi \I}
        \int_{\boundary \mathbb{D}}
        \frac{f(v_0' + zv)}{z^{n+1}}
        \D z
        \qquad
        \textrm{for all }
        v \in U_{0}',
    \end{equation}
    as
    \begin{equation}
        v_0' + zv
        \in
        U_0' + U_0'
        \subseteq U_0
        \qquad
        \textrm{for all }
        z \in \boundary \mathbb{D}.
    \end{equation}
    In view of the continuity of $f$, this implies continuity of
    \begin{equation}
        P_{n,\bullet}
        \colon
        U_0''
        \times
        U_0'
        \longrightarrow
        \widehat{W}.
    \end{equation}
    But $P_{n,\bullet}$ is polynomial with respect to its second argument and thus
    Proposition~\ref{prop:PolynomialContinuity} coupled with variation of $v_0$ yields
    the claim.
\end{proof}

Under which conditions pointwisely convergent series of continuous polynomials yield a
Fréchet holomorphic limit is much trickier, as the Taylor series
\eqref{eq:GateauxTaylor} derived from Gateaux holomorphy has no reason to converge
uniformly on any \emph{open} zero neighbourhood. We will return to this problem in
Corollary~\ref{cor:FrechetPowerSeries} and Proposition~\ref{prop:Bastiani}, once we
have established more technology. The latter provides a converse to
Corollary~\ref{cor:TaylorFrechetContinuity}, where the crucial assumption is the
continuity with respect to the \emph{expansion} point.

Combining the continuity of Taylor polynomials for Fréchet differentiable mappings
with Theorem~\ref{thm:GateauxTaylor} and Corollary~\ref{cor:GateauxPowerSeries}
moreover allows us to prove that Fréchet differentiability is, up to some technical
subtleties, stable under composition.
\begin{theorem}
    \index{Fréchet!Compositions}
    \index{Chain rule}
    \label{thm:FréchetCompositions}
    Let $V_1,V_2,V_3$ be locally convex spaces, $U_1
    \subseteq V_1$ be open and $v_0 \in U$. Let moreover $f \colon U_1 \longrightarrow
    V_2$ be Fréchet differentiable at~$v_0$, $U_2 \subseteq V_2$ an open neighbourhood
    of~$f(v_0)$ and $g \in \Holomorphic(U_2,V_3)$. Then the composition $g \circ f$ is
    Fréchet differentiable at $v_0$.
\end{theorem}
\begin{proof}
    By Theorem~\ref{thm:GateauxTaylor} there exist balanced neighbourhoods of zero
    $U_{0,1} \subseteq V_1$ and~$U_{0,2} \subseteq V_2$ as well as polynomials
    $P_n
    \colon V_1 \longrightarrow \widehat{V}_2$ and $Q_n \colon V_2 \longrightarrow
    \widehat{V}_3$ for all $n \in \N_0$ such that
    \begin{equation}
        f(v_0 + v)
        =
        \sum_{n=0}^\infty
        P_n(v)
        \quad \textrm{and} \quad
        g
        \bigl(
        f(v_0) + w
        \bigr)
        =
        \sum_{n=0}^\infty
        Q_n(w)
    \end{equation}
    for $v \in U_{0,1}$ and $w \in U_{0,2}$. By
    Proposition~\ref{prop:PolynomialContinuity}, each of the $Q_n$ is actually
    locally uniformly continuous and thus maps Cauchy nets to Cauchy nets.
    Consequently, it
    admits a unique locally uniformly continuous extension to a polynomial $Q_n
    \colon
    \widehat{V}_2 \longrightarrow \widehat{V}_3$ and this in turn yields an extension
    of $\widecheck{Q}_n$ to $\widehat{V}_2 \times \cdots \times \widehat{V}_2$
    for all $n \in \N_0$. The auxiliary function
    \begin{equation}
        h
        \colon
        U_{0,1}
        \longrightarrow
        V_2, \quad
        h(v)
        \coloneqq
        f(v_0 + v) - f(v_0)
        =
        \sum_{n=1}^{\infty}
        P_n(v)
    \end{equation}
    is continuous at zero by continuity of $f$ at $v_0$ and thus $U_0 \coloneqq
    U_{0,1} \cap h^{-1}(U_{0,2})$ constitutes an open neighbourhood of $v_0$.
    Composing the Taylor expansions of $f$ and $g$ then yields
    \begin{equation}
        \bigl(
        g \circ f
        \bigr)(v_0 + v)
        =
        g
        \bigl(
        f(v_0 + v)
        \bigr)
        =
        g
        \biggl(
        f(v_0)
        +
        \sum_{n=1}^{\infty}
        P_n(v)
        \biggr)
        =
        \sum_{m=0}^{\infty}
        Q_m
        \biggl(
        \sum_{n=1}^{\infty}
        P_n(v)
        \biggr)
    \end{equation}
    for $v \in U_{0} \subseteq U_{0,1}$, as $h(v) = \sum_{n=1}^{\infty} P_n(v) \in
    U_{0,2}$ by construction. Note that, for fixed~$v \in U_0$, this series
    converges
    \emph{unconditionally}, as it arises from a Taylor series in finite dimensions.
    By Corollary~\ref{cor:TaylorFrechetContinuity}, each of the polynomials $Q_n$ is
    continuous, and by Corollary~\ref{cor:PolynomialContinuityVsMultilinear} the same is
    true for the corresponding multilinear symmetric mapping $\widecheck{Q}_n$.
    Consequently,
    \begin{equation}
        Q_m
        \biggl(
        \sum_{n=1}^{\infty}
            P_n(v)
        \biggr)
        =
        \sum_{k_1, \ldots, k_n = 1}^{\infty}
        \widecheck{Q}_m
        \bigl(
            P_{k_1}(v),
            \ldots,
            P_{k_n}(v)
        \bigr)
        \qquad
        \textrm{for all }
        v \in U_0.
    \end{equation}
    By the Polarization identity \eqref{eq:Polarization} and the fact that the
    composition of two polynomials is once again a polynomial, each of the
    compositions
    \begin{equation}
        v
        \mapsto
        \widecheck{Q}_m
        \bigl(
        P_{k_1}(v),
        \ldots,
        P_{k_n}(v)
        \bigr)
    \end{equation}
    is thus a linear combination of homogeneous polynomials itself. By unconditional
    convergence, we may thus rearrange the series into one indexed by their degrees,
    which still converges pointwisely. Each term is then a pointwisely convergent
    sequence of homogeneous polynomials and thus converges to a homogeneous
    polynomial itself. By Corollary~\ref{cor:GateauxPowerSeries},
    this in turn implies that $g \circ f$ is Gateaux differentiable at $v_0$. Finally,
    as the continuity of $g \circ f$ at $v_0$ is clear, this completes the proof.
\end{proof}

In the situation that $f$ is a polynomial, i.e. the setting of
Proposition~\ref{prop:PullbackPolynomials}, the complexity of the argument reduces
significantly. Indeed, then the inner sum is \emph{finite}, and thus one does not
need to reorder the series.

Our next goal is a sufficient criterion for Fréchet holomorphy based on a
boundedness condition, which generalizes
Lemma~\ref{lem:PolynomialContinuityVsBoundedness} to Fréchet holomorphic
functions. Before establishing the abstract result, we convince ourselves that general
Fréchet holomorphic functions need, unlike continuous polynomials, not be bounded.
\begin{example}[An unbounded holomorphic function]
    \label{ex:HolomorphicUnbounded}
    \index{Fréchet!Unbounded}
    Let $V \coloneqq \gls{AbsolutelySummable}$ be the space of absolutely
    summable sequences indexed by $\N_0$. We write $e_n \in \ell^1$ for the
    sequences with $\gls{UnitSequence}(k) \coloneqq \delta_{n,k}$ for $n,k \in
    \N_0$. The dual
    vectors
    \begin{equation}
        \varphi_n
        \colon
        \ell^1 \longrightarrow \C, \quad
        \varphi_n(a)
        \coloneqq
        a_n
    \end{equation}
    are continuous linear functionals with operator norm $\norm{\varphi_n} = 1$ for
    all $n \in \N_0$. Define
    \begin{equation}
        \label{eq:HolomorphicUnbounded}
        f \colon \ell^1 \longrightarrow \C, \quad
        f(a)
        \coloneqq
        \sum_{n=0}^\infty
        \varphi_n(a)^n.
    \end{equation}
    Note that $f$ is well defined, as for $a \in \ell^1$, there exists an index $N \in \N_0$ with $\abs{a_n} \le 1$ for all $n \ge N$. Consequently,
    \begin{equation*}
        \abs[\big]
        {f(a)}
        \le
        \sum_{n=0}^\infty
        \abs{a_n}^n
        \le
        \sum_{n=0}^{N-1}
        \abs{a_n}^n
        +
        \sum_{n=N}^\infty
        \abs{a_n}
        \le
        \sum_{n=1}^{N-1}
        \abs{a_n}^n
        +
        \norm{a}_1
        <
        \infty.
    \end{equation*}
    By construction, the Taylor polynomials of $f$ at the origin are given by
    $P_n = \varphi_n^n$, each of which is $n$-homogeneous and continuous. Thus $f
    \in \Holomorphic_G(\ell^1,\C)$ by Corollary~\ref{cor:GateauxPowerSeries}.
    However, $f$ is unbounded on every ball with radius $R > 1$. Indeed, if $1 <
    r < R$, then
    \begin{equation*}
        \abs[\big]
        {f(r e_n)}
        =
        r^n
        \qquad
        \textrm{for all }
        n \in \N_0,
    \end{equation*}
    which is unbounded for $n \rightarrow \infty$. Conversely, the function $f$ is bounded
    by $1$ on the closed unit ball. This also illustrates that -- in contrast to
    polynomials -- a
    holomorphic function may behave wildly different on a set $S$ and its rescalings $rS$.

    Another related observation is that the Taylor series
    \eqref{eq:HolomorphicUnbounded} converges uniformly on the open ball
    $\Ball_r(0) \subseteq
    \ell^1$ for any $r \le 1$, but not on any ball with radius $R > 1$. Indeed, every
    sequence $a \in \Ball_r(0)$ with $r \le 1$ fulfils $\abs{a_n} < r$ for all $n \in
    \N_0$. This yields the estimate
    \begin{equation}
        \abs[\bigg]
        {
            \sum_{n=N}^{\infty}
            \varphi_n(a)^n
        }
        \le
        \sum_{n=N}^{\infty}
        \abs[\big]{a_n}^n
        <
        \sum_{n=N}^{\infty}
        r^n
        =
        \frac{r^N}{1 - r}
        \overset{N \rightarrow \infty}{\longrightarrow}
        0
    \end{equation}
    for the remainder. As a uniform limit of continuous and bounded functions, our function $f$ is thus continuous and bounded on $\Ball_1(0)$ for abstract reasons. That is, $f \in \Holomorphic(\Ball_1(0),\C)$. We will later see that $f$ is Fréchet holomorphic on all of $\ell^1$. Conversely, the unboundedness of $f$ on $\Ball_R(0)$ for all $R > 1$ implies that the Taylor series does not converge uniformly there. Regarding the situation on the unit sphere, we note
    \begin{equation}
        \abs[\bigg]
        {
            f(e_{N+1})
            -
            \sum_{n=0}^{N}
            a_n(e_{N+1})^n
        }
        =
        1
    \end{equation}
    for any $N \in \N_0$. In particular, the convergence is not uniform on the closed unit ball.
\end{example}

We proceed with a more or less obvious Lemma to develop an intuition on how to
establish a notion of ``local boundedness''.
\begin{lemma}[Local Boundedness]
    \index{Boundedness!Local}
    \label{lem:LocalBoundedness}
    Let $V,W$ be locally convex spaces.
    \begin{lemmalist}
        \item \label{item:LocalBoundednessVariants}
        Let $U \subseteq V$ be open and $f \colon U \longrightarrow W$ be such that there exists an open set $U' \subseteq U$ such that $f(U') \subseteq W$ is bounded. Then the restriction $f \at{U'} \colon U' \longrightarrow W$ is bounded.
        \item \label{item:LocalBoundednessImpliesBoundednessPolynomial}
        Let $P \colon V \longrightarrow W$ be a polynomial such that there exists an
        open zero neighbourhood~$U$ within $V$ such that the restriction $P \at{U}$ is
        bounded. Then $P$ is bounded.
    \end{lemmalist}
\end{lemma}
\begin{proof}
    Take objects as described in \ref{item:LocalBoundednessVariants} and let $B \subseteq U'$ be bounded. Then also the image $f(B) \subseteq f(U') \subseteq W$ is bounded as a subset of a bounded set. For \ref{item:LocalBoundednessImpliesBoundednessPolynomial} we may additionally assume that $P$ is $n$-homogeneous in the sense of \eqref{eq:NHomogeneity}. If now $B \subseteq V$ is bounded, then there exists some $r > 0$ such that $rB \subseteq U$ and clearly $rB$ is still bounded. Hence,
    \begin{equation}
        P(B)
        =
        \frac{1}{r^n}
        P(rB)
    \end{equation}
    is bounded by boundedness of $P \at{U}$. We have proved boundedness of $P$.
\end{proof}

In the sequel, we call a function $f \colon U \longrightarrow W$ \emph{locally
bounded at $v_0 \in U$} if the more restrictive condition in
\ref{item:LocalBoundednessVariants} is fulfilled, i.e. $v_0$ has an open
neighbourhood $U' \subseteq U$ such that the image $f(U') \subseteq W$ is
bounded. Analogously, $f$ is called \emph{locally bounded} if it is locally
bounded at all $v_0 \in U$. This yields the following sufficient criterion for
Fréchet holomorphy, which can be found in \cite[Prop.~3.7]{dineen:1999a} and
formalizes some of our observations from Example~\ref{ex:HolomorphicUnbounded}.
\begin{proposition}
    \label{prop:FrechetVsLocalBoundedness}
    Let $V,W$ be locally convex spaces, $U \subseteq V$ open and $f \colon U \longrightarrow W$.
    \begin{propositionlist}
        \item \label{item:FrechetVsLocalBoundedness}
        If $f \in \Holomorphic_G(U,W)$ is locally bounded, then it is Fréchet holomorphic.
        \item \label{item:FrechetVsLocalBoundednessNormed}
        If $W$ is normed and $f$ is continuous\footnote{Say, by virtue of
            Fréchet differentiability of $f$ at $v_0$.} at $v_0 \in U$, then $f$
            locally bounded at $v_0$.
    \end{propositionlist}
\end{proposition}
\begin{proof}
    Let $f$ be locally bounded and $v_0 \in U$. By assumption, we find an open neighbourhood $U' \subseteq U$ of zero such that $f(v_0 + U') \subseteq W$ is bounded. By local convexity, we may assume that $U'$ is absolutely convex. That is, for every $\seminorm{q} \in \cs(W)$ we have
    \begin{equation}
        C_{\seminorm{q}}
        \coloneqq
        \sup_{v \in U'}
        \seminorm{q}
        \bigl(
            f(v_0 + v)
        \bigr)
        <
        \infty.
    \end{equation}
    We prove that the Taylor series \eqref{eq:GateauxTaylor} of $f$ around $v_0$
    converges uniformly on $v_0 + U'$, i.e. within the complete space $\Continuous(v_0 +
    U',\widehat{W})$. Invoking the locally convex Cauchy
    estimates~\eqref{eq:CauchyEstimatesLocallyConvex}, we get
    \begin{equation}
        \sup_{v \in U'}
        \seminorm{q}
        \bigl(
            P_n(v)
        \bigr)
        \le
        r^{-n}
        \cdot
        \sup_{v \in v_0 + rU'}
        \seminorm{q}
        \bigl(
            f(v)
        \bigr)
        \le
        r^{-n}
        \cdot
        \sup_{v \in U'}
        \seminorm{q}
        \bigl(
            f(v_0 + v)
        \bigr)
        \le
        r^{-n}
        \cdot
        C_{\seminorm{q}}
    \end{equation}
    for $0 < r < 1$ by absolute convexity of $U'$. Hence,
    \begin{equation}
        \sum_{n=0}^{\infty}
        \seminorm{q}
        \bigl(
            P_n(v)
        \bigr)
        \le
        C_{\seminorm{q}}
        \cdot
        \sum_{n=0}^{\infty}
        r^{-n}
        =
        \frac{C_{\seminorm{q}}}{1 - r}
        <
        \infty
        \qquad
        \textrm{for }
        0 < r < 1.
    \end{equation}
    Thus, \eqref{eq:GateauxTaylor} converges in $\Continuous(v_0 + U', \widehat{W})$
    and its limit $f$ is continuous around $v_0$. Variation of $v_0$ proves $f \in
    \Holomorphic(U,V)$.

    Let now conversely $W$ be normed and assume $f$ is continuous at $v_0 \in U$.
    The crucial point is that open norm balls $\Ball_r(v_0)$ are bounded open subsets of
    $W$. By continuity at~$v_0$, the preimage
    \begin{equation}
        U'
        \coloneqq
        f^{-1}
        \bigl(
            \Ball_1(f(v_0))
        \bigr)
    \end{equation}
    is an open neighbourhood of $v_0$ in $U$. Consequently,
    \begin{equation}
        \norm[\big]
        {f(v)}
        \le
        \norm[\big]
        {f(v) - f(v_0)}
        +
        \norm[\big]
        {f(v_0)}
        <
        1
        +
        \norm[\big]
        {f(v_0)}
    \end{equation}
    for all $v \in U'$ and thus $f(U') \subseteq W$ is indeed bounded.
\end{proof}

The condition that $W$ is normed is essential for \ref{item:FrechetVsLocalBoundednessNormed}. Indeed, general locally convex spaces lack open bounded sets\footnote{The existence of an open bounded set is quite close to the topology coming from some norm. We will make this precise in Proposition~\ref{prop:BoundedOpenVsNormability}.} and so the condition may already fail for the \emph{identity mapping}.
\begin{example}[Fréchet holomorphic mappings need not be locally bounded]
    \index{Fréchet!Locally unbounded}
    \label{ex:FrechetLocallyUnbounded}
    For a concrete example, consider the space
    \begin{equation}
        \gls{Sequences}
        =
        \Map(\N,\C)
        =
        \prod_{n \in \N}
        \C
    \end{equation}
    of \emph{all} complex sequences, endowed with the topology of pointwise convergence, i.e. the Cartesian product topology inherited from countably many copies $\C$. A subset
    \begin{equation}
        B
        =
        \prod_{n \in \N}
        B_n
        \subseteq
        \C^\N
    \end{equation}
    is bounded if and only if $B_n \subseteq \C$ is bounded for all $n \in \N$. On the
    other hand, we are dealing with the product topology: A nonempty set
    \begin{equation}
        U
        =
        \prod_{n \in \N}
        U_n
        \subseteq
        \C^\N
    \end{equation}
    is open if and only if almost all\footnote{All but finitely many.} $U_n = \C$
    and the remaining $U_n \subseteq \C$ are open. As $\C$ itself is unbounded,
    these two conditions are clearly incompatible. Hence, the identity mapping
    \begin{equation}
        \id
        \colon
        \C^\N
        \longrightarrow
        \C^\N,
    \end{equation}
    which is Fréchet holomorphic as a continuous linear map by
    Proposition~\ref{prop:Frechet}, \ref{item:FrechetPolynomials}, is not locally bounded.
\end{example}

Using the ideas from Proposition~\ref{prop:FrechetVsLocalBoundedness},
\ref{item:FrechetVsLocalBoundedness} we may prove a variant of
Corollary~\ref{cor:GateauxPowerSeries} for Fréchet holomorphy.
\begin{corollary}
    \label{cor:FrechetPowerSeries}
    \index{Fréchet!Power series}
    Let $P_n \in \Pol^n(V,W)$ be $n$-homogeneous polynomials between locally convex
    spaces $V$ and $W$ for all $n \in \N_0$.
    Assume moreover that for every $\seminorm{q} \in \cs(W)$
    there exists a zero neighbourhood $U \subseteq V$ such that
    \begin{equation}
        \label{eq:FrechetPowerSeriesCondition}
        C
        \coloneqq
        \sup_{n \in \N_0}
        \sup_{v \in U}
        \seminorm{q}
        \bigl(
            P_n(v)
        \bigr)
        <
        \infty.
    \end{equation}
    Then the power series
    \begin{equation}
        \label{eq:FrechetPowerSeries}
        f
        \colon
        U \longrightarrow \widehat{W}, \quad
        f(v)
        \coloneqq
        \sum_{n=0}^{\infty}
        P_n(v)
    \end{equation}
    defines a Fréchet holomorphic function.
\end{corollary}
\begin{proof}
    Let $0 < r < 1$ and note that, by homogeneity,
    \begin{equation}
        \sum_{n=0}^{\infty}
        \sup_{v \in rU}
        \seminorm{q}
        \bigl(
            P_n(v)
        \bigr)
        =
        \sum_{n=0}^{\infty}
        r^n
        \cdot
        \sup_{v \in U}
        \seminorm{q}
        \bigl(
        P_n(v)
        \bigr)
        \le
        \frac{C}{1 - r}.
    \end{equation}
    By assumption, we get such an estimate for all $\seminorm{q} \in \cs(W)$. Thus, by
    continuity of each~$P_n$, this implies that the series
    \eqref{eq:FrechetPowerSeries}
    converges absolutely within the space $\Continuous(rU,\widehat{W})$. Taking the
    limit $r \uparrow 1$ shows that this yields a continuous function $f \colon U
    \longrightarrow \widehat{W}$. As the series~\eqref{eq:FrechetPowerSeries} converges
    in particular pointwisely, Corollary~\ref{cor:GateauxPowerSeries} implies that $f$
    is Gâteaux holomorphic, as well. Thus, $f \in \Holomorphic(U,\widehat{W})$
    as desired.
\end{proof}

Recall that a topological space $X$ is called $T_1$-space\footnote{The $T$ alludes to
the German ``Trennungseigenschaft'', which translates to ``separation axiom'', and was
coined in \cite{tietze:1923a}. The idea
behind numbering them was that axiom $T_k$ implies $T_\ell$ whenever $k \ge \ell$. In
practice, this did not work out all too well due to a number competing inequivalent
definitions. For instance, the modern notion of $T_3$-space does not imply property
$T_2$ without also assuming $T_1$. Moreover, this has lead to amusing notions such as
$T_{2 \; 1/2}$, and the full beauty of the arising taxonomy can be marvelled upon
in the wonderful
\cite[Sec.~I.2]{steen.seebach:1995a}.} if for every pair of
distinct
points~$x,y$ within $X$, there exists a neighbourhood $U$ of $x$ such that $y \notin
U$.
\begin{proposition}[Kolmogorov's normability criterion, {\cite{kolmogorov:1934a}}]
    \index{Kolmogorov, Andrey}
    \index{Boundedness!Open Set}
    \label{prop:BoundedOpenVsNormability}
    Let $V$ be a topological vector space. Then there is a norm that induces the topology of $V$ if and only if it is a $T_1$-space and admits a bounded convex zero neighbourhood.
\end{proposition}
\begin{proof}
    We only prove the non-trivial implication. Let $V$ be $T_1$ and $U \subseteq V$ a
    bounded convex neighbourhood of zero. Consider the balanced hull
    \begin{equation}
        B
        \coloneqq
        \bigl\{
            \lambda v
            \colon
            v \in U, \,
            \lambda \in \C, \,
            \abs{\lambda} = 1
        \bigr\}
        =
        \bigcup_{\abs{\lambda} = 1}
        \lambda \cdot U
    \end{equation}
    of $U$. By continuity of the multiplication with scalars, each of
    the sets $\lambda \cdot U$ is open again, and hence also $B$ is open. By
    Lemma~\ref{lem:BoundednessVsAbsoluteConvexity} the set $B$ is moreover
    bounded, as subsets of bounded sets are bounded. We study the associated
    Minkowski functional
    \begin{equation}
        \label{eq:MinkowskiProof}
        \seminorm{p}_B
        \colon
        V \longrightarrow [0,\infty), \quad
        \seminorm{p}_B(v)
        \coloneqq
        \inf
        \bigl\{
            r > 0
            \colon
            v \in rB
        \bigr\}.
    \end{equation}
    By openness of $B$, the Minkowski functional is well defined. As $B$ is an absolutely
    convex neighbourhood of zero, $\seminorm{p}_B$ is a continuous seminorm by
    \cite[Thm.~3.7]{osborne:2014a}. We check that~$\seminorm{p}_B$ actually
    defines a
    norm. To this end, let $v \in V \setminus \{0\}$. As $V$ is a $T_1$-space, we find a
    zero neighbourhood $U' \subseteq V$ with $v \notin U'$. By boundedness of $B$,
    there moreover exists a radius~$r > 0$ such that $r B \subseteq U'$, i.e. $v \notin r'B$
    for all $0 < r' \le r$. Taking another look at~\eqref{eq:MinkowskiProof},
    this means
    $\seminorm{p}_B(v) \ge r > 0$. Thus $\seminorm{p}_B$ is indeed a norm and the
    closure of $B$ is the associated closed unit ball. It remains to establish that
    $\seminorm{p}_B$ generates the topology of~$V$. By continuity of
    translations, it
    suffices to prove that the open balls $\Ball_r(0)$ associated to~$\seminorm{p}_B$
    constitute a basis of zero neighbourhoods. Indeed, if $\tilde{U} \subseteq V$ is a zero
    neighbourhood, then by boundedness of $B$, there exists a radius $r > 0$ such that
    $rB \subseteq \tilde{U}$. By what we have shown, this implies $\Ball_r(0) \subseteq
    \tilde{U}$. This completes the proof.
\end{proof}

Using Kolmogorov's\footnote{Andrey Nikolaevich Kolmogorov (1903-1987) was a Russian
mathematician and the father of the modern formulation of probability theory by means
of measure spaces. The first of his many major discoveries was an almost everywhere
divergent Fourier series \cite{kolmogorov:1923a}.} normability criterion,
we may infer that bounded subsets of
topological spaces are metrizable.
\begin{corollary}
    \index{Boundedness!Metrizability}
    \label{cor:BoundedIsMetrizable}
    Let $V$ be a topological vector space with the $T_1$-property and $B
    \subseteq V$ be
    bounded. Then $B$, endowed with the subspace topology, is metrizable.
\end{corollary}
\begin{proof}
    Let $W \subseteq V$ be the vector space generated by $B$, endowed with the subspace topology. Then $W$ is still $T_1$ and $B \subseteq W$ is a bounded open subset of $W$. By Proposition~\ref{prop:BoundedOpenVsNormability}, the topological vector space $W$ is thus normable, say with norm $\norm{\argument}$. The restriction of the corresponding metric to $B$ provides the desired metric.
\end{proof}

Our proof of Proposition~\ref{prop:BoundedOpenVsNormability} moreover provides an
explicit formula for a norm if $B$ is absolutely convex, namely its
Minkowski\footnote{Hermann Minkowski (1864-1909) was a Polish geometer. His
introduction of Minkowski spacetime provides an efficient and precise mathematical
framework for special relativity. This, in particular, builds a bridge between physics
and the field
of hyperbolic geometry.}
functional. This is sometimes useful, hence we make it precise here.
\begin{corollary}
    \label{cor:BoundedAbsconvNormableByMinkowski}
    \index{Minkowski!Functional}
    \index{Seminorms!Minkowski}
    \index{Minkowski!Hermann}
    Let $V$ be a topological vector space with the $T_1$-property and an open bounded absolutely convex subset $B \subseteq V$. Then the Minkowski functional
    \begin{equation}
        \label{eq:Minkowski}
        \seminorm{p}_B
        \colon
        V \longrightarrow [0,\infty), \quad
        \seminorm{p}_B(v)
        \coloneqq
        \inf
        \bigl\{
        r > 0
        \colon
        v \in rB
        \bigr\}
    \end{equation}
    associated to $B$ constitutes a norm on $V$, which induces its topology.
\end{corollary}

After this topological detour, we return to the regularly scheduled study of Fréchet
holomorphic functions.
\begin{remark}[Hartogs' Theorem on separate holomorphy: Redux]
    \index{Hartogs' Theorem!Fréchet holomorphy}
    \label{rem:HartogRedux}
    We revisit Remark~\ref{rem:Hartog}. In terms of Gâteaux and Fréchet
    holomorphy, we may rephrase Hartogs' Theorem as
        \begin{equation}
            \label{eq:Hartog}
            \Holomorphic_G(U, \C^m)
            =
            \Holomorphic(U, \C^m)
        \end{equation}
        for any $n,m \in \N$ and open $U \subseteq \C^n$. Remarkably, this still
        holds true if the codomain is infinite dimensional. Early versions of this
        observation are due to Dunford~\cite{dunford:1938a}, who studied the
        problem in terms of the Gâteaux derivative from
        Corollary~\ref{cor:GateauxVsDirectional}.
\end{remark}

\begin{proposition}
    \index{Dunford, Nelson}
    \label{prop:HartogRedux}
    Let $n \in \N$, $U \subseteq \C^n$ open and $W$ a locally convex space. Then
    \begin{equation}
        \Holomorphic_G(U, W)
        =
        \Holomorphic(U, W).
    \end{equation}
\end{proposition}
\begin{proof}
    Let $f \in \Holomorphic_G(U,W)$. We check local boundedness of $f$, as we may
    then infer Fréchet holomorphy by means of
    Proposition~\ref{prop:FrechetVsLocalBoundedness}. To this end, fix $v \in U$
    and consider some norm ball~$U_v \coloneqq \Ball_r(v)$ such that its closure
    $K$ is contained within $U$. By \eqref{eq:Hartog}, the composition $\varphi
    \circ f$ is continuous for all
    $\varphi \in W'$. Consequently, the images
    \begin{equation}
        \bigl(
            \varphi
            \circ
            f
        \bigr)(K)
        \subseteq
        \C
    \end{equation}
    are compact and in particular bounded for all $\varphi \in W'$. By the Hahn-Banach
    Theorem, this is equivalent to the boundedness of $f(K) \subseteq W$. As
    $f(\Ball_r(v)) \subseteq f(K)$, this implies that $f$ is locally bounded at $v$. As $v$
    was arbitrary, we may thus apply Proposition~\ref{prop:FrechetVsLocalBoundedness}
    to complete the proof.
\end{proof}

As a simple application, we note that the function from
Example~\ref{ex:GateauxInfiniteDimensionalRange} is continuous and thus Fréchet
holomorphic. Before endowing $\Holomorphic(U,W)$ with a locally convex topology of
its own, we investigate the connection between Fréchet holomorphy and \emph{total
differentiability}.\footnote{Historically, this is also referred to as Fréchet
differentiability.
We have, however, given this phrase an ostensibly different meaning in
Definition~\ref{def:Frechet}.
In Proposition~\ref{prop:FrechetVsTotal} we shall prove that both approaches are
consistent whenever both concepts are applicable.} Recall that a mapping $f
\colon U
\longrightarrow W$ defined on an open
subset $U$ of a normed space $B$ with values in a locally convex space $W$ is called
totally differentiable at $v_0 \in U$ if there exists a continuous \emph{complex}
linear mapping
$L \colon B \longrightarrow W$ such that
\begin{equation}
    \label{eq:TotalDifferentiability}
    \frac{f(v_0 + v) - f(v_0) - Lv}{\norm{v}}
    \overset{v \rightarrow 0}{\longrightarrow}
    0.
\end{equation}
If one wants to generalize this condition to domains $U$ within locally convex spaces, one runs into the principal problem that seminorms may possess non-trivial kernels. In some spaces, one may choose a continuous auxiliary norm, but these need not exist:
\begin{example}[No continuous norms]
    \index{Auxiliary norms}
    Consider the space $\C^\N$ of all complex sequences from
    Example~\ref{ex:FrechetLocallyUnbounded}, again endowed with the locally convex
    topology of pointwise convergence. Let $\seminorm{p} \in \cs(\C^\N)$ be a continuous
    seminorm. Then, by definition of the topology of $\C^\N$, there exists a finite set
    $\{n_1, \ldots, n_k\} \subseteq \N$ and some $c_1, \ldots, c_k > 0$ such that
    \begin{equation}
        \seminorm{p}(\gamma)
        \le
        c_1
        \cdot
        \abs[\big]{\gamma_{n_1}}
        + \cdots +
        c_k
        \cdot
        \abs[\big]{\gamma_{n_k}}
        \qquad
        \textrm{for all }
        \gamma \in \C^\N.
    \end{equation}
    In particular, if $n \notin \{n_1, \ldots, n_k\}$, then the sequence $e_n$
    with $e_n(m) = \delta_{m,n}$ fulfils
    \begin{equation}
        0
        \le
        \seminorm{p}(e_n)
        \le
        c_1
        \cdot
        \abs[\big]{e_n(n_1)}
        + \cdots +
        c_k
        \cdot
        \abs[\big]{e_n{n_k}}
        =
        \seminorm{p}(\gamma)
        \le
        c_1
        \cdot
        \delta_{n,n_1}
        + \cdots +
        c_k
        \cdot
        \delta_{n,n_k}
        =
        0.
    \end{equation}
    Hence, $e_n \in \gls{Kernel} \seminorm{p}$ and thus $\seminorm{p}$ is not a
    norm.
\end{example}

Without a straightforward generalization of \eqref{eq:TotalDifferentiability} in sight, we nevertheless note the following classical result for Banach spaces, which we may extend to locally convex codomains without too much trouble.
\begin{proposition}[Fréchet Holomorphy vs. Total Differentiability, {\cite{zorn:1946a}}]
    \index{Fréchet!vs. total differentiability}
    \index{Zorn, Max}
    \label{prop:FrechetVsTotal} \;\\
    Let $B$ be normed, $U \subseteq B$ be open, $W$ locally convex and $f \colon
    U \longrightarrow W$ be a mapping. Then~$f$ is Fréchet differentiable at $v_0
    \in U$ if and only if it is totally differentiable at $v_0$.
\end{proposition}
\begin{proof}
    Assume first that $f$ is Fréchet differentiable at $v_0 \in U$, i.e. Gâteaux
    differentiable and continuous at $v_0$. As the statement is local, we may assume
    that $U$ is some open ball around $v_0$. By Theorem~\ref{thm:GateauxTaylor}, we
    then have
    \begin{equation}
        \frac{f(v_0 + v) - f(v_0) - P_1(v)}{\norm{v}}
        =
        \frac{1}{\norm{v}}
        \sum_{n=2}^{\infty}
        P_n(v)
        =
        \sum_{n=2}^{\infty}
        \norm{v}^{n-1}
        \cdot
        P_n
        \biggl(
            \frac{v}{\norm{v}}
        \biggr)
    \end{equation}
    for all $v \in V \setminus \{0\}$ with $v_0 + v \in U$, where $P_n$ are the Taylor
    polynomials of $f$ at~$v_0$. Let~$\seminorm{q} \in \cs(W)$. By continuity of
    $\seminorm{q} \circ f$ at $v_0$, there exists a radius $r > 0$ such that
    \begin{equation}
        \Ball_r(v_0)
        \subseteq U
        \quad \textrm{and} \quad
        C \coloneqq \sup_{v \in \Ball_r(v_0)} \seminorm{q}(v) <
        \infty.\footnote{This is essentially
            Proposition~\ref{prop:FrechetVsLocalBoundedness},
            \ref{item:FrechetVsLocalBoundednessNormed} for the local Banach space
            $\pi_\seminorm{q} \colon W \longrightarrow W_\seminorm{q}$ at
            $\seminorm{q}$
            applied to the composition $\pi_\seminorm{q} \circ f \colon U
            \longrightarrow
            W_\seminorm{q}$.}
    \end{equation}
    Invoking the locally convex Cauchy estimates
    from Proposition~\ref{prop:CauchyEstimatesLocallyConvex}, we get
    \begin{equation}
        P_n
        \biggl(
            \frac{v}{\norm{v}}
        \biggr)
        \le
        r^{-n}
        \cdot
        \sup_{v' \in \Ball_r(v_0)}
        \seminorm{q}
        \bigl(
            f(v')
        \bigr)
        =
        \frac{C}{r^n}
    \end{equation}
    for $v \in \Ball_r(0)$ and $n \in \N_0$. Consequently, if $v \in \Ball_{r}(0)$, then
    \begin{equation}
        \seminorm{q}
        \biggl(
            \frac{f(v_0 + v) - f(v_0) - P_1(v)}{\norm{v}}
        \biggr)
        \le
        \frac{C}{r}
        \sum_{n=1}^{\infty}
        \biggl(
            \frac{\norm{v}}{r}
        \biggr)^n
        =
        \norm{v}
        \cdot
        \frac{C}{r^2}
        \frac{r}{r - \norm{v}}
        \overset{v \rightarrow 0}{\longrightarrow}
        0
        \cdot
        \frac{C}{r^2}
        =
        0.
    \end{equation}
    As this holds for all $\seminorm{q} \in \cs(W)$, we have proved
    \begin{equation}
        \lim_{v \rightarrow 0}
        \frac{f(v_0 + v) - f(v_0) - P_1(v)}{\norm{v}}
        =
        0
    \end{equation}
    and thus $f$ is totally differentiable at $v_0$ with total derivative $P_1$, which is
    indeed linear and continuous by a combination of Theorem~\ref{thm:GateauxTaylor}
    with Corollary~\ref{cor:TaylorFrechetContinuity}.

    Assume conversely that $f$ is totally differentiable at $v_0 \in U$ with total derivative
    \begin{equation}
        L
        \in
        \gls{Linear}
        =
        \Pol^1(V,W).
    \end{equation}
     Then also all restrictions $f \at{F}$ of $f$ to finite dimensional subspaces $F
     \subseteq V$ are totally differentiable with total derivative $L \at{F}$ and thus
     holomorphic. Hence, $f \in \Holomorphic_G(U,W)$ by Remark~\ref{rem:Hartog}. The
     condition \eqref{eq:TotalDifferentiability} moreover implies for any $\seminorm{q} \in
     \cs(W)$ that
    \begin{equation}
        \seminorm{q}
        \bigl(
            f(v_0 + v)
            -
            f(v_0)
        \bigr)
        \le
        \seminorm{q}
        \bigl(
            Lv
        \bigr)
        +
        \seminorm{q}
        \bigl(
            f(v_0 + v) - f(v_0) - Lv
        \bigr)
        \overset{v \rightarrow 0}{\longrightarrow}
        0
    \end{equation}
    by continuity of $L$ and $\norm{v} \le 1$ for $v \rightarrow 0$. Thus the mapping $f$
    is also continuous at $v_0$, that is to say $f \in \Holomorphic(U,W)$.
\end{proof}

Using the same techniques, one may establish that $\widecheck{P_n}$ are
the higher order total derivatives of $f$, i.e. $P_2$ is the total derivative of $v_0
\mapsto P_{1,v_0}$ and so on. The symmetry of the inducing maps $\widecheck{P_n}$
-- or equivalently, that we may pass to
$P_n$ in \eqref{eq:GateauxTaylor} -- reflects the symmetry of mixed higher derivatives
with respect to their ordering. We proceed with another characterization of Fréchet
holomorphy by means of only the first derivative.
\begin{proposition}
    \label{prop:Bastiani}
    \index{Fréchet!vs. Bastiani}
    Let $V,W$ be locally convex spaces such that $W$ is sequentially complete, $U
    \subseteq V$ open and $f \in \Holomorphic_G(U,W)$. Moreover assume that there
    is some some convex open set $U_{0} \subseteq U$ such
    that the first Taylor polynomial is jointly continuous as a mapping
    \begin{equation}
        \label{eq:BastianiContinuityDirectional}
        P_{1,\bullet}
        \colon
        U_{0} \times V
        \longrightarrow
        W.
    \end{equation}
    Then the restriction $f \at{U_{0}}$ is Fréchet holomorphic.
\end{proposition}
In fact, there is a version of this result for functions that are only \emph{real}
Gâteaux differentiable and $P_{1,\bullet}$ is replaced by the continuous
Gâteaux derivative $df \colon U \times V \longrightarrow \widehat{W}$. We
provide a proof that works in this generality, but refrain from making the
additional terminology precise.
\begin{proof}
    Fix $v_0 \in U_0$ and let $v \in U_0$ such that the line segment
    \begin{equation}
        \gls{LineSegment}
        \coloneqq
        \bigl\{
            tv + (1-t)(v - v_0)
            \colon
            t \in [0,1]
        \bigr\}
        \end{equation}
    between $v_0$ and $v$ is contained in $U_0$. We consider the auxiliary
    function
    \begin{equation}
        g
        \colon
        [0,1] \longrightarrow W, \quad
        g(t)
        \coloneqq
        f
        \bigl(
            v_0 + t v
        \bigr),
    \end{equation}
    which is differentiable by virtue of $f \in \Holomorphic_G(U,W)$ with
    derivative
    \begin{equation}
        g'(t)
        =
        \lim_{h \rightarrow 0}
        \frac{g(t + h) - g(t)}{h}
        =
        \lim_{h \rightarrow 0}
        \frac{f(v_0 + tv + hv) - g(v_0 + tv)}{h}
        =
        P_{1,v_0 + tv}(v)
    \end{equation}
    for $t \in [0,1]$. By the vector-valued fundamental theorem of calculus as
    it can be found e.g. in \cite[2.6.~Corollary~(6)]{kriegl.michor:1997a},
    this yields
    \begin{equation}
        f(v_0 + v)
        -
        f(v_0)
        =
        g(1)
        -
        g(0)
        =
        \int_0^1
        g'(t)
        \D t
        =
        \int_0^1
        P_{1,v_0 + tv}(v)
        \D t
    \end{equation}
    and thus
    \begin{equation}
        \label{eq:BastianiProof}
        \seminorm{q}
        \bigl(
            f(v_0 + v)
            -
            f(v_0)
        \bigr)
        \le
        \int_0^1
        \seminorm{q}
        \bigl(
            P_{1,v_0 + tv}(v)
        \bigr)
        \D t.
        \tag{$\heartsuit$}
    \end{equation}
    Now, $P_{1,v_0}(0) = 0$ by linearity, and the continuity of
    \eqref{eq:BastianiContinuityDirectional} implies that given $\epsilon > 0$
    there exists a $\delta > 0$ and an open convex neighbourhood $U_0'
    \subseteq U_0$ such that
    \begin{equation}
        \seminorm{q}
        \bigl(
            P_{1,v_0 + tv}(v)
        \bigr)
        \le
        \epsilon
        \qquad
        \textrm{for all $t \in [0,1]$ and $v \in U_0'$ with $[v_0,v] \subseteq
        U_0'$.}
    \end{equation}
    Plugging this into \eqref{eq:BastianiProof} yields
    \begin{equation}
        \seminorm{q}
        \bigl(
            f(v_0 + v)
            -
            f(v_0)
        \bigr)
        \le
        \epsilon
        \qquad
        \textrm{for all $v \in U_0'$ with $[v_0,v] \subseteq U_0'$.}
    \end{equation}
    As $U_0'$ is an open convex neighbourhood of~$v_0$, this implies
    continuity of $f$ at $v_0$. Variation of $v_0$ completes the proof.
\end{proof}

Along the way, we have established the following well-known version of the
fundamental theorem of calculus for Fréchet holomorphic functions and Riemann
integrals.
\begin{corollary}[Fundamental Theorem of Calculus]
    \index{Fundamental Theorem of Calculus}
    Let $V,W$ be locally convex spaces such that $W$ is sequentially
    complete, $U \subseteq V$ be open and $f \in \Holomorphic(U,W)$. Then
    \begin{equation}
        f(v_0 + v)
        -
        f(v_0)
        =
        \int_0^1
        P_{1,v_0 + tv}(v)
        dt
    \end{equation}
    for all $v_0,v \in U$ with $[v_0,v] \subseteq U$.
\end{corollary}

\begin{remark}[Bastiani holomorphy]
    \index{Bastiani, Andrée}
    \index{Ehresmann, Andrée}
    \index{Holomorphy!Bastiani}
    \index{Open questions!Bastiani Holomorphy in aQFT}
    Recent work such as \cite{brunetti.fredenhagen.lauridsen:2019a,
    hawkins.rejzner.visser:2023a, brunetti.moro:2024a} in algebraic
    quantum
    field theory have identified the notion of \emph{Bastiani\footnote{Andrèe Ehresmann
    (born as Andrèe Bastiani in 1935) is a French mathematician specialized in infinite
    dimensional differential calculus, category theory and differential geometry. Since the
    death of her husband, the differential geometer Charles
    Ehresmann (1905--1979), she has been the director of the mathematical journal
    \emph{Cahiers de Topologie et Géométrie Différentielle Catégoriques}.} smoothness}
    as
    a both sufficiently strong and general notion of differentiability for the
    observable algebras: for a function $f \colon U \longrightarrow W$, one
    demands the existence of directional
    derivatives $d^n f$ of \emph{all} orders and their continuity as mappings
    \begin{equation}
        d^n f
        \colon
        U
        \times
        V^n
        \longrightarrow
        W
    \end{equation}
    with respect to the product topology. In Proposition~\ref{prop:Bastiani}
    we have proved that Bastiani smooth maps are, in particular, continuous
    and fulfil a version of the fundamental theorem of calculus. The idea first
    appeared in \cite{michal:1938a} as $M$-differential and was later refined
    significantly in \cite{bastiani:1964a}. Replacing real with complex
    differentiability to define \emph{Bastiani holomorphy}, a combination of
    Corollary~\ref{cor:TaylorFrechetContinuity} with
    Proposition~\ref{prop:Bastiani} proves that this ostensibly new notion
    actually coincides with Fréchet holomorphy. It would be interesting to
    study whether the functions that appear in the context of algebraic
    quantum field theory admit for holomorphic extensions.
\end{remark}

\section{Uniform Convergence on Bounded Sets}
\label{sec:UniformConvergenceOnBoundedSets}
\epigraph{Multiple exclamation marks,' he went on, shaking his head, 'are a sure sign
of a diseased mind.}{\emph{Eric} -- Terry Pratchett}
% !TeX root = ../Dissertation.tex

In this section, we study the spaces of polynomials of a fixed degree $n \in \N_0$ and
of bounded Fréchet holomorphic functions, each endowed with the topology of uniform
convergence on bounded sets. In Lemma~\ref{lem:PolynomialContinuityVsBoundedness}, we
have seen that every continuous polynomial is bounded and that both notions are
equivalent if the domain is bornological. Taking a step back, we first establish that
the space of all bounded mappings
\begin{equation}
    \index{Bounded mapping!Space}
    \gls{Bounded}
    \coloneqq
    \big\{
        f
        \colon
        U \longrightarrow W
        \;\big|\;
        \forall_{B \subseteq U \textrm{ bounded}}
        \colon
        f(B)
        \subseteq
        W
        \textrm{ is bounded}
    \big\}
\end{equation}
is complete if endowed with the locally convex topology generated by the seminorms
$\seminorm{p}_{B,\seminorm{q}}$ from \eqref{eq:BoundedSeminorms}.\footnote{The author
would like to thank Praful Rahangdale for asking about and discussing this.}
Here, $U \subseteq V$ is some fixed open set of a locally convex space, $B
\subseteq U$ runs through all of its bounded subsets and
$\seminorm{q}$ through $\cs(W)$. The results of this section are well known, but
references are sparse as analysts tend to favour coarser topologies. That being said,
the material we are going to present may also be derived from
Mackey's\footnote{George Whitelaw Mackey (1916-2006) was an American representation
theorist and early noncommutative geometer. His principal insight concerning duality
theory was that different topologies may produce the same topological dual space,
which lead him to the notion of vector bornologies, which focus on bounded sets
instead of open ones. Remarkably, this also provides a suitable framework for
homological algebra on locally convex spaces.} \cite{mackey:1946a}
more algebraic theory of bornologies, a
comprehensive exposition of which is the monograph \cite{hogbe-nlend:1977a}. In view
of our goals, we shall stick to our more analytic approach.
\begin{proposition}
    \label{prop:BoundedCompleteness}
    \index{Bounded mapping!Completeness}
    Let $V,W$ be locally convex spaces and $U \subseteq V$ be open.
    \begin{propositionlist}
        \item \label{item:BoundedCompleteness}
        If $W$ is complete Hausdorff, then so is $\Bounded(U,W)$.
        \item \label{item:BoundedCompletenessSequential}
        If $W$ is sequentially complete Hausdorff, then so is $\Bounded(U,W)$.
    \end{propositionlist}
\end{proposition}
\begin{proof}
    Let $(f_\alpha)_{\alpha \in J} \subseteq \Bounded(U,W)$ be a Cauchy net. As finite
    sets are always bounded, the topology of pointwise convergence is coarser
    than the topology of uniform convergence on bounded subsets. Consequently,
    $\Bounded(U,W)$ is Hausdorff and the net $(f_\alpha)$ converges pointwisely to its
    pointwise limit
    \begin{equation}
        f
        \colon
        U \longrightarrow W, \quad
        f(v)
        \coloneqq
        \lim_{\alpha \in J}
        f_\alpha(v).
    \end{equation}
    We prove that this convergence is actually uniform on bounded subsets of $U$. To
    this end, fix $\epsilon > 0$, $\seminorm{q} \in
    \cs(W)$ and a bounded subset $B \subseteq U$. As $(f_\alpha)_{\alpha \in J}
    \subseteq \Bounded(U,W)$ is a Cauchy net, we find an index $\alpha_0 \in J$ such
    that
    \begin{equation}
        \sup_{v \in B}
        \seminorm{q}
        \bigl(
            f_\alpha(v) - f_\beta(v)
        \bigr)
        =
        \seminorm{p}_{B,\seminorm{q}}
        \bigl(
            f_\alpha - f_\beta
        \bigr)
        \le
        \epsilon
        \qquad
        \textrm{for all }
        \alpha, \beta \later \alpha_0.
    \end{equation}
    Consequently, the continuity of $\seminorm{q}$ implies
    \begin{equation}
        \seminorm{q}
        \bigl(
        f(v) - f_\beta(v)
        \bigr)
        =
        \lim_{\alpha \in J}
        \seminorm{q}
        \bigl(
        f_\alpha(v) - f_\beta(v)
        \bigr)
        \le
        \epsilon
        \qquad
        \textrm{for all }
        v \in B
        \textrm{ and }
        \beta \later \alpha_0.
    \end{equation}
    Taking the supremum over $v \in B$ yields
    \begin{equation}
        \label{eq:BoundedCompletenessProof}
        \seminorm{p}_{B,\seminorm{q}}
        \bigl(
            f
            -
            f_\beta
        \bigr)
        \le
        \epsilon
        \qquad
        \textrm{for all }
        \beta \later \alpha_0.
        \tag{$\diamondsuit$}
    \end{equation}
    This establishes the uniform convergence $f_\alpha \rightarrow f$ on bounded
    subsets of $V$. It remains to show that $f$ is bounded itself. To prove this,
    keep the objects from as above, say for the particular choice $\epsilon \coloneqq
    1$. Then \eqref{eq:BoundedCompletenessProof} gives
    \begin{equation}
        \sup_{v \in B}
        \seminorm{q}
        \bigl(
            f(v)
        \bigr)
        =
        \seminorm{p}_{B,\seminorm{q}}
        \bigl(
            f
        \bigr)
        \le
        \seminorm{p}_{B,\seminorm{q}}
        \bigl(
            f - f_{\alpha_0}
        \bigr)
        +
        \seminorm{p}_{B,\seminorm{q}}
        \bigl(
            f_{\alpha_0}
        \bigr)
        \le
        1
        +
        \seminorm{p}_{B,\seminorm{q}}
        \bigl(
        f_{\alpha_0}
        \bigr)
        <
        \infty
    \end{equation}
    by boundedness of $f_{\alpha_0}$. Varying $B$ and $\seminorm{q}$ proves the
    boundedness of $f$, which establishes~\ref{item:BoundedCompleteness}. Replacing
    the nets with sequences throughout our considerations provides a proof
    of~\ref{item:BoundedCompletenessSequential}.
\end{proof}

Taking another look at the upper part of our argument, we get the following
somewhat technical, but nevertheless useful, lemma.
\begin{lemma}
    \index{Uniform convergence!Lemma}
    Let $V,W$ be locally convex spaces, $U \subseteq V$ be open and assume that $W$ is
    complete Hausdorff. Moreover, let $(f_\alpha)_{\alpha \in J}$ be a net of
    mappings $f_\alpha \colon U \longrightarrow W$ that is Cauchy with respect to the
    topology $\tau$ of uniform convergence on some collection of sets containing the
    finite subsets of $U$. Then the net $(f_\alpha)_{\alpha \in J}$ converges to
    its pointwise limit with respect to $\tau$.
\end{lemma}

Of course, there is an analogous statement for sequences if $W$ is merely sequentially
complete. Two examples for $\tau$ are the topologies of uniform convergence on all
compact or all bounded subsets of $U$. Here we use that finite subsets $\{v_1, \ldots,
v_n\} \subseteq V$ are compact, and thus in particular bounded, regardless of the
topology of $V$. One should also be mindful that, given any collection~$S$ of subsets
of $U$, one typically has to include additional sets to ensure that uniform
convergence on all members of $S$ induces a topology. A~comprehensive list of the
most important constructions within the linear theory can be found in
\cite[Sec.~3.6]{osborne:2014a}.

We have already seen that boundedness and continuity are related for certain types of
maps. In the case of $V$ being a Montel\footnote{Paul Montel (1876-1975) was a French
complex analyst. He introduced and systematized the notion of \emph{normal families}
of holomorphic functions, which notably predates the abstract concept of compactness.
His comprehensive
exposition on the field of complex dynamics as pioneered by Pierre Fatou and Gaston
Julia was the first of its kind.}
space\footnote{Montel's Theorem
holds: every bounded and closed subset of $V$ is compact. Some authors call such
spaces semi-Montel. Roughly speaking, the
Montel property ensures that there is a rich
supply of compact sets, which intertwines topology and bornology.}, we get
relations between both notions for arbitrary globally defined functions.
\begin{proposition}
    \index{Montel!Paul}
    \index{Montel!Boundedness vs. continuity}
    \index{Montel!Uniform convergence}
    \index{Uniform convergence!Montel}
    \label{prop:Montel}
    Let $V,W$ be locally convex spaces such that $V$ is Montel.
    \begin{propositionlist}
        \item \label{item:MontelUniformConvergence}
        The topologies of uniform convergence on bounded and uniform convergence on
        compact subsets of $V$ coincide on the space of bounded mappings
        $\Bounded(V,W)$.
        \item \label{item:MontelContinuousAreBounded}
        We have
        \begin{equation}
            \Continuous(V,W)
            \subseteq
            \Bounded(V,W),
        \end{equation}
        where \gls{Continuous} denotes the set of continuous functions between
        topological spaces $X$ and $Y$.
    \end{propositionlist}
\end{proposition}
\begin{proof}
    Every compact set $K \subseteq V$ is bounded: Indeed,
    \begin{equation}
        \bigcup_{n \in \N}
        \Ball_{\seminorm{p},n}(0)
        =
        V
    \end{equation}
    provides an open covering of $K$ for every $\seminorm{p} \in \cs(V)$. Passing to a
    finite subcover means that there is some $n \in \N_0$ with $K \subseteq
    \Ball_{\seminorm{p},n}(0)$. This proves the boundedness of $K$ by varying
    $\seminorm{p} \in \cs(V)$. Conversely, let $B \subseteq V$ be bounded. The Montel
    property implies that its closure $B^\cl \subseteq V$ is compact. Consequently,
    uniform convergence on all compact subsets of $V$ implies uniform convergence on
    all bounded subsets of $V$. We have shown~\ref{item:MontelUniformConvergence}.

    To see \ref{item:MontelContinuousAreBounded}, let again $B \subseteq V$ be
    bounded
    and $f \in \Continuous(V,W)$. Invoking the Montel property as before, its closure
    $B^\cl$ is compact. Hence, the same is true for its image $f(B^\cl)$ under the
    continuous mapping $f$. By the first part of the proof, this in particular means
    $f(B^\cl)$ is a bounded superset of $f(B)$, proving its boundedness. This
    establishes $f \in \Bounded(V,W)$ and completes the proof.
\end{proof}

The upshot of Proposition~\ref{prop:Montel} is that our results on the
$\beta$-topology on the polynomial functions $\Pol^n(V,W)$ pass to the topology of
uniform convergence on compact subsets under the additional assumption that the domain
$V$ is Montel.
\begin{remark}
    Note that we have restricted ourselves to globally defined bounded mappings in
    Proposition~\ref{prop:Montel}. The boundedness ensures that both topologies are
    locally convex. Having global domains of definitions reflects that a
    compact subset $K$ of $V$ needs not have a compact intersection with an open
    subset $U$. For a simple example, consider the vector space $V \coloneqq \R$ with
    the $U \coloneqq (-1,1)$ and $K \coloneqq [0,2]$. The problem here is, of course,
    that $K$ is not contained within $U$.
\end{remark}

We return to the study of polynomials and note the following continuity estimates
regarding translations on $\Pol^n_\beta(V,W)$, which can be found in
\cite[Lem.~1.10]{dineen:1999a} and will be useful in various contexts.
\begin{proposition}
    \label{prop:PolynomialsTranslations}
    \index{Polynomials!Translations}
    Let $V,W$ be locally convex spaces, $B \subseteq V$ be a bounded and balanced
    subset, $\seminorm{q} \in \cs(W)$ and $P \in \Pol^n(V,W)$ for some $n \in \N_0$.
    Then
    \begin{equation}
        \label{eq:PolynomialsTranslationForward}
        \seminorm{p}_{B,\seminorm{q}}
        (P)
        \le
        \seminorm{p}_{v_0 + B,\seminorm{q}}
        (P)
        \qquad
        \textrm{for all }
        v_0 \in V.
    \end{equation}
    If $B$ is moreover convex and $r > 0$ as well as $v_0 \in V$ are such that $r v_0 \in
    B$, then conversely
    \begin{equation}
        \label{eq:PolynomialsTranslationBack}
        \seminorm{p}_{B + v_0,\seminorm{q}}
        (P)
        \le
        \bigl(
        1 + 1/r
        \bigr)^n
        \cdot
        \seminorm{p}_{B,\seminorm{q}}
        (P).
    \end{equation}
\end{proposition}
\begin{proof}
    As $B$ is balanced and $P$ is $n$-homogeneous, we have on the one hand
    \begin{equation}
        \seminorm{p}_{v_0 + B,\seminorm{q}}
        (P)
        =
        \sup_{v \in B}
        \seminorm{q}
        \bigl(
        P(v_0 + v)
        \bigr)
        =
        \sup_{v \in B, t \in \R}
        \seminorm{q}
        \bigl(
        P(v_0 + e^{\I t} v)
        \bigr)
        =
        \sup_{v \in B, t \in \R}
        \seminorm{q}
        \bigl(
        P(e^{-\I t} v_0 + v)
        \bigr).
    \end{equation}
    On the other hand, utilizing the mean value property \eqref{eq:MeanValue} around $v
    \in B$ in direction~$v_0 \in V$, we get
    \begin{equation}
        P(v)
        =
        \frac{1}{2 \pi \I}
        \int_{\boundary \mathbb{D}}
        \frac{P(zv_0 + v)}{z}
        \D z
    \end{equation}
    and thus
    \begin{equation}
        \seminorm{q}
        \bigl(
            P(v)
        \bigr)
        \le
        \sup_{t \in \R}
        \seminorm{q}
        \bigl(
            P(e^{\I t} v_0 + v)
        \bigr).
    \end{equation}
    Taking the suprema over $B$ on both sides and using the first equality proves
    \eqref{eq:PolynomialsTranslationForward}. For the second part, we additionally
    assume convexity of $B$ and fix $r > 0$ such that $rv_0 \in B$. Then
    \begin{equation}
        v_0 + B
        \in
        B/r
        +
        B
        =
        \bigl(
          1 + 1/r
        \bigr)
        B.
    \end{equation}
    By homogeneity of $P$, we arrive at desired estimate
    \begin{equation}
        \seminorm{p}_{B + v_0,\seminorm{q}}
        (P)
        \le
        \seminorm{p}_{(1+1/r)B,\seminorm{q}}
        (P)
        =
        \bigl(
        1 + 1/r
        \bigr)^n
        \cdot
        \seminorm{p}_{B,\seminorm{q}}
        (P).
    \end{equation}
\end{proof}

\begin{remark}
    \label{rem:PolynomialTranslation}
    One should think of Proposition~\ref{prop:PolynomialsTranslations} as a
    quantitative incarnation of the algebraic fact that a polynomial is uniquely
    determined by its values on any balanced zero neighbourhood.
    Taking another look at our proof, we see that both the estimates
    \eqref{eq:PolynomialsTranslationForward} and \eqref{eq:PolynomialsTranslationBack}
    also hold for unbounded subsets $B$ with seminorms defined as in
    \eqref{eq:BoundedSeminorms}, but may degenerate into vacuous statements about
    infinities.
\end{remark}

Next, we establish completeness of $\Pol^n_\beta(V,W)$ for any fixed $n \in \N_0$
under mild assumptions, slightly generalizing \cite[Prop.~1.30]{dineen:1981a}.
\begin{theorem}
    \label{thm:CompletenessHomogeneousPolynomials}%
    \index{Polynomials!Completeness $\beta$}
    \index{Completeness!Polynomials $\beta$}
    Let $V$ be a bornological locally convex space and $W$
    be complete Hausdorff. Then $\Pol^n_\beta(V,W)$ is complete Hausdorff for every $n
    \in \N_0$.
\end{theorem}
\begin{proof}
    The case $n=0$ is just the completeness of $W$. Thus assume $n \in \N$. By
    Lemma~\ref{lem:PolynomialContinuityVsBoundedness},
    \ref{item:PolynomialContinuityImpliesBoundedness} we may view $\Pol^n_\beta(V,W)$
    as endowed with the subspace topology inherited from the inclusion
    \begin{equation}
        \Pol^n_\beta(V,W)
        \subseteq
        \Bounded(V,W).
    \end{equation}
    As the space of bounded mappings $\Bounded(V,W)$ is complete Hausdorff by
    Proposition~\ref{prop:BoundedCompleteness},~\ref{item:BoundedCompleteness} this
    means that $\Pol^n_\beta(V,W)$ is Hausdorff itself and that it suffices to
    establish the closedness of $\Pol^n_\beta(V,W)$ within $\Bounded(V,W)$. To this
    end, let
    \begin{equation}
        (P_\alpha)_{\alpha \in J}
        \subseteq
        \Pol^n_\beta(V,W)
    \end{equation}
    be a convergent net with limit $P \in \Bounded(V,W)$. As the pointwise convergence
    of $(P_\alpha)_{\alpha \in J}$
    corresponds to pointwise convergence of the corresponding net
    $(\widecheck{P}_\alpha)$ of $n$-linear mappings by
    Proposition~\ref{prop:Polarization}, it is clear that $P$ is an $n$-homogeneous
    polynomial itself. Finally, invoking
    Lemma~\ref{lem:PolynomialContinuityVsBoundedness},
    \ref{item:PolynomialBoundednessImpliesContinuity}, which we may apply by the
    assumption that $V$ is bornological, the boundedness of $P$ implies its continuity.
\end{proof}

Dropping the assumption that $V$ is bornological, our proof shows that the spaces of
\emph{bounded} $n$-homogeneous polynomials are complete Hausdorff for all $n \in
\N_0$. In other words, having a bornological domain is only required to pass back from
boundedness to continuity. Invoking Proposition~\ref{prop:Montel},
\ref{item:MontelUniformConvergence}, we get the completeness of $\Pol^n(V,W)$ also
with respect to the topology of uniform convergence on compact subsets under the
additional assumption that $V$ is Montel.

This completes our considerations on polynomials. In
Example~\ref{ex:FrechetLocallyUnbounded}, we have seen that boundedness typically
singles out a proper subspace of the space of
all Fréchet holomorphic functions. We make this additional assumption, as the
functions of interest in Section~\ref{sec:RTopologiesPolynomial} are always bounded.
Once again, our goal is to establish completeness, this time under suitable
assumptions. We begin by fixing some notation.
\begin{definition}[Bounded Fréchet holomorphic functions]
    \label{def:HolomorphicBounded}
    \index{Bounded mapping!Holomorphic}
    \index{Fréchet!Bounded}
    Let $V,W$ be locally convex spaces and $U \subseteq V$ be open. We write
    \begin{equation}
        \gls{HolomorphicBounded}
        \coloneqq
        \bigl\{
        f \in \Holomorphic(U,W)
        \colon
        f
        \; \textrm{is bounded}
        \bigr\}
        =
        \Holomorphic(U,W)
        \cap
        \Bounded(U,W)
    \end{equation}
    for the space of bounded Fréchet holomorphic functions endowed with the subspace
    topology inherited from $\Bounded(U,W)$.
\end{definition}

This topology is, of course, the topology of uniform convergence on bounded
subsets. By Proposition~\ref{prop:Montel}, the following is immediate.
\begin{corollary}
    Let $V,W$ be locally convex spaces such that $V$ is Montel, and let $U
    \subseteq V$ be open. Then
    \begin{equation}
        \Holomorphic_\beta(U,V)
        =
        \Holomorphic(U,V).
    \end{equation}
\end{corollary}

As indicated in the end of \cite[Sec.~3.1]{dineen:1999a}, there is an alternative
defining set of seminorms for
$\Holomorphic_\beta(V,W)$ based on summability of the Taylor polynomials at the origin.
We generalize this to functions defined on open subsets $U \subseteq V$. As the
correspondence is ultimately based on the existence of Taylor expansions
\eqref{eq:GateauxTaylor} and the locally convex Cauchy estimates
\eqref{eq:CauchyEstimatesLocallyConvex}, we introduce a more precise notation for
their typical domains of convergence. Indeed, by what we have
shown\footnote{After carefully reverting the additional assumption of $v_0 = 0$.}, the
Taylor series
\begin{equation}
    f(v_0 + v)
    =
    \sum_{n=0}^{\infty}
    P_n(v)
\end{equation}
of some $f \in \Holomorphic_G(U,W)$ converges for all $v \in U$ such that the line segment
\begin{equation}
    [v_0, v]
    \coloneqq
    \bigl\{
        v_0 + t(v-v_0)
        \colon
        t \in [0,1]
    \bigr\}
\end{equation}
between $v_0$ and $v$ is contained in $U$. That is, on any starlike subset of $U$ with
starpoint~$v_0$. As we are interested in bounded sets, we consider
\begin{equation}
    \label{eq:Polydisks}
    \index{Polydisks}
    \gls{Polydisks}
    \coloneqq
    \bigl\{
        B \subseteq V
        \colon
        B \;
        \textrm{is an absolutely convex bounded set with} \;
        v_0 + B
        \subseteq
        U
    \bigr\}
\end{equation}
for any open $U \subseteq V$ and $v_0 \in U$. Indeed, if $B \in \Bdd(U,v_0)$ is
non-empty, then $B$ is starlike with starpoint zero by absolute convexity, and
thus $v_0 + B$ is starlike with starpoint $v_0$. Note that, by
Lemma~\ref{lem:BoundednessVsAbsoluteConvexity} there is finally always a wealth
of absolutely convex bounded sets. One can regard the
collection \eqref{eq:Polydisks} as a
locally convex adaptation of polydisks. The locally convex Cauchy estimates
\eqref{eq:CauchyEstimatesLocallyConvex} now yield the following:
\begin{proposition}
    \label{prop:BetaAlternativeSeminorms}
    \index{Seminorms!$\beta$ power series}
    \index{Gâteaux!Boundedness}
    \index{Fréchet!Bounded}
    Let $V,W$ be locally convex spaces, $U \subseteq V$ be open and $f \in \Holomorphic_G(U,W)$ with Taylor polynomials $P_{n,v_0}$ at $v_0 \in U$.
    \begin{propositionlist}
        \item \label{item:AlternativeBounded}
        The function $f$ is bounded if and only if
        \begin{equation}
            \gls{SeminormsBoundedAlternative}(f)
            \coloneqq
            \sum_{n=0}^{\infty}
            \seminorm{p}_{B, \seminorm{q}}
            \bigl(
                P_{n,v_0}
            \bigr)
            <
            \infty
        \end{equation}
        for all $v_0 \in U$, $B \in \Bdd(U,v_0)$ and $\seminorm{q} \in \cs(W)$.
        \item \label{item:AlternativeBoundedAreSeminorms}
        The functions $\seminorm{r}_{B, \seminorm{q}, v_0}$ are seminorms on
        $\Holomorphic_\beta(U,W)$ for all vectors $v_0 \in U$, bounded sets $B
        \in \Bdd(U,v_0)$ and continuous seminorms $\seminorm{q} \in \cs(W)$.
        \item \label{item:BoundedAlternativeSeminorms}
        The locally convex topology generated by the system
        \begin{equation}
            \label{eq:BoundedAlternativeSeminorms}
            \bigl\{
            \seminorm{r}_{B, \seminorm{q}, v_0}
            \colon
            \seminorm{q} \in \cs(W), \,
            v_0 \in U, \,
            B \in \mathfrak{B}(U,v_0)
            \bigr\}
        \end{equation}
        is the $\Holomorphic_\beta(U,W)$ topology. More precisely, we have the estimates
        \begin{equation}
            \seminorm{r}_{B, \seminorm{q}, v_0}
            \le
            \frac{\seminorm{p}_{v_0 + rB, \seminorm{q}}}{1 - r}
            \qquad \textrm{and} \qquad
            \seminorm{p}_{B, \seminorm{q}}
            \le
            \seminorm{r}_{B, \seminorm{q}}
        \end{equation}
        for any $B \in \Bdd(U,v_0)$ and $r > 0$ with $v_0 + rB \subseteq U$.
        \item \label{item:BoundedAlternativeEntireSeries}
        If $U = V$, then the seminorms corresponding to $\Bdd(U,v_0)$ in
        \eqref{eq:BoundedAlternativeSeminorms} generate the topology of
        $\Holomorphic_\beta(U,W)$ for every $v_0 \in V$.
        \item \label{item:BoundedAlternativeEntire}
        If $U = V$, then the seminorms corresponding to $\Bdd(U,v_0)$ in
        \eqref{eq:BoundedSeminorms} generate the topology of
        $\Holomorphic_\beta(U,W)$
        for every $v_0 \in V$.
    \end{propositionlist}
\end{proposition}
\begin{proof}
    Assume first that $f$ is bounded and let $v_0 \in U$, $B \in \Bdd(U,v_0)$ and $\seminorm{q} \in \cs(W)$. Applying the locally convex Cauchy estimates \eqref{eq:CauchyEstimatesLocallyConvex}, we get
    \begin{equation}
        \seminorm{p}_{B, \seminorm{q}}
        (P_{n,v_0})
        \le
        \frac{\seminorm{p}_{v_0 + rB, \seminorm{q}}(f)}{r^n}
    \end{equation}
    for any $r > 0$ with $v_0 + rB \subseteq U$ and $n \in \N_0$. This yields the estimate
    \begin{equation}
        \seminorm{r}_{B, \seminorm{q}, v_0}(f)
        =
        \sum_{n=0}^{\infty}
        \seminorm{p}_{B, \seminorm{q}}
        (P_{n,v_0})
        \le
        \frac{\seminorm{p}_{v_0 + rB, \seminorm{q}}(f)}{1 - r}.
    \end{equation}
    This implies $\seminorm{r}_{B, \seminorm{q}, v_0}(f) < \infty$, as with $B$ also
    $v_0 + rB$ is bounded. Conversely, we know that the series
    \eqref{eq:GateauxTaylor} converges on all of $B$ by our definition of
    $\Bdd(U,v_0)$ and thus we get
    \begin{equation}
        \seminorm{p}_{B,\seminorm{q}}(f)
        \le
        \sum_{n=0}^{\infty}
        \seminorm{p}_{B,\seminorm{q}}
        (P_{n,v_0})
        =
        \seminorm{r}_{B, \seminorm{q}, v_0}(f).
    \end{equation}
    This completes the proof of \ref{item:AlternativeBounded} and
    \ref{item:BoundedAlternativeSeminorms}. Moreover,
    \ref{item:AlternativeBoundedAreSeminorms} is now clear, as pointwisely convergent
    series over seminorms are seminorms. Finally, to see
    \ref{item:BoundedAlternativeEntireSeries} and \ref{item:BoundedAlternativeEntire},
    note that
    \begin{equation}
        \mathfrak{B}(V,v_0)
        =
        \bigl\{
        B \subseteq V
        \colon
        B \;
        \textrm{is absolutely convex and bounded,} \;
        v_0 + B
        \subseteq
        V
        \bigr\}
    \end{equation}
    is simply the collection of \emph{all} absolutely convex and bounded subsets of $V$,
    as the second condition is vacuous. In particular, we get
    \begin{equation}
        \mathfrak{B}(V,v_0)
        =
        \mathfrak{B}(V,v_0')
        \qquad
        \textrm{for all}
        \quad
        v_0,v_0' \in V.
        \tag*{$\qed$}
    \end{equation}
\end{proof}

By Corollary~\ref{cor:TaylorFrechetContinuity}, the Taylor polynomials of a Fréchet
holomorphic function are continuous, which by
Lemma~\ref{lem:PolynomialContinuityVsBoundedness} in particular implies their
boundedness.
\begin{corollary}
    \label{cor:FrechetToTaylorContinuous}
    \index{Taylor series!Continuity}
    Let $V,W$ be locally convex spaces, $U \subseteq V$ open and $v_0 \in U$.
    \begin{corollarylist}
        \item \label{item:FrechetToTaylorContinuous}
        The linear mapping
        \begin{equation}
            P_{n, v_0}
            \colon
            \Holomorphic_\beta(U,W)
            \longrightarrow
            \Pol_\beta^n(V,\widehat{W})
        \end{equation}
        is continuous with respect to the topology of uniform convergence on bounded
        subsets.
        \item \label{item:TaylorOfConvergentBounded}
        If $(f_\alpha) \subseteq \Holomorphic_\beta(U,W)$ converges within
        $\Bounded(U,W)$ to some function
        \begin{equation}
            f \colon U \longrightarrow W,
        \end{equation}
        then $f$ is bounded as well as Gâteaux holomorphic. Furthermore, the Taylor
        polynomials $P_{n,v_0,\alpha}$ of $f_\alpha$ at $v_0$ converge uniformly
        on bounded
        subsets to the Taylor polynomials~$P_{n,v_0}$ of $f$ at $v_0$ for all $n
        \in \N_0$.
    \end{corollarylist}
\end{corollary}
\begin{proof}
    The first part is a direct consequence of the locally convex Cauchy estimates.
    Indeed, let $B
    \in \mathfrak{B}(U,v_0)$, $\seminorm{q} \in \cs(W)$ and $f \in \Holomorphic(U,W)$.
    Then, by \eqref{eq:CauchyEstimatesLocallyConvex} with $r=1$,
    \begin{equation}
        \seminorm{p}_{B, \seminorm{q}}
        (P_{n,v_0})
        \le
        \seminorm{p}_{v_0 + B, \seminorm{q}}
        (f)
        \qquad
        \textrm{for all $\seminorm{q} \in \cs(W)$ and $n \in \N_0$.}
    \end{equation}
    If now $B \subseteq V$ is an arbitrary bounded, absolutely convex subset of $V$,
    then there exists some $r > 0$ such that
    \begin{equation}
        rB \subseteq U - v_0,
        \qquad
        \textrm{that is}
        \quad
        rB \in
        \mathfrak{B}(U,v_0).
    \end{equation}
    By what we have already shown, this yields the continuity estimate
    \begin{equation}
        \seminorm{p}_{B, \seminorm{q}}
        (P_{n,v_0})
        =
        \frac{1}{r^n}
        \cdot
        \seminorm{p}_{rB, \seminorm{q}}
        (P_{n,v_0})
        \le
        \frac{1}{r^n}
        \cdot
        \seminorm{p}_{v_0 + rB, \seminorm{q}}
        (f),
    \end{equation}
    completing the proof of \ref{item:FrechetToTaylorContinuous}. Note that the second
    statement would be an immediate consequence of the first if we assumed the
    continuity of the limit
    $f$, which we do not. Let $(f_\alpha)_{\alpha \in J}$ and~$f$ as described
    in~\ref{item:TaylorOfConvergentBounded} be given. As the topology of uniform
    convergence on bounded subsets is finer than the topology of uniform convergence
    on intersections of $U$ with finite dimensional subspaces, we know that $f$ is
    Gâteaux holomorphic by Proposition~\ref{prop:GateauxTopology}.
    Proposition~\ref{prop:BoundedCompleteness},~\ref{item:BoundedCompleteness}
    moreover implies the boundedness of $f$. Taking another look at the proof of the
    first part, we have not used the continuity of~$f$, but only its boundedness. Thus,
    \ref{item:TaylorOfConvergentBounded} follows.
\end{proof}

Another remarkable and important consequence of our estimates is that Taylor series of
bounded holomorphic functions converge uniformly on bounded subsets.
\begin{corollary}
    \index{Taylor series!Fréchet $\beta$}
    Let $V$ and $W$ be Hausdorff locally convex spaces, $U_0$ a balanced
    neighbourhood of
    zero, $v_0 \in V$ and {$f \in \Holomorphic_\beta(v_0 + U_0,W)$}. Then the Taylor
    series
    \begin{equation}
        \label{eq:GateauxTaylorOnBounded}
        \sum_{n=0}^{\infty}
        P_{n,v_0}
    \end{equation}
    converges towards $f(v_0 + \argument)$ in the space
    $\Holomorphic_\beta(U_0,\widehat{W})$.
\end{corollary}
\begin{proof}
    By Theorem~\ref{thm:GateauxTaylor}, the series \eqref{eq:GateauxTaylorOnBounded}
    converges pointwisely to $f(v_0 + \argument)$ on all of~$U_0$. As the topology
    of uniform convergence on bounded sets is finer than the topology of pointwise
    convergence, it thus suffices to establish that the series
    \eqref{eq:GateauxTaylorOnBounded} converges absolutely. To this end, fix
    $\seminorm{q} \in \cs(W)$ and
    \begin{equation}
        B
        \in
        \mathfrak{B}(U_0, 0)
        =
        \mathfrak{B}
        \bigl(
            v_0
            +
            U_0,v_0
        \bigr).
    \end{equation}
    Then, by the locally convex Cauchy estimates
    \eqref{eq:CauchyEstimatesLocallyConvex} and \eqref{eq:GateauxTaylor}, we have
    \begin{equation}
        \seminorm{p}_{B,\seminorm{q}}
        \biggl(
            f
            -
            \sum_{n=0}^{N}
            P_{n,v_0}
        \biggr)
        =
        \sum_{n=N+1}^{\infty}
        \seminorm{p}_{B,\seminorm{q}}
        \bigl(
            P_{n,v_0}
        \bigr)
        \le
        \seminorm{p}_{v_0 + rB, \seminorm{q}}(f)
        \cdot
        \frac{r^{N+1}}{1-r}
        \overset{N \rightarrow \infty}{\longrightarrow}
        0
    \end{equation}
    for any $0 < r < 1$ with $rB \subseteq U_0$ as before.
\end{proof}

Note that we have to assume $f \in \Holomorphic_\beta(v_0 + U_0,W)$ to make this
argument work: we do not yet know whether the $\Holomorphic_\beta$-topology is
complete. Which brings us to our next goal: to work out sufficient conditions to
ensure this. As a preliminary consideration, we prove another continuity property of
taking Taylor
polynomials of Fréchet holomorphic functions, which is worth comparing with
Corollary~\ref{cor:TaylorFrechetContinuity} and
Corollary~\ref{cor:FrechetToTaylorContinuous}.
\begin{proposition}
    \label{prop:FrechetExpansionToTaylorContinuous}
    Let $V,W$ be locally convex spaces, $U \subseteq V$ be open, $f \in
    \Holomorphic(U,W)$ and $n \in \N_0$. Then the mapping
    \begin{equation}
        \label{eq:FrechetExpansionToTaylorContinuous}
        P_{n,\bullet}
        \colon
        U
        \longrightarrow
        \Pol_\beta(U,\widehat{W}),
    \end{equation}
    which sends $v_0 \in U$ to the $n$-th Taylor polynomial
    $P_{n,\bullet}$ of $f$, is continuous.
\end{proposition}
\begin{proof}
    By Corollary~\ref{cor:TaylorFrechetContinuity}, the mapping
    \eqref{eq:FrechetExpansionToTaylorContinuous} is well defined. Let
    $(v_\alpha)_{\alpha \in J} \subseteq U$ be a convergent net with limit $v_0 \in
    U$. Fix moreover $\epsilon > 0$, $\seminorm{q} \in \cs(W)$ and $B \in
    \mathfrak{B}(U,v_0)$. By continuity of $f$ at $v_0$, there exists an absolutely
    convex zero neighbourhood $U_0 \subseteq V$ such that
    \begin{equation}
        \seminorm{q}
        \bigl(
            f(v_0 + v)
            -
            f(v_0)
        \bigr)
        \le
        \epsilon
        \qquad
        \textrm{for all $v \in U_0$.}
    \end{equation}
    The continuity of vector addition provides another absolutely convex zero
    neighbourhood
    \begin{equation}
        U_0' \subseteq U_0
        \qquad
        \textrm{such that}
        \qquad
        U_0' + U_0'
        \subseteq
        U_0.
    \end{equation}
    By convergence of $(v_\alpha)$ to $v_0$ we find an index $\alpha_0 \in J$ such that
    $v_\alpha \in U_0'$ for all $\alpha \later \alpha_0$. As~$B$ is bounded, there
    moreover exists some $0 < r < 1$ such that the rescaled set $rB \subseteq U_0'$
    and \eqref{eq:GateauxTaylorPolynomials} yields
    \begin{equation}
        P_{n,v_0}(v)
        =
        \frac{1}{2\pi \I}
        \int_{r \boundary \mathbb{D}}
        \frac{f(v_0 + zv)}{z^n}
        \D z
        \quad \textrm{and} \qquad
        P_{n,v_\alpha}(v)
        =
        \frac{1}{2\pi \I}
        \int_{r \boundary \mathbb{D}}
        \frac{f(v_\alpha + zv)}{z^n}
        \D z
    \end{equation}
    for all $v_\alpha$ with $\alpha \later \alpha_0$ and $v \in U_0$. Putting
    everything together, this provides the estimate
    \begin{align}
        &\seminorm{q}
        \bigl(
            P_{n,v_0}(v)
            -
            P_{n,v_\alpha}(v)
        \bigr) \\
        &\le
        \frac{1}{2\pi r^n}
        \int_0^{2\pi}
        \seminorm{q}
        \bigl(
            f(v_0 + r e^{\I t}v)
            -
            f(v_\alpha + r e^{\I t}v)
        \bigr)
        \D t \\
        &\le
        \frac{1}{2\pi r^n}
        \biggl(
            \int_0^{2\pi}
            \seminorm{q}
            \bigl(
                f(v_0 + e^{\I t} \underbrace{r v}_{\in rB \subseteq U_0'})
                -
                f(v_0)
            \bigr)
            +
            \seminorm{q}
            \bigl(
                f(v_\alpha + e^{\I t} \underbrace{r v}_{\in rB \subseteq U_0'})
                -
                f(v_0)
            \bigr)
            \D t
        \biggr) \\
        &\le
        \frac{1}{2\pi r^n}
        \biggl(
        \int_0^{2\pi}
        \seminorm{q}
        \bigl(
            f(v_0 + \underbrace{r e^{\I t} v}_{\in U_0' \subseteq U_0})
            -
            f(v_0)
        \bigr)
        +
        \seminorm{q}
        \bigl(
            f(v_0 + \underbrace{\underbrace{v_\alpha - v_0}_{\in U_0'} + \underbrace{r
            e^{\I t}v}_{\in U_0'}}_{\in U_0})
            -
            f(v_0)
        \bigr)
        \D t
        \biggr) \\
        &\le
        \frac{2\epsilon}{r^n}
    \end{align}
    for all $v \in B$ and $\alpha
    \later \alpha_0$, where we have used the absolute convexity of $U_0'$ several
    times.
\end{proof}

Building on Theorem~\ref{thm:CompletenessHomogeneousPolynomials}, we add another
assumption on the domain space $V$ to facilitate completeness.
Recall that a topological space $X$ is called \emph{Baire} space if countable
intersections of dense open subsets of $X$ are dense. Or, equivalently, if any
countable
union~$A$ of closed subsets $A_n$ of $X$ has an inner point, then so does one of the
$A_n$. Crucially, complete metric spaces are always Baire by Baire's first category
theorem\footnote{René-Louis Baire (1874-1932) was a French mathematician, whose
study of pointwise limits of continuous functions lead to him to what the mathematical
community nowadays calls his category theorems, which are at the core of modern
functional analysis. Surprising many of his contemporaries, Baire proved that there is no
function $f \colon \R \longrightarrow \R$, which is continuous \emph{exactly} at all
rational numbers.}
and thus the same is true for Banach and Fréchet spaces. We refer to
\cite[§4.6]{koethe:1969a} or \cite[Ch.~IX~§5.]{bourbaki:1998d} for a systematic
treatment, but cannot resist to mention the
following characterization of the Baire property for topological vector spaces due to
Saxon.\footnote{Stephen Apollos Saxon is a Professor Emeritus of Florida State
University and a specialist in locally convex analysis and various notions related to
barrels.}
\begin{proposition}[Baire property in topological vector spaces
    {\cite[Thm.~1]{saxon:1974a}}]
    \label{prop:BaireTVS}
    \index{Baire space}
    \index{Saxon, Stephen}
    \index{Baire, René-Louis}
    Let $V$ be a topological vector space not endowed with the indiscrete
    topology.\footnote{The topology $\{\emptyset,V\}$.}
    Then $V$ is Baire if and only if every closed balanced and absorbing subset $A
    \subseteq V$ has non-empty interior.
\end{proposition}
\begin{proof}[Of the trivial implication]
    Let $A \subseteq V$ be closed balanced and absorbing. We denote the interior of a
    set $S$ by $S^\interior$. As the multiplication with a fixed scalar is a
    homeomorphism of $V$, we have
    \begin{equation*}
        A^\interior \neq \emptyset
        \iff
        \exists_{n \in \N}
        \colon
        (nA)^\interior \neq \emptyset
        \iff
        \forall_{n \in \N}
        \colon
        (nA)^\interior \neq \emptyset.
    \end{equation*}
    Moreover, as $A$ is absorbing, we may write
    \begin{equation}
        V
        =
        \bigcup_{n \in \N}
        (nA),
    \end{equation}
    which certainly has an inner point and is a countable union of closed sets, again
    as the scaling with $n \in \N$ is a homeomorphism and $A$ is closed by assumption.
    Thus, if $V$ is Baire, then there exists some $n \in \N$ with $(nA)^\interior \neq
    \emptyset$ and thus $A^\interior \neq \emptyset$.
\end{proof}

Recall that a \emph{barrel} in $V$ is a closed, balanced, absorbing and \emph{convex}
subset and that~$V$ is called \emph{barrelled} if every barrel is a neighbourhood of
zero.
The condition in Proposition~\ref{prop:BaireTVS} is thus reminiscent of $V$ being
barrelled: every closed, absolutely convex, absorbing subset of $V$ is a zero
neighbourhood. Comparing this with
Proposition~\ref{prop:BaireTVS}, we see that one trades convexity for asserting that
zero is an inner point -- rather than just having nonempty interior. Adapting the proof
of the easy implication yields the following well known lemma.
\begin{lemma}[Baire implies barrelled]
    \index{Baire space!Barrels}
    \label{lem:BaireBarrels}
    Let $V$ be a topological vector space with the Baire property. Then $V$ is barrelled.
\end{lemma}
\begin{proof}
    By the easy implication of Proposition~\ref{prop:BaireTVS}, every barrel $B \subseteq
    V$ possesses an inner point $v_0$, say with open neighbourhood $U \subseteq B$.
    By balancedness of $B$ and the continuity of the multiplication with scalars, this
    implies that $-U \subseteq B$ is an open neighbourhood of $-v_0$. But then $0 = v_0
    + (-v_0)$ is a linear combination of inner points and thus inner itself, with open
    neighbourhood given by
    \begin{equation}
        U + (-U)
        =
        \bigl\{
            u_1 - u_2
            \colon
            u_1, u_2 \in U
        \bigr\}.
    \end{equation}
    Thus $V$ is indeed barrelled.
\end{proof}

The converse of Lemma~\ref{lem:BaireBarrels} does not hold, as every LF
space\footnote{Countable strict inductive limits of Fréchet spaces $V_n$ with
$V_{n+1}
\setminus V_n \neq \emptyset$ for all $n \in \N_0$. One prominent example of such a
space is the space of test functions $\Cinfty_c(\R^n)$.} is
barrelled, but never Baire.\footnote{Indeed, as a vector space, a strict inductive limit
$\varinjlim V_n$ is given by the union over all $V_n$ and as proper subspaces, each of
them has empty interior.} The Baire property supplies a useful sufficient
criterion for when the continuity of the Taylor polynomials passes to the
corresponding Gâteaux holomorphic function. In view of \eqref{eq:GateauxTaylor}, this
is not too surprising, as Baire's main objective within \cite{baire:1899a,baire:1905a} was
the study of pointwise limits of continuous mappings, which ultimately resulted in the
isolation of the Baire property as a useful tool. The following proposition is a
generalization of \cite[Example~3.8(a)]{dineen:1999a}.
\begin{proposition}
    \label{prop:FrechetBaire}
    \index{Fréchet!Baire}
    Let $V,W$ be locally convex spaces such that $V$ is Baire and $W$ is complete.
    Furthermore let $U \subseteq V$ be open as well as connected and $v_0 \in U$.
    Moreover assume that $f \in \Holomorphic_G(U,W)$ possesses continuous Taylor
    polynomials $P_{n,v_0} \in \Pol^n(V,W)$ for all $n \in \N_0$. Then $f \in
    \Holomorphic(U,W)$.
\end{proposition}
\begin{proof}
    We have to prove continuity of $f$, which is equivalent to continuity of the
    compositions $\pi_\seminorm{q} \circ f \colon U \longrightarrow
    \gls{LocalBanach}$
    for all $\seminorm{q} \in \cs(W)$, where $\pi_\seminorm{q} \colon W
    \longrightarrow W_\seminorm{q}$ denotes the local Banach
    space at $\seminorm{q}$.\footnote{The completion of
    $W/\ker\seminorm{q}$ with respect to the norm $\norm{[w]}_{\seminorm{q}} \coloneqq
    \seminorm{q}(w)$ for $w \in W$.} By translating and shrinking $U$ if necessary, we
    may assume that~$v_0 = 0$ and that $U$ is absolutely convex. Consider
    \begin{equation*}
        U_n
        \coloneqq
        \bigcap_{m=0}^\infty
        \bigl\{
            v \in U
            \colon
            \seminorm{q}
            \bigl(
                P_m(v)
            \bigr)
            \le
            n
        \bigr\},
    \end{equation*}
    which is closed subset of $U$ as an intersection of such, where we use the
    continuity of the Taylor polynomials $P_m$. By convergence of the Taylor series
    \eqref{eq:GateauxTaylor} for fixed $v \in U$, we know that
    $(\seminorm{q}(P_m(v)))_m \subseteq W$ is a zero sequence and thus
    $\bigcup_{n \in \N_0} U_n = U$.

    By the Baire property, there thus exists an index
    $N \in \N_0$ and $v_0' \in U_N$ such that $v_0'$ is an interior point of $U_N$. Let
    $U_0 \subseteq V$ be an absolutely convex zero neighbourhood such that $v_0' +
    U_0 \subseteq V_N$. Using homogeneity and
    \eqref{eq:PolynomialsTranslationForward} in the extended setting as explained in
    Remark~\ref{rem:PolynomialTranslation}, we estimate
    \begin{equation}
        \sum_{n=0}^\infty
        \seminorm{p}_{U_0/2, \seminorm{q}}(P_n)
        =
        \sum_{n=0}^\infty
        2^{-n}
        \cdot
        \seminorm{p}_{U_0, \seminorm{q}}(P_n)
        \le
        \sum_{n=0}^\infty
        2^{-n}
        \cdot
        \seminorm{p}_{v_0' + U_0, \seminorm{q}}(P_n)
        \le
        2N.
    \end{equation}
    Hence, the Taylor series \eqref{eq:GateauxTaylor} converges uniformly on the open
    neighbourhood $v_0 + U_0/2$ of $v_0$ and consequently $\pi_\seminorm{q} \circ f$
    is continuous on $v_0 + U_0/2$. We have shown
    \begin{equation}
        \pi_\seminorm{q}
        \circ f
        \in
        \Holomorphic
        \bigl(
            U_0/2,W_{\seminorm{q}}
        \bigr).
    \end{equation}

    Let now $v_0'' \in U$. Then, by Corollary~\ref{cor:GateauxTaylorContinuity}, and
    indicating the function dependence in square brackets, we have
    \begin{equation}
        P_{n,v_0''}[\pi_\seminorm{q} \circ  f]
        =
        P_{n,v_0''}
        \biggl[
            \sum_{m=0}^{\infty}
            \pi_\seminorm{q} \circ P_m
        \biggr]
        =
        \sum_{m=0}^{\infty}
        P_{n,v_0''}
        [\pi_\seminorm{q} \circ P_m]
        =
        \sum_{m=n}^\infty
        P_{n,v_0''}
        [\pi_\seminorm{q} \circ P_m],
    \end{equation}
    where we have also used \eqref{eq:PolynomialDerivative}, and that the series
    converges uniformly on intersections of $U_0/2$ with finite dimensional subspaces,
    see again Theorem~\ref{thm:GateauxTaylor}, \ref{item:GateauxTaylorConvergence}.
    By continuity of the $P_m$ and
    Corollary~\ref{cor:TaylorFrechetContinuity} -- or by virtue of the explicit
    formula -- each of the Taylor polynomials
    $P_{n,v_0''}[\pi_\seminorm{q} \circ P_m]$ is continuous. Invoking Baire's Theorem
    on points of continuity~\cite{baire:1905a}, we get that the pointwise limit
    $P_{n,v_0''}[\pi_\seminorm{q}
    \circ f]$ of continuous functions has a point of continuity. By
    Proposition~\ref{prop:PolynomialContinuity}, this implies continuity of~$P_{n,v_0''}$
    on all of $V$. By the first part of the proof, we get continuity of $\pi_\seminorm{q} \circ
    f$ on some neighbourhood of $v_0''$. As $v_0''$ was
    arbitrary, this yields $\pi_\seminorm{q} \circ f \in \Holomorphic(U,W_\seminorm{q})$
    and variation of~$\seminorm{q}$ completes the proof.
\end{proof}

Taking another look at Example~\ref{ex:GateauxInfiniteDimensionalRange},
Example~\ref{ex:HolomorphicUnbounded} and
Example~\ref{ex:FrechetLocallyUnbounded}, all
functions we have considered are Fréchet
holomorphic on the entirety of their domains. By a slight variation of
\cite[Example~3.8(d)]{dineen:1999a}, there is another class of
examples with the pleasant property described in Proposition~\ref{prop:FrechetBaire}.
\begin{corollary}
    \label{cor:DirectLimitBaire}
    Let $V = \gls{DirectLimit} B_\alpha$ be a direct limit of Baire locally convex
    spaces $B_\alpha$ both in the category of topological spaces and continuous maps
    and in the category of topological vector spaces. If $W$ is
    a complete Hausdorff locally convex space and $f \in \Holomorphic_G(V,W)$
    such that the
    Taylor polynomials $P_{n,0} \in \Pol^n(V,W)$, then $f \in \Holomorphic(V,W)$.
\end{corollary}
\begin{proof}
    By the characteristic property of the \emph{topological} inductive limit, the
    mapping
    \begin{equation}
        f \colon V \longrightarrow W
    \end{equation}
    is continuous if and only if all pullbacks
    \begin{equation}
        f_\alpha \coloneqq f \circ \phi_\alpha \colon B_\alpha \longrightarrow W
    \end{equation}
    are continuous for all $\alpha \in J$. Now, by Gâteaux holomorphy of~$f$, we
    may
    Taylor expand~$f$ around $0 = \phi_\alpha(0)$ by
    Theorem~\ref{thm:GateauxTaylor}.
    This yields unique $n$-homogeneous polynomials
    \begin{equation}
        P_{n,0} \colon V
        \longrightarrow W
        \qquad
        \textrm{for all }
        n \in \N_0
    \end{equation}
    such that
    \begin{equation}
        \bigl(
            f \circ \phi_\alpha
        \bigr)(v)
        =
        f
        \bigl(
            \phi_\alpha(v)
        \bigr)
        =
        \sum_{n=0}^{\infty}
        P_{n,0}
        \bigl(
            \phi_\alpha(v)
        \bigr)
        =
        \sum_{n=0}^{\infty}
        \bigl(
            P_{n,0} \circ \phi_\alpha
        \bigr)(v)
    \end{equation}
    for all $v \in B_\alpha$. By $1$-homogeneity of $\phi_\alpha$\footnote{Here we
    use that $V$ is the direct limit of the $B_\alpha$ in the category of topological
    vector spaces.}, the mappings
    $P_{n,0} \circ \phi_\alpha$ are $n$-homogeneous for all $n \in \N_0$. As Taylor
    expansions are unique, this implies that $f_\alpha$ is Gâteaux holomorphic with
    Taylor polynomials $Q_n \coloneqq P_n \circ \phi_\alpha$ at $0$. By assumption,
    every $Q_n$ is continuous and as $B_\alpha$ is Baire, this implies continuity of
    $f_\alpha$ by Proposition~\ref{prop:FrechetBaire}. Hence, $f_\alpha \in
    \Holomorphic(B_\alpha,W)$ and by the characteristic property of the topological
    direct limit we are done.
\end{proof}

\begin{theorem}[Completeness of $\Holomorphic_\beta$]
    \label{thm:CompletenessBetaFrechet}
    \index{Completeness!Fréchet holomorphic $\beta$}
    \index{Fréchet!Completeness $\beta$}
    Let $V$ and $W$ be locally convex spaces such that $V$ is bornological as well as
    Baire and $W$ is complete Hausdorff. Then $\Holomorphic_\beta(U,W)$ is complete
    Hausdorff for any open $U \subseteq V$.
\end{theorem}
\begin{proof}
    As in the proof of Theorem~\ref{thm:CompletenessHomogeneousPolynomials}, we view
    $\Holomorphic_\beta(U,W)$ as a subspace of $\Bounded(U,W)$. Consequently, the
    Hausdorff property is obvious and it suffices to prove the closedness of the
    subspace $\Holomorphic_\beta(U,W)$ within $\Bounded(U,W)$ by its completeness,
    which we have established in Proposition~\ref{prop:BoundedCompleteness},
    \ref{item:BoundedCompleteness}.

    To this end, let
    $(f_\alpha)_{\alpha \in J} \subseteq \Holomorphic_\beta(U,W)$ be a convergent net
    with limit $f \in \Bounded(U,W)$. Then~$f$ is Gâteaux holomorphic by virtue of
    Proposition~\ref{prop:GateauxTopology}. It remains to establish the continuity of
    $f$. Fix $v_0 \in U$ and
    first note that by Corollary~\ref{cor:FrechetToTaylorContinuous},
    \ref{item:TaylorOfConvergentBounded}, the Taylor polynomials $P_{n,v_0,\alpha}$ of
    $f_\alpha$ at $v_0$ converge uniformly on bounded subsets to the Taylor
    polynomials $P_{n,v_0}$ of $f$ at $v_0$. The latter indeed exist by a
    combination of $f \in \Holomorphic_G(U,W)$ and Theorem~\ref{thm:GateauxTaylor}.
    By Fréchet holomorphy of $f_\alpha$, each of the Taylor polynomials
    \begin{equation}
        P_{n,v_0,\alpha}
        \colon
        V
        \longrightarrow
        W
    \end{equation}
    is continuous by virtue of Corollary~\ref{cor:TaylorFrechetContinuity}. Invoking the
    completeness of $\Pol_\beta^n(V,W)$ for every~$n \in \N_0$ from
    Theorem~\ref{thm:CompletenessHomogeneousPolynomials} this, in turn, implies that
    $P_{n,v_0}$ is continuous for all~$n \in \N_0$. As $V$ is Baire by assumption,
    this is sufficient to ensure the continuity of $f$ by
    Proposition~\ref{prop:FrechetBaire}. Thus, $f \in \Holomorphic(U,W)$ and we are
    done.
\end{proof}

As we have noted before, if $V$ is Banach and $W$ is Fréchet, then the
$\beta$-topology may be
induced by countably many seminorms. Combining this with
Theorem~\ref{thm:CompletenessBetaFrechet} yields the following:
\begin{corollary}
    Let $V$ be Banach, $W$ be Fréchet and $U \subseteq V$ be open. Then the space
    \begin{equation}
        \Holomorphic_\beta
        (U,W)
    \end{equation}
    is Fréchet.
\end{corollary}

Weakening completeness to sequential completeness once again yields an analogous
statement, essentially by the very same methods.
\begin{corollary}
    Let $V,W$ be locally convex spaces such that $V$ is bornological as well as Baire
    and $W$ is a sequentially complete Hausdorff space. Then
    $\Holomorphic_\beta(U,W)$ is sequentially complete for any open $U \subseteq V$.
\end{corollary}

\begin{remark}[Oka principle]
    \index{Oka!Principle}
    \index{Oka!Kiyoshi}
    \label{rem:Oka}
    Taking another look at our proof of Theorem~\ref{thm:CompletenessBetaFrechet}, we
    note that the only problematic aspect is the continuity of the pointwise limit. Indeed,
    by virtue of constant functions being holomorphic, the completeness of $W$ is
    certainly a necessary condition. One may view this as an instance of the Oka
    principle\footnote{Named in honour of the Japanese complex geometer Kiyoshi Oka
    (1901-1978), whose seminal coherence Theorems are at the heart of modern
    complex analysis of several variables, at least for its sheaf theoretic
    incarnation.}: it
    constitutes a holomorphic problem that only
    possesses a topological
    obstruction.
\end{remark}

We conclude the section with studying an explicit example of a Fréchet holomorphic
function defined on the LF space of continuous test functions.
\begin{example}[Test functions]
    \index{Test functions}
    Consider the LF space of continuous test functions
    \begin{equation}
        \gls{TestFunctions}(\R^n)
        \coloneqq
        \varinjlim
        \Continuous_K(\R^n)
    \end{equation}
    with
    \begin{equation}
        \Continuous_K(\R^n)
        \coloneqq
        \bigl\{
            f \in \Continuous(\R^n)
            \colon
            f \at{\R^n \setminus K}
            \equiv
            0
        \bigr\}
    \end{equation}
    for compact $K \subseteq \R^n$ and the natural inclusions $\Continuous_K(\R^n)
    \subseteq \Continuous_{K'}(\R^n)$ whenever $K' \subseteq K$. We study the
    variant of
    Dirac's\footnote{Paul Dirac (1902-1984) was a British theoretical and mathematical
    physicist, who won the Nobel prize in 1933. He may be regarded as one of the
    founding fathers of quantum mechanics and, subsequently, quantum field theory.}
    comb given by
    \begin{equation}
        f
        \colon
        \Continuous_c(\R) \longrightarrow \C, \quad
        f(\phi)
        \coloneqq
        \sum_{n \in \N_0}
        \phi(n)^n.
    \end{equation}
    As the series terminates for every $\phi \in \Continuous_c(\R)$, the function $f$
    is well
    defined. As $n$-th term is a continuous polynomial of degree $n$, we have thus
    defined $f$ by its Taylor series and thus it is Gâteaux holomorphic by
    Corollary~\ref{cor:GateauxPowerSeries}.

    If $B \subseteq \Continuous_c(\R)$ is bounded, then there exists a compact
    interval $[-N,N] \subseteq \R$ such that $B$ is contained within
    $\Continuous_K(\R)$ and bounded there by \cite[§19.5~(5)]{koethe:1969a}. That
    is to say
    \begin{equation}
        \sup_{\phi \in B}
        \supnorm{\phi}
        \coloneqq
        \sup_{\phi \in B}
        \sup_{x \in \R}
        \abs[\big]
        {\phi(x)}
        <
        \infty.
    \end{equation}
    Hence,
    \begin{equation}
        \sup_{\phi \in B}
        \abs[\big]
        {f(\phi)}
        \le
        \sup_{\phi \in B}
        \sum_{n=0}^{N-1}
        \phi(n)^n
        \le
        N
        \cdot
        \max
        \biggl\{
            \sup_{\phi \in B}
            \supnorm{\phi},
            \sup_{\phi \in B}
            \supnorm{\phi}^N
        \biggr\},
    \end{equation}
    and consequently $f$ is bounded.

    Next, we prove the continuity of $f$. As already noted before, we may not rely on
    Proposition~\ref{prop:FrechetBaire} as LF spaces are never Baire and the condition
    in Corollary~\ref{cor:DirectLimitBaire} is impractical. Instead, we will be using
    the explicit description of the continuous seminorms for $\Continuous_c(\R)$ as it
    can be found in \cite[Thm.~2.1.5]{hoermander:1990a} for
    $\gls{TestFunctionsSmooth}(\R)$ with the
    obvious modification of dropping the derivatives. Indeed, using the same ideas it is not
    hard to show that if $V = \varinjlim V_n$ is a countable strict inductive limit
    realized as the union $V = \bigcup_{n \in \N} V_n$, then the seminorms
    \begin{equation}
        \index{Seminorms!LF}
        \seminorm{p}_a
        \coloneqq
        \sup_{n \in \N}
        a_n
        \cdot
        \seminorm{p_n}
    \end{equation}
    with arbitrary sequences $a = (a_n)_n \subseteq [0,\infty)$ of weights and
    $\seminorm{p}_n \in \cs(V_n)$ for all $n \in \N$ provide a defining system of
    seminorms for $V$. For $V = \Continuous_c(\R)$, one may take the Banach spaces
    \begin{equation}
        V_n
        \coloneqq
        \Continuous
        \bigl(
            [-n,n]
        \bigr)
        \quad \textrm{with norms} \quad
        \norm{\phi}_n
        \coloneqq
        \max_{x \in [-n,n]}
        \abs[\big]
        {\phi(x)}.
    \end{equation}
    Interpolating the sequence $a$ to a continuous -- or smooth -- function on $\R$
    results in the version of the seminorms provided in \cite{hoermander:1990a}.

    Let $a_n \coloneqq 2^n$ and $\phi \in \Ball_{\seminorm{p}_a,1}(0)$. Then
    $\norm{\phi}_n \le 1$, i.e. $\norm{\phi}_n^n \le \norm{\phi}_n$, for all $n \in
    \N$, and we we may estimate
    \begin{equation}
        \sum_{n=N+1}^{\infty}
        \abs[\big]
        {\phi(n)}^n
        \le
        \sum_{n=N+1}^{\infty}
        2^{-n}
        \cdot
        2^n
        \cdot
        \norm{\phi}_n^n
        \le
        \sum_{n=N+1}^{\infty}
        2^{-n}
        \cdot
        \seminorm{p}_a(\phi)
        \le
        2^N
    \end{equation}
    for all $N \in \N_0$. Hence, the Taylor series of $f$ converges \emph{uniformly}
    on the open set $\Ball_{\seminorm{p}_a,1}(0)$ and $f$ is continuous on this
    neighbourhood of zero.

    Let now $\phi_0 \in \Continuous_c(\R)$ and $C \coloneqq \max\{\supnorm{\phi},
    1\}$. We compute
    \begin{equation}
        f(\phi_0 + \phi)
        =
        \sum_{n = 0}^\infty
        \bigl(
            \phi_0(n) + \phi(n)
        \bigr)^n
        =
        \sum_{n = 0}^\infty
        \sum_{k=0}^n
        \binom{n}{k}
        \phi_0(n)^{n-k}
        \cdot
        \phi(n)^{k},
    \end{equation}
    for $\phi \in \Continuous_c(\R)$. If now $b_n \coloneqq (4C)^n$ for $n \in \N$ and
    $\phi \in \Ball_{\seminorm{p}_{b},1}(0)$, then this provides the estimate
    \begin{equation}
        \sum_{n = N+1}^\infty
        \abs[\big]
        {
            \phi_0(n) + \phi(n)
        }^n
        \le
        \sum_{n = N+1}^\infty
        C^n
        \sum_{k=0}^n
        \binom{n}{k}
        \norm[\big]{\phi}_n
        =
        \sum_{n = N+1}^\infty
        (2C^n)
        \cdot
        \norm[\big]{\phi}_n
        \le
        2^N
    \end{equation}
    for $N \in \N$, proving the uniform convergence of the Taylor series of $f$ on the
    open neighbourhood $\Ball_{\seminorm{p}_{b},1}(\phi_0)$ of $\phi_0$. Hence, $\phi$
    is continuous at $\phi_0$ and thus continuous everywhere.

    Summarizing, we have established that $f \in
    \Holomorphic_\beta(\Continuous_c(\R))$
    and have thus found an example of a continuous and bounded holomorphic
    function on a non-Baire space. It is straightforward to generalize this to
    series over powers of distributions $\varphi_n$ with supports, say, contained
    in $[-n,n] \setminus [n-1,n-1]$ for all $n \in \N$.
\end{example}

\section{Symmetric Tensors as Polynomials}
\label{sec:RTopologiesPolynomial}
\epigraph{``The truth has got its boots on,'' he said. ``It's going to start
kicking.''}{\emph{The Truth} -- Terry Pratchett}
% !TeX root = ../Dissertation.tex

Having established a firm control over holomorphic differential calculus on
infinite dimensional locally convex spaces and the topology of uniform
convergence on bounded sets, we return to strict deformation quantization. More
precisely, we investigate the complex analytic properties of the spaces of
tensors central to \cite{waldmann:2014a}, when interpreted as polynomials on the
strong dual space. We begin by fixing our notation. Let $V$ be a locally convex
space and denote its $n$-th symmetric tensor power by
\begin{equation}
    V^{\vee n}
    \coloneqq
    \gls{SymmetricPower}
    \coloneqq
    \underbrace{V \vee \cdots \vee V}_{n-\text{times}}
    \qquad
    \textrm{for }
    n \in \N
\end{equation}
and $V^0 \coloneqq \Sym^0(V) \coloneqq \C$. We endow each of the
spaces $\Sym^n(V)$ with the \emph{projective tensor product topology}, which is
generated by the seminorms
\begin{equation}
    \index{Seminorms!Projective tensor power}
    \label{eq:ProjectiveSeminorms}
    \gls{SymmetricSeminormPower}(v)
    \coloneqq
    \inf
    \biggl\{
        \sum_k
        \seminorm{p}
        \bigl(v_1^{(k)}\bigr)
        \cdots
        \seminorm{p}
        \bigl(v_n^{(k)}\bigr)
        \colon
        v
        =
        \sum_k
        v_1^{(k)}
        \tensor \cdots \tensor
        v_n^{(k)}
    \biggr\},
\end{equation}
where the summation over $k$ is finite, and $\seminorm{p}$ varies over $\cs(V)$. The
infimum accounts for the non-uniqueness of the decomposition of a given tensor $v$ into
factorizing ones. For~$n = 0$, one sets $\seminorm{p}^0 \coloneqq
\abs{\argument}$. We indicate this topology with a subscript $\pi$. That is, we
write $\Sym^n_\pi(V)$, and call $\gls{SymmetricPowerProjective}$ the $n$-th
\emph{projective symmetric tensor power} of $V$.

Note that \eqref{eq:ProjectiveSeminorms} makes sense regardless of the symmetry of $v$
and that one may define tensor product of different locally convex spaces and
seminorms in an analogous fashion.\footnote{We will make this precise in
\eqref{eq:ProjectiveSeminormMixed}.} We collect some well known properties of the
projective tensor powers.
\begin{proposition}[Projective tensor powers]
    \index{Tensor product!Projective}
    \index{Projective tensor product}
    \index{Projective tensor power}
    \label{prop:ProjectiveTensor}
    Let $V$ be a locally convex space.
    \begin{propositionlist}
        \item The space $\Sym^1_\pi(V)$ coincides with $V$ as a locally convex space.
        \item If $v_1, \ldots, v_n \in V$ and $\seminorm{p} \in \cs(V)$, then
        \begin{equation}
            \index{Projective tensor power!Factorizing tensors}
            \label{eq:ProjectiveOnFactorizing}
            \seminorm{p}^n(v_1 \vee \cdots \vee v_n)
            =
            \seminorm{p}(v_1) \cdots \seminorm{p}(v_n).
        \end{equation}
        \item The locally convex space $V$ is Hausdorff if and only if $\Sym^n_\pi(V)$
        is for some $n \in \N$. In this case, $\Sym^n_\pi(V)$ is Hausdorff for
        all $n
        \in
        \N_0$.
        \item If $\seminorm{p}_1,\ldots,\seminorm{p}_n \in \cs(V)$, then
        \begin{equation}
            \index{Seminorms!Projective tensor product}
            \label{eq:ProjectiveSeminormMixed}
            \bigl(
                \seminorm{p}_1 \tensor \cdots \tensor \seminorm{p}_n
            \bigr)(v)
            \coloneqq
            \inf
            \biggl\{
            \sum_k
            \seminorm{p}_1
            \bigl(v_1^{(k)}\bigr)
            \cdots
            \seminorm{p}_n
            \bigl(v_n^{(k)}\bigr)
            \colon
            v
            =
            \sum_k
            v_1^{(k)}
            \tensor \cdots \tensor
            v_n^{(k)}
            \biggr\}
        \end{equation}
        defines a continuous seminorm on $\Sym^n_\pi(V)$.
        \item \label{item:ProjectiveAssociativity}
        If $\seminorm{p}_1,\ldots,\seminorm{p}_{n} \in \cs(V)$, then
        \begin{equation}
            \bigl(
                \seminorm{p}_1 \tensor \cdots \tensor \seminorm{p}_m
            \bigr)
            \tensor
            \bigl(
                \seminorm{p}_{m+1} \tensor \cdots \tensor \seminorm{p}_{n}
            \bigr)
            =
            \seminorm{p}_1 \tensor \cdots \tensor \seminorm{p}_{n}
        \end{equation}
        for any $m,n \in \N$ with $1 \le m < n$.\footnote{Indeed, this is what makes
        speaking of tensor powers feasible in the first place.}
        \item The symmetric tensor product is continuous as a bilinear mapping
        \begin{equation}
            \vee
            \colon
            \Sym^n_\pi(V) \times \Sym_\pi^m(V)
            \longrightarrow
            \Sym^{n+m}_\pi(V).
        \end{equation}
        More precisely, if $v \in \Sym^n_\pi(V)$ and $w \in \Sym^m_\pi(V)$, then
        \begin{equation}
            \seminorm{p}^{n+m}(v \vee w)
            \le
            \seminorm{p}^n(v)
            \cdot
            \seminorm{p}^m(w).
        \end{equation}
    \end{propositionlist}
\end{proposition}
\begin{proof}
    All presented statements and a systematic treatment of projective tensor products
    can be found in the textbooks \cite[Ch.~43]{treves:2006a}, \cite[§41]{koethe:1979a} or
    \cite[Ch.~15]{jarchow:1981a}.
\end{proof}

In the context of star products we will actually need a slightly more general construction,
namely the projective tensor product of two different locally convex spaces $V$
and~$W$. The topology of $V \tensor W$ is then generated by the seminorms
\eqref{eq:ProjectiveSeminormMixed} with continuous seminorms~$\seminorm{p}_1 \in
\cs(V)$ and $\seminorm{p}_2 \in \cs(W)$. We indicate this topology
by writing~$V \tensor_\pi W$, whenever there is any room for confusion, and otherwise
suppress the index. Unsurprisingly, the resulting space has properties
analogous to the projective symmetric tensor powers. Moreover, the following rather
convenient characterization of continuity for $n$-multilinear mappings holds, which
provides significant simplifications for certain technical arguments within strict
deformation quantization. This is exemplified within proof-sketch of
Theorem~\ref{Thm:StdOrderedStarProduct} and the proof of
Theorem~\ref{thm:StarProductHolomorphicContinuity}.
\begin{proposition}[Infimum argument, {\cite[Prop.~43.4]{treves:2006a}}]
    \index{Projective tensor product!Infimum argument}
    \index{Infimum Argument}
    \label{prop:InfimumArgument}
    Let $V_1, \ldots, V_n, W$ be locally convex spaces and
    \begin{equation}
        \phi
        \colon
        V_1 \times \cdots \times V_n
        \longrightarrow
        W
    \end{equation}
    be $n$-linear with corresponding linear map
    \begin{equation}
        \Phi
        \colon
        V_1 \tensor \cdots \tensor V_n
        \longrightarrow
        W.
    \end{equation}
    Endow $V_1 \times \cdots \times V_n$ with the Cartesian product topology and
    $V_1 \tensor \cdots \tensor V_n$ with the projective tensor product topology.
    Then $\phi$ is continuous if and only if $\Phi$ is. More precisely, if for
    a continuous seminorm $\seminorm{q} \in \cs(W)$ there are $\seminorm{p}_1 \in
    \cs(V_1)$, $\ldots$,
    $\seminorm{p_n} \in \cs(V_n)$ such that
    \begin{equation}
        \seminorm{q}
        \bigl(
            \phi
            (v_1, \ldots, v_n)
        \bigr)
        \le
        \seminorm{p}_1(v_1)
        \cdots
        \seminorm{p}_n(v_n)
        \qquad
        \textrm{for all }
        v_1 \in V_1,
        \ldots,
        v_n \in V_n,
    \end{equation}
    then
    \begin{equation}
        \seminorm{q}
        \bigl(
            \Phi(v)
        \bigr)
        \le
        \bigl(
            \seminorm{p}_1
            \tensor \cdots \tensor
            \seminorm{p}_n
        \bigr)(v)
        \qquad
        \textrm{for all }
        v \in V_1 \tensor \cdots \tensor V_n,
    \end{equation}
    and vice versa.
\end{proposition}

This means that it suffices to prove continuity estimates for multilinear mappings,
such as star products, on \emph{factorizing tensors} only. In view of the infimum in
\eqref{eq:ProjectiveOnFactorizing}, this constitutes a drastic simplification.
\begin{remark}[Injective tensor product]
    \index{Projective tensor product!Vs. injective}
    \index{Nuclearity}
    \index{Injective tensor product}
    \index{Tensor product!Injective}
    \label{rem:InjectiveTensor}
    There is another natural way of endowing the algebraic tensor product $V \tensor
    W$ with a topology, namely by means of \emph{injective tensor products} of
    seminorms. Indeed, if $\seminorm{p} \in \cs(V)$ and $\seminorm{q} \in \cs(W)$,
    then their injective tensor product is defined by
    \begin{equation}
        \index{Seminorms!Injective tensor product}
        \label{eq:InjectiveSeminorms}
        \bigl(
            \seminorm{p} \, \gls{InjectiveTensor} \, \seminorm{q}
        \bigr)(x)
        \coloneqq
        \sup
        \bigl\{
            \bigl(
                v' \tensor w'
            \bigr)(x)
            \colon
            v'
            \in
            \Ball_{\seminorm{p},1}(0)^\polar, \;
            w'
            \in
            \Ball_{\seminorm{\seminorm{q}},1}(0)^\polar
       \bigr\},
    \end{equation}
    where
    \begin{equation}
        \index{Polar}
        \label{eq:Polar}
        \gls{Polar}
        \coloneqq
        \bigl\{
            \varphi \in V'
            \colon
            \abs[\big]{\varphi(v)}
            \le
            1
            \textrm{ for all }
            v \in A
        \bigr\}
    \end{equation}
    denotes the \emph{polar} of $A \subseteq V$. Remarkably, an analogue of
    \eqref{eq:ProjectiveOnFactorizing} still holds for the injective tensor product, as we
    shall prove in a moment. However, there is no injective version of
    Proposition~\ref{prop:InfimumArgument}, a fact we shall make precise in
    Remark~\ref{rem:InjectiveNoInfimum} after some preliminary
    considerations. Having to work with the supremum and the presence of polars
    make injective tensor products less convenient to work with than projective
    tensor products in the context of strict deformation quantization.

    That being said, ultimately, all locally convex spaces $V,W$ of observables
    utilized in
    the strict deformation quantizations \cite{waldmann:2014a,
    esposito.stapor.waldmann:2017a, schoetz.waldmann:2018a,
    heins.roth.waldmann:2023a} turn out to be \emph{nuclear}: this implies that
    the
    identity mapping
    \begin{equation}
        V \gls{ProjectiveTensor} W
        \longrightarrow
        V \tensor_\epsilon W
    \end{equation}
    constitutes an isomorphism of locally convex spaces. In general, it is at least
    continuous.
\end{remark}

\begin{lemma}[Injective vs. projective $\tensor$, {\cite[Cor. of
Prop.~43.4]{treves:2006a}}]
    \label{lem:ProjectiveVsInjective}
    Let $V$ and~$W$ be locally convex spaces. Then the identity mapping
    \begin{equation}
        \label{eq:ProjectiveVsInjectiveIdentity}
        V \tensor_\pi W
        \longrightarrow
        V \tensor_\epsilon W
    \end{equation}
    is continuous. More precisely, if $\seminorm{p} \in \cs(V)$ and $\seminorm{q} \in
    \cs(W)$, then
    \begin{equation}
        \label{eq:ProjectiveVsInjective}
        \seminorm{p} \tensor_\epsilon \seminorm{q}
        \le
        \seminorm{p} \tensor_\pi \seminorm{q}.
    \end{equation}
\end{lemma}
\begin{proof}
    Let $x \in V \tensor W$ and consider a finite decomposition
    \begin{equation}
        x
        =
        \sum_{k}
        v_k \tensor w_k
    \end{equation}
    into factorizing tensors as well as $v' \in \Ball_{\seminorm{p},1}(0)^\polar$
    and $w' \in \Ball_{\seminorm{q},1}(0)^\polar$. By definition of the polars,
    see again \eqref{eq:Polar}, we have
    \begin{equation}
        \abs[\bigg]
        {
            v'
            \biggl(
                \frac{v}{\seminorm{p}(v)}
            \biggr)
        }
        \le
        1, \qquad
        \textrm{that is} \quad
        \abs[\big]
        {v'(v)}
        \le
        \seminorm{p}(v)
    \end{equation}
    for all $v \in V \setminus \ker \seminorm{p}$. If $\seminorm{p}(v) = 0$, then
    $\seminorm{p}(rv) = 0 \le 1$ for all $r > 0$, as well. By definition of the polar
    this in turn means that
    \begin{equation}
       r
       \cdot
       \abs[\big]
       {v'(v)}
       =
       \abs[\big]
       {v'(rv)}
       \le
       1
       \qquad
       \textrm{for all }
       r > 0,
    \end{equation}
    which implies $v'(v) = 0$. Hence, $\abs{v'} \le
    \seminorm{p}$ holds and analogously we get $\abs{w'} \le
    \seminorm{q}$.\footnote{In passing, we note that this proves ``$\le$'' in
    \eqref{eq:ProjectiveOnFactorizing} for injective tensor products. To see the other
    inequality, fix vectors $v \in V$, $w \in W$ and extend the continuous linear
    functionals
        \begin{equation}
            \varphi_v(\lambda v)
            \coloneqq
            \lambda \cdot \seminorm{p}(v)
            \quad
            \textrm{and}
            \quad
            \varphi_w(\lambda w)
            \coloneqq
            \lambda \cdot \seminorm{q}(w)
        \end{equation}
        from $\C \cdot v$ resp. $\C \cdot w$ to $V$ resp.~$W$ by means of the
        Hahn-Banach Theorem.}
    Consequently,
    \begin{equation}
        \abs[\big]
        {
            \bigl(
                v' \tensor w'
            \bigr)(x)
        }
        \le
        \sum_{k}
        \abs[\big]
        {v'(v_k) \cdot w'(w_k)}
        \le
        \sum_{k}
        \seminorm{p}(v_k)
        \cdot
        \seminorm{q}(w_k).
    \end{equation}
    Note that the left-hand side is independent of the chosen decomposition,
    whereas the right-hand side does not depend on the choice of $v'$ and $w'$. Thus,
    taking the infimum over all decompositions of $x$ and the supremum over all $v'$
    and $w'$ in the respective polars proves the continuity estimate
    \eqref{eq:ProjectiveVsInjective}.
\end{proof}

Along the way, we showed the following alternative description of polars, which
is conceptually pleasing and will be useful later.
\begin{corollary}
    \label{cor:Polars}
    \index{Polar}
    Let $V$ be a locally convex space and $p \in \cs(V)$. Then
    \begin{equation}
        \label{eq:PolarAlternative}
        \Ball_{\seminorm{p},1}(0)^\polar
        =
        \bigl\{
            \varphi
            \in
            V'
            \colon
            \abs{\varphi}
            \le
            \seminorm{p}
        \bigr\}.
    \end{equation}
\end{corollary}

By \cite[Thm.~50.1(f)]{treves:2006a} nuclearity of a locally convex space $V$ may be
characterized as
\begin{equation}
    V \tensor_\pi W
    \cong
    V \tensor_\epsilon W
\end{equation}
as locally convex spaces via the canonical mapping
\eqref{eq:ProjectiveVsInjectiveIdentity} for all locally convex spaces~$W$,
which in view of Lemma~\ref{lem:ProjectiveVsInjective} is in turn equivalent to the
continuity of the identity mapping
\begin{equation}
    V \tensor_\epsilon W
    \longrightarrow
    V \tensor_\pi W
\end{equation}
for all locally convex spaces $W$. Indeed, this was also
Grothendieck's\footnote{Alexander Grothendieck (1928-2014) was the father of all
of
algebraic geometry and received the Field's medal in 1966. Before turning to geometry,
he pursued and largely laid to rest locally convex analysis under his doctoral
advisors
Laurent Schwartz and
Jean Dieudonné.} original definition in \cite{grothendieck:1955a}. We may use this to
prove that having an injective version of Proposition~\ref{prop:InfimumArgument} for
bilinear mappings already implies nuclearity in the following sense.
\begin{remark}
    \label{rem:InjectiveNoInfimum}
    \index{Injective tensor product!Infimum argument}
    Let $V,W$ be locally convex spaces and consider the bilinear mapping
    \begin{equation}
        \tensor
        \colon
        V \times W
        \longrightarrow
        V \tensor_\pi W, \quad
        \tensor(v, w)
        \coloneqq
        v \tensor w.
    \end{equation}
    By \eqref{eq:ProjectiveOnFactorizing}, the mapping $\tensor$ is continuous. If now
    Proposition~\ref{prop:InfimumArgument} were to hold for bilinear mappings, then the
    same would be true for the induced mapping
    \begin{equation}
        V \tensor_\epsilon W
        \longrightarrow
        V \tensor_\pi W,
    \end{equation}
    which is just the identity. Thus we could infer nuclearity of $V$ by varying $W$.
\end{remark}

Note that the utilization of the seminorms \eqref{eq:InjectiveSeminorms} is rather
non-standard in the older literature, which usually defines the topology by
providing a basis of zero neighbourhoods. This has the advantage of capturing the
higher generality of topological vector spaces instead of only locally convex ones,
but is arguably less intuitive to work with. For a comprehensive discussion of
injective and projective tensor products and in particular their interpretation as
continuous $\tensor$-$\Hom$-adjunctions, we refer to the textbooks
\cite[Ch.~43~\&~44]{treves:2006a}, \cite[Ch.~II~§3-4]{bourbaki:1998b}, \cite[Ch.~15
\& 16]{jarchow:1981a} and in particular \cite[Sec.~6]{vogt:2000a} for a more modern
exposition based on seminorms.

After this short digression, we return to the goal of this section, which is to
investigate the mapping
\begin{equation}
    \label{eq:Iota}
    \gls{CanonicalEmbedding}
    \colon
    \gls{SymmetricAlgebra}
    \coloneqq
    \bigoplus_{n=0}^\infty
    \Sym^n(V)
    \longrightarrow
    \Pol(V'), \quad
    \iota(v_1 \vee \cdots \vee v_n)
    \at[\Big]{\varphi}
    \coloneqq
    \prod_{j=1}^n
    \varphi(v_j)
\end{equation}
and its components $\iota_n \coloneqq \iota \at{\Sym^n(V)}$. Our considerations are
thus similar to the classical study of duality as it can e.g. be found in the textbook
\cite[Sec.~1.2 \& Ch.~2]{dineen:1999a}, but with
a notably different outlook and focus: the application to strict deformation
quantization. As a preliminary consideration, we prove the injectivity of $\iota$.
\begin{lemma}
    \label{lem:IotaInjectivity}%
    Let $V$ be a Hausdorff locally convex space. Then the
    mapping~$\iota$ from \eqref{eq:Iota} and its components $\iota_n$ are injective for all
    $n \in \N_0$.
\end{lemma}
\begin{proof}
    As polynomials of different homogeneous degrees are linearly
    independent, it suffices to prove the injectivity of the maps
    \begin{equation}
        \iota_n
        \colon
        \Sym^n(V)
        \longrightarrow
        \Pol(V'_\beta)
    \end{equation}
    for all $n \in \N_0$. For $n = 0$, we simply have
    \begin{equation}
        \label{eq:IotaZeroIsIdentity}
        \iota_0
        =
        \id_\C
    \end{equation}
    and thus there is nothing to be shown. Thus assume $n \in \N$ and let
    \begin{equation}
        0 \neq v_n \in \Sym^n(V).
    \end{equation}
    We commit the cardinal sin\footnote{Sorry!} and choose a linear algebraic basis
    $\{\basis{e}_\alpha \colon \alpha \in I \} \subseteq V$ to expand
    \begin{equation}
        v_n = v^{\alpha_1 \cdots \alpha_n}
        \cdot
        \basis{e}_{\alpha_1} \vee \cdots \vee \basis{e}_{\alpha_n},
    \end{equation}
    where -- and in the sequel -- we also adhere to Einstein's summation
    convention. As~$v_n$ is nonzero,
    there are indices $\alpha_1, \ldots, \alpha_n \in I$ such that the
    corresponding
    coefficient $v^{\alpha_1 \cdots \alpha_n}$ is nonzero, as well. Moreover, if
    $\beta_1, \ldots,
    \beta_n \in I$ is a \emph{different choice}\footnote{That is, $(\beta_1,
    \ldots, \beta_n)$ is not a permutation of $(\alpha_1, \ldots, \alpha_n)$.} of
    indices such that the coefficient
    $v^{\beta_1 \cdots \beta_n} \neq 0$, then
    \begin{equation}
        \beta_j \notin \{\alpha_1, \ldots,
        \alpha_n\}
        \qquad
        \textrm{for some}
        \quad
        j \in \{1, \ldots, n\}.
    \end{equation}
    Choosing one such index for every other non-vanishing coefficient defines a
    finite index set $I_0 \subseteq I$, as we are working with a liner algebraic
    basis. Consider the linear functional
    \begin{align}
        &\psi
        \colon
        \Span
        \bigl(
        \bigl\{
        \basis{e}_{\alpha_1}, \ldots, \basis{e}_{\alpha_n}
        \bigr\}
        \cup
        \{\basis{e}_\beta \colon \beta \in I_0\}
        \bigr)
        \longrightarrow
        \field{C}, \\
        &\psi
        \biggl(
        \sum_{k=1}^n
        \eta_k
        \basis{e}_{\alpha_k}
        +
        \sum_{\beta \in I_0}
        \lambda_\beta
        \basis{e}_{\beta}
        \biggr)
        \coloneqq
        \sum_{k=1}^{n}
        \eta_k .
    \end{align}
    As the various $\basis{e}_\beta$ are linearly independent, our functional $\psi$
    is well
    defined. Moreover, as its domain is finite dimensional and Hausdorff, $\psi$ is
    continuous. By the Hahn-Banach Theorem, we find a continuous extension $\varphi
    \in V'$ of $\psi$. By construction,
    \begin{equation*}
        \iota_n(v_n)
        \at{\varphi}
        =
        v^{\beta_1 \cdots \beta_n}
        \cdot
        \varphi(\basis{e}_{\beta_1})
        \cdots
        \varphi(\basis{e}_{\beta_n})
        =
        v^{\alpha_1 \cdots \alpha_n}
        \neq
        0,
    \end{equation*}
    as every other term contains a factor $0$. This proves that the linear mapping
    \begin{equation}
        \iota_n
        \colon
        \Sym^n(V) \longrightarrow \Pol^n(V'_\beta)
    \end{equation}
    is injective and thus the same is true for $\iota$ itself, as we are dealing with
    direct sums throughout.
\end{proof}

Notice that $\iota_1
\colon V \longrightarrow \Pol(V') = V''$ is the canonical embedding of $V$ into its
bidual.\footnote{To speak of the continuous bidual $V''$ or continuous polynomials on
$V'$, one has to choose a topology for $V'$. This is an oversight, which we shall rectify
momentarily.}
Consequently,
the mapping $\iota$ may be viewed as the algebra morphism
obtained from the universal property of $\Sym^\bullet(V)$ as the \emph{free
commutative algebra} over $V$, see \cite[Ch.~III,
§6]{bourbaki:1998b} or \cite[Sec.~10.13]{stacks-project} for a proper discussion of
this
algebraic point of view. This has the pleasant consequence that
\begin{equation}
    \label{eq:IotaNFactorization}
    \iota_n(v_1 \vee \cdots \vee v_n)
    =
    \iota_1(v_1)
    \cdots
    \iota_1(v_n),
\end{equation}
where we use the pointwise product $V'$ inherits from $\C$ on the right-hand side. As
we are first and foremost interested in the continuity of $\iota$, we have to
endow the
symmetric tensor algebra $\Sym^\bullet(V)$ and its direct summands with a topology; we
have already achieved the latter by means of projective tensor products.
\index{Free commutative algebra}

For the purposes of strict deformation quantization \cite[Sec.~3.1]{waldmann:2014a}
has endowed $\Sym^\bullet(V)$ with a topology finer than the locally convex direct
sum, but coarser than the Cartesian product topology, the so called $R$-topology with
respect to a real parameter $R$. We have encountered this topology for $V = \R$ and $R
= 1/2$ in Example~\ref{ex:StdOrdII}. Both extremal cases may be ruled out a priori,
which has lead to this intermediate choice:

On the one hand, the star products fail to be continuous for the
full Cartesian product, which essentially corresponds to the fact that not every
formal power series converges. And, on the other hand, one wishes to incorporate
certain formal power series, for example
observables of Gaussian decay, in completions of the observable algebras. As the
locally convex direct sum is already complete, this necessitates the passage to a
finer
topology. We refer to the survey \cite{waldmann:2019a} for a comprehensive discussion
of these ideas.

The main result of this section is that these power series and their topology possess
a conceptual interpretation as bounded Fréchet holomorphic functions. To make this
precise, we first recall the definition of the $R$-topologies by means of a defining
system of seminorms. Fix $R \ge 0$. If $\seminorm{p} \in \cs(V)$, then we define
\begin{equation}
    \index{Seminorms!R-Topology}
    \label{eq:RTopologySeminorms}
    \gls{SeminormsRTopology}
    \coloneqq
    \sum_{n=0}^{\infty}
    c^n
    \cdot
    n!^R
    \cdot
    \seminorm{p}^n
\end{equation}
for any $c \ge 0$, where $\seminorm{p}^n$ denotes the $n$-th projective tensor power
from
\eqref{eq:ProjectiveSeminorms}. Note that the series in \eqref{eq:RTopologySeminorms}
is only formally infinite, as $\Sym^\bullet(V)$ is a direct sum. This means that
every element $v \in \Sym^\bullet(V)$ is a finite linear combination of
factorizing tensors. We then write
\begin{equation}
    \label{eq:RTopology}
    \gls{RTopology}
    \coloneqq
    \bigoplus_{n=0}^\infty
    \Sym^n(V)
\end{equation}
for the symmetric algebra of $V$ endowed with the locally convex topology generated by
the seminorms \eqref{eq:RTopologySeminorms} with $\seminorm{p} \in \cs(V)$ and $c \ge
0$. In the sequel, we shall refer to this topology as the \emph{R-topology}. Note
that even in the
case that $V \neq \{0\}$ is normable,
$\Sym_R^\bullet(V)$ is \emph{not}. Different constants~$c \ge 0$
result in inequivalent seminorms.\footnote{In view of the presence of canonical
commutation relations in many important observable algebras, see again
\eqref{eq:CanonicalCommutation}, this is a feature and not a bug if
one wants to achieve the continuity of the star product.} We collect the
properties of~$\Sym_R^\bullet(V)$, which will play a role in our further
considerations.
\begin{proposition}[$R$-Topologies, {\cite[Lem.~3.6]{waldmann:2014a} or
\cite[Prop.~3.1 \& 3.2]{heins.roth.waldmann:2023a}}]
    \index{R-Topologies}
    \label{prop:RTopologies} \; \\
    Let $V$ be a locally convex space and $R \ge 0$.
    \begin{propositionlist}
        \item If $T \le R$, then $\Sym_R^\bullet(V) \subseteq \Sym_T^\bullet(V)$ and
        the inclusion is continuous.
        \item \label{item:RTopologiesSubspace}
        The subspace topology inherited from the inclusion $\Sym^n(V) \subseteq
        \Sym_R^\bullet(V)$ is the projective tensor power topology for all $n \in
        \N_0$.
        \item The space $\Sym_R^\bullet(V)$ is Hausdorff if and only if $V$ is.
        \item \label{item:RTopologyCompletion}
        The completion $\widehat{\Sym}_R^\bullet(V)$ of $\Sym_R^\bullet(V)$ is
        given by
        \begin{equation}
            \widehat{\Sym}_R^\bullet(V)
            =
            \biggl\{
                v
                =
                \sum_{n=0}^{\infty}
                v_n
                \in
                \gls{CartesianProduct}
                \widehat{\Sym}^n(V)
                \colon
                \seminorm{p}_{R,c}(v)
                <
                \infty
                \textrm{ for all }
                \seminorm{p} \in \cs(V),
                c \ge 0
            \biggr\}.
        \end{equation}
        Moreover, the infinite series converges in the $R$-topology.
        \item The symmetric tensor product $\vee$ endows $\Sym^\bullet_R(V)$ with the
        structure of a locally convex algebra. If $R = 0$, then $\vee$ is
        submultiplicative.
        \item For every $\varphi \in V'$, the evaluation
        functional
        \begin{equation}
            \delta_\varphi
            \colon
            \Sym_R^\bullet(V)
            \longrightarrow
            \C, \quad
            \delta_\varphi
            (v)
            \coloneqq
            \sum_{n=0}^{\infty}
            \varphi^k(v_k)
        \end{equation}
        is a continuous character, i.e. constitutes an algebra morphism.
    \end{propositionlist}
\end{proposition}

One may use the very same seminorms to define a topology on the full tensor
algebra~$\Tensor^\bullet(V)$. This is sometimes useful, as the symmetrizer
\begin{equation}
    \gls{Symmetrizer}
    \colon
    \Tensor_R^\bullet(V)
    \longrightarrow
    \Sym_R^\bullet(V), \quad
    \Symmetrizer(v_1 \tensor \cdots \tensor v_n)
    \coloneqq
    \frac{1}{n!}
    \sum_{\sigma \in S_n}
    v_{\sigma(1)}
    \tensor \cdots \tensor
    v_{\sigma(n)}
\end{equation}
then turns out to be a continuous linear projection by
\eqref{eq:ProjectiveOnFactorizing}, where \gls{PermutationGroup}
denotes the permutation group in the letters $\{1, \ldots, n\}$. This more general
point of view is sometimes useful, for instance for our considerations within
Section~\ref{sec:HilberTensorProduct}.

In analogy to indicating the topology of uniform convergence on bounded sets with a
subscript $\beta$, we shall in the sequel use the subscript
\gls{PointwiseConvergence} to indicate the
topology of pointwise convergence -- that is to say, the topology of uniform
convergence on finite sets. For example, we write $V'_\sigma$ and
$\Pol^\bullet_\sigma$. This yields the
following continuity properties of \eqref{eq:Iota}.
\begin{theorem}[Tensors as polynomials on the dual space]
    \label{thm:TensorsAsPolynomials}%
    \index{Polynomials!Induced by tensors}
    \index{Tensors as polynomials}
    Let $V$ be a locally convex space and $n \in \field{N}_0$.
    \begin{theoremlist}
        \item \label{item:IotaV}
        Let $v_1, \ldots, v_n \in V$. The polynomial
        \begin{equation}
            \iota(v_1 \vee \cdots \vee v_n)
            \colon
            V'_\sigma
            \longrightarrow
            \field{C}
        \end{equation}
        is a homogeneous polynomial of degree $n$ of finite type\footnote{A
        polynomial $P$ is
        of \emph{finite type} if there are
        finitely
            many $a_1, \ldots, a_N \in (V'_\beta)'$ and corresponding coefficients
            $\lambda_1,
            \ldots, \lambda_N \in \C$ such that
            \begin{equation}
                P(\varphi)
                =
                \sum_{k=0}^{N} \lambda_k
                a_k(\varphi)^n
                \qquad
                \textrm{for all }
                \varphi \in V'.
            \end{equation}
            Note that this linear combination takes place within fixed
            homogeneity and is thus different from the decomposition of $\iota$ into
            its homogeneous degrees~$\iota_n$.}
            and continuous. In particular, it is continuous as a polynomial on $V'_\beta$.
        \item \label{item:IotaNSigmaContinuous}%
        The restrictions
        \begin{equation}
            \iota_n
            \coloneqq
            \iota
            \at[\Big]{\Sym^n(V)}
            \colon
            \Sym^n_\pi(V)
            \longrightarrow
            \Pol_\sigma^n
            (V'_\beta)
        \end{equation}
        are continuous for all $n \in \N_0$.
        \item \label{item:IotaNBetaContinuous}%
        If $V$ is barrelled, then the restrictions
        \begin{equation}
            \iota_n
            =
            \iota
            \at[\Big]{\Sym^n}
            \colon
            \Sym^n_\pi(V)
            \longrightarrow
            \Pol_\beta^n(V'_\beta)
        \end{equation}
        are continuous.
        \item \label{item:IotaContinuous}%
        If $V$ is barrelled and $R \ge 0$, then the mapping
        \begin{equation}
            \iota
            \colon
            \Sym_R(V)
            \longrightarrow
            \Pol^\bullet_\beta
            (V'_\beta)
        \end{equation}
        is continuous.
    \end{theoremlist}
\end{theorem}
\begin{proof}
    The case $n=0$ is once again trivial by \eqref{eq:IotaZeroIsIdentity}. Let $n
    \in \field{N}$ and $v_1, \ldots, v_n \in V$. The homogeneity of $\iota(v_1
    \vee \cdots \vee v_n)$ is clear by \eqref{eq:Iota}. Note that the seminorm
    \begin{equation}
        \seminorm{p}(\varphi)
        \coloneqq
        \max_{j=1, \ldots, n}
        \abs[\big]{\varphi(v_j)}
    \end{equation}
    is $\sigma$-continuous. Thus the estimate
    \begin{equation}
        \abs[\big]
        {\iota(v_1 \vee \cdots \vee v_n)\varphi}
        =
        \prod_{j=1}^n
        \abs[\big]
        {\varphi(v_j)}
        \le
        \prod_{j=1}^n
        \seminorm{p}(\varphi)
        =
        \seminorm{p}(\varphi)^n
    \end{equation}
    proves the $\sigma$-continuity of $\iota(v_1 \vee \cdots \vee v_n)$ by virtue of
    Proposition~\ref{prop:PolynomialContinuity},
    \ref{item:PolynomialContinuityEstimate}. As the strong topology is finer than
    the $\sigma$-topology, the polynomial $\iota(v_1 \vee \cdots \vee v_n)$ is thus
    also continuous on $V_\beta'$. We have shown the continuity assertions in
    \ref{item:IotaV}. Let now
    \begin{equation}
        v
        =
        \sum_k v_1^{(k)} \vee \cdots \vee v_n^{(k)}
        \in
        \Sym^n(V)
    \end{equation}
    and $\varphi \in V'$ with corresponding seminorm
    $\seminorm{p}_\varphi(Q) \coloneqq \abs{Q(\varphi)}$ for $Q \in
    \Pol^\bullet_\sigma(V'_\beta)$. By continuity of $\varphi$, we moreover note
    $\seminorm{q}_\varphi \coloneqq \abs{\varphi} \in \cs(V)$. We arrive at the
    estimate
    \begin{equation}
        \label{eq:IotaProof1}
        \seminorm{p}_\varphi
        \bigl(
        \iota_n(v)
        \bigr)
        =
        \abs[\bigg]
        {
            \sum_k
            \prod_{j=1}^n
            \varphi(v_j^{(k)})
        }
        \le
        \sum_k
        \prod_{j=1}^n
        \abs[\big]
        {\varphi(v_j^{(k)})}
        =
        \sum_k
        \prod_{j=1}^n
        \seminorm{q}_\varphi
        \bigl(
            v_j^{(k)}
        \bigr).
    \end{equation}
    Taking the infimum over all decompositions of $v$ into factorizing tensors thus yields the continuity estimate
    \begin{equation}
        \label{eq:IotaTaylorEstimate}
        \abs[\big]
        {
            \bigl(
                \iota_n(v)
            \bigr)(\varphi)
        }
        =
        \seminorm{p}_\varphi
        \bigl(
        \iota_n(v)
        \bigr)
        \le
        \seminorm{q}_\varphi^n(v),
    \end{equation}
    where $\seminorm{q}^n_\varphi$ denotes the $n$-th projective tensor power of
    $\seminorm{q}_\varphi$. This proves \ref{item:IotaNSigmaContinuous}. Let now $V$
    be barrelled and $B' \subseteq V'_\beta$ a bounded subset. By what we have already
    argued, we get~\eqref{eq:IotaProof1} for every $\varphi \in B'$. Consider
    $\seminorm{q} \coloneqq \sup_{\varphi \in B'} \seminorm{q}_\varphi$, which is a
    well defined, continuous seminorm on $V$ by the Banach-Steinhaus Theorem. We
    arrive at the continuity estimate
    \begin{align}
        \seminorm{p}_{B'}
        \bigl(
        \iota_n(v)
        \bigr)
        &=
        \sup_{\varphi \in B'}
        \seminorm{p}_\varphi
        \bigl(
        \iota_n(v)
        \bigr) \\
        &\le
        \sup_{\varphi \in B'}
        \sum_k
        \prod_{j=1}^n
        \seminorm{q}_\varphi
        \bigl(
            v_j^{(k)}
        \bigr) \\
        &\le
        \sum_k
        \prod_{j=1}^n
        \sup_{\varphi \in B'}
        \seminorm{q}_\varphi(v_j^{(k)}) \\
        &=
        \sum_k
        \prod_{j=1}^n
        \seminorm{q}(v_j^{(k)}).
    \end{align}
    This proves \ref{item:IotaNBetaContinuous} by once again taking the infimum over
    all decompositions of $v$ into factorizing tensors. Next, we prove that
    $\iota_n(v_1 \vee \cdots \vee v_n)$ is of finite type. As noted before, the
    mapping $\iota_1 \colon V \longrightarrow V''$ is just the usual canonical
    injection, whence there is nothing to be shown in this case. Turning to
    higher degrees, let $v_1, \ldots, v_n \in V$ be given. By the polarization
    identity \eqref{eq:Polarization}, we may write
    \begin{equation}
        \varphi(v_1) \cdots \varphi(v_n)
        =
        \frac{1}{2^n \cdot n!}
        \sum_{\epsilon_j = \pm 1}
        \epsilon_1 \cdots \epsilon_n
        \cdot
        \varphi
        \biggl(
        \sum_{k=1}^{n}
        \epsilon_k v_k
        \biggr)^n
    \end{equation}
    for $\varphi \in V'$, as the left-hand side is linear with respect to $v_1 \tensor
    \cdots \tensor v_n$. Consequently, taking the index set $J \coloneqq \{-1,1\}^n$,
    \begin{equation}
        a_j
        \coloneqq
        \sum_{k=1}^{n}
        j(k)
        v_k
        \qquad \textrm{and} \qquad
        \lambda_j
        \coloneqq
        \frac{j(1) \cdots j(n)}{2^n \cdot n!}
    \end{equation}
    for all $j \in J$ does the job. Finally, it suffices to consider the
    case $R = 0$ in \ref{item:IotaContinuous}, as it is the coarsest of the Hausdorff
    $R$-topologies and biggest as a set by Proposition~\ref{prop:RTopologies}. By
    what we have already shown, we have
    \begin{equation}
        \seminorm{p}_{B'}
        \bigl(
        \iota(v)
        \bigr)
        \le
        \sum_{n=0}^\infty
        \seminorm{p}_{B'}
        \bigl(
        \iota_n(v_n)
        \bigr)
        \le
        \sum_{n=0}^\infty
        \seminorm{q}^n
        \bigl(v_n)
        =
        \seminorm{q}_{0,1}
        \bigl(\iota(v)\bigr)
    \end{equation}
    for $v = \sum_{n=0}^{\infty} v_n \in \Sym_0^\bullet(V)$ with homogeneous
    components $v_n \in \Sym^n(V)$ for all $n \in \N_0$.
\end{proof}

Polynomials of finite type have been studied extensively, beginning with the work of
Aron-Schottenloher \cite{aron.schottenloher:1976} and we refer to
\cite[Sec.~2.6]{dineen:1999a} for a comprehensive overview over the literature.

Our next goal is to investigate the flavour of functions the unique continuous linear
extension
\begin{equation}
    \gls{CanonicalEmbeddingExtension}
    \colon
    \widehat{\Sym}_0(V)
    \longrightarrow
    \Map(V_\beta')
\end{equation}
of \eqref{eq:Iota} produces for locally convex Hausdorff spaces $V$ and $R \ge 0$. To
this end, we first observe that
\begin{equation}
    \label{eq:IotaExtension}
    \hat{\iota}
    \biggl(
    \sum_{n=0}^\infty
    v_n
    \biggr)
    =
    \sum_{n=0}^\infty
    \hat\iota(v_n)
    =
    \sum_{n=0}^\infty
    \hat\iota_n(v_n).
\end{equation}
for vectors $v = \sum_{n=0}^\infty v_n \in
\widehat{\Sym}_0^\bullet(V)$ with components $v_n \in \widehat{\Sym}^n_\pi(V)$.
To see this, recall that the partial sums of the series $\sum_{n=0}^\infty v_n$ actually
converge to $v$ in the~$\Sym_R$-topology and each of the restrictions $\hat{\iota}_n$
is continuous, as well, and given by extending the restrictions~$\iota_n$
of~$\iota$.
For this it is crucial that the subspace topology induced by the inclusion
$\Sym^n(V) \subseteq \Sym_0^\bullet(V)$ is the~$\pi$-topology by
Proposition~\ref{prop:RTopologies},~\ref{item:RTopologiesSubspace}. Moreover, we
observe that already the homogeneous
components~$v_n$ might not actually be linear combinations of factorizing tensors any
more.

Thus there are two limiting procedures to be understood: a \emph{horizontal}
completion within each homogeneous degree and a \emph{vertical} completion of the
symmetric algebra, where the summability condition comes into play. By completeness of
$\Pol_\beta(V_\beta')$ for bornological $V_\beta'$, see again
Theorem~\ref{thm:CompletenessHomogeneousPolynomials}, we get that the extended
components $\hat{\iota}_n$ map again into $\Pol_\beta(V_\beta')$. Moreover, it turns
out that Gâteaux holomorphy is preserved without any assumptions.
\begin{corollary}
    \label{cor:Gateaux}
    Let $V$ be a Hausdorff locally convex space and $v \in \widehat{\Sym}_0(V)$. Then
    \begin{equation}
        \hat\iota(v)
        \colon
        V_\beta'
        \longrightarrow
        \C
    \end{equation}
    is Gâteaux holomorphic with Taylor polynomials at the origin given by
    $\hat\iota_n(v_n)$ for $n \in \N_0$.
\end{corollary}
\begin{proof}
    This is clear by Corollary~\ref{cor:GateauxPowerSeries} and
    Theorem~\ref{thm:TensorsAsPolynomials}, \ref{item:IotaNSigmaContinuous}, where
    we once again use that multilinearity is stable under pointwise limits.
\end{proof}

Our next goal is to understand when $\iota$ is an embedding. We begin with a
preliminary lemma on polars, which is remarkable in its own right.
\begin{lemma}
    \label{lem:PolarsOfBallsBounded}
    \index{Polar!Boundedness}
    Let $V$ be a locally convex space and $\seminorm{p} \in \cs(V)$. Then the polar
    \begin{equation}
        \Ball_{\seminorm{p},1}(0)^\polar
        \subseteq
        V_\beta'
    \end{equation}
    is bounded with respect to the topology of uniform convergence on bounded sets.
\end{lemma}
\begin{proof}
    Recall that a subset $S$ of a locally convex space $W$ is bounded if and only if
    \begin{equation}
        \sup_{v \in S}
        \seminorm{q}(v)
        <
        \infty
        \qquad
        \textrm{for all }
        \seminorm{q}
        \in
        \cs(W).
    \end{equation}
    For $W = V'_\beta$ the continuous seminorms are given by $\seminorm{p}_B$ from
    \eqref{eq:BoundedSeminorms} for all bounded subsets $B \subseteq V$. Thus let $B
    \subseteq V$ be bounded and $\seminorm{p} \in \cs(V)$. By boundedness, there
    exists some radius $r > 0$ with $B \subseteq \Ball_{\seminorm{p},r}(0)$. Thus,
    also investing
    \eqref{eq:PolarAlternative}, we get
    \begin{equation}
        \sup_{\varphi \in \Ball_{\seminorm{p},1}(0)^\polar}
        \seminorm{p}_B(\varphi)
        =
        \sup_{\varphi \in \Ball_{\seminorm{p},1}(0)^\polar}
        \sup_{v \in B}
        \abs[\big]
        {\varphi(v)}
        \le
        \sup_{v \in B}
        \seminorm{p}(v)
        \le
        r
        <
        \infty,
    \end{equation}
    proving the boundedness of $\Ball_{\seminorm{p},1}(0)^\polar$.
\end{proof}

The final ingredient we need to guarantee that $\iota$ is an embedding turns out to
be nuclearity. We will use this property mostly as a blackbox. For us, its crucial
consequence is that, on the tensor product of nuclear spaces, the locally convex
topologies generated by projective and injective tensor products coincide, see again
Remark~\ref{rem:InjectiveTensor}.
\begin{theorem}[Embedding I]
    \label{thm:IotaEmbedding}
    \index{Tensors as polynomials!Embedding I}
    Let $V$ be a barrelled Hausdorff locally convex space. Assume
    furthermore that $V$ is nuclear. Then the mapping
    \begin{equation}
        \label{eq:IotaEmbedding}
        \iota
        \colon
        \Sym_0^\bullet(V)
        \longrightarrow
        \Pol_\beta^\bullet(V'_\beta)
    \end{equation}
    is a grading preserving embedding, i.e. a linear injection such that the subspace
    topology inherited from the inclusion of its image into
    $\Pol^\bullet_\beta(V'_\beta)$ coincides with the $0$-topology.
\end{theorem}
\begin{proof}
    We have just proved the injectivity of $\iota$ in Lemma~\ref{lem:IotaInjectivity}.
    By Theorem~\ref{thm:TensorsAsPolynomials},~\ref{item:IotaContinuous}, our
    mapping \eqref{eq:IotaEmbedding} is continuous. It remains to assert that the
    subspace topology inherited from the inclusion $\iota(\Sym_0^\bullet(V))
    \subseteq \Pol_\beta^\bullet(V'_\beta)$ is finer than the $\Sym_0$-topology.

    To this end, let $\seminorm{p} \in \cs(V) \setminus \{0\}$. By nuclearity, there exists
    another seminorm $\seminorm{p}' \in \cs(V)$ with $\seminorm{p} \le \seminorm{p}'$
    such that
    \begin{equation}
        \seminorm{p} \tensor_\pi \seminorm{q}
        \le
        \nu
        \cdot
        \bigl(
        \seminorm{p} \tensor_\epsilon \seminorm{q}
        \bigr)
        \qquad
        \textrm{for all }
        \seminorm{q} \in \cs(V),
    \end{equation}
    where $\nu$ corresponds to the nuclear norm of the canonical mapping
    $V_\seminorm{p'} \longrightarrow V_{\seminorm{p}}$
    between the local Banach spaces, see \cite[Thm.~6.38]{vogt:2000a}. Crucially,
    $\nu$ thus \emph{only} depends on the seminorm on the left-hand side. For tensor
    powers, we claim that this implies
    \begin{equation}
        \label{eq:NuclearProjectiveVsInjective}
        \seminorm{p}^n
        \le
        \nu^{n-1}
        \cdot
        \seminorm{p}
        \tensor_\epsilon \cdots \tensor_\epsilon
        \seminorm{p}
        \qquad
        \textrm{for all }
        n \in \N.
%        \tag{\Radiation}
    \end{equation}
    Indeed, for $n = 1$, this is clear and, in view of
    Proposition~\ref{prop:ProjectiveTensor}, \ref{item:ProjectiveAssociativity} and the
    analogous statement for injective tensor powers\footnote{Which, taking another look
    at \eqref{eq:InjectiveSeminorms}, is indeed clear.}, we get by induction
    \begin{equation}
        \seminorm{p}^n
        =
        \seminorm{p}
        \tensor_\pi
        \seminorm{p}^{n-1}
        \le
        \nu^{n-2}
        \cdot
        \seminorm{p}
        \tensor_\pi
        \bigl(
            \seminorm{p}
            \tensor_\epsilon \cdots \tensor_\epsilon
            \seminorm{p}
        \bigr)
        \le
        \nu^{n-1}
        \cdot
        \seminorm{p}
            \tensor_\epsilon \cdots \tensor_\epsilon
        \seminorm{p},
    \end{equation}
    proving \eqref{eq:NuclearProjectiveVsInjective}. Consider now the polar $B'
    \coloneqq \Ball_{\seminorm{p},1}(0)^\polar$ as defined
    in \eqref{eq:Polar}. By Corollary~\ref{cor:Polars}, every functional $\varphi \in
    B'$ fulfils
    $\abs{\varphi} \le p$, and by Lemma~\ref{lem:PolarsOfBallsBounded} the set $B'$ is
    bounded within $V'_\beta$. Let
    \begin{equation}
        v = \sum_{n=0}^{\infty} v_n \in
        \Sym^\bullet_0(V)
        \qquad
        \textrm{with homogeneous components} \quad
        v_n \in \Sym_\pi^n(V)
    \end{equation}
    and write $P_n \coloneqq \iota_n(v_n) \in \Pol^n(V'_\beta)$. Using
    \eqref{eq:NuclearProjectiveVsInjective} and \eqref{eq:PolarizationEstimate} for
    $B'$, where we use that polars are always absolutely convex by the triangle
    inequality, leads to the estimate
    \begin{align}
        \seminorm{p}^n
        \bigl(
            v_1 \vee \cdots \vee v_n
        \bigr)
        &\le
        \nu^{n-1}
        \cdot
        \seminorm{p}
        \tensor_\epsilon \cdots \tensor_\epsilon
        \seminorm{p} \\
        &=
        \nu^{n-1}
        \sup_{\varphi_1, \ldots, \varphi_n \in B'}
        \abs[\Big]
        {
            \bigl(
                \varphi_1 \tensor \cdots \tensor \varphi_n
            \bigr)(v_n)
        } \\
        &=
        \nu^{n-1}
        \sup_{\varphi_1, \ldots, \varphi_n \in B'}
        \abs[\Big]
        {
            \bigl(
            \varphi_1 \vee \cdots \vee \varphi_n
            \bigr)(v_n)
        } \\
        &=
        \nu^{n-1}
        \sup_{\varphi_1, \ldots, \varphi_n \in B'}
        \abs[\Big]
        {
            \widecheck{P}_n
            \bigl(
                \varphi_1, \ldots, \varphi_n
            \bigr)
        } \\
        &\le
        \frac{\nu^{n-1} \cdot n^n}{n!}
        \sup_{\varphi \in B'}
        \abs[\big]
        {P_n(\varphi)} \\
        &=
        \frac{\nu^{n-1} \cdot n^n}{n!}
        \cdot
        \seminorm{p}_{B'}(P_n)
    \end{align}
    for all $n \in \N$. Putting everything together, setting
    $C \coloneqq \max\{\nu,1\}$ and
    using \eqref{eq:CauchyEstimateFactorial} yields the desired inequality
    \begin{align}
        \seminorm{p}_{0,c}(v)
        &=
        \sum_{n=0}^{\infty}
        c^n
        \cdot
        \seminorm{p}^n(v_n) \\
        &\le
        \abs{v_0}
        +
        \sum_{n=1}^{\infty}
        (C \cdot c)^n
        \cdot
        \frac{n^n}{n!}
        \cdot
        \seminorm{p}_{B'}(P_n) \\
        &\le
        \sum_{n=0}^{\infty}
        (C \cdot c \cdot e)^n
        \cdot
        \seminorm{p}_{B'}(P_n) \\
        &=
        \seminorm{r}_{mB',\abs{\argument},0}
        \bigl(
            \iota(v)
        \bigr)
    \end{align}
    with $m \coloneqq C \cdot c \cdot e$ and the continuous seminorm
    $\seminorm{r}_{mB',\abs{\argument},0}$
    from \eqref{eq:BoundedAlternativeSeminorms}.
\end{proof}

In the preceding proof we have used the nuclearity of $V$ to pass from the projective to
the injective tensor product within each symmetric tensor degree. This raises the
question whether one could have worked with injective tensor powers instead of
projective ones from the start to construct an alternative $\Sym_R$-topology compatible
with $\iota$ also without nuclearity. This turns out to be feasible.

In analogy to \eqref{eq:RTopologySeminorms}, given $\seminorm{p} \in \cs(V)$, we write
$\seminorm{p}^{\tensor_\epsilon n}$ for its $n$-th injective tensor power and define for
any $c \ge 0$ the seminorm
\begin{equation}
    \index{Seminorms!Injective $0$ topology}
    \label{eq:InjectiveSym0Seminorms}
    \tilde{\seminorm{p}}_{0,c}
    \coloneqq
    \sum_{n=0}^{\infty}
    c^n
    \cdot
    \seminorm{p}^{\tensor_\epsilon n}
\end{equation}
on $\Sym^\bullet(V)$, which like \eqref{eq:RTopologySeminorms} is only formally an
infinite series. Finally, we denote the corresponding locally convex space by
\gls{RTopologyInjective}.\footnote{Of course, one can in principle introduce the
parameter $R$ as before, but not going to this generality allows us to simplify the
bookkeeping considerably.} By Lemma~\ref{lem:ProjectiveVsInjective}, this results in a
coarser topology than the $0$-topology. Taking another look at
Theorem~\ref{thm:TensorsAsPolynomials}, we thus have to revisit the continuity of
$\iota$, now as a mapping defined on $\Sym_\epsilon^\bullet(V)$.
\begin{theorem}[Embedding II]
    \label{thm:IotaEmbedding2}
    \index{Tensors as polynomials!Embedding II}
    Let $V$ be a barrelled Hausdorff locally convex space.
    Then the mapping
    \begin{equation}
        \label{eq:IotaEmbedding2}
        \iota
        \colon
        \Sym_\epsilon^\bullet(V)
        \longrightarrow
        \Pol_\beta^\bullet(V'_\beta)
    \end{equation}
    is a grading preserving embedding, i.e. a linear injection such that the subspace
    topology inherited from the inclusion of its image into
    $\Pol^\bullet_\beta(V'_\beta)$ coincides with the $0$-topology.
\end{theorem}
\begin{proof}
    It is instructive to proceed as in the proof of
    Theorem~\ref{thm:TensorsAsPolynomials}, i.e. by first considering each component
    mapping $\iota_n$ and $\Pol_\sigma^n(V'_\beta)$. To this end, let $\varphi \in V'$ with
    corresponding seminorms $\seminorm{q}_\varphi \coloneqq \abs{\varphi} \in
    \cs(V)$ and $\seminorm{p}_\varphi \in \cs(\Pol_\sigma(V'_\beta))$
    given by $\seminorm{p}_\varphi(Q) = \abs{Q(\varphi)}$ for $Q \in
    \Pol(V'_\beta)$. Then we have $\varphi \in
    \Ball_{\seminorm{q}_\varphi,1}(0)^\polar$ by \eqref{eq:PolarAlternative} and thus
    \begin{align}
        \seminorm{p}_\varphi
        \bigl(
            \iota_n(v)
        \bigr)
        &=
        \abs[\Big]
        {
            \bigl(
                \varphi \tensor \cdots \tensor \varphi
            \bigr)
            (v)
        } \\
        &\le
        \sup_{\varphi_1, \ldots, \varphi_n \in \Ball_{\seminorm{q}_\varphi,1}(0)^\polar}
        \abs[\Big]
        {
            \bigl(
            \varphi_1 \tensor \cdots \tensor \varphi_n
            \bigr)
            (v)
        } \\
        &=
        \bigl(
            \seminorm{q}_\varphi
            \tensor_\epsilon \cdots \tensor_\epsilon
            \seminorm{q}_\varphi
        \bigr)(v)
    \end{align}
    for all  $v \in \Sym^n(V)$. This means that
    \begin{equation}
        \iota_n
        \colon
        \Sym_\epsilon^n(V)
        \longrightarrow
        \Pol_\sigma^n(V'_\beta)
    \end{equation}
    is continuous. If now $V$ is barrelled and $B' \subseteq V'$ some bounded set, we
    may replace $\seminorm{q}_\varphi$ with $\seminorm{q} \coloneqq \sup_{\varphi \in
    B'} \seminorm{q}_\varphi$ to prove the continuity also of
    \begin{equation}
        \iota_n
        \colon
        \Sym_\epsilon^n(V)
        \longrightarrow
        \Pol_\beta^n(V'_\beta).
    \end{equation}
    By \eqref{eq:InjectiveSym0Seminorms}, the corresponding continuity estimates then
    glue to a continuity estimate of \eqref{eq:IotaEmbedding2}. As we have already
    argued that the remaining statements hold also in this alternative setting, this
    completes the proof.
\end{proof}

\index{Questions!Star products vs. injective tensor products}
It would be interesting to investigate the continuity of the star products in
\cite{waldmann:2014a, schoetz.waldmann:2018a, heins.roth.waldmann:2023a} after
replacing $\Sym_0^\bullet(V)$ with the appropriately parametrized version of
$\Sym_\epsilon^\bullet(V)$. In view of the absence of the infimum argument, see again
Remark~\ref{rem:InjectiveTensor}, this might however be prohibitively difficult.
\begin{question}
    \index{Open questions!Injective $R$-topology}
    Are the continuous star products from \cite{waldmann:2014a} continuous with
    respect to the~$R$-topologies constructed from injective tensor products?
\end{question}

We return to the investigation of the vertical completion of
$\Sym_0^\bullet(V)$ in terms of bounded entire Fréchet holomorphic functions.
\begin{corollary}
    Let $V$ be a bornological and barrelled Hausdorff locally convex space. Assume
    furthermore that $V$ is nuclear and that $V'_\beta$ is Baire. Then each element of
    $\widehat{\Sym}_0(V)$ may be identified with a unique entire bounded Fréchet
    holomorphic function on $V'_\beta$.
\end{corollary}
\begin{proof}
    Composing \eqref{eq:IotaEmbedding} with the further embedding
    \begin{equation}
        \Pol_\beta^\bullet(V'_\beta)
        \subseteq
        \Holomorphic_\beta(V'_\beta,\C)
    \end{equation}
    allows us to identify $\Sym_0^\bullet(V)$ with a subspace $S$ of the complete
    Hausdorff space $\Holomorphic_\beta(V'_\beta,\C)$. Thus, one may realize the
    completion of $\Sym_0^\bullet(V)$ as the closure within
    $\Holomorphic_\beta(V'_\beta,\C)$.

    Alternatively, one may argue as follows. Let $v \in \widehat{\Sym}^\bullet_0(V)$
    with homogeneous components~$v_n \in \widehat{\Sym}_\pi^n(V)$. By
    Corollary~\ref{cor:Gateaux}, we know that $\hat{\iota}(v)$ is Gâteaux holomorphic
    with Taylor polynomials at the origin given by $\hat{\iota}_n(v_n)$.
    Theorem~\ref{thm:TensorsAsPolynomials}, \ref{item:IotaV} and
    \ref{item:IotaNBetaContinuous} in combination with
    Theorem~\ref{thm:CompletenessHomogeneousPolynomials} asserts that
    $\hat{\iota}_n(v_n)$ is continuous for every $n \in \N_0$. By
    Proposition~\ref{prop:FrechetBaire}, this implies Fréchet holomorphy of
    $\hat{\iota}(v)$ on all of $V'_\beta$.
\end{proof}

Our alternative proof shows that the nuclearity is not necessary for the Fréchet
holomorphy of $\hat{\iota}(v)$. However, without nuclearity of $V$, it may in
principle happen
that there are non-convergent nets in $\Sym^\bullet_0(V)$ such that the image net under
$\iota$ converges in $\Holomorphic_\beta(V'_\beta)$. We have seen in
Theorem~\ref{thm:CompletenessBetaFrechet} that the Baire property on
the domain ensures completeness of the space of bounded Fréchet holomorphic functions.
To apply this criterion we run into the complication that $V'_\beta$ is a dual
space, which is at odds with the usual sufficient criteria for the Baire property.
\begin{remark}[Duality]
    \label{rem:Duality}
    \index{Tensors as polynomials!Duality}
    One way around this is to switch the roles of $V$ and $V'$: As they are in
    duality, every $\varphi \in \Sym^\bullet(V')$ with homogeneous components
    $\varphi_n \in \Sym^n_\pi(V')$ defines a polynomial on $V$ by
    means of
    \begin{equation}
        \label{eq:IotaDual}
        v
        \quad \mapsto \quad
        \sum_{n=0}^{\infty}
        \varphi_n
        \bigl(
            \underbrace{v \tensor \cdots \tensor v}_{n-\textrm{times}}
        \bigr).
    \end{equation}
    As the resulting maps are nothing else than the restrictions of maps obtained by
    further dualizing, i.e. by considering
    \begin{equation}
        \iota
        \colon
        \Sym^\bullet_0
        \bigl(
            V'_\beta
        \bigr)
        \longrightarrow
        \Pol
        \bigl(
            (V_\beta')_\beta'
        \bigr),
    \end{equation}
    analogous statements to the ones we have proved hold for
    \eqref{eq:IotaDual}. Here one runs into the slight complication that, a priori,
    the bounded subsets of $V$ and its bidual $(V'_\beta)'_\beta$ may not be the same:
    If $B \subseteq V$ is a bounded subset of a barrelled and bornological space, then
    \begin{equation}
        \iota_1(B)
        \subseteq
        \Pol_\beta(V_\beta')
        =
        (V_\beta')_\beta'
    \end{equation}
    is bounded by Theorem~\ref{thm:TensorsAsPolynomials},
    \ref{item:IotaNBetaContinuous} and
    Lemma~\ref{lem:PolynomialContinuityVsBoundedness},
    \ref{item:PolynomialBoundednessImpliesContinuity}. However, in general there are
    many more bounded subsets of $(V'_\beta)'_\beta$. Nevertheless, as the properties
    of $V$ and $V'_\beta$ are typically vastly different e.g. in the setting of Fréchet
    spaces, swapping roles might be useful. Another particular case arises if $V$
    admits for a pre-dual. Then one does not even have to restrict. We will see these
    duality phenomenona in action in Example~\ref{ex:SequenceSpace} and
    Example~\ref{ex:Continuous}.
\end{remark}

Even in the situation of Theorem~\ref{thm:IotaEmbedding}, the immediate
question is however the following:
\begin{question}
    How large is the closure of $\iota(\Sym_0(V))$ within
    $\Holomorphic_\beta(V_\beta')$?
\end{question}

\begin{remark}[Reflexivity]
    \index{Reflexivity}
    \index{Tensors as polynomials!Reflexivity}
    Approaching this question, it is notable that it is interesting already for the
    component mappings $\iota_n$. As we have noted and used before, $\iota_1$ is
    nothing else than the canonical inclusion of $V$ into its strong bidual
    $(V'_\beta)'_\beta$. Surjectivity of~$\iota_1$ for barrelled space $V$ thus
    corresponds to \emph{reflexivity} of $V$, as we have already proved the continuity
    of $\iota_1$ at least for barrelled $V$ in Theorem~\ref{thm:TensorsAsPolynomials},
    \ref{item:IotaNBetaContinuous}. As such, the continuity of $\iota_1$ is well-known
    and can be found in \cite[5.5 \& 5.6]{schaefer:1999a}.
\end{remark}

Combining reflexivity with our considerations on duality from Remark~\ref{rem:Duality}
gives rise to a remarkable class of examples.\footnote{I would like to thank Stefan
Waldmann for pointing this out.}
\begin{example}[Reflexive nuclear Fréchet spaces]
    \index{Reflexivity!Nuclear Fréchet spaces}
    Let $F$ be a nuclear reflexive Fréchet space such as:
    \begin{examplelist}
        \item The algebra of smooth functions $\Cinfty(\Omega)$ defined on some
        open subset~$\emptyset \neq \Omega \subseteq \R^n$ with values in $\C$, see
        \cite[Example~28.9~(1)]{meise.vogt:1992a} and \cite[Prop~36.10]{treves:2006a}.
        \item The Schwartz space $\Schwartz(\R^n)$ of complex valued rapidly
        decreasing smooth functions on $\R^n$, which we will discuss in the context of
        smooth vectors in Example~\ref{ex:Schroedinger}, see
        \cite[Prop~36.10 \& Thm.~51.5]{treves:2006a}.
        \item The algebra of holomorphic functions $\Holomorphic(\Omega)$ defined on
        some complex nonempty domain~$\Omega \subseteq \C^n$ with values in
        $\C$, see
        \cite[Example~28.9~(4)]{meise.vogt:1992a} and \cite[Corollary of
        Prop.~36.9 or Prop.~36.10]{treves:2006a}, both of which are applicable by
        virtue of Montel's Theorem.
    \end{examplelist}
    Consider $V \coloneqq F'_\beta$. Then, by reflexivity, $V' =
    (F'_\beta)_\beta' \cong F$ is Baire as a
    Fréchet space. Moreover, $V$ is nuclear as the strong dual of a nuclear space by
    \cite[Prop.~50.6]{treves:2006a}, and bornological as the strong dual of a
    reflexive Fréchet space. Finally, $F$ is Montel as a nuclear Fréchet space by
    \cite[Prop.~50.2]{treves:2006a} and thus its strong dual $F'$ is also Montel by
    \cite[Prop.~36.10]{treves:2006a}, which in particular
    implies that $F'$ is barrelled itself by \cite[Cor.~of Prop~36.9 \&
    Prop.~36.4]{treves:2006a}. Thus all assumptions of
    Theorem~\ref{thm:TensorsAsPolynomials} and
    Theorem~\ref{thm:IotaEmbedding} are fulfilled.
\end{example}

    Going to higher orders, surjectivity of
    \begin{equation}
        \label{eq:IotaNExtensions}
        \hat\iota_n
        \colon
        \widehat{\Sym}^n_\pi(V)
        \longrightarrow
        \Pol^n_\beta(V)
    \end{equation}
    may then be seen as a natural definition of $n$-reflexivity. This immediately
    opens several avenues for further research such as:
    \index{Questions!Reflexivity}
    \begin{question}
        Does $n$-reflexivity imply $k$-reflexivity for all smaller (greater) $k \in
        \N$?
    \end{question}
    \begin{question}
        Does $n$-reflexivity for all $n \in \N$ imply surjectivity of $\hat\iota$?
    \end{question}
    \begin{question}
        Does surjectivity of $\hat\iota$ imply $n$-reflexivity?
    \end{question}

    As we have seen in the abstract setting in Section~\ref{sec:HolomorphicFrechet}
    and as
    we will again see in Section~\ref{sec:EntireVectorSpacesExamples}, Taylor series
    need not converge in the $\beta$-topology. Thus the relationship between
    $\hat\iota$ and its components \eqref{eq:IotaNExtensions} not clear: the
    polynomials in some approximating net may change completely in every step. One
    should think of this like the distinction between having a \emph{convergent power
    series} and \emph{approximation by polynomials}. Already on compact subsets of the
    real line $\R$, the former implies holomorphy, the latter only continuity.

\begin{remark}[Positive $R$ and Gelfand-Shilov]
    \index{Gelfand-Shilov!Spaces}
    \index{Gelfand, Isreal}
    \index{Shilov, Georgiy}
    In \cite[Sec.~5.4]{waldmann:2014a} the particular case of finite dimensional $V$
    was investigated. Similarly to what we have done for $\Sym_0^\bullet(V)$ in
    Theorem~\ref{thm:IotaEmbedding}, one may identify $\Sym_R^\bullet(V)$ with the
    space of entire functions of finite order~$1/R$ and minimal type for any $R \ge
    0$. The order constitutes a growth condition on the Taylor coefficients, which is
    most apparent in \eqref{eq:RTopologySeminorms}, and may be rephrased as a growth
    condition
    of the full function by the methods we have employed in
    Proposition~\ref{prop:BetaAlternativeSeminorms}.
\end{remark}

\section{A Collection of Examples}
\epigraph{Ook!}{The Librarian}
\label{sec:EntireVectorSpacesExamples}
% !TeX root = ../Dissertation.tex

It is instructive to study how one can understand \eqref{eq:Iota} in examples, where
the continuous dual $V'_\beta$ may be described explicitly. Of particular interest are
the questions of surjectivity of $\hat{\iota}_n$ and $\hat{\iota}$ itself. We begin
with a sanity check.
\begin{example}[Finite dimensions]
    Let $V = \C^n$. Then every $n$-linear mapping on the dual space $V' \cong V$ is
    continuous. Consequently, $\Pol^n(V') \cong \C[z_1, \ldots, z_n]$ with
    \begin{equation}
        \iota_n
        \bigl(
            e_{j_1}
            \vee
            \cdots
            \vee
            e_{j_n}
        \bigr)
        =
        z_{j_1} \cdots z_{j_n}.
    \end{equation}
    Moreover, the $\beta$-topology reduces to the topology of locally uniform
    convergence. Thus the subspace $\iota(\Sym \C^n) \subseteq \Holomorphic(\C^n)$ is
    dense by virtue of Taylor expansions.
\end{example}

As our next examples will show, approximation by Taylor series does not always work.
\begin{example}[Non reflexive sequence spaces]
    \label{ex:SequenceSpace}
    \index{Zero sequences}
    \index{Sequence space!Zero sequences}
    \index{Sequence space!Absolutely summable sequences}
    We consider the Banach space of zero sequences
    \begin{equation}
        \gls{ZeroSequences}
        \coloneqq
        \bigl\{
            \gamma
            \colon
            \field{N}
            \longrightarrow
            \field{C}
            \colon
            \lim_{n \rightarrow \infty}
            \gamma_n
            =
            0
        \bigr\}
    \end{equation}
    endowed with the supremum norm $\gls{Supnorm} \coloneqq \sup_{n \in \N}
    \abs{\gamma_n}$. By \cite[Satz~II.2.3]{werner:2002a}, we have
    \begin{equation}
        c_0'
        \cong
        \gls{AbsolutelySummable}
        \coloneqq
        \biggl\{
            \gamma
            \colon
            \field{N}
            \longrightarrow
            \field{C}
            \colon
            \gls{Ell1Norm}
            \coloneqq
            \sum_{n=1}^{\infty}
            \abs[\big]{\gamma_n}
            <
            \infty
        \biggr\}
    \end{equation}
    and
    \begin{equation}
        (\ell^1)'
        \cong
        \gls{BoundedSequences}
        \coloneqq
        \bigl\{
        \gamma
        \colon
        \field{N}
        \longrightarrow
        \field{C}
        \colon
        \supnorm{\gamma}
        <
        \infty
        \bigr\}
    \end{equation}
    as Banach spaces by means of the natural pairings
    \begin{equation}
        \langle
            \gamma,
            s
        \rangle
        \coloneqq
        \sum_{n=1}^{\infty}
        \gamma_n \cdot s_n
    \end{equation}
    with $\gamma \in c_0$ resp. $\gamma \in \ell^\infty$ and $s \in \ell^1$. In
    particular, $c_0$ is not reflexive, as the canonical mapping
    $\iota_1 \colon c_0 \longrightarrow \ell^\infty$ is simply given by the set
    inclusion $c_0 \subseteq \ell^\infty$. Consequently, $\iota_1$ is not surjective.

    Investigating the second order, we claim that the completed projective tensor
    square
    \begin{equation}
        c_0 \gls{TensorCompletion} c_0
    \end{equation}
    may be viewed as a dense, but proper, subspace
    of the space $c_0(c_0)$ of zero sequences with entries in $c_0$. An element
    $\Lambda \in c_0(c_0)$ is a sequence $(\Lambda_n)_n$ of zero sequences, which
    converges uniformly to the zero sequence. The latter translates to the condition
    \begin{equation}
        \label{eq:ZeroSequenceOfSequences}
        0
        =
        \lim_{n \rightarrow \infty}
        \supnorm{\Lambda_n}
        =
        \lim_{n \rightarrow \infty}
        \sup_{m \in \N}
        \abs[\big]
        {\Lambda_n(m)},
    \end{equation}
    whereas the former implies
    \begin{equation}
        \label{eq:SequenceOfZeroSequences}
        \lim_{m \rightarrow \infty}
        \Lambda_n(m)
        =
        0
        \qquad
        \textrm{for all }
        n \in \N.
    \end{equation}
    In particular, both limits commute and one may visualize $\Lambda$ as the infinite
    matrix
    \begin{equation}
        \begin{pmatrix}
            \Lambda_1(1) & \Lambda_1(2) & \ldots & 0 \\
            \Lambda_2(1) & \Lambda_2(2) & \ldots & 0 \\
            \ldots & \ldots & \ldots & \ldots \\
            0 & 0 & \ldots & 0
        \end{pmatrix}
    \end{equation}
    such that any path ending on the right or bottom edge corresponds to a zero
    sequence. This also means that, a posteriori, both indices exhibit the same
    behaviour and thus we may view $\Lambda$ as a doubly indexed sequence without
    losing information. In particular, the norm of $c_0(c_0)$ is by slight abuse of
    notation simply given by
    \begin{equation}
        \supnorm{\Lambda}
        \coloneqq
        \sup_{n \in \N}
        \supnorm{\Lambda_n}
        =
        \sup_{n,m \in \N}
        \abs[\big]
        {\Lambda_{n,m}}.
    \end{equation}
    We consider the injective linear mapping
    \begin{equation}
        \label{eq:ZeroSequenceTensorSquare}
        \chi
        \colon
        c_0 \tensor_\pi c_0
        \longrightarrow
        c_0(c_0), \quad
        \chi(\gamma \tensor \eta)
        \at[\Big]{(n,m)}
        \coloneqq
        \gamma(n) \cdot \eta(m).
    \end{equation}
    Then $\chi(\gamma \tensor \eta)$ fulfils both \eqref{eq:ZeroSequenceOfSequences}
    and \eqref{eq:SequenceOfZeroSequences}, as writing $\Lambda \coloneqq \chi(\gamma
    \tensor \eta)$ we get
    \begin{equation}
        \lim_{n \rightarrow \infty}
        \supnorm{\Lambda_n}
        =
        \lim_{n \rightarrow \infty}
        \sup_{m \in \N}
        \abs{\Lambda_n(m)}
        =
        \lim_{n \rightarrow \infty}
        \sup_{m \in \N}
        \abs[\big]
        {\gamma_n}
        \cdot
        \abs[\big]
        {\eta_m}
        =
        \lim_{n \rightarrow \infty}
        \abs[\big]
        {\gamma_n}
        \cdot
        \supnorm{\eta}
        =
        0
    \end{equation}
    and
    \begin{equation}
        \lim_{m \rightarrow \infty}
        \Lambda_n(m)
        =
        \lim_{m \rightarrow \infty}
        \gamma_n
        \cdot
        \eta_m
        =
        \gamma_n
        \cdot
        \lim_{m \rightarrow \infty}
        \eta_m
        =
        \gamma_n
        \cdot
        0
        =
        0
    \end{equation}
    for all $n \in \N$. Thus, \eqref{eq:ZeroSequenceTensorSquare} is well defined. We
    moreover have the continuity estimate
    \begin{equation}
        \supnorm[\big]
        {\chi(\gamma \tensor \eta)}
        =
        \sup_{n,m \in \N}
        \abs[\big]
        {\gamma_n}
        \cdot
        \abs[\big]
        {\eta_m}
        =
        \supnorm{\gamma}
        \cdot
        \supnorm{\eta}
    \end{equation}
    on factorizing tensors and Proposition~\ref{prop:InfimumArgument} implies the
    continuity of $\chi$.

    Next, we check that the image of $\chi$ is dense. To this end, fix $\Lambda \in
    c_0(c_0)$ and $\epsilon > 0$. By our preliminary considerations, there exists an
    index $N \in \N$ such that $\abs{\Lambda_{n,m}} \le \epsilon$ for all $n,m \ge N$.
    Using the standard unit vectors $\basis{e}_n \in c_0$ with entries
    $(\basis{e}_n)(m) \coloneqq \delta_{n,m}$ we moreover note
    \begin{equation}
        \label{eq:SequencesDenseImage}
        \chi
        \bigl(
            \basis{e}_n
            \tensor
            \Lambda_{n,\argument}
        \bigr)
        \at[\Big]{(k,\ell)}
        =
        \delta_{n,k}
        \cdot
        \Lambda_{n,\ell}
    \end{equation}
    for all $n,k,\ell \in \N$. In terms of infinite matrices, one may visualize
    $\chi(\basis{e}_n \tensor \Lambda_{n,\argument})$ as the~$n$-th row of $\Lambda$
    with zeros otherwise. Consequently,
    \begin{equation}
        \Gamma
        \coloneqq
        \sum_{n=1}^{N}
        \chi
        \bigl(
            \basis{e}_n
            \tensor
            \Lambda_{n,\argument}
        \bigr)
        =
        \begin{pNiceArray}{c|c}
            \Lambda & \Lambda \\
            \hline
            0 & 0
        \end{pNiceArray}
    \end{equation}
    fulfils $\supnorm{\Lambda - \Gamma} \le \epsilon$ by choice of $N$.
    Hence, the image of \eqref{eq:ZeroSequenceTensorSquare} is dense.

    Recall now the explicit description of the completion $c_0 \tensor_\pi c_0$ as it
    can be found in \cite[Thm.~45.1]{treves:2006a}. That is, for every $\Gamma
    \in c_0 \widehat{\tensor}_\pi c_0$ there are two sequences $\gamma \in
    \ell^1(c_0)$ and~$\eta \in \ell^\infty(c_0)$
    such that
    \begin{equation}
        \Gamma
        =
        \sum_{k=1}^{\infty}
        \gamma_k
        \tensor
        \eta_k.
    \end{equation}
    This has two consequences: Firstly, the continuous extension of $\chi$ is to $c_0
    \widehat{\tensor}_\pi c_0$ is still injective. Secondly, the series
    \begin{equation}
        \sum_{m,n=1}^{\infty}
        \chi
        \bigl(
            \basis{e}_m
            \tensor
            \Lambda_{n,\argument}
        \bigr)
    \end{equation}
    is \emph{not} convergent in $c_0 \widehat{\tensor}_\pi c_0$ unless
    $\Lambda_{n,\argument}$ is actually absolutely summable. Indeed, replacing zero
    sequences with absolutely summable sequences throughout, i.e. dualizing, proves
    the well-known identities
    \begin{equation}
        \label{eq:Ell1TensorSquare}
        \ell^1
        \widehat{\tensor}_\pi
        \ell^1
        \cong
        \ell^1(\ell^1)
        \cong
        \ell^1(\N^2)
        \cong
        \ell^1.
    \end{equation}
    We will meet some of these identifications from a rather different point of view
    and in higher detail in Lemma~\ref{lem:ProjectiveTensorOfEll1} and refer to
    \cite[Ch.~46]{treves:2006a} and \cite[Sec.~15.7]{jarchow:1981a} for a discussion
    of completed projective tensor products with $\Lone$-spaces associated to some
    measure.
    Our discussion may then be recovered by considering the particular case of
    the counting measure.

    Moreover, \eqref{eq:Ell1TensorSquare}
    gives an abstract reason for $\chi(c_0 \widehat{\tensor}_\pi c_0) \subsetneq
    c_0(c_0)$: strong duals essentially turn injective into projective
    tensor products by \cite[Ch.~49 \& 50]{treves:2006a} or \cite[§45.3]{koethe:1979a}
    and infinite dimensional normed spaces are never nuclear, see again
    Remark~\ref{rem:InjectiveTensor} and \cite[Ch.~50, Cor.~2]{treves:2006a}.
    Either way, using our identifications, $\iota_2$ takes the form
    \begin{equation}
        \iota_2(a)
        \at[\Big]{\gamma}
        =
        \sum_{n,m=1}^{\infty}
        \gamma(n) \Lambda(n,m) \gamma(m)
        \qquad
        \textrm{for all }
        \Lambda \in c_0(c_0), \;
        \gamma \in \ell^1.
    \end{equation}
    Once again, this series still converges for $a \in \ell^\infty$, so $\iota_2$ is
    also not surjective. As noted before, this does not imply that the full mapping $\iota$
    is not surjective, as one might be able to approximate by mixing degrees. It is clear
    that these considerations generalize to higher degrees without complications.
\end{example}

\begin{example}[Continuous Functions]
    \label{ex:Continuous}%
    \index{Continuous functions}
    Let $V = \Continuous([0,1])$ equipped with the norm topology induced by the
    supremum norm $\supnorm{\argument}$. By the Riesz-Markov Theorem, every
    functional $\varphi \in V'$ is given by a unique regular complex Radon
    measure $\mu_\varphi$ on $[0,1]$, i.e.
    \begin{equation}
        \label{eq:RieszMarkov}
        \varphi(f)
        =
        \int_{[0,1]}
        f
        \D \mu_\varphi
        \qquad
        \textrm{for all }
        f \in \Continuous([0,1]).
    \end{equation}
    Moreover, it is a straightforward consequence of the Stone-Weierstraß Theorem, as
    it can e.g. be found in \cite[5.7, p.~122]{rudin:1991a}, that
    \begin{equation}
        \Continuous\bigl([0,1]\bigr)
        \widehat{\otimes}_\pi
        \Continuous\bigl([0,1]\bigr)
        \cong
        \Continuous\bigl([0,1]^2\bigr)
    \end{equation}
    is just the space of continuous functions on the unit square, again equipped with the
    norm topology induced by the corresponding supremum norm: Indeed, the set
    \begin{equation}
        \algebra{A}
        \coloneqq
        \biggl\{
            \Bigl(
                (x,y)
                \mapsto
                \sum_{k=1}^N
                f_k(x) \cdot g_k(y)
            \Bigr)
            \colon
            N \in \N, \;
            f_k, g_k \in \Continuous([0,1])
        \biggr\}
    \end{equation}
    of finite sums of factorizing functions constitutes a point separating
    $^*$-subalgebra of the continuous functions on the unit square
    $\Continuous([0,1] \times [0,1])$. Given a mapping $f \in
    \Continuous([0,1]^2)$ with $f(x,y) = f(y,x)$ for all $x,y \in [0,1]$, it is
    thus reasonable to consider the function~$\iota(f)$. Taking another look at
    \eqref{eq:RieszMarkov} suggests that
    \begin{equation}
        \label{eq:ProductMeasure}
        \iota(f)
        \at[\Big]{\varphi}
        =
        \int_{[0,1]^2}
        f
        \D (\mu_\varphi \times \mu_\varphi),
        \qquad
        \varphi \in V'.
    \end{equation}
    Here, $\mu_\varphi \times \mu_\varphi$ denotes the product measure of
    $\mu_\varphi$ with itself. Note that \eqref{eq:ProductMeasure} makes sense
    regardless of the symmetry properties of our function $f$, so we ignore them in the
    sequel. If $f$ is factorizing, i.e. $f(x,y) = \sum_k g_k(x) h_k(y)$ for all points $x,y \in
    [0,1]$ and suitable functions $g_k, h_k \in \Continuous([0,1])$, then
    \eqref{eq:ProductMeasure} reproduces the expected result $\sum_k \varphi(g_k)
    \varphi(h_k)$ by Fubini's Theorem. By density of such factorizing functions and the
    continuity of \eqref{eq:ProductMeasure} as a linear functional on
    $\Continuous([0,1]^2)$, this is thus the correct description also for non-factorizing
    elements. In particular, all of our continuity estimates carry over unharmed.

    Notably, one may conversely consider
    \eqref{eq:ProductMeasure} as a linear functional $\Xi$ on $\Continuous([0,1]^2)$ by
    fixing functionals $\varphi, \xi \in \Continuous([0,1])'$ and setting
    \begin{equation}
        \Xi(f)
        \coloneqq
        \int_{[0,1]^2}
        f
        \D (\mu_\varphi \times \mu_\xi)
        \qquad
        \textrm{for all }
        f \in \Continuous([0,1]^2).
    \end{equation}
    This is, of course, just the duality of $V$ and $V'$ in action. In particular, $\Xi$
    extends to the algebraic bidual $V^{**}$. Finally, this also showcases the manner in
    which the continuous dual of the projective tensor product $V \otimes_\pi V$ is
    isomorphic to the space of quadratic polynomials on $V$, see
    \cite[Sec.~1.2]{dineen:1999a} for a comprehensive abstract discussion.
\end{example}

\begin{example}[Sequence Spaces]
    \label{ex:Sequences}%
    \index{Sequence spaces!All sequences}
    This is motivated by \cite[Example~1.27]{dineen:1999a}. Throughout, we equip the
    complex numbers $\field{C}$ with their unique Hausdorff locally convex topology.
    That is to say, any seminorm on $\field{C}$ is continuous. Let $J$ be an
    index set and
    \begin{equation}
        V
        \coloneqq
        \gls{CartesianPower}
        \coloneqq
        \prod_{j \in J} \field{C}
    \end{equation}
    be endowed with the product topology. We
    conjecture that the strong dual is given by
    \begin{equation}
        \label{eq:SequenceSpaceStrongDual}
        \bigl(\field{C}^{J}\bigr)'_\beta
        \cong
        \gls{CartesianDirectPower}
        \coloneqq
        \bigoplus_{j \in J}
        \field{C},
    \end{equation}
    where the finite sequences $\field{C}^{(J)}$ carry the locally convex direct sum
    topology. Summing over $\field{C}$, this results in the finest locally convex topology
    on $\field{C}^{(J)}$, i.e. also here any seminorm is continuous.

    To prove \eqref{eq:SequenceSpaceStrongDual} for countable\footnote{We stick
        to a general index set $J$ until we actually use the countability, namely for
        the continuity of the inverse mapping.} index sets $J$, we consider the linear
        mapping
    \begin{equation}
        \Phi
        \colon
        \field{C}^{(J)}
        \longrightarrow
        \bigl(\field{C}^{J}\bigr)'_\beta, \quad
        \Phi(\gamma)
        \at[\Big]{c}
        \coloneqq
        \sum_{j \in J}
        \gamma_j
        \cdot
        c_j.
    \end{equation}
    As $\gamma_j = 0$ for almost all $j \in J$, any choice of linear ordering on $J$
    results in the same~$\Phi$, so we suppress this choice in the sequel. The map
    $\Phi$ is clearly injective. Let $\varphi \in (\field{C}^{J})'_\beta$ be given and define
    $\gamma_j \coloneqq \varphi(\basis{e}_j)$, where $(\basis{e}_j)_{k} =
    \delta_{j,k}$ for $j,k \in J$. All
    but finitely many of the coefficients $\gamma_j$ vanish: Consider the sequence $c \in
    \field{C}^J$ with
    \begin{equation*}
        c_j
        \coloneqq
        \begin{cases}
            1/\gamma_j
            &\; \textrm{if} \;
            \gamma_j \neq 0, \\
            0
            &\; \textrm{if} \;
            \gamma_j = 0,
        \end{cases}
    \end{equation*}
    for all $j \in J$ and let $J_0 \coloneqq \{j \in J \colon c_j = 0\}$. Note
    moreover that, ordering
    the finite subsets $F$ of $J$ by inclusion, the net limit
    \begin{equation}
        \lim_{F \subseteq J}
        \sum_{j \in F}
        c_j \basis{e}_j
        \eqqcolon
        \sum_{j \in J}
        c_j \basis{e}_j
    \end{equation}
    exists and reproduces $c$. Thus, by continuity and linearity of $\varphi$, we have
    \begin{equation*}
        \infty
        >
        \varphi(c)
        =
        \sum_{j \in J}
        \varphi(c_j \basis{e}_j)
        =
        \sum_{j \in J}
        c_j
        \varphi(\basis{e}_j)
        =
        \sum_{j \in J \setminus J_0}
        1.
    \end{equation*}
    Hence $J \setminus J_0$ is indeed finite, i.e. $\gamma \in \field{C}^{(J)}$. An analogous computation yields $\Phi(\gamma) = \varphi$ at once. That is, the mapping $\Phi$ is bijective, and we have constructed the inverse explicitly. By the characteristic property of the final topology, it suffices to check the continuity of the mappings
    \begin{equation*}
        \phi_j
        \colon
        \field{C}
        \longrightarrow
        \bigl(\field{C}^{J}\bigr)'_\beta, \quad
        \phi_j(\lambda)
        \at{c}
        \coloneqq
        \lambda
        \cdot
        c_j
    \end{equation*}
    for all $j \in J$ to infer the continuity of $\Phi$. But this is trivial, as any seminorm on
    $\field{C}$ is continuous.

    We turn towards the continuity of the inverse mapping $\Phi^{-1}$. To facilitate
    our proof we make the additional assumption that the index set $J$ is
    countably infinite. That is, we
    assume $J = \N$. Let $\seminorm{q} = \sum_{j \in \N} \seminorm{q}_j$ be a
    formal sum of seminorms $\seminorm{q}_j$ on $\field{C}$. By continuity, there are
    scalars $\lambda_j \ge 0$ such that $\seminorm{q}_j(\lambda) \le \lambda_j
    \abs{\lambda}$ for all $\lambda \in \field{C}$ and $j \in \N$. Recall that the seminorms
    $\abs{\delta_j}$ -- with the usual Dirac \gls{Dirac} functional $\delta_j(c)
    \coloneqq c_j$ --
    form a defining system of seminorms for~$\field{C}^\N$. Hence, the set
    \begin{equation*}
        B
        \coloneqq
        \bigl\{
        2^j
        \cdot
        \lambda_j
        \cdot
        \basis{e}_j
        \colon
        j \in \N
        \bigr\}
        \subseteq
        \field{C}^\N
    \end{equation*}
    is bounded.\footnote{This is where the countability of the index set $J$ is
    crucial.} Putting both
    observations together yields the estimate
    \begin{equation}
        \seminorm{q}
        \bigl(
            \Phi^{-1}(\varphi)
        \bigr)
        =
        \sum_{j \in \N}
        \seminorm{q}_j
        \bigl(
        \varphi(\basis{e}_j)
        \bigr)
        \le
        \sum_{j \in \N}
        \lambda_j
        \abs[\big]
        {\varphi(\basis{e}_j)}
        \le
        \sum_{j \in \N}
        2^{-j}
        \cdot
        \abs[\big]
        {\varphi(2^j \cdot \lambda_j \cdot \basis{e}_j)}
        \le
        \seminorm{p}_B(\varphi)
    \end{equation}
    for all $\varphi \in \bigl(\field{C}^{\N}\bigr)'_\beta$. As formal sums of seminorms are a
    defining system for the direct sum topology, this finally proves
    \eqref{eq:SequenceSpaceStrongDual}, albeit only for $J = \N$.

    Next, we show that
    \begin{equation}
        \label{eq:SequenceSpacePower}
        \field{C}^J \widehat{\otimes}_\pi \field{C}^J
        \cong
        \field{C}^{J \times J}.
    \end{equation}
    To this end, we consider the linear extension of the mapping
    \begin{equation}
        \Psi
        \colon
        \field{C}^J \otimes_\pi \field{C}^J
        \longrightarrow
        \field{C}^{J \times J}, \quad
        \Psi(c \tensor d)_{jk}
        \coloneqq
        c_j \cdot d_k
    \end{equation}
    and first prove that it is continuous with dense image. By the characteristic property of the product topology, the continuity of the linear mappings $\psi_{jk} \colon \field{C}^J \otimes_\pi \field{C}^J \longrightarrow \field{C}$,
    \begin{equation*}
        \psi_{jk}(c \tensor d)
        \coloneqq
        \bigl(
            \pr_{jk} \circ \Psi
        \bigr)
        (c \tensor d)
        =
        c_j \cdot d_k
        =
        \bigl(
            \delta_j \tensor \delta_k
        \bigr)
        (c \tensor d)
    \end{equation*}
    for all $j,k \in J$ implies continuity of $\Psi$. Here we have used that tensor
    product of continuous linear mappings are themselves continuous on the projective
    tensor product. This general fact is readily proved by combining
    Proposition~\ref{prop:InfimumArgument} with \eqref{eq:ProjectiveOnFactorizing}.

    To see the density of the image, we shift our perspective, disregard the involved
    topologies for a moment and consider $\field{C}^{J \times J}$ as the continuous
    functions $\Continuous(J \times J)$, where we endow the index set $J
    \times J$ with the discrete topology. Notice now that the image $\Psi(\field{C}^J
    \otimes_\pi \field{C}^J)$ is a unital $^*$-subalgebra of $\Continuous(J \times J)$. By
    the Stone-Weierstraß Theorem, it is thus dense in $\Continuous(J
    \times J)$ with respect to the topology of convergence on compact subsets
    of $J \times J$. That is, by discreteness, the topology of pointwise convergence.
    Passing back to $\field{C}^{J \times J}$, this is exactly the product topology. Hence
    our map $\Psi$ indeed has dense image.

    Alternatively and more explicitly, given $a \in \field{C}^{J \times J}$ one may
    consider the sequences of one index given by $(a_{j \bullet})_k \coloneqq a_{jk}$
    for all $k \in J$. Notice that, in
    analogy with \eqref{eq:SequencesDenseImage}, we have
    \begin{equation*}
        \biggl(
        \sum_{j \in J}
        \Psi(\basis{e}_j \tensor a_{j \bullet})
        \biggr)_{mn}
        =
        \sum_{j \in J}
        \delta_{j,m} \cdot a_{j n}
        =
        a_{mn}
        \qquad
        \textrm{for all }
        m,n \in J.
    \end{equation*}
    Making once again use of the silently agreed on linear ordering $\succeq$ of $J$,
    we may view the series $\sum_{j \in J} \basis{e}_j \tensor a_{j
    \bullet}$ as a net in $\field{C}^J \tensor_\pi \field{C}^J$. Notice now that,
    fixing $j \in J$, the seminorm
    \begin{equation*}
        \bigl(
            \abs{\delta_n}
            \tensor
            \abs{\delta_m}
        \bigr)
        (
            \basis{e}_j \tensor a_{j \bullet}
        )
        =
        \delta_{n,j}
        \cdot
        \abs{a_{jm}}
    \end{equation*}
    simply vanishes for $j \succ n$. Thus, taking $k \succ \ell \succ n \in J$,
    \begin{equation*}
        \bigl(
            \abs{\delta_n}
            \tensor
            \abs{\delta_m}
        \bigr)
        \biggl(
        \sum_{k \succeq j \succeq \ell}
        \basis{e}_j \tensor a_{j \bullet}
        \biggr)
        \le
        \sum_{k \succeq j \succeq \ell}
        \delta_{n,j}
        \cdot
        \abs{a_{jm}}
        =
        0,
    \end{equation*}
    i.e. $\sum_{k \succeq j \succeq \ell} \basis{e}_j \tensor a_{j \bullet}$ is
    Cauchy in
    $\field{C}^J \otimes_\pi \field{C}^J$ and thus defines an element of $\field{C}^J
    \widehat{\otimes}_\pi \field{C}^J$.

    We extend $\Psi$ to $\field{C}^J \widehat{\otimes}_\pi \field{C}^J$ by continuity
    and
    denote the resulting map again by $\Psi$. It remains to establish continuity of the
    inverse mapping. Using the concrete description of open subsets of the product
    topology, one may prove  openness of $\Psi$ directly. This, however, disregards most
    of the linear structure, making the argument quite cumbersome. Having constructed
    the inverse, we may instead show continuity by direct estimation. Indeed, this boils
    down to the familiar
    \begin{equation*}
        \bigl(
            \abs{\delta_n}
            \tensor
            \abs{\delta_m}
        \bigr)
        \biggl(
        \sum_{j \in J}
        \basis{e}_j \tensor a_{j \bullet}
        \biggr)
        =
        \sum_{j \in J}
        \delta_{j,n}
        \cdot
        \abs{a_{j m}}
        =
        \abs{a_{nm}}
        =
        \abs[\big]{\delta_{n,m}}(a)
    \end{equation*}
    for all $a \in \field{C}^{J \times J}$ and $n,m \in J$.

    Having established these preliminaries, guessing the behaviour of our map
    $\iota$ is straightforward: Let $a \in \field{C}^{J \times J}$. Suppressing
    $\Psi$, we claim that
    \begin{equation}
        \label{eq:SequenceSpaceQuadratic}
        \iota(a)
        \at[\Big]{\gamma}
        =
        \sum_{j,k \in J}
        \gamma_j
        \cdot
        a_{jk}
        \cdot
        \gamma_k
        \qquad
        \textrm{for all }
        \gamma \in \field{C}^{(J)}
    \end{equation}
    is just the obvious quadratic form. Indeed, using once again our various constructions yields
    \begin{equation*}
        \iota(a)
        \at[\Big]{\gamma}
        =
        \sum_{j \in J}
        \iota
        (\basis{e}_j \tensor a_{j \bullet})
        \at[\Big]{\gamma}
        =
        \sum_{j \in J}
        \sum_{\ell \in J}
        \gamma_j \delta_{j,\ell}
        \sum_{k \in J}
        \gamma_k
        a_{j k}
        =
        \sum_{j,k \in J}
        \gamma_j
        a_{jk}
        \gamma_k.
    \end{equation*}
    Notably, the bilinear mapping $L_a(\gamma, \delta) \coloneqq \sum_{j,k \in J}
    \gamma_j a_{jk} \delta_k$, which induces the quadratic polynomial $\iota(a)$,
    remains well defined if only one argument is a finite sequence. By induction, our
    discussion extends to arbitrary finite tensor powers.

    Let now $a = \sum_{n=0}^\infty a^{(n)} \in
    \widehat{\Sym}_0^\bullet(\field{C}^J)$ with
    $a^{(n)} \in
    \widehat{\Sym}^n_\pi(\field{C}^J)$. The symmetry condition simply means that the
    $\Sym_n$-action $(\sigma \acts a^{(n)})_{j_1 \cdots j_n} = a^{(n)}_{j_{\sigma(1)}
    \cdots j_{\sigma(n)}}$ is trivial for all $n \in \field{N}$. By what we have shown,
    \begin{equation}
        \label{eq:SequenceSpaceIotaSeries}
        \iota(a)
        \at[\Big]{\gamma}
        =
        \sum_{n=0}^\infty
        \sum_{j_1, \ldots, j_n \in J}
        a_{j_1 \cdots j_n}^{(n)}
        \cdot
        \gamma_{j_1} \cdots \gamma_{j_n}
        \qquad
        \textrm{for all }
        \gamma \in \field{C}^{(J)},
    \end{equation}
    which already very much looks like a Taylor series. In this language, the vertical
    summability condition in $\Sym_0^\bullet(\field{C}^J)$ takes the form
    \begin{equation}
        \label{eq:SequenceSpaceSummability}
        \abs{\delta_j}_{0,c}(a)
        =
        \sum_{n=0}^\infty
        c^n
        \abs{\delta_j}^n
        \bigl(a^{(n)}\bigr)
        =
        \sum_{n=0}^\infty
        c^n
        \abs[\big]{a^{(n)}_{j \cdots j}}
        <
        \infty
    \end{equation}
    for all $c \ge 0$ and $j \in J$. That is to say, the functions $z \mapsto
    \sum_{n=0}^\infty z^n a^{(n)}_{j \cdots j}$ of one complex variable have to be entire for
    every $j \in J$. As the absolute values of the Dirac functionals are a defining system
    of seminorms for $\field{C}^J$, this indeed suffices. As there are continuous
    seminorms beyond those, there is actually more structure. Fixing finitely many
    indices~$j_1, \ldots, j_n \in J$, we consider $\gamma \coloneqq z_1
    \basis{e}_{j_1} + \cdots + z_n \basis{e}_{j_n}$ with $z_1, \ldots, z_n \in
    \field{C}$. Then~\eqref{eq:SequenceSpaceIotaSeries} gives the convergent
    series
    \begin{equation}
        \label{eq:SequenceSpaceGateaux}
        \iota(a)
        \at[\Big]{\gamma}
        =
        \sum_{m=0}^\infty
        \sum_{\ell_1, \ldots, \ell_m =1}^n
        a_{j_{\ell_1} \cdots j_{\ell_m}}^{(m)}
        z_{\ell_1} \cdots z_{\ell_m}.
    \end{equation}
    That is, restrictions of $\iota(a)$ to finite dimensional subspaces are entire functions
    in the usual sense. That is exactly the notion of Gâteaux holomorphy.
    By Theorem~\ref{thm:TensorsAsPolynomials},
    \ref{item:IotaContinuous} we know that $\iota(a)$ is bounded on bounded sets $B
    \subseteq \field{C}^{(J)}$. We claim that any such set is contained in some finite
    dimensional subspace and bounded there. The boundedness is clear by continuity of
    $\abs{\delta_j}$ for every $j \in J$. We prove the contraposition of the first statement.
    To this end, let $S \subseteq \field{C}^{(J)}$ generate an infinite dimensional
    subspace. Then we find a countably infinite subset $\tilde{J} \subseteq J$ with
    corresponding linearly independent $\gamma_j \in S$ for every $j \in \tilde{J}$. We
    complete this linearly independent set to some linear algebraic basis of
    $\field{C}^{(J)}$ and write $\tilde{\delta}_j$ for the corresponding
    coefficient functionals.
    Finally, choose some enumeration of $\tilde{J} = \{j_n \;|\; n \in \field{N}\}$ and
    consider the formal sum of seminorms
    \begin{equation*}
        \seminorm{p}(\gamma)
        \coloneqq
        \sum_{n=1}^\infty
        j_n
        \abs[\big]{\tilde{\delta}_{j_n}(\gamma_{j_n})}
        \qquad
        \textrm{for all }
        \gamma \in \field{C}^{(J)}.
    \end{equation*}
    By construction, $\seminorm{p}(\gamma_{j_n}) = n$, i.e. $\seminorm{p}$ is
    unbounded on $S$. Hence $S$ is not bounded. Consequently, having
    \eqref{eq:SequenceSpaceGateaux} recovers the boundedness of $\iota(a)$ on
    bounded sets at once, as holomorphic functions on finite dimensional spaces are
    always continuous. Thus the example is consistent with our claims in
    Theorem~\ref{thm:TensorsAsPolynomials}.
\end{example}

\section{Hilbert Tensor Product}
\label{sec:HilberTensorProduct}
\epigraph{Time is a drug. \\ Too much of it kills you.}{\emph{Small Gods} -- Terry
Pratchett}
% !TeX root = ../Dissertation.tex

We conclude the chapter by developing a Hilbert\footnote{We all know David
Hilbert (1862-1943) and no finite footnote could ever do him justice.} space version of the
considerations from Section~\ref{sec:RTopologiesPolynomial}, replacing projective and
injective tensor products with the Hilbert space tensor product, as it is
discussed in \cite[Ch.~11]{kadison.ringrose:1997b}. In this setting, many of
the technical problems we have faced so far turn out to be simply absent.

Recall that given Hilbert spaces $(H_1, \langle \argument, \argument \rangle_1)$ and
$(H_2, \langle \argument, \argument \rangle_2)$ we may endow their tensor product
$H_1 \tensor H_2$ with the structure of a Pre-Hilbert space by
sesquilinearly\footnote{We follow the physicist's convention of constituting
anti-linearity of \gls{Pairing} in the \emph{first} argument.} extending
\begin{equation}
    \index{Hilbert space tensor product}
    \index{Tensor product!Hilbert}
    \label{eq:HilbertSP}
    \langle v_1 \tensor v_2, w_1 \tensor w_2 \rangle
    \coloneqq
    \langle v_1, w_1 \rangle_1
    \cdot
    \langle v_2, w_2 \rangle_2,
\end{equation}
where $v_1, w_1 \in H_1$ and $v_2, w_2 \in H_2$. The resulting space is not complete,
whenever both~$H_1$ and $H_2$ are infinite dimensional. In the sequel, we write
$H_1
\widehat{\tensor} H_2$ for the completion of the algebraic tensor product $H_1 \tensor
H_2$ with respect to the norm \gls{Norm} induced by \eqref{eq:HilbertSP}.
\begin{lemma}
    \index{Hilbert-Schmidt operator}
    Let $H_1$ and $H_2$ be Hilbert spaces. The linear extension of
    \begin{equation}
        \Xi
        \colon
        H_1 \tensor H_2
        \longrightarrow
        \overline{\HS}(H_1, H_2),
        \quad
        \Xi(v_1 \tensor v_2)v
        \coloneqq
        \langle
            v, v_1
        \rangle_1
        \cdot
        v_2
    \end{equation}
    is an isometry with dense image, where $\overline{\gls{HilbertSchmidt}}(H_1, H_2)$
    denotes the space of antilinear Hilbert-Schmidt operators from $H_1$ to $H_2$. In
    particular,
    \begin{equation}
        H_1
        \gls{TensorCompletion}
        H_2
        \cong
        \overline{\HS}(H_1, H_2) .
    \end{equation}
\end{lemma}
\begin{proof}
    Let $v_1 \in H_1$ and $v_2 \in H_2$. The operator $\Xi(v_1 \tensor v_2)$ is
    linear with finite rank, so $\Xi$ is well defined. Let $\{e_j\}_{j \in J_1}
    \subseteq H_1$ be
    a Hilbert basis. Invoking Parseval's identity and~\eqref{eq:HilbertSP}, we get
    \begin{equation*}
        \norm[\big]
        {\Xi(v_1 \tensor v_2)}_{\HS}^2
        =
        \sum_{j \in J_1}
        \norm[\big]
        {\Xi(v_1 \tensor v_2)e_j}_2^2
        =
        \sum_{j \in J_1}
        \abs[\big]
        {\langle e_j, v_1 \rangle_1}^2
        \cdot
        \norm{v_2}^2
        =
        \norm{v_1}^2_1 \cdot \norm{v_2}_2^2
        =
        \norm{v_1 \tensor v_2}^2
    \end{equation*}
    for any $v_1 \in H_1$ and $v_2 \in H_2$.
    Thus, our mapping $\Xi$ is indeed isometric. Moreover, given a Hilbert basis
    $\{f_j\}_{j \in
    J_2} \subseteq H_2$, we claim that the operators $\Xi(e_j \tensor f_k)$ with $j \in
    J_1$ and~$k \in J_2$ constitute a Hilbert basis of $\overline{\HS}(H_1,H_2)$: Indeed,
    if $T
    \in
    \overline{\HS}(H_1, H_2)$ fulfils
    \begin{equation}
        0
        =
        \big\langle
            \Xi(e_j \tensor f_k), T
        \big\rangle_{\HS}
        =
        \sum_{\ell \in J_1}
        \big\langle
        \Xi(e_j \tensor f_k) e_\ell,
        T e_\ell
        \big\rangle_2
        =
        \sum_{\ell \in J_1}
        \langle
            e_\ell, e_j
        \rangle_1
        \cdot
        \big\langle
        f_k,
        T e_\ell
        \big\rangle_2
        =
        \big\langle
        f_k,
        T e_j
        \big\rangle_2
    \end{equation}
    for all $j \in J_1$ and $k \in J_2$, then by completeness of $\{f_j\}_{j \in J_2}$ it
    follows
    that $T e_j = 0$ for all~$j \in J_1$. That is, $T = 0$. Hence, $\Xi$ does indeed
    have dense image. Finally, we note
    \begin{equation}
        \big\langle
            e_j \tensor f_k,
            e_{\ell} \tensor f_m
        \big\rangle
        =
        \langle
            e_j, e_\ell
        \rangle
        \cdot
        \langle
            f_k, f_m
        \rangle
        =
        \delta_{j,\ell}
        \cdot
        \delta_{k, m}
        \qquad
        \textrm{for $j,\ell \in J_1$ and $k,m \in J_2$},
    \end{equation}
    proving also the orthonormality by virtue of our first computation.
\end{proof}

As a by-product, we have proved the following Hilbert basis description of the Hilbert
tensor product, which in particular demystifies the structure of the completion.
\begin{corollary}
    \index{Hilbert-Schmidt operator!Orthonormal basis}
    Let $\{e_j\}_{j \in J_1} \subseteq H_1$ and $\{f_j\}_{j \in J_2} \subseteq H_2$ be
    Hilbert bases of Hilbert spaces~$H_1$ and $H_2$. Then the set
    \begin{equation}
        \big\{
            e_j \tensor f_k
            \colon
            j \in J_1,
            k \in J_2
        \big\}
    \end{equation}
    constitutes a Hilbert basis of $H_1 \widehat{\tensor} H_2$.
\end{corollary}

Reducing the Hilbert spaces $H_1$ and $H_2$ to the index sets $J_1$ and $J_2$ of their
Hilbert bases, their completed Hilbert tensor product $H_1 \widehat{\tensor} H_2$ thus
corresponds to the Cartesian product~$J_1 \times J_2$ of the index sets. In other words,
\begin{equation}
    \ell^2(J_1)
    \; \widehat{\tensor} \;
    \ell^2(J_2)
    \cong
    \ell^2(J_1 \times J_2)
\end{equation}
as Hilbert spaces. In particular, tensorizing with yet another Hilbert space $H_3$ yields
\begin{equation*}
    (H_1 \; \widehat{\tensor} \; H_2) \; \widehat{\tensor} \; H_3
    \cong
    H_1 \; \widehat{\tensor} \; (H_2 \; \widehat{\tensor} \; H_3)
\end{equation*}
canonically via the usual identification. Here, continuity is trivially fulfilled, as both ways
of bracketing yield the same norm by \eqref{eq:HilbertSP}. In the sequel, we are
interested in $n$-th tensor powers of a single Hilbert space $H$, that is the case $H_1 =
\ldots = H_n = H$. By what we have shown and the uniqueness of vector space
completions, any order of tensorizing results, up to unique isometric
isomorphism, in the same Hilbert space
\begin{equation}
    \gls{TensorPower}(H)
    \coloneqq
    \underbrace{H \widehat{\tensor} \cdots \widehat{\tensor} H}_{n \; \text{times}}.
\end{equation}
By \eqref{eq:HilbertSP}, its norm $\norm{\argument}$ on factorizing tensors is given by
\begin{equation}
    \label{eq:HilbertNormFactorizing}
    \norm{v_1 \tensor \cdots \tensor v_n}
    =
    \norm{v_1} \cdots \norm{v_n}
    \qquad
    \textrm{for }
    v_1, \ldots, v_n \in H.
\end{equation}
Moreover, we set $\Tensor^0(H) \coloneqq \field{C}$ endowed with the Euclidean scalar
product. Unwrapping the definition of a Hilbert basis yields the following coordinate
based description of the completion $\widehat{\Tensor}^n(H)$ of $\Tensor^n(H)$ by
means of square summable sequences \gls{SquareSummable}.
\begin{corollary}
    \label{cor:HilbertCoordinates}%
    \index{Hilbert-Schmidt operator!Coordinates}
    Let $H$ be a Hilbert space with Hilbert basis $\{e_j \colon j \in J\}$ and $v \in
    \widehat{\Tensor}^n(H)$. Then there exist unique coefficients $(v^{j_1 \cdots j_n}) \in
    \ell^2(J^n)$ such that
    \begin{equation}
        \label{eq:HilbertBasisDecomposition}
        v
        =
        \sum_{j_1, \ldots, j_n \in J}
        v^{j_1 \cdots j_n}
        \cdot
        e_{j_1} \tensor \cdots \tensor e_{j_n}.
    \end{equation}
    The series converges absolutely, at most countably many $\alpha^{j_1 \cdots j_n}$
    are nonzero, and they are explicitly given by
    \begin{equation}
        v^{j_1 \cdots j_n}
        =
        \langle
        e_{j_1} \tensor \cdots \tensor e_{j_n},
        v
        \rangle.
    \end{equation}
    Conversely, every sequence $\alpha \in \ell^2(J^n)$ defines an element $v \in
    \widehat{\Tensor}^n(H)$ by \eqref{eq:HilbertBasisDecomposition} and
    \begin{equation}
        \index{Seminorms!Hilbert tensor powers}
        \label{eq:HilbertNorms}
        \norm{v}_n^2
        =
        \sum_{j_1, \ldots, j_n \in J}
        \abs[\big]{v^{j_1 \cdots j_n}}^2.
    \end{equation}
\end{corollary}

As for locally convex spaces, an element $v \in \Tensor^n(H)$ is called symmetric if
it is invariant under the the symmetrization operator
\begin{equation}
    \label{eq:SymmetrizerHilbert}
    \index{Symmetrizer}
    \gls{Symmetrizer}
    \colon
    \Tensor^n(H) \longrightarrow \Tensor^n(H), \quad
    \Symmetrizer(v_1 \tensor \cdots \tensor v_n)
    \coloneqq
    \frac{1}{n!}
    \sum_{\sigma \in S_n}
    v_{\sigma(1)}
    \tensor
    \cdots
    \tensor
    v_{\sigma(n)},
\end{equation}
which defines subspaces $\Sym^n(H) \coloneqq \Symmetrizer(\Tensor^n(H)) \subseteq
\Tensor^n(H)$ for all $n \in \N$ and the corresponding graded vector space
\begin{equation}
    \Sym^\bullet(H)
    \coloneqq
    \bigoplus_{n =0}^\infty
    \Sym^n(H),
\end{equation}
where we once again set $\Sym^0(H) \coloneqq \C$. The symmetrizer turns out to be
compatible with the Hilbert space structure.
\begin{lemma}
    \label{lem:Symmetrizer}%
    \index{Symmetrizer!Hilbert}
    Let $H$ be a Hilbert space and $n \in \N$. Then the symmetrizer extends to an
    orthogonal projection
    \begin{equation}
        \Symmetrizer
        \colon
        \widehat{\Tensor}^n(H)
        \longrightarrow
        \widehat{\Sym}^n(H)
        \subseteq
        \widehat{\Tensor}^n(H).
    \end{equation}
    In particular, we may endow $\widehat{\Sym}^n(H)$ with the structure of a Hilbert
    space with scalar product given by the restriction of \eqref{eq:HilbertSP}.
\end{lemma}
\begin{proof}
    To facilitate explicit computations, we begin by studying
    \eqref{eq:SymmetrizerHilbert}. The principal difficulty here is that $\Tensor^n(H)$ is
    only pre-Hilbert. The algebraic insight that the symmetrizer $\Symmetrizer$ is a linear
    projection is standard and shall not be repeated here. Note
    \begin{align}
        \big\langle
        v_1 \tensor \cdots \tensor v_n,
        \Symmetrizer
        (w_1 \tensor \cdots \tensor w_n)
        \big\rangle
        &=
        \frac{1}{n!}
        \sum_{\sigma \in S_n}
        \langle v_1, w_{\sigma(1)} \rangle
        \cdots
        \langle v_n, w_{\sigma(n)} \rangle \\
        &=
        \frac{1}{n!}
        \sum_{\sigma \in S_n}
        \langle v_{\sigma(1)}, w_1 \rangle
        \cdots
        \langle v_{\sigma(1)}, w_n \rangle \\
        &=
        \big\langle
        \Symmetrizer
        (v_1 \tensor \cdots \tensor v_n),
        w_1 \tensor \cdots \tensor w_n
        \big\rangle
    \end{align}
    for $v_1, \ldots, v_n, w_1, \ldots, w_n \in H$. Hence, $P^* = P$ is a symmetric,
    densely defined operator. It remains to assert its continuity, which then means that its
    extension to $\widehat{\Tensor}^n(H)$ is an orthogonal projection. The
    Cauchy-Schwarz inequality facilitates the continuity estimate
    \begin{align}
        \norm[\big]
        {\Sym(v_1 \tensor \cdots \tensor v_n)}^2
        &=
        \frac{1}{n!^2}
        \sum_{\sigma, \tau \in S_n}
        \big\langle
            v_{\sigma(1)}, v_{\tau(1)}
        \big\rangle
        \cdots
        \big\langle
        v_{\sigma(n)}, v_{\tau(n)}
        \big\rangle \\
        &=
        \frac{1}{n!}
        \sum_{\sigma, \tau \in S_n}
        \big\langle
        v_{\sigma(1)}, v_1
        \big\rangle
        \cdots
        \big\langle
        v_{\sigma(n)}, v_n
        \big\rangle \\
        &\le
        \frac{1}{n!}
        \sum_{\sigma, \tau \in S_n}
        \abs[\big]
        {
            \big\langle
            v_{\sigma(1)}, v_1
            \big\rangle
        }
        \cdots
        \abs[\big]
        {
            \big\langle
            v_{\sigma(n)}, v_n
            \big\rangle
       } \\
       &\le
       \frac{1}{n!}
       \sum_{\sigma, \tau \in S_n}
       \norm[\big]{v_{\sigma(1)}}
       \cdot
       \norm[\big]{v_1}
       \cdots
       \norm[\big]{v_{\sigma(n)}}
       \cdot
       \norm[\big]{v_n} \\
       &=
       \norm[\big]{v_1}^2
       \cdots
       \norm[\big]{v_n}^2
    \end{align}
    for $v_1, \ldots, v_n \in H$, completing the proof.
\end{proof}

Using the Riesz Representation Theorem, the Hilbert space analogue of our map
\eqref{eq:Iota} may now be decomposed into
\begin{equation}
    \label{eq:IotaNHilbert}
    \iota_n
    \bigl(
        v_1 \tensor \cdots \tensor v_n
    \bigr)
    \at[\Big]{w}
    =
    \langle v_1, w \rangle
    \cdots
    \langle v_n, w \rangle
    =
    \prod_{j=1}^n
    \langle v_j, w \rangle
    =
    \iota_n
    \bigl(
        v_1 \vee \cdots \vee v_n
    \bigr)
    \at[\Big]{w}
\end{equation}
for $v_1, \ldots, v_n, w \in H$, where we use the symmetric tensor product
\begin{equation}
    v_1
    \vee \cdots \vee
    v_n
    \coloneqq
    \Symmetrizer
    \bigl(
        v_1 \tensor \cdots \tensor v_n
    \bigr)
\end{equation}
as in Section~\ref{sec:RTopologiesPolynomial}. Note that $\iota_n$ is \emph{antilinear}
with respect to $v_1 \tensor \cdots \tensor v_n$ in this description. Having
established these prelimiaries, we arrive at the following Hilbertian analogue of
our Theorem~\ref{thm:TensorsAsPolynomials} relating tensors to polynomials.
\begin{theorem}
    \label{thm:TensorsAsPolynomialsHilbert}
    \index{Tensors as polynomials!Hilbert}
    Let $H$ be a Hilbert space and $n \in \field{N}_0$.
    \begin{theoremlist}
        \item \label{item:HilbertIotaNVContinuity}
        Let $v_1, \ldots, v_n \in H$. Then the mapping
        \begin{equation}
            \iota(v_1 \vee \cdots \vee v_n)
            \colon
            H \longrightarrow \field{C}
        \end{equation}
        is a continuous $n$-homogeneous polynomial of finite type.
        \item \label{item:HilbertIotaNContinuity}
        The antilinear mapping $\iota_n \colon \Tensor^n(H) \longrightarrow
        \Pol_\beta^n(H)$ is continuous and fulfils
        \begin{equation}
            \label{eq:HilbertIotaNContinuity}
            \abs[\big]{\iota_n(v)w}
            \le
            \norm{v}_n
            \cdot
            \norm{w}^n
            \qquad
            \textrm{for all }
            v \in \Tensor^n(H), \,
            w \in H .
        \end{equation}
    \end{theoremlist}
\end{theorem}
\begin{proof}
    Let $v \in \Tensor^n(H)$ and $w \in H$. By Corollary~\ref{cor:HilbertCoordinates} we
    have
    \begin{align*}
        \iota_n(v)w
        &=
        \sum_{j_1, \ldots, j_n \in J}
        \cc{v^{j_1 \ldots j_n}}
        \cdot
        \iota_n(e_{j_1} \tensor \cdots \tensor e_{j_n})w \\
        &=
        \sum_{j_1, \ldots, j_n \in J}
        \cc{v^{j_1 \ldots j_n}}
        \cdot
        \langle e_{j_1}, w \rangle
        \cdots
        \langle e_{j_n}, w \rangle \\
        &=
        \sum_{j_1, \ldots, j_n \in J}
        \cc{v^{j_1 \ldots j_n}}
        \cdot
        w^{j_1} \cdots w^{j_n} \\
        &=
        \langle
        v,
        W
        \rangle ,
    \end{align*}
    where the final scalar product is taken in $\ell^2(J^n)$ and $W^{j_1 \ldots j_n}
    \coloneqq w^{j_1} \cdots w^{j_n}$. By the Cauchy-Schwarz inequality this implies
    \begin{equation*}
        \abs[\big]{\iota_n(v)w}
        \le
        \norm{v}_n
        \cdot
        \norm{W}_n
        =
        \norm{v}_n
        \cdot
        \norm{w}^n ,
    \end{equation*}
    proving the continuity in \ref{item:HilbertIotaNVContinuity} and taking a supremum
    over $w \in B$ for bounded subsets $B \subseteq H$ yields
    \ref{item:HilbertIotaNContinuity} at once. The statement about finite type is purely
    algebraic and independent of the topology, so it follows from
    Theorem~\ref{thm:TensorsAsPolynomials}, \ref{item:IotaV}.
\end{proof}

Passing to the full tensor algebra, we introduce the seminorms
\begin{equation}
    \label{eq:HilbertSeminormsTensor}
    \index{Seminorms!Hilbert tensor algebra}
    \seminorm{p}_{R,c}
    \colon
    \Tensor^\bullet(H) \longrightarrow [0,\infty), \quad
    \seminorm{p}_{R,c}(v)
    \coloneqq
    \sum_{n=0}^{\infty}
    n!^{R}
    \cdot
    c^n
    \cdot
    \norm{v_n}_n
\end{equation}
for $R, c \ge 0$, where $v_n \in \Tensor^n(H)$ are the homogeneous components of $v$
and $\norm{\argument}_n$ is defined as in \eqref{eq:HilbertNorms}. Fixing $R \ge 0$,
we write $\Tensor_R^\bullet(H)$ for the tensor algebra $\Tensor^\bullet(H)$
endowed with the locally
convex topology generated by the seminorms $\seminorm{p}_{R,c}$ with $c \ge 0$.
While this is technically conflicting with already established notation, we shall
not need the projective or injective $R$-topologies for the remainder of this
section, and should thus cause a bounded amount of confusion.
\begin{corollary}
    Let $H$ be a Hilbert space and $R \ge 0$. Then the mapping
    \begin{equation}
        \iota
        \colon
        \Tensor^\bullet_R(H)
        \longrightarrow
        \Pol^\bullet_\beta(H)
    \end{equation}
    and its components $\iota_n \colon \Tensor^n(H) \longrightarrow \Pol^n_\beta(H)$ are
    continuous.
\end{corollary}
\begin{proof}
    It suffices to consider the case $R = 0$, as the inclusions $\Tensor_R^\bullet(H)
    \subseteq \Tensor_S^\bullet(H)$ for parameters $S \le R$ are continuous. Let
    $v \in
    \Tensor^\bullet(H)$. By \eqref{eq:HilbertIotaNContinuity}, we have
    \begin{equation*}
        \abs[\big]{\iota(v)w}
        \le
        \sum_{n=0}^{\infty}
        \abs[\big]{\iota_n(v)w_n}
        \le
        \sum_{n=0}^{\infty}
        \norm{v}_n
        \cdot
        \norm{w}^n
        =
        \seminorm{p}_{0,\norm{w}}(v)
        \qquad
        \textrm{for }
        w \in H.
    \end{equation*}
    Taking suprema over bounded subsets of $H$ yields the continuity of
    $\iota$, which in turn implies continuity of its components at once.
\end{proof}

By continuity, we extend the mappings $\iota_n$ and $\iota$ to the completions
$\widehat{\Tensor}^n(H)$ and~$\widehat{\Tensor}_R^\bullet(H)$ as~$\hat{\iota}_n$
for all
$n \in \N_0$ and~$\hat{\iota}$, respectively. By \cite[Lemma~3.6,
\textit{iii.)}]{waldmann:2014a}, we have
\begin{equation}
    \widehat{\Tensor}_R^\bullet(H)
    =
    \bigg\{
    v
    =
    \sum_{n=0}^\infty
    v_n
    \in
    \prod_{n=0}^\infty
    \widehat{\Tensor}^n(H)
    \colon
    \seminorm{p}_{R,c}(v)
    <
    \infty
    \textrm{ for all }
    c \ge 0
    \bigg\}
\end{equation}
in analogy with Proposition~\ref{prop:RTopologies}, \ref{item:RTopologyCompletion}. As
in the locally convex situation, the extended mapping $\hat{\iota}$ and its components
map into the space of Fréchet
holomorphic functions.
\begin{theorem}
    \index{Tensors as polynomials!Hilbert space}
    Let $H$ be a Hilbert space and $v \in \widehat{\Tensor}_0^\bullet(H)$ with
    homogeneous
    components~$v_n \in \widehat{\Tensor}^n(H)$ for all $n \in \N_0$.
    \begin{theoremlist}
        \item The mappings $\hat{\iota}_n(v_n)$ and $\hat{\iota}(v)$ are Fréchet
        holomorphic and
        bounded.
        \item The $n$-th Taylor polynomial of $\hat{\iota}(v)$ is given by
        $\hat{\iota}_n(v_n)$.
    \end{theoremlist}
\end{theorem}
\begin{proof}
     As for projective tensor products, the infinite series $\sum_{n=0}^{\infty} v_n$
     actually converges in the $\Tensor_0$-topology, essentially by the definitions of the
     seminorms \eqref{eq:HilbertSeminormsTensor}. This proves the Gâteaux holomorphy
     of $\hat{\iota}(v)$, the second statement and thus also the Gâteaux holomorphy of
     $\hat{\iota}(v)$ by Corollary~\ref{cor:GateauxPowerSeries} and
     Theorem~\ref{thm:CompletenessHomogeneousPolynomials}. As each Taylor
     polynomial $\hat{\iota}_n(v_n)$ is continuous and thus bounded by
     Theorem~\ref{thm:TensorsAsPolynomialsHilbert} and
     Lemma~\ref{lem:PolynomialContinuityVsBoundedness}, this in
     turn implies Fréchet holomorphy and boundedness of $\iota(v)$ by
     Proposition~\ref{prop:FrechetBaire}.
\end{proof}

Finally, one may wonder whether an analogue of Theorem~\ref{thm:IotaEmbedding} and
Theorem~\ref{thm:IotaEmbedding2} holds for symmetric Hilbert tensor products.
This turns out to be remarkably delicate. Firstly, one has to restrict to symmetric tensors
to make $\iota_n$ and $\iota$ injective. This then makes the norm more complicated, as
\begin{equation}
    \norm[\big]
    {v_1 \vee \cdots \vee v_n}_n^2
    =
    \frac{1}{n!}
    \sum_{\sigma \in S_n}
    \big\langle
        v_1, v_{\sigma(1)}
    \big\rangle
    \cdots
    \big\langle
        v_n, v_{\sigma(n)}
    \big\rangle
    \qquad
    \textrm{for all }
    v_1, \ldots, v_n \in H.
\end{equation}
While for factorizing tensors, this provides the estimate
\begin{equation}
    \norm[\big]
    {v_1 \vee \cdots \vee v_n}_n^2
    \le
    \norm{v_1}^2
    \cdots
    \norm{v_n}^2
\end{equation}
by Lemma~\ref{lem:Symmetrizer}, mixed tensors are much more complicated as a
result. To see this, we consider $v \in \Sym^n(H)$. The polarization estimate
\eqref{eq:PolarizationEstimate} yields
\begin{equation}
    \label{eq:HilbertPolarizationEstimate}
    \sup_{w \in \Ball_1(0)}
    \abs[\big]
    {\iota_n(v)w}
    \ge
    \frac{n!}{n^n}
    \sup_{w_1, \ldots, w_n \in \Ball_1(0)}
    \abs[\big]
    {
        \big\langle
        v, w_1 \tensor \cdots \tensor w_n
        \big\rangle_n
    }.
\end{equation}
One would now like to achieve $w_1 \tensor \cdots \tensor w_n = v$, as then
\begin{equation}
    \abs[\big]
    {
        \big\langle
            v, w_1 \tensor \cdots \tensor w_n
        \big\rangle_n
    }
    =
    \langle
        v,v
    \rangle_n
    =
    \norm{v}_n^2.
\end{equation}
However, unless $v = v_1 \tensor v_1 \tensor \cdots \tensor v_1$ for some $v_1 \in H$,
this is not possible\footnote{The symmetry of $v$ was essential to apply the polarization
estimate.} and thus the supremum in
\eqref{eq:HilbertPolarizationEstimate} may be strictly smaller than $\norm{v}$.
\begin{example}
    For a
    simple example, consider $H = \C^2 = \Span(\basis{e}_1, \basis{e}_2)$ with
    \begin{equation}
        v
        \coloneqq
        2 \basis{e}_1 \vee (\basis{e}_1 + \basis{e}_2)
        =
        2 \basis{e}_1 \tensor \basis{e}_1
        +
        \basis{e}_1 \tensor \basis{e}_2
        +
        \basis{e}_2 \tensor \basis{e}_1.
    \end{equation}
    Making the ansatz $w \coloneqq (\alpha \basis{e}_1 + \sqrt{1 - \alpha^2}
    \basis{e}_2) \tensor (\beta \basis{e}_1 + \sqrt{1 - \beta^2} \basis{e}_2)$
    with $0 \le \abs{\alpha}, \abs{\beta} \le 1$ it is not hard to see that
    \begin{equation}
        \sup_{w_1, w_2 \in \Ball_1(0)^\cl}
        \abs[\big]
        {
            \langle
            v, w_1 \tensor w_2
            \rangle_2
        }
        =
        1 + \sqrt{2}
        <
        \sqrt{6}
        =
        \norm{v}_2,
    \end{equation}
    where the maximum is attained for $\alpha = \beta = \frac{\sqrt{2 + \sqrt{2}}}{2}$.
\end{example}

This does, of course, not rule out the existence of a universal estimate, but it seems
unlikely. We leave this as another question for posterity.
\begin{question}
    \label{q:Embedding}
    \index{Questions!Hilbert tensor embedding}
    Let $H$ be a Hilbert space and $n \in \N \setminus \{1\}$. Is the continuous map
    \begin{equation}
        \iota_n
        \colon
        \Sym^n_0(H)
        \longrightarrow
        \Pol_\beta^n(H)
    \end{equation}
    an embedding?
\end{question}

\chapter{Entire Functions on Lie groups}
\label{ch:Lie}
\epigraph{God does not play dice with the universe; He plays an ineffable game of His
own devising, which might be compared, from the perspective of any of the other
players [i.e. everybody], to being involved in an obscure and complex variant of poker
in a pitch-dark room, with blank cards, for infinite stakes, with a Dealer who won't
tell you the rules, and who smiles all the time.}{\emph{Good Omens} -- \\Neil Gaiman
\& Terry Pratchett}

% !TeX root = ../Dissertation.tex

In this chapter, the letter \gls{LieGroup} always denotes a finite dimensional
real Lie group with unit~\gls{GroupUnit} and of
dimension~$n \in \N_0$.\footnote{Which, in retrospect, introduces a persistent need
for mindfulness regarding the naming of indices and turns out to be a surprisingly
easy way to lose access to a perfectly functional letter. The author would not
recommend this and, indeed, would not wish this fate for his enemies. Yet, here we
are.} We
follow the tradition of denoting the corresponding Lie
algebra by the matching lowercase letter in fractur font. For instance, the letter
\gls{LieAlgStandalone} corresponds to~$G$, $\liealg{h}$ to $H$ and $\liealg{k}$ to $K$.
If the situation is more complicated, we use the more precise functorial notation
\gls{LieAlg} to denote $\liealg{g}$ instead.
When speaking of Lie groups, we always refer to \emph{real}\footnote{Considering the
general theme of the text, this constitutes a major diversion.} Lie groups unless
specified differently. The Lie exponential of $G$ is denoted by \gls{Exponential} or
simply $\exp$ if there is no room for ambiguity. All group morphisms are assumed to be
smooth, unless stated otherwise, and we refer to them as Lie group morphisms.
We define the left and right multiplications -- or viewed from an
outside perspective the translations -- by a group element~$g$ as
\begin{equation}
    \index{Translations}
    \label{eq:LeftAndRightMultiplication}
    \text{\gls{LeftMultiplication}}, \text{\gls{RightMultiplication}}
    \colon
    G \longrightarrow G, \quad
    \ell_g(h)
    \coloneqq
    gh
    \quad \textrm{and} \quad
    r_g(h)
    \coloneqq
    hg.
\end{equation}
Occasionally, we will also need the conjugation
\begin{equation}
    \label{eq:Conjugation}
    \index{Conjugation}
    \gls{Conjugation}
    \coloneqq
    r_g^{-1} \circ \ell_g \colon G \longrightarrow G,
\end{equation}
the adjoint action
\begin{equation}
    \index{Adjoint action}
    \label{eq:AdjointActionGroup}
    \gls{AdjointGroup}
    \coloneqq
    \gls{Tangent}_\E
    \Conj_g
    \colon
    \liealg{g} \longrightarrow \liealg{g}
\end{equation}
of a group element $g \in G$ on $\liealg{g}$ and the -- legally distinct -- adjoint action
\begin{equation}
    \gls{AdjointAlgebra}
    \coloneqq
    T_\E
    \Ad_{\bullet}(\xi)
    \colon
    \liealg{g} \longrightarrow \liealg{g}
\end{equation}
of a Lie algebra element $\xi \in \liealg{g}$ on $\liealg{g}$. Note that varying $g$,
respectively $\xi$, indeed yields smooth group actions of $G$, respectively a Lie
algebra action of $\liealg{g}$. Given $\xi \in \liealg{g}$, we moreover denote the
corresponding left invariant vector field by
\begin{equation}
    \index{Left invariant vector fields}
    X(\xi)
    \coloneqq
    \gls{LeftInvariantVectorField}
    \coloneqq
    T_e
    \ell_g
    \xi
    \in
    \gls{Sections}(TG).
\end{equation}
Finally, we once and for all fix a basis $(\basis{e}_1, \ldots,
\basis{e}_n)$ of $\liealg{g}$, and denote the corresponding left invariant vector
fields by
\begin{equation}
    X_j \coloneqq X_{\basis{e}_j}
    \qquad
    \textrm{for }
    j=1, \ldots, n.
\end{equation}

Returning to the complex world\footnote{Which, as far as we are concerned, is really the
\emph{real world}.}, we understand a complex Lie group as a real Lie group,
whose underlying manifold is complex and with holomorphic group multiplication.
As for smoothness within
the real setting, the holomorphicity of the group inversion is then automatic by virtue of
the holomorphic incarnation of the implicit function theorem. Note that the textbook
\cite{hilgert.neeb:2012a}, whose comprehensive and thorough nature makes it our most
powerful and dependable ally throughout the chapter\footnote{Expressed in customarily
available numbers: Twenty-nine -- that is, 29 -- times.}, uses a ostensibly
different definition of complex Lie groups: for their more group theoretic purposes, a
complex Lie group $G$ is a real Lie group, whose Lie algebra $\liealg{g}$ is a complex
vector space and its adjoint representation \eqref{eq:AdjointActionGroup} acts by
\emph{complex} linear mappings
\begin{equation}
    \index{Lie!Complex}
    \index{Complex!Lie group}
    \Ad_g
    \colon
    \liealg{g}
    \longrightarrow
    \liealg{g}
    \qquad
    \textrm{for all }
    g \in G.
\end{equation}
This alternative approach indeed turns out to be equivalent to our definition by an
application of the celebrated Newlander\footnote{August Newlander
was a doctoral student of Nirenberg. The Newlander-Nirenberg theorem was his
doctoral thesis.}-Nirenberg\footnote{Louis
Nirenberg (1925-2020) was an extraordinarily productive Canadian-American
mathematician specialized in the study of partial differential equations. His results have
had major impacts in geometric and complex analysis as well as real and complex
differential geometry. Already in his doctoral thesis, he resolved a well-known open
question known as the \emph{Weyl problem}.} Theorem
\cite{newlander.nirenberg:1957a}.

The structure of this chapter is as follows. In Section~\ref{sec:Gutt}, we review a
classical construction of the standard ordered star product $\star_\std$ for the
polynomial functions on the cotangent bundle $T^*G$ of a Lie group $G$, which
culminates
in the \emph{explicit} factorization \eqref{eq:StarProductFactorization}.
Along the way, we obtain a left invariant symbol calculus for the algebra~$\Cinfty(G)$
and the Gutt star product for the polynomial functions on the dual of the Lie algebra
$\liealg{g}^*$.
Section~\ref{sec:LieTaylor} is concerned with establishing and studying the notion of
entire functions, resulting in the algebras
\begin{equation}
    \Entire_R(G)
    \subseteq
    \Comega(G),
    \qquad
    \textrm{where $R \ge 0$},
\end{equation}
which are modelled as the Lie theoretic analogues of the algebras $\Entire$ from
\eqref{eq:EntireIntro} and $\widehat{\Sym}_R^\bullet$, which we have discussed in
Section~\ref{sec:RTopologiesPolynomial}. To this end, we first review the Lie-Taylor
formula in one and then multiple variables. As most of the material was already
presented in \cite{heins.roth.waldmann:2023a}, we are going to focus on the underlying
concepts and ideas instead of replicating all details. Of course, there will be precise
references
throughout. Whenever there is an interesting variant of a proof available, we present
it. Our efforts culminate in Theorem~\ref{thm:LieGroupStarProductContinuity}, where the
continuity of $\star_\std$ is established.

Sections~\ref{sec:UniversalComplexification},
\ref{sec:UniversalComplexificationExamples} and \ref{sec:ExtensionAndRestriction} then
serve the overarching goal of characterizing the algebra~$\Entire_0(G)$ by means of
holomorphic functions on the universal complexification $G_\C$ of $G$. That is to say,
we pass from the real analytic picture painted in Section~\ref{sec:QuantizationStrict}
to the holomorphic one in Section~\ref{sec:QuantizationHolomorphic}. This endeavour
requires careful study of the geometric properties of the universal complexification,
which we conduct in Section~\ref{sec:UniversalComplexification}. Everything we discuss
seems to be well known, at least as folklore. Having established the abstract theory,
we discuss a number of examples for universal complexifications in
Section~\ref{sec:UniversalComplexificationExamples} with a focus on possible
pathologies that may occur in the context of holomorphic extension. With these
preparations, we then prove Theorem~\ref{thm:Extension} and
Theorem~\ref{thm:Restriction}, which together provide a one-to-one correspondence
between entire functions on $G$ and holomorphic functions on the universal
complexification as locally convex algebras in the setting of linear Lie groups.
Afterwards, we establish a complexified version of the standard
ordered star product in Theorem~\ref{thm:StarProductHolomorphicContinuity} and prove
its continuity both by real techniques and by direct estimation within the holomorphic
world. Our key auxiliary result throughout is provided by the Lie theoretic
incarnation of the Cauchy estimates
\eqref{eq:CauchyEstimates}.

In Section~\ref{sec:EntireVectors}, we recast these results in the language of entire and
strongly entire vectors of the translation representation of $G$ on its algebra of
continuous functions. This shift in perspective allows for another proof of
Theorem~\ref{thm:Extension} by means of methods from
Chapter~\ref{ch:InfiniteDimensions} and a group theoretic incarnation of the monodromy
theorem. To put these considerations into their natural and conceptual
context, we first give an overview over the most important properties of continuous
and smooth vectors for representations on locally convex spaces in
Section~\ref{sec:ContinuousAndSmoothVectors}. Here, our aims beyond collecting the
hallmarks of the abstract theory are two-fold.

On the one hand, we develop practical criteria to determine the space of smooth
vectors in concrete situations. We showcase the resulting methods in
Example~\ref{ex:Schroedinger}, where we discuss the Schrödinger
representation of the Heisenberg group, which is at the heart of basic quantum
mechanics. On the other hand, we establish
differentiation of the representation of $G$ to a Lie algebra representation of
$\liealg{g}$. This then naturally leads to the converse problem of integration. That
is, the reconstruction of the group action from the Lie algebraic one.

Looking for a resolution, in turn, motivates the notions of analytic, entire and
strongly entire vectors, whose characterization is the content of
Section~\ref{sec:EntireVectors}. Our considerations culminate in
Corollary~\ref{cor:StronglyEntireIntegration}, which gives sufficient conditions for
the resolvability of the integration problem, and in a
generalization of \cite[Thm.~4.23]{heins.roth.waldmann:2023a} in the form of
Theorem~\ref{thm:Universality}, which establishes the universality of the
translation representation.

\section{Standard Ordered Quantization}
\label{sec:Gutt}
\epigraph{The presence of those seeking the truth is infinitely to be preferred to the
presence of those who think they’ve found it.}{\emph{Monstrous Regiment} -- \\ Terry
Pratchett}
% !TeX root = ../Dissertation.tex

In this section, we discuss a brief and concrete construction of the standard ordered
quantization map and the Gutt star product on the cotangent bundle $\gls{Cotangent}G$,
which goes
back to \cite{gutt:1983a}. This recovers explicit formulas already appearing in
\cite{gutt:1983a} and \cite{bordemann.neumaier.waldmann:1998a}, but adapted to our
notation and in a manner convenient for estimation. Note that the approach followed in
\cite[Sec.~2 \& Appendix~A]{heins.roth.waldmann:2023a} is different from what follows
and closer to \cite{bordemann.neumaier.waldmann:1998a} than to \cite{gutt:1983a}. Both
methods do however share that they first establish a left invariant symbol
calculus, which is utilized as a quantization map in the sense of
\eqref{eq:QuantizationMap} afterwards.

Recall that the natural manner of differentiation on a Lie group $G$ consists in left
invariant derivatives implemented by taking Lie derivatives in direction of left
invariant vector fields. Geometrically, this quantifies the rate of change of a
function along the flow of a left invariant vector field, that is to say along the image of the
Lie exponential. One may encode this as the Lie algebra morphism
\begin{equation}
    \index{Lie!Derivative}
    \label{eq:LieDerivativeOnLieAlg}
    \gls{LieDer}
    \colon
    \liealg{g} \longrightarrow
    \Diffop^1
    \bigl(
        \Cinfty(G,\R)
    \bigr) \cong \Sec(TG), \quad
    \xi
    \mapsto
    \Lie(\xi)
    =
    \Lie_{X_\xi},
\end{equation}
where we realize tangent vector fields $\gls{Sections}(TG)$ as linear derivations of
the real-valued smooth functions $\Cinfty(G,\R)$ on $G$. Before continuing, we make
the notion of linear differential operator precise.
\begin{remark}[Differential operators]
    \index{Differential operator}
    \index{Grothendieck, Alexander}
    \label{rem:DiffOps}
    Let $\algebra{A}$ be an associative and commutative algebra over some field
    $\field{k}$. We use Grothendieck's \cite{grothendieck:1967a} recursive
    definition of differential operators as it is discussed in
    \cite[Ch.~15]{mcconell.robson:2001a}. That
    is, $\gls{DiffOps}(\algebra{A})$ is the filtered algebra\footnote{A filtration of
    a vector space $V^\bullet$ indexed by $\N_0$ is a family of subspaces $\{V^m\}_{m
    \in \N_0}$ such that
    \begin{equation}
        V^\bullet
        =
        \bigcup_{m \in \N_0} V^m
        \qquad \textrm{and} \qquad
        V^m \subseteq V^{m+1}
        \quad
        \textrm{for all }
        m \in \N_0.
    \end{equation}
    If $V^\bullet$ moreover carries a multiplication such that $V^m \cdot
    V^n \subseteq V^{n+m}$ for all $n,m \in \N_0$, then $V^\bullet$ is called filtered
    algebra.}
    \begin{align}
        \index{Filtered Algebra}
        \index{Graded Algebra}
        \index{Algebra!Filtered}
        \index{Algebra!Graded}
        &\Diffop^0(\algebra{A})
        \coloneqq
        \bigl\{
        M_a
        \colon
        a \in \algebra{A}
        \bigr\}, \\
        &\Diffop^k(\algebra{A})
        \coloneqq
        \bigl\{
        D \in \Linear(\algebra{A})
        \;\big|\;
        \forall_{a \in \algebra{A}}
        \colon
        [D, M_a] \in \Diffop^{k-1}(\algebra{A})
        \bigr\}
        \qquad
        \textrm{for }
        k > 0,
    \end{align}
    where $\Linear(\algebra{A})$ is the space of $\field{k}$-linear maps from
    $\algebra{A}$ to $\algebra{A}$, the \gls{MultiplicationOperator} are multiplication
    operators with $a \in \algebra{A}$ and
    \begin{equation}
        \label{eq:Commutator}
        \index{Commutator}
        \gls{LieBracket}
        \colon
        \Linear(\algebra{A}) \times \Linear(\algebra{A})
        \longrightarrow
        \Linear(\algebra{A}), \quad
        [D,D']
        \coloneqq
        D \circ D'-D'\circ D
    \end{equation}
    is the commutator induced from the associative algebra structure of
    $\Linear(\algebra{A})$. If $\algebra{A}$ is unital, then
    \begin{equation}
        \Diffop^0(\algebra{A})
        \cong
        \End_{\algebra{A}}(\algebra{A})
        \cong
        \algebra{A},
    \end{equation}
    where $\End_{\algebra{A}}(\algebra{A})$ denotes the set of $\algebra{A}$-linear
    endomorphisms of $\algebra{A}$. Moreover, note that by definition,
    \begin{equation}
        \Diffop^k(\algebra{A})
        \subseteq
        \Diffop^{k+1}(\algebra{A})
        \qquad
        \textrm{for all }
        k \in \N_0.
    \end{equation}
    Hence, we do not get a grading\footnote{A grading of
        a vector space $V^\bullet$ indexed by $\N_0$ is a family of subspaces
        $\{V^m\}_{m \in \N_0}$ such that
        \begin{equation}
            \label{eq:Grading}
            V
            =
            \bigoplus_{m \in \N_0} V^m.
        \end{equation}
        If $V^\bullet$ moreover carries a
        multiplication such that $V^m \cdot V^n \subseteq V^{n+m}$ for all $n,m \in
        \N_0$, then $V^\bullet$ is called graded algebra. Setting $W^m \coloneqq
        \bigcup_{n=1}^m V^n$ for all $m \in \N_0$ moreover endows $V^\bullet$ with the
        structure of a filtered algebra, meaning that every graded algebra possesses
        a canonical filtration.}, but only a
        filtration. Nevertheless, we refer to $k$ as the \emph{order} of the
        differential operator. In the sequel, we will need the cases~$\algebra{A} =
        \Cinfty(G,\R)$ with $\field{k} = \R$ and $\algebra{A} = \Cinfty(G)$ with
        $\field{k} = \C$.
\end{remark}

The left invariance of $X_\xi$ translates via \eqref{eq:LieDerivativeOnLieAlg} to the
invariance of $\gls{LieDer}(\xi)$ under the natural action of $G$ on
$\Diffop^\bullet(\Cinfty(G,\R))$ given by
\begin{equation}
    \index{Differential operator!Invariance}
    \label{eq:DiffOpGAction}
    g
    \acts
    D
    \coloneqq
    \ell_{g^{-1}}^*
    \circ
    D
    \circ
    \ell_g^*
    \qquad
    g \in G, \,
    D
    \in
    \Diffop^\bullet
    \bigl(
        \Cinfty(G,\R)
    \bigr).
\end{equation}
Note that, for $\phi \in \Cinfty(G,\R)$, we have
\begin{equation}
    g \acts M_\phi
    =
    M_{\psi}
    \qquad
    \textrm{with}
    \quad
    \psi \coloneqq \ell_{g^{-1}}^* \phi.
\end{equation}
Consequently, $g \acts D$ is a differential operator of the same order as $D$. As the
action \eqref{eq:DiffOpGAction} is moreover compatible with composition, it also
respects the induced Lie bracket given by the commutator. This yields the filtered
subalgebra
\begin{equation}
    \label{eq:DiffOpInvariant}
    \Diffop^\bullet
    \bigl(
        \Cinfty(G,\R)
    \bigr)^G
    \coloneqq
    \bigl\{
        D \in \Diffop^\bullet(\Cinfty(G,\R))
        \colon
        g \acts D = D
        \textrm{ for all }
        g \in G
    \bigr\}
\end{equation}
of $G$-invariant differential operators, which is moreover stable under taking commutators.

The inverse of \eqref{eq:LieDerivativeOnLieAlg} on the $G$-invariant first order
differential operators
\begin{equation}
    \Diffop^1
    \bigl(
        \Cinfty(G,\R)
    \bigr)^G
\end{equation}
is constructed by choosing a linear algebraic basis of $\liealg{g}$, which in turn
induces a global trivialization of the tangent bundle $TG$. In particular, we then have
\begin{equation}
    \label{eq:LeftTrivialization}
    \Sec(TG)
    \cong
    \Cinfty(G,\R)
    \tensor
    \liealg{g}
    \quad \textrm{and} \quad
    \Lie(\liealg{g})
    \cong
    \Sec(TG)^G
    \cong
    \{1\}
    \tensor
    \liealg{g}
    \cong
    \liealg{g}.
\end{equation}
Notice that we may complexify all of the above structures and $\C$-linearly extend
\eqref{eq:LieDerivativeOnLieAlg}, which yields an isomorphism
\begin{equation}
    \label{eq:LieDerivativeOnLieAlgComplexified}
    \Lie
    \colon
    \liealg{g}_\C
    \longrightarrow
    \Diffop^1
    \bigl(\Cinfty(G)\bigr)
    \cong
    \Sec(\gls{TangentComplexified} G)
\end{equation}
of complex Lie algebras. As we are ultimately interested in the algebra $\Cinfty(G)$ of
complex-valued smooth functions, we shall complexify all objects in the sequel without
further comment. The fine properties of $\liealg{g}_\C \coloneqq \liealg{g} \tensor_\R
\C$ will only become important in Section~\ref{sec:UniversalComplexification}, so we
postpone its systematic discussion until then.

For now, we only need that any $\R$-linear map $\phi \colon \liealg{g} \longrightarrow
V$ into another complex vector space $V$ has a unique $\C$-linear extension $\phi_\C
\colon \liealg{g}_\C \longrightarrow V$. Doing this in every fiber of a vector bundle
then naturally leads to complexifications of vector bundles. Indeed, this universal
property characterizes $\liealg{g}_\C$ up to unique $\C$-linear isomorphism and its
natural generalization to Lie groups is precisely Hochschild's\footnote{Gerhard
Hochschild
(1915-2010) was a American mathematician, who worked in Lie theory and homological
algebra. He introduced the Hochschild cohomology, which formalizes and governs the
deformation theory of associative algebras.} universal complexification, which we
shall encounter in Definition~\ref{def:UniversalComplexification}.
\begin{remark}[Lie brackets]
    \index{Lie!Bracket}
    \index{Complex!Lie bracket}
    \index{Complexification!Lie bracket}
    \label{rem:LieBrackets}
    Note that our considerations did not require the Lie bracket of the Lie algebra
    $\liealg{g}$ whatsoever. Going back to \eqref{eq:LieDerivativeOnLieAlg}, one may
    thus reverse the logic and \emph{define} the Lie bracket on $\liealg{g} = T_\E G$
    by declaring the isomorphism $\Lie$ to be bracket preserving. In fact, this is how
    the Lie algebra structure on $\liealg{g}$ is typically defined in the first place,
    see e.g. \cite[Def.~9.1.7]{hilgert.neeb:2012a}.

    This is also a good opportunity to review the construction of the \emph{complex}
    Lie algebra of a complex Lie group $G$ as it can be found in
    \cite[Ch.~X~6.]{kobayashi.nomizu:1969a}. Indeed, using the complex structure, we
    may define the holomorphic tangent spaces $T^\C_g G^{(1,0)}$ as the set of
    $\C$-linear derivations of the space of complex valued holomorphic germs
    $\Holomorphic_g(G)$ at $g \in G$. Their disjoint union
    \begin{equation}
        \index{Holomorphic!Tangent bundle}
        T^\C G^{(1,0)}
        \coloneqq
        \gls{DisjointUnion}
        T^\C_g G^{(1,0)}
    \end{equation}
    then constitutes the holomorphic tangent bundle of $G$, which is involutive. Its
    complex structure is obtained from the one of $G$ by applying the
    tangent functor.

    By our definition of a complex Lie group, the left translations $\ell_g$ from
    \eqref{eq:LeftAndRightMultiplication} are holomorphic. In particular, they are
    almost holomorphic, which means that
    \begin{equation}
        T_\E \ell_g
        \circ
        J_\E
        =
        J_g
        \circ
        T_\E \ell_g
        \qquad
        \textrm{for all }
        g \in G.
    \end{equation}
    Here, $J \in \Sec[\infty](\End(G))$ denotes the almost complex structure induced
    by the complex structure of $G$. This implies that the left invariant vector
    fields fulfil
    \begin{equation}
        X_\xi
        =
        T_\E \ell_\bullet \xi
        \in
        \Sec[\infty]
        \bigl(
             T^\C G^{(1,0)}
        \bigr)
        \qquad
        \textrm{for all }
        \xi \in T_\E^\C G^{(1,0)}.
    \end{equation}
    The involutivity of $T^\C G^{(1,0)}$ now means that, building on what is known
    from the real situation, the $G$-invariant vector fields form a Lie
    subalgebra of $\Sec[\infty](T^\C G^{(1,0)})$, which is linearly isomorphic to the
    \emph{complex Lie algebra}
    \begin{equation}
        \label{eq:LieAlgebraComplex}
        \gls{LieAlgStandaloneComplex}
        \coloneqq
        T_\E^\C G^{(1,0)}
    \end{equation}
    by means of left translation. Finally, we endow $\hat{\liealg{g}}$ with a Lie
    bracket by demanding that this linear isomorphism is even an isomorphism of complex
    Lie algebras.
%    It is now a straightforward exercise to prove that
%    \begin{equation}
%        \hat{\liealg{g}}
%        \cong
%        \liealg{g}
%    \end{equation}
%    as \emph{real} Lie algebras. Here, $\liealg{g}$ is the usual Lie algebra of $G$,
%    when viewed as a real Lie group.
\end{remark}

\index{Universal!Enveloping algebra}
Returning to the setting of a real Lie group $G$ for now, we essentially follow ideas
of \cite{berezin:1967a}. Starting out, we model the universal
enveloping algebra~$\Universal^\bullet(\liealg{g_\C})$ as the complexification of
the tensor
algebra~$\Tensor^\bullet(\liealg{g})$ over $\liealg{g}$ modulo the relation
\begin{equation}
    \label{eq:UniversalEnvelopeRelation}
    \xi \tensor \chi - \chi \tensor \xi
    \sim
    [\xi,\chi]
    \qquad
    \textrm{for all }
    \xi,\chi \in \liealg{g}.
\end{equation}
As $\Diffop^\bullet(\Cinfty(G))$ is a unital associative algebra, the universal
property of the universal enveloping algebra $\Universal^\bullet(\liealg{g_\C})
\cong
\gls{UniversalEnvelope}_\C$ implies that the Lie algebra morphism
\eqref{eq:LieDerivativeOnLieAlgComplexified} has a unique extension to a unital
algebra morphism $\Lie \colon \Universal^\bullet(\liealg{g}_\C) \longrightarrow
\Diffop^\bullet(\Cinfty(G))$.\footnote{For a systematic discussion of the universal
enveloping
algebra, we refer to the textbook \cite[Sec.~7.1]{hilgert.neeb:2012a}.} In fact,
\begin{equation}
    \label{eq:LieDerivativeOnEnveloping}
    \Lie(\xi_1 \tensor \cdots \tensor \xi_k)
    =
    \Lie_{X_{\xi_1}}
    \circ
    \cdots
    \circ
    \Lie_{X_{\xi_k}}
\end{equation}
holds for all $\xi_1, \ldots, \xi_k \in \liealg{g}_\C$. Hence, $\Lie$ even provides a morphism of filtered algebras
\begin{equation}
    \label{eq:EnvelopingAsDiffops}
    \Lie
    \colon
    \Universal^\bullet(\liealg{g}_\C)
    \longrightarrow
    \Diffop^\bullet\bigl(\Cinfty(G)\bigr)^G.
\end{equation}
Here, we define the filtration of $\Universal^\bullet(\liealg{g}_\C)$ by declaring
$\Universal^0(\liealg{g}_\C) \coloneqq \C$, $\Universal^1(\liealg{g}_\C) \coloneqq
\liealg{g}_\C$ and
\begin{equation}
    \Universal^{k+\ell}(\liealg{g}_\C)
    \coloneqq
    \Span
    \bigl(
        \Universal^k(\liealg{g}_\C) \tensor \Universal^\ell(\liealg{g}_\C)
    \bigr)
    \qquad
    \textrm{for all }
    k,\ell \in \N_0.
\end{equation}
As $\liealg{g}_\C$ generates $\Universal^\bullet(\liealg{g}_\C)$, this
indeed suffices to define a filtration. Alternatively, one may inherit the canonical
filtration induced by the grading of the tensor algebra via
\eqref{eq:UniversalEnvelopeRelation}, as the relation identifies degree two elements
with degree one elements. Note that this does not quite respect the grading, but we
nevertheless get a filtration. In the end, both constructions yield the same filtered
algebra $\Universal^\bullet(\liealg{g}_\C)$.

The mapping \eqref{eq:EnvelopingAsDiffops} is the higher order analogue of \eqref{eq:LeftTrivialization}. We will use both notations occurring in \eqref{eq:LieDerivativeOnEnveloping} interchangeably in the sequel, depending on which leads to more comprehensible formulas.

The idea is now that we after applying \eqref{eq:EnvelopingAsDiffops}, we may take the
leading symbol $\sigma$ as discussed in \cite[Sec.~1.1]{schapira:1985a} and then
evaluate at the group unit. This provides symbol maps
\begin{equation}
    \index{Symbol calculus}
    \label{eq:SymbolMap}
    \tilde{\sigma}_k
    \colon
    \Universal^k(\liealg{g}_\C)
    \longrightarrow
    \Sym^k(\liealg{g}_\C)
\end{equation}
with values in the symmetric algebra $\Sym^\bullet(\liealg{g}_\C)$ over
$\liealg{g}_\C$ for all $k \in \N_0$. By construction,
\begin{align}
    \tilde{\sigma}_k
    (\xi_1 \tensor \cdots \tensor \xi_k)
    &=
    \sigma_k
    \bigl(
        \Lie_{X_{\xi_1}}
        \circ
        \cdots
        \circ
        \Lie_{X_{\xi_k}}
    \bigr) \\
    &=
    \gls{SymbolLeading}(\Lie_{X_{\xi_1}})
    \vee
    \cdots
    \vee
    \sigma_1(\Lie_{X_{\xi_k}}) \\
    &=
    \xi_1 \vee \cdots \vee \xi_k
    \label{eq:SymbolMapOnGenerators}
\end{align}
for all $\xi_1, \ldots, \xi_k \in \liealg{g}$ by virtue of
\eqref{eq:LieDerivativeOnEnveloping} and the fact that leading symbols are always
symmetric. Hence, \eqref{eq:SymbolMap} is compatible with products. Taking another
look at \eqref{eq:SymbolMapOnGenerators}, we see that $\sigma_k$ is moreover
surjective for all $k \in \N_0$. Somewhat surprisingly, the formal sum
\begin{equation}
    \tilde{\sigma}
    \coloneqq
    \bigoplus_{k=0}^{\infty}
    \tilde{\sigma}_k
    \colon
    \Universal^\bullet(\liealg{g}_\C)
    \longrightarrow
    \Sym^\bullet(\liealg{g}_\C)
\end{equation}
turns out to be a bijection, which may be viewed as a particular instance of the
celebrated Poincaré-Birkhoff-Witt Theorem\footnote{Of which Capelli proved a
special case, a fact which Poincaré was unaware of. In turn, neither Birkhoff
nor Witt refer to Poincaré's work, which was only widely acknowledged after
appearing within Bourbaki's book \cite{bourbaki:1975a}.}
\cite{capelli:1890a, poincare:1900a,
birkhoff:1937a, witt:1937a}, a modern Lie algebraic incarnation of which can be
found in the textbook \cite[Thm.~7.1.9]{hilgert.neeb:2012a}.
\index{Poincaré-Birkhoff-Witt Theorem}

We provide a proof for this for our particular situation within
Proposition~\ref{prop:SymbolCalculusLeftInvariant}, which facilitates the construction
of the standard ordered quantization map as $\Lie \circ \tilde{\sigma}^{-1}$. This
leads us to the following educated guess, where we additionally incorporate a
combinatorial prefactor and a suitable power of the parameter $\hbar \in \C$, see
again the discussion in Section~\ref{sec:QuantizationStrict}.
\begin{definition}[Standard ordered quantization map]
    \index{Standard ordered!Quantization map}
    Let $G$ be a Lie group with corresponding Lie algebra $\liealg{g}$ and~$\hbar \in
    \C$. The standard ordered quantization map is
    \begin{align}
        \label{eq:StdOrdering}
        \begin{split}
        &\gls{StandardOrdering}
        \colon \Cinfty(G)
        \tensor
        \Sym^\bullet(\liealg{g}_\C)
        \longrightarrow
        \Diffop^\bullet
        \bigl(\Cinfty(G)\bigr), \\
        &\varrho_\std
        \bigl(
            \phi \tensor \xi_1 \vee \cdots \vee \xi_k
        \bigr)
        \coloneqq
        \biggl(
            \frac{\hbar}{\I}
        \biggr)^k
        \frac{1}{k!}
        \cdot
        M_\phi
        \circ
        \sum_{\sigma \in S_k}
        \Lie
        \bigl(
            \sigma
            \acts
            \xi_1 \tensor \cdots \tensor \xi_k
        \bigr),
    \end{split}
    \end{align}
    where $M_\phi \colon \Cinfty(G) \longrightarrow \Cinfty(G)$ is the multiplication
    operator with symbol $\phi \in \Cinfty(G)$.
\end{definition}

Going back to \eqref{eq:LeftTrivialization}, we note the isomorphisms
\begin{equation}
    \label{eq:PolynomialsCotangent}
    \index{Polynomials!Cotangent bundle $T^*G$}
    \Cinfty(G)
    \tensor
    \Sym^\bullet(\liealg{g}_\C)
    \cong
    \Sec
    \bigl(
        \Sym^\bullet_\C(TG)
    \bigr)
    \cong
    \gls{PolynomialsVectorBundle}(T^*G),
\end{equation}
where $\Pol(T^*G)$ is the set of smooth maps $P \colon T^*G \longrightarrow \C$
that are polynomial in fiber directions. Indeed, restriction to the zero section
$\gls{CanonicalEmbedding} \colon G \longrightarrow T^*G$ yields $\iota^* P \in
\Cinfty(G)$ and, provided $\iota^* P(g) \neq 0$, the map
\begin{equation}
    P
    \bigl(
        g, T_\E \ell_g \argument
    \bigr)
    \colon
    \liealg{g}
    \cong T_\E G
    \longrightarrow
    \C
\end{equation}
is in $\Sym^\bullet(\liealg{g})$. If there is no such $g$, we may infer $P = 0$. The
inverse process is given by the pullback with the bundle projection $\pi \colon T^*G
\longrightarrow G$ and left translation. Consequently, the standard ordering
\eqref{eq:StdOrdering} provides a symbol calculus for the polynomials on the cotangent
bundle $T^*G$.
\begin{proposition}[Standard ordered symbol calculus,
{\cite[Prop.~2.3]{heins.roth.waldmann:2023a}}]
    \index{Standard ordered!Symbol calculus}
    \label{prop:SymbolCalculusLeftInvariant}
    Let $G$ be a Lie group and $\hbar \in {\C \setminus \{0\}}$. Then the
    standard ordered quantization map \eqref{eq:StdOrdering} constitutes an
    isomorphism of filtered vector spaces, i.e. $\varrho_\std$ is a linear bijection
    such that
    \begin{equation}
        \varrho_\std
        \bigl(
            \Cinfty(G)
            \tensor
            \Sym^m(\liealg{g}_\C)
        \bigr)
        =
        \Diffop^m
        \bigl(
            \Cinfty(G)
        \bigr)
        \qquad
        \textrm{for all }
        m \in \N_0.
    \end{equation}
    Here, we endow $\Cinfty \tensor \Sym^\bullet(\liealg{g}_\C)$ with the filtration
    corresponding to the grading of $\Sym^\bullet(\liealg{g}_\C)$ and by declaring
    $\Cinfty(G) \tensor 1$ to be degree zero. Moreover, the same is true for its
    restriction
    \begin{equation}
        \label{eq:StandardOrderingRestriction}
        \varrho_\std
        \colon
        \Sym^\bullet(\liealg{g}_\C)
        \cong
        \Sym^\bullet_\C(\liealg{g})
        \longrightarrow
        \Diffop^\bullet
        \bigl(
            \Cinfty(G)
        \bigr)^G.
    \end{equation}
\end{proposition}
\begin{proof}
     We begin with the proof of the second statement. By \eqref{eq:SymbolMapOnGenerators}, we have
     \begin{equation}
         \label{eq:SymbolOfStd}
         \bigl(
            \sigma_k
            \circ
            \varrho_\std
         \bigr)
         (\xi_1 \vee \cdots \vee \xi_k)
         =
         \biggl(
            \frac{\hbar}{\I}
         \biggr)^k
         \xi_1 \vee \cdots \vee \xi_k
         \tag{$\sim$}
     \end{equation}
     for all $\xi_1, \ldots, \xi_k \in \liealg{g}$ and thus
     \eqref{eq:StandardOrderingRestriction} is injective,
     matching the surjectivity of the symbol maps~$\sigma_k$. Taking another look
     at \eqref{eq:StdOrdering}, we also see that the operator
     \begin{equation}
         \varrho_\std
         \bigl(
            \xi_1
            \vee \cdots \vee
            \xi_k
         \bigr)
         =
         \biggl(
            \frac{\hbar}{\I}
         \biggr)^k
         \frac{1}{k!}
         \sum_{\sigma \in S_k}
         \Lie(\xi_{\sigma(1)})
         \circ \cdots \circ
         \Lie(\xi_{\sigma(k)})
     \end{equation}
     is $G$-invariant, as all of its factors $\Lie(\xi_1), \ldots,
     \Lie(\xi_k)$ are. To establish its surjectivity, let
     \begin{equation}
         D
         \in
         \Diffop^k
         \bigl(
            \Cinfty(G)
         \bigr)^G
         \qquad \textrm{with} \quad
         \tilde{\sigma}_k(D)
         \in
         \Sym^k(\liealg{g}_\C)
     \end{equation}
     be given. By the left invariance of $D$, also its leading symbol
     $\sigma_k(D)$ is left invariant, which means that
     \begin{equation}
         \sigma(D)
         \at[\Big]{g}
         =
         (T_\E \ell_g)^{\tensor k}
         \tilde{\sigma}_k(D)
         \qquad
         \textrm{for all }
         g \in G.
     \end{equation}
    Consequently,
    \begin{equation}
        D'
        \coloneqq
        D
        -
        \biggl(
            \frac{\I}{\hbar}
        \biggr)^k
        \varrho_\std
        \bigl(
            \sigma(D)
        \bigr)
    \end{equation}
    is another $G$-invariant differential operator with
    \begin{equation}
        \sigma_k(D')
        =
        \sigma_k(D)
        -
        \biggl(
            \frac{\I}{\hbar}
        \biggr)^k
        \Bigl(
            \sigma_k \circ \varrho_\std
        \Bigr)
        \bigl(
            \sigma(D)
        \bigr)
        =
        0
    \end{equation}
    by virtue of \eqref{eq:SymbolOfStd}, where we also use the assumption $\hbar
    \neq 0$.
    That is to say, $D'$ has strictly lower degree than $D$. Iterating
    this procedure proves the surjectivity of
    \eqref{eq:StandardOrderingRestriction}. Taking another look at
    \eqref{eq:StdOrdering}, it is left $\Cinfty(G)$-linear. Consequently, it
    remains to assert that
    \begin{equation}
        \label{eq:DiffOpFactorization}
        \Diffop^\bullet
        \bigl(
            \Cinfty(G)
        \bigr)
        \cong
        \Cinfty(G) \tensor
        \Diffop^\bullet
        \bigl(
            \Cinfty(G)
        \bigr)^G.
    \end{equation}
    By what we have already shown, we know that $1 \in \C$ and the operators of the
    form $\Lie(\xi)$ generate the algebra $\Diffop^\bullet(\Cinfty(G))^G$. Thus it
    suffices to check the claim in filtration degree one. But this is just the
    bijectivity of \eqref{eq:LieDerivativeOnLieAlgComplexified}. We have completed the
    proof.
\end{proof}

\index{Star product!Standard ordered}
\index{Star product!Cotangent bundle $T^*G$}
Having established the bijectivity of the quantization map $\varrho_\std$, we may
define the standard ordered star product by pulling back the product $\circ$ of
$\Diffop^\bullet(\Cinfty(G))$, which is the composition of differential
operators. More precisely, we define a linear map
\begin{equation}
    \gls{StarProductStandardOrdered}
    \colon
    \Pol(T^*G) \tensor \Pol(T^*G)
    \longrightarrow
    \Pol(T^*G)
\end{equation}
by setting
\begin{equation}
    \label{eq:StarProductStd}
    P \star_\std Q
    \coloneqq
    \varrho_\std^{-1}
    \bigl(
        \varrho_\std(P)
        \circ
        \varrho_\std(Q)
    \bigr).
\end{equation}
By construction, the associativity of $\star_\std$ is clear, as the composition of
functions is always associative. Taking another look at \eqref{eq:StdOrdering}, we
moreover see that there are no convergence issues. As the degrees of the polynomials
$P$ and $Q$ are finite, the same is true for the degree of the operator
$\varrho_\std(P) \circ \varrho_\std(Q)$ and thus $P \star_\std Q$ is a polynomial of
degree at most $\deg P + \deg Q$. Of course, we ultimately want to go beyond the
purely polynomial situation regardless, at which point we will have to deal with
convergence questions. The following corollary is immediate from
Proposition~\ref{prop:SymbolCalculusLeftInvariant}.
\begin{corollary}[Gutt star product]
    \index{Gutt, Simone}
    \index{Star product!Gutt}
    \index{Star product!Lie algebra}
    Let $G$ be a Lie group. Then the mapping
    \begin{equation}
        \label{eq:StarProductGutt}
        \gls{StarProductGutt}
        \colon
        \Sym^\bullet(\liealg{g}_\C) \tensor \Sym^\bullet(\liealg{g}_\C)
        \longrightarrow
        \Sym^\bullet(\liealg{g}_\C), \quad
        p \star_\Gutt q
        \coloneqq
        (\mathbb{1} \tensor p)
        \star_\std
        (\mathbb{1} \tensor q)
    \end{equation}
    constitutes an associative product on $\Sym^\bullet(\liealg{g}_\C)$, where
    \gls{ConstantFunction} denotes the function on $G$ with constant value equal
    to $1 \in \C$.
\end{corollary}

We call $\star_\Gutt$ the \emph{Gutt star product} on $\liealg{g}$ or the \emph{Lie
algebra star product}. By virtue of {Lie-Cartan's~Theorem}, which we will review in
Theorem~\ref{thm:Lie3}, one may forget about the underlying Lie group~$G$ and work
with an abstract real or complex Lie algebra $\liealg{g}$ from the start. This is
also clear from the construction we have presented.
\begin{remark}[Gutt star product]
    Calling $\star_\Gutt$ the Gutt\footnote{Simone Gutt (born 1956) is a Belgian
    differential geometer and professor emeritus at Université libre de Bruxelles. She
    has made numerous contributions to symplectic geometry and symplectic formal
    deformation quantization, including early existence and classification results.}
    star product is historically somewhat inaccurate.
    The construction in \cite{gutt:1983a} goes beyond the purely Lie algebraic
    situation and
    provides a Weyl\footnote{Weyl-Wigner-Groenewold-Moyal.} ordering on $T^*G$,
    which is different from -- but equivalent as a
    star product to -- the standard ordering we have presented so far. We will come
    back to the ordering problem in Remark~\ref{rem:Ordering}. However, also the
    Weyl-ordered star product descends to the Lie algebra -- and yields a star product
    equivalent to the Gutt star product -- which warrants our choice of terminology.
    Moreover, it
    should be noted that the star product was independently constructed by
    Drinfel'd\footnote{Vladimir Drinfel'd (born 1954) is a Ukrainian mathematical physicist,
    who moved to the United States of America in 1999 and is currently working at the
    University of Chicago. His main interest consists in algebraic geometry over finite
    fields and their connection to number theory. In the year 1990, he received the Fields
    Medal for his contributions to the geometric
    Langlands correspondence and the theory of quantum groups.}
    \cite{drinfeld:1983a} in the setting of quantum groups. In the surrounding literature,
    many authors speak of $h$-deformations instead of star products.

    Unwrapping the definitions, one may think of the combinatorics of the Lie algebra star
    product of $p,q \in \Sym(\liealg{g}_\C)$ in the following manner. First, one associates
    the invariant and symmetric differential operators $\varrho_\std(p)$ and
    $\varrho_\std(q)$ to $p$ and $q$. Composition results in another invariant and
    symmetric differential operator, which is however no longer of the form in
    \eqref{eq:StdOrdering}. The star product $p \star q$ then constitutes the
    \emph{symmetric polynomial} inducing the composition $\varrho_\std(p) \circ
    \varrho_\std(q)$.

    An instructive example is the following: Let $\xi, \chi \in \liealg{g}$ and $k,\ell \in \N_0$. Then
    \begin{equation}
        \varrho_\std
        \bigl(
            \underbrace{\xi \vee \cdots \vee \xi}_{k-\textrm{times}}
        \bigr)
        =
        \biggl(
            \frac{\hbar}{\I}
        \biggr)^k
        \Lie_{X_\xi}^k
        \quad \textrm{and} \quad
        \varrho_\std
        \bigl(
        \underbrace{\chi \vee \cdots \vee \chi}_{\ell-\textrm{times}}
        \bigr)
        =
        \biggl(
        \frac{\hbar}{\I}
        \biggr)^\ell
        \Lie_{X_\chi}^\ell.
    \end{equation}
    From here, it is not at all clear what $\xi^k \star_\Gutt \chi^\ell$ might be, even for $k = \ell = 1$. The principal issue is that, while $\xi^k$ and $\chi^\ell$ are homogeneous, their product will typically longer be. In the degenerate case $\xi = \chi$ we of course get that $\xi \star_\Gutt \chi = \xi \vee \chi$.

    \index{Baker-Campbell-Hausdorff formula}
    A comprehensive discussion of the Lie algebra star product and its intimate
    connection to the Baker-Campbell-Hausdorff formula can be found in
    \cite[Sec.~3.2 \& 3.4]{stapor:2016a} and
    \cite[Sec.~2]{esposito.stapor.waldmann:2017a}. The inclined reader may
    moreover find a proper analysis of our example in terms of the
    Baker-Campbell-Hausdorff series within \cite[Lemma~4.1.9]{stapor:2016a}.

    The main result of \cite{esposito.stapor.waldmann:2017a} is the following continuity
    result, which the authors even extended to certain infinite dimensional Lie algebras
    with the asymptotic estimate property or fulfilling a nilpotency condition. We confine
    ourselves to the finite dimensional situation and use the notation from
    Section~\ref{sec:RTopologiesPolynomial} and in particular \eqref{eq:RTopology} by
    using the fact that a Lie algebra is, in particular, a vector space.
\end{remark}

\begin{theorem}[Continuity of the Gutt star product, {\cite[Sec.~3.1]{esposito.stapor.waldmann:2017a}}]
    \index{Continuity!Gutt star product}
    \label{thm:LieAlgebraStarProductContinuity}
    Let $\liealg{g}$ be a finite dimensional Lie algebra, $\hbar \in \C$ and $R \ge 1$. The Lie algebra star product
    \begin{equation}
        \star_\Gutt
        \colon
        \Sym_R^\bullet(\liealg{g}_\C)
        \otimes
        \Sym_R^\bullet(\liealg{g}_\C)
        \longrightarrow
        \Sym_R^\bullet(\liealg{g}_\C)
    \end{equation}
    is well defined and continuous. Moreover,
    \begin{equation}
        \C \ni \hbar
        \mapsto
        p \star_\Gutt q
        \in
        \Sym_R^\bullet(\liealg{g}_\C)
    \end{equation}
    is Fréchet holomorphic for all $p,q \in \Sym_R(\liealg{g}_\C)$. The condition
    $R \ge 1$ is necessary in general.
\end{theorem}

Before investigating the Poisson brackets to warrant calling $\star_\std$ and
$\star_\Gutt$ star products in the sense of Definition~\ref{def:StarProduct}, we
collect some explicit formulas. As before, $\mathbb{1} \in \Cinfty(G)$ denotes the
constant function with value $1 \in \field{C}$, i.e. the unit of $(\Cinfty(G),
\cdot)$. Recall that the standard ordered quantization map \eqref{eq:StdOrdering}
is left $\Cinfty(G)$-linear by definition. Thus,
\begin{equation}
    \label{eq:StdOrderedLeftLinearity}
    (\phi \tensor 1)
    \star_\std
    (\mathbb{1} \tensor q)
    =
    \frac{\hbar}{\I}
    \varrho_\std^{-1}
    \bigl(
        M_\phi
        \circ
        \Lie(q)
    \bigr)
    =
    \phi \tensor q
    =
    (\phi \tensor 1)
    \cdot
    (\mathbb{1} \tensor q)
\end{equation}
for $\phi \in \Cinfty(G)$ and $q \in \Sym^\bullet_\C(\liealg{g})$ by
\eqref{eq:StarProductStd}. Here, $\cdot$ acts as pointwise multiplication in the first
factor and as the symmetric tensor product $\vee$ in the second. The remaining
combination is more interesting.
\begin{theorem}[Standard ordered star product,
{\cite[Prop.~2.4]{heins.roth.waldmann:2023a}}]
    \index{Star product!Standard ordered}
    \index{Standard ordered!Star product}
    \label{Thm:StdOrderedStarProduct}
    Let $G$ be a Lie group, $\hbar \in \C \setminus \{0\}$, $\phi, \psi \in
    \Cinfty(G)$, $\xi_1, \ldots, \xi_k \in \liealg{g}_\C$ and $q \in
    \Sym^\bullet_\C(\liealg{g})$. Then
    \begin{multline}
        \label{eq:StarProductFactorization}
        (\phi \tensor \xi_1 \vee \cdots \vee \xi_k)
        \star_\std
        (\psi \tensor q)
        =
        (\phi \tensor 1)
        \cdot
        (\mathbb{1} \tensor \xi_1 \vee \cdots \vee \xi_k)
        \star_\std
        (\psi \tensor 1)
        \star_\std
        (\mathbb{1} \tensor q) \\
        =
        \sum_{p=0}^{k}
        \biggl(
            \frac{\hbar}{\I}
        \biggr)^p
        \frac{\phi}{p! (k-p)!}
        \sum_{\sigma \in S_k}
        \Lie
        \bigl(
            \xi_{\sigma(1)} \tensor \cdots \tensor \xi_{\sigma(p)}
        \bigr)
        \psi
        \tensor
        (\xi_{\sigma(p+1)} \vee \cdots \vee \xi_{\sigma(k)}) \star_\Gutt q.
    \end{multline}
\end{theorem}
\begin{proof}
    One may obtain the presented formulae from
    \cite[Prop.~11]{bordemann.neumaier.waldmann:1998a} by making
    \eqref{eq:PolynomialsCotangent} explicit. We provide another proof based on
    \eqref{eq:StdOrdering}. Let $p \coloneqq \xi_1 \vee
    \cdots \vee \xi_k$ for ease of notation. By associativity of the star product, the
    neutrality of $\mathbb{1}$ and $1$ as well as~\eqref{eq:StdOrderedLeftLinearity},
    we first note
    \begin{equation}
        (\phi \tensor p)
        \star_\std
        (\psi \tensor q)
        =
        (\phi \tensor 1)
        \cdot
        (\mathbb{1} \tensor p)
        \star_\std
        (\psi \tensor 1)
        \star_\std
        (\mathbb{1} \tensor q).
    \end{equation}
    To compute the middle product, we work with a basis $(\basis{e}_1, \ldots,
    \basis{e}_n)$ of the Lie algebra. By bilinearity of $\star_\std$, we may assume
    that
    \begin{equation}
        p
        =
        \basis{e}_{j_1}
        \vee \cdots \vee
        \basis{e}_{j_k}
        \qquad
        \textrm{for some}
        \qquad
        j_1, \ldots, j_k
        \in
        \{1,\ldots,n\}.
    \end{equation}
    The trick is to apply the isomorphism $\varrho_\std$
    on both sides of \eqref{eq:StarProductFactorization} after setting $\phi =
    \mathbb{1}$ and $q = 1$. By definition of the star product in
    \eqref{eq:StarProductStd} and
    \eqref{eq:StdOrdering}, the left-hand side yields on the one hand
    \begin{align}
       \varrho_\std
       \bigl(
            (\mathbb{1} \tensor p)
            \star_\std
            (\psi \tensor 1)
       \bigr)
       &=
       \varrho_\std
       (\mathbb{1} \tensor \basis{e}_{j_1} \vee \cdots \vee \basis{e}_{j_k})
       \circ
       \varrho_\std
       (\psi \tensor 1) \\
       &=
       \biggl(\frac{\hbar}{\I}\biggr)^{k}
       \frac{1}{k!}
       \sum_{\sigma \in S_k}
       \Lie
       \bigl(
        \basis{e}_{j_{\sigma(1)}}
        \tensor \cdots \tensor
        \basis{e}_{j_{\sigma(k)}}
       \bigr)
       \circ
       M_\psi \\
       &=
       \biggl(\frac{\hbar}{\I}\biggr)^{k}
       \frac{1}{k!}
       \sum_{\sigma \in S_k}
       \Lie\bigl(\basis{e}_{j_{\sigma(1)}}\bigr)
       \circ \cdots \circ
       \Lie\bigl(\basis{e}_{j_{\sigma(k)}}\bigr)
       \circ
       M_\psi.
    \end{align}
    To arrive back at an expression of the form \eqref{eq:StdOrdering}, we need to pull the
    multiplication operator~$M_\psi$ to the left. As
    $\Lie(\basis{e}_{j_{\sigma(1)}}), \ldots, \Lie(\basis{e}_{j_{\sigma(k)}})$ are
    derivations, we may do so at the cost of a Leibniz rule. Each time, this produces
    two terms: one where the derivative
    acts on $\psi$ and another where it does not. Due to the additional symmetrization,
    this leads to the noncommutative higher Leibniz rule
    \begin{equation}
        \index{Leibniz rule}
        \sum_{\sigma \in S_k}
        D_{\sigma(1)}
        \circ \cdots \circ
        D_{\sigma(k)}
        (ab)
        =
        \sum_{\sigma \in S_k}
        \sum_{p=0}^{k}
        \binom{k}{p}
        \bigl(
            D_{\sigma(1)}
            \circ \cdots \circ
            D_{\sigma(p)}
            a
        \bigr)
        \bigl(
            D_{\sigma(p+1)}
            \circ \cdots \circ
            D_{\sigma(k)}
            b
        \bigr)
    \end{equation}
    for any derivations $D_1, \ldots, D_n$ of some -- not necessarily associative --
    ring
    $\ring{R}$ with $\field{Q} \subseteq \ring{R}$ and $a,b \in \ring{R}$, see
    \cite[Appendix.~B]{heins.roth.waldmann:2023a}. In our situation, this yields
    \begin{align}
        &\varrho_\std
        \bigl(
        (\mathbb{1} \tensor p)
        \star_\std
        (\psi \tensor 1)
        \bigr) \\
        &=
        \biggl(\frac{\hbar}{\I}\biggr)^{k}
        \frac{1}{k!}
        \sum_{\sigma \in S_k}
        \sum_{p=0}^{k}
        \binom{k}{p}
        M_{
            \Lie(\basis{e}_{j_{\sigma(1)}}
            \tensor \cdots \tensor
            \basis{e}_{j_{\sigma(p)}})
            \psi
        }
        \circ
        \Lie\bigl(\basis{e}_{j_{\sigma(p+1)}}\bigr)
        \circ \cdots \circ
        \Lie\bigl(\basis{e}_{j_{\sigma(k)}}\bigr) \\
        &=
        \sum_{p=0}^{k}
        \biggl(\frac{\hbar}{\I}\biggr)^{k}
        \frac{1}{p! (k-p)!}
        \sum_{\sigma \in S_k}
        M_{
            \Lie(\basis{e}_{j_{\sigma(1)}}
            \tensor \cdots \tensor
            \basis{e}_{j_{\sigma(p)}})
            \psi
        }
        \circ
        \Lie\bigl(\basis{e}_{j_{\sigma(p+1)}}\bigr)
        \circ \cdots \circ
        \Lie\bigl(\basis{e}_{j_{\sigma(k)}}\bigr).
    \end{align}
    Using the left $\Cinfty(G)$-linearity of $\varrho_\std$ gives on the other hand
    \begin{align}
        \varrho_\std
        &\biggl(
            \sum_{p=0}^k
            \biggl(\frac{\hbar}{\I}\biggr)^{p}
            \frac{1}{p! (k-p)!}
            \sum_{\sigma \in S_k}
            \Lie
            \bigl(
                \basis{e}_{j_{\sigma(1)}}
                \tensor \cdots \tensor
                \basis{e}_{j_{\sigma(p)}}
            \bigr)
            \psi
            \tensor
            \basis{e}_{j_{\sigma(p+1)}} \vee \cdots \vee \basis{e}_{j_{\sigma(k)}}
        \biggr) \\
        &=
        \sum_{p=0}^k
        \biggl(\frac{\hbar}{\I}\biggr)^{p}
        \frac{1}{p! (k-p)!}
        \sum_{\sigma \in S_k}
        M_{
            \Lie
            (   \basis{e}_{j_{\sigma(1)}}
                \tensor \cdots \tensor
                \basis{e}_{j_{\sigma(p)}}
            )
            \psi
        }
        \circ
        \varrho_\std
        \bigl(
            \mathbb{1}
            \tensor
            \basis{e}_{j_{\sigma(p+1)}} \vee \cdots \vee \basis{e}_{j_{\sigma(k)}}
        \bigr).
    \end{align}
    Dropping the matching prefactors and using the auxiliary variables
    \begin{equation}
        \chi^{\sigma}_{j}
        \coloneqq
        \basis{e}_{j_{\sigma(p+j)}}
        \qquad
        \textrm{for }
        j=1,\ldots,k-p
        \quad \textrm{and} \quad
        \sigma \in S_k,
    \end{equation}
    we finally compute
    \begin{align}
        \sum_{\sigma \in S_k}
        \varrho_\std
        &\bigl(
            \mathbb{1}
            \tensor
            \basis{e}_{j_{\sigma(p+1)}}
            \vee \cdots \vee
            \basis{e}_{j_{\sigma(k)}}
        \bigr) \\
        &=
        \sum_{\sigma \in S_k}
        \varrho_\std
        \bigl(
            \mathbb{1}
            \tensor
            \chi_1^{\sigma}
            \vee \cdots \vee
            \chi_{k-p}^{\sigma}
        \bigr) \\
        &=
        \sum_{\sigma \in S_k}
        \biggl(
            \frac{\hbar}{\I}
        \biggr)^{k-p}
        \frac{1}{(k-p)!}
        \sum_{\tau \in S_{k-p}}
        \Lie
        \bigl(
            \chi_{\tau(1)}^{\sigma}
            \tensor \cdots \vee
            \chi_{\tau(k-p)}^{\sigma}
        \bigr) \\
        &=
        \sum_{\sigma \in S_k}
        \biggl(
        \frac{\hbar}{\I}
        \biggr)^{k-p}
        \frac{1}{(k-p)!}
        \sum_{\tau \in S_{k-p}}
        \Lie
        \bigl(
            \basis{e}_{j_{\sigma(p+\tau(1))}}
            \tensor \cdots \tensor
            \basis{e}_{j_{\sigma(p+\tau(k-p))}}
        \bigr) \\
        &=
        \frac{1}{(k-p)!}
        \sum_{\tau \in S_{k-p}}
        \sum_{\sigma \in S_k}
        \biggl(
        \frac{\hbar}{\I}
        \biggr)^{k-p}
        \Lie
        \bigl(
        \basis{e}_{j_{\sigma(p+\tau(1))}}
        \tensor \cdots \tensor
        \basis{e}_{j_{\sigma(p+\tau(k-p))}}
        \bigr) \\
        &=
        \frac{1}{(k-p)!}
        \sum_{\tau \in S_{k-p}}
        \sum_{\sigma' \in S_k}
        \biggl(
        \frac{\hbar}{\I}
        \biggr)^{k-p}
        \Lie
        \bigl(
        \basis{e}_{j_{\sigma'(p+1)}}
        \tensor \cdots \tensor
        \basis{e}_{j_{\sigma'(k))}}
        \bigr) \\
        &=
        \sum_{\sigma \in S_k}
        \biggl(\frac{\hbar}{\I}\biggr)^{k-p}
        \Lie
        \bigl(
            \basis{e}_{j_{\sigma(p+1)}}
            \tensor \cdots \tensor
            \basis{e}_{j_{\sigma(k)}}
        \bigr),
    \end{align}
    as $\sigma' \coloneqq \sigma \circ (\id \times \tau)$ yields all of $S_k$ if we fix $\tau \in S_{k-p}$ and vary $\sigma \in S_k$. Summarizing, we have now shown \eqref{eq:StarProductFactorization} under the restraint $q = 1$. But, going back to the factorization and once more using associativity, we get
    \begin{align}
        &(\phi \tensor 1)
        \star_\std
        (\mathbb{1} \tensor p)
        \star_\std
        (\psi \tensor 1)
        \star_\std
        (\mathbb{1} \tensor q) \\
        &=
        \Bigg(
        \sum_{p=0}^{k}
        \underbrace{
        \biggl(
            \frac{\hbar}{\I}
        \biggr)^p
        \frac{\phi}{p! (k-p)!}
        \sum_{\sigma \in S_k}
        \Lie
        \bigl(
            \xi_{\sigma(1)} \tensor \cdots \tensor \xi_{\sigma(p)}
        \bigr)
        \psi}_{\eqqcolon \psi_p}
        \tensor
        \bigl(
            \xi_{\sigma(p+1)} \vee \cdots \vee \xi_{\sigma(k)}
        \bigr)
        \Bigg)
        \star_\std
        (\mathbb{1} \tensor q) \\
        &=
        \sum_{p=0}^{k}
        \bigl(
            \psi_p
            \tensor
            \xi_{\sigma(p+1)} \vee \cdots \vee \xi_{\sigma(k)}
        \bigr)
        \star_\std
        (\mathbb{1} \tensor q) \\
        &=
        \sum_{p=0}^{k}
        \bigl(
        \psi_p
        \tensor
        1
        \bigr)
        \star_\std
        \bigl(
            \mathbb{1}
            \tensor
            \xi_{\sigma(p+1)} \vee \cdots \vee \xi_{\sigma(k)}
        \bigr)
        \star_\std
        (\mathbb{1} \tensor q) \\
        &=
        \sum_{p=0}^{k}
        \bigl(
        \psi_p
        \tensor
        1
        \bigr)
        \star_\std
        \bigl(
            \mathbb{1}
            \tensor
            \xi_{\sigma(p+1)} \vee \cdots \vee \xi_{\sigma(k)}
            \star_\Gutt
            q
        \bigr) \\
        &=
        \sum_{p=0}^{k}
        \psi_p
        \tensor
        \bigl(
            \xi_{\sigma(p+1)} \vee \cdots \vee \xi_{\sigma(k)}
            \star_\Gutt
            q
        \bigr)
    \end{align}
    by virtue of \eqref{eq:StarProductGutt}. This completes the proof.
\end{proof}

The factorization \eqref{eq:StarProductFactorization} facilitated the continuity
estimates in \cite[Thm.~6.3]{heins.roth.waldmann:2023a}, which build on
Theorem~\ref{thm:LieAlgebraStarProductContinuity}. We will use similar techniques in
Theorem~\ref{thm:StarProductHolomorphicContinuity} to provide a holomorphic version of
this result. Note that \eqref{eq:StarProductFactorization} makes sense for $\hbar =
0$, which gives
\begin{align}
    (\phi \tensor \xi_1 \vee \cdots \vee \xi_k)
    \star_\std
    (\psi \tensor q)
    &=
    \frac{1}{k!}
    \sum_{\sigma \in S_k}
    \psi
    \tensor
    \bigl(
        \xi_{\sigma(1)}
        \vee \cdots \vee
        \xi_{\sigma(k)}
        \vee
        q
    \bigr) \\
    &=
    \psi
    \tensor
    \bigl(
    \xi_1
    \vee \cdots \vee
    \xi_k
    \vee
    q
    \bigr).
\end{align}
This means that, despite our construction breaking down for $\hbar = 0$, naively
setting $\hbar = 0$ within the explicit formula \eqref{eq:StarProductFactorization}
provides an associative product. Taking a closer look, it is nothing else than the
\emph{undeformed product} of $\Pol(T^*G)$ and thus our construction is consistent with
the classical limit \eqref{eq:ClassicalLimit}. In the sequel, this is what we mean
when speaking of the case $\hbar = 0$.
\begin{remark}[Interpretation of $\star_\std$]
    \label{rem:InterpretationStarStd}
    \index{Standard ordered!Interpretation}
    One may think of \eqref{eq:StarProductFactorization} for $\phi = \mathbb{1}$ and
    $q = 1$ as arising from the \emph{commutator} of the operators $\varrho_\std(p)$
    and $M_\psi$ up to some combinatorial factors. In particular, the natural pairing
    of $\Diffop^\bullet(G)$ and $\Cinfty(G)$ appears as the term with $p=k$ in the
    summation.

    In \cite[Thm.~6.3 \& Thm.~6.4]{heins.roth.waldmann:2023a} an extension of
    Theorem~\ref{thm:LieAlgebraStarProductContinuity} to the full standard ordered
    star product $\star_\std$ on a suitable subalgebra was proved. In view of
    Theorem~\ref{thm:LieAlgebraStarProductContinuity}, one expects operators
    associated to elements of the form
    \begin{equation}
        \label{eq:ExponentialElement}
        \sum_{k=0}^{\infty}
        \frac{\xi^{\vee k}}{k!}
        =
        \exp_\vee(\xi)
    \end{equation}
    for some fixed $\xi \in \liealg{g}$ to appear in the completion of
    $\Sym_1^\bullet(\liealg{g}_\C)$, see again Proposition~\ref{prop:RTopologies},~
    \ref{item:RTopologyCompletion}. That is,
    actions of operators such as $\exp(\Lie_{X_\xi})$ will be crucial. The first half
    of Section~\ref{sec:LieTaylor} will thus be dedicated to the study of such
    operators, as for sufficiently nice functions, the operator $\exp(\Lie_{X_\xi})$
    simply constitutes a pullback with the right translation by $\exp \xi$
    on~$G$. We
    will phrase this as the \emph{Lie-Taylor formula} on~$G$ in
    Theorem~\ref{thm:LieTaylor}.

    \index{Goodman, Roe}
    Note that \eqref{eq:ExponentialElement} is \emph{not quite} contained in
    $\widehat{\Sym}_1^\bullet(\liealg{g}_\C)$ for general Lie algebras, see
    \cite[Ex.~3.3]{esposito.stapor.waldmann:2017a} and the discussion in
    \cite[Sec.~5.4]{stapor:2016a}. Indeed, this is only possible if the
    Baker-Campbell-Hausdorff series converges on the entirety of $\liealg{g}_\C$,
    which limits the geometric complexity and generality severely.
    Nevertheless, as elements of slightly slower growth are still relevant, this
    provides good
    intuition regardless. In particular, we see that $\Cinfty(G)$ is much too large of an
    algebra if we want to be able to take star products with the entirety of
    $\Sym_1^\bullet(\liealg{g}_\C)$. This
    type of question for natural pairings was investigated by Goodman\footnote{Roe
    Goodman is a harmonic analyst and distinguished Professor Emeritus at Rutgers
    University since 2015. He is a specialist on analytic, entire and strongly entire
    vectors of strongly continuous unitary representations.} in
    \cite{goodman:1969a,goodman:1970a,goodman:1971a}, and we will return to the
    accompanying concepts in Section~\ref{sec:ContinuousAndSmoothVectors} and
    Section~\ref{sec:EntireVectors}.

    The geometric interpretation of operators associated to non-real Lie algebra
    elements is unclear at this point, but we will resolve this in
    Theorem~\ref{thm:Extension}: They provide translations within the universal
    complexification $G_\C$ of $G$.
\end{remark}

We proceed by taking a closer look at the associated Poisson brackets. That is,
we enforce the semiclassical limit \eqref{eq:SemiclassicalLimitContinuous} by
defining
\begin{equation}
    \label{eq:PoissonStd}
    \index{Standard ordered!Poisson bracket}
    \index{Poisson bracket!Standard ordered}
    \bigl\{
        P, Q
    \bigr\}
    \coloneqq
    \I
    \frac{\D}{\D \hbar}
    \bigl(
        P \star_\std Q
        -
        Q \star_\std P
    \bigr)
    \at[\bigg]{\hbar = 0}´
\end{equation}
and
\begin{equation}
    \index{Poisson bracket!Gutt}
    \label{eq:PoissonGutt}
    \bigl\{
    p,q
    \bigr\}_\Gutt
    \coloneqq
    \I
    \frac{\D}{\D \hbar}
    \bigl(
    p \star_\Gutt q
    -
    q \star_\Gutt p
    \bigr)
    \at[\bigg]{\hbar = 0}
\end{equation}
for $P,Q \in \Pol(T^*G)$ and $p,q \in \Sym^\bullet(\liealg{g}_\C)$. Note that
$\star_\Gutt$ does indeed have a non-trivial dependence on $\hbar$. Thus
Theorem~\ref{Thm:StdOrderedStarProduct} yields the following:
\begin{corollary}[Poisson bracket,
{\cite[Prop.~8.2]{bordemann.neumaier.waldmann:1998a}}]
    \index{Poisson bracket}
    Let $G$ be a Lie group with Lie algebra $\liealg{g}$ and let $\phi, \psi \in \Cinfty(G)$
    as well as $\xi_1, \ldots, \xi_k \in \liealg{g}_\C$.
    \begin{corollarylist}
        \item \label{item:StdPoissonBracket}
        If the second factor is constant in fiber direction, the Poisson bracket is given by
        \begin{equation}
            \bigl\{
                \phi \tensor \xi_1 \vee \cdots \vee \xi_k,
                \psi \tensor 1
            \bigr\}
            =
            \phi
            \cdot
            \sum_{p=1}^{k}
            \Lie(\xi_p) \psi
            \tensor
            \xi_{1} \vee \overset{\xi_p}{\cdots} \vee \xi_{k},
        \end{equation}
        where \gls{AbsentFactor} indicates that the factor $\xi_p$ is absent.
        \item \label{item:StdPoissonBracketVsKKS}
        \index{Kostant-Kirillov-Souriau bracket}
        The Poisson bracket $\{\argument,\argument\}_\Gutt$ on $\Pol(\liealg{g}^*)$
        coincides with the Kostant-Kirillov-Souriau bracket\footnote{For a
        comprehensive discussion, we refer to \cite[Ch.~1]{kirillov:2004a}.}
        \gls{KostantKirillovSouriau} on $\Pol(\liealg{g}^*) \subseteq
        \Cinfty(\liealg{g}^*)$ if one identifies $\Sym^\bullet(\liealg{g})$ with
        polynomials on the dual
        $\liealg{g}^*$ by natural pairing for linear polynomials and
        compatibility with $\vee$.\footnote{Just like we did earlier in
        Section~\ref{sec:RTopologiesPolynomial}.}
    \end{corollarylist}
\end{corollary}
\begin{proof}
    The idea for the first part is that differentiation by $\hbar$ and then plugging in $\hbar = 0$ singles out the term with $p=1$ in \eqref{eq:StarProductFactorization}. Consequently, we have
    \begin{align}
        \I
        \frac{\D}{\D \hbar}
        \bigl(
            (\phi \tensor \xi_1 \vee \cdots \vee \xi_k)
            \star_\std
            (\psi \tensor 1)
        \bigr)
        \at[\Big]{\hbar = 0}
        &=
        \frac{\phi}{(k-1)!}
        \sum_{\sigma \in S_k}
        \Lie(\xi_{\sigma(1)})
        \psi
        \tensor
        \xi_{\sigma(2)} \vee \cdots \vee \xi_{\sigma(k)} \\
        &=
        \phi
        \cdot
        \sum_{p=1}^{k}
        \Lie(\xi_p) \psi
        \tensor
        \xi_1 \vee \overset{\xi_p}{\cdots} \vee \xi_k,
    \end{align}
    where we have used the symmetry of $\vee$. The other contribution vanishes, as
    left multiplication with $\psi \tensor 1$
    reduces to the pointwise product, which is indepedent of $\hbar$. We have proved
    \ref{item:StdPoissonBracket}. It suffices to
    check \ref{item:StdPoissonBracketVsKKS} on generators, i.e. constant functions and
    linear
    monomials.\footnote{Here, we use that in order $\hbar^k$ of the star product there
    is a bi-differential operator of order $(k,k)$ acting on the inputs. For $k=1$,
    this guarantees that \eqref{eq:PoissonStd} and \eqref{eq:PoissonGutt} define
    Poisson brackets, as $(1,1)$ bi-differential operators are exactly bi-derivations
    and \eqref{eq:PoissonStd} as well as \eqref{eq:PoissonGutt} have built-in
    antisymmetry.} For the former, there is nothing to be shown, as both sides simply
    vanish. Thus let $\xi_1, \xi_2 \in \liealg{g}$. By \cite[Lemma.~4.1]{stapor:2016a}
    with $z = \hbar/\I$, we have
    \begin{equation}
        \xi_1 \star_\Gutt \xi_2
        =
        \xi_1 \vee \xi_2
        +
        \frac{\hbar}{2 \I}
        [\xi_1,\xi_2]
    \end{equation}
    and thus by antisymmetry of the Lie bracket
    \begin{equation}
        \bigl\{
            \xi_1, \xi_2
        \bigr\}_G
        =
        [\xi_1, \xi_2].
    \end{equation}
    This yields the statement, as suppressing the identification, we simply have
    $\{\xi,\eta\}_{\KKS} = [\xi,\eta]$.
\end{proof}

We conclude this section with some more remarks.
\begin{remark}[Complex Lie algebras]
    \index{Star product!Gutt complex}
    \index{Star product!Standard ordered complex}
    \index{Lie!Algebra complex}
    \label{rem:ComplexLieAlgebras}
    Let $G$ be a complex Lie group. Replacing the complexified Lie algebra
    $\liealg{g}_\C$ with
    the complex Lie algebra $\hat{\liealg{g}}$ from Remark~\ref{rem:LieBrackets}
    and~$\Cinfty(G)$ with $\Holomorphic(G)$ throughout provides a complex version of
    the
    Gutt and standard ordered star products. To make this precise, one either
    establishes a holomorphic version of the left invariant symbol calculus discussed
    in Proposition~\ref{prop:SymbolCalculusLeftInvariant} or observes that the
    explicit formulas from \cite[Sec.~2.3]{esposito.stapor.waldmann:2017a} and
    \eqref{eq:StarProductFactorization} respect holomorphy. Either way, one is lead to
    defining the \emph{the holomorphic polynomial algebra}
    \begin{equation}
        \label{eq:HolomorphicPolynomialAlgebra}
        \gls{HolomorphicPolynomials}
        \coloneqq
        \Holomorphic(G)
        \tensor_\pi
        \Sym^\bullet_1
        (\hat{\liealg{g}}),
    \end{equation}
    where $\tensor_\pi$ denotes the projective tensor product\footnote{See again the
    discussion around \eqref{eq:ProjectiveSeminorms}.} and
    $\Sym^\bullet_R(\hat{\liealg{g}})$ is the symmetric algebra endowed with topology
    induced by the seminorms \eqref{eq:RTopologySeminorms} with $R=1$. As the
    continuity result from Theorem~\ref{thm:LieAlgebraStarProductContinuity} is
    agnostic regarding the base field, the resulting complex Gutt star product is
    still continuous. We shall prove that the same is true for the resulting complex
    standard ordered star product in
    Theorem~\ref{thm:StarProductHolomorphicContinuity}.
\end{remark}

\begin{remark}[Fedosov's construction]
    \index{Fedosov, Boris Vasilievich}
    \label{rem:Fedosov}
    The principal idea of \cite[Sec.~8]{bordemann.neumaier.waldmann:1998a} is to
    identify the star product constructed in \cite{gutt:1983a} as the
    Fedosov star product \cite{fedosov:1994a} associated to the half-commutator
    connection
    \begin{equation}
        \index{Half commutator connection}
        \nabla_{X_\xi}
        X_\chi
        \coloneqq
        \frac{1}{2}
        X_{[\xi,\chi]}
    \end{equation}
    for the tangent bundle $TG$. From Fedosov's classification results it is clear that
    there exists some connection \gls{Connection} for $TG$ with this property,
    but the abstract
    machinery does not provide its exact form. In our formulation, the choice of this
    connection is virtually invisible, but may be recovered from \eqref{eq:StdOrdering} by
    comparing with the abstract formulas involving symmetrized covariant derivatives.
    We refer to \cite[Sec.~6.4]{waldmann:2007a} for a comprehensive discussion of
    Fedosov's construction and to \cite{collini:2016a} for a reformulation and
    subsequent generalization within the framework of algebraic quantum field theory.
\end{remark}

\begin{remark}[$\kappa$-orderings and Weyl star product]
    \label{rem:Ordering}
    \index{Orderings}
    \index{Star product!Weyl}
    \index{Neumaier, Nicolai}
    In \eqref{eq:StdOrdering}, we have chosen an ordering prescription. Namely, we have
    put the coefficient function to the very left as a multiplication operator and all
    polynomials to the right as derivatives. This leads to particularly simple formulas, as
    the derivatives would act on their coefficient functions otherwise. However,
    this is not
    the only choice, and has the serious flaw of being incompatible with the
    $^*$-involution on $\Pol(T^*G)$ given by pointwise complex conjugation. By
    additionally choosing a smooth positive density $\mu$ on $T^*G$ -- such as the
    density associated to the left invariant volume form $\theta^1 \wedge \cdots \wedge
    \theta^n$ -- one may pass from the standard ordering to other ordering
    prescriptions
    by means of a Neumaier operator, which provides a linear bijection of
    $\Pol^\bullet(T^*G)$
    intertwining the products. Introducing an additional real parameter $\kappa$,
    this
    yields a continuum of equivalent star products, both in the formal and the
    strict setting. In the case of $\kappa = 1/2$,
    this results in the so called Weyl\footnote{No, we are not doing this again.}
    star product, which is compatible with the
    $^*$-involution. Roughly speaking, it arises in the same fashion as
    \eqref{eq:StdOrdering}, but with full symmetrization over all possible orderings.
    Note
    that these considerations are not specific to $T^*G$ and work for any cotangent
    bundle, leading to an astonishing explicit formula for the adjoint operator of
    \emph{any} differential
    operator acting on a suitable $\Ltwo$-space. We refer the interested reader to
    \cite{bordemann.neumaier.waldmann:1998a, bordemann.neumaier.waldmann:1999a}
    for the algebraic details and \cite[Appendix~A, Sec.~2.2 \&
    Prop.~6.7]{heins.roth.waldmann:2023a}, which contains a continuity result for the
    Neumaier operators.
\end{remark}

\section{Lie-Taylor Expansions and Entire Functions}
\label{sec:LieTaylor}
\epigraph{Tiffany jumped when she saw a balloon sail up above the trees, catch the
wind, and swoop away, but it turned out to be just a balloon and not a lump of excess
Brian. She could tell this because it was followed by a long scream of rage mixed with a
roar of complaint: ``AAaargwannawannaaaagongongonaargggaaaaBLOON!'' which is
the traditional sound of a very small child learning that with balloons, as with life itself, it
is important to know when not to let go of the string. The whole point of balloons is to
teach small children this.}{\emph{A Hat Full of Sky} -- Terry Pratchett}
% !TeX root = ../Dissertation.tex

This section serves as a review of the most important properties of the algebra of
\emph{entire functions} on $G$. That is to say, we are in the realm of strict
deformation quantization of the cotangent bundle $T^*G$, as discussed in
Section~\ref{sec:QuantizationStrict}. As such, almost everything we present is already
discussed in \cite{heins.roth.waldmann:2023a} to vastly greater detail. The inclusion
of this section is thus purely for the reader's convenience and to keep the text as
self-contained as possible.

We begin with the classical Lie-Taylor formulas, as they can be found in the
beginning of \cite[Sec.~2.1.4]{helgason:2001a} for Lie groups and
\cite[(1.48)]{forstneric:2011a} or \cite[Prop.~4.5]{heins.roth.waldmann:2023a} for
arbitrary analytic manifolds. We already specialize to the setting to Lie groups, as
that
is the only situation we are ultimately interested in. In the sequel, we write
$\Comega(G)$ for
the real analytic functions on a Lie group $G$, where we use the fact that the
exponential atlas endows $G$ with the structure of an analytic manifold, see for
instance
\cite[(1.6.3)~Theorem]{duistermaat.kolk:2000a}.
%\newpage
\begin{proposition}[One variable Lie-Taylor formula]
    \index{Lie-Taylor!One variable}
    \index{Taylor!Lie one variable}
    \label{prop:LieTaylorOneVariable}
    Let $G$ be a Lie group and fix a real analytic $\phi \in \Comega(G)$, an expansion
    point $g \in G$ and a direction $\xi \in \liealg{g}$.
    \begin{propositionlist}
        \item There exists a radius $R > 0$ such that the Lie-Taylor formula
        \begin{equation}
            \label{eq:LieTaylorOneVariable}
            \phi
            \bigl(
                g \cdot \gls{ExponentialSeries}(t\xi)
            \bigr)
            =
            \bigl(
                \exp(t \Lie_{X_\xi})
                \phi
            \bigr)(g)
            =
            \sum_{k=0}^{\infty}
            \frac{t^k}{k!}
            \Lie_{X_\xi}^k
            \phi
            (g)
        \end{equation}
        holds whenever $t \in \R$ fulfils $\abs{t} < R$.
        \item The series \eqref{eq:LieTaylorOneVariable} is convergent in the Fréchet
        space $\Cinfty(-R,R)$ with respect to~$t$. That is to say, each individual
        derivative with respect to $t$ converges uniformly on compact subsets of the
        open interval $\gls{Interval} \subseteq \R$.
    \end{propositionlist}
\end{proposition}

One should think of $g \in G$ as the expansion point of the power series, $\xi \in
\liealg{g}$ the direction and $t \in \R$ as the magnitude of the perturbation or,
alternatively, a time parameter. For the additive group $G = \R^n$, the Lie-Taylor
formula \eqref{eq:LieTaylorOneVariable} takes the familiar form
\begin{equation}
    \index{Taylor!Classical}
    \phi(x + t\xi)
    =
    \sum_{k=0}^{\infty}
    \frac{t^k}{k!}
    \frac{\partial^k \phi}{\partial \xi^k}
    (x)
    \qquad
    \textrm{for all $x,\xi \in \R^n$},
\end{equation}
where we suppress the Lie exponential as usual and $\frac{\partial}{\partial \xi}$
denotes the directional derivative in direction of $\xi$.
\begin{proof}[of Proposition~\ref{prop:LieTaylorOneVariable}]
    Let $U \subseteq \liealg{g}$ be an open neighbourhood of $0$ such that
    \begin{equation}
        \exp
        \colon
        U \longrightarrow \exp(U)
    \end{equation}
    is a diffeomorphism. The idea is to consider the left-hand side of
    \eqref{eq:LieTaylorOneVariable} for $t\xi \in U$ as our function $\phi$
    precomposed with the exponential chart $\ell_g \circ \exp$ around $g$. Or in
    other words, we may view it as a coordinate expression for $\phi$. By
    analyticity of $\phi$ and by
    shrinking $U$ if necessary, we may thus assume that $\phi \circ \ell_g \circ
    \exp
    \colon U \longrightarrow \C$ is globally given by its Taylor series around zero.
    Moreover, we may assume that $\xi = \basis{e}_1$ by choosing an appropriate basis
    of $\liealg{g}$. As the origin is an inner point of $U$, we find some $R > 0$ such
    that the open~\gls{SupnormVector}-ball $\Ball_{R}(0)$ is contained within $U$. For
    $\norm{t\xi} = \abs{t} < R$, the classical Taylor formula around the origin then
    takes the simple form
    \begin{equation}
        \label{eq:LieTaylorProof}
        \phi
        \bigl(
            g \exp(t \xi)
        \bigr)
        =
        \sum_{k=0}^{\infty}
        \frac{t^k}{k!}
        \frac{\D^k}{\D s^k}
        \phi
        \bigl(
            g \exp(s \xi)
        \bigr)
        \at[\Big]{s=0}.
        \tag{$\ast$}
    \end{equation}
    By left invariance, we have
    \begin{align}
        \Bigl(
            \Lie(\xi)
            \phi
        \Bigr)
        \bigl(
            g \exp(s\xi)
        \bigr)
        &=
        \bigl(
            \ell_{g \exp(s\xi)}^*
            \circ
            \Lie(\xi)
        \bigr)
        \phi
        \at[\Big]
        {\E} \\
        &=
        \bigl(
            \Lie(\xi)
            \circ
            \ell_{g \exp(s\xi)}^*
        \bigr)
        \phi
        \at[\Big]
        {\E} \\
        &=
        \frac{\D}{\D t}
        \phi
        \bigl(
            g \exp(s\xi) \exp(t\xi)
        \bigr)
        \at[\Big]{t=0} \\
        &=
        \frac{\D}{\D t}
        \phi
        \bigl(
        g \exp( (s + t) \xi)
        \bigr)
        \at[\Big]{t=0} \\
        &=
        \frac{\D}{\D s}
        \phi
        \bigl(
            g \exp(s\xi)
        \bigr)
    \end{align}
    and thus
    \begin{equation}
        \Lie_{X_\xi}^k
        \phi
        \at[\Big]
        {g \exp(s\xi)}
        =
        \frac{\D^k}{\D s^k}
        \phi
        \bigl(
            g \exp(s\xi)
        \bigr)
        \qquad
        \textrm{for all $k \in \N_0$.}
    \end{equation}
    Plugging this into \eqref{eq:LieTaylorProof} and once more invoking left invariance, we obtain
    \begin{equation}
        \phi
        \bigl(
            g \exp(t \xi)
        \bigr)
        =
        \sum_{k=0}^{\infty}
        \frac{t^k}{k!}
        \Lie_{X_\xi}^k
        \phi
        \bigl(
            g \exp(s\xi)
        \bigr)
        \at[\Big]{s=0}
        =
        \sum_{k=0}^{\infty}
        \frac{t^k}{k!}
        \Lie_{X_\xi}^k
        \phi
        (g)
    \end{equation}
    whenever $\abs{t} < R$, which is \eqref{eq:LieTaylorOneVariable}.
    The remaining statement follows from the fact that the partial sums of
    \eqref{eq:LieTaylorOneVariable} correspond exactly to the partial sums
    of the classical Taylor series of $\phi \circ \ell_g \circ \exp$ and the latter is
    even convergent in the Fréchet space $\Cinfty(U)$.
\end{proof}

We proceed with some remarks on the construction.
\begin{remark}[A wrong conclusion]
    \label{rem:ExponentialOperators}
    \index{Exponential operators}
    Notably, Proposition~\ref{prop:LieTaylorOneVariable} does not imply that the
    exponential operator $\exp(\Lie_{X_\xi})$
    acts as a well defined endomorphism of $\Comega(G)$, contrary to what was stated
    in~\cite[Prop.~4.5, \textit{iii.)}]{heins.roth.waldmann:2023a}. This is already
    wrong in the simplest non-trivial case of the additive reals $G = (\R,+)$. Indeed,
    consider the family of functions
    \begin{equation}
        f_\epsilon
        \colon
        \R \longrightarrow \C, \quad
        f_\epsilon(x)
        \coloneqq
        \frac{1}{x - \I \epsilon}
        \qquad
        \textrm{for }
        \epsilon > 0,
    \end{equation}
    which may be meromorphically extended to $\C$ with a simple pole at $x = \I
    \epsilon$. By virtue of
    \begin{equation}
        \exp(t \Lie_{X_1})
        f_\epsilon
        \at[\Big]{x}
        =
        \exp
        \biggl(
            t \frac{\D}{\D x}
        \biggr)
        f_\epsilon
        \at[\Big]{x}
        =
        \sum_{k=0}^{\infty}
        \frac{t^k}{k!}
        \frac{\D^k f_\epsilon}{\D x^k}
        (x)
    \end{equation}
    the Lie-Taylor series coincides with the usual Taylor series; this also follows
    from the Lie exponential of $G$ being given by the identity mapping. Consequently,
    the maximal choice of the radius $R$ for $f_\epsilon$ is $\epsilon$ and thus
    $\exp(t \Lie_{X_1})$ does not act on all $f_\epsilon$ for any $t \neq 0$.
\end{remark}

\begin{remark}[Uniform choices, {\cite[Cor.~4.9]{heins.roth.waldmann:2023a}}]
    \label{rem:UniformChoices}
    \index{Lie-Taylor!Uniform choices}
    Let $g \in G$. By what we have shown, the function $R_g
    \colon \liealg{g} \longrightarrow (0,\infty) \cup \{\infty\}$,
    \begin{equation}
        R_g(\xi)
        \coloneqq
        \sup
        \bigl\{
        R > 0
        \colon
        \eqref{eq:LieTaylorOneVariable}
        \; \textrm{holds whenever} \;
        \abs{t} < R
        \bigr\}
    \end{equation}
    is well defined. Taking another look at our construction for fixed $\xi \in \liealg{g}$, we see that the zero neighbourhood $U$ was independent of $\xi$ and thus
    \begin{equation}
        R_g
        \coloneqq
        \min_{\supnorm{\chi} \le 1}
        R_g(\chi)
        \ge
        R_g(\xi)
        >
        0,
    \end{equation}
    where we take the supnorm with respect to some linear algebraic basis of
    $\liealg{g}$ with $\xi = \basis{e}_1$. Indeed, given $\chi \in \liealg{g}$
    with $\supnorm{\chi} \le 1$ that is not a
    multiple of $\xi$, we may change basis such that $\basis{e}_2 = \chi$. Having a
    uniform choice $R_g > 0$ reflects the fact that power series converge in open
    sets. We will incorporate this fact in our formulation of the multivariable
    Lie-Taylor formula. Note that $R_g$ depends on the choice of our basis and thus
    has no intrinsic geometric meaning. However, $R_g > 0$ with respect to some basis
    does imply that $R_g > 0$ for \emph{any} basis. This, in turn, is equivalent to
    $\phi$ being real analytic at the point $g$. Indeed, this observation will be at
    the core of Proposition~\ref{prop:AnalyticVectors}.

    Having noted this, we may return to Remark~\ref{rem:ExponentialOperators}.
    By Proposition~\ref{prop:LieTaylorOneVariable}, the operator
    \begin{equation}
        \exp
        \bigl(
        \Lie(\xi)
        \bigr)
        \colon
        \Comega_{\abs{\xi}}(G)
        \longrightarrow
        \Comega(G)
    \end{equation}
    is well defined and corresponds to the pullback with the right translation
    $r_{\exp(\xi)}$, where
    \begin{equation}
        \label{eq:ComegaR}
        \gls{RealAnalyticMinumumRadius}
        =
        \bigl\{
        \phi \in \Comega(G)
        \;\big|\;
        \forall_{g \in G}
        \colon
        R_g
        \ge
        r
        \bigr\}
        \qquad
        \textrm{for any }
        r > 0.
    \end{equation}
    As already noted, the spaces $\Comega_r(G)$ for finite values of $r$ depend on a
    choice of basis and is thus not intrinsically interesting. Consequently and up to
    some technical subtleties regarding the placement of absolute values, our
    and the interest of \cite{heins.roth.waldmann:2023a} lies in the extreme case
    $r = \infty$, i.e. the intersection
    \begin{equation}
        \bigcap_{r > 0}
        \Comega_r(G),
    \end{equation}
    where all such operators are well
    defined. One may think of the Lie-Taylor formula \eqref{eq:LieTaylorOneVariable}
    as the Lie algebra action \eqref{eq:LieDerivativeOnLieAlg} integrating to the Lie
    group action
    \begin{equation}
        g
        \acts
        \phi
        =
        \phi \circ r_g
        =
        r_g^* \phi
    \end{equation}
    and as the intersection $\bigcap_{r > 0} \Comega_r(G)$ as the
    space of \emph{entire vectors} for the action. We refer to
    \cite[Rem.~4.16]{heins.roth.waldmann:2023a} and the references
    therein for a comprehensive discussion of these concepts, and will return to this
    point of view to resolve the integration problem in
    Section~\ref{sec:EntireVectors}. Note however that the
    spaces \eqref{eq:ComegaR} and also the notion of entire vectors already appeared
    and were studied in the series of papers \cite{goodman:1969a, goodman:1970a,
    goodman:1971a}.
\end{remark}

To simplify bookkeeping for the multivariable version and the rest of this chapter, we
fix some more notation, which was devised by Oliver Roth and introduced in
\cite[Sec.~4]{heins.roth.waldmann:2023a}. Recall that we have once and for all fixed a
basis $(\basis{e}_1, \ldots, \basis{e}_n)$ of $\liealg{g}$. Expanding on
\eqref{eq:LieDerivativeOnEnveloping}, we write
\begin{equation}
    \label{eq:LieDerivativeOfMulti}
    \Lie_{X_\alpha}
    \coloneqq
    \Lie(\alpha)
    \coloneqq
    \Lie_{X_{\alpha_1}}
    \circ
    \cdots
    \circ
    \Lie_{X_{\alpha_k}}
\end{equation}
as well as
\begin{equation}
    \underline{z}^\alpha
    \coloneqq
    z^{\alpha_1}
    \cdots
    z^{\alpha_k}
    \quad \textrm{for} \quad
    \underline{z}
    =
    (z^1, \ldots, z^n)
    \in
    \C^n
\end{equation}
for any $k$-tuple $\alpha \in \{1,\ldots,n\}^k \eqqcolon \N_n^k$. For $\phi \in \Cinfty(G)$ and $g \in G$, we call the formal series
\begin{equation}
    \label{eq:LieTaylorFormal}
    \index{Lie-Taylor!Formal}
    \gls{LieTaylor}
    \coloneqq
    \sum_{k=0}^{\infty}
    \frac{1}{k!}
    \sum_{\alpha \in \N_n^k}
    \bigl(
        \Lie(\alpha)\phi
    \bigr)(g)
    \cdot
    \underline{z}^\alpha
    \in
    \C\formal{\underline{z}}
\end{equation}
the Lie-Taylor series of $\phi$ at the point $g \in G$. This is warranted by the
following theorem. Recall that we make use of Einstein's summation
convention\footnote{Which is a marvellous gadget, whose self-correcting properties
make it an indispensable tool for any calculation involving coordinates.} to
implicitly sum over repeated indices.
\begin{theorem}[Multivariable Lie-Taylor formula,
{\cite[Cor.~4.6]{heins.roth.waldmann:2023a}}]
    \index{Lie-Taylor!Multivariable}
    \index{Taylor!Lie multivariable}
    \label{thm:LieTaylor}
    Let $G$ be a Lie group, $g \in G$ and $\phi \in \Comega(G)$. Fix a basis $(\basis{e}_1, \ldots, \basis{e}_n)$ of $\liealg{g}$. Then there is a radius $R_g > 0$ such that
    \begin{equation}
        \label{eq:LieTaylor}
        \phi
        \bigl(
            g \exp(x^j \basis{e}_j)
        \bigr)
        =
        \Taylor_\phi(\underline{x};g)
    \end{equation}
    for all $\underline{x} = (x^1, \ldots, x^n) \in \R^n$ with $\supnorm{\underline{x}} < R_g$.
\end{theorem}
\begin{proof}
    We define $R_g$ as discussed in Remark~\ref{rem:UniformChoices} corresponding to
    the basis \gls{Basis}. Consequently, we may apply
    \eqref{eq:LieTaylorOneVariable} to the direction $\xi \coloneqq x^j \basis{e}_j$
    for all $\underline{x} = (x^1, \ldots, x^n) \in \R^n$ with
    $\supnorm{\underline{x}} < R_g$. The formula \eqref{eq:LieTaylor} is now purely a
    matter of multiplying out
    within the universal enveloping algebra $\Universal^\bullet(\liealg{g})$, see
    again \eqref{eq:UniversalEnvelopeRelation}. Indeed,
    we have
    \begin{equation}
        \xi^k
        =
        (x^j \basis{e}_j)^k
        =
        \sum_{\alpha \in \N_n^k}
        \basis{e}_{\alpha_k}
        \cdots
        \basis{e}_{\alpha_1}
        \cdot
        \underline{x}^\alpha
    \end{equation}
    and thus
    \begin{equation}
        \Lie(\xi)^k
        =
        \sum_{\alpha \in \N_n^k}
        \Lie(\alpha)
        \cdot
        \underline{x}^\alpha.
        \tag*{\qed}
    \end{equation}
\end{proof}

Motivated by Remark~\ref{rem:InterpretationStarStd}, we now pull the absolute values
all the way inside of \eqref{eq:LieTaylor} to define the \emph{Lie-Taylor majorant} of
$\phi \in \Cinfty(G)$ at $g \in G$ as the formal power series
\begin{equation}
    \label{eq:LieTaylorMajorant}
    \index{Lie-Taylor!Majorant}
    \gls{Majorant}
    \coloneqq
    \sum_{k=0}^{\infty}
    c_{k}(\phi;g)
    \cdot
    \underline{z}^k
    \in
    \C\formal{\underline{z}},
\end{equation}
where
\begin{equation}
    \label{eq:LieTaylorMajorantCoefficients}
    \gls{MajorantCoefficients}
    \coloneqq
    \frac{1}{k!}
    \sum_{\alpha \in \N_n^k}
    \abs[\Big]
    {
        \bigl(
            \Lie(\alpha)\phi
        \bigr)(g)
    }
    \qquad
    \textrm{for all }
    k \in \N_0.
\end{equation}
For $g = \E$, we simply write $\gls{MajorantAtUnit} \coloneqq
\Majorant_{\phi}(\underline{z};\E)$
and $\gls{MajorantAtUnitCoefficients} \coloneqq c_k(\phi;\E)$. Note that we have
suppressed the choice of basis in our notation. This is warranted by the following
lemma.
\begin{lemma}
    \label{lem:MajorantVsBasis}
    Let $G$ be a Lie group, $g \in G$, $\phi \in \Comega(G)$ and $\mathcal{B},
    \mathcal{B}' \subseteq \liealg{g}$ be bases.
    \begin{lemmalist}
        \item If $\Majorant^\mathcal{B}_\phi(\argument;g)$ converges on all of $\C$, then also $\Majorant^{\mathcal{B}'}_\phi(\argument;g)$ converges on all of $\C$.
        \item If $\Majorant^\mathcal{B}_\phi(\argument;g)$ has positive radius of convergence, then so does $\Majorant^{\mathcal{B}'}_\phi(\argument;g)$.
    \end{lemmalist}
\end{lemma}
\begin{proof}
    Let $\mathcal{B} = (\basis{e}_1, \ldots, \basis{e}_n)$ and $\mathcal{B}' = (\basis{e}_1', \ldots, \basis{e}_n')$. Then $\basis{e}_j' = \sum_{\ell=1}^{n} \lambda_j^\ell \basis{e}_\ell$ and we set
    \begin{equation}
        M
        \coloneqq
        \max_{j,\ell=1,\ldots,n}
        \abs[\big]
        {\lambda_j^\ell},
    \end{equation}
    which is simply the modulus of the largest entry of the linear change of basis
    from $\mathcal{B}'$ to~$\mathcal{B}$. This yields the estimate
    \begin{equation}
        c_k^{\mathcal{B}'}
        (\phi;g)
        =
        \frac{1}{k!}
        \sum_{\alpha \in \N_n^k}
        \abs[\Big]
        {
            \bigl(
                \Lie
                (
                    \basis{e}_{\alpha_k}'
                    \tensor \cdots \tensor
                    \basis{e}_{\alpha_1}'
                )
                \phi
            \bigr)(g)
        }
        \le
        (Mn)^k
        \cdot
        c_k(\phi;g)
    \end{equation}
    for all $k \in \N_0$. Crucially, this grows like a constant raised to the power of
    $k$. Consequently, this factor may be absorbed into the complex variable
    $\underline{z}$, resulting in a scaling of the radius of convergence, and both
    claims follow.
\end{proof}

The upshot is that the choice of basis is inconsequential in the extremal situations
of infinite and positive but unspecified radius of convergence. Moving on, let
$\phi \in \Comega(G)$. Then
\begin{equation}
    \abs[\big]
    {\Taylor_\phi(\underline{z};g)}
    \le
    \Majorant_\phi
    \bigl(
        \supnorm{\underline{z}};g
    \bigr)
    \qquad
    \textrm{for all }
    \underline{z}
    \in
    \C^n,
\end{equation}
where either side might be equal to $\infty$ and we set $x \le \infty$ for all $x \in
\R \cup \{\infty\}$. This warrants our choice of terminology. In particular, if
the
Lie-Taylor majorant $\Majorant_\phi$ defines an entire function, then the formal
Lie-Taylor
series $\Taylor_\phi(\underline{z};g)$ converges for all $\underline{z} \in \C^n$. As
$\phi$ was assumed to be real analytic, Theorem~\ref{thm:LieTaylor} then
implies that
the
multivariable Lie-Taylor formula \eqref{eq:LieTaylor} holds for all $\underline{x} \in
\R^n$. More precisely, we get the following:
\begin{proposition}[Translation invariance,
{\cite[Thm.~4.17,~\textit{ii.)}]{heins.roth.waldmann:2023a}}]
    \label{prop:TranslationInvariance}
    \index{Entire!Translation invariance}
    Let $G$ be a connected Lie group and $\phi \in \Comega(G)$ such that
    $\Majorant_\phi(\argument;g_0) \in \Holomorphic(\C)$ for some $g_0 \in G$. Then
    also~$\Majorant_\phi(\argument;g) \in \Holomorphic(\C)$ for all $g \in G$.
\end{proposition}
\begin{proof}
    We sketch the argument. By left translation, we may assume $g_0 = \E$ by
    replacing $\phi$ with $\phi \circ \ell_{g_0}$. Writing $g = \exp(\xi_1) \cdots
    \exp(\xi_k)$ for suitable $\xi_1, \ldots, \xi_k \in \liealg{g}$ adhering the
    normalization $\supnorm{\xi} \le 1$, one may then apply
    \cite[Prop.~4.15,~\textit{iv.)}]{heins.roth.waldmann:2023a} $k$-times, which
    yields the claim.
\end{proof}

For this reason, we shall in the sequel always suppose connectedness of the Lie group
$G$. Another
pleasant consequence of this assumption is that
\begin{equation}
    \phi
    =
    0
    \quad \iff \quad
    \exists_{g \in G}
    \colon
    \Majorant_{\phi}(\argument; g)
    =
    0,
\end{equation}
as all coefficients \eqref{eq:LieTaylorMajorantCoefficients} are nonnegative. This
leads us to the notion of \emph{entire functions} on $G$. That is to say, we demand
that the Lie-Taylor majorant $\Majorant_\phi$ at the group unit defines an entire
function of one complex variable.
\begin{definition}[Entire functions, {\cite[Def.~4.10]{heins.roth.waldmann:2023a}}]
    \label{def:EntireFunction}
    \index{Entire functions}
    \index{Entire!Functions}
    \index{Lie-Taylor!Entire}
    Let $G$ be a connected Lie group.
    \begin{definitionlist}
        \item A real analytic function $\phi \in \Comega(G)$ is called entire
        if its Lie-Taylor majorant~$\Majorant_\phi$ at the group
        unit is holomorphic on all of $\C$.
        \item We denote the set of all entire functions on $G$ by
        \begin{equation}
            \gls{Entire}
            \coloneqq
            \bigl\{
                \phi \in \Comega(G)
                \colon
                \Majorant_{\phi}
                \in
                \Holomorphic(\C)
            \bigr\}.
        \end{equation}
        \item For $c \ge 0$, we define a family of norms on $\Entire_0(G)$ by
        \begin{equation}
            \label{eq:EntireSeminorms}
            \gls{SeminormsEntire}(\phi)
            \coloneqq
            \Majorant_{\phi}(c)
        \end{equation}
        and endow $\Entire_0(G)$ with the corresponding locally convex topology.
    \end{definitionlist}
\end{definition}

The advantage of working with $\Majorant_\phi$ instead of with the norms
\eqref{eq:EntireSeminorms} is that $\Majorant_\phi$ is a holomorphic function. This
simple shift of perspective allows for the utilization of many useful complex analytic
techniques, such as the Cauchy estimates and Montel's Theorem. We collect the most
important properties and symmetries of $\Entire_0(G)$. A systematic discussion and
detailed proofs may be found in \cite[Sec.~4]{heins.roth.waldmann:2023a}.
\begin{theorem}[Symmetries, {\cite[Thm.~4.17,
4.18]{heins.roth.waldmann:2023a}}]
    \index{Entire functions!Symmetries}
    \label{thm:Symmetries}
    Let $G$ be a connected Lie group.
    \begin{theoremlist}
        \item The set $\Entire_0(G)$ is a nuclear Fréchet algebra that is moreover
        a Montel space\footnote{Montel's Theorem holds: every bounded and
        closed subset of $\Entire_0(G)$ is compact.}, separable as well as reflexive.

        \item The group inversion $\gls{Inversion} \colon G \longrightarrow G$ acts
        isometrically on $\Entire_0(G)$ by pullback. More precisely,
        \begin{equation}
            \Majorant_{\phi \circ \inv}
            =
            \Majorant_\phi
            \qquad
            \textrm{for all }
            \phi \in \Comega(G).
        \end{equation}

        \item \label{item:TranslationInvariance}
        The group $G$ acts on $\Entire_0(G)$ via pullbacks with left and right
        translations in a strongly continuous fashion, i.e. the functions
        \begin{equation}
            \label{eq:MultiplicationActionStronglyContinuous}
            G
            \ni
            g
            \mapsto
            \ell_g^*
            \phi
            \in
            \Entire_0(G)
            \quad \textrm{and} \quad
            G
            \ni
            g
            \mapsto
            r_g^*
            \phi
            \in
            \Entire_0(G)
        \end{equation}
        are continuous for all $\phi \in \Entire_0(G)$.

        \item \label{item:DifferentiationInvariance}
        The universal enveloping algebra $\Universal^\bullet(\liealg{g}_\C)$ acts on
        $\Entire_0(G)$ via \eqref{eq:LieDerivativeOnEnveloping} by continuous
        linear mappings.

        \item \label{item:LieTaylorAbsoluteConvergence}
        The Lie-Taylor series \eqref{eq:LieTaylorFormal} is an absolutely convergent power series in $\Entire_0(G)$, i.e. converges absolutely in $\Entire_0(G)$ for any $\underline{z} \in \C^n$.
        \item The $\Entire_0(G)$-topology is finer than the $\Cinfty(G)$-topology. In particular, the evaluation functionals
        \begin{equation}
            \delta_{g,\xi}
            \colon
            \Entire_0(G) \longrightarrow \C, \quad
            \delta_{g,\xi}(\phi)
            \coloneqq
            \bigl(
                \Lie(\xi)
                \phi
            \bigr)
            (g)
        \end{equation}
        are continuous for all $g \in G$ and $\xi \in
        \Universal^\bullet(\liealg{g}_\C)$.
    \end{theoremlist}
\end{theorem}

Note that the statement \ref{item:LieTaylorAbsoluteConvergence} relies on
\ref{item:DifferentiationInvariance}, which in particular establishes that every
partial sum corresponding to the formal Lie-Taylor series
$\Taylor(\underline{z},\argument)$ is an element of $\Entire_0(G)$.
\begin{remark}[$R$-Entire functions, {\cite[Def.~4.13]{heins.roth.waldmann:2023a}}]
    \label{rem:GelfandShilov}
    \index{Gelfand-Shilov!Entire}
    \index{Entire functions!$R$-entire}
    Rescaling the coefficients $c_k(\phi)$ from
    \eqref{eq:LieTaylorMajorantCoefficients} with $k!^R$ for a fixed parameter $R \in
    \R$ yields a spectrum of locally convex spaces
    \begin{equation}
        \Entire_R(G)
        \coloneqq
        \biggl\{
            \phi
            \in
            \Comega(G)
            \;\Big|\;
            \forall_{c \ge 0}
            \colon
            \sum_{k=0}^{\infty}
            k!^R
            \cdot
            c_{k}(\phi)
            \cdot
            c^k
            <
            \infty
        \biggr\}.
    \end{equation}
    Note that $\Entire_R(G) \subseteq \Entire_S(G)$ whenever $S \le R$. They may be
    viewed as
    the Lie theoretic analogue of \emph{entire functions of finite order and minimal
    type} or Gelfand\footnote{Israïl Moyseyovich Gel'fand (1913-2009) was a
    Soviet-American functional analyst specialized in representation theory. For most
    of his life, he worked at the Moscow State University, until he emigrated to the
    United States in 1989, where he remained until his death.}-Shilov\footnote{Georgi
    Evgen'evich Shilov (1917-1975) was a Soviet functional analyst and a doctoral
    student of Gelfand. He continued to collaborate with him throughout his entire
    career and both proved numerous theorems on various flavours of generalized
    functions.}
    spaces. This is of course motivated by the case $G = \field{R}^n$, which we have
    discussed in Section~\ref{sec:RTopologiesPolynomial} at least for $R = 0$, and the
    case $R=1/2$ appears in Example~\ref{ex:StdOrdII} and Example~\ref{ex:StdOrdIII}.
    Indeed, in this case it is straightforward to see that
    \begin{equation}
        \index{Entire functions!Vector space}
        \Entire_R(\field{R}^n)
        =
        \widehat{\Sym}_R(\field{R}^n)
    \end{equation}
    as Fréchet algebras, see \cite[Lemma~5.6]{heins.roth.waldmann:2023a}. The absence
    of a notion of polynomials on~$G$ has lead us to model the completion directly.
    The case of $R = 0$ recovers our prior considerations and this is why we included
    the subscript $0$. For $R \ge 0$, the algebras~$\Entire_R(G)$ share all discussed
    properties with $\Entire_0(G)$. For $R < 0$, they degenerate in the sense that the
    evaluation functionals $\delta_g(\phi) \coloneqq \phi(g)$ are no longer continuous
    for $g \neq \E$ and the same is true for most of the actions discussed in
    Theorem~\ref{thm:Symmetries}, see \cite[Rem.~4.20]{heins.roth.waldmann:2023a}. The
    underlying problem is that the Lie-Taylor series of $\phi \in \Entire_R(G)$ with
    $R < 0$ does not necessarily have infinite radius of convergence any more and
    thus the Taylor data at the group unit might not describe its global behaviour. In
    the sequel, we will mostly be interested in the case~$R = 0$, as it results in the
    largest algebra with reasonable properties.
\end{remark}

Another important observation is that $\Entire_R$ constitutes a contravariant functor
from the category of connected Lie groups into the category of Fréchet algebras. That
is to say, there is the following compatibility with group morphisms.
\begin{proposition}[Functoriality of $\Entire_R$, {\cite[Prop.~4.21]{heins.roth.waldmann:2023a}}]
    \index{Entire functions!Functoriality}
    \label{prop:EntireFunctor}%
    Let $\Phi \colon G \longrightarrow H$ be a Lie group morphism between connected
    Lie groups and $R \in \R$. Then the pullback with $\Phi$ is a morphism of Fréchet
    algebras
    \begin{equation}
        \Phi^*
        \colon
        \Entire_R(H)
        \longrightarrow
        \Entire_R(G).
    \end{equation}
\end{proposition}

\index{Matrix elements}
\index{Representative functions}
Finally, we note that there always is a wealth of entire functions, at least for
linear\footnote{Lie groups admitting for a faithful representation on a finite
dimensional vector space.} Lie groups, namely the \emph{representative functions} or
\emph{matrix elements} corresponding to representations $\pi \colon
G \longrightarrow \gls{GeneralLinear}$ on finite dimensional vector spaces $V$. The
idea is that combining a vector $v \in V$ with a functional $\varphi \in V'$ yields
a function
\begin{equation}
    \label{eq:RepresentativeFunction}
    \gls{RepresentativeFunction}
    \colon
    G \longrightarrow \C, \quad
    \pi_{v,\varphi}(g)
    \coloneqq
    \varphi
    \bigl(
        \pi(g)v
    \bigr).
\end{equation}
For a systematic discussion of representative functions, we refer to
\cite[Sec.~4.3]{broecker.tomdieck:1985a}, and for their remarkable relationship with
joint eigenfunctions of the Laplacian and other invariant differential operators to
\cite[Sec.~10]{hall:2003a}. The choice of a basis facilitates a proof of the following.
\begin{theorem}[Representative functions, {\cite[Thm.~4.23]{heins.roth.waldmann:2023a}}]
    \index{Matrix elements!Entirety}
    \index{Representative functions!Entirety}
    \index{Entire functions!Representative functions}
    \label{thm:RepresentativeFunctions}
    Let
    \begin{equation}
        \pi \colon G \longrightarrow \gls{GeneralLinearComplex}
    \end{equation}
    be a representation of a connected Lie group $G$ and $R < 1$. Then
    \begin{equation}
        \label{eq:RepresentativeFunctionPairing}
        \pi_{k,\ell}
        \colon
        G \longrightarrow \C, \quad
        \pi_{k,\ell}(g)
        \coloneqq
        \big\langle
            \basis{e}_\ell,
            \pi(g) \basis{e}_k
        \big\rangle
    \end{equation}
    is in $\Entire_R(G)$ for all $k,\ell = 1, \ldots, d$. Here, $(\basis{e}_1, \ldots,
    \basis{e}_d)$ denotes the standard basis of $\C^d$ and~\gls{Pairing} is the
    euclidean scalar product on $\C^d$.
\end{theorem}

Notably, the restriction $R < 1$ is necessary already for the exponential
representation of the additive reals $(\R,+)$ on themselves by
\cite[Ex.~4.24]{heins.roth.waldmann:2023a}.
\begin{remark}[Representation theory and entire functions]
    \index{Entire functions!Representation theory}
    The fact that representative functions are entire is not a coincidence. It turns
    out that their Lie derivatives are in direct correspondence with the ones of the
    \emph{orbit mappings}
    \begin{equation}
        \gls{OrbitMap}
        \colon
        G \longrightarrow V, \quad
        \pi_v(g)
        \coloneqq
        \pi(g)v,
    \end{equation}
    where $v \in \C^n$ is a fixed vector.
    Extending Theorem~\ref{thm:RepresentativeFunctions} to suitable representations on
    infinite dimensional spaces then naturally leads to the notions of \emph{entire}
    and \emph{strongly entire vectors}. We will make this precise in the form of the
    universality Theorem~\ref{thm:Universality}.
\end{remark}
\begin{remark}[Representative functions in physics]
    \index{RepresentativeFunctions!Physics}
    Physically, Theorem~\ref{thm:RepresentativeFunctions} has the following
    interpretation: Going into a representation~$\pi$ of $G$ is the same as choosing a
    concrete realization of our system, as it lifts to the cotangent bundle $T^*G$ as
    a point transformation, which then also respects the symplectic structure. That is
    to say, it preserves the dynamics. The representative functions
    \begin{equation}
        \pi_{v,w}(g)
        \coloneqq
        \big\langle
            v,
            \pi(g)w
        \big\rangle
    \end{equation}
    then carry the interpretation of \emph{correlations} between states described by
    normalized vectors
    $v, w \in V$. In the case of $v = w$, one moreover obtains the \emph{expectation
    value} of $\pi(g)$ in the state $v$.
    Theorem~\ref{thm:RepresentativeFunctions} therefore asserts that these
    quantities are always observable, both in the classical and the quantum situation.
    This also makes it clear that for serious applications, we are going to have to
    generalize our result to representations on infinite dimensional spaces, which
    will ultimately lead us to Theorem~\ref{thm:Universality}.
\end{remark}

The principal purpose of the algebra of entire functions $\Entire_0(G)$ is to serve as
a natural partner of $\Sym_{1}^\bullet(\liealg{g}_\C)$ with respect to $\star_\std$,
which yields the following extension of
Theorem~\ref{thm:LieAlgebraStarProductContinuity}.
\begin{theorem}[Continuity of $\star_\std$, {\cite[Thm.~6.3]{heins.roth.waldmann:2023a}}]
    \index{Continuity!Standard ordered star product}
    \index{Standard ordered!Continuity}
    \label{thm:LieGroupStarProductContinuity}
    Let $G$ be a connected Lie group and fix $\hbar \in \C$, $R \ge 0$ and $R' \ge 1$.
    Then the standard ordered star product \eqref{eq:StarProductStd} is a well
    defined and
    continuous product
    \begin{equation}
        \star_\std
        \colon
        \bigl(
            \Entire_R(G)
            \tensor
            \Sym_{R'}^\bullet(\liealg{g}_\C)
        \bigr)
        \tensor
        \bigl(
        \Entire_R(G)
        \tensor
        \Sym_{R'}^\bullet(\liealg{g}_\C)
        \bigr)
        \longrightarrow
        \Entire_R(G)
        \tensor
        \Sym_{R'}^\bullet(\liealg{g}_\C),
    \end{equation}
    where $\tensor$ denotes the projective tensor product. Moreover,
    \begin{equation}
        \C \ni \hbar
        \; \mapsto \;
        P \star_\std Q
        \in
        \Entire_R(G) \tensor \Sym_{R'}^\bullet(\liealg{g}_\C)
    \end{equation}
    is an entire Fréchet holomorphic function for all $P,Q \in \Entire_R(G) \tensor
    \Sym_{R'}^\bullet(\liealg{g}_\C)$.
\end{theorem}
\begin{proof}
    We discuss the concepts underlying the proof and profitable technical
    simplifications one should employ. The missing details can be found in
    \cite[Thm.~6.3]{heins.roth.waldmann:2023a}.

    The first crucial observation is that working with projective tensor products means
    that it suffices to estimate products of \emph{factorizing tensors}. This is the
    \emph{infimum argument} from Proposition~\ref{prop:InfimumArgument}.

    The factorization \eqref{eq:StarProductFactorization}, which we have obtained from
    the associativity of the star product, then reduces the continuity of $\star_\std$
    to the continuity of the pointwise product, the Gutt star product~$\star_\Gutt$
    and the mixed product
    \begin{equation}
        \star_\std
        \colon
        \bigl(
        \Sym_{R'}(\liealg{g}_\C) \tensor 1
        \bigr)
        \tensor
        \bigl(
            \mathbb{1} \tensor \Entire_R(G)
        \bigr)
        \longrightarrow
        \Entire_R(G) \tensor \Sym_{R'}^\bullet(\liealg{g}_\C).
    \end{equation}
    The continuity of the pointwise product is a straightforward exercise, which was
    spelled out in \cite[Prop.~4.15,~\textit{iii.)}]{heins.roth.waldmann:2023a} and
    the continuity of $\star_\Gutt$ for $R' \ge 1$ is the main result of
    \cite[Sec.~3.1]{esposito.stapor.waldmann:2017a}, see again
    Theorem~\ref{thm:LieAlgebraStarProductContinuity}. Setting $\phi = \mathbb{1}$ and
    $q = 1$ in \eqref{eq:StarProductFactorization}, we are left with estimating the
    norm of the expression
    \begin{multline}
        (\mathbb{1} \tensor \xi_1 \vee \cdots \vee \xi_k)
        \star_\std
        (\psi \tensor 1) \\
        =
        \sum_{p=0}^{k}
        \biggl(
        \frac{\hbar}{\I}
        \biggr)^p
        \frac{1}{p! (k-p)!}
        \sum_{\sigma \in S_k}
        \Lie
        \bigl(
            \xi_{\sigma(1)} \tensor \cdots \tensor \xi_{\sigma(p)}
        \bigr)
        \psi
        \tensor
        (\xi_{\sigma(p+1)} \vee \cdots \vee \xi_{\sigma(k)}),
    \end{multline}
    for $\xi_1, \ldots, \xi_k \in \liealg{g}_\C$ and $\psi \in \Entire_R(G)$ by means of seminorms of both factors. As $\Entire_R(G)$ is a complex vector space, it suffices to consider real $\xi$.

    The next useful observation is that we may \emph{choose} which seminorms we use
    for $\liealg{g}$, as it constitutes a finite dimensional vector space. Here, a
    good choice turns out to be the~$\ell^1$-norm
    $\seminorm{p}$ associated to a fixed basis $(\basis{e}_1, \ldots, \basis{e}_n)$ of
    $\liealg{g}$, as the projective tensor product of two~$\ell^1$-norms, say
    associated to index sets $I$ and $J$, is simply given by the $\ell^1$-norm
    associated to $I \times J$. A proof of this fact can be found in
    \cite[Lem.~A.1]{cahen.gutt.waldmann:2020a}. In our situation, this leads to
    ``orthogonality within homogeneous degrees'', which manifests in \emph{equality}
    in the triangle inequality for certain seminorms. As we will need the resulting
    estimates later for the continuity of the holomorphic version of the star product, which
    we are going to establish within
    Theorem~\ref{thm:StarProductHolomorphicContinuity}, we make this precise in the
    form of Lemma~\ref{lem:ProjectiveTensorOfEll1}.

    The estimation of the mixed product then roughly goes as follows: First, one uses
    the triangle inequality to pull the seminorms $\seminorm{q}_{R,c} \tensor
    \seminorm{p}_{R',c'}$ all the way inside. By \eqref{eq:ProjectiveOnFactorizing},
    we get the factorization
    \begin{equation*}
        \bigl(
            \seminorm{p} \tensor \seminorm{q}
        \bigr)
        (v \tensor w)
        =
        \seminorm{p}(v) \cdot \seminorm{q}(w),
    \end{equation*}
    which in our case leads to $\seminorm{q}_{R,c}$ of a derivative and terms of the
    form \eqref{eq:ProjectivePowersEll1}. Spelling everything out, one is left with a
    shifted series, which may be estimated by the full series by virtue of the
    nonnegativity of all terms; here it is crucial that \eqref{eq:LieTaylorMajorant}
    converges and not just the one-variable Lie-Taylor series. We will come back to
    this problem in Section~\ref{sec:EntireVectors}. Finally, having \emph{equalities}
    in \eqref{eq:TruncatedPolynomialEstimate} allows to extend the estimate to
    arbitrary polynomials in the first factor.

    The holomorphic dependence on $\hbar$ is then comparatively easy to obtain.
    Throughout the estimates, one may treat $\hbar$ in a locally uniform way and each
    of the partial sums is simply a polynomial in $\hbar$. Our estimate moreover
    implies \emph{absolute} convergence of the series in~$\star_\std$, which overall
    yields the claim by Corollary~\ref{cor:FrechetPowerSeries}. We will later see
    that, in a purely
    complex setting, one may
    view $(\hbar, P, Q) \mapsto P \star_\std Q$ as a single Fréchet holomorphic
    mapping in the sense of Definition~\ref{def:Frechet}, matching with the programme
    outlined in Section~\ref{sec:QuantizationHolomorphic}.
\end{proof}

As we are going to need it, we conclude the section by discussing a precise
formulation of the ``orthogonality within fixed homogeneous degree'' from the proof of
\cite[Lemma~6.2]{heins.roth.waldmann:2023a}.
\begin{lemma}
    \label{lem:ProjectiveTensorOfEll1}
    \index{Projective tensor power!Of $\ell^1$}
    Let $G$ be a complex Lie group and fix a basis $(\basis{e}_1, \ldots, \basis{e}_n)$ of its Lie algebra. Let moreover $\seminorm{p}$ be the $\ell^1$-norm with respect to this basis and denote its $k$-th projective tensor power by $\seminorm{p}^k$. Fix $k \in \N_0$ and $1 \le j_1,\ldots, j_k \le n$.
    Then
    \begin{equation}
        \label{eq:ProjectivePowersEll1}
        \seminorm{p}^{\ell}
        \bigl(
        \basis{e}_{j_{\sigma(1)}}
        \vee \cdots \vee
        \basis{e}_{j_{\sigma(\ell)}}
        \bigr)
        =
        1
        =
        \seminorm{p}^{k}
        \bigl(
        \basis{e}_{j_1}
        \vee \cdots \vee
        \basis{e}_{j_k}
        \bigr)
    \end{equation}
    and
    \begin{equation}
        \label{eq:TruncatedPolynomialEstimate}
        \seminorm{p}_{R,c}
        \bigl(
        \basis{e}_{j_{\sigma(1)}}
        \vee \cdots \vee
        \basis{e}_{j_{\sigma(\ell)}}
        \bigr)
        =
        \biggl(
        \frac{\ell!}{k!}
        \biggr)^R
        c^{\ell-k}
        \cdot
        \seminorm{p}_{R,c}
        \bigl(
        \basis{e}_{j_1}
        \vee \cdots \vee
        \basis{e}_{j_k}
        \bigr)
    \end{equation}
    hold for all parameters $R, c \ge 0$, indices $1 \le \ell \le k$ and permutations
    $\sigma \in S_k$.
\end{lemma}
\begin{proof}
    The equality \eqref{eq:ProjectivePowersEll1} is immediate from
    \cite[Lem.~A.1]{cahen.gutt.waldmann:2020a}. Consequently,
    \begin{align*}
        \seminorm{p}_{R,c}
        \bigl(
        \basis{e}_{j_{\sigma(1)}}
        \vee \cdots \vee
        \basis{e}_{j_{\sigma(\ell)}}
        \bigr)
        &=
        \ell!^R
        c^\ell
        \cdot
        \seminorm{p}^{k-j}
        \bigl(
        \basis{e}_{j_{\sigma(1)}}
        \vee \cdots \vee
        \basis{e}_{j_{\sigma(\ell)}}
        \bigr) \\
        &=
        \biggl(
        \frac{\ell!}{k!}
        \biggr)^R
        c^{\ell-k}
        \cdot
        k!^R c^k
        \cdot
        \seminorm{p}^{k}
        \bigl(
        \basis{e}_{j_1}
        \vee \cdots \vee
        \basis{e}_{j_k}
        \bigr) \\
        &=
        \biggl(
        \frac{\ell!}{k!}
        \biggr)^R
        c^{\ell-k}
        \cdot
        \seminorm{p}_{R,c}
        \bigl(
        \basis{e}_{j_1}
        \vee \cdots \vee
        \basis{e}_{j_k}
        \bigr),
    \end{align*}
    which is exactly \eqref{eq:TruncatedPolynomialEstimate}.
\end{proof}

As noted before, the crucial point is that \eqref{eq:TruncatedPolynomialEstimate} is an \emph{equality} and not just an estimate.

\section{Universal Complexification of Lie Groups}
\label{sec:UniversalComplexification}
\epigraph{Thunder rolled. It rolled a six.}{\emph{Guards! Guards!} -- Terry Pratchett}
% !TeX root = ../Dissertation.tex
The principal idea behind complexification of a real Lie group $G$ is to find a complex
Lie group $G_\C$, whose properties are in controlled correspondence with the ones of
$G$.
\index{Hochschild, Gerhard}
\index{Chevalley, Claude}
\index{Complexification!Chevalley}
Historically, there have been quite a few different notions and constructions of
complexifications. The textbook \cite[Ch.~X~§13-14]{naimark.stern:1982} advocates for
a notion of complexification based on the Lie algebra complexification,
\cite[Ch.~VI~§VII-IX]{chevalley:2018a}
discusses the Chevalley\footnote{Claude Chevalley (1909-1984) was a French algebraic
group theorist,
a founding member of the Bourbaki group and the doctoral advisor of Gerhard
Hochschild. He introduced and studied what we nowadays refer to as Chevalley groups,
which are essentially Lie groups over finite fields.}
complexification of compact Lie groups\footnote{The construction is based on the Hopf
algebra of representative functions. In analogy to Milnor's exercise, one then
identifies the real characters of this algebra as the points of the underlying real
group by considering evaluation functionals. The complex characters then provide the
complexification.} and
\cite[Sec.~15.1]{hilgert.neeb:2012a} discusses the
universal complexification introduced by Hochschild \cite{hochschild:1966a}, which is
most suited for our needs.

The common denominator is that the Lie algebra of the
complexification should be isomorphic to the complexification of $\LieAlg(G)$ as a
vector space. While this holds in many important cases also for
universal
complexifications, we shall however see in Example~\ref{ex:DimensionDrop} that this
fails in general. This will result in some quite interesting technical complications, and
may be viewed as the price one pays for its guaranteed existence.

\index{Complexification!Vector Space}
Recall that given a -- not necessarily finite dimensional -- real vector space $V$,
its complexification may be defined as $\gls{ComplexifiedVectorSpace} \coloneqq V
\tensor_\R \C$ with its usual
addition and the linearly extended multiplication by scalars
\begin{equation}
    w \cdot (v \tensor z)
    \coloneqq
    v \tensor (zw)
    \eqqcolon
    (v \tensor z) \cdot w
    \qquad
    \textrm{for $v \in V$ and $z,w \in \C$.}
\end{equation}
By slight abuse of notation, one then writes $zv \coloneqq v \tensor z$ and uses
complex multiples of real vectors. An alternative model is given by $V \oplus V$ with
\begin{equation}
    z
    \cdot
    (v,v')
    \coloneqq
    \bigl(
        \RE(z) v
        -
        \IM(z)v',
        \RE(z)v'
        +
        \IM(z)v
    \bigr)
    \qquad
    \textrm{for $v,v' \in V$ and $z \in \C$.}
\end{equation}
This corresponds to the decomposition into real \gls{RealPart} and imaginary parts
\gls{ImaginaryPart}, which provides the
complex linear isomorphism
\begin{equation}
    \label{eq:VectorSpaceComplexificationVariants}
    V_\C
    \ni
    v \tensor z
    \mapsto
    \bigl(
        \RE(z)v,
        \IM(z)v
    \bigr)
    \in
    V \oplus V.
\end{equation}
Note that the original vector space $V$ may be identified with the real subspaces $V
\tensor \{1\}$ and $V \oplus \{0\}$, respectively. If now $\Phi \colon V
\longrightarrow W$ is a linear map into a complex vector space $W$, we may set
\begin{equation}
    \label{eq:VectorSpaceComplexification}
    \Phi_\C
    \colon
    V_\C \longrightarrow W, \quad
    \Phi_\C(v \tensor z)
    \coloneqq
    \Phi(v) \tensor z,
\end{equation}
which is the unique $\C$-linear extension of $\Phi$ to $V_\C$. It turns out that this
universal property characterizes the pair $(V_\C, V \hookrightarrow V_\C)$
uniquely up to unique isomorphism of complex vector spaces. As we are ultimately
interested
in relating entire functions on~$G$ with holomorphic functions on $G_\C$ and in view
of the functoriality of $\Entire_0$ from Proposition~\ref{prop:EntireFunctor}, having a
Lie group morphism from the group into its complexification will be crucial. This
leads us to the notion of universal complexification.
\begin{definition}[Universal Complexification of a Lie group,
{\cite{hochschild:1966a}}]
    \index{Complexification!Universal}
    \index{Complexification!Hochschild}
    \index{Universal!Complexification}
    \label{def:UniversalComplexification}%
    \; \\
    Let $G$ be a Lie group. A pair $(\gls{Complexification},
    \gls{ComplexificationMorphism})$ of a
    complex Lie group $G_\C$ and a Lie group morphism
    \begin{equation}
        \eta_G
        \colon G \longrightarrow G_\C
    \end{equation}
    is called the universal
    complexification of the Lie group $G$ if for every complex Lie group~$H$ and every
    Lie group morphism $\Phi \colon G \longrightarrow H$ there exists a unique
    holomorphic group morphism
    \begin{equation}
        \label{eq:HolomorphicExtension}
        \Phi_\C \colon
        G_\C \longrightarrow H
        \qquad
        \textrm{such that}
        \qquad
        \Phi_\C \circ
        \eta_G = \Phi.
    \end{equation}
\end{definition}

We suppress the subscript and write $\eta$ instead when only one complexification is involved. In terms of a commutative diagram, the universal complexification looks like this:
\begin{equation}
    \label{eq:UniversalComplexification}
    \begin{tikzcd}[column sep = huge, row sep = huge]
        G
        \arrow[d, "\eta"]
        \arrow[r, "\Phi"]
        &H \\
        G_\C
        \ar[swap,ru, dashed, "\exists_1 \; \Phi_\C"]
    \end{tikzcd}
\end{equation}
Note that there are two different kinds of arrows: While all arrows
are morphisms of real Lie groups, the mapping $\Phi_\C$ is additionally
holomorphic. The usual argument for universal objects
shows that the universal complexification of a given Lie group is
unique up to a unique biholomorphic group morphism, i.e. an isomorphism of complex Lie
groups. In the sequel, we shall thus speak of \emph{the} universal complexification by
slight abuse of terminology. We confirm that our preliminary considerations
are consistent with Definition~\ref{def:UniversalComplexification}.
\begin{example}[Universal Complexification of vector space groups]
    \index{Universal Complexification!Vector space group}
    \label{ex:ComplexificationVectorSpace}
    \; \\
    Let $G = (V,+)$ be a finite dimensional vector space, viewed as a real Lie group.
    Then the pair $(V_\C, V \hookrightarrow V_\C)$ is the universal complexification
    of $V$ as a Lie group. Indeed, let
    \begin{equation}
        \Phi
        \colon
        V \longrightarrow H
    \end{equation}
    be a Lie group morphism into a complex Lie group $H$. Choosing a basis
    $(\basis{e}_1, \ldots, \basis{e}_n)$ of $V$, we get real one parameter subgroups
    \begin{equation}
        \R
        \ni
        t
        \mapsto
        \Phi(t \basis{e}_k)
        \in
        H
    \end{equation}
    of $H$ for $k=1, \ldots, n$. By \cite[Thm.~9.2.15]{hilgert.neeb:2012a} there exist
    unique corresponding $\xi_1, \ldots, \xi_k \in \liealg{h}$ with
    \begin{equation}
        \Phi(t \basis{e}_k)
        =
        \exp_H(t \xi_k)
        \qquad
        \textrm{for all }
        t \in \R
        \textrm{ and }
        k=1, \ldots, n.
    \end{equation}
    As $H$ is a complex Lie group, we may then define
    \begin{equation}
        \Phi_\C
        (\basis{e}_k \tensor z)
        \coloneqq
        \exp_H(z \xi_k)
        \qquad
        \textrm{for all }
        z \in \C
        \textrm{ and }
        k=1, \ldots, n,
    \end{equation}
    which we extend to all of $V$ by demanding additivity. By construction, this
    provides a holomorphic extension of $\Phi$ to $V_\C$. This extension is unique by
    uniqueness of $\xi_1, \ldots, \xi_k$. Thus we have shown that the pair $(V_\C, V
    \hookrightarrow V_\C)$ is indeed the universal complexification of $V$, also as a
    Lie group.
\end{example}

The existence of universal complexifications is more involved, but as we shall see,
always guaranteed by Theorem~\ref{thm:UniversalComplexification}. To prove this, we
are going to utilize Lie's second and third Theorems. Taking its existence for
granted for a
moment, we collect some first abstract properties of universal complexification.
\begin{lemma}
    \label{lem:UniversalComplexificationMorphisms}
    Let $G,H$ be real Lie groups and $\Phi \colon G \longrightarrow H$ a Lie group
    morphism. Then there exists a unique holomorphic group morphism $\Phi_\C \colon
    G_\C \longrightarrow H_\C$ with
    \begin{equation}
        \Phi_\C \circ \eta_G
         =
         \eta_H \circ \Phi.
    \end{equation}
\end{lemma}
\begin{proof}
    The idea is that $\eta_H \circ \Phi$ constitutes a Lie group morphism from $G$ into the complex Lie group $H_\C$. By the universal property of $G_\C$, there exists a a unique holomorphic group morphism $(\eta_H \circ \Phi)_\C \colon G_\C \longrightarrow H_\C$ such that the diagram
    \begin{equation}
        \label{eq:UniversalComplexificationMorphisms}
        \begin{tikzcd}[column sep = huge, row sep = huge]
            G
            \arrow[d, "\eta_G"]
            \arrow[r, "\Phi"]
            \arrow[dr, "\eta_H \circ \Phi"]
            &H
            \arrow[d, "\eta_H"] \\
            G_\C
            \ar[r, dashed, swap, "\exists_1 \; (\eta_H \circ \Phi)_\C"]
            &H_\C
        \end{tikzcd}
    \end{equation}
    is commutative. That is, $(\eta_H \circ \Phi)_\C \circ \eta_G = \eta_H \circ
    \Phi$, which means that $\Phi_\C \coloneqq (\eta_H \circ \Phi)_\C$ is the
    desired holomorphic
    group morphism and it is uniquely determined by $\Phi$, as every solution makes
    the diagram \eqref{eq:UniversalComplexification} for $\eta_H \circ \Phi$ commute.
\end{proof}

We may rephrase Lemma~\ref{lem:UniversalComplexificationMorphisms} as functoriality.
Here we bump into a slight technical complication\footnote{I would like to thank Stefan
Waldmann for telling me about a similar problematic, which made me aware of what
follows and how to elegantly resolve it in the first place.}: To
define a functor, we
need to \emph{specify} a concrete realization of the universal complexification for
every given group $G$. That is to say, one needs to pick a representative within the
singular equivalence class. We get around invoking a suitably powerful version
of the
axiom of
choice by defining the functor by means of the \emph{construction} we will provide in
Theorem~\ref{thm:UniversalComplexification} to prove the existence of universal
complexifications.

Having dealt with this subtlety, we denote the category of real Lie
groups with Lie group morphisms as morphisms by \gls{LieCat} and the category of
complex Lie groups with holomorphic group morphisms as morphisms by
\gls{LieCatComplex}.
\begin{corollary}
    \label{cor:UniversalComplexificationFunctor}
    \index{Complexification!Functoriality}
    Universal complexification constitutes a covariant functor
    \begin{equation}
        \gls{ComplexificationFunctor}
        \colon
        \categoryname{Lie}
        \longrightarrow
        \categoryname{Lie}_\C.
    \end{equation}
    In particular, universal complexification $\argument_\C$ preserves isomorphisms.
\end{corollary}
\begin{proof}
    Let $G$ be a Lie group. Taking another look at \eqref{eq:UniversalComplexificationDiagram}, it is clear that $\Phi = \id_G$ induces $(\id_G)_\C = \id_{G_\C}$. Let moreover be $\Phi \colon G \longrightarrow H$ and $\Psi \colon H \longrightarrow K$ be Lie group morphisms. By Lemma~\ref{lem:UniversalComplexificationMorphisms}, we get the commutative diagram
    \begin{equation}
        \label{eq:UniversalComplexificationComposition}
        \begin{tikzcd}[column sep = huge, row sep = huge]
            G
            \arrow[d, "\eta_G"]
            \arrow[r, "\Phi"]
            \ar[rr,bend left = 25, "\Psi \circ \Phi"]
            &H
            \arrow[d, "\eta_H"]
            \arrow[r, "\Psi"]
            &K
            \arrow[d, "\eta_K"] \\
            G_\C
            \ar[r, dashed, "\exists_1 \; \Phi_\C"]
            \ar[rr,dashed,bend right = 25, "\exists_1 \; (\Psi \circ \Phi)_\C"]
            &H_\C
            \ar[r, dashed, "\exists_1 \; \Psi_\C"]
            &K_\C
        \end{tikzcd}
    \end{equation}
    with holomorphic group morphisms $\Phi_\C, \Psi_\C$ and $(\Psi \circ \Phi)_\C$.
    By uniqueness,
    \begin{equation}
        (\Psi \circ \Phi)_\C
        =
        \Psi_\C \circ \Phi_\C
    \end{equation}
    follows. The additional statement is clear.
\end{proof}

The next question we investigate is the following: Given Lie groups $G$ and $H$ with
$H \subseteq G$, is the universal complexification $H_\C$ of $H$ a
subgroup of the universal complexification $G_\C$ of $G$? Encoding the subgroup as an
injective group morphism $\iota \colon H \longrightarrow G$,
Lemma~\ref{lem:UniversalComplexificationMorphisms} yields a unique holomorphic group
morphism
\begin{equation}
    \label{eq:ComplexificationSubgroup}
    \iota_\C \colon H_\C \longrightarrow G_\C
    \qquad \textrm{with} \quad
    \iota_\C \circ \eta_H
    =
    \iota.
\end{equation}
Thus $H_\C$ is a subgroup of $G_\C$ if $\iota_\C$ is injective. In
Example~\ref{ex:DimensionDrop}, we will see that
non-injectivity of $\eta_G$ may lead to a negative answer to our initial question.
Conversely, assuming injectivity of $\eta_G$ ensures injectivity of $\eta_H$.
\begin{proposition}
    \label{prop:ComplexificationOfSubgroupsInjectivity}%
    Let $G$ be a Lie group with universal complexification $\eta_G \colon G \longrightarrow G_\C$ and $\iota \colon H \longrightarrow G$ a group morphism.
    \begin{propositionlist}
        \item \label{item:ComplexificationOfSubgroupsInjectivityLocal}
        If $\eta_G$ and $\iota$ are locally injective at the group unit, then so is
        $\eta_H \colon H \longrightarrow H_\C$.
        \item \label{item:ComplexificationOfSubgroupsInjectivity}
        If $\eta_G$ and $\iota$ are injective, then so is $\eta_H \colon H
        \longrightarrow H_\C$.
    \end{propositionlist}
\end{proposition}
\begin{proof}
    By Lemma~\ref{lem:UniversalComplexificationMorphisms} there exists a unique holomorphic group morphism $\iota_\C$ such that the diagram
    \begin{equation}
        \label{eq:ComplexificationOfSubgroups}
        \begin{tikzcd}[column sep = huge, row sep = huge]
            H
            \arrow[d, "\eta_H"]
            \arrow[r, "\iota"]
            &G
            \arrow[d, "\eta_G"] \\
            H_\C
            \ar[r, dashed, swap, "\iota_\C"]
            &G_\C
        \end{tikzcd}
    \end{equation}
    commutes. By assumption in \ref{item:ComplexificationOfSubgroupsInjectivityLocal},
    the edge $\eta_G \circ \iota$ is locally injective at the group unit, as
    $\iota$ preserves the unit as a group morphism. Thus the same is true for the
    other edge $\iota_\C \circ \eta_H$, necessitating local injectivity at $\E$.
    Alternatively, one may also
    apply the tangent functor to the diagram \eqref{eq:ComplexificationOfSubgroups}
    and evaluate at the group unit to argue with injectivity instead of local
    injectivity. The proof of part~\ref{item:ComplexificationOfSubgroupsInjectivity}
    is
    analogous.
\end{proof}

Thus, if $\eta_G$ is injective, one is within the situation $H \subseteq G
\subseteq G_\C$
and $H \subseteq H_\C$. The idea is now that the smallest complex closed subgroup
containing $H$ within $G_\C$ should be isomorphic to $H_\C$. Under certain additional
assumptions this is indeed true; to check the universal property directly, we take a
brief detour to group theory. Recall that a \emph{complement} of some subgroup $H
\subseteq G$ is another subgroup $K \subseteq G$ with $G = HK$ and $H \cap K =
\{\E\}$. That is to say, every element $g \in G$ may be uniquely written as a
product
\begin{equation}
    hk
    \qquad
    \textrm{with} \quad
    h \in H
    \quad \textrm{and} \quad
    k \in K.
\end{equation}
Note that, if $G = H \gls{Semidirect} N$ is a semidirect product,
then $N$ is always a complement of $H$ that is in addition a normal subgroup of $G$.
We will discuss this situation in detail and from a rather different angle as
Example~\ref{ex:ComplexificationOfProducts}.
\begin{proposition}[Complemented groups, {\cite{hall:1937a}}]
    \label{prop:Complements}
    \index{Group complements}
    \index{Group complements!Characterization}
    Let $G$ be a group and $H \subseteq G$ be a subgroup. Then the following are equivalent:
    \begin{propositionlist}
        \item \label{item:ComplementNormal}
        The subgroup $H$ admits a complement that is normal in $G$.
        \item \label{item:ComplementQuotient}
        We have $H \cong G / N$ for some normal subgroup $N \subseteq G$.
        \item \label{item:ComplementExtension}
        Every group morphism $\phi \colon H \longrightarrow L$ into another group $L$ admits an extension to a group morphism $\Phi \colon G \longrightarrow L$.
    \end{propositionlist}
\end{proposition}
\begin{proof}
    The equivalence of \ref{item:ComplementNormal} and \ref{item:ComplementQuotient} is clear. Let $N \subseteq G$ be a normal subgroup that is also a complement of $H$ and let $\phi \colon H \longrightarrow L$ be a group morphism into another group~$L$. We define
    \begin{equation}
        \label{eq:ComplementExtension}
        \Phi
        \colon
        G
        =
        HN
        \longrightarrow L, \quad
        \Phi(hn)
        \coloneqq
        \phi(h),
    \end{equation}
    where $h \in H$ and $n \in N$. If now $h_1, h_2 \in H$ and $n_1, n_2 \in N$, then
    there exists some element $n_3 \in N$ with
    \begin{equation}
        h_1 n_1 h_2 n_2
        =
        h_1 h_2 n_3
    \end{equation}
    as by normality $N h_2 = h_2 N$. Consequently, we have
    \begin{equation}
        \Phi
        \bigl(
            h_1 n_1 h_2 n_2
        \bigr)
        =
        \Phi(h_1 h_2 n_3)
        =
        \phi(h_1 h_2)
        =
        \phi(h_1) \phi(h_2)
        =
        \Phi(h_1 n_1)
        \Phi(h_2 n_2)
    \end{equation}
    and thus $\Phi$ is a group morphism with $\Phi \at{H} = \phi$. Alternatively, one
    may use \ref{item:ComplementQuotient} with the quotient projection $\pi \colon G
    \longrightarrow H$ and define $\Phi \coloneqq \phi \circ \pi$ as the pullback of
    $\phi$ along $\pi$. This recovers the $\Phi$ we have defined. Assume conversely
    that \ref{item:ComplementExtension} holds and let $\phi \colon H \longrightarrow
    H$ be $\id_H$. By assumption, there exists a group morphism $\Phi \colon G
    \longrightarrow H$ with $\Phi(h) = h$ for all $h \in H$. We set $N \coloneqq \ker
    \Phi$, which yields a normal subgroup of $G$. If $g \in H \cap N$, then
    \begin{equation}
        g
        \in
        \ker \phi
        =
        \{\E_H\},
        \qquad
        \textrm{i.e.}
        \quad
        H \cap N
        =
        \{\E_H\}.
    \end{equation}
    If $g \in G$, then $n \coloneqq
    \Phi(g)^{-1}g$ fulfils
    \begin{equation}
        \Phi(n)
        =
        \E_H
        \quad \textrm{and} \quad
        \Phi(g)
        n
        =
        \Phi(g)
        \Phi(g)^{-1}
        g
        =
        g
    \end{equation}
    and thus $G = HN$. This completes the proof.
\end{proof}

Taking $H = \{\E\}$ as the trivial group, it is clear that the extension
\ref{item:ComplementExtension} is highly non-unique in general: Any group
morphism $\Phi \colon G \longrightarrow L$ restricts to the trivial morphism
\begin{equation}
    \phi
    \colon
    \{\E\}
    \longrightarrow
    L
\end{equation}
The preceding considerations pass to topological groups and continuous group
morphisms without complications and thus in particular to Lie groups.
\begin{corollary}
        \index{Group complements!Topological}
        \label{prop:ComplementsTopological}
        Let $G$ be a topological group and $H \subseteq G$ be a closed subgroup. Then the following are equivalent:
        \begin{corollarylist}
            \item \label{item:ComplementTopologicalNormal}
            The subgroup $H$ admits a complement that is normal and closed in $G$.
            \item \label{item:ComplementTopologicalQuotient}
            We have $H \cong G / N$ as topological groups for some normal subgroup $N \subseteq G$.
            \item \label{item:ComplementTopologicalExtension}
            Every continuous group morphism $\phi \colon H \longrightarrow L$ into another group $L$ admits an extension to a continuous group morphism $\Phi \colon G \longrightarrow L$.
        \end{corollarylist}
\end{corollary}

Having a normal and closed complement now yields a sufficient condition for $H_\C \subseteq G_\C$.
\begin{proposition}
    \index{Universal Complexification!Subgroups}
    Let $G$ be a connected Lie group and $H \subseteq G$ be a closed subgroup with
    \begin{equation}
        0 < \gls{Dimension} H < \dim G
    \end{equation}
    admitting a normal as well as closed complement. Assume moreover that $\eta_G
    \colon G \longrightarrow G_\C$ is injective. Then $H_\C$ is isomorphic to a
    subgroup of $G_\C$.
\end{proposition}
\begin{proof}
    As already noted, injectivity of $\eta_G$ lets us view $H$ as a subgroup of
    $G_\C$.
    This allows us to consider the smallest complex closed subgroup $H_\C$
    containing
    $H$ within $G_\C$
    and $\eta_H \colon H \longrightarrow H_\C$ be the inclusion of $H$ into
    $H_\C$.\footnote{Which is, in accordance with what we expect in view of
    Proposition~\ref{prop:ComplexificationOfSubgroupsInjectivity}, is indeed
    injective.} We verify that $(H_\C, \eta_H)$ possesses the universal property of
    the universal complexification. To this end, let $\phi \colon H \longrightarrow L$
    be a Lie group morphism into a complex Lie group $L$. Then, by assumption and
    Proposition~\ref{prop:ComplementsTopological}, there exists an extension of $\phi$
    to a Lie group morphism $\Phi \colon G \longrightarrow L$. By the universal
    property of $G_\C$, we may extend $\Phi$ even further to a unique holomorphic
    group morphism~$\Phi_\C \colon G_\C \longrightarrow L$. Setting $\phi_\C
    \coloneqq
    \Phi_\C \at{H_\C}$ yields a holomorphic group morphism with
    \begin{equation*}
        \phi_\C(h)
        =
        \Phi_\C(h)
        =
        \Phi(h)
        =
        \phi(h)
        \qquad
        \textrm{for all }
        h \in H.
    \end{equation*}
    It remains to prove the uniqueness of $\phi_\C$, which turns out to be a bit more
    involved. If the subgroup $N \subseteq G$ is a closed normal complement of $H$, we
    may define $N_\C$ in the same manner as $H_\C$. Moreover,
    \cite[Exercise~15.1.3(d)]{hilgert.neeb:2012a} asserts that the resulting closed
    complex subgroup $N_\C \subseteq G_\C$ is normal again. We first establish that
    $H_\C$ and $N_\C$ are complementary within $G_\C$. To prove this, we tacitly use
    two -- rather reasonable -- facts about universal complexification, which we
    are going to demonstrate in
    Theorem~\ref{thm:UniversalComplexification} and
    Lemma~\ref{lem:ComplexificationEtaLocallyInjective} without using what we are
    proving right now.
    \begin{lemmalist}
        \item If $G$ is connected, then so is $G_\C$.
        \item If the morphism $\eta$ is locally injective, then the Lie algebra of the
        universal complexification $G_\C$ is given by the complexification
        $\liealg{g}_\C$ of the Lie algebra $\liealg{g}$ of $G$.
    \end{lemmalist}

    Assuming this and recalling the injectivity of $\eta_H$ and $\eta_N$ by
    Proposition~\ref{prop:ComplexificationOfSubgroupsInjectivity}, we see that the
    condition $H \cap
    N = \{\E\}$ implies $\liealg{g} \cap \liealg{n} = \{0\}$ and thus $\liealg{g}_\C
    \cap \liealg{n}_\C = \{0\}$, as well. Exponentiating, this means that there exist
    zero neighbourhoods $U_H \subseteq H_\C$ and $U_N \subseteq N_\C$ with $U_H \cap
    U_N = \{\E\}$, which implies $H_\C \cap N_\C = \{\E\}$.

    Similarly, investing moreover $G = HN$ implies that $\liealg{g} = \liealg{h}
    \oplus \liealg{n}$ as vector spaces and thus $\liealg{g}_\C = \liealg{h}_\C \oplus
    \liealg{n}_\C$. Consequently, the product $\chi \coloneqq \exp_{H_\C} \times
    \exp_{N_\C}$ maps $\liealg{g}_\C$ onto a neighbourhood of the group unit
    within $G_\C$, as a
    suitable restriction constitutes an exponential chart of the second kind by
    \cite[Lem.~9.2.6~(ii.)]{hilgert.neeb:2012a}. By connectedness of $G$,
    also~$G_\C$
    is connected and thus every element of $G_\C$ is a finite product of elements in
    the image of $\chi$, proving $G_\C = H_\C N_\C$ by normality of $N_\C$ as
    desired.

    Let now $\psi_\C \colon H_\C \longrightarrow L$ be another holomorphic
    group morphism with $\psi_\C(h) = \phi(h)$ for all $h \in H$. By what we have
    shown and Proposition~\ref{prop:ComplementsTopological}, we may extend
    $\psi_\C$ to a group morphism $\Psi_\C \colon G_\C \longrightarrow L$;
    taking another look at \eqref{eq:ComplementExtension}, we moreover note
    that our construction even yields a holomorphic extension $\Psi_\C$. But then
    \begin{equation}
        \Psi_\C(hn)
        =
        \psi_\C(h)
        =
        \phi(h)
        =
        \Phi(hn)
    \end{equation}
    holds for all $h \in H$ and $n \in N$. As $\Phi_\C$ was the unique holomorphic
    extension of the morphism ${\Phi \colon G \longrightarrow L}$ to $G_\C$, this
    implies $\Psi_\C = \Phi_\C$ and thus $\phi_\C = \psi_\C$ by restriction to the
    subgroup $H_\C$.
\end{proof}

\begin{corollary}
    Let $G$ be a connected Lie group and $H \subseteq G$ be a closed subgroup with
    \begin{equation}
        0 < \dim H < \dim G
    \end{equation}
    admitting a normal as well as closed complement $N \subseteq G$. Assume moreover
    that $\eta_G$ is locally injective. Then $H_\C$ and $N_\C$ are complementary
    subgroups of~$G_\C$ and $N_\C$ is normal.
\end{corollary}

This leaves us with the problem of when $\eta_G \colon G \longrightarrow G_\C$ is
injective. We postpone this investigation until we have a concrete construction
of $G_\C$. Instead, we proceed with a simple example that
covers the case of discrete subgroups and will prove useful for the understanding of
universal complexifications of disconnected groups.
\begin{example}[Discrete Groups]
    \label{ex:DiscreteGroups}%
    \index{Complexification!Discrete groups}
    Let $G$ be a discrete Lie group, i.e. a countable group endowed
    with the discrete topology. As $\{g\}$ is an open subset of
    $G$ for all $g \in G$, we get an atlas of $G$ by means of the maps
    \begin{equation}
        z_g
        \colon
        \{g\}
        \longrightarrow
        \{0\}, \quad
        z_g(g) \coloneqq 0.
    \end{equation}
    The idea is now to view $\{0\} \cong \C^0$ as a
    zero dimensional \emph{complex} vector space, so we even
    have a holomorphic atlas for $G$. Consequently, we claim that $G =
    G_\C$ for similar reasons: any Lie group morphism $\Phi \colon G
    \longrightarrow H$ into some complex Lie group $H$ is holomorphic for this atlas:
    the compositions $\Phi \circ z_g^{-1}$ are constant.

    While this example is trivial, we are
    going encounter this situation in the course of the construction of the universal
    complexification. The component groups $G / G_0$ of any Lie group $G$
    are discrete, where~\gls{UnitComponent} denotes the connected component of the
    group unit. The
    crucial idea is then going to be that they should not change under
    complexification beyond the alteration of perspective just described. These
    ideas are
    at the core of Remark~\ref{rem:ComplexificationDisconnected}.
\end{example}

Our next goal is to prove that every Lie group admits a universal complexification. We begin by restating Lie's seminal theorems in the incarnations we will need.
\begin{theorem}[Lie's second Theorem, {\cite[Thm.~9.5.9]{hilgert.neeb:2012a}} or {\cite[Thm.~1.8.3]{duistermaat.kolk:2000a}}]
    \index{Lie's Theorems!Second Theorem}
    \label{thm:Lie2}%
    Let $G$ and $H$ be Lie groups. Assume furthermore
    that $G$ is connected as well as simply connected and let $\phi
    \colon \liealg{g} \longrightarrow \liealg{h}$ be a morphism of Lie
    algebras. Then there exists a unique Lie group morphism $\Phi \colon G
    \longrightarrow H$ such that $T_\E \Phi = \phi$. If $G$ and $H$ are complex Lie groups and $\phi$ is $\C$-linear, then $\Phi$ is holomorphic.
\end{theorem}

Diagrammatically, the situation is thus the following:
\begin{equation}
    \begin{tikzcd}[column sep = huge, row sep = huge]
    G
    \arrow[d, "\exp_G", swap, leftarrow]
    \arrow[r, dashed, "\exists_1 \; \Phi"]
    &H
    \arrow[d, "\exp_H", leftarrow] \\
    \liealg{g}
    \arrow[r, "\phi"]
    &\liealg{h}
    \end{tikzcd}
\end{equation}
%\begin{proof}[Idea]
%    Integrate the graph of $\phi$ to an (immersed) subgroup $K
%    \hookrightarrow G \times H$ and define $\Phi \colon G \cong K
%    \overset{\pr_2}{\longrightarrow} H$. Here the simply connectedness
%    of $G$ is required to obtain $K \cong G$ from the invertibility of
%    $\liealg{g} \hookrightarrow \liealg{k} = \graph(\phi) \subseteq
%    \liealg{g} \times \liealg{h}$.
%\end{proof}

As a first application, we provide a simple criterion to identify universal
complexifications in practice.
\begin{proposition}
    \label{prop:ClassicalComplexifications}
    Let $H$ be a connected and simply connected complex Lie group and let $G \subseteq
    H$ be a real closed connected subgroup such that $\LieAlg(G)_\C \cong \LieAlg(H)$.
    Then $G_\C \coloneqq H$ with the inclusion $\eta \colon G \hookrightarrow H$ is
    the universal complexification of $G$.
\end{proposition}
\begin{proof}
    Let $\Phi \colon G \longrightarrow H'$ be a Lie group morphism into another complex Lie group $H'$. By the universal property of the algebra complexification, its tangent map
    \begin{equation}
        \phi \coloneqq T_\E \Phi \colon \LieAlg(G) \longrightarrow \LieAlg(H')
    \end{equation}
    has a unique extension $\phi_\C \colon \LieAlg(G)_\C = \LieAlg(H) \longrightarrow \LieAlg(H')$ to a morphism of complex Lie algebras. By Theorem~\ref{thm:Lie2}, we may integrate $\phi_\C$ to a unique holomorphic group morphism $\Phi_\C \colon H \longrightarrow H'$ with tangent map $\phi_\C$. By construction, we have
    \begin{equation}
        \label{eq:ClassicalComplexificationsProof}
        T_\E
        \bigl(
            \Phi_\C \circ \eta
        \bigr)
        =
        \phi_\C
        \circ
        \iota'
        =
        \phi
        =
        T_\E
        \Phi
        \tag{$\star$}
    \end{equation}
    where $\iota'$ is the inclusion $\LieAlg(G) \hookrightarrow \LieAlg(H)$. By
    connectedness of $G$, this implies
    \begin{equation*}
        \Phi_\C \circ \eta
        =
        \Phi.
    \end{equation*}
    If $\Psi \colon H \longrightarrow H'$ constitutes another holomorphic group
    morphism with
    $\Psi \circ \eta = \Phi$, then its tangent map fulfils
    \begin{equation*}
        T_\E \Psi
        \at{\LieAlg(G)}
        =
        T_\E \Phi
    \end{equation*}
    by \eqref{eq:ClassicalComplexificationsProof}. As $\LieAlg(G)$
    generates $\LieAlg(H) = \LieAlg(G)_\C$ as a complex vector space, this in turn
    implies $T_\E \Psi = T_\E \Phi_\C$ by $\C$-linearity and thus $\Psi = \Phi_\C$ by
    connectedness.
\end{proof}

We resume our review of Lie's theorems.
\begin{theorem}[Cartan-Lie's Theorem, {\cite[Thm.~9.4.11]{hilgert.neeb:2012a}},
{\cite[Thm.~1.14.3]{duistermaat.kolk:2000a}}]
    \index{Lie's Theorems!Cartan-Lie's Theorem}
    \index{Lie's Theorems!Lie's third Theorem}
    \label{thm:Lie3}%
    Let $\liealg{g}$ be a finite dimensional Lie algebra.
    \begin{theoremlist}
        \item \label{item:Lie3Existence}
        There exists a connected and simply connected Lie
        group $G$ with $\LieAlg(G) = \liealg{g}$.
        \item The group $G$ from \ref{item:Lie3Existence} is unique up to Lie group
        isomorphism.
        \item If $\liealg{g}$ is complex, then $G$ from \ref{item:Lie3Existence} may
        be endowed with the structure of a complex Lie group. In this case, $G$ is
        unique up to a biholomorphic group isomorphism.
    \end{theoremlist}
\end{theorem}

Putting both of Lie's theorems together and demanding $H$ to be simply
connected gives now the following diagram:
\begin{equation}
    \label{eq:Lie2And3}
    \begin{tikzcd}[column sep = huge, row sep = huge]
    G
    \arrow[d, "\exp_G", swap, leftarrow]
    \arrow[r, dashed, "\exists_1 \; \Phi"]
    &\exists_1 \; H
    \arrow[d, dashed, "\exists_1 \; \exp_H", leftarrow] \\
    \liealg{g}
    \arrow[r, "\phi"]
    &\liealg{h}
    \end{tikzcd}
\end{equation}
Going even further, one can integrate any Lie algebra morphism to a
Lie group morphism between the corresponding connected and simply
connected Lie groups. Note that to apply Theorem~\ref{thm:Lie2},
we needed $G$ to be simply connected regardless.

The idea to construct a
universal complexification in the sense of
Definition~\ref{def:UniversalComplexification} is now to
integrate the embedding $\iota \colon \liealg{g} \hookrightarrow \liealg{g}_\C$
of vector spaces in the sense of \eqref{eq:Lie2And3}. Here,
\gls{ComplexifiedLieAlgebra} is the vector space complexification of $\liealg{g}$
endowed with the
Lie bracket
\begin{equation}
    \label{eq:ComplexificationLieBracketTensor}
    \big[
        \xi \tensor z,
        \chi \tensor w
    \big]
    =
    [\xi, \chi]
    \tensor
    (z \cdot w)
    \quad \iff \quad
    \big[
        z \xi,
        w \chi
    \big]
    =
    (z´w)
    \cdot
    [\xi, \chi]
\end{equation}
for $\xi, \eta  \in \liealg{g}$ and $z, w \in \C$ or in the other model
\begin{equation}
    \label{eq:ComplexificationLieBracketSum}
%    \big[
%        (\xi_1, \xi_2),
%        (\chi_1, \chi_2)
%    \big]
%    =
%    \bigl(
%        [\xi_1, \chi_1],
%        [\xi_1, \chi_2]
%    \bigr)
%    \quad \iff \quad
    \big[
        \xi_1 + \I \xi_2,
        \chi_1 + \I \chi_2
    \big]
    =
    [\xi_1, \chi_1]
    +
    \I [\xi_2, \chi_2]
\end{equation}
for $\xi_1, \xi_2, \chi_1, \chi_2 \in \liealg{g}$. Note that this is consistent both
with \eqref{eq:VectorSpaceComplexificationVariants} and with extending the adjoint
mappings
\begin{equation*}
    \gls{AdjointAlgebra}
    \colon
    \liealg{g} \longrightarrow \liealg{g}, \quad
    \ad_\xi(\chi)
    \coloneqq
    [\xi,\chi]
\end{equation*}
as linear maps to $\liealg{g}_\C$ for all $\xi \in \liealg{g}$ by the universal
property.

\index{Covering}
However, we only get the diagram \eqref{eq:Lie2And3} if $G$ is simply connected, which is a rather strong
assumption we do not want to make. To resolve this, we will first complexify the universal
covering group $\widetilde{G}$ of $G$ as an intermediate step and quotient out the
redundancy afterwards. Recall that a covering is a locally trivial fiber bundle
\begin{equation*}
    \gls{Covering}
    \colon
    X \longrightarrow Y
\end{equation*}
in the category of topological spaces, whose fibers are discrete. In this situation,
we say that $X$ covers $Y$. Unwrapping the definition, this means that $X, Y$ are
topological spaces, the projection $p$ is continuous and that the following
holds:

Every $y \in Y$ has an open neighbourhood $V$, whose preimage $p^{-1}(V)$ is the
disjoint union of finitely many open sets $U_1, \ldots, U_k \subseteq X$ such that $p
\at{U_j} \colon U_j \longrightarrow Y$ is a homeomorphism
for $j=1, \ldots, k$. One calls $V$ is an evenly covered neighbourhood of $y$ in this
case. Note that the integer $k$ typically depends on the point $y$, but is locally
constant. We are going to utilize both the discreteness of the fibers
and the local bijectivity in the sequel. For a comprehensive discussion of covering
spaces for topological spaces, we refer to the textbooks
\cite[Sec.~53]{munkres:1985a}, \cite[Sec.~1.3]{hatcher:2002a} and
\cite[Appendix~A]{hilgert.neeb:2012a}.

\index{Universal!Covering}
\index{Covering!Universal}
Under mild connectedness assumptions\footnote{Namely, path connectedness, local path
connectedness and semi-locally simple connectedness. Many authors refer to this
situation as \emph{well connected}.} on $Y$, there is always a simply connected
topological space $\widetilde{Y}$ covering $Y$ by
\cite[Thm.~A.2.12]{hilgert.neeb:2012a}. It is called universal covering, as it covers
all other covering spaces of $Y$ by a straightforward application of the Lifting
Theorem \cite[Thm.~A.2.9]{hilgert.neeb:2012a}, and as it is unique up to
fiber preserving homeomorphism by \cite[Cor.~A.2.10]{hilgert.neeb:2012a}. Taking $Y
\coloneqq G$ as a connected Lie group suffices to guarantee the existence of a
universal covering space \gls{CoveringGroup}, which turns out to be a Lie group itself
by \cite[Cor.~9.4.7]{hilgert.neeb:2012a}. Moreover, the covering morphism $p \colon
\widetilde{G} \longrightarrow G$ is then a group morphism. By local injectivity, we
moreover have $\LieAlg(G) = \LieAlg(\widetilde{G})$.
\index{Fundamental Group}
Finally, for our purposes the fundamental group \gls{FundamentalGroup} of $G$ is
defined as the Lie subgroup $\pi_1(G) \coloneqq \ker p$. It is always a central
discrete
subgroup of $\tilde{G}$ by \cite[Thm.~9.5.4]{hilgert.neeb:2012a}. We have gathered all
ingredients to prove the existence of universal complexifications.
\begin{theorem}[Universal Complexification, {\cite[Thm.~15.1.4]{hilgert.neeb:2012a}}]
    \label{thm:UniversalComplexification}%
    \index{Complexification!Construction}
    Let $G$ be a Lie group. Then there exists a universal complexification
    $(G_\C, \eta)$ of $G$ and it is unique up to a unique biholomorphic group morphism.
\end{theorem}
\begin{proof}
    We have already noted the uniqueness. Assume $G$ to be connected. Denoting the universal covering group
    of $G$ by $p \colon \widetilde{G} \longrightarrow G$, we have the
    commutative diagram
    \begin{equation*}
        \label{eq:UniversalComplexificationProof0}
        \begin{tikzcd}
%            [column sep = huge]
        G
        \arrow[dr, "\exp_G", leftarrow,swap]
        \arrow[rr, "p", leftarrow]
        &\,
        &\widetilde{G}
        \arrow[dl, "\exp_{\widetilde{G}}", leftarrow]
        \arrow[r, dashed, "\widetilde{\eta}"]
        &\widetilde{G}_\C
        \arrow[d, dashed, "\exp_{\widetilde{G}_\C}", leftarrow] \\
        \,
        &\liealg{g}
        \arrow[rr, "\iota", hookrightarrow, swap]
        &\,
        &\liealg{g}_\C
        \end{tikzcd}
        \tag{$\#$}
    \end{equation*}
    where the right half is \eqref{eq:Lie2And3} with $\phi \coloneqq
    \iota$, i.e. $\iota(\xi) = \xi$. Note that the universal covering group is simply
    connected, so the assumptions of Theorem~\ref{thm:Lie2} are
    indeed fulfilled and we get $\widetilde{\eta}, \iota$ and $\widetilde{G}_\C$ such that \eqref{eq:UniversalComplexificationProof0} commutes.

    By local bijectivity of the covering $p$, its
    kernel $\pi_1(G) \subseteq \widetilde{G}$ is a discrete normal
    subgroup. Indeed, we even have $\pi_1(G) \subseteq
    \Center(\widetilde{G})$, where
    \begin{equation}
        \gls{Center}
        \coloneqq
        \bigl\{
            g \in G
            \;|\;
            \forall_{h \in G}
            \colon
            gh = hg
        \bigr\}
    \end{equation}
    denotes the center of a group $G$, which always constitutes a normal
    subgroup:

    Given $h \in \pi_1(G)$ and $g \in \widetilde{G}$ we
    find a path $\gamma$ from $g$ to $\E$ in $\widetilde{G}$. This induces a path
    \begin{equation}
        \sigma(t) \coloneqq \gamma(t) h \gamma^{-1}(t)
    \end{equation}
    from $ghg^{-1}$ to $h$, which is entirely contained in $\pi_1(G)$. But $\pi_1(G)$
    is discrete, so this
    implies $\sigma \equiv \E$, i.e. $h \in \Center(\widetilde{G})$.

    Let now
    $g \in \pi_1(G) \subseteq \Center(\widetilde{G})$, then we find a Lie algebra
    element $\xi
    \in \liealg{g}$ with $g = \exp_{\widetilde{G}} \xi$ and $\spec(\ad_\xi)
    \subseteq 2\pi \I \field{Z}$ by
    \cite[Theorem~14.2.8]{hilgert.neeb:2012a}. Consequently,
     \begin{equation*}
         \Ad_{\widetilde{\eta}(g)}
         \chi
         =
         \Ad_{\widetilde{\eta}(\exp_{\widetilde{G}}\xi)}
         \chi
         \overset{\eqref{eq:UniversalComplexificationProof0}}{=}
         \Ad_{\exp_{\widetilde{G}_\C}(\iota(\xi))}
         \chi
         =
         \exp_{\End(\liealg{g}_\C)}
         \bigl(
            \ad_{\iota(\xi)}
         \bigr)
         \chi
         =
         \chi
     \end{equation*}
    holds for all $\chi \in \liealg{g}_\C$. In the final step we have used that the
    only nilpotent matrix exponentiating to the identity is the zero matrix.
    Summarizing, we have shown
    \begin{equation}
        \widetilde{\eta}
        \bigl(
            \pi_1(G)
        \bigr)
        \subseteq
        \ker
        \bigl(
            \Ad_{\widetilde{G}_\C}
        \bigr)
        =
        \Center(\widetilde{G}_\C),
    \end{equation}
    where the equality holds by \cite[Theorem~14.2.8]{hilgert.neeb:2012a}. We denote
    the closure of the
    complex group span of some subset $S$ of a complex Lie group by
    \gls{GroupSpanComplex}, which constitutes a normal subgroup whenever
    $S$ was by continuity of the group operations. By
    \cite[Theorem~11.1.5]{hilgert.neeb:2012a} we can quotient by such a group to
    obtain a complex Lie group. We define
    \begin{equation}
        \label{eq:ComplexificationFromUniversalCovering}
        G_\C
        \coloneqq
        \widetilde{G}_\C
        \Big/
        \big<
            \widetilde{\eta}(\pi_1(G))
        \big>_{\C}
        \quad \textrm{and} \quad
        \eta
        \colon
        G
        \longrightarrow
        G_\C, \quad
        \eta
        \bigl(
            p(g)
        \bigr)
        \coloneqq
        \pi
        \bigl(
            \widetilde{\eta}(g)
        \bigr)
    \end{equation}
    with the quotient projection $\pi \colon \tilde{G}_\C \longrightarrow G_\C$. In
    terms of a commutative diagram, the situation is
    \begin{equation}
        \label{eq:ComplexificationFromUniversalCoveringDiagram}
        \begin{tikzcd}[column sep = huge]
            \widetilde{G}
            \arrow[d, "p", swap]
            \arrow[r, "\tilde{\eta}"]
            &\tilde{G}_\C
            \arrow[d, "\pi"] \\
            G
            \arrow[r, "\eta"]
            &G_\C
        \end{tikzcd}.
    \end{equation}
    We prove that $(G_\C, \eta)$ is the universal complexification of $G$, and first
    check that $\eta$ is smooth. By definition, we have $\eta \circ p = \pi \circ
    \widetilde{\eta}$. As $\pi$ and $\widetilde{\eta}$ are smooth, the same is true
    for $\eta \circ p$. This in turn implies smoothness of $\eta$, as $p$ is a local
    diffeomorphism. Existence of local inverses of $p$, i.e. of smooth sections, also
    yields that $\eta$ is a
    group morphism.

    It remains to check the universal property. To this end, let
    $\Phi \colon G \longrightarrow H$ be a Lie group morphism into a
    complex Lie group $H$ with tangent map $\phi \coloneqq T_{\E_G} \Phi
    \colon \liealg{g} \longrightarrow \liealg{h}$. This induces a
    $\C$-linear Lie algebra morphism $\phi_\C \colon
    \liealg{g}_\C \longrightarrow \liealg{h}$ by universal complexification of vector spaces. This puts us into the position to once again apply
    Theorem~\ref{thm:Lie2} to integrate $\phi_\C$, i.e. we
    get the commutative diagram
    \begin{equation*}
        \begin{tikzcd}[column sep = huge]
        \widetilde{G}_\C
        \arrow[d, "\exp_{\widetilde{G}_\C}", leftarrow, swap]
        \arrow[r, dashed, "\exists! \; \widetilde{\Phi}_\C"]
        &H
        \arrow[d, "\exp_H", leftarrow] \\
        \liealg{g}_\C
        \arrow[r, "\phi_\C"]
        &\liealg{h}
        \end{tikzcd}.
    \end{equation*}
    Notice now that
    \begin{equation}
        \label{eq:ComplexificationProof0}
        T_{\E_{\widetilde{G}}}
        \bigl(
            \Phi \circ p
        \bigr)
        =
        \phi
        \circ
        \id_{\liealg{g}}
        =
        \phi
        =
        \phi_\C
        \circ
        \iota
        =
        T_{\E_{\widetilde{G}}}
        \bigl(
            \widetilde{\Phi}_\C
            \circ
            \widetilde{\eta}
        \bigr),
        \tag{$\star$}
    \end{equation}
    which implies $\Phi \circ p = \widetilde{\Phi}_\C \circ \widetilde{\eta}$
    by connectedness of $\widetilde{G}$, see
    \cite[Cor.~9.2.12]{hilgert.neeb:2012a}. In particular, assuming $g \in
    \pi_1(G)$, this yields
    \begin{equation}
        0
        =
        (\widetilde{\Phi}_\C
        \circ
        \widetilde{\eta})(g),
        \qquad \textrm{i.e.} \quad
        \widetilde{\eta}
        \bigl(
            \pi_1(G)
        \bigr)
        \subseteq
        \ker
        \widetilde{\Phi}_\C.
    \end{equation}
    Consequently, $\widetilde{\Phi}_\C$ descends to the quotient $G_\C$, inducing a
    holomorphic group morphism~$\Phi_\C \colon G_\C \longrightarrow H$. We moreover
    compute
    \begin{align*}
        T_{\E_G}
        (\Phi_\C \circ \eta)
        =
        T_{\E_{\widetilde{G}}}
        (\Phi_\C \circ \eta \circ p)
        =
        T_{\E_{\widetilde{G}}}
        (\Phi_\C \circ \pi \circ \widetilde{\eta})
        =
        T_{\E_{\widetilde{G}}}
        (\widetilde{\Phi}_\C \circ \widetilde{\eta})
        =
        T_{\E_{\widetilde{G}}}
        (\Phi \circ p)
        =
        T_{\E_{G}} \Phi,
    \end{align*}
    where we have used $T_{\E_{\widetilde{G}}}p = \id_{\liealg{g}}$, the definition of
    $\eta$, the characteristic property of the quotient and
    \eqref{eq:ComplexificationProof0}. By connectedness, $\Phi_\C \circ \eta = \Phi$
    follows, as required in
    Definition~\ref{def:UniversalComplexification}. As complex multiples of any basis of
    $\liealg{g}$ generate $\LieAlg(G_\C)$, the $\C$-linear tangent mapping~$T_\E
    \Phi_\C$ is uniquely determined by $\phi$. By connectedness this also determines
    $\Phi_\C$ uniquely. This completes the proof for connected $G$.

    Let now $G$ be disconnected and write $G_0$ for the connected
    component of the group unit. As $G_0$ is connected, it possesses a
    universal complexification $(G_{0, \C}, \eta_0)$ by what we have already
    established. For fixed $g \in G$ this gives the commutative diagram
    \begin{equation*}
        \begin{tikzcd}[column sep = huge]
        G_0
        \arrow[d, "\eta_0", rightarrow]
        \arrow[r, "\Conj_g \at{G_0}"]
        &G_0
        \arrow[d, "\eta_0"] \\
        G_{0, \C}
        \arrow[r, dashed, "\Conj_{g, \C}"]
        &G_{0, \C}
        \end{tikzcd},
    \end{equation*}
    by Lemma~\ref{lem:UniversalComplexificationMorphisms} applied to the restriction
    $\Conj_g \at{G_0}$ of the conjugation from \eqref{eq:Conjugation}. The uniqueness
    gives us the explicit formula
    \begin{equation*}
        \Conj_{g, \C}(h)
        =
        \eta_0(g) h \eta_0(g^{-1}),
        \qquad
        h \in G_{0,\C},
    \end{equation*}
    as this definition makes \eqref{eq:UniversalComplexificationMorphisms} commute due to $\eta_0$ being
    a group morphism. In particular, we may view this as a
    group morphism
    \begin{equation*}
        \label{eq:UniversalComplexificationProof1}
        c
        \colon
        G_0
        \longrightarrow
        \Aut(G_{0, \C}), \quad
        c_g(h)
        \coloneqq
        \Conj_{g, \C}(h)
        =
        \eta_0(g) h \eta_0(g^{-1}).
        \tag{$*$}
    \end{equation*}
    Consider now the outer semidirect product $S \coloneqq
    G_{0, \C} \ltimes_c G$ together with the subset
    \begin{equation*}
        \label{eq:UniversalComplexificationProof2}
        K
        \coloneqq
        \bigl\{
            \bigl( \eta_0(g), g^{-1} \bigr)
            \in
            S
            \colon
            g \in G_0
        \bigr\}.
        \tag{$\dagger$}
    \end{equation*}
    We claim that $K$ is a normal subgroup of $S$. To show this, let
    first $g, h \in G_0$ and compute
    \begin{align*}
        \bigl( \eta_0(g), g^{-1} \bigr)
        \cdot
        \bigl( \eta_0(h), h^{-1} \bigr)
        &=
        \bigl(
            \eta_0(g)c_{g^{-1}}(\eta_0(h)),
            g^{-1} h^{-1}
        \bigr) \\
        &=
        \bigl(
            \eta_0(g) \eta_0(g^{-1}) \eta_0(h) \eta_0(g),
            (hg)^{-1}
        \bigr) \\
        &=
        \bigl(
            \eta_0(hg),
            (hg)^{-1}
        \bigr)
        \in
        K,
    \end{align*}
    where we have used \eqref{eq:UniversalComplexificationProof1} and
    that $G_0$ is a subgroup. Let $g \in G_0$ and $(z, h) \in
    S$. Then
    \begin{align*}
        (z, h)
        \cdot
        \bigl( \eta_0(g), g^{-1} \bigr)
        \cdot
        \bigl(
            c_{h^{-1}}(z^{-1}),
            h^{-1}
        \bigr)
        &=
        (z, h)
        \cdot
        \bigl( \eta_0(g), g^{-1} \bigr)
        \cdot
        \bigl(
            \eta_0(h^{-1}) z^{-1} \eta_0(h),
            h^{-1}
        \bigr) \\
        &=
        (z, h)
        \cdot
        \bigl(
            \eta_0(g)
            c_{g^{-1}}
            (\eta_0(h^{-1}) z^{-1} \eta_0(h)),
            g^{-1} h^{-1}
        \bigr) \\
        &=
        (z, h)
        \cdot
        \bigl(
           \eta_0(h^{-1}) z^{-1} \eta_0(hg),
           g^{-1} h^{-1}
        \bigr) \\
        &=
        \bigl(
            z c_h
            (\eta_0(h^{-1}) z^{-1} \eta_0(hg)),
            h g^{-1} h^{-1}
        \bigr) \\
        &=
        \bigl(
            z \eta_0(h)
            (\eta_0(h^{-1}) z^{-1} \eta_0(hg))
            \eta_0(h^{-1}),
            h g^{-1} h^{-1}
        \bigr) \\
        &=
        \bigl(
            \eta_0(hgh^{-1}),
            h g^{-1} h^{-1}
        \bigr)
        \in
        K,
    \end{align*}
    which
    proves the normality. Note furthermore that $K$ is closed by continuity of group inversion and $\eta_0$. We are thus in a position to define
    \begin{equation}
        \label{eq:ComplexificationNonConnected}
        G_\C
        \coloneqq
        S \big/ K
        \quad \textrm{and} \quad
        \eta
        \colon
        G
        \longrightarrow
        G_\C, \quad
        \eta(g)
        \coloneqq
        \big[
            (\E_{G_{0, \C}}, g)
        \big]
    \end{equation}
    as a real Lie group and claim that $(G_\C, \eta)$ is the universal complexification of $G$. To this end, we compute
    \begin{equation*}
        \eta(g) \eta(h)
        =
        \big[
            (\E_{G_{0, \C}}, g)
            \cdot
            (\E_{G_{0, \C}}, h)
        \big]
        =
        \big[
            (
                \E_{G_{0, \C}} c_g(\E_{G_{0, \C}}),
                gh
            )
        \big]
        =
        [(\E_{G_{0, \C}}, gh)]
        =
        \eta(gh)
    \end{equation*}
    for $g, h \in G$, i.e. $\eta$ is a group morphism. As the map $G \ni g \mapsto
    (\E_{G_{0, \C}}, g) \in S$ is continuous and induces
    $\eta$, the smoothness of $\eta$ follows from the quotient projection being a surjective submersion.

    Next, we construct a holomorphic atlas for $G_\C$. To this end, let $(U_\alpha, z_\alpha)_{\alpha \in J}$ be a holomorphic atlas of $G_{0,\C}$. Note that setting $h \coloneqq g^{-1}$ in \eqref{eq:UniversalComplexificationProof3} gives $(z,g) \sim (z \eta_0(g^{-1}), \E_G)$, i.e. every equivalence class $[(z,g)]$ contains a unique representative of the form $(z',e_G)$. Consequently, the sets $V_\alpha \coloneqq [U_\alpha \times \{\E_G\}]$ for $\alpha \in J$ cover $G_\C$. Moreover, the mappings
    \begin{equation*}
        w_\alpha
        \colon
        V_\alpha \longrightarrow z_\alpha(U_\alpha), \quad
        w_\alpha
        \bigl(
            [z,\E_G]
        \bigr)
        \coloneqq
        z_\alpha(z)
    \end{equation*}
    are well defined biholomorphisms. This provides the desired holomorphic atlas for
    $G_\C$.

    Let now $\Phi \colon G \longrightarrow H$ be a Lie
    group morphism into a complex Lie group $H$. Its restriction~$\Phi_0 \coloneqq
    \Phi \at{G_0}$ induces a unique holomorphic group morphism
    $\Phi_{0, \C} \colon G_{0, \C} \longrightarrow H$
    satisfying~$\Phi = \Phi_{0, \C} \circ \eta_0$ by the universal property of $(G_0,
    \eta_0)$. We consider
    \begin{equation*}
        \Phi_\C
        \colon
        G_\C
        \longrightarrow
        H, \quad
        \Phi_\C
        \bigl( [(z, g)] \bigr)
        \coloneqq
        \Phi_{0, \C}(z)
        \Phi(g),
    \end{equation*}
    which is well defined: varying $h \in G_0$ we get all other representatives of $(z,g) \in S$ as
    \begin{equation*}
        \label{eq:UniversalComplexificationProof3}
        (z, g)
        \cdot
        (\eta_0(h), h^{-1})
        =
        \bigl(
            z c_g(\eta_0(h)),
            g h^{-1}
        \bigr)
        \overset{\eqref{eq:UniversalComplexificationProof1}}{=}
        \bigl(
            z \eta_0(ghg^{-1}),
            g h^{-1}
        \bigr).
        \tag{$\ddag$}
    \end{equation*}
    Now, as $\Phi = \Phi_{0, \C} \circ \eta_0$, we get
    \begin{align}
       \Phi_{0, \C}
       \bigl(z \eta_0(ghg^{-1})\bigr)
       \cdot
       \Phi(gh^{-1})
       &=
       \Phi_{0, \C}(z)
       \cdot
       \Phi(ghg^{-1})
       \cdot
       \Phi(gh^{-1}) \\
       &=
       \Phi_{0, \C}(z)
       \cdot
       \Phi(g) \\
       &=
       \Phi_\C(z,g).
    \end{align}
    Thus the value of $\Phi_\C$ is indeed independent of
    the chosen representative. For $z, z' \in G_{0, \C}$ and
    $g, g' \in G$ we compute
    \begin{align*}
        \Phi_\C
        \bigl(
            [(z, g)]
            \cdot
            [(z', g')]
        \bigr)
        &=
        \Phi_\C
        \bigl(
            [(
                z \eta_0(g) z' \eta_0(g^{-1}),
                gg'
            )]
        \bigr) \\
        &=
        \Phi_{0, \C}
        \bigl(
            z \eta_0(g) z' \eta_0(g^{-1})
        \bigr)
        \Phi(gg') \\
        &=
        \Phi_{0, \C}(z)
        \cdot
        \Phi(g)
        \cdot
        \Phi_{0, \C}(z')
        \cdot
        \Phi(g^{-1})
        \cdot
        \Phi(gg') \\
        &=
        \Phi_\C
        \bigl(
            [(z, g)]
        \bigr)
        \cdot
        \Phi_\C
        \bigl(
            [(z', g')]
        \bigr),
    \end{align*}
    showing $\Phi_\C$ is a group morphism. Moreover, the description of $\Phi_\C$ in the atlas $(V_\alpha, w_\alpha)_{\alpha \in J}$ coincides with the description of $\Phi_{0,\C}$ in the atlas $(U_\alpha, z_\alpha)_{\alpha \in J}$, i.e.
    \begin{equation*}
        \Phi_\C
        \circ
        w_\alpha^{-1}
        \at[\Big]{w_\alpha([z,\E_G])}
        =
        \Phi_{0,\C}
        \circ
        z_\alpha^{-1}
        \at[\Big]{z_\alpha(z)}
    \end{equation*}
    for all $\alpha \in J$ and $z \in U_\alpha$. This implies holomorphicity of $\Phi_\C$ at once.
    By definition it is also true that
    \begin{equation*}
        \Phi_\C
        \bigl(
            \eta(g)
        \bigr)
        =
        \Phi_{0,\C}
        (\E_{G_{0, \C}})
        \cdot
        \Phi(g)
        =
        \Phi(g).
    \end{equation*}
    Thus $\Phi_{0, \C}$ makes the diagram
    \eqref{eq:UniversalComplexification} commute. Conversely, the morphism
    property and the preceding condition already determine
    $\Phi_\C$ uniquely. This completes the proof.
\end{proof}

Taking a closer look at our construction and keeping our notation, we obtain the
following additional well known properties of the universal complexification, all of
which can be found in \cite[Sec.~15.1 \& Exercises]{hilgert.neeb:2012a}.
\begin{corollary}
    \label{cor:UniversalComplexificationProperties}%
    \index{Complexification!Properties}
    Let $G$ be a Lie group with universal complexification
    $(G_\C, \eta)$.
    \begin{corollarylist}
        \item \label{item:ComplexificationDimension}
        We have $\dim_\C G_\C \le \dim_\R G$.
        \item \label{item:UniversalComplexificationFaithful}%
        If $G$ is linear, then $\eta$ is injective.
    \end{corollarylist}
    If $G$ is connected, we moreover have the following:
    \begin{corollarylist}
        \setcounter{enumi}{2}
        \item \label{item:ComplexificationKernelCentral}
        The subgroup $\eta(\pi_1(G))$ is central.
        \item \label{item:UniversalComplexificationSimplyConnected1}
        The pair $(\widetilde{G}_\C, \widetilde{\eta})$ constitutes the
          universal complexification of the universal covering group $\widetilde{G}$
          of $G$.
        \item \label{item:UniversalComplexificationSimplyConnected2}%
          If $G$ is simply connected, the same is true for its
          universal complexification $G_\C$. Moreover,
          $T_\E \eta \colon \liealg{g} \longrightarrow \liealg{g}_\C$
          is injective, $\LieAlg(G_\C) \cong \liealg{g}_\C$ and $\dim_\C G_\C = \dim_\R G$.
        \item \label{item:ComplexificationSemisimple}
        If $G$ is semisimple, then $\dim_\C G_\C = \dim_\R G$.
    \end{corollarylist}
\end{corollary}
\begin{proof}
    To see the first part, recall first that $\dim_\C(\widetilde{G}_\C) =
    \dim_\R(\widetilde{G}) = \dim_\R(G)$. Thus, in the connected case, the
    quotient $G_\C$ of $\widetilde{G}$ as defined in
    \eqref{eq:ComplexificationFromUniversalCovering} fulfils
    \begin{equation}
        dim_\C
        (G_\C)
        \le
        \dim_\R(G).
    \end{equation}
    For the disconnected case we first note
    \begin{align}
        \dim_\R
        \bigl(
            G_{0,\C} \rtimes G
        \bigr)
        &=
        \dim_\R
        (G_{0,\C})
        \cdot
        \dim_\R(G) \\
        &=
        2
        \cdot
        \dim_\C(G_{0,\C})
        \cdot
        \dim_\R(G) \\
        &\le
        2
        \dim_\R(G)^2,
    \end{align}
    as semidirect products and direct products coincide as manifolds and the unit
    component~$G_0$ fulfils $\dim_\R(G_0) = \dim_\R(G)$. Taking another look at
    \eqref{eq:UniversalComplexificationProof2}, we see that
    \begin{equation}
        \dim_\R(K)
        =
        \dim_\R(G).
    \end{equation}
    Putting everything together and using \eqref{eq:ComplexificationNonConnected} thus
    yields $\dim_\C(G_\C) \le \dim_\R(G)$ also in the disconnected case.

    Given a faithful representation $\rho \colon G \longrightarrow
    \operatorname{GL}_n(\C)$, we have the commutative diagram
    of group morphisms
    \begin{equation*}
        \begin{tikzcd}[column sep = huge]
            G
            \arrow[d, "\eta"]
            \arrow[r, "\rho"]
            &\operatorname{GL}_n(\C) \\
            G_\C
            \ar[ru, dashed, swap, "\exists! \; \rho_\C"]
        \end{tikzcd}.
    \end{equation*}
    As $\rho$ is injective, the same is true for $\eta$ regardless of the behaviour of $\rho_\C$.

    We have already shown \ref{item:ComplexificationKernelCentral} during the proof of Theorem~\ref{thm:UniversalComplexification}. Going through the construction with $G = \widetilde{G}$, we note $\pi_1(G)$ and thus $\widetilde{\eta}(\pi_1(G))$ are trivial, i.e. we quotient out nothing in \eqref{eq:ComplexificationFromUniversalCovering} and get $G_\C = \widetilde{G}_\C$ as well as $T_\E \eta = \iota$, proving \ref{item:UniversalComplexificationSimplyConnected1} and \ref{item:UniversalComplexificationSimplyConnected2}.

    For the final statement recall that a semisimple group $G$ has discrete center
    $\Center(G)$. By~\ref{item:ComplexificationKernelCentral}, this implies that
    $\tilde{\eta}(\pi_1(G)) \subseteq \tilde{G}_\C$ is discrete and as such both
    complex and closed, see again Example~\ref{ex:DiscreteGroups}. Thus the
    quotient in
    \eqref{eq:ComplexificationFromUniversalCovering} does not lose dimensions and from
    \ref{item:ComplexificationDimension} we infer
    \ref{item:ComplexificationSemisimple}.
\end{proof}

The observation from \ref{item:ComplexificationSemisimple} is due to \cite{hochschild:1966a}.
Our next goal is to describe holomorphicity on a universal complexification in a
left invariant manner. To this end, we once again choose a basis $(\basis{e}_1,
\ldots, \basis{e}_n)$ of its Lie algebra $\liealg{g}$. This induces the basis
$(\basis{e}_1, \ldots, \basis{e}_n, \I \basis{e}_1, \ldots, \I \basis{e}_n)$ of the
complexification $\liealg{g}_\C$ and the canonical complex structure $J_\canonical$
becomes multiplication with the imaginary unit $\I$. Alternatively, one may consider
the Wirtinger\footnote{Wilhelm Wirtinger (1865-1945) was an Austrian function
theorist, who first proposed generalizing eigenvalues to spectra. His
doctoral students include Wilhelm Blaschke, Kurt Gödel and Leopold Vietoris.} linear
combinations
\begin{equation}
    \label{eq:Wirtinger}
    \basis{f}_k
    \coloneqq
    \frac{\basis{e}_k - \I \basis{e}_k}{2}
    \quad \textrm{and} \quad
    \cc{\basis{f}}_k
    \coloneqq
    \frac{\basis{e}_k + \I \basis{e}_k}{2}
    \qquad
    \textrm{for }
    k=1, \ldots, n,
\end{equation}
which constitute another basis of $\liealg{g}_\C$. Note that the left translations of
$f_k$ span the holomorphic tangent bundle, whereas the ones of $\cc{f}_k$ generate the
antiholomorphic one.
\begin{proposition}[Left invariant Cauchy-Riemann equations]
    \index{Cauchy-Riemann equations}
    \index{Left invariant complex analysis!Cauchy-Riemann}
    \label{prop:CauchyRiemann}
    Let $G$ be a Lie group such that $\LieAlg(G_\C) = \LieAlg(G)_\C$ and let $\phi \colon G_\C \longrightarrow \C$ be a function on the universal complexification $G_\C$ of $G$. Then the following are equivalent:
    \begin{propositionlist}
        \item The function $\phi$ is holomorphic.
        \item \label{item:CauchyRiemann}
        We have $\Lie(\I \basis{e}_k) \phi = \I \cdot \Lie(\basis{e}_k) \phi$ for
        $k=1, \ldots, n$.
        \item \label{item:CauchyRiemannWirtinger}
        We have $\Lie(\cc{\basis{f}}_k) \phi = 0$ for $k=1, \ldots, n$.
    \end{propositionlist}
    If one -- and thus all -- of the above conditions are fulfilled, then
    \begin{equation}
        \Lie(\basis{e}_k)
        \phi
        \in
        \Holomorphic(G)
        \qquad
        \textrm{for }
        k=1, \ldots, n.
    \end{equation}
\end{proposition}
\begin{proof}
    The equivalence of \ref{item:CauchyRiemann} and \ref{item:CauchyRiemannWirtinger} is clear by \eqref{eq:Wirtinger}. As the Lie exponential is the flow of the left invariant vector field, we have
    \begin{equation}
        \Lie_{X_{\basis{e}_k}}
        \phi
        (g)
        =
        \frac{\D}{\D t}
        \phi
        \bigl(
            g \exp_{G_\C}(t \basis{e}_k)
        \bigr)
        \at[\Big]{t=0}
    \end{equation}
    for all $g \in G$. Now, \ref{item:CauchyRiemann} and \ref{item:CauchyRiemannWirtinger} are exactly the Cauchy-Riemann equations for
    \begin{equation}
        \phi
        \circ
        \exp_{G_\C}
        \colon
        \LieAlg(G)_\C
        \longrightarrow
        \C
    \end{equation} and as we may check holomorphicity in the exponential atlas, this establishes the remaining equivalences. The additional statement is now clear.
\end{proof}

It turns out that the almost complex structure of a complex Lie group arising from universal complexification $\eta \colon G \longrightarrow G_\C$ always comes with a complex conjugation.
\begin{proposition}[Complex conjugation,
{\cite[Thm.~15.1.4~\textit{iv.)}]{hilgert.neeb:2012a}}]
    \index{Complex!Conjugation}
    \label{prop:ComplexConjugation}
    Let $G$ be a Lie group with universal complexification $(G_\C, \eta)$. Then there
    exists a unique antiholomorphic group morphism
    \begin{equation}
        \gls{ComplexConjugationGroup}
        \colon G_\C \longrightarrow G_\C
    \end{equation}
    such that $\sigma^2 = \id_{G_\C}$ and fixed point group
    \begin{equation}
        \label{eq:ComplexConjugationFixed}
        \bigl\{
            g \in G_{\C,0}
            \colon
            \sigma(g)
            =
            g
        \bigr\}
        =
        \eta(G_0),
    \end{equation}
    where $G_{\C,0} \subseteq G_\C$ denotes the connected component of the group unit.
\end{proposition}
\index{Complex!Conjugated group}
We call $\sigma \colon G_\C \longrightarrow G_\C$ the \emph{complex conjugation} of $G_\C$ and refer to \gls{ComplexConjugatedGroup} as the \emph{complex conjugated group} to $G$. As usual, we indicate the corresponding group with an additional subscript when several complex conjugations are involved. For $G = \R^n$ this reproduces the usual complex conjugation.
\begin{proof}[of Proposition~\ref{prop:ComplexConjugation}]
    We define a complex Lie algebra $\cc{\liealg{g}}_\C$ by endowing $\liealg{g}_\C$ with the same addition and Lie bracket, but a new multiplication by scalars
    \begin{equation*}
        \label{eq:ComplexConjugationProof1}
        \gls{ScalarConjugated}
        \colon
        \C
        \times
        \cc{\liealg{g}}_\C
        \longrightarrow
        \cc{\liealg{g}}_\C, \quad
        z
        \odot
        \cc{\xi}
        =
        \cc{ \cc{z} \xi},
        \tag{$\flat$}
    \end{equation*}
    where we distinguish the elements of $\cc{\liealg{g}}_\C$ by writing
    $\cc{\xi}$.\footnote{I would like to thank Stefan Waldmann for making me aware of
    this extraordinarily useful notational convention.} As $\cc{\liealg{g}}_\C$ and
    $\liealg{g}_\C$ are
    isomorphic as real vector spaces, we may replace the latter with the former in our
    construction in Theorem~\ref{thm:UniversalComplexification}. This yields another
    pair $(\cc{G}_\C, \cc{\eta})$ of complex Lie group and Lie group morphism
    $\cc{\eta} \colon G \longrightarrow \cc{G}_\C$. As real Lie groups, $\cc{G}_\C$
    and $G_\C$ coincide, but $\cc{G}_\C$ does not fulfil the universal property of the
    universal complexification. The universal property of $G_\C$ yields the diagram
    \begin{equation*}
        \begin{tikzcd}[column sep = huge]
            G
            \arrow[d, "\eta"]
            \arrow[r, "\overline{\eta}"]
            &\overline{G}_\C \\
            G_\C
            \ar[ru, dashed, "\exists! \; \overline{\sigma}"]
        \end{tikzcd}
    \end{equation*}
    with a unique holomorphic group morphism $\cc{\sigma}$. Let now
    \begin{equation*}
        \sigma
        \colon
        G_\C \longrightarrow G_\C, \quad
        \sigma
        \coloneqq
        \cc{\id}
        \circ
        \overline{\sigma},
    \end{equation*}
    where $\cc{\id} \colon \overline{G}_\C \longrightarrow G_\C$ is given by
    $\cc{\id}(g) \coloneqq g$. Note that $\cc{\id} \circ \cc{\eta} = \eta$. By
    \eqref{eq:ComplexConjugationProof1}, the mapping~$\cc{\id}$ is antiholomorphic and
    thus the same is true for $\sigma$ as a composition of a holomorphic with an
    anthiholomorphic map. By construction, we moreover have
    \begin{equation}
        \label{eq:ComplexConjugationProof2}
        \sigma
        \circ
        \eta
        =
        \cc{\id}
        \circ
        \overline{\sigma}
        \circ
        \eta
        =
        \cc{\id}
        \circ
        \cc{\eta}
        =
        \eta.
        \tag{$\sharp$}
    \end{equation}
    Thus $\sigma^2
    \colon G_\C \longrightarrow G_\C$ is holomorphic and
    fulfils $\sigma^2 \circ \eta = \eta$. The
    universal property asserts that the unique holomorphic group
    morphism with this feature is given by $\id_{G_\C}$, proving that $\sigma$ is an
    involution.
    Moreover, \eqref{eq:ComplexConjugationProof2} proves the inclusion ``$\supseteq$''
    in \eqref{eq:ComplexConjugationFixed}. Let conversely be $g \in G_{\C,0}$ be such
    that $\sigma(g) = g$. As $g$ is in the unit component of the group unit, we find
    $\xi_1, \ldots, \xi_n \in \liealg{g}_\C$ such that $g = \exp(\xi_1) \cdots
    \exp(\xi_n)$. Now, $\sigma(\exp(\xi)) = \exp(\xi)$ for some Lie algebra element
    $\xi \in \liealg{g}_\C$ holds if and only if $\xi \in \liealg{g}$, as
    \begin{equation}
        \sigma
        \circ
        \exp
        =
        \exp
        \circ
        T_\E \sigma
        =
        \exp
        \circ
        \operatorname{cc}
    \end{equation}
    with the complex conjugation \gls{ComplexConjugation}$(\xi) = \cc{\xi}$.
    Consequently, we get $\xi_1, \ldots, \xi_n \in \liealg{g}$, which implies $g \in
    \eta(G_0)$.
\end{proof}

\begin{remark}[Fixed point groups]
    \index{Fixed point groups}
    \index{Cartan, Élie}
    \index{Borel, Armand}
    Let $G$ be simply connected and semisimple Lie group. By
    Corollary~\ref{cor:UniversalComplexificationProperties},
    \ref{item:UniversalComplexificationSimplyConnected2}, $G_\C$ is then also simply
    connected and Cartan's\footnote{Élie Cartan (1869-1951) was a French differential
    geometer, who shaped much of our modern understanding of analysis on manifolds
    and, in particular, Lie groups. His oldest son Henri Cartan (1904-2008) also left
    significant marks on mathematics, especially in the field of algebraic topology,
    and was a founding member of Bourbaki.}
    semisimplicity criterion
    \cite[Thm.~5.5.9]{hilgert.neeb:2012a} yields its semisimplicity. Under these
    assumptions, the fixed point group
    \begin{equation}
        \bigl\{
            g \in G_\C
            \colon
            \sigma(g)
            =
            g
        \bigr\}
    \end{equation}
    is always connected by a theorem of Borel\footnote{Armand Borel (1923-2003) was a
    Swiss mathematician specialized in algebraic topology and Lie groups. Not to be
    confused with the French mathematician Èmile Borel (1871-1956), whose Lemma on
    Taylor coefficients we have come across within
    Chapter~\ref{ch:StrictDeformationQuantization}.}
    \cite[Thm.~3.4.]{borel:1961a}. Thus, the
    restriction to $G_{\C,0}$ that facilitated our proof of
    \eqref{eq:ComplexConjugationFixed} was not necessary in this case.
\end{remark}

Along the way we have shown an antiholomorphic version of the universal property.
\begin{corollary}
    For every group morphism $\Phi \colon G \longrightarrow
    H$ into a complex Lie group $H$ there is a unique
    antiholomorphic group morphism $\cc{\Phi}_\C \colon
    G_\C \longrightarrow H$ such that
    \begin{equation}
        \label{eq:AntiholomorphicExtension}
        \cc{\Phi}_\C \circ \eta
        =
        \Phi.
    \end{equation}
\end{corollary}
\begin{proof}
    Setting $\cc{\Phi}_\C \coloneqq \sigma \circ \Phi_\C$ with $\Phi_\C$ induced by
    \eqref{eq:UniversalComplexification} does the job. Conversely, every
    antiholomorphic solution $\cc{\Phi}_\C$ of~\eqref{eq:AntiholomorphicExtension}
    induces a holomorphic solution $\Phi_\C \coloneqq \cc{\Phi}_\C \circ \sigma$
    of~\eqref{eq:HolomorphicExtension}, proving the uniqueness of
    $\cc{\Phi}_\C$.
\end{proof}

By continuity of $\sigma$, the following useful observation is now immediate.
\begin{corollary}
    \label{cor:ImageOfEtaSubmanifold}
    Let $G$ be a Lie group with universal complexification $(G_\C, \eta)$.
    \begin{corollarylist}
        \item \label{item:ImageOfEtaSubmanifold}
        The image $\eta(G) \subseteq G_\C$ is a closed real subgroup and submanifold.
        \item \label{item:ImageOfEtaComplexification}
        The embedding $\iota \colon \eta(G) \hookrightarrow G_\C$ constitutes a universal complexification.
    \end{corollarylist}
\end{corollary}
\begin{proof}
    By Proposition~\ref{prop:ComplexConjugation}, there exists a complex conjugation
    $\sigma \colon G_\C \longrightarrow G_\C$. Consider the smooth auxiliary mapping
    \begin{equation}
        \Phi
        \colon
        G_\C
        \longrightarrow
        G_\C, \quad
        \Phi(g)
        \coloneqq
        \sigma(g)g^{-1}.
    \end{equation}
    By \eqref{eq:ComplexConjugationFixed}, $\eta(G)$ is the set of fixed points for
    $\sigma$ intersected with the connected component of the group unit, both of which
    are closed. Consequently, $\eta(G)$ is closed, as well. As images of group
    morphisms are subgroups, this implies that $\eta(G)$ is a real Lie group itself
    and in particular a submanifold of $G_\C$ with exponential charts as submanifold
    charts by the closed subgroup Theorem \cite[Thm.~9.3.7]{hilgert.neeb:2012a}. This
    completes the proof of the first part. We proceed with
    \ref{item:ImageOfEtaComplexification}. To this end, let $\Phi \colon \eta(G)
    \longrightarrow H$ be a group morphism into a complex Lie group $H$. Then
    \begin{equation}
        \Psi
        \coloneqq
        \Phi
        \circ
        \eta
        \colon
        G
        \longrightarrow
        H
    \end{equation}
    is a Lie group morphism from $G$ into a complex Lie group. By the universal
    property of $(G_\C, \iota \circ \eta)$,\footnote{In this proof, $\eta$ was
    co-restricted to its image $\eta(G)$, so the usual universal complexification
    map for~$G$ is given by $\iota \circ \eta \colon G \longrightarrow G_\C$
    instead.}
    there thus exists a unique holomorphic group morphism $\Psi_\C \colon G_\C
    \longrightarrow H$ with
    \begin{equation}
        \Psi_\C
        \circ
        \iota
        \circ
        \eta
        =
        \Psi
        =
        \Phi
        \circ
        \eta.
    \end{equation}
    By surjectivity of $\eta \colon G \longrightarrow \eta(G)$, this is equivalent to $\Psi_\C \circ \iota = \Phi$ and $\Psi_\C \colon G_\C \longrightarrow H$ is the only holomorphic group morphism with this property.
\end{proof}

Our next goal is to prove that
\begin{equation}
    \LieAlg(G_\C)
    =
    \LieAlg
    \bigl(
        \eta(G)
    \bigr)_\C,
\end{equation}
which matches with the intuition behind the existence of unique holomorphic extensions
in \ref{item:ImageOfEtaComplexification}. Before we begin working towards the
realization of our goal in the form of
Lemma~\ref{lem:ComplexificationEtaLocallyInjective}, we take a short detour. The
complex conjugation allows us to characterize those holomorphic group morphisms
between universal complexifications, which arise from group morphisms between the
underlying real groups.
\begin{proposition}
    \label{prop:InducedMorphismReal}
    Let $G, H$ be connected Lie groups and $\Phi \colon G_\C \longrightarrow H_\C$ a group morphism. Then the following are equivalent:
    \begin{propositionlist}
        \item \label{item:MorphismInducedReal}
        We have $\Phi = \Psi_\C$ in the sense of \eqref{eq:UniversalComplexificationMorphisms} with $\Psi \coloneqq \Phi \at{\eta_G(G)} \colon \eta_G(G) \longrightarrow \eta_H(H)$.
%        \item \label{item:MorphismInducedComplex}
%        We have $\Phi = \Psi_\C$ in the sense of the universal property with $\Psi \coloneqq \Phi \at{\eta_G(G)} \colon \eta_G(G) \longrightarrow H_\C$.
        \item \label{item:MorphismCommutesConjugation}
        The group morphism $\Phi$ fulfils
        \begin{equation}
            \label{eq:MorphismCommutesConjugation}
            \Phi \circ \sigma_G
            =
            \sigma_H \circ \Phi.
        \end{equation}
        \item \label{item:MorphismCommutesConjugationTangent}
        The tangent map $\phi \coloneqq T_\E \Phi$ fulfils
        \begin{equation}
            \label{eq:MorphismCommutesConjugationTangent}
            \phi \circ \operatorname{cc}
            =
            \operatorname{cc} \circ \phi.
        \end{equation}
    \end{propositionlist}
    If any -- and thus all -- of the above conditions are fulfilled, then $\Phi$ is
    holomorphic.
\end{proposition}
\begin{proof}
    As for connected group also its universal complexification is connected and thus the equivalence of \ref{item:MorphismCommutesConjugation} and \ref{item:MorphismCommutesConjugationTangent} is immediate. If $G$ is not connected, the latter of course still follows from the former by differentiating. Assuming \ref{item:MorphismInducedReal}, we compute for $g \in G$
    \begin{equation}
        \sigma_H
        \circ
        \Phi
        \circ
        \sigma_G
        \at[\Big]{\eta_G(g)}
        =
        \sigma_H
        \circ
        \Phi
        \at[\Big]{\eta_G(g)}
        =
        \sigma_H
        \at[\Big]
        {\Phi(\eta_G(g))}
        =
        \Phi
        \at[\Big]{\eta_G(g)},
    \end{equation}
    where we have used \eqref{eq:ComplexConjugationFixed} and $\Phi(\eta(G)) \subseteq \eta_H(H)$. By Corollary~\ref{cor:IdentityPrinciple}, \eqref{eq:MorphismCommutesConjugation} follows. Assume conversely \ref{item:MorphismCommutesConjugation} holds. Then we get
    \begin{equation}
        \Phi
        \bigl(
            \eta_G(g)
        \bigr)
        =
        \Phi
        \at[\Big]{\eta_G(g)}
        =
        \Phi
        \circ
        \sigma_G
        \at[\Big]{\eta_G(g)}
        =
        \sigma_H
        \circ
        \Phi
        \at[\Big]{\eta_G(g)}
        =
        \sigma_H
        \at[\Big]{\Phi(\eta_G(g))}
        =
        \sigma
        \bigl(
            \Phi
            \bigl(
            \eta(g)
            \bigr)
        \bigr)
    \end{equation}
    for all $g \in G$. Consequently, $\Phi \at{\eta(g)} \subseteq \eta_H(H)$ by
    connectedness of $H$ and \eqref{eq:ComplexConjugationFixed}. It remains to
    establish that $\Phi$ is holomorphic, then the uniqueness in
    Lemma~\ref{lem:UniversalComplexificationMorphisms}
    proves~\ref{item:MorphismInducedReal}. We do this by checking $\C$-linearity of
    the tangent map $\phi$. By what we have already shown, we know that $\phi$ maps
    $\liealg{g} \coloneqq \LieAlg{\eta_G(G)}$ into $\liealg{h} \coloneqq
    \LieAlg{\eta_H(H)}$. Consequently, using
    \eqref{eq:MorphismCommutesConjugationTangent}, we get
    \begin{equation}
        \phi(\I \xi)
        =
        \phi \circ \operatorname{cc}
        \at[\Big]
        {-\I \xi}
        =
        \operatorname{cc} \circ \phi
        \at[\Big]
        {-\I \xi}
        =
        -
        \cc{\phi(\I \xi)}
    \end{equation}
    for all $\xi \in \liealg{g}$. That is, $\phi$ maps $\I \liealg{g}$ into $\I
    \liealg{h}$ and thus splits into $\phi = \phi_1 \oplus \I \phi_2$ with $\R$-linear
    functions $\phi_1 \colon \liealg{g} \longrightarrow \liealg{h}$ and $\phi_2 \colon
    \liealg{g} \longrightarrow \liealg{h}$. But such mappings are $\C$-linear.
\end{proof}

Note that \eqref{eq:MorphismCommutesConjugation} is a stronger condition than
holomorphicity. There typically exist holomorphic group morphisms between universal
complexifications that do not restrict nicely to $\eta(G)$. Indeed, this already
happens in the simplest case $G = \R$. While a generic complex linear -- and thus
holomorphic -- mapping $\Phi \colon \C \longrightarrow \C$ is given by $\Phi(z) = az$
for some $a \in \C$, only those with $a \in \R$ fulfil
\eqref{eq:MorphismCommutesConjugation}, which in this case coincides
with~\eqref{eq:MorphismCommutesConjugationTangent}. The reason for this is
that \emph{any}
one dimensonal real subspace of~$\C$ complexifies back to~$\C$. That is to
say, real forms are
typically not unique. We will see a more drastic instance of this phenomenon in
Example~\ref{ex:SpecialLinear}, where two non-isomorphic groups have the same
universal complexification. Strictly speaking, this is impossible, since the
morphism~$\eta$ is part of the data for a universal complexification. However,
the
involved complex Lie groups \emph{do} indeed coincide. Keeping these
phenomenona in mind disproves a number of reasonably sounding abstract
statements.

Next, we investigate the fundamental groups of $G_\C$ for connected Lie groups $G$. In Corollary~\ref{cor:UniversalComplexificationProperties}, \ref{item:UniversalComplexificationSimplyConnected2}, we have already seen that $\pi_1(G) = \{\E\}$ implies that $\pi_1(G_\C) = \{\E\}$. It thus remains to cover the multiply connected situation.

Taking another look at the construction of the universal complexification we have
discussed in Theorem~\ref{thm:UniversalComplexification}, one might expect
$\tilde{G}_\C$ to be the universal covering group of~$G_\C$ with the quotient
projection $\pi \colon \tilde{G}_\C \longrightarrow G_\C$ from
\eqref{eq:ComplexificationFromUniversalCovering} being the universal covering map.
This then would imply the desired equalities of fundamental groups at once. Note that
by uniqueness $\pi = (\eta \circ p)_\C$, which is the other natural candidate for a
universal covering projection. As a first ingredient, we use the existence of complex
conjugations to infer that injectivity of $\tilde{\eta}$ guarantees it is actually an
embedding of manifolds, i.e. a diffeomorphism when endowing $\tilde{\eta}(\tilde{G})
\subseteq \tilde{G}_\C$ with the subspace topology.
\index{Embedding}
\begin{lemma}
    \label{lem:Embedding}
    Let $\tilde{G}$ be a connected and simply connected Lie group such that
    \begin{equation}
        \tilde{\eta}
        \colon
        \tilde{G}
        \longrightarrow
        \tilde{G}_\C
    \end{equation}
    is injective. Then $\tilde{\eta}$ is an embedding of manifolds.
\end{lemma}
\begin{proof}
    By
    Corollary~\ref{cor:ImageOfEtaSubmanifold}, $\tilde{\eta}(\tilde{G})$ constitutes a
    closed real subgroup and thus a closed submanifold of $\tilde{G}_\C$.
    Consequently, $\tilde{\eta} \colon \tilde{G} \longrightarrow
    \tilde{\eta}(\tilde{G})$ is a bijective immersion by virtue of
    Corollary~\ref{cor:UniversalComplexificationProperties},
    \ref{item:UniversalComplexificationSimplyConnected2}. Applying the smooth inverse
    function theorem, we see that immersivity and smoothness imply that $\tilde{\eta}$
    is a local diffeomorphism. As smoothness of the inverse may be checked locally,
    this proves that $\tilde{\eta}$ is a diffeomorphism.
\end{proof}

Recall that a subset $D \subseteq X$ of a topological space $X$ is called discrete if the subspace topology the subset $D$ inherits from $X$ is the discrete topology.
\begin{lemma}
    \index{Discrete}
    \label{lem:Discrete}
    Let $\iota \colon M \longrightarrow N$ be an embedding of topological spaces and let $D \subseteq M$ be a discrete subset. Then $\iota(D) \subseteq N$ is discrete.
\end{lemma}
\begin{proof}
    Let $x \in D$. By discreteness of the set $D \subseteq M$, there thus exists an
    open neighbourhood $U_x$ of $x$ such that $U_x \cap D = \{x\}$. As the mapping
    $\iota$ is an embedding, the image $\iota(U_x) \subseteq \iota(D)$ is an open
    neighbourhood of
    $\iota(x)$ with respect to the subspace topology $\iota(D) \subseteq N$. By
    injectivity of $\iota$, we moreover have
    \begin{equation}
        \iota(U_x)
        \cap
        \iota(D)
        =
        \iota
        \bigl(
            U_x \cap \iota(D)
        \bigr)
        =
        \iota
        \bigl(
            \{x\}
        \bigr)
        =
        \bigl\{
            \iota(x)
        \bigr\}.
    \end{equation}
    This completes the proof.
\end{proof}

\index{Properly discontinuous action}
Recall that an action $\Phi \colon H \times M \longrightarrow M$ of a topological
group $H$ on a topological space~$M$ is called properly discontinuous if every $x \in
M$ has a neighbourhood $U_x \subseteq M$ such that
\begin{equation}
    (h \acts U_x)
    \cap
    U_x
    \neq
    \emptyset
    \qquad
    \textrm{implies}
    \quad
    h = \E
    \quad
    \textrm{for all }
    h \in H.
\end{equation}
Notably, such an action is always free. We borrow the following elementary lemma from
topology.
\begin{lemma}
    \index{Discrete group}
    \label{lem:QuotientIsUniversal}
    Let $\Phi \colon H \times M \longrightarrow M$ be a continuous and properly discontinuous action of a topological group $H$ on a connected topological space $M$. Then $H$ carries the discrete topology and the quotient projection $\pi \colon M \longrightarrow M/H$ is a covering of topological spaces.
\end{lemma}
\begin{proof}
    Given $x \in M$, proper discontinuity yields a neighbourhood $U_x \subseteq M$ of
    $x$ such that $(h \acts U_x) \cap U_x = \emptyset$ for all $H \setminus \{\E\}$.
    As the map $\Phi$ is an action, this in particular implies
    \begin{equation}
        (h \acts U_x)
        \cap
        (g \acts U_x)
        =
        \emptyset
        \qquad
        \textrm{for all }
        h,g \in H
        \textrm{ with }
        h \neq g.
    \end{equation}
    As $\Phi$ acts by homeomorphisms, $h \acts U_x$ is moreover an open neighbourhood
    of $h \acts x$. Consequently, the fiber over $\pi(U_x)$ is given by
    \begin{equation*}
        \pi^{-1}
        \bigl(
            \pi(U_x)
        \bigr)
        =
        \bigsqcup_{h \in H}
        h \acts U_x,
    \end{equation*}
    where $\pi(U_x)$ is an open neighbourhood of $\pi(x)$ by openness of quotient
    projections. Next, we assert that $H$ carries the discrete topology. Assuming the
    contrary, we find a net~$(h_\alpha)_{\alpha \in J} \subseteq H$ converging to some
    $h \in H$ with $h_\alpha \neq h$ for all $\alpha \in J$. Left translation with
    $h^{-1}$ yields a convergent net $(h^{-1}h_\alpha)_{\alpha \in J}$ with limit
    $\E$. By continuity of $\Phi$, we get
    \begin{equation}
        \Phi
        \bigl(
            h^{-1} h_\alpha, x
        \bigr)
        \rightarrow
        \Phi(\E,x)
        =
        x.
    \end{equation}
    Consequently, given an open neighbourhood $U_x$ of $x$, there exists an $\alpha_0
    \in J$ such that we have $\Phi(h^{-1} h_\alpha, x) \in U_x$ for all $\alpha \later
    \alpha_0$. In particular, $(h^{-1} h_{\alpha_0} \acts U_x) \cap U_x \neq
    \emptyset$ and the quotient $h^{-1} h_{\alpha_0} \neq \E$, so the action of $H$ is
    not properly
    discontinuous. Hence $H$ indeed carries the discrete topology. Now, note that the
    restrictions $\pi \at{h \acts U_x}$ are injective for any fixed $h
    \in H$: If $\pi(h \acts x_1) = \pi(h \acts x_2)$ for some $x_1, x_2 \in U_x$, then
    there exists a $g \in H$ with
    \begin{equation}
         g
         \acts
         (h \acts x_2) = h \acts x_1
    \end{equation}
    and thus
    \begin{equation}
        x_1
        =
        (h^{-1}gh)
        \acts
        x_2
        \in
        U_x
        \cap
        (h^{-1}gh
        \acts
        U_x).
    \end{equation}
    By proper discontinuity, this means $h^{-1}gh = \E$, which implies $g = \E$ by
    left multiplication with $h$ and right multiplication with $h^{-1}$. But then $x_1
    = x_2$ and thus the restriction $\pi \at{h \acts U_x}$ is injective. Using once
    more the openness of $\pi$ and putting everything together, we obtain that the
    restrictions $\pi \at{h \acts U_x}$ are homeomorphisms onto their images for all
    $h \in H$.
\end{proof}

\index{Covering!Action}
For this reason and to avoid the arguably awkward terminology, some modern literature
such as \cite[Ch.~11 \& 12]{lee:2000a}, which also contains the above result, refers
actions that are both continuous and properly discontinuous as \emph{covering
actions}
instead.
In the special case of a subgroup $H \subseteq G$ acting on a topological group $M =
G$, it suffices to check proper discontinuity at $x = \E$:
\begin{lemma}
    \label{lem:ProperlyDiscontinuousOnGroups}
    Let $G$ be a topological group and $H \subseteq G$ be a subgroup such that there exists a zero neighbourhood $U_0 \subseteq G$ with $(h \acts U_0) \cap U_0 = \emptyset$ for all $h \in H \setminus \{\E\}$. Then the left action of $H$ on $G$ by left multiplications is properly discontinuous.
\end{lemma}
\begin{proof}
    Let $g \in G$ and define $U_g \coloneqq U_0 g$, which is a neighbourhood of $g$. Then
    \begin{equation}
        (h \acts U_g)
        \cap
        U_g
        =
        h U_0 g
        \cap
        U_0 g
        =
        (h U_0 \cap U_0)g
        =
        \emptyset \cdot g
        =
        \emptyset.
        \tag*{$\,$}
    \end{equation}
\end{proof}

Next, we prove that discrete subgroups always act properly discontinuously on the ambient group and thus induce covering projections by Lemma~\ref{lem:QuotientIsUniversal}.
\begin{lemma}
    \label{lem:QuotientsByDiscreteIsCovering}
    \index{Discrete!Covering}
    Let $H \subseteq G$ be a discrete subgroup of a Lie group $G$. Then the quotient projection $\pi \colon G \longrightarrow G/H$ is a covering of Lie groups.
\end{lemma}
\begin{proof}
     By continuity of the group multiplication and discreteness of $H$, there exists a
     zero neighbourhood $U_0 \subseteq \tilde{G}_\C$ with $(U_0 \cdot U_0) \cap H =
     \{\E\}$ and $U_0^{-1} = U_0$. Let $h \in H$ be such that
    \begin{equation*}
        (h
        \acts
        U_0)
        \cap
        U_0
        \neq
        \emptyset.
    \end{equation*}
    Taking a closer look, this means that there exist $g_1, g_2 \in U_0$ with
    \begin{equation*}
        h
        g_1
        =
        g_2,
        \quad \textrm{i.e.} \quad
        h
        =
        g_2 g_1^{-1}
        \in
        U_0 \cdot U_0^{-1}
        =
        U_0 \cdot U_0.
    \end{equation*}
    By choice of $U_0$, this means $h = \E$. Consequently, the action of $H$ is
    properly discontinuous in view of
    Lemma~\ref{lem:ProperlyDiscontinuousOnGroups} and thus the quotient
    projection $\pi \colon G \longrightarrow G/H$ is a covering of Lie groups by
    Lemma~\ref{lem:QuotientIsUniversal}.
\end{proof}

Putting everything together, we arrive at the following, a variant of which can be found in \cite[Exercise~15.1.2(b)]{hilgert.neeb:2012a}:
\begin{proposition}
    \label{prop:ComplexificationVsCovering}
    \index{Complexification!Fundamental group}
    Let $G$ be a connected Lie group such that $\tilde{\eta} \colon \tilde{G} \longrightarrow \tilde{G}_\C$ is injective. Then the quotient projection $\pi \colon \tilde{G}_\C \longrightarrow G_\C$ is the universal covering projection of the universal complexification $G_\C$ of $G$. In particular,
    \begin{equation}
        \label{eq:ComplexificationVsCovering}
        \pi_1(G_\C)
        \cong
        \tilde{\eta}
        \bigl(
            \pi_1(G)
        \bigr)
        \cong
        \pi_1(G).
    \end{equation}
\end{proposition}
\begin{proof}
    By Lemma~\ref{lem:Embedding}, $\tilde{\eta}$ is an embedding. Thus the image $\tilde{\eta}(\pi_1(G))$ of the discrete subgroup $\pi_1(G) \subseteq \tilde{G}$ is a discrete complex subgroup of $\tilde{G}_\C$ by Lemma~\ref{lem:Discrete}. Lemma~\ref{lem:QuotientsByDiscreteIsCovering} implies then that the quotient projection $\pi \colon \tilde{G}_\C \longrightarrow G_\C$ is a covering projection. As $\tilde{G}_\C$ is simply connected, $\pi$ is universal and \eqref{eq:ComplexificationVsCovering} follows immediately.
\end{proof}

Summarizing, if $\tilde{\eta}$ injective, passing to universal complexifications and universal coverings are commutative processes, i.e. the diagram
\begin{equation}
    \label{eq:UniversalComplexificationVsCovering}
    \begin{tikzcd}[column sep = huge, row sep = huge]
        \tilde{G}
        \arrow[r, "\tilde{\eta}"]
        \arrow[d, "p"]
        &\tilde{G}_C
        \arrow[d, "\pi"] \\
        G
        \ar[r, "\eta"]
        &G_\C
    \end{tikzcd}
\end{equation}
commutes and its vertical arrows are universal covering projections. In particular,
there are no new topological obstructions to analytic continuation to $G_\C$ beyond
what is already present within $G$. This will be crucial for our considerations in
Section~\ref{sec:ExtensionAndRestriction}, where we extend entire functions defined on
$G$, i.e. functions already abiding by all topological obstructions, to holomorphic
ones
on $G_\C$. As another aside, we once again see that $\eta$ is locally injective in
this situation, as all other arrows in \eqref{eq:UniversalComplexificationVsCovering}
are. Truncating our proof of Proposition~\ref{prop:ComplexificationVsCovering}, we
have moreover shown the following.
\begin{corollary}
    \label{cor:EtaOfPi1Discrete}%
    Let $G$ be a connected Lie group such that $\tilde{\eta}(\pi_1(G)) \subseteq
    \tilde{G}_\C$ is discrete. Then the quotient projection $\pi \colon \tilde{G}_\C
    \longrightarrow G_\C$ is the universal covering morphism and
    \begin{equation}
        \pi_1(G_\C)
        \cong
        \tilde{\eta}
        \bigl(
            \pi_1(G)
        \bigr).
    \end{equation}
    In particular,
    \begin{equation}
        \pi_1(G)
        \supseteq
        \pi_1(G_\C)
        \quad \textrm{and} \quad
        \LieAlg(G_\C)
        =
        \LieAlg(\tilde{G}_\C).
    \end{equation}
\end{corollary}

In Example~\ref{ex:SpecialLinear} we will see that dropping the assumption that
$\tilde{\eta} \colon \tilde{G} \longrightarrow \tilde{G}_\C$ is injective in
Proposition~\ref{prop:ComplexificationVsCovering} may result in $\pi_1(G_\C)
\subsetneq \pi_1(G)$. We are now in a position to revisit
Corollary~\ref{cor:ImageOfEtaSubmanifold}, for which we first note the following
lemmas, which are also of independent interest. Recall that a Lie group morphism $\Phi
\colon G \longrightarrow H$ is called locally injective at a group element $g \in
G$ if there
exists a neighbourhood $U \subseteq G$ of $g$ such that the restriction $\Phi
\at{U} \colon U \longrightarrow H$ is injective.
\begin{lemma}
    Let $\Phi \colon G \longrightarrow H$ be a Lie group morphism. Then $\Phi$ is
    locally injective at the group unit if and only if its tangent map $T_\E \Phi$ is
    injective.
\end{lemma}
\begin{proof}
    The interesting implication is going from local injectivity at $\E$ to injectivity
    of the tangent map, as the converse statement is just the usual inverse
    function
    Theorem. Thus let $U \subseteq G$ be a neighbourhood of $\E$ such that
    $\Phi
    \at{U} \colon U \longrightarrow H$ is injective. Moreover, we choose a zero
    neighbourhood $V \subseteq \liealg{g}$ such that $\exp_G \colon V \longrightarrow
    G$ is a diffeomorphism with image $\exp(V)$ contained within $U$. By
    injectivity of $\Phi$, this implies injectivity of
    \begin{equation}
        \exp_{H} \circ T_\E \Phi
        =
        \Phi \circ \exp_G
        \colon
        V
        \longrightarrow
        H,
    \end{equation}
    which implies injectivity of $T_\E \Phi$ on $V$. But $T_\E \Phi$ is a linear map,
    so injectivity on a zero neighbourhood implies its global injectivity.
\end{proof}

Note that even global injectivity of a generic smooth map does not necessarily lead to
injectivity of the tangent map. For a simple example, consider the mapping $x \mapsto
x^3$ on the real line. The underlying problem is that this type of map is not a
diffeomorphism. As bijective Lie group morphisms are always diffeomorphisms, this
yields another way to prove the preceding lemma.
\begin{lemma}
    \label{lem:ComplexificationEtaLocallyInjective}
    Let $G$ be a connected Lie group with universal complexification $(G_\C, \eta)$
    such that $\eta \colon G \longrightarrow G_\C$ is locally injective at the group
    unit. Then $\LieAlg(G_\C) = \liealg{g}_\C$.
\end{lemma}
\begin{proof}
    Taking another look at \eqref{eq:ComplexificationFromUniversalCoveringDiagram}, we
    see that by our assumption, the composition~$\pi \circ \tilde{\eta}$ is locally
    injective at $\E$. Let $U \subseteq G$ be a neighbourhood of $\E$ such that
    \begin{equation}
        \pi
        \circ
        \tilde{\eta}
        \at{U}
        \colon
        U \longrightarrow G_\C
    \end{equation}
    is injective. By \eqref{eq:ComplexificationFromUniversalCovering}, this means that
    \begin{equation}
        \tilde{\eta}(U)
        \cap
        \tilde{\eta}
        \bigl(
            \pi_1(G)
        \bigr)
        =
        \{\E\}.
    \end{equation}
    That is, $\tilde{\eta}(\pi_1(G))$ is discrete in $\tilde{G}_\C$. By
    Corollary~\ref{cor:EtaOfPi1Discrete}, this implies $\LieAlg(G_\C) = \liealg{g}_\C$.
\end{proof}

As a special case we obtain the following result about the fine structure of $\eta(G)$ for a generic connected Lie group $G$, which will be useful for our considerations in Section~\ref{sec:ExtensionAndRestriction}.
\begin{corollary}
    \label{cor:ImageOfEtaLieAlgebra}
    Let $G$ be a connected Lie group and $k \coloneqq \dim_\R(\eta(G)))$. Every basis
    \begin{equation}
        \bigl(
            \basis{e}_1, \ldots, \basis{e}_k
        \bigr)
        \subseteq
        \LieAlg\bigl(\eta(G)\bigr)
    \end{equation}
    induces a basis of $\LieAlg(G_\C)$ via $(\basis{e}_1, \ldots, \basis{e}_k, \I
    \basis{e}_1, \ldots, \I \basis{e}_k)$. In particular,
    \begin{equation}
        \dim_\C(G_\C)
        =
        \dim_\R
        \bigl(
            \eta(G)
        \bigr).
    \end{equation}
\end{corollary}
\begin{proof}
    By Corollary~\ref{cor:ImageOfEtaSubmanifold}, \ref{item:ImageOfEtaComplexification}, the embedding $\eta(G) \hookrightarrow G_\C$ constitutes a universal complexification. As injective mappings are in particular locally injective at every point, Lemma~\ref{lem:ComplexificationEtaLocallyInjective} yields $\LieAlg(G_\C) = \LieAlg(\eta(G))_\C$, from which the assertions are immediate.
\end{proof}

Combining Corollary~\ref{cor:ImageOfEtaLieAlgebra} with Proposition~\ref{prop:CauchyRiemann} yields the following version of the identity principle, which will prove useful in many places.
\begin{corollary}[Identity Principle]
    \label{cor:IdentityPrinciple}
    \index{Identity principle}
    \index{Left invariant complex analysis!Identity principle}
    Let $G$ be a connected Lie group and let
    \begin{equation}
        \phi \colon G_\C \longrightarrow M
    \end{equation}
    be a holomorphic map into a complex manifold $M$. If $\phi \at{\eta(G)}$ is
    constant, then $\phi$ is constant.
\end{corollary}
\begin{proof}
    Choosing a holomorphic chart $(U,z)$ of $M$, we may consider the component mappings $z^j \circ \phi \colon G_\C \longrightarrow \C$. By assumption, each of them is constant on $\eta(G)$, which by Proposition~\ref{prop:CauchyRiemann} implies they are constant on the connected component of the group unit. By connectedness of $G$, the assertion follows.
\end{proof}

We round out our abstract considerations about universal complexification by taking a closer look at the disconnected case. The idea is that one may recast the construction as a Lie group extension problem.
\begin{remark}[Interpretation of the disconnected case]
    \label{rem:ComplexificationDisconnected}
    \index{Universal Complexification!Disconnected groups}
    \index{Complexification!Disconnected groups}
    Consider the short exact sequence
    \begin{equation}
        \label{eq:ComponentGroupDiagram}
        \begin{tikzcd}[column sep = huge]
            0
            \ar[r]
            &G_0
            \ar[r, "\iota"]
            &G
            \ar[r, "\pi"]
            &G \big/ G_0
            \ar[r]
            &0
        \end{tikzcd}
    \end{equation}
    with the inclusion $\iota \colon G_0 \hookrightarrow G$ and the
    quotient projection $\pi \colon G \longrightarrow G / G_0$ to the
    component group $G / G_0$. This group should also
    be the component group of the complexification, i.e.
    \begin{equation}
        G_\C
        /
        G_{0, \C}
        \cong
        G / G_0
    \end{equation}
    By Example~\ref{ex:DiscreteGroups} this equality makes sense
    both as real and as complex Lie groups if~$G$ has positive
    dimension. Otherwise, the respective unit components are just
    given by a point and we obtain $G_\C \cong G$ once again in
    accordance to Example~\ref{ex:DiscreteGroups}. Note now that
    by the connected case, we do have a universal complexification
    $(G_{0, \C}, \eta_0)$ of the connected component $G_0$ of
    the group unit. Thus the situation looks like
    \begin{equation}
        \label{eq:UniversalComplexificationDiagram}
        \begin{tikzcd}[column sep = huge]
            0
            \ar[r]
            &G_0
            \ar[r, "\iota"]
            \ar[d, "\eta_0"]
            &G
            \ar[r, "\pi"]
            &G \big/ G_0
            \ar[r]
            \ar[d, "\id", leftrightarrow]
            &0 \\
            0
            \ar[r]
            &G_{0, \C}
            \ar[r, dashed, "\iota_\C"]
            &G_\C
            \ar[r, dashed, "\pi_\C"]
            \ar[u, leftarrow, dashed, "\eta"]
            &G \big/ G_0
            \ar[r]
            &0,
        \end{tikzcd}
    \end{equation}
    where we have to construct $G_\C$ and the dashed arrows in the appropriate
    categories
    such that the diagram commutes and the second row is exact. Doing this for the
    second line only is known as a \emph{group extension} problem,
    which typically has a multitude of inequivalent solutions. It is
    now the additional conditions derived from the commutativity of
    \eqref{eq:UniversalComplexificationDiagram} that make the pair
    $(G_\C, \eta)$ unique. This type of problem gets discussed
    e.g. in \cite[Chapter~3§6, Lemma~7]{bourbaki:1975a} in some
    more generality. In the proof of
    Theorem~\ref{thm:UniversalComplexification} the complexified
    Lie group was then realized as the quotient $(G_{0, \C}
    \rtimes_c G) / K$ with respect to an action by conjugations $c$
    and $K$ as in \eqref{eq:UniversalComplexificationProof2}. To
    obtain the remaining Lie group morphisms, the crucial observation
    is that we have as a set $G_{0, \C} \rtimes_c G = G_{0,
      \C} \times G$. Notably, products with only one non
    trivial factor in $G_{0, \C} \rtimes_c G$ just reduce to
    the products of $G_{0, \C}$ and $G$, respectively. Consequently,
    simply embedding $G_{0,\C}$ or $G$ yields Lie group morphisms.
    In particular, we can map from $G_{0, \C}$ and $G$ into
    $G_\C$ by taking equivalence classes afterwards, that is
    \begin{equation*}
        \iota_\C
        (g)
        \coloneqq
        \big[
            (g, \E_G)
        \big]
        \quad \textrm{and} \quad
        \eta(g)
        \coloneqq
        \big[
            (\E_{G_{0, \C}}, g)
        \big],
    \end{equation*}
    which are exactly the maps we used in the proof of
    Theorem~\ref{thm:UniversalComplexification}. We have
    \begin{equation*}
        \iota_\C
        \circ
        \eta_0
        \at{g}
        =
        \big[
            (\eta_0(g), \E_G)
        \big]
        =
        \big[
            (\E_{G_{0, \C}}, \iota(g))
        \big]
        =
        \eta
        \circ
        \iota
        \at{g}
    \end{equation*}
    for $g \in G_0$ by \eqref{eq:UniversalComplexificationProof3}, i.e. the left square
    in \eqref{eq:UniversalComplexificationDiagram} commutes. Note that
    we only used elements $g \in G_0$ in the definition of $K$ and
    \emph{suppressed} the map $\iota$ there. As every equivalence class $[(z,g)]$ contains a unique representative of the form $(z',e_G)$, the morphism $\iota_\C$ is
    injective and thus the diagram
    \eqref{eq:UniversalComplexificationDiagram} is exact at $G_{0,
      \C}$. An alternative way of obtaining the morphism
    $\iota_\C$ is by means of the universal property of $(G_{0,
      \C}, \eta_0)$: the map $\eta \circ \iota \colon G_0
    \longrightarrow G_\C$ is a Lie group morphism into a
    complex Lie group. Taking another look at the proof, we then used
    Lie's second Theorem~\ref{thm:Lie2} to integrate the
    complexified tangent map, which is not very explicit. By the
    uniqueness statement in the universal property, this however
    reproduces our $\iota_\C$. Similarly, we can use the
    universal property of $(G_\C, \eta)$ to obtain the final
    map $\pi_\C$ by once again using
    Example~\ref{ex:DiscreteGroups} to view $G / G_0$ as
    a complex Lie group: as our notation already suggests, the map
    $\pi_\C$ really is the complexification of the quotient
    projection $\pi \colon G \longrightarrow G / G_0$. The right
    square in \eqref{eq:UniversalComplexificationDiagram} thus
    commutes by construction, which then means also the outer square
    commutes. Taking another look at the construction of
    $\pi_\C$, we first note that $\pi \circ \iota = [\E_G]$ by
    exactness of~\eqref{eq:ComponentGroupDiagram}. Thus
    $\pi \at{G_0}$ complexifies to
    \begin{equation*}
        \pi_\C
        \bigl(
            [(z, g)]
        \bigr)
        =
        \pi(g),
    \end{equation*}
    which is well defined, as we have for $(\eta_0(h) z,
    h^{-1}g) \in [(z, g)]$ also $h^{-1} g \in G_0 g$,
    i.e. multiplication with an element of $G_0$ does not change the
    connected component. As $\pi$ is surjective, the same is thus true
    for $\pi_\C$,
    i.e. \eqref{eq:UniversalComplexificationDiagram} is exact also at
    the lower $G / G_0$. Finally, we compute
    \begin{equation*}
        \pi_\C
        \circ
        \iota_\C
        \at{g}
        =
        \pi_\C
        \bigl(
            [(g, \E_G)]
        \bigr)
        =
        \pi(\E_G)
        =
        [\E_G]
        =
        \E_{G / G_0}
        \qquad
        \textrm{for all }
        g \in G_{0, \C},
    \end{equation*}
    which also proves exactness at $G_\C$.
\end{remark}

Summarizing, we have the following:
\begin{corollary}
    Let $G$ be a Lie group with universal complexification
    $(G_\C, \eta)$.
    \begin{corollarylist}
        \item The component groups of $G$ and $G_\C$ are
          isomorphic as both real and complex Lie groups.
        \item The universal complexification $G_\C$ is a group
          extension of $G_{0, \C}$ relative to the component
          group $G / G_0$.
        \item We have $G_{0, \C} = G_{\C, 0}$, i.e. the
          universal complexification of the unit component $G_0$ of
          $G$ is given by the unit component of the universal
          complexification $G_\C$ of $G$.
        \item The diagram of real Lie groups
        \begin{equation}
            \begin{tikzcd}[column sep = huge]
                0
                \ar[r]
                &G_0
                \ar[r, "\iota"]
                \ar[d, "\eta_0"]
                &G
                \ar[r, "\pi"]
                &G \big/ G_0
                \ar[r]
                \ar[d, "\id", leftrightarrow]
                &0 \\
                0
                \ar[r]
                &G_{0, \C}
                \ar[r, "\iota_\C"]
                &G_\C
                \ar[r, "\pi_\C"]
                \ar[u, leftarrow, "\eta"]
                &G \big/ G_0
                \ar[r]
                &0,
            \end{tikzcd}
        \end{equation}
        is commutative with exact lines and $\iota_\C$ as well as $\pi_\C$ are holomorphic.
    \end{corollarylist}
\end{corollary}

\section{Examples of Universal Complexifications}
\label{sec:UniversalComplexificationExamples}
\epigraph{Moist was, by inclination, a stranger to the concept of two in the morning, a
time that happened to other people.}{\emph{Raising Steam} -- Terry Pratchett}
% !TeX root = ../Dissertation.tex
The purpose of this section is to fill the abstract machinery we have developed with life
in the form of concrete constructions and counterexamples. While our proof of
Theorem~\ref{thm:UniversalComplexification} is constructive, it is typically difficult to
make all the required structures concrete for a given example. The circle group turns out
to be sufficiently simple.
\begin{example}[Circle group, {\cite[Ex.~4.19]{heins.roth.waldmann:2023a}}]
    \label{ex:CircleGroupUniversalComplexification}%
    \index{Universal Complexification!Circle group}
    \index{Complexification!Circle group}
    Let
    \begin{equation}
        G
        \coloneqq
        \gls{OneSphere}
        \coloneqq
        \bigl\{
            z \in \C
            \colon
            \abs{z}
            =
            1
        \bigr\}
        =
        \boundary \disk
        =
        \liegroup{U}(1),
    \end{equation}
    which is connected, but not simply
    connected. Its universal covering group $\widetilde{G}$ is given by
    \begin{equation}
        \label{eq:CircleGroupUniversalCovering}%
        p
        \colon
        \widetilde{\mathbb{S}}^1 = \R
        \longrightarrow
        G, \quad
        p(t)
        =
        e^{2 \pi \I t}.
    \end{equation}
    Any interval of length less than $1$ gives a evenly covered
    neighbourhood of all its points in $G$. The complexification of $\R$ is
    \begin{equation}
        \widetilde{\eta}
        \colon
        \R
        \longrightarrow
        \C = \widetilde{G}_\C, \quad
        \widetilde{\eta}(x)
        =
        x + 0\I,
    \end{equation}
    see again Example~\ref{ex:ComplexificationVectorSpace}.
    Taking another look at \eqref{eq:CircleGroupUniversalCovering}, we see
    that $\pi_1(G) = \field{Z} \subseteq \R$ and thus also
    $\widetilde{\eta}(\pi_1(G)) = \field{Z} \subseteq \C$. This is
    already a discrete and thus closed subgroup. That is to say, $\field{Z} =
    \left<\field{Z}\right>_{\C}$. Consequently, the universal complexification of the circle
    group is
    \begin{equation}
        \mathbb{S}^1_\C
        =
        \C
        \big/
        \field{Z}
        \quad \textrm{and} \quad
        \eta
        \colon
        \mathbb{S}^1
        \longrightarrow
        \mathbb{S}^1_\C, \quad
        \eta
        \bigl(
        \E^{2 \pi \I t}
        \bigr)
        =
        \big[
        \widetilde{\eta}(t)
        \big]
        =
        [t],
    \end{equation}
    which is indeed well defined. Geometrically, $\mathbb{S}^1_\C$ is thus given by a
    cylinder $\mathbb{S}^1 \times \R$. Notably, the universal complexification~
    $\mathbb{S}^1_\C$ is not compact, even though the circle
    group~$\mathbb{S}^1$ was. This is good news for the prospects of interesting
    globally defined holomorphic functions. We shall see in
    Section~\ref{sec:ExtensionAndRestriction} that this non-compactness is a
    general feature, which makes holomorphic extension of entire functions on
    $\mathbb{S}^1$ to holomorphic ones on $\mathbb{S}^1_\C$ feasible. Ultimately, this
    results in Theorem~\ref{thm:CompactExtension}.

    Analogously, one obtains the universal complexification for higher tori as
    products of cylinders. Notice that the group morphism
    \begin{equation}
        \sigma
        \colon
        \mathbb{S}^1_\C
        \longrightarrow
        \mathbb{S}^1_\C, \quad
        \sigma\bigl( e^{2\pi \I t}, r \bigr)
        =
        (e^{2 \pi \I t}, -r)
    \end{equation}
    fixes all $g \in \eta(\mathbb{S}^1) = \mathbb{S}^1 \times \{ 0 \}$.
    Moreover, we have
    \begin{equation}
        T_{(1, 0)} \sigma
        =
        \id_{\R}
        \times
        \bigl(
        - \id_{\R}
        \bigr),
    \end{equation}
    which is $\C$-antilinear. Thus $\sigma$ is the complex conjugation
    of $\mathbb{S}^1_\C$ from Proposition~\ref{prop:ComplexConjugation}.
\end{example}

\begin{example}[Semidirect products, {\cite[Exercise~15.1.7]{hilgert.neeb:2012a}}]
    \label{ex:ComplexificationOfProducts}
    \index{Universal Complexification!Semidirect products}
    \index{Complexification!Semidirect products}
    Let $G$ and $N$ be Lie groups and $\gamma \colon G \times N \longrightarrow N$
    be a smooth action of $G$ by automorphisms of $N$. We write
    \begin{equation}
        \gamma_g
        \coloneqq
        \gamma(g,\argument)
        \in
        \gls{AutomorphismGroup}(N)
        \qquad
        \textrm{for }
        g \in G.
    \end{equation}
    By the universal property of $N_\C$, the compositions $\eta_N \circ \gamma_g \colon
    N \longrightarrow N_\C$
    induce unique holomorphic group morphisms $\widetilde{\gamma}_g \colon N_\C
    \longrightarrow N_\C$ such that
    \begin{equation}
        \widetilde{\gamma}_g
        \circ
        \eta_N
        =
        \eta_N
        \circ
        \gamma_g
        \qquad
        \textrm{for all }
        g \in G.
    \end{equation}
    This yields a smooth action $\widetilde{\gamma} \colon
    G \times N_\C \longrightarrow N_\C$ by holomorphic automorphisms of $N_\C$ as
    \begin{equation*}
        \widetilde{\gamma}_{gh}
        \circ
        \eta_N
        =
        \eta_N
        \circ
        \gamma_{gh}
        =
        \eta_N
        \circ
        \gamma_g
        \circ
        \gamma_h
        =
        \widetilde{\gamma}_g
        \circ
        \eta_N
        \circ
        \gamma_h
        =
        \widetilde{\gamma}_g
        \circ
        \widetilde{\gamma}_h
        \circ
        \eta_N
    \end{equation*}
    holds for all $g,h \in G$. That is, $\eta_N \circ \gamma_g \circ \gamma_h$ and $\eta_N \circ \gamma_{gh}$ are simply the same morphisms and thus also induce the same holomorphic extensions. For every $n \in N_\C$, we now consider the group morphisms
    \begin{equation*}
        \psi_n
        \colon
        G
        \longrightarrow
        N_\C, \quad
        \psi_n(g)
        \coloneqq
        \widetilde{\gamma}_g(n).
    \end{equation*}
    By the universal property of $G_\C$, there exist unique holomorphic group morphisms
    \begin{equation}
        \widetilde{\psi}_n
        \colon
        G_\C
        \longrightarrow
        N_\C
        \quad \textrm{with} \quad
        \widetilde{\psi}_n
        \circ
        \eta_G
        =
        \psi_n.
    \end{equation}
    Repeating the arguments above, we see that the $\widetilde{\psi}_n$ induce a
    smooth action
    \begin{equation}
        \gamma_\C
        \colon
        G_\C \times N_\C
        \longrightarrow
        N_\C
    \end{equation}
    by biholomorphic morphisms of $N_\C$. By construction, both of the partial maps
    \begin{equation}
        \gamma_\C(g,\argument)
        =
        \widetilde{\gamma}_g
        \colon
        N_\C
        \longrightarrow
        N_\C
        \quad \textrm{and} \quad
        \gamma_\C(\argument,n)
        =
        \widetilde{\psi}_n
        \colon
        G_\C
        \longrightarrow
        N_\C
    \end{equation}
    are holomorphic for all $g \in G_\C$ and $n \in N_\C$. By Hartogs' Theorem, see
    again Remark~\ref{rem:Hartog}, this implies that the full map $\gamma_\C$ is
    holomorphic, as well. Consequently, the semidirect product $N_\C
    \rtimes_{\gamma_\C} G_\C$ endowed with the product atlas carries the structure of a
    complex Lie group: its multiplication may be decomposed into the multiplications of
    the complex Lie groups $G_\C$, $N_\C$ and applications of $\gamma_\C$. Finally,
    we prove that
    \begin{equation}
        \label{eq:UniversalComplexificationProduct}
        \eta
        \colon
        N \rtimes_\gamma G
        \longrightarrow
        N_\C \rtimes_{\gamma_\C} G_\C,
        \quad
        \eta
        \coloneqq
        \eta_N \times \eta_G
    \end{equation}
    constitutes a universal complexification. To show this, let $\Phi \colon N
    \rtimes_\gamma G \longrightarrow H$ be a group morphism into a complex Lie group
    $H$. Denoting the canonical inclusions by
    \begin{equation}
        \iota_G
        \colon
        G \longrightarrow
        N \rtimes_\gamma G
        \quad \textrm{and} \quad
        \iota_N
        \colon
        N
        \longrightarrow
        N \rtimes_\gamma G,
    \end{equation}
    we apply the universal properties of $G_\C$ resp. $N_\C$ to
    \begin{equation}
        \Phi_G
        \coloneqq
        \iota_G^*
        \Phi
        \colon
        G \longrightarrow H
        \quad \textrm{resp.} \quad
        \Phi_N
        \coloneqq
        \iota^*_N \Phi
        \colon
        N
        \longrightarrow
        H.
    \end{equation}
    This provides holomorphic group morphisms $\Phi_{G,\C} \colon
    G_\C \longrightarrow H$ fulfilling $\Phi_{G,\C} \circ \eta_G = \Phi_G$ and $\Phi_{N,\C}
    \colon N_\C \longrightarrow H$ fulfilling $\Phi_{N,\C} \circ \eta_N = \Phi_N$. Using
    both maps, we define
    \begin{equation}
        \Phi_\C
        \colon
        N_\C \rtimes_{\gamma_\C} G_\C
        \longrightarrow
        H, \quad
        \Phi_\C(n,g)
        \coloneqq
        \Phi_{N,\C}(n)
        \cdot
        \Phi_{G,\C}(g),
    \end{equation}
    which is holomorphic. By construction, we have for $g \in G$ and $n \in N$ that
    \begin{equation}
        \Phi_\C
        \circ
        \eta
        \at[\Big]{(g,n)}
        =
        \Phi_{N,\C}
        \bigl(
        \eta_N(n)
        \bigr)
%        \cdot
        \Phi_{G,\C}
        \bigl(
        \eta_G(g)
        \bigr)
        =
        \Phi_N(n)
%        \cdot
        \Phi_G(g)
        =
        \Phi(e_G,n)
%        \cdot
        \Phi(g,e_N)
        =
        \Phi(g,e).
    \end{equation}
    As $\eta(N \rtimes_\gamma G) = \eta_N(N) \times \eta_G(G)$, the identity principle
    from Corollary~\ref{cor:IdentityPrinciple} implies that~$\Phi_\C$ is the group
    morphism uniquely determined by $\Phi_\C \circ \eta = \Phi$.

    Note that \eqref{eq:UniversalComplexificationProduct} implies that $\eta$ is
    (locally) injective if and only if both $\eta_N$ and $\eta_G$ are. Thus twisting
    with a group action can never lose dimensions beyond what was already lost upon
    complexifying the factors.
\end{example}

As a special case, namely taking $\gamma$ as the trivial action, we see that universal
complexification distributes over finite Cartesian products. Note also that, if the
automorphism group~$\Aut(N)$ carries a Lie group structure, then the complex
automorphism group~$\Aut_\C(N_\C)$ may be endowed with the structure of a complex Lie
group by \cite[Prop.~15.4.1]{hilgert.neeb:2012a}. In this case, one may
construct $\gamma_\C$ from $\widetilde{\gamma}$ directly by applying the universal
property of $G_\C$, which simplifies the preceding construction somewhat. We continue
with the principal example within Lie theory.
\begin{example}[Special linear group]
    \label{ex:SpecialLinear}
    \index{Universal Complexification!Special linear group}
    \index{Complexification!Special linear group}
    \index{Special linear group!Complexification}
    \index{Special linear group!Universal Covering}
    We consider the special linear group
    \begin{equation}
        G
        \coloneqq
        \gls{SpecialLinear}
        =
        \bigl\{
            M \in \gls{Matrices}(\R)
            \colon
            \det(M) = 1
        \bigr\},
    \end{equation}
    which is connected, but not simply connected. Its universal covering group is
    $\tilde{G} = \widetilde{\liegroup{SL}}_2(\R)$, which is a \emph{non linear} group
    by \cite[Example~9.5.18]{hilgert.neeb:2012a}. We denote the universal covering
    projection by $p \colon \widetilde{\liegroup{SL}}_2(\R) \longrightarrow
    \liegroup{SL}_2(\R)$.
    The Lie algebras of both groups are given by the traceless matrices
    \begin{equation}
         \gls{SpecialLinearAlgebra}
        \coloneqq
        \bigl\{
            M
            \in
            \Mat_2(\R)
            \colon
            \tr(M)
            =
            0
        \bigr\},
    \end{equation}
    whose complexification is simply
    \begin{equation}
        \liealg{sl}_2(\R)_\C
        =
        \bigl\{
            M
            \in
            \Mat_2(\C)
            \colon
            \tr(M)
            =
            0
        \bigr\}
        \eqqcolon
        \liealg{sl}_2(\C).
    \end{equation}
    Consequently, by our construction within
    Theorem~\ref{thm:UniversalComplexification}, we have
    \begin{equation}
        \widetilde{\liegroup{SL}}_2(\R)_\C
        =
        \liegroup{SL}_2(\C)
        =
        \bigl\{
            M
            \in
            \Mat_2(\C)
            \colon
            \det(M)
            =
            1
         \bigr\},
    \end{equation}
    as $\liegroup{SL}_2(\C)$ is simply connected by
    \cite[Exercise~15.1.1]{hilgert.neeb:2012a} with Lie algebra $\liealg{sl}_2(\C)$. Next,
    we prove that $\tilde{\eta} = \iota \circ p$, where $\iota \colon \liegroup{SL}_2(\R)
    \hookrightarrow \liegroup{SL}_2(\C)$ simply reinterprets real matrices as complex
    ones. Indeed, by construction its tangent map is given by
    \begin{equation}
        T_\E
        \tilde{\eta}
        =
        T_\E
        \iota
        =
        T_\E
        \iota
        \circ
        T_\E p ,
    \end{equation}
    where $T_\E \iota$ is the inclusion $\liealg{sl}_2(\R) \hookrightarrow \liealg{sl}_2(\C)$
    and by connectedness the claimed form of~$\tilde{\eta}$ follows. Notably,
    $\tilde{\eta}$ is not injective, but locally injective in accordance with
    Corollary~\ref{cor:UniversalComplexificationProperties},~\ref{item:UniversalComplexificationSimplyConnected2}.
     Moreover, we get $\tilde{\eta}(\pi_1(G)) = \{\E\}$, as $\pi_1(G) = \ker p$.
    Consequently,~\eqref{eq:ComplexificationFromUniversalCovering} yields that
    $\liegroup{SL}_2(\R)_\C = \liegroup{SL}_2(\C)$ with $\eta \coloneqq \iota$ in
    accordance with Proposition~\ref{prop:ClassicalComplexifications}. This example is
    remarkable in at least two regards: firstly, we see that non isomorphic Lie groups may
    have isomorphic universal complexifications. Secondly, the conclusion in
    Proposition~\ref{prop:ComplexificationVsCovering} may fail if $\tilde{\eta}$ is not
    injective. Indeed, we in particular have
    \begin{equation}
        \pi_1(G_\C)
        =
        \pi_1
        (\liegroup{SL}_2(\C))
        =
        \{\E\}
        \subsetneq
        \Z
        =
        \pi_1
        \bigl(
            \liegroup{SL}_2(\R)
        \bigr)
        =
        \pi_1(G).
    \end{equation}
    On the other hand, the example matches what we expect from
    Corollary~\ref{cor:EtaOfPi1Discrete}.
\end{example}

Keeping some of our notation and combining all three of the preceding
examples allows us to show that the inequality $\dim_\C G_\C \le \dim_\R G$
may be strict.
\begin{example}[Dimension drop, {\cite[Exercise~15.1.4]{hilgert.neeb:2012a}}]
    \label{ex:DimensionDrop}
    \index{Universal Complexification!Dimension drop}
    \index{Complexification!Dimension drop}
    Consider the quotient
    \begin{equation}
        G
        \coloneqq
        \bigl(
            \widetilde{\liegroup{SL}}_2(\R)
            \times
            \R
        \bigr)
        \big/
        D
    \end{equation}
    by the subgroup
    \begin{equation}
        \label{eq:DimensionDropD}
        D
        \coloneqq
        \big<
            (g_0, 1),
            (\E, \sqrt{2})
        \big>
        \subseteq
        \widetilde{\liegroup{SL}}_2(\R)
        \times
        \R,
    \end{equation}
    where $\left<g_0\right> = \pi_1(\liegroup{SL}_2(\R))$ is a generator. Explicitly, one may choose
    $u \coloneqq (\begin{smallmatrix}
        0 & -1 \\
        1 & 0
    \end{smallmatrix}) \in \liealg{sl}_2(\R)$
    and $g_0 \coloneqq \exp_{\widetilde{\liegroup{SL}}_2(\R)}(2\pi u)$. Crucially,
    $\widetilde{\liegroup{SL}}_2(\R)$ is an infinite cyclic covering\footnote{That is
    to say, the
    fibers $p^{-1}(\{A\})$ are countably infinite cyclic groups for all $A \in
    \liegroup{SL}_2(\R)$. To see this, consider the polar decomposition into the circle
    group $\mathbb{S}^1$ times positive matrices, which form a cone and are as such
    contractible.} and
    thus the powers of $g_0$ are distinct. Consequently,
    \begin{equation}
        D
        =
        \bigl\{
        \bigl(
            g_0^k,
            k+\ell \sqrt{2}
        \bigr)
        \colon
        k,\ell \in \Z
        \bigr\}
    \end{equation}
    is indeed discrete, second countable and in particular closed. Thus $G$ carries the
    structure of a Lie group itself. Its universal covering group is given by the quotient
    projection
    \begin{equation}
        \pi
        \colon
        \tilde{G}
        =
        \widetilde{\liegroup{SL}}_2(\R)
        \times
        \R
        \longrightarrow
        G
    \end{equation}
    by Lemma~\ref{lem:QuotientsByDiscreteIsCovering}. By the compatibility of universal complexification with direct products,
    \begin{equation}
        \tilde{\eta}
        \colon
        \tilde{G}
        \longrightarrow
        \tilde{G}_\C
        \coloneqq
        \liegroup{SL}_2(\C)
        \times
        \C, \quad
        \tilde{\eta}
        \coloneqq
        p \times \id
    \end{equation}
    is the universal complexification of $\tilde{G}$, see Example~\ref{ex:SpecialLinear}. By construction, we arrive at
    \begin{equation}
        \tilde{\eta}
        \bigl(
            \pi_1(G)
        \bigr)
        =
        \tilde{\eta}(D)
        =
        \Bigl\{
            \bigl(
            \mathbb{1}_2,
            k + \ell \sqrt{2}
            \bigr)
            \colon
            k,\ell \in \Z
        \Bigr\}
    \end{equation}
    with $\left<\tilde{\eta}(\pi_1(G))\right>_{\C} = \{\mathbb{1}_2\} \times \C$. We have shown
    \begin{equation}
        G_\C
        =
        \tilde{G}_\C
        /
        \left<\tilde{\eta}(\pi_1(G))\right>_{\C}
        =
        \liegroup{SL}_2(\C)/\{\mathbb{1}_2\}
        \times
        \C/\C
        =
        \liegroup{SL}_2(\C)
        \times
        \{\E\}
        \cong
        \liegroup{SL}_2(\C)
    \end{equation}
    with universal complexification map inherited from $\liegroup{SL}_2(\C)$. In particular, we see that
    \begin{equation}
        \dim_\C G_\C
        <
        \dim_\R G
        \quad \textrm{and} \quad
        \pi_1(G_\C)
        =
        \{\E\}
        \subsetneq
        D
        =
        \pi_1(G).
    \end{equation}
\end{example}

\begin{example}[Complexifying complex groups,
{\cite[Exercise~15.2.3]{hilgert.neeb:2012a}}]
    \label{ex:ComplexificationOfComplex}
    \index{Universal Complexification!Complex groups}
    \index{Complexification!Complex groups}
    Let $\widetilde{G}$ be a connected and simply connected -- an assumption we shall
    drop in a second step -- \emph{complex} Lie group. The following provides a model
    of
    $\widetilde{G}$ within $\widetilde{G}_\C$ as the diagonal. That is, we define
    \begin{equation}
        \label{eq:ComplexificationOfComplex}
        \tilde{\eta}
        \colon
        \widetilde{G} \longrightarrow \widetilde{G} \times \cc{\widetilde{G}}, \quad
        \tilde{\eta}(g)
        \coloneqq
        (g,g),
    \end{equation}
    where $\cc{\widetilde{G}}$ is the complex conjugated group to $\tilde{G}$. The idea is to solve the extension problem posed by the universal property on the level of the Lie algebras and then to integrate by means of Theorem~\ref{thm:Lie2}. We thus begin with some preliminary considerations pertaining the linear algebraic situation.

    Note that one should not expect complexification to be idempotent. This already
    occurs for vector spaces, as there are $\R$-linear mappings that fail to be
    $\C$-linear, e.g. the complex conjugation $\operatorname{cc} \colon \C
    \longrightarrow \C$, see again our considerations surrounding
    Proposition~\ref{prop:InducedMorphismReal}.
    Here, the two instances of $\C$ play different roles. We consider the domain of
    $\operatorname{cc}$ as the \emph{real} group we want to complexify and the image
    as a generic \emph{complex} group.

    We double the dimension and consider $\C_\C \coloneqq \C \oplus \cc{\C}$, where
    $\cc{\C} = \C$ as an additive group endowed with the opposite scalar multiplication
    $z \odot \cc{\xi} \coloneqq \cc{\cc{z}\xi}$ for scalars $z \in \C$ and vectors~$\cc{\xi} \in
    \C$ as we have done in the proof of Proposition~\ref{prop:ComplexConjugation}. Then
    we may embed~$\C$ as the diagonal
    \begin{equation}
        \bigl\{
            (\xi,\xi)
            \colon
            \xi \in \C
        \bigr\}
        \subseteq
        \C_\C,
    \end{equation}
    which yields a candidate for the universal complexification, namely $\eta(\xi) \coloneqq (\xi,\xi)$. And indeed, for the specific example of $\operatorname{cc}$, setting
    \begin{equation}
        \operatorname{cc}_\C
        \colon
        \C_\C
        \longrightarrow
        \C, \quad
        \operatorname{cc}_\C(\xi,\eta)
        \coloneqq
        \cc{\eta}
    \end{equation}
    yields a $\C$-linear and thus holomorphic extension of $\operatorname{cc}$. Analogously,  the holomorphic extension a of $\C$-linear mapping $\phi \colon \C \longrightarrow \C$ is simply given by $\phi_\C(\xi,\eta) \coloneqq \phi(\xi)$.

    For general $\R$-linear maps, i.e. maps that are neither $\C$-linear nor
    $\C$-antilinear, the situation is more complicated. We denote the spaces of
    $\R$-linear, $\C$-linear, respectively $\C$-antilinear mappings between two complex
    vector spaces $V$ and $W$ by \gls{Linear},  \gls{CLinear} and \gls{CALinear}.
    Then
    \begin{equation}
        \Linear(V,W)
        \cong
        \Linear_\C(V,W)
        \oplus
        \cc{\Linear}_\C(V,W)
    \end{equation}
    forms a direct decomposition with projections
    \begin{align}
        &\pi_\C
        \colon
        \Linear(V,W)
        \longrightarrow
        \Linear_\C(V,W), \quad
        \bigl(
            \pi_\C \phi
        \bigr)(v)
        \coloneqq
        \frac{1}{2}
        \bigl(
        \phi(v)
        -
        \I
        \phi(\I v)
        \bigr), \\
        &\cc{\pi}_\C
        \colon
        \Linear(V,W)
        \longrightarrow
        \cc{\Linear}_\C(V,W), \quad
        \bigl(
            \cc{\pi}_\C \phi
        \bigr)(v)
        \coloneqq
        \frac{1}{2}
        \bigl(
        \phi(v)
        +
        \I
        \phi(\I v)
        \bigr).
    \end{align}
    Note that we use that both $V$ and $W$ are complex vector spaces, so the
    multiplications with the imaginary unit are well defined both within and outside of
    $\phi$. By construction, we moreover have $\pi_\C + \cc{\pi}_\C = \id_{\Linear(V,W)}$,
    \begin{equation}
        \bigl(
        \pi_\C
        \circ
        \cc{\pi}_\C
        \bigr)
        \phi(v)
        =
        \frac{1}{2}
        \bigl(
        (\cc{\pi}_\C\phi)(v)
        -
        \I
        (\cc{\pi}_\C \phi)(\I v)
        \bigr)
        =
        \frac{1}{4}
        \bigl(
        (\phi(v) + \I \phi(\I v))
        -
        \I
        (\phi(\I v) + \I \phi(- v))
        \bigr)
        =
        0
    \end{equation}
    and
    \begin{equation}
        \bigl(
        \cc{\pi}_\C
        \circ
        \pi_\C
        \bigr)
        \phi(v)
        =
        \frac{1}{2}
        \bigl(
        (\pi_\C\phi)(v)
        +
        \I
        (\pi_\C \phi)(\I v)
        \bigr)
        =
        \frac{1}{4}
        \bigl(
        (\phi(v) - \I \phi(\I v))
        +
        \I
        (\phi(\I v) - \I \phi(- v))
        \bigr)
        =
        0.
    \end{equation}
    We return to the more specific problem of extending an $\R$-linear mapping
    \begin{equation}
        \phi
        \colon
        V
        \longrightarrow
        W,
    \end{equation}
    between complex vector spaces $V$ and $W$ to a $\C$-linear mapping defined on
    $V \oplus \cc{V}$. Our considerations yield the extension
    \begin{equation}
        \label{eq:ComplexificationOfComplexMorphism}
        \phi_\C
        \colon
        V \oplus \cc{V}
        \longrightarrow
        W, \quad
        \phi_\C(\xi, \chi)
        \coloneqq
        \pi_\C \phi(\xi, \chi)
        =
        \frac{1}{2}
        \bigl(
        \phi(\xi+\chi)
        -
        \I
        \phi(\I \xi - \I \chi)
        \bigr),
    \end{equation}
    where the final minus sign arises from the conjugated scalar multiplication $\odot$ on $\cc{V}$, which is defined exactly as in the case of $\cc{\C}$.

    Let now $\Phi \colon \widetilde{G} \longrightarrow H$ be a group morphism into
    another complex Lie group $H$. We study its $\R$-linear tangent map $\phi
    \coloneqq T_\E \Phi \colon \liealg{g} \longrightarrow \liealg{h}$ and define $\phi_\C$
    by \eqref{eq:ComplexificationOfComplexMorphism}. By our preliminary
    considerations, this yields a $\C$-linear mapping with
    \begin{equation}
        \label{eq:ComplexificationOfComplexTangent}
        \bigl(
        \phi_\C \circ T_\E \tilde{\eta}
        \bigr)(\xi)
        =
        \frac{1}{2}
        \bigl(
        \phi(\xi+\xi)
        -
        \I
        \phi(\I \xi - \I \xi)
        \bigr)
        =
        \phi(\xi)
    \end{equation}
    for $\xi \in \liealg{g}$ by $\R$-linearity of $\phi$. It remains to check that $\phi_\C$
    respects the Lie bracket, which by $\C$-linearity of $\phi$ and $\C$-bilinearity of the
    bracket suffices to be done on generators of the form $(\xi,0)$ and $(0,\chi)$ with
    $\xi,\chi \in \liealg{g}$. To this end, we first note that
    \begin{equation}
        \phi_\C
        \bigl(
        [(\xi, 0), (0, \chi)]
        \bigr)
        =
        \phi_\C
        \bigl(
        [\xi,0], [0,\chi]
        \bigr)
        =
        \phi_\C(0)
        =
        0
    \end{equation}
    as well as
    \begin{align}
        \big[
        \phi_\C(\xi,0),
        \phi_\C(0,\chi)
        \big]
        &=
        \frac{1}{4}
        \big[
        \phi(\xi) - \I \phi(\I \xi),
        \phi(\chi) + \I \phi(\I \chi)
        \big] \\
        &=
        \frac{1}{4}
        \bigl(
        [\phi(\xi), \phi(\chi)]
        +
        [\phi(\I \xi), \phi(\I \chi)]
        +
        \I
        (
        -
        [\phi(\I \xi),\phi(\chi)]
        +
        [\phi(\xi), \phi(\I \chi)]
        )
        \bigr) \\
        &=
        \frac{1}{4}
        \bigl(
        \phi([\xi, \phi])
        +
        \phi([\I \xi, \I \chi])
        +
        \I
        (
        -
        \phi
        ([\I \xi,\chi])
        +
        \phi([\xi, \I \chi])
        )
        \bigr) \\
        &=
        0,
    \end{align}
    where we have used that $\phi$ is a Lie algebra morphism and that the brackets of $\liealg{g}$ and $\liealg{h}$ are $\C$-bilinear. Similarly, we obtain
    \begin{align}
        \phi_\C
        \bigl(
        [(\xi, 0), (\chi,0)]
        \bigr)
        &=
        \phi_\C
        \bigl(
        ([\xi,\chi],0)
        \bigr) \\
        &=
        \frac{1}{2}
        \bigl(
        \phi([\xi,\chi])
        -
        \I
        \phi(\I[\xi,\chi])
        \bigr) \\
        &=
        \frac{1}{4}
        \bigl(
        2[\phi(\xi),\phi(\chi)]
        -
        \I
        [\phi(\I \xi),\phi(\chi)]
        -
        \I
        [\phi(\xi),\phi(\I \chi)]
        \bigr) \\
        &=
        \frac{1}{4}
        \big[
        \phi(\xi) - \I \phi(\I \xi),
        \phi(\chi) - \I \phi(\I \chi)
        \big] \\
        &=
        \big[
        \phi_\C(\xi,0),
        \phi_\C(\chi,0)
        \big].
    \end{align}
    The antisymmetry of the Lie bracket yields the remaining combinations, so we have shown that $\phi_\C \colon \liealg{g} \oplus \cc{\liealg{g}} \longrightarrow \liealg{h}$ is a $\C$-linear Lie algebra morphism.

    Invoking Lie's second Theorem~\ref{thm:Lie2}, there thus exists a unique
    holomorphic Lie group
    morphism
    \begin{equation}
        \Phi_\C \colon \widetilde{G} \longrightarrow H
    \end{equation}
    such that $T_\E \Phi_\C = \phi_\C$. By connectedness of $\widetilde{G}$ and
    \eqref{eq:ComplexificationOfComplexTangent}, we have $\Phi_\C \circ \tilde{\eta} =
    \Phi$. As in the proof of Theorem~\ref{thm:UniversalComplexification}, $\Phi_\C$
    is uniquely determined by its tangent map $\phi_\C$, which in turn is the unique
    $\C$-linear extension of $\phi$. That is, $\Phi_\C$ is unique. Consequently
    \eqref{eq:ComplexificationOfComplex} is a universal complexification of
    $\widetilde{G}$.

    Let now $G$ a be connected, but not necessarily simply connected Lie group. Using
    what we have already shown, we may mimic our proof of
    Theorem~\ref{thm:UniversalComplexification} to obtain a description of the universal
    complexification of $G$ as a quotient of the universal complexification of its universal
    covering group $p \colon \widetilde{G} \longrightarrow G$. In the notation of
    Theorem~\ref{thm:UniversalComplexification}, we get that~$\widetilde{\eta}(g) =
    (g,g)$
    and thus
    \begin{equation}
        \widetilde{\eta}
        \bigl(
            \pi_1(G)
        \bigr)
        =
        \bigl\{
            (g,g)
            \colon
            g \in \pi_1(G)
        \bigr\}
        \subseteq
        \widetilde{G}
        \times
        \cc{\widetilde{G}}
    \end{equation}
    is already a complex subgroup. Consequently,
    \begin{equation}
        \label{eq:UniversalComplexificationOfComplex}
        \eta
        \colon
        G
        \longrightarrow
        \bigl(
            \widetilde{G} \times \cc{\widetilde{G}}
        \bigr)
        \big/
        \{(g,g) \colon g \in \pi_1(G)\}, \quad
        \eta
        \bigl(
            p(h)
        \bigr)
        =
        \bigl[
            (h,h)
        \bigr]
    \end{equation}
    is a universal complexification of $G$ by
    \eqref{eq:ComplexificationFromUniversalCovering}. One may thus consider repeated
    universal complexification as a tower of complex Lie groups obtained by starting with
    a complex group, such as the complexification of some Lie group. Note that the
    mapping
    $\tilde{\eta}$ from~\eqref{eq:ComplexificationOfComplex} is injective, so we may
    apply
    Proposition~\ref{prop:ComplexificationVsCovering}. Indeed, taking another look
    at~\eqref{eq:UniversalComplexificationOfComplex}, we see that
    \begin{equation*}
        \pi_1(G_\C)
        =
        \bigl\{
            (g,g)
            \colon
            g \in \pi_1(G)
        \bigr\}
        \cong
        \pi_1(G)
    \end{equation*}
    in accordance with \eqref{eq:ComplexificationVsCovering}, as direct products of
    simply connected groups are simply connected. These considerations match with the
    form of the homogeneous space \eqref{eq:Omega}.
\end{example}

\begin{example}[The tower of cylinders]
    \index{Universal Complexification!Tower of cylinders}
    \index{Complexification!Tower of cylinders}
    Let $G \coloneqq (\C^\times, \cdot)$ be the cylinder, which is the universal complexification of the circle group $\mathbb{S}^1$ by Example~\ref{ex:CircleGroupUniversalComplexification}. Its universal covering projection is given by the complex exponential function
    \begin{equation}
        \exp
        \colon
        \C \longrightarrow \C^\times
    \end{equation}
    and $\pi_1(G) \cong \{2\pi \I k \colon k \in \Z\}$.
    By \eqref{ex:ComplexificationOfComplex}, the universal complexification is given by
    \begin{equation}
        \eta
        \colon
        \C^\times
        \longrightarrow
        G_\C
        \coloneqq
        \bigl(
            \C \times \cc{\C}
        \bigr)
        \big/
        \{(2\pi \I k, 2\pi \I k) \colon k \in \Z\}
        , \quad
        \eta
        \bigl(
        \exp(z)
        \bigr)
        =
        \bigl[(z,z)\bigr]
    \end{equation}
    Two points $(z_1, w_1), (z_2, w_2) \in G_\C$ are thus identified if and only if
    \begin{equation}
        z_1 - z_2
        =
        2\pi \I k
        =
        w_1 - w_2
        \qquad
        \textrm{for some }
        k \in \Z.
    \end{equation}
    That is, $\RE(z_1) = \RE(z_2)$, $\RE(w_1) = \RE(w_2)$ and
    \begin{equation}
        \bigl(
            \IM(z_1), \IM(w_1)
        \bigr)
        =
        \bigl(
            \IM(w_1) + 2\pi \I k,
            \IM(w_2) + 2\pi \I k
       \bigr).
    \end{equation}
    Grouping together real and imaginary parts,
    $G_\C$ is thus geometrically given by $\C$ (the real parts) times a cylinder (the
    imaginary parts), which matches nicely with
    Proposition~\ref{prop:ComplexificationVsCovering}. Using the compatibility of direct
    products from Example~\ref{ex:ComplexificationOfProducts} with universal
    complexification, we see that the $n$-fold complexification of $G = \C^\times$ is
    given by
    \begin{equation}
        \C^{2n-1}
        \times
        \C^\times
        \cong
        \widetilde{G}^{2n-1}
        \times
        G.
    \end{equation}
\end{example}

\section{Extension and Restriction of Entire Functions}
\label{sec:ExtensionAndRestriction}
\epigraph{Wir müssen wissen. \\ Wir werden wissen.}{David Hilbert}
% !TeX root = ../Dissertation.tex
The purpose of this section is to establish a one-to-one correspondence between entire
functions on $G$ and holomorphic ones on $G_\C$ for connected linear Lie groups.
These are
the problems of \emph{extension} and \emph{restriction}. We begin with the former,
which we approach by treating two somewhat redundant but nevertheless instructive
special cases first: Simply connected nilpotent and compact groups.

Understanding these in many regards simpler situations first lets us develop machinery
suitable for the proof of Theorem~\ref{thm:Extension} and its generalization
Proposition~\ref{prop:ExtensionGeneral}, which then subsume
Theorem~\ref{thm:NilpotentExtension} and Theorem~\ref{thm:CompactExtension}.
Starting out, we are guided by the general assumption that
the universal complexification morphism $\eta \colon G \longrightarrow G_\C$ is
injective. This means that $G$ constitutes a closed Lie subgroup of $G_\C$ by virtue
of Corollary~\ref{cor:ImageOfEtaSubmanifold}. Having a set inclusion is of
course necessary to speak of ``extensions'' of functions in a strict way. However, we
shall convince ourselves shortly that pullback with $\eta$ provides a natural
replacement in the general situation. Revisiting the particular case of representative
functions motivates the indicated program quite nicely.
\begin{example}[Representative functions]
    \index{Representative functions!Extension}
    \index{Matrix elements!Extension}
    \label{ex:RepresentativeExtension}
    Let $\pi \colon G \longrightarrow \GL(V)$ be a representation of a connected Lie
    group $G$ on a finite dimensional vector space $V$. By
    Theorem~\ref{thm:RepresentativeFunctions}, every choice of $v \in V$ and $\varphi
    \in V'$ induces an entire function
    \begin{equation}
        \gls{RepresentativeFunction}
        \colon
        G \longrightarrow \C, \quad
        \pi_{v,\varphi}(g)
        \coloneqq
        \varphi
        \bigl(
            \pi(g)v
        \bigr),
    \end{equation}
    which is even $R$-entire for $0 \le R < 1$ by
    Theorem~\ref{thm:RepresentativeFunctions}. Complexifying $V$, we may extend
    $\pi$ to a Lie group morphism taking values in the complex Lie group $\GL(V_\C)$.
    We denote this extension also as $\pi$. By the universal property of $G_\C$,
    there in turn exists a unique holomorphic representation $\pi_\C \colon G_\C
    \longrightarrow \GL(V_\C)$ with $\pi_\C \circ \eta = \pi$. Setting
    \begin{equation}
        \Pi_{v,\varphi}
        \colon
        G_\C \longrightarrow \C, \quad
        \Pi_{v,\varphi}(g)
        \coloneqq
        \varphi
        \bigl(
            \pi_\C(g)v
        \bigr),
    \end{equation}
    where we also extend the linear mapping $\varphi$ to $V_\C$, thus yields a
    holomorphic function with
    \begin{equation}
        \label{eq:ExtensionRepresentative}
        \Pi_{v,\varphi}
        \circ
        \eta
        =
        \pi_{v,\varphi}.
    \end{equation}
    Note that asking for the condition \eqref{eq:ExtensionRepresentative} makes sense
    regardless of injectivity of $\eta$ and thus yields the aforementioned natural
    generalization of the extension problem. In Example~\ref{ex:SpecialLinear}, we will
    see that the failure of injectivity of $\eta$ then enforces certain periodicity
    conditions for the extensions. We shall encounter a quite drastic instance of this
    phenomenon in Example~\ref{ex:AutomaticPeriodicity}.
\end{example}

The first class of geometrically quite simple groups we want to investigate are
connected and simply connected groups, which are in addition nilpotent. A concrete
example is then provided by the Heisenberg groups, which we will investigate in
Example~\ref{ex:Schroedinger} and a comprehensive discussion of which can be found
in the form of the monograph \cite{binz.pods:2008a}. By
\cite[Cor.~11.2.7]{hilgert.neeb:2012a} these assumptions imply that the Lie
exponential
\begin{equation}
    \exp
    \colon
    \liealg{g}
    \longrightarrow
    G
\end{equation}
is a diffeomorphism, which is the property we are interested in. It turns out that
this feature, along with many others, passes to the universal complexification.
\begin{proposition}
    \index{Nilpotent!Complexification}
    \index{Complexification!Nilpotent}
    \index{Universal Complexification!Nilpotent}
    \label{prop:ComplexificationNilpotent}
    Let $G$ be a nilpotent connected and simply connected Lie group with universal complexification $(G_\C, \eta)$.
    \begin{propositionlist}
        \item \label{item:ComplexificationNilpotentSimplyConnected}
        The universal complexification $G_\C$ is connected and simply connected.
        \item \label{item:ComplexificationNilpotentInjective}
        The mapping $\eta \colon G \longrightarrow G_\C$ is injective.
        \item \label{item:ComplexificationNilpotentLieAlgebra}
        We have $\LieAlg(G_\C) \cong \liealg{g}_\C$ as Lie algebras.
        \item \label{item:ComplexificationNilpotent}
        The universal complexification $G_\C$ is a nilpotent group.
        \item \label{item:ComplexificationNilpotentExpDiffeo}
        The exponential mapping $\exp_{G_\C} \colon \liealg{g}_\C \longrightarrow G_\C$ is a diffeomorphism.
    \end{propositionlist}
\end{proposition}
\begin{proof}
    By Corollary~\ref{cor:UniversalComplexificationProperties},
    \ref{item:UniversalComplexificationSimplyConnected2}, the claims
    \ref{item:ComplexificationNilpotentSimplyConnected} and
    \ref{item:ComplexificationNilpotentLieAlgebra} hold and $T_\E \eta$ is injective. But
    then
    \begin{equation}
        \eta
        =
        \exp_{G_\C}
        \circ
        T_\E \eta
        \circ
        \log_G
    \end{equation}
    is injective itself as a composition of injective maps, where we write $\log_G$ for the
    inverse of $\exp_G$. We have thus
    shown~\ref{item:ComplexificationNilpotentInjective}. To prove
    \ref{item:ComplexificationNilpotent}, we check that $\liealg{g}_\C$ is nilpotent, which
    by \cite[Thm.~11.2.5]{hilgert.neeb:2012a} then implies that $G_\C$ is. Conversely, we
    may use that $\liealg{g}$ is nilpotent by assumption. Invoking
    Engel's\footnote{Friedrich Engel (1861-1941) was a German
    mathematician. Together with Sophus Lie, he wrote the series of books ``Theory of
    transformation groups'' \cite{engel.lie:1888a}, which may be viewed as the
    inception of modern Lie theory.
    Not to be confused with the philosopher Friedrich Engels.} Theorem
    \cite[Thm.~5.2.8]{hilgert.neeb:2012a} it suffices to check that the linear mappings
    $\ad_\xi \colon \liealg{g}_\C \longrightarrow \liealg{g}_\C$ are nilpotent for all $\xi \in
    \liealg{g}_\C$. That is, for every $\xi \in \liealg{g}_\C$ there exists some integer $k \in
    \N$ such that $\ad_\xi^k = 0$. But
    this is clear by either \eqref{eq:ComplexificationLieBracketTensor}
    or~\eqref{eq:ComplexificationLieBracketSum} and the nilpotency of $\liealg{g}$.
    Finally,
    \ref{item:ComplexificationNilpotentExpDiffeo} follows from
    \cite[Cor.~11.2.7]{hilgert.neeb:2012a} in the same manner as it did for $G$ itself.
\end{proof}

The upshot is that $G_\C$ is, just like $G$, geometrically a vector space.
\begin{theorem}[Extension from nilpotent groups]
    \label{thm:NilpotentExtension}
    \index{Nilpotent!Extension}
    \index{Extension!Nilpotent}
    Let $G$ be a connected, simply connected and nilpotent Lie group with universal
    complexification $(G_\C, \eta)$. Moreover, fix a parameter $R \ge 0$ and $f \in
    \Entire_R(G)$. Then there exists a unique
    \begin{equation}
        F
        \in
        \Holomorphic(G_\C)
        \cap
        \Entire_R(G_\C)
        \qquad
        \textrm{such that}
        \qquad
        \eta^* F = f.
    \end{equation}
    More precisely,
    \begin{equation}
        \label{eq:NilpotentExtension}
        F
        \bigl(
            \eta(g) \exp \xi
        \bigr)
        =
        \Taylor_f
        \bigl(
        \underline{\xi}; g
        \bigr)
        =
        \Taylor_F
        \bigl(
        \underline{\xi}; \eta(g)
        \bigr)
    \end{equation}
    holds for all $g \in G$ and $\xi \in \liealg{g}_\C$.
\end{theorem}
\begin{proof}
    By Proposition~\ref{prop:ComplexificationNilpotent},
    \ref{item:ComplexificationNilpotentExpDiffeo}, the exponential map $\exp_{G_\C}
    \colon \liealg{g}_\C \longrightarrow G_\C$ constitutes a global chart of $G_\C$ with
    $\exp_{G_\C}(\liealg{g}) = \eta(G)$. Moreover, the formal Lie-Taylor series
    $\Taylor_f(\underline{\xi}; \E)$ converges for all $\xi \in \liealg{g}_\C$ by
    Theorem~\ref{thm:LieTaylor}. We define $F \colon G_\C \longrightarrow \C$ by
    \begin{equation}
        F
        \bigl(
        \exp_{G_\C}(\xi)
        \bigr)
        \coloneqq
        \Taylor_f
        \bigl(
        \underline{\xi}; \E
        \bigr)
        =
        \sum_{k=0}^{\infty}
        \frac{1}{k!}
        \sum_{\alpha \in \N_n^k}
        \bigl(
        \Lie(\alpha)
        f
        \bigr)
        (\E)
        \cdot
        \underline{\xi}^\alpha,
    \end{equation}
    which is holomorphic with respect to $\xi$ and thus holomorphic on $G_\C$. For $g =
    \exp_G(\xi) \in G$, we get
    \begin{equation}
        F
        \bigl(
        \eta(g)
        \bigr)
        =
        F
        \bigl(
        \exp_{G_\C}(\xi)
        \bigr)
        =
        \Taylor_f
        \bigl(
        \underline{\xi}; \E
        \bigr)
        =
        f
        \bigl(
        \exp_G(\xi)
        \bigr)
        =
        f(g).
    \end{equation}
    Finally, \eqref{eq:NilpotentExtension} holds by construction and the fact that $\E \in G_\C \cap \eta(G)$.
\end{proof}

\begin{remark}[Nilpotent groups]
    \index{Nilpotent!Group}
    \index{Dixmier, Jacques}
    \index{Saitō, Morihiko}
    To prove Theorem~\ref{thm:NilpotentExtension}, it would suffice to demand that
    \begin{equation}
        \exp_{G_\C} \colon \liealg{g}_\C \longrightarrow G_\C
    \end{equation}
    is a diffeomorphism. By a combination of \cite[Thm.~14.4.3]{hilgert.neeb:2012a}
    and Dixmier-Saito's Theorem\footnote{Morihiko Saitō (born 1961) is a Japenese
    algebraic geometer and algebraic analyst working at Kyoto University. Dixmier and
    Saitō discovered their theorem independently.}~\cite{dixmier:1957a, saito:1957a}, a
    modern formulation of which can be found in
    \cite[Thm.~14.4.8]{hilgert.neeb:2012a}, this however
    already implies that $G_\C$ and thus in particular any of its subgroups are
    nilpotent. As our general assumption is the injectivity of $\eta \colon G
    \longrightarrow G_\C$, we have thus not lost any generality by dealing with the
    nilpotent situation directly.
\end{remark}

\begin{figure}[H]
    \begin{center}
        \includegraphics[width=10cm]{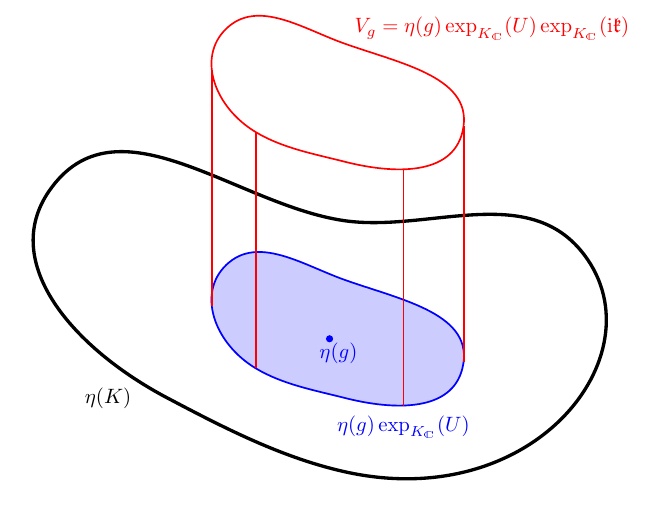}
        \caption{Schematic picture of $V_g$}
        \label{fig:ExtensionFromCompact}
    \end{center}
\end{figure}

\begin{figure}[h]
    \begin{center}
        \includegraphics[width=12cm]{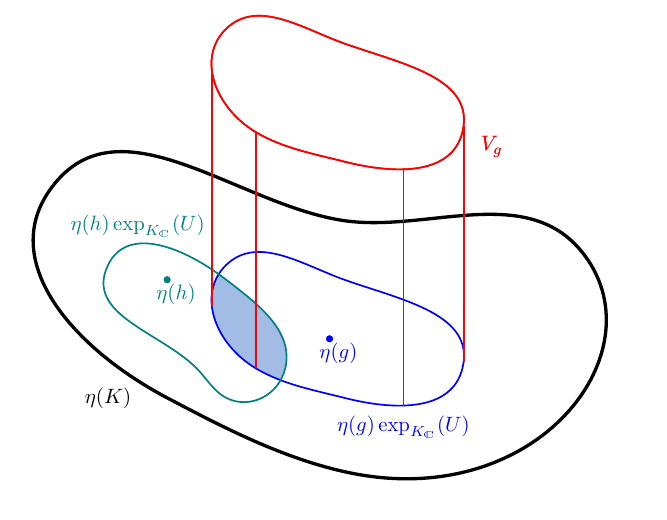}
        \caption{Schematic picture of $V_g \cap V_h$}
        \label{fig:ExtensionFromCompactIntersection}
    \end{center}
\end{figure}

Next, we investigate connected compact Lie groups $K$ by applying the
techniques from the nilpotent case locally. By
\cite[Prop.~15.2.1(c)]{hilgert.neeb:2012a}, the polar decomposition
\begin{equation}
    \label{eq:CompactPolar}
    \Phi
    \colon
    K \times \liealg{k}
    \longrightarrow
    K_\C, \quad
    \Phi(g,\xi)
    \coloneqq
    \eta(g) \cdot \exp_{K_\C}(\I \xi)
\end{equation}
is a diffeomorphism. Geometrically, one may thus view $K_\C$ as a cylinder over
$\eta(K) \cong K$, see also \cite[Prop.~15.2.1(a)]{hilgert.neeb:2012a}. In particular, this
implies $\pi_1(K_\C) \cong \pi_1(K)$. Consequently, all topological obstructions for
analytic continuation are already present within $K$.
\begin{theorem}[Extension from compact groups]
    \label{thm:CompactExtension}
    \index{Extension!Compact}
    Let $K$ be a connected as well as compact Lie group with universal
    complexification $(K_\C, \eta)$, $R \ge 0$ and $f \in \Entire_R(K)$. Then there
    exists a unique
    \begin{equation}
        F
        \in
        \Holomorphic(K_\C)
        \cap
        \Entire_R(K_\C)
        \qquad
        \textrm{such that}
        \qquad
        \eta^* F = f.
    \end{equation}
    More precisely,
    \begin{equation}
        \label{eq:CompactExtension}
        F
        \bigl(
            \eta(g) \exp_{K}(\chi) \exp_{K_\C}(\I \xi)
        \bigr)
        =
        \Taylor_f
        \bigl(
        \underline{\chi + \I \xi}; g
        \bigr)
        =
        \Taylor_F
        \bigl(
            \underline{\chi + \I \xi}; \eta(g)
        \bigr)
    \end{equation}
    holds for all $g \in K$ and $\chi, \xi \in \liealg{k}$.
\end{theorem}
\begin{proof}
    By \cite[Prop.~9.2.5]{hilgert.neeb:2012a} there exists a zero neighbourhood
    $U \subseteq \liealg{k}$ such that $\exp_K \at{U}$ is a diffeomorphism onto
    its image. Combining this with \eqref{eq:CompactPolar} and the injectivity of
    $\eta$ yields a holomorphic atlas $(z_g, V_g)_{g \in K}$ of $K_\field{C}$
    with $V_g \coloneqq \eta(g) \exp_{K_\C}(U)\exp_{K_\C}(\I \liealg{k})$ and
    \begin{equation}
        z_g^{-1}(\chi, \xi)
        \coloneqq
        \eta(g) \exp_{K_\C}(\chi) \exp_{K_\C}(\I \xi).
    \end{equation}
    A schematic picture of $V_g$ can be found in
    Figure~\ref{fig:ExtensionFromCompact}, where only a finite section of the
    vertical lines is pictured and $\exp_{K_\C}(\I \xi)$ identified with $\I \xi$
    by means of \eqref{eq:CompactPolar} for simplicity.
    By Theorem~\ref{thm:Symmetries}, \ref{item:LieTaylorAbsoluteConvergence}
    the Lie-Taylor series $\Taylor_f(\underline{\xi}; g)$ of $f$ from
    \eqref{eq:LieTaylorFormal} around $g$ converges for all $\xi \in
    \liealg{k}_\field{C}$. For every $g \in G$, we define
    functions $F_g \colon V_g \longrightarrow \C$ by
    \begin{equation}
        F_g
        \bigl(
            \eta(g) \exp_{K_\C}(\chi) \exp_{K_\C}(\I \xi)
        \bigr)
        \coloneqq
        \Taylor_f
        \bigl(
            \underline{\chi + \I \xi}; g
        \bigr)
        =
        \sum_{k=0}^{\infty}
        \frac{1}{k!}
        \sum_{\alpha \in \N_n^k}
        \bigl(
        \Lie(\alpha)
        f
        \bigr)
        (g)
        \cdot
        (\underline{\chi + \I \xi})^\alpha,
    \end{equation}
    which is holomorphic with respect to the complex variable $\chi + \I \xi$. Moreover, we
    note that
    \begin{equation}
       V_g
       \cap
       \eta(K)
       =
       \eta(g)
       \exp_{K_\C}(U)
       =
       \eta
       \bigl(
        g \exp_K(U)
       \bigr)
    \end{equation}
    and thus
    \begin{equation}
        \eta^* F_g
        \bigl(
            g \exp_K \chi
        \bigr)
        =
        F_g
        \bigl(
            \eta(g)
            \exp_{K_\C} \chi
        \bigr)
        =
        T_f(\underline{\chi};g)
        =
        f
        \bigl(
            g \exp_K \chi
        \bigr).
    \end{equation}
    Consequently, $F_g$ constitutes a holomorphic extension of $f \at{g \exp(U)}$.
    Moreover, if $V_g \cap V_h \neq \emptyset$ for some $g,h \in K$ as sketched in
    Figure~\ref{fig:ExtensionFromCompactIntersection}, then
    \begin{equation}
        F_g
        \at[\Big]
        {\eta(g\exp U \cap h \exp U)}
        =
        f
        \at[\Big]
        {g \exp U \cap h \exp U}
        =
        F_h
        \at[\Big]
        {\eta(g\exp U \cap h \exp U)}
    \end{equation}
    and thus $F_g = F_h$ by virtue of the identity principle from
    Corollary~\ref{cor:IdentityPrinciple}. Here we use that
    $\liealg{k}_\C$ is the $\C$-span of $\liealg{k}$. Consequently and as $\bigcup_{g \in
    G} V_g = K_\C$, the family $(F_g)_{g \in G}$ glues to a globally defined holomorphic
    function $F$ on $K_\C$ fulfilling $\eta^* F = f$ and \eqref{eq:CompactExtension}.
    Once again invoking the identity principle from
    Corollary~\ref{cor:IdentityPrinciple} establishes the uniqueness of $F$. Finally,
    \eqref{eq:CompactExtension} implies that the Lie-Taylor coefficients of $F$ at
    $\eta(g)$ simply coincide with those of $f$ at $g$ for all $g \in G$. This means that $F
    \in \Entire_R(K_\C)$ by virtue of $\E_{K_\C} = \E_K$.
\end{proof}

Notably, using the notation of Theorem~\ref{thm:CompactExtension} and considering
\begin{equation}
    \widetilde{F}
    \bigl(
        \eta(g) \exp_{K_\C}(\chi) \exp_{K_\C}(\I \xi)
    \bigr)
    \coloneqq
    \Taylor_f
    \bigl(
        \underline{\I \xi}; g \exp_K(\chi)
    \bigr)
\end{equation}
also yields an entire function $\widetilde{F} \in \Entire_R(K_\C)$ with $\eta^*
\widetilde{F} = f$. However, $\widetilde{F}$ fails to be holomorphic in general. It
coincides with \eqref{eq:CompactExtension} whenever $[\chi, \xi] = 0$, as in this case
\begin{equation}
    \exp_{K_\C}(\chi)
    \exp_{K_\C}(\I \xi)
    =
    \exp_{K_\C}
    (\chi + \I \xi)
\end{equation}
and thus by Theorem~\ref{thm:LieTaylor} and \eqref{eq:CompactExtension}
\begin{align}
    \Taylor_f
    \bigl(
        \underline{\I \xi}; g \exp(\chi)
    \bigr)
    &=
    \Taylor_F
    \bigl(
        \underline{\chi + \I \xi}; \eta(g)
    \bigr) \\
    &=
    F
    \bigl(
        \eta(g)
        \exp_{K_\C}(\chi + \I \xi)
    \bigr) \\
    &=
    F
    \bigl(
        \eta(g)
        \exp_{K_\C}(\chi)
        \exp_{K_\C}(\I \xi)
    \bigr) \\
    &=
    F
    \bigl(
        \eta(g \exp_K \chi)
        \exp_{K_\C}(\I \xi)
    \bigr) \\
    &=
    \sum_{k=0}^{\infty}
    \frac{1}{k!}
    \sum_{\alpha \in \N_n^k}
    \bigl(
    \Lie(\alpha)
    F
    \bigr)
    \bigl(
        \eta(g \exp_K \chi)
    \bigr)
    \cdot
    (\underline{\I \xi})^\alpha \\
    &=
    \sum_{k=0}^{\infty}
    \frac{1}{k!}
    \sum_{\alpha \in \N_n^k}
    \bigl(
    \Lie(\alpha)
    f
    \bigr)
    (g \exp_K \chi)
    \cdot
    (\underline{\I \xi})^\alpha \\
    &=
    \Taylor_f
    \bigl(
        \underline{\I \xi}; g \exp_K \chi
    \bigr) \\
    &=
    \widetilde{F}
    \bigl(
        \eta(g) \exp_{K_\C}(\chi) \exp_{K_\C}(\I \xi)
    \bigr).
\end{align}

This completes our preliminary -- and ultimately redundant -- considerations on
extension. Guided by these particular cases, we will now first use the Lie-Taylor series
\eqref{eq:LieTaylor} to obtain a holomorphic extension of a given entire function around
the group unit. In a second step, we then analytically continue along any path in $G_\C$
using essentially Steiner\footnote{Jakob Steiner (1796-1863) was a Swiss synthetic
geometer infamous for never preparing his lectures beforehand, and then suffering the
consequences. His treatment of geometry was entirely algebraic, as he harboured a
deep disdain for analysis and, in particular, coordinates.} chains of balls within charts.
Finally, our considerations from
Section~\ref{sec:UniversalComplexification} will allow us to turn these local statements
into global ones, which may be seen as an instance of the Oka principle, see again
Remark~\ref{rem:Oka}: In the absence of topological obstructions, analytic
continuation is always possible.

We slightly alter our notation for the original function and its extension from $f$
and $F$ to $\phi$ and $\Phi$, respectively, as we shall work with the Lie-Taylor
Majorant $\Majorant$ in the sequel. Before putting our plan into practice, we
recall the notion of analytic continuations along paths, which notably need not
be differentiable but only continuous. Indeed, let~$M$ and~$N$ be complex
manifolds and $\gamma$ be a path within $M$. We say that a holomorphic
mapping~$g_0 \colon U_0 \longrightarrow N$, where $U_0 \subseteq M$ is a
connected
open neighbourhood of $\gamma(0)$, may be analytically continued along
$\gamma$ if there are $0 = t_0 < t_1 < \cdots < t_n = 1$, open connected
neighbourhoods $U_j \subseteq M$ of $\gamma(t_j)$ and
holomorphic mappings~$g_j \colon U_j \longrightarrow N$ with
\begin{equation}
    \gamma(t_{j})
    \in
    U_{j-1}
    \qquad
    \textrm{and}
    \qquad
    g_{j-1}
    \at[\Big]{U_{j-1} \cap U_{j}}
    =
    g_{j}
    \at[\Big]{U_{j-1} \cap U_{j}}
    \qquad
    \textrm{for }
    j=1,\ldots,n.
\end{equation}
Note that this does not mean that the $f_j$ glue to a globally defined
holomorphic function, as e.g. the principal branch of the holomorphic logarithm
may be analytically continued along any path within $\field{C}^\times$.
\begin{proposition}
    \label{prop:ExtensionAlongPaths}
    \index{Extension!Along paths}
    Let $G$ be a connected Lie group such that the universal complexification morphism
    $\eta \colon G \longrightarrow G_\C$ is injective, and let $\phi \in \Entire_0(G)$.
    \begin{propositionlist}
        \item There exists zero neighbourhood $U \subseteq G_\C$ and a holomorphic
        function
        \begin{equation}
            \phi_0 \colon U \longrightarrow \C
        \end{equation}
        such that $\phi_0 \at{U \cap G} = \phi$.
        \item The function $\phi_0$ may be analytically continued along any path in
        $G_\C$ starting at $\E$.
    \end{propositionlist}
\end{proposition}
\begin{proof}
    Throughout, we use the notation from \eqref{eq:LieTaylorMajorant} and \eqref{eq:LieTaylorMajorantCoefficients}.
    By Theorem~\ref{thm:LieTaylor} we may set $U \coloneqq \exp(\liealg{g}_\C)$ and
    $\phi_0 \coloneqq \Taylor_\phi(\argument; \E)$. Let $\gamma \colon [0,1]
    \longrightarrow G_\C$ be a path starting at the group unit. To show the second
    statement, we first replace $U = \exp(V)$ by a sufficiently small zero neighbourhood
    such that $\exp \colon V \longrightarrow U$ is a diffeomorphism. By continuity
    of~$\gamma$, its trace $\gamma([0,1]) \subseteq G_\C$ is compact. Hence, there are
    $0
    = t_0 < t_1 < t_2 < \cdots < t_k = 1$ with
    \begin{equation}
        \gamma(t_{\ell+1})
        \in
        \gamma(t_\ell)
        \exp(\liealg{g}_\C)
        \qquad
        \textrm{for}
        \quad
        \ell = 1, \ldots, k-1
    \end{equation}
    and
    \begin{equation}
        \gamma\bigl([0,1]\bigr)
        \subseteq
        \bigcup_{\ell=0}^k
        \gamma(t_\ell)
        U
        =
        U
        \cup
        \bigcup_{\ell=1}^k
        \gamma(t_\ell)U.
    \end{equation}
    We prove that the Lie-Taylor series $\Taylor_{\phi_0}(\argument; \gamma(t_1))$ yields
    an entire function on $\C^n$ by generalizing Theorem~\ref{thm:Symmetries},
    \ref{item:TranslationInvariance} to \emph{complex translations}. By construction,
    $\gamma(t_1) \in U$, so there exists a $\xi \in \liealg{g}_\C$ with $\gamma(t_1) =
    \exp \xi$. Recall that $\LieAlg(G_\C) = \liealg{g}_\C$ by construction or
    Lemma~\ref{lem:ComplexificationEtaLocallyInjective} and thus choosing a basis of
    $\liealg{g}$ yields a basis of $\LieAlg(G_\C)$, see
    Corollary~\ref{cor:ImageOfEtaLieAlgebra}. Unwrapping the definitions and plugging in
    \eqref{eq:LieTaylorOneVariable} for $\Lie(\alpha)\phi_0$ leads to
    \begin{align}
        \Majorant_{\phi_0}
        \bigl(\abs{z}; \gamma(t_1)\bigr)
        &=
        \sum_{k=0}^{\infty}
        \frac{1}{k!}
        \sum_{\alpha \in \N_n^k}
        \abs[\Big]
        {
            \bigl(
                \Lie(\alpha)
                \phi_0
            \bigr)(\exp \xi )
        }
        \cdot
        \abs{z}^k \\
        \label{eq:ExtensionProof}
        &\le
        \sum_{k,\ell=0}^{\infty}
        \frac{\abs{z}^k}{k! \ell!}
        \sum_{\alpha \in \N_n^k}
        \abs[\Big]
        {
            \bigl(
            \Lie(\xi)^\ell
            \Lie(\alpha)
            \phi
            \bigr)(\E)
        }.
    \end{align}
    Decomposing $\xi = \xi_1 + \I \xi_2$ with $\xi_1,\xi_2 \in \liealg{g}$
    yields by Proposition~\ref{prop:CauchyRiemann}
    \begin{equation}
        \Lie(\xi)\psi
        =
        \Lie(\xi_1)\psi
        +
        \I
        \cdot
        \Lie(\xi_2)\psi
        \qquad
        \textrm{for all }
        \psi \in \Holomorphic(U).
    \end{equation}
    Plugging in $\psi \coloneqq \Lie(\alpha) \phi_0$ for any $\alpha \in \N_n^k$, we
    estimate
    \begin{align*}
        \abs[\bigg]
        {
            \sum_{\alpha \in \N_n^k}
            \bigl(
            \Lie(\xi)^\ell
            \Lie(\alpha)
            \phi
            \bigr)(\E)
        }
        &\le
        \sum_{\alpha \in \N_n^k}
        \abs[\Big]
        {
            \Bigl(
            \bigl(
                \Lie(\xi_1)
                +
                \I
                \cdot
                \Lie(\xi_2)
            \bigr)^\ell
            \Lie_{X_\alpha}
            \phi
            \Bigr)(\E)
        } \\
        &\le
        \sum_{\alpha \in \N_n^k}
        \sum_{\epsilon_1, \ldots, \epsilon_\ell = 1}^2
        \abs[\Big]
        {
            \Bigl(
                \Lie
                \bigl(
                    \xi_{\epsilon_1}
                    \tensor \cdots \tensor
                    \xi_{\epsilon_\ell}
                \bigr)
                \Lie(\alpha)
                \phi
            \Bigr)(\E)
        } \\
        &=
        \sum_{\alpha \in \N_n^k}
        \sum_{\epsilon_1, \ldots, \epsilon_\ell = 1}^2
        \abs[\big]
        {
            \xi_{\epsilon_1}^{m_1}
            \cdots
            \xi_{\epsilon_\ell}^{m_\ell}
        }
        \cdot
        \abs[\Big]
        {
            \Bigl(
                \Lie
                \bigl(
                    \basis{e}_{m_1}
                    \tensor \cdots \tensor
                    \basis{e}_{m_\ell}
                \bigr)
                \Lie(\alpha)
                \phi
            \Bigr)(\E)
        } \\
        &\le
        \supnorm{2\xi}^\ell
        \sum_{\alpha \in \N_n^k}
        \sum_{m_1, \ldots, m_\ell = 1}^n
        \abs[\Big]
        {
            \Bigl(
            \Lie
            \bigl(
                \basis{e}_{m_1}
                \tensor \cdots \tensor
                \basis{e}_{m_\ell}
            \bigr)
            \Lie(\alpha)
            \phi
            \Bigr)(\E)
        } \\
        &=
        \supnorm{2\xi}^\ell
        \sum_{\alpha \in \N_n^{k+\ell}}
        \abs[\big]
        {
            \bigl(
                \Lie(\alpha)
                \phi
            \bigr)
            (\E)
        } \\
        &=
        \supnorm{2\xi}^\ell
        \cdot
        c_{k+\ell}(\phi)
        \cdot
        (k+\ell)!
    \end{align*}
    for all $k,\ell \in \N_0$. Plugging this into \eqref{eq:ExtensionProof} with $z \in \C$, we arrive at
    \begin{align}
        \abs[\big]
        {
            \Majorant_{\phi_0}
            \bigl(z; \gamma(t_1)\bigr)
        }
        &\le
        \sum_{k,\ell=0}^{\infty}
        \frac{\abs{z}^k}{k! \ell!}
        \abs[\bigg]
        {
            \sum_{\alpha \in \N_n^k}
            \bigl(
            \Lie(\xi)^\ell
            \Lie(\beta)
            \phi
            \bigr)(\E)
        } \\
        &\le
        \sum_{\ell=0}^{\infty}
        \frac{(2 \supnorm{\xi})^\ell}{\ell!}
        \sum_{k=0}^{\infty}
        \frac{(k+\ell)!}{k!}
        \cdot
        c_{k+\ell}(\phi)
        \cdot
        \abs{z}^k \\
        &=
        \sum_{\ell=0}^{\infty}
        \frac{(2 \supnorm{\xi})^\ell}{\ell!}
        \Majorant_\phi^{(\ell)}
        \bigl(
            \abs{z}
        \bigr) \\
        &=
        \Majorant_\phi
        \bigl(
            \abs{z} + 2 \supnorm{\xi}
        \bigr),
    \end{align}
    where $\Majorant_\phi^{(\ell)}$ denotes the $\ell$-th derivative of the Lie-Taylor
    majorant $\Majorant_\phi \in \Holomorphic(\C)$. In the final two steps, we have used
    that the Lie-Taylor majorant $\Majorant_\phi$ is entire. Thus
    $\Taylor_{\phi_0}(\argument; \gamma(t_1))$ indeed constitutes an entire function on
    $\C^n$. With other words, using that $\exp \colon V \longrightarrow U$ is bijective, we
    may define
    \begin{equation}
        \phi_1
        \colon
        \gamma(t_1)
        U
        \longrightarrow
        \C, \quad
        \phi_1
        \bigl(
            \gamma(t_1)
            \exp \chi
        \bigr)
        \coloneqq
        \Taylor_{\phi_0}
        \bigl(
            \underline{\chi}; \gamma(t_1)
        \bigr).
    \end{equation}
    Note that $U \cap \gamma(t_1)U \neq \emptyset$, as $\gamma(t_1)$ is within
    by construction. Hence, the intersection even contains an open connected
    neighbourhood $W$ of $\gamma(t_1)$, on which $\phi_0$ and $\phi_1$
    coincide. By the identity principle, we thus have $\phi_0 = \phi_1$ on the
    intersection of their domains. We may now iterate the procedure and define
    $\phi_2 \coloneqq \Taylor_{\phi_1}(\argument; \gamma(t_2))$ and so on.
    Crucially, we have shown that the Lie-Taylor majorant
    $\Majorant_{\phi_1}(\argument; \gamma(t_1))$ at $\gamma(t_1)$ is entire, as
    well. This yields the desired analytic continuation of $\phi_0$ along
    $\gamma$ after finitely many steps.
\end{proof}

\begin{remark}
    \label{rem:ExtensionFromGerm}
    Note that it suffices to know the values of $\phi \in \Entire_0(G)$ in some
    neighbourhood of the group unit to set off this process. This observation will be
    crucial to drop a somewhat unnatural assumption we are going to make in our initial
    version of Theorem~\ref{thm:Extension} to simplify the proof. Moreover, a combination
    of our construction, the uniqueness of the extension, the identity principle and
    Theorem~\ref{thm:LieTaylor} implies that an analogue of
    \eqref{eq:CompactExtension} holds for $\phi$ and $\Phi$.
\end{remark}

Investing the monodromy theorem, as it can be found in the textbook
\cite[Sec.~10.3.4]{krantz:1999a}, yields the existence of a globally defined holomorphic
extensions under the additional assumption that $G_\C$ is simply connected.
\begin{corollary}
    \label{cor:ExtensionSimplyConnected}
    Let $G$ be a connected Lie group such that $\eta \colon G \longrightarrow G_\C$ is
    injective and $G_\C$ is simply connected. Then, for every $\phi \in \Entire_0(G)$,
    there is a unique holomorphic function $\Phi \in \Holomorphic(G_\C)$ such that $\Phi
    \circ \eta = \phi$.
\end{corollary}

Taking another look at \eqref{eq:ComplexificationVsCovering}, the idea is now that
this describes the situation on the universal covering group, at least when
$\tilde{\eta}$ is injective.
\begin{theorem}[Extension]
    \index{Entire functions!Extension}
    \index{Extension}
    \label{thm:Extension}
    Let $G$ be a connected Lie group such that the universal complexification morphism
    $\tilde{\eta} \colon \widetilde{G} \longrightarrow \widetilde{G}_\C$ of the universal
    covering group $\widetilde{G}$ is injective. Then, for every $\phi \in
    \Entire_R(G)$, where $R \ge 0$,
    there exists a unique holomorphic function
    \begin{equation}
        \Phi
        \in
        \Holomorphic(G_\C)
        \cap
        \Entire_R(G_\C)
        \qquad
        \textrm{such that}
        \qquad
        \eta^* \Phi
        =
        \phi.
    \end{equation}
    Moreover,
    \begin{equation}
        \label{eq:Extension}
        \Phi
        \bigl(
            \eta(g) \exp_{G_\C}(\chi) \exp_{G_\C}(\I \xi)
        \bigr)
        =
        \Taylor_\phi
        \bigl(
            \underline{\chi + \I \xi}; g
        \bigr)
        =
        \Taylor_\Phi
        \bigl(
            \underline{\chi + \I \xi}; \eta(g)
        \bigr)
    \end{equation}
    holds for all $g \in G_\C$, $\chi,\xi \in \LieAlg(\eta(G)) \subseteq \LieAlg(G_\C)
    = \LieAlg(\eta(G))_\C$.
\end{theorem}
\begin{proof}
    The idea is to combine Proposition~\ref{prop:ComplexificationVsCovering} with
    \eqref{eq:ComplexificationVsCovering}. Denoting the universal covering projection by
    $p \colon \tilde{G} \longrightarrow G$, we define $\psi \coloneqq \phi \circ p$, which is
    $\pi_1(G)$-invariant and an element of $\Entire_0(G)$ by
    Proposition~\ref{prop:EntireFunctor}. Applying
    Corollary~\ref{cor:ExtensionSimplyConnected}, we obtain a unique holomorphic
    function $\Psi \in \Holomorphic(\tilde{G}_\C)$ with $\Psi \circ \tilde{\eta} = \psi$. Next,
    we show that $\Psi$ is invariant under the action of $\tilde{\eta}(\pi_1(G))$. To
    this end, fix $h \in
    \pi_1(G)$ and consider the auxiliary function $\Psi_h(g) \coloneqq \Psi(gh)$ for
    some $g \in G_\C$. Notice that
    \begin{equation*}
        \Psi
        \bigl(
            \tilde{\eta}(g)
            \tilde{\eta}(h)
        \bigr)
        =
        \psi
        (gh)
        =
        \psi(g)
        =
        \Psi
        \bigl(
            \tilde{\eta}(g)
        \bigr)
        \qquad
        \textrm{for all }
        g \in \tilde{G}.
    \end{equation*}
    By Proposition~\ref{prop:CauchyRiemann}, which we may apply by virtue of
    $\LieAlg(G_\C) = \liealg{g}_\C$, this means that the Lie-Taylor series of the
    difference $\Psi-\Psi_h$
    around the group unit vanishes and thus we get $\Psi = \Psi_h$ by a
    combination of the Lie-Taylor formula \eqref{eq:LieTaylorFormal} with the
    usual identity principle.
    Consequently, $\Psi$ is indeed invariant under the action of
    $\tilde{\eta}(\pi_1(G))$ and thus factors
    through the quotient projection $\pi \colon \tilde{G} \longrightarrow G_\C$, yielding the
    desired extension~$\Phi$ with $\Psi = \Phi \circ \pi$: Indeed, recall that $\eta
    \circ p = \pi
    \circ \tilde{\eta}$ by virtue of \eqref{eq:ComplexificationFromUniversalCovering}
    and thus
    \begin{equation*}
        \Phi
        \circ
        \eta
        \circ
        p
        =
        \Phi
        \circ
        \pi
        \circ
        \tilde{\eta}
        =
        \Psi
        \circ
        \tilde{\eta}
        =
        \psi
        =
        \phi
        \circ
        p.
    \end{equation*}
    Invoking the local bijectivity of $p$, this implies $\Phi \circ \eta = \phi$ as
    desired. The uniqueness of $\Phi$ is clear by
    Corollary~\ref{cor:IdentityPrinciple} and we have already noted the validity of
    \eqref{eq:Extension} in Remark~\ref{rem:ExtensionFromGerm}. This also proves the
    remaining statement.
\end{proof}

Being more careful, we can do even better.
\begin{proposition}
    \index{Extension!General}
    \index{Oka!Principle}
    \label{prop:ExtensionGeneral}
    The assertions of Theorem~\ref{thm:Extension} also hold if $\tilde{\eta}$ is not
    injective.
\end{proposition}
\begin{proof}
    Let $\phi \in \Entire_0(G)$ and $\psi \coloneqq \phi \circ p$ as before. By
    Corollary~\ref{cor:UniversalComplexificationProperties},
    \ref{item:UniversalComplexificationSimplyConnected2}, the universal complexification
    morphism $\tilde{\eta} \colon \tilde{G} \longrightarrow \tilde{G}_\C$ is locally injective.
    Thus, there exists a neighbourhood $\tilde{U} \subseteq \tilde{G}$ of the group unit
    such that $\tilde{\eta} \at{\tilde{U}}$ is a diffeomorphism with inverse $\Xi$ defined in
    the neighbourhood $\tilde{\eta}(\tilde{U})$ of the group unit in the image group
    $\tilde{\eta}(\tilde{G})$. Here we use the fact that $\tilde{\eta}(\tilde{G}) \subseteq
    \tilde{G}_\C$ is a closed submanifold by
    Corollary~\ref{cor:ImageOfEtaSubmanifold}. Consider the function
    \begin{equation}
        \tilde{\psi}
        \colon
        \tilde{\eta}(\tilde{U})
        \longrightarrow
        \C, \quad
        \tilde{\psi}
        \coloneqq
        \psi
        \circ
        \Xi.
    \end{equation}
    By continuity of the group multiplication, we find another neighbourhood $V \subseteq
    \tilde{\eta}(\tilde{U})$ of the group unit in $\tilde{\eta}(\tilde{G})$ such that $V \cdot V
    \subseteq \tilde{\eta}(\tilde{U})$. As $\Xi$ is a group morphism as the inverse of a
    group morphism, we thus have $\Xi(gh) = \Xi(g) \Xi(h)$ for all $g,h \in V$.
    Consequently, for any $\xi \in \liealg{g}$, we have
    \begin{equation}
        \Lie(\xi)
        \tilde{\psi}
        \at[\Big]{\E}
        =
        \Lie_{X_\xi}
        \bigl(\Xi^* \psi\bigr)
        \at[\Big]{\E}
        =
        \Xi^*
        \bigl(
            \Lie(T_\E \Xi \xi)
            \psi
        \bigr)
        \at[\Big]{\E}
        =
        \Lie\bigl(T_\E \Xi \xi\bigr)
        \psi
        \at[\Big]{\E}.
    \end{equation}
    As in \cite[Prop.~4.3, \textit{ii.)}]{heins.roth.waldmann:2023a}, it follows that the
    Lie-Taylor majorant $T_{\tilde{\psi}}(\argument;\E)$ is entire. By
    Remark~\ref{rem:ExtensionFromGerm}, we may thus use the techniques from
    Proposition~\ref{prop:ExtensionAlongPaths} to obtain another zero neighbourhood
    $V_0 \subseteq \tilde{G}_\C$ and a holomorphic
    \begin{equation}
        \psi_0 \colon V_0 \longrightarrow
        \C
        \quad \textrm{with} \quad
        \psi_0
        \at[\Big]{V_0 \cap V}
        =
        \tilde{\psi},
    \end{equation}
    which may be analytically continued along any path in $\tilde{G}_\C$ starting at the
    group unit. By the simple connectedness of $\tilde{G}_\C$, which follows from
    Corollary~\ref{cor:UniversalComplexificationProperties},
    \ref{item:UniversalComplexificationSimplyConnected2}, and the monodromy theorem,
    there thus exists a unique globally defined holomorphic extension
    \begin{equation}
        \tilde{\psi}
        \colon
        G_\C
        \longrightarrow
        \C
    \end{equation}
    of $\psi_0$. From here, one may proceed as in the proof of
    Theorem~\ref{thm:Extension} to check that $\psi$ factors through the quotient
    projection $\pi \colon \tilde{G}_\C \longrightarrow G$ and that this provides the
    unique
    \begin{equation}
        \Phi \in \Holomorphic(G_\C)
        \qquad
        \textrm{with} \quad
        \Phi \circ \eta = \phi.
    \end{equation}
    Here, we use that $\tilde{\psi}$ is \emph{locally}
    $\tilde{\eta}(\pi_1(G))$-invariant: To see this, we fix $h
    \in \pi_1(G)$ and if necessary shrink $V$ further to achieve $\tilde{\eta}(h) \cdot V \in
    \tilde{\eta}(\tilde{U})$. By virtue of
    \begin{equation}
        \tilde{\psi}
        \at[\big]{g \tilde{\eta}(h)}
        =
        \psi
        \circ
        \Xi
        \at[\big]{g \tilde{\eta}(h)}
        =
        \psi
        \bigl(
            \Xi(g) h
        \bigr)
        =
        \psi(g)
    \end{equation}
    for $g \in V$ this then implies $\pi_1(G)$-invariance of $\psi_0$ and thus of $\psi$ by
    the identity principle.
\end{proof}

\begin{example}[Automatic periodicity]
    \index{Special linear group!Automatic periodicity}
    \index{Extension!Automatic periodicity}
    \label{ex:AutomaticPeriodicity}
    We revisit the special linear group $\liegroup{SL}_2(\R)$ and its
    universal covering group $p \colon \widetilde{\liegroup{SL}}_2(\R) \longrightarrow
    \liegroup{SL}_2(\R)$
    from Example~\ref{ex:SpecialLinear}. Recall that
    \begin{equation}
        \liegroup{SL}_2(\R)_\C
        \cong
        \liegroup{SL}_2(\C)
        \cong
        \widetilde{\liegroup{SL}}_2(\R)_\C
    \end{equation}
    and let $f \in \Entire_0(\widetilde{\liegroup{SL}}_2(\R))$. By
    Theorem~\ref{thm:Extension}, there
    exists a unique holomorphic extension
    \begin{equation}
        F
        \in
        \Holomorphic
        \bigl(
            \liegroup{SL}_2(\C)
        \bigr)
        \cap
        \Entire_0
        \bigl(
            \liegroup{SL}_2(\C)
        \bigr)
    \end{equation}
    such that $F \circ p = f$. This means that $F$ factors through $p$ and is thus
    necessarily periodic with respect to the fundamental group
    $\pi_1(\liegroup{SL}_2(\R))$. Consequently, the algebra
    $\Entire_0(\widetilde{\liegroup{SL}}_2(\R))$ does not separate points in this
    case. In Theorem~\ref{thm:Restriction} we will see that
    \begin{equation}
        \Holomorphic
        \bigl(
            \liegroup{SL}_2(\C)
        \bigr)
        \subseteq
        \Entire_0
        \bigl(
            \liegroup{SL}_2(\C)
        \bigr),
    \end{equation}
    where we view $\liegroup{SL}_2(\C)$ as a real Lie group to define entire functions
    as before. This matches nicely with the well-known fact that every representation
    of
    $\widetilde{\liegroup{SL}}_2(\R)$ factors over the fundamental group
    $\pi_1(\liegroup{SL}_2(\R)$, which in particular implies the same for all
    representative functions, which are entire by
    Theorem~\ref{thm:RepresentativeFunctions} and always possess holomorphic
    extensions by Example~\ref{ex:RepresentativeExtension}. Indeed,
    the analogous statement holds for any covering group
    of~$\liegroup{SL}_2(\R)$. For
    a comprehensive discussion and a proof of these observations we refer to
    \cite[Ex.~9.5.18]{hilgert.neeb:2012a}.
\end{example}

This completes the resolution of the extension problem: Every entire function
defined on a Lie group possesses -- for better and for worse -- a unique
holomorphic extension to the universal complexification.
We turn towards the problem of restriction. That is, we investigate, whether
\begin{equation}
    \Entire_0(G)
    \cap
    \Holomorphic(G)
    =
    \Holomorphic(G)
\end{equation}
for a complex Lie group $G$ and how the corresponding topologies are related. In
particular, we may consider this problem for general complex groups, not just groups
arising as universal complexifications.

As it has been the theme throughout, the key towards a resolution of the problem
consists in a suitable version of the
Cauchy estimates. To formulate them, we fix some notation. For
continuous functions $\phi \colon G \longrightarrow \C$ and compact subsets
$K \subseteq G$, we write
\begin{equation}
    \gls{SupnormCompact}
    \coloneqq
    \max_{g \in K}
    \abs[\big]
    {\phi(g)}.
\end{equation}
To facilitate the estimate, we make the additional assumption of working with a matrix
Lie group. That is, we assume that $G$ is linear and choose a faithful finite
dimensional representation as matrices. In this context, we denote the operator
norm~\gls{Norm} of a matrix~$M \in \Mat_n(\C)$ simply by $\norm{M}$. Moreover, it makes
sense to speak \emph{both} of operator
norms of group elements and of Lie algebra elements. Thus there is a direct way to
compare sizes of infinitesimal and exponentiated elements. Moreover, as the Lie
exponential is simply given by the exponential series for matrix Lie groups, we have
the estimate
\begin{equation}
    \label{eq:OperatorNormExponential}
    \norm
    [\big]
    {\exp \xi}
    \le
    e^{\norm{\xi}}
    \qquad
    \textrm{for all }
    \xi \in \liealg{g},
\end{equation}
where $e$ denotes Euler's number.\footnote{Leonhard Euler (1707-1783) was an
extraordinarily
productive Swiss polymath. So much so that many of his discoveries were not named
after him, but rather after the second scientist to observe them.}
This may be thought of as a continuity estimate for the Lie exponential.
Keeping this in mind,
we arrive at the following Lie theoretic incarnation of the Cauchy estimates.
\begin{lemma}[Lie theoretic Cauchy estimates]
    \index{Cauchy estimates!Lie group}
    \label{lem:CauchyEstimates}
    Let $G$ be a complex connected matrix Lie group, $\phi \in \Holomorphic(G)$, $K_0
    \subseteq \liealg{g}$ a compact neighbourhood of zero, $r > 0$ and $k \in \N_0$.
    Then, writing $K_r \coloneqq \{g \in G \colon \norm{g} \le \gls{Euler}^r \}$, we
    have
    \begin{equation}
        \label{eq:CauchyEstimates}
        \abs[\Big]
        {
            \bigl(
                \Lie(\xi_1 \tensor \cdots \tensor \xi_k)
                \phi
            \bigr)
            (g \exp \xi)
        }
        \le
        \norm{\phi}_{g\exp(K_0)K_r}
        \cdot
        \frac{k^k}{r^k}
    \end{equation}
    for all $g \in G$, $\xi \in K_0$ and $\xi_1, \ldots, \xi_k \in \liealg{g}$ with $\norm{\xi_j} \le 1$ for $j=1, \ldots, k$.
\end{lemma}
\begin{proof}
    As the Lie exponential $\exp(z\xi_j)$ is the flow of the left invariant vector field $X_{\xi_j}$, we have
    \begin{equation}
        \Lie(\xi_j)
        \phi
        (g)
        =
        \frac{\partial}{\partial z}
        \phi
        \bigl(
            g \exp(\xi) \exp(z \xi_j)
        \bigr)
        \at[\Big]{z = 0}
    \end{equation}
    for $j = 1, \ldots, k$, $g \in G$ and $\xi \in K_0$. Alternatively, this readily
    follows from \eqref{eq:LieTaylorOneVariable}. As the Lie exponential induces
    a holomorphic atlas, the composition
    $\phi \circ \ell_{g \exp(\xi)} \circ \exp \colon \liealg{g} \longrightarrow \C$ is
    holomorphic. By induction, we get
    \begin{equation*}
        \Lie(\xi_1 \tensor \cdots \tensor \xi_k)
        \phi
        (g \exp(\xi))
        =
        \frac{\partial}{\partial z^1}
        \cdots
        \frac{\partial}{\partial z^k}
        \phi
        \bigl(
            g \exp(\xi) \exp(z^k \xi_k) \cdots \exp(z^1 \xi_1)
        \bigr)
        \at[\Big]{z^1 = \ldots = z^k = 0},
    \end{equation*}
    where
    \begin{equation}
        \C^k
        \ni
        z
        \mapsto
        \phi
        \bigl(
            g
            \exp(\xi)
            \exp(z^k \xi_k)
            \cdots
            \exp(z^1 \xi_1)
        \bigr)
    \end{equation}
    is holomorphic. Notably, there are only first derivatives involved. Applying the
    usual
    Cauchy estimates for first derivatives with all radii equal to $r/k$ yields
    \begin{equation*}
        \label{eq:CauchyEstimatesProof}
        \abs[\Big]
        {
            \bigl(
            \Lie(\xi_1 \tensor \cdots \tensor \xi_k)
            \phi
            \bigr)
            (g \exp \xi)
        }
        \le
        \frac{k^k}{r^k}
        \max_{\abs{z^1} = \cdots = \abs{z^k} = r/k}
        \abs[\Big]
        {\phi(g \exp(\xi) \exp(z^k \xi_k) \cdots \exp(z^1 \xi_1))}.
        \tag{$\star$}
    \end{equation*}
    By submultiplicativity of the operator norm and \eqref{eq:OperatorNormExponential},
    we note
    \begin{align*}
        \norm[\big]
        {\exp(z^k \xi_k) \cdots \exp(z^1 \xi_1)}
        &\le
        \norm[\big]
        {\exp(z^k \xi_k)}
        \cdots
        \norm[\big]
        {\exp(z^1 \xi_1)} \\
        &\le
        \exp
        \bigl(
            \abs{z^k} \cdot \norm{\xi_k}
        \bigr)
        \cdots
        \exp
        \bigl(
            \abs{z^1} \cdot \norm{\xi_1}
        \bigr) \\
        &=
        \exp
        \bigl(
            \abs{z^k} \cdot \norm{\xi_k}
            +
            \cdots
            +
            \abs{z^k} \cdot \norm{\xi_1}
        \bigr) \\
        &\le
        \exp(r)
    \end{align*}
    for $\abs{z^1} = \cdots = \abs{z^k} = r/k$. That is, we have
    \begin{equation}
        \exp(z^k \xi_k)
        \cdots
        \exp(z^1 \xi_1)
        \in
        K_r
        \qquad
        \textrm{whenever}
        \quad
        \abs{z^1} = \ldots = \abs{z^k} = r/k.
    \end{equation}
    Consequently, \eqref{eq:CauchyEstimates} follows from
    \eqref{eq:CauchyEstimatesProof} by varying $z$ and $\xi$.
\end{proof}

This simple inequality allows us to resolve the restriction problem with a positive
answer for linear groups. Here, we tacitly -- and in the sequel without further
mention -- use the functoriality of $\Entire_0$ and $\Holomorphic$, to pass from a
linear group to a concrete realization as matrices to apply the Lie theoretic Cauchy
estimates.
\begin{theorem}[Restriction]
    \index{Entire functions!Restriction}
    \index{Restriction}
    \label{thm:Restriction}
    Let $G$ be a complex connected linear Lie group. Then
    \begin{equation}
        \Holomorphic(G)
        \subseteq
        \Entire_0(G)
    \end{equation}
    induces an embedding of Fréchet algebras, where we endow $\Holomorphic(G)$ with the topology of locally uniform convergence.
\end{theorem}
\begin{proof}
    Let $\phi \in \Holomorphic(G)$ and fix a basis $(\basis{e}_1, \ldots, \basis{e}_n)$ of $\liealg{g}$ such that $\norm{\basis{e}_j} = 1$ for $j=1,\ldots,n$. We once again use the Lie-Taylor majorants from \eqref{eq:LieTaylorMajorant}. By Lemma~\ref{lem:CauchyEstimates}, we have
    \begin{align}
        \Majorant_\phi
        \bigl(
            \abs{z}
        \bigr)
        &=
        \sum_{k=0}^{\infty}
        \frac{\abs{z}^k}{k!}
        \sum_{\alpha \in \N_n^k}
        \abs[\Big]
        {
            \bigl(
                \Lie(\alpha) \phi
            \bigr)(\E)
        } \\
        &\le
        \norm{\phi}_{K_r}
        \cdot
        \sum_{k=0}^{\infty}
        \frac{\abs{z}^k}{k!}
        \sum_{\alpha \in \N_n^k}
        \frac{k^k}{r^k} \\
        &=
        \norm{\phi}_{K_r}
        \cdot
        \sum_{k=0}^{\infty}
        \biggl(
            \frac{n}{r}
        \biggr)^k
        \frac{k^k}{k!}
        \cdot
        \abs{z}^k
    \end{align}
    for \emph{any} $r > 0$ and $z \in \C$. The power series on the right-hand side has a positive radius of convergence $R$ for any $r > 0$, namely
    \begin{equation*}
        R
        =
        \frac{r}{n}
        \lim_{k \rightarrow \infty}
        \frac{\sqrt[k]{k!}}{k}
        =
        \frac{r}{ne},
    \end{equation*}
    where $e$ denotes again Euler's number. Taking the limit $r \rightarrow
    \infty$, we
    see that $\Majorant_\phi$
    is indeed entire. This moreover proves that the seminorms $\seminorm{q}_{0,c}$ are
    continuous with respect to the topology of locally uniform convergence on
    $\Holomorphic(G)$; here one first fixes $c > 0$ and then chooses $r$ sufficiently
    large. Conversely, let $K_0 \subseteq \liealg{g}$ be a compact zero neighbourhood
    such that the restriction $\exp \at{K_0}$ is bijective and $\norm{\xi} \le 1$ for all $\xi \in
    K$. The Lie-Taylor formula \eqref{eq:LieTaylor} yields then
    \begin{equation*}
        \norm[\big]
        {\phi}_{g\exp(K_0)}
        \le
        \max_{\xi \in K_0}
        \abs[\big]
        {\phi(g \exp(\xi))}
        \le
        \sum_{k=0}^{\infty}
        \frac{1}{k!}
        \sum_{\alpha \in \N_n^k}
        \abs[\Big]
        {
            \bigl(
            \Lie(\alpha) \phi
            \bigr)(g)
        }
        \eqqcolon
        \seminorm{q}_{0,1,g}(\phi)
    \end{equation*}
    for any $g \in G$. By the translation invariance from Theorem~\ref{thm:Symmetries},
    \ref{item:TranslationInvariance}, the seminorm~$\seminorm{q}_{0,1,g}$ is
    continuous on $\Entire_0(G)$, completing the proof.
\end{proof}

Having established both an extension and a restriction result, the natural question is
what happens during a roundtrip. Astonishingly, the ultimately rather pleasing answer
is \emph{nothing whatsoever}, including the topologies.
\begin{corollary}
    \label{cor:EntireVsHolomorphic}
    \index{Entire!vs. Holomorphic}
    Let $G$ be a connected Lie group such that its universal
    complexification~$(G_\C,\eta)$ is a linear Lie group.\footnote{Note that, by the
    universal property from Definition~\ref{def:UniversalComplexification}, every
    linear Lie
    group complexifies to a linear complex Lie group.} Then the pullback
    \begin{equation}
        \label{eq:EntireVsHolomorphicIso}
        \eta^*
        \colon
        \Holomorphic(G_\C)
        \longrightarrow
        \Entire_0(G), \quad
        \eta^*
        \Phi
        \coloneqq
        \Phi \circ \eta
    \end{equation}
    with the universal complexification morphism $\eta$ is well defined and an isomorphism of Fréchet algebras, where the products are defined pointwisely.
\end{corollary}
\begin{proof}
    We check that extension and restriction, which is modelled by
    \eqref{eq:EntireVsHolomorphicIso}, are inverse processes. Given $\Phi \in
    \Holomorphic(G_\C)$, Theorem~\ref{thm:Restriction} yields $\Phi \in
    \Entire_0(G)$ and thus $\eta^* \Phi \in \Entire_0(G)$ by
    Proposition~\ref{prop:EntireFunctor}. Invoking Theorem~\ref{thm:Extension},
    the function $\eta^* \Phi$ induces a unique holomorphic function $\Psi \in
    \Holomorphic(G_\C)$ such that
    \begin{equation}
        \eta^*
        \Psi
        =
        \Psi
        \circ
        \eta
        =
        \eta^* \Phi.
    \end{equation}
    But $\Phi$ also has this property, so $\Psi = \Phi$ by uniqueness. Conversely, let $\phi \in \Entire_0(G)$ be given. By Theorem~\ref{thm:Extension}, there exists a unique holomorphic function $\Phi \in \Holomorphic(G_\C) \cap \Entire_0(G)$ such that $\eta^* \Phi = \phi$. Consequently, \eqref{eq:EntireVsHolomorphicIso} is well defined and bijective. It is clear that \eqref{eq:EntireVsHolomorphicIso} preserves sums and pointwise products.
\end{proof}

On the level of the observable algebra $\Entire_0(G)$, we have thus put the
program of holomorphic deformation quantization as
outlined in Section~\ref{sec:QuantizationHolomorphic} into practice
successfully. The observable algebra
allows for a holomorphic description and the conditions on the Lie-Taylor
coefficients
simply translate to locally uniform convergence within the domains of the
holomorphic extensions.

Using our results, one may now transfer the continuous star products on
$\Entire_0(G) \tensor \Sym_1^\bullet(\liealg{g}_\C)$ to star products on the
holomorphic polynomial algebra
\begin{equation}
    \Holomorphic_\Pol(T^*G)
    \coloneqq
    \Holomorphic(G_\C)
    \tensor_\pi
    \Sym_1^\bullet(\liealg{g}_\C)
\end{equation}
from \eqref{eq:HolomorphicPolynomialAlgebra} by declaring
\eqref{eq:EntireVsHolomorphicIso} to be an
isomorphism of Fréchet algebras, now with respect to the star products. In view
of the following remark this is however not satisfactory in general:
\begin{remark}[Real forms]
    \label{rem:RealForms}
    \index{Real form}
    Not every complex Lie group $G$ is the universal complexification of a real Lie group
    \gls{RealForm}. Real subgroups $H \subseteq G$ with $H_\C \cong G$ are called
    \emph{real
    forms}. The obstructions already arise on the Lie algebraic level. For an example of
    minimal dimension -- namely, three -- we refer to
    \cite[Ex.~5.1.24]{hilgert.neeb:2012a}. In view of
    Corollary~\ref{cor:ImageOfEtaSubmanifold} and
    Lemma~\ref{lem:ComplexificationEtaLocallyInjective}, one may not obtain the
    corresponding connected and simply connected Lie group by means of universal
    complexification.

    Indeed, assume $H_\C \cong G$ for some Lie group $H$. Then $H$
    is connected and its image~$\eta(H)$ complexifies to $H_\C$, as well, by
    Corollary~\ref{cor:ImageOfEtaSubmanifold}. But $\LieAlg(\eta(H))_\C = \LieAlg(G)$ by
    Lemma~\ref{lem:ComplexificationEtaLocallyInjective}, which is a contradiction.
\end{remark}

However, by Remark~\ref{rem:ComplexLieAlgebras} we may follow the general philosophy
outlined
in Section~\ref{sec:QuantizationHolomorphic} and use our explicit formulas
\eqref{eq:StarProductFactorization} for the standard ordered star product to treat
complex linear Lie groups $G$ and their \emph{complex} Lie algebras
\gls{LieAlgStandaloneComplex}, see again Remark~\ref{rem:LieBrackets}
directly.\footnote{The author would like to sincerely
thank Matthias Schötz for discussing until the construction was no longer both
unnecessarily complicated and eminently inelegant.} By virtue
of the Lie theoretic Cauchy estimates from Lemma~\ref{lem:CauchyEstimates}, this
turns out to be considerably simpler than in \cite[Sec.~6]{heins.roth.waldmann:2023a},
see also the general strategy sketched after
Theorem~\ref{thm:LieGroupStarProductContinuity}. Note that by
Corollary~\ref{cor:EntireVsHolomorphic}, working with all holomorphic functions in
\eqref{eq:HolomorphicPolynomialAlgebra} corresponds to $R = 0$.\footnote{This loses
some generality for the mixed products,
but if $G$ arises as a universal complexification, we ultimately recover the biggest
of the algebras from \cite{heins.roth.waldmann:2023a} for the full product
regardless.} We use Lemma~\ref{lem:CauchyEstimates} and
Lemma~\ref{lem:ProjectiveTensorOfEll1} to establish the continuity of mixed products
first.
\begin{proposition}
    \label{prop:StarProduct}
    \index{Holomorphic!Standard ordered star product}
    Let $G$ be a complex connected linear Lie group.
    \begin{propositionlist}
        \item The restricted standard ordered star product
        \begin{equation}
            \star_\std
            \colon
            \bigl(
                \mathbb{1}
                \tensor
                \Sym_1^\bullet(\hat{\liealg{g}})\bigr)
            \tensor
            \bigl(
                \Holomorphic(G)
                \tensor
                1
            \bigr)
            \longrightarrow
            \Holomorphic_\Pol(T^*G)
        \end{equation}
        is well defined and continuous.
        \item The mapping
        \begin{equation}
            \C
            \times
            \Holomorphic(G)
            \tensor
            \Sym_1^\bullet(\hat{\liealg{g}})
            \ni
            (\hbar, \phi, P)
            \mapsto
            (\mathbb{1}
            \tensor
            P)
            \star_\std
            (\phi
            \tensor
            1)
            \in
            \Holomorphic_\Pol(T^*G)
        \end{equation}
        is Fréchet holomorphic.
    \end{propositionlist}
\end{proposition}
\begin{proof}
    Fix a complex basis $(\basis{e}_1, \ldots, \basis{e}_n)$
    of $\hat{\liealg{g}}$ and let $\seminorm{p}$ be the corresponding~$\ell^1$-norm.
    Let moreover $K_0
    \subseteq \hat{\liealg{g}}$ be a compact neighbourhood of zero, $g \in G$, $\phi
    \in
    \Holomorphic(G)$, $c \ge 0$ and indices
    $1 \le j_1, \ldots, j_k \le n$. Motivated by the explicit formula
    \eqref{eq:StarProductFactorization}, we use
    Lemma~\ref{lem:CauchyEstimates} and \eqref{eq:TruncatedPolynomialEstimate}
    with $R = 1$ to estimate
    \begin{align}
        &\sum_{p=0}^{k}
        \frac{\abs{\hbar}^p}{p! (k-p)!}
        \sum_{\sigma \in S_k}
        \norm[\Big]
        {
            \Lie(f_{j_{\sigma(p)}} \tensor \cdots \tensor f_{j_{\sigma(1)}})
            \phi
        }_{g \exp(K_0)}
        \seminorm{p}_{R,c}
        \bigl(
            \basis{e}_{j_{\sigma(p+1)}}
            \vee \cdots \vee
            \basis{e}_{j_{\sigma(k)}}
        \bigr) \\
        \le
        &\sum_{p=0}^{k}
        \frac{\abs{\hbar}^p}{p! (k-p)!}
        \sum_{\sigma \in S_k}
        \norm{\phi}_{g \exp(K_0) K_r}
        \frac{p^p}{r^p}
        \biggl(
            \frac{(k-p)!}{k!}
        \biggr)^R
        c^{-p}
        \cdot
        \seminorm{p}_{R,c}
        \bigl(
            \basis{e}_{j_1}
            \vee \cdots \vee
            \basis{e}_{j_k}
        \bigr) \\
        \le
        &\norm{\phi}_{g \exp(K_0) K_r}
        \cdot
        \seminorm{p}_{R,c}
        \bigl(
            \basis{e}_{j_1}
            \vee \cdots \vee
            \basis{e}_{j_k}
        \bigr)
        \sum_{p=0}^{k}
        \biggl(
            \frac{(k-p)!}{k!}
        \biggr)^{R-1}
        \frac{\abs{\hbar}^p}{p!}
        \frac{p^p}{r^p \cdot c^p} \\
        \le
        &\norm{\phi}_{g \exp(K_0) K_r}
        \cdot
        \seminorm{p}_{R,c}
        \bigl(
        \basis{e}_{j_1}
        \vee \cdots \vee
        \basis{e}_{j_k}
        \bigr)
        \sum_{p=0}^{\infty}
        \frac{\abs{\hbar}^p}{p!}
        \frac{p^p}{r^p \cdot c^p}
    \end{align}
    for any $r > 0$. Taking $r \coloneqq  e \abs{\hbar}/c$ with Euler's number $e$ yields a
    convergent series
    \begin{equation}
        M_r
        \coloneqq
        \sum_{p=0}^{\infty}
        \frac{\abs{\hbar}^p}{p!}
        \frac{p^p}{r^p \cdot c^p}
        <
        \infty,
    \end{equation}
    which may be chosen locally uniformly with respect to $\hbar \in \C$. This proves
    both claims by a combination of the triangle inequality with
    Corollary~\ref{cor:FrechetPowerSeries}.
\end{proof}

Mutatis mutandis, this should also work for holomorphic functions of fixed order
tensorized with $\Sym^\bullet_R(\hat{\liealg{g}})$ with parameters chosen as in
\cite[Lem.~6.2]{heins.roth.waldmann:2023a}. However, one would not gain anything for
the full product, as $R = 1$ is optimal in
Theorem~\ref{thm:LieAlgebraStarProductContinuity}. This is also why we
have defined \eqref{eq:HolomorphicPolynomialAlgebra} with $R = 1$ in the first
place.
\begin{theorem}[Complex Standard Ordered Star Product]
    \index{Standard ordered star product!Complexified}
    \index{Continuity!Complexified Gutt star product}
    \index{Star product!Complex standard ordered}
    \label{thm:StarProductHolomorphicContinuity}
    Let $G$ be a complex connected linear Lie group.
    \begin{theoremlist}
        \item The standard ordered star product
        \begin{equation}
            \star_\std
            \colon
            \Holomorphic_\Pol(T^*G)
            \tensor
            \Holomorphic_\Pol(T^*G)
            \longrightarrow
            \Holomorphic_\Pol(T^*G)
        \end{equation}
        is well defined and continuous.
        \item The mapping
        \begin{equation}
            \C
            \times
            \Holomorphic_\Pol(T^*G)
            \otimes
            \Holomorphic_\Pol(T^*G)
            \ni
            (\hbar, \Phi, \Psi)
            \mapsto
            \Phi
            \star_\std
            \Psi
            \in
            \Holomorphic_\Pol(T^*G)
        \end{equation}
        is Fréchet holomorphic.
    \end{theoremlist}
\end{theorem}
\begin{proof}
    The proof is almost identical to \cite[Proof of
    Thm.~6.3]{heins.roth.waldmann:2023a}, where our Proposition~\ref{prop:StarProduct}
    replaces \cite[Lem.~6.2]{heins.roth.waldmann:2023a}. The key idea is to use the
    associativity of $\star_\std$ to separate the full product into the pointwise
    product, mixed products and the Lie algebra star product, which we have noted in
    \eqref{eq:StarProductFactorization}. Each may then be estimated with what we have
    already shown, where we also use the elementary estimate
    \begin{equation}
        \norm{\phi \cdot \psi}_K
        \le
        \norm{\phi}_K
        \cdot
        \norm{\psi}_K
    \end{equation}
    for compact subsets $K \subseteq G$ and $\phi,\psi \in \Holomorphic(G)$.
\end{proof}

\index{Questions!Gelfand-Shilov for $G_\C$}
\index{Gelfand, Isreal}
\index{Shilov, Georgiy}
From here, one may in principle use Corollary~\ref{cor:EntireVsHolomorphic} to
establish a star product on real forms, reversing our earlier logic. This of course
yields nothing new, but the circle closes nicely. It would be interesting to describe
the holomorphic analogues of the spaces $\Entire_R(G)$ with $R > 0$. By what we have
shown, they are tautologically given by
\begin{equation}
    \label{eq:GelfandShilovLie}
    \Entire_R(G)
    =
    \Holomorphic(G_\C)
    \cap
    \Entire_R(G_\C)
\end{equation}
for Lie groups $G$ such that the universal complexification $\eta \colon G
\longrightarrow G_\C$ is a linear Lie group. As in Example~\ref{ex:StdOrdIII}, the
condition on the Lie-Taylor majorant
may be rephrased as its holomorphic extensions being \emph{of finite order $1/R$ and
minimal type}. The precise correspondence can be found in
\cite[Prop.~4.15]{heins.roth.waldmann:2023a}. In the case of $\R^n$, this construction
yields the classical
Gelfand-Shilov spaces and it would be interesting to characterize the
spaces~\eqref{eq:GelfandShilovLie} without referring to their Lie-Taylor
majorants, say by means of an intrinsic growth condition on the universal
complexification $G_\C$. See also
Remark~\ref{rem:GelfandShilov} for our prior discussion on the algebra of
$R$-entire functions.

\section{Continuous and Smooth Vectors}
\label{sec:ContinuousAndSmoothVectors}
\epigraph{Slake was one of those places, Moist thought, that you put on the map
because it was embarrassing to to have a map with holes in it.}{\emph{Raising Steam}
-- Terry Pratchett}
% !TeX root = ../Dissertation.tex

This section serves as an elementary and self contained discussion of the notions of
continuous and smooth vectors for Lie algebra and Lie group representations on locally
convex spaces.
Everything we cover is standard, but the presentation is catered towards our ultimate
goal: Understanding analytic, entire and strongly entire vectors in
Section~\ref{sec:EntireVectors}. As these notions turn out to be stronger than
continuity and smoothness, we thus get to distinguish between smooth and holomorphic
phenomena more easily. Moreover, we put particular emphasis in developing
criteria to
decide the regularity of a particular vector, whose power we showcase in
Example~\ref{ex:Schroedinger} at the end of the section.

Throughout, we assume that $V$ is a \emph{complex} vector space by complexifying if
necessary, as this does not complicate our considerations whatsoever and thus
essentially comes along for free. Indeed, once again, some arguments even simplify.

\index{Strong continuity}
\index{Orbit map}
In the sequel, we shall use the terminology and symbols from
Section~\ref{sec:UniformConvergenceOnBoundedSets}.
Moreover, we call a group representation $\pi \colon G \longrightarrow \GL(V)$ of a
Lie group $G$ on a topological vector space $V$ \emph{strongly continuous} if
it is continuous with respect to the $\sigma$-topology -- that is, the weak
topology\footnote{This is a rather unfortunate clash of terminology, which
historically arose from the presence of even weaker topologies in the setting of
representations on Hilbert spaces $V$. Instead of fighting these particular windmills,
we join them and rejoice in rotational bliss.} -- on
\begin{equation}
    \GL(V)
    \coloneqq
    \bigl\{
        L
        \in
        \Linear(V)
        \colon
        L
        \textrm{ is invertible}
    \bigr\}
    \subseteq
    \Linear_\sigma(V).
\end{equation}
This means that the \emph{orbit} maps
\begin{equation}
    \label{eq:OrbitMapGroup}
    G
    \ni
    g
    \; \mapsto \;
    \gls{OrbitMap}(g)
    \coloneqq
    \pi(g)v
    \in
    V
\end{equation}
are continuous for all $v \in V$. Notably, having a topological group would
suffice for the study of continuity. However, many of our arguments require $G$
to be a topological manifold, and thus working with Lie groups from the start is
no loss of generality by the positive resolution of Hilbert's fifth problem. For a
comprehensive discussion, we refer to the monograph \cite[§1]{tao:2014a}.

The analogous notion for Lie algebra representations $\varrho \colon \liealg{g}
\longrightarrow \Linear(V)$ would demand the continuity of the corresponding
orbit maps
\begin{equation}
    \label{eq:OrbitMapAlgebra}
    \liealg{g}
    \ni
    \xi
    \; \mapsto \;
    \varrho(\xi)v
    \in
    V
\end{equation}
for all $v \in V$. However, by finite dimensionality of $\liealg{g}$ and linear
dependence of \eqref{eq:OrbitMapAlgebra} on $\xi$, this is always fulfilled.
Hence, we refrain from using this terminology and speak of Lie algebra
representations instead. The non-trivial assumption is then that $\varrho(\xi)$ is
continuous for every $\xi \in \liealg{g}$. Similarly, we implicitly assume the
continuity of the linear mappings~$\pi(g)$ for all $g \in G$ by
assuming~$\pi$ maps into $\GL(V)$. Whenever both types of orbit maps are
involved, we sometimes refer to \eqref{eq:OrbitMapAlgebra} as
\emph{infinitesimal} orbit maps to distinguish both concepts. If $V$ is barrelled,
the Banach-Steinhaus Theorem yields the following by
\cite[Beginning of Sec.~4.1]{warner:1972a}.
\begin{lemma}
    \label{lem:StrongContinuityContinuityOfAction}
    Let $\pi \colon G \longrightarrow V$ be a representation of a Lie group
    $G$ on a barrelled topological space $V$. Then $\pi$ is strongly continuous
    if and only if the action mapping
    \begin{equation}
        \label{eq:ActionMapping}
        \Phi
        \colon
        G \times V
        \longrightarrow
        V, \quad
        \Phi(g,v)
        \coloneqq
        \pi(g)v
    \end{equation}
    is continuous with respect to the product topology.
\end{lemma}
\begin{proof}
    By virtue of $\pi_v = \Phi(\argument, v)$ it is clear that continuity of
    \eqref{eq:ActionMapping} implies the strong continuity of $\pi$. To see the
    converse implication, let $(g_n)_{n \in \N} \subseteq G$ and
    $(v_\alpha)_{\alpha \in J}$ be convergent
    with limits $g \in G$ and $v \in V$. Fix a compact neighbourhood $K \subseteq G$
    of $g$. Then there exists an index $N \in \N$ such that~$g_n \in K$ for all $n \ge
    N$. Consider the set
    \begin{equation}
        \mathcal{T}
        \coloneqq
        \bigl\{
            \pi(k)
            \colon
            k \in K
        \bigr\}
        \subseteq
        \Linear(V).
    \end{equation}
    By continuity of \eqref{eq:OrbitMapGroup}, the set $\mathcal{T}$ is weakly bounded, i.e. bounded in $\Linear_\sigma(V)$. By the Banach-Steinhaus Theorem, $\mathcal{T}$ is thus equicontinuous. Let now $O \subseteq V$ be some zero neighbourhood and use the continuity of the vector space addition to choose another zero neighbourhood $U \subseteq V$ such that $U + U \subseteq V$. By equicontinuity of $\mathcal{T}$, there exists an index $\alpha_0 \in J$ such that
    \begin{equation}
        \Phi(k,v)
        -
        \Phi(k,v_\alpha)
        =
        \pi(k)
        \bigl(
            v - v_\alpha
        \bigr)
        \in
        U
        \qquad
        \textrm{for all }
        k \in K
        \textrm{ and }
        \alpha \later \alpha_0.
    \end{equation}
    By continuity of $\pi_v$, we moreover find an index $M \ge N$ such that
    \begin{equation}
        \Phi(g_n,v)
        -
        \Phi(g,v)
        \in U
        \qquad
        \textrm{for all }
        n \ge M.
    \end{equation}
    Putting everything together, we get
    \begin{align}
        \Phi(g,v)
        -
        \Phi(g_n, v_\alpha)
        =
        \Phi(g,v)
        -
        \Phi(g_n,v)
        +
        \Phi(g_n,v)
        -
        \Phi(g_n,v_\alpha)
        \in
        U + U
        \subseteq
        O
    \end{align}
    for all $n \ge M$ and $\alpha \later \alpha_0$. This completes the proof.
\end{proof}

We proceed by establishing some more or less obvious, but useful technical simplifications for checking strong continuity. Having a group representation first of all implies that it suffices to check continuity of \eqref{eq:OrbitMapGroup} at the group unit.
\begin{lemma}
    \label{lem:OrbitMapContinuityAtE}
    \index{Strong continuity!At group unit}
    Let $G$ be a Lie group and $\pi \colon G \longrightarrow \GL(V)$ be a
    representation on some locally convex space $V$ and $v \in V$. Then the orbit map
    $\pi_v \colon G \longrightarrow V$ is continuous if and only if it is continuous
    at the group unit.
\end{lemma}
\begin{proof}
    Let $v \in V$ be such that $\pi_v$ is continuous at the group unit $\E$. If $g_n \rightarrow g$ in $G$, then $g^{-1} g_n \rightarrow \E$ and thus
    \begin{equation*}
        \pi_v(g_n)
        =
        \pi(g_n)v
        =
        \pi(g)
        \pi(g^{-1} g_n)
        v
        =
        \pi(g)
        \pi_v(g^{-1} g_n)
        \overset{n \rightarrow \infty}{\longrightarrow}
        \pi(g)
        \pi_v(\E)
        =
        \pi(g)v
        =
        \pi_v(g).
    \end{equation*}
    Consequently, the orbit map $\pi_v$ is continuous at $g$.
\end{proof}

Moreover, strong continuity may be checked on generators.
\begin{lemma}
    \label{lem:SemidirectProductContinuousVectors}
    \index{Strong continuity!Generators}
    Let $G$ and $N$ be Lie groups and $\varphi \colon G \times N \longrightarrow
    N$ be a smooth action of~$G$ by automorphisms of $N$. Assume moreover
    that
    \begin{equation*}
        \pi
        \colon
        N \rtimes_\varphi G
        \longrightarrow
        \GL(V)
    \end{equation*}
    is a group representation on some locally convex space $V$ such that the
    restrictions $\pi \at{N}$ and $\pi \at{G}$ are strongly continuous. Then $\pi$ is
    strongly continuous.
\end{lemma}
\begin{proof}
    Let $v \in V$. By Lemma~\ref{lem:OrbitMapContinuityAtE} it suffices to check
    continuity of $\pi_v$ at the group unit. To this end, let $(h_n, g_n)_n \in N
    \rtimes_\varphi G$ be a sequence converging to the group unit~$(\E_N, \E_G)$. As
    $N \rtimes_\varphi G$ carries the product topology, this implies that $h_n
    \rightarrow \E_N$ and $g_n \rightarrow \E_G$. Now, writing $\pi^N \coloneqq \pi
    \at{N}$ and $ \pi^G \coloneqq \pi \at{G}$, we note
    \begin{equation}
        \pi_v(h_n, g_n)
        =
        \pi(h_n, g_n)v
%        =
%        \pi(h_n, \E_G)
%        \pi(\E_N, g_n)v
        =
        \pi^N(h_n)
        \pi^G(g_n)v
        =
        \pi^N_{\pi^G(g_n)v}(h_n)
        =
        \pi^N_{\pi_v^G(g_n)}(h_n),
    \end{equation}
    which implies $\pi_v(h_n, g_n) \rightarrow v = \pi_v(\E_N,\E_G)$.
\end{proof}

Along the way, we have shown the following factorization of the orbit maps, which is occasionally useful.
\begin{corollary}
    \index{Orbit map!Factorization}
    Let $G$ and $N$ be Lie groups and $\varphi \colon G \times N
    \longrightarrow N$ be a smooth action of~$G$ by automorphisms of $N$. If
    \begin{equation*}
        \pi
        \colon
        N \rtimes_\varphi G
        \longrightarrow
        \GL(V)
    \end{equation*}
    is a group representation on some locally convex space $V$, then the orbit
    maps fulfil
    \begin{equation}
        \label{eq:OrbitMapSemidirectFactorization}
        \pi_v(h,g)
        =
        \pi^N_{\pi_v^G(g)}(h)
        \qquad
        \textrm{for all }
        v \in V, \,
        g \in G, \,
        h \in N,
    \end{equation}
    where $\pi^N \coloneqq \pi \at{N}$ and $ \pi^G \coloneqq \pi \at{G}$.
\end{corollary}

Beyond Lemma~\ref{lem:StrongContinuityContinuityOfAction}, the Banach-Steinhaus
Theorem moreover yields closedness of the set of all continuous vectors within
barrelled locally convex spaces $V$.
\begin{proposition}
    \label{prop:ContinuousVectorsClosed}
    Let $G$ be a Lie group and $\pi \colon G \longrightarrow \GL(V)$ be a
    representation on some barrelled locally convex space $V$. Assume
    moreover that for every compact set $K \subseteq G$, the image~$\pi(K)
    \subseteq \GL(V)$ is weakly bounded.\footnote{The sets $\pi(K)v \subseteq
    V$ are bounded for all $v \in V$.} Then
    \begin{equation}
        \label{eq:ContinuousVectors}
        \Continuous(\pi)
        \coloneqq
        \bigl\{
            v \in V
            \colon
            \pi_v
            \in
            \Continuous(G,V)
        \bigr\}
        \subseteq
        V
    \end{equation}
    is closed.
\end{proposition}
\begin{proof}
    Let $(v_\alpha)_{\alpha}  \subseteq \Continuous(\pi)$ be a convergent net
    with limit $v \in V$. Consider the corresponding net $(\pi_{v_\alpha})_{\alpha}
    \subseteq \Continuous(G,V)$. We prove that this net converges to the orbit
    map~$\pi_v$ in the space $\Continuous(G,V)$ -- that is, uniformly on
    compact subsets of $G$. To
    this end, let $K \subseteq G$ be compact
    and~$\seminorm{p} \in \cs(V)$. By assumption,~$\pi(K) \subseteq \GL(V)$ is then
    weakly bounded and thus equicontinuous by the Banach-Steinhaus Theorem.
    Thus, there is a $\seminorm{q} \in \cs(V)$ with
    \begin{equation*}
        \seminorm{q}
        \bigl(
            \pi(g)w
        \bigr)
        \le
        \seminorm{q}(w)
        \qquad
        \textrm{for all }
        g \in K
        \textrm{ and }
        w \in V.
    \end{equation*}
    Consequently,
    \begin{equation*}
        \max_{g \in K}
        \seminorm{p}
        \bigl(
            \pi_{v_\alpha}(g)
            -
            \pi_v(g)
        \bigr)
        =
        \max_{g \in K}
        \seminorm{p}
        \bigl(
            \pi(g)
            (v - v_\alpha)
        \bigr)
        \le
        \seminorm{q}
        \bigl(
            v_\alpha
            -
            v
        \bigr)
        \qquad
        \textrm{for all $\alpha$},
    \end{equation*}
    which implies the convergence of the net $(\pi_{v_\alpha})_{\alpha}$ to $\pi_v$
    within $\Continuous(G,V)$. In particular, the limit $\pi_v$ is continuous and thus
    $v \in \Continuous(\pi)$ by \eqref{eq:ContinuousVectors}.
\end{proof}

An important special instance is given by unitary representations $\pi$ on some
Hilbert space $V$. Indeed, in this case each of the operators $\pi(g)$ is an isometry,
which implies the continuity estimate $\norm{\pi(g)v} = \norm{v}$ for all $v \in V$.
The closedness of $\Continuous(\pi)$ now means that it suffices to check strong
continuity
on a dense subspace.
\begin{corollary}
    \label{cor:StrongContinuityOnDense}
    \index{Strong continuity!Density}
    Let $G$ be a Lie group and $\pi \colon G \longrightarrow \GL(V)$ be a
    representation on some barrelled locally convex space $V$. Assume
    moreover that for every compact subset $K \subseteq G$, the image~$\pi(K)
    \subseteq \GL(V)$ is weakly bounded. If $\Continuous(\pi)$ is dense in $V$,
    then $\pi$ is strongly continuous.
\end{corollary}
\begin{proof}
    Closed and dense subsets coincide with the whole space.
\end{proof}

%\index{$\Fun$-vector}
\index{Vector!Smooth}
\index{Smooth vectors}
Going one step further, we use the smooth structure of $G$ to ask whether the
orbit mappings \eqref{eq:OrbitMapGroup} and \eqref{eq:OrbitMapAlgebra} are
smooth, i.e. have
Fréchet derivatives of all orders with values in the completion of $V$. In
Section~\ref{sec:HolomorphicFrechet} we have seen that sequential completeness
suffices
to take Fréchet derivatives. In the sequel we shall therefore always assume that $V$
is sequentially complete and Hausdorff to avoid certain
technical complications.\footnote{For instance, without the Hausdorff property of $V$
the mere definition of the first derivative as a function could already entail a
\emph{choice}. And even after making this choice, one immediately runs into questions
of the following unpleasant flavour: If $\phi \colon V \longrightarrow V$ is
continuously differentiable, can the existence of the second derivative depend on
which limit points one uses to define the first?} We denote the subspace of all such
$v \in V$ by \gls{SmoothVectorGroup} and \gls{SmoothVectorAlgebra}, respectively, and
call them \emph{smooth vectors}. As the domains of $\pi_v$ and
\eqref{eq:OrbitMapAlgebra} are finite dimensional, it suffices to consider left
invariant or partial derivatives to check their smoothness. For a comprehensive
discussion of smooth vectors, we refer to \cite[Sec.~4.4]{warner:1972a}.
\begin{remark}[Gårding space]
    \index{Gårding!Space}
    \index{Gårding!Lars}
    \index{Dixmier, Jacques}
    \index{Malliavin, Paul}
    \index{Pettis, Billy James}
    \index{Gelfand, Isreal}
    \label{rem:Garding}
    Let $G$ be a Lie group and $\pi \colon G \longrightarrow \GL(V)$ a representation on
    a sequentially complete Hausdorff locally convex space $V$. Then space of
    smooth
    vectors $\Cinfty(\pi) \subseteq V$ is typically not sequentially complete itself.
    In fact, essentially the opposite is true, as
    Gårding's Theorem \cite{garding:1947a} asserts that $\Cinfty(\pi)$ is always
    \emph{sequentially dense} in $V$ -- and thus in particular dense. He proved this
    by considering the
    Gårding\footnote{Lars Gårding (1919-2014) was a Swedish mathematician specialized
    in partial differential equations. His doctoral advisor was Marcel Riesz and he
    supervised Lars Hörmander.} space
    \begin{equation}
        \label{eq:Garding}
        \gls{Garding}
        \coloneqq
        \Span
        \biggl\{
        \int_G
        \phi(g)
        \pi(g)v
        \D g
        \;\Big|\;
        \phi \in \Cinfty_c(G),
        v \in V
        \biggr\},
    \end{equation}
    where $\Cinfty_c(G)$ denotes the space of smooth functions with compact support
    and the integral is most easily handled as a weak integral in the sense of
    Gelfand-Pettis\footnote{Billy James Pettis (1913-1979) was an American functional
    analyst specialized on infinite dimensional integration and measure theory. He
    generalized Gelfand's work on weak integrals from the Lebesgue measure on an
    interval to general measure spaces.}
    \cite{gelfand:1936a, pettis:1938a}, see also
    \cite[Ch.~3.5]{rudin:1991a}. Alternatively, assuming $V$ to be Fréchet, one may
    treat it as a vector-valued Riemann integral by \cite[Ch.~3,
    Exercise~23]{rudin:1991a}. The integrals in \eqref{eq:Garding} are closely related
    to the convolution algebra of $G$\footnote{In fact, taking $\pi$ as the
    translation representation on the Banach space $\Lone(G)$, the integral reduces to
    the usual convolution product. In general, sending $\phi$ to
    $\int_G\phi(g)\pi(g)v\D g$ yields a homomorphism from the convolution algebra
    $\Lone(G)$ to $\Linear(V)$.} and using standard convolution techniques it is not
    hard to prove that $\mathfrak{G}$ consists entirely of smooth vectors and is
    sequentially dense in~$V$. Astonishingly, Dixmier\footnote{Jacques Dixmier (born
    1924) is a French operator algebraist, who introduced numerous fundamental objects
    within the theory such as the Dixmier trace and the Dixmier mapping. He also
    appears to be impervious to the passage of time.} and
    Malliavin\footnote{Paul Malliavin (1925-2010) was a French harmonic and stochastic
    analyst. He devised what the Malliavin calculus, which is central
    to Monte Carlo simulations and mathematical finance.}
    \cite{dixmier.malliavin:1978a} proved that one even has $\mathfrak{G} =
    \Cinfty(\pi)$ for strongly continuous representations on Fréchet spaces.
\end{remark}

\begin{example}[Translation Representation I: Continuous \& smooth vectors]
    \index{Translation representation!Continuous vectors}
    \index{Translation representation!Smooth vectors}
    \label{ex:TranslationRepContinuousSmooth}
    We shall mostly discuss right translations, but of course, analogous statements hold
    for left translations.\footnote{The opposite group is, after all, itself a Lie group.}
    Recall our notation $r_g$ for the
    left translation by $g$ from \eqref{eq:LeftAndRightMultiplication}. They induce a
    representation
    \begin{equation}
        \label{eq:RightTranslationPullback}
        r^*
        \colon
        G
        \longrightarrow
        \GL
        \bigl(
            \gls{Mappings}(G)
        \bigr), \quad
        r^*_g\phi(h)
        \coloneqq
        \phi(hg)
    \end{equation}
    on the space $\Map(G)$ of all complex valued functions defined on $G$,
    endowed with the
    locally convex and complete topology of pointwise convergence. To investigate the
    strong continuity of $r^*$ at some $\phi \in \Map(G)$, we have to consider the
    continuity of the mappings
    \begin{equation}
        G
        \ni
        g
        \; \mapsto \;
        r^* \phi
        \in
        \Map(G),
    \end{equation}
    which by first countability of $G$ is equivalent to sequential continuity. If now $(g_n) \subseteq G$ is a sequence with limit $g \in G$, then we have to check that
    \begin{equation}
        \phi(h g_n)
        =
        r_{g_n}^* \phi(h)
        \longrightarrow
        r_g^* \phi(h)
        =
        \phi(hg)
        \qquad
        \textrm{for all }
        h \in G.
    \end{equation}
    That is, the function $\phi$ needs to be continuous itself.
    Thus $r^*$ restricts to a strongly continuous representation on the space of
    complex valued continuous functions~$\Continuous(G)$, endowed with the topology of
    pointwise convergence.

    However, this space is not even sequentially complete and its natural locally
    convex topology is instead given by uniform convergence on compact subsets of $G$.
    This raises the question whether this more natural topology yields the same
    continuous vectors. Proceeding as before, we have to check when $\phi(h g_n)
    \rightarrow \phi(hg)$, but now uniformly on compact subsets of $G$. This is
    equivalent to the \emph{uniform continuity} of $\phi$ on every compact subset,
    which is automatic if $\phi$ is continuous at all. Thus $r^*$ is also a strongly
    continuous representation on the space $\Continuous(G)$, endowed with the topology
    of uniform convergence on compact sets.

    Let now $\phi \in \Cinfty(r^*)$. By smoothness, the difference quotient fulfils
    \begin{equation}
        \label{eq:RightTranslationDifferenceQuotient}
        \frac{r^*_{\exp(t\xi)} \phi - \phi}{t}
        \at[\bigg]{g}
        =
        \frac{\phi(g \exp(t\xi)) - \phi(g)}{t}
        \overset{t \rightarrow 0}{\longrightarrow}
        \Lie(\xi)
        \phi(g)
    \end{equation}
    for all $g \in G$ and $\xi \in \liealg{g}$. Again, we may either consider the limits pointwisely or locally uniformly on compact subsets by means of uniform continuity, this time also of the derivatives. Consequently, every smooth vector of $r^*$ is continuously differentiable as a function. Applying this reasoning inductively yields that smooth vectors are necessarily smooth as functions and conversely it is clear that every smooth function constitutes a smooth vector. Thus,
    \begin{equation}
        \label{eq:TranslationSmoothVectors}
        \Cinfty(r^*)
        =
        \Cinfty(G).
    \end{equation}

    Note that \eqref{eq:RightTranslationDifferenceQuotient} recovers the Lie algebra action \eqref{eq:LieDerivativeOnLieAlg} on $\Cinfty(G)$ and its higher orders yield \eqref{eq:LieDerivativeOnEnveloping}. In particular, we recover
    \begin{equation}
        \label{eq:TranslationVsDifferentiation}
        \Lie(\xi)
        \circ
        r^*_g
        =
        r^*_g
        \circ
        \Lie
        \bigl(
            \Ad_{g^{-1}} \xi
        \bigr)
        \qquad
        \textrm{for all }
        g \in G
        \textrm{ and }
        \xi \in \Universal(\liealg{g}_\C)
    \end{equation}
    with the adjoint action $\Ad_g \colon \liealg{g} \longrightarrow \liealg{g}$ from
    \eqref{eq:AdjointActionGroup}. Indeed, on the one hand, we have
    \begin{equation}
        \Lie(\xi)
        r_g^*
        \phi
        \at[\Big]{h}
        =
        \frac{\D}{\D t}
        \at[\Big]{t=0}
        \bigl(
            r^*_{\exp(t\xi)}
            (r_g^* \phi)
        \bigr)(h)
        =
        \frac{\D}{\D t}
        \at[\Big]{t=0}
        r_{\exp(t\xi)g}^*
        \phi(h)
        =
        \frac{\D}{\D t}
        \at[\Big]{t=0}
        \phi
        \bigl(
            h \exp(t\xi) g
        \bigr)
    \end{equation}
    for all $g,h \in G$ and $\xi \in \liealg{g}$. Now,
    \begin{equation}
        \label{eq:ConjVsAdjoint}
        \exp(t\xi) g
        =
        g g^{-1} \exp(t\xi) g
        =
        g \Conj_g^{-1}(\exp(t\xi))
        =
        g
        \exp(\Ad_g^{-1} t\xi)
    \end{equation}
    and thus
    \begin{equation}
        \Lie(\xi)
        r_g^*
        \phi
        \at[\Big]{h}
        =
        \frac{\D}{\D t}
        \at[\Big]{t=0}
        \phi
        \bigl(
            hg
            \exp(\Ad_g^{-1} t\xi)
        \bigr)
        =
        \Lie(\Ad_g^{-1})
        \phi
        \at[\Big]{hg}
        =
        r_g^*
        \Lie(\Ad_g^{-1})
        \phi
        \at[\Big]{h}.
    \end{equation}
    This proves \eqref{eq:TranslationVsDifferentiation} by
    \eqref{eq:LieDerivativeOnEnveloping} and an induction.

    The space of smooth vectors thus provides a natural domain for the
    \emph{infinitesimal} counterpart of the translation action. However,
    \eqref{eq:LieDerivativeOnLieAlg} is \emph{not} continuous with respect to the
    subspace topology induced by the inclusion $\Cinfty(G) \subseteq \Continuous(G)$.
    For $G = \R$, this is the well-known issue of discontinuity of infinitessimal
    generators of semigroups of operators, say acting on a Hilbert space. That being
    said, it \emph{is} continuous if we endow $\Cinfty(G)$ with its natural Fréchet
    topology.

    \index{Convolution product}

    Taking another look a Remark~\ref{rem:Garding} and in particular
    \eqref{eq:Garding} for $v = \psi \in \Continuous(G)$, we may evaluate $\pi(g)v =
    \ell_{g^{-1}}^* \phi$ at some $h \in G$, which yields the left convolution product
    \begin{equation*}
        \int_G
        \phi(g)
        \pi(g)v
        \D g
        \at[\Big]{h}
        =
        \int_G
        \phi(g)
        \psi(g^{-1} h)
        \D g
        =
        \phi
        \ast
        \psi
        \at[\Big]{h}.
    \end{equation*}
    The aforementioned compatibility of the Gårding space $\mathfrak{G}$ with the
    convolution algebra~$(\Cinfty_c(G), \ast)$ reduces to the associativity of
    $\ast$ in
    this case. Moreover, Gårding's Theorem~\cite{garding:1947a} reflects the fact that
    $\phi \ast \psi$ is smooth for all $\phi \in \Cinfty_c(G)$ and $\psi \in
    \Continuous(G)$. Coupled with the Dixmier-Malliavin Theorem
    \cite{dixmier.malliavin:1978a}, this finally has the not at all obvious
    consequence that every smooth function on $G$ may be written as a \emph{finite}
    linear combination of convolution products, even though the convolution algebra
    lacks a unit if $G$ is not discrete.
\end{example}

As we will see in Theorem~\ref{thm:ActionDifferentiation}, the preceding observations
on the differentiated action are a general feature of strongly continuous group
representations. To make this precise, including its continuity, we give
$\Cinfty(\pi)$ its own locally convex topology, where we follow the discussion in
\cite[Sec.~4.4.1]{warner:1972a}. The idea is to identify $v$ with the mapping $\pi_v$
from \eqref{eq:OrbitMapGroup}, i.e. we realize $\Cinfty(\pi)$ as a subspace of
$\Cinfty(G,V)$. We have used this idea in the continuous setting to prove
Proposition~\ref{prop:ContinuousVectorsClosed}. Recall that the seminorms
\begin{equation}
    \index{Smooth vectors!Seminorms}
    \index{Seminorms!Smooth vectors}
    \label{eq:SmoothVectorsSeminorms}
    \gls{SeminormsSmooth}(\phi)
    \coloneqq
    \max_{g \in K}
    \seminorm{p}
    \bigl(
        \Lie(\xi)
        \phi
        (g)
    \bigr)
\end{equation}
for $\xi \in \Universal^\bullet(\liealg{g}_\C)$, compact sets $K \subseteq G$ and
$\seminorm{p} \in \cs(V)$ yield a defining system of seminorms for~$\Cinfty(G,V)$ that
endow it with the structure of a Hausdorff locally convex space. It is
Fréchet if and only if $V$ is, as $V \hookrightarrow \Cinfty(G,V)$ as the constant
functions.
\begin{lemma}
    \label{lem:SmoothVectorsTopology}
    \index{Smooth vectors!Topology}
    Let $\pi \colon G \longrightarrow \GL(V)$ be a strongly continuous
    representation of a Lie group $G$ on a complete Hausdorff locally convex
    space $V$. Then
    \begin{equation}
        \label{eq:OrbitMapSpace}
        \gls{SmoothVectorGroupAsFunctions}
        \coloneqq
        \bigl\{
            \pi_v
            \colon
            G \longrightarrow V
            \;\big|\;
            v \in \Cinfty(\pi)
        \bigr\}
        \subseteq
        \Cinfty(G,V)
    \end{equation}
    is a closed subspace, which is linearly isomorphic to $\Cinfty(\pi)$ by means of
    \begin{equation}
        \label{eq:OrbitMapSpaceIsomorphism}
        \gls{VectorToOrbitMap}
        \colon
        \Cinfty(\pi) \longrightarrow \Cinfty_\pi(G,V), \quad
        \Pi(v)
        \coloneqq
        \pi_v.
    \end{equation}
    If $V$ is merely sequentially complete, then the same is true for $\Cinfty(\pi)$.
\end{lemma}
\begin{proof}
    Let $v \in \Cinfty(\pi)$ be such that $\pi_v = 0$. By \eqref{eq:OrbitMapGroup} for
    $g = \E$, this implies
    \begin{equation}
        0
        =
        \pi(\E)v
        =
        v.
    \end{equation}
    This establishes the linear isomorphy of $\Cinfty(\pi)$ and
    \eqref{eq:OrbitMapSpace} by means of \eqref{eq:OrbitMapSpaceIsomorphism}. Let now
    $(v_\alpha)_{\alpha \in J} \subseteq \Cinfty(\pi)$ be such that $\pi_{v_\alpha}
    \rightarrow \phi \in \Cinfty(G,V)$. Notice that
    \begin{equation}
        \pi_{v_\alpha}(g)
        =
        \pi(g)v_\alpha
        =
        \pi(g)
        \pi(\E)
        v_\alpha
        =
        \pi(g)
        \pi_{v_\alpha}(\E)
        \qquad
        \textrm{for all }
        \alpha \in J
        \textrm{ and }
        g \in G.
    \end{equation}
    Taking limits on both sides yields
    \begin{equation}
        \phi(g)
        =
        \pi(g)
        \phi(\E)
        =
        \pi_{\phi(\E)}(g)
        \qquad
        \textrm{for all }
        g \in G.
    \end{equation}
    By smoothness of $\phi$, this in particular implies $\phi(\E) \in \Cinfty(\pi)$.
    Thus we have shown that
    \begin{equation}
        \phi
        =
        \pi_{\phi(\E)}
        \qquad
        \textrm{with}
        \quad
        \phi(\E) \in \Cinfty(\pi),
    \end{equation}
    which completes the proof.
\end{proof}

In the sequel, we identify vectors $v$ with their orbit maps $\pi_v$ and endow
$\Cinfty(\pi)$ with the subspace topology inherited from the inclusion
$\Cinfty_\pi(G,V) \subseteq \Cinfty(G,V)$. That is to say, we demand that
\eqref{eq:OrbitMapSpaceIsomorphism} is an isomorphism of locally convex spaces. Using
$\C \subseteq \Universal^\bullet(\liealg{g}_\C)$ and $K \coloneqq \{\E\}$ in
\eqref{eq:SmoothVectorsSeminorms} yields that this topology is finer than the one
inherited from the inclusion $\Cinfty(\pi) \subseteq V$.
Lemma~\ref{lem:SmoothVectorsTopology} moreover states that we obtain a Fréchet space
in this manner whenever $V$ is Fréchet itself.
\begin{proposition}[Invariance of $\Cinfty(\pi)$]
    \index{Smooth vectors!Invariance under $\pi$}
    \index{Invariance!Smooth vectors}
    \label{prop:SmoothVectorsRestriction}
    Let $\pi \colon G \longrightarrow \GL(V)$ be a strongly continuous
    representation of a Lie group $G$ on a complete Hausdorff locally convex
    space $V$.
    \begin{lemmalist}
        \item \label{item:SmoothVectorsRestrictionPiInvariance}
        We have $\pi(g)v \in \Cinfty(\pi)$ for all $g \in G$ and $v \in \Cinfty(\pi)$.
        \item The representation $\pi$ restricts to a strongly continuous representation on $\Cinfty(\pi)$, endowed with the subspace topology inherited from $V$.
        \item The representation $\pi$ restricts to a strongly continuous
        representation on $\Cinfty(\pi)$, endowed with the subspace topology inherited
        from $\Cinfty(G,V)$.
        \item The smooth vectors of both restrictions coincide with $\Cinfty(\pi)$.
    \end{lemmalist}
\end{proposition}
\begin{proof}
    Let $v \in \Cinfty(\pi)$, $g \in G$ and write $w \coloneqq \pi(g)v$. As $\pi$ is a representation, we have
    \begin{equation}
        \label{eq:SmoothVectorInvarianceProof}
        \pi_{w}(h)
        =
        \pi(h)w
        =
        \pi(h)
        \pi(g)
        v
        =
        \pi(hg)v
        =
        \pi
        \bigl(
            r_g(h)
        \bigr)v
        =
        \bigl(
            \pi_v
            \circ
            r_g
        \bigr)(h)
    \end{equation}
    for all $h \in G$. Thus $\pi_w$ is smooth as composition of the smooth maps
    $\pi_v$ and $r_g$, which means $w \in \Cinfty(\pi)$. The strong continuity of $\pi
    \at{\Cinfty(\pi)}$ with respect to the inclusion~$\Cinfty(\pi) \subseteq V$ is
    then clear. Next, we study the topology inherited from $\Cinfty(V,G)$. Indeed,
    combining \eqref{eq:SmoothVectorInvarianceProof} and
    \eqref{eq:TranslationVsDifferentiation} we compute
    \begin{align}
        \seminorm{q}_{\xi_1 \tensor \cdots \tensor \xi_k, K, \seminorm{p}}
        \bigl(
            \pi(g)v
        \bigr)
        &=
        \max_{h \in K}
        \seminorm{p}
        \Bigl(
            \bigl(
                \Lie(\xi_k \tensor \cdots \tensor \xi_1)
                \pi_w
            \bigr)(h)
        \Bigr) \\
        &=
        \max_{h \in K}
        \seminorm{p}
        \Bigl(
            \bigl(
                \Lie(\xi_k \tensor \cdots \tensor \xi_1)
                r_g^*
                \pi_v
            \bigr)(h)
        \Bigr) \\
        &=
        \max_{h \in K}
        \seminorm{p}
        \Bigl(
        \bigl(
            r_g^*
            \Lie(\Ad_g^{-1} \xi_k \tensor \cdots \tensor \Ad_g^{-1} \xi_1)
            \pi_v
        \bigr)(h)
        \Bigr) \\
        &=
        \max_{h' \in Kg}
        \seminorm{p}
        \Bigl(
        \bigl(
            \Lie(\Ad_g^{-1} \xi_k \tensor \cdots \tensor \Ad_g^{-1} \xi_1) \pi_v
        \bigr)(h')
        \Bigr) \\
        &=
        \seminorm{q}_{\Ad_{g}^{-1} \xi_1 \tensor \cdots \tensor \Ad_{g}^{-1} \xi_k, Kg, \seminorm{p}}(v)
    \end{align}
    for all $\xi_1, \ldots, \xi_k \in \liealg{g}$, compact $K \subseteq G$ and
    $\seminorm{p} \in \cs(V)$. Thus $\pi(g)$ is continuous also in the topology
    inherited from $\Cinfty(G,V)$. As strong continuity and smoothness are both based
    on the properties of the same map $\pi_v$, it suffices to prove that every $v \in
    \Cinfty(\pi)$ induces a smooth orbit map $\pi_v \colon G \longrightarrow V$, i.e.
    $\pi_v \in \Cinfty(G, \Cinfty(\pi))$. Equivalently, by means of currying, it
    suffices to assert
    \begin{equation}
        \Pi(\pi_v)
        \in
        \Cinfty
        \bigl(
            G, \Cinfty_\pi(G,V)
        \bigr)
        \subseteq
        \Cinfty
        \bigl(
            G, \Cinfty(G,V)
        \bigr)
        \cong
        \Cinfty(G \times G, V).
    \end{equation}
    Now, keeping the ordering of arguments, $\Pi(\pi_v)(g,h) = \pi(h)\pi(g)v$ holds for all $g,h \in G$ and thus we indeed have $\Pi(\pi_v) \in \Cinfty(G \times G, V)$. Finally, $\Pi(\pi_v)$ comes from $v \in \Cinfty(\pi)$, which means that $v$ is a smooth vector also for the restriction of $\pi$ to $\Cinfty(\pi)$.
\end{proof}

With these preliminaries, we may now realize our hopes and dreams of differentiating
$\pi$ to a strongly continuous Lie algebra representation.
\begin{theorem}[Differentiation on $\Cinfty(\pi)$]
    \label{thm:ActionDifferentiation}
    \index{Smooth vectors!Differentiation}
    Let $G$ be a Lie group and
    \begin{equation}
        \pi
        \colon
        G \longrightarrow \GL(V)
    \end{equation}
    a strongly continuous representation on a sequentially complete Hausdorff
    locally convex space $V$. Then
    \begin{equation}
        \label{eq:LieAlgebraRep}
        T_\E \pi
        \colon
        \liealg{g}
        \longrightarrow
        \Linear
        \bigl(
        \Cinfty(\pi)
        \bigr), \quad
        T_\E \pi \xi \at[\Big]{v}
        \coloneqq
        \frac{\D}{\D t}
        \pi
        \bigl(
        \exp(t \xi)
        \bigr)
        v
        \at[\Big]{t=0}
        =
        \bigl(
            \Lie(\xi)
            \pi_v
        \bigr)
        (\E)
    \end{equation}
    is a Lie algebra representation with $\Cinfty(T_\E \pi) =
    \Cinfty(\pi)$. In particular, the tangent map~$T_\E \pi$ extends to an action of
    the universal enveloping algebra
    $\Universal^\bullet(\liealg{g}_\C)$ on $\Cinfty(\pi)$.
\end{theorem}
\begin{proof}
    The crucial point is that the limit in \eqref{eq:LieAlgebraRep} exists for all $v \in \Cinfty(\pi)$ and yields another smooth vector, as $t \mapsto \pi(\exp(t\xi))v$ is a coordinate expression of \eqref{eq:OrbitMapGroup}. To see the latter, note that
    \begin{equation}
        \label{eq:ActionsCommuteProof}
        \pi_{T_\E \pi \xi v}(g)
        =
        \pi(g)
        \bigl(
        (T_\E \pi \xi)(v)
        \bigr)
        =
        \frac{\D}{\D t}
        \pi
        \bigl(
        g \exp(t \xi)
        \bigr)
        v
        \at[\Big]{t=0}
        =
        \Lie(\xi)
        \pi_v(g)
    \end{equation}
    for all $g \in G$ by continuity of $\pi(g)$ and the fact that $\pi$ is a group representation. Having established this, it follows that $T_\E \pi$ is a Lie algebra morphism at once, as $\Lie$ is. By the universal property of the universal enveloping algebra $\Universal^\bullet(\liealg{g})$, the Lie algebra representation $T_\E \pi$ has a unique extension to $\Universal^\bullet(\liealg{g})$ as an algebra morphism, which we denote by the same symbol. Iterating \eqref{eq:LieAlgebraRep}, the extension is explicitly given by
    \begin{equation*}
        \label{eq:LieAlgebraRepProof}
        T_\E \pi
        \bigl(
        \xi_1 \tensor \cdots \tensor \xi_k
        \bigr)
        v
        =
        \Lie(\xi_1 \tensor \cdots \tensor \xi_k)
        \pi_v
        (\E)
        \tag{$\ast$}
    \end{equation*}
    for $\xi_1, \ldots, \xi_k \in \liealg{g}$, which is indeed well defined by
    smoothness of $v$. As $T_\E \pi$ is a homomorphism, this moreover implies
    $\Cinfty(T_\E \pi) = \Cinfty(\pi)$. Let now $\xi \in \liealg{g}$. Taking another
    look at~\eqref{eq:SmoothVectorsSeminorms} and \eqref{eq:LieAlgebraRep}, we get the
    continuity estimate
    \begin{equation}
        \seminorm{q}_{\chi, K, \seminorm{p}}
        \bigl(
            T_\E \pi \xi v
        \bigr)
        =
        \max_{g \in K}
        \seminorm{p}
        \Bigl(
            \Lie(\chi)
            \Lie(\xi)
            \pi_v(g)
        \Bigr)
        =
        \seminorm{q}_{\xi \tensor \chi, K, \seminorm{p}}
        (v)
    \end{equation}
    for all $\chi \in \Universal^\bullet(\liealg{g}_\C)$, compact sets $K \subseteq
    G$ and $\seminorm{p} \in \cs(V)$. Thus the linear map $T_\E \pi \xi$ is a
    continuous self-map of $\Cinfty(\pi)$.
\end{proof}

We proceed with several corollaries of our considerations.
\begin{corollary}
    \index{Smooth vectors!Intertwiners}
    Let $G$ be a Lie group and $\pi \colon G \longrightarrow \GL(V)$ a strongly continuous representation on a sequentially complete Hausdorff locally convex space $V$. Then
    \begin{equation}
        \label{eq:Intertwining}
        \pi(g)
        \circ
        T_\E \pi \xi
        \at[\Big]{v}
        =
        \Lie(\xi)
        \pi_v(g)
        \quad \textrm{and} \quad
        \Pi
        \circ
        T_\E \pi \xi
        =
        \Lie(\xi) \circ \Pi
    \end{equation}
    for all $g \in G$, $\xi \in \Universal^\bullet(\liealg{g}_\C)$ and $v \in \Cinfty(\pi)$.
\end{corollary}
\begin{proof}
    The first equality is \eqref{eq:ActionsCommuteProof}. For the second, note that
    \begin{equation}
        \Pi
        \circ
        T_\E \pi \xi
        \in
        \Linear
        \bigl(
            \Cinfty(\pi),\Cinfty(G,V)
        \bigr)
    \end{equation}
    and so we may evaluate both sides at first at $v \in \Cinfty(\pi)$ and
    subsequently at $g \in G$. Starting with the left-hand side in
    \eqref{eq:Intertwining} and using the first identity, this yields
    \begin{equation}
        \Pi
        \bigl(
            T_\E \pi \xi v
        \bigr)
        \at[\Big]{g}
        =
        \pi
        (g)
        T_\E \pi \xi v
        =
        \Lie(\xi)
        \pi_v(g)
        =
        \bigl(
            \Lie(\xi)
            \circ
            \Pi
        \bigr)
        (v)
        \at[\Big]{g}.
        \tag*{\qed}
    \end{equation}
\end{proof}

\index{Intertwiner}
The interpretation of \eqref{eq:Intertwining} is that the isomorphism $\Pi$ acts as an \emph{intertwiner} between the representations \eqref{eq:LieAlgebraRep} and \eqref{eq:LieDerivativeOnEnveloping} of $\Universal^\bullet(\liealg{g}_\C)$. Next, we generalize \eqref{eq:TranslationVsDifferentiation}.
\begin{corollary}
    \index{Smooth vectors!Intertwiners}
    Let $G$ be a Lie group and $\pi \colon G \longrightarrow \GL(V)$ a strongly continuous representation on a sequentially complete Hausdorff locally convex space $V$. Then
    \begin{equation}
        \label{eq:ActionsDoNotQuiteCommute}
        T_\E \pi \xi
        \circ
        \pi(g)
        =
        \pi(g)
        \circ
        T_\E \pi
        \bigl(
            \Ad_g^{-1} \xi
        \bigr)
    \end{equation}
    as operators on $\Cinfty(\pi)$ for all $g \in G$ and $\xi \in \Universal^\bullet(\liealg{g}_\C)$.
%    If $G$ is unimodular, then
%    \begin{equation}
%        \label{eq:ActionsCommute}
%        T_\E \pi \xi
%        \circ
%        \pi(g)
%        =
%        \pi(g)
%        \circ
%        T_\E \pi \xi
%    \end{equation}
\end{corollary}
\begin{proof}
    The reasoning is essentially identical to the proof of
    \eqref{eq:TranslationVsDifferentiation}. Let $g \in G$, $\xi \in \liealg{g}$
    and~$v \in \Cinfty(\pi)$. Using \eqref{eq:LieAlgebraRep}, \eqref{eq:ConjVsAdjoint}
    and \eqref{eq:Intertwining} we compute
    \begin{align}
        T_\E \pi \xi
        \circ
        \pi(g)
        \at[\Big]{v}
        &=
        \Lie(\xi)
        \pi_{\pi(g)v}
        (\E) \\
        &=
        \frac{\D}{\D t}
        \pi
        \bigl(
            \exp(t\xi)
        \bigr)
        \pi(g)v
        \at[\bigg]{t=0} \\
        &=
        \frac{\D}{\D t}
        \pi
        \bigl(
            \exp(t\xi)g
        \bigr)v
        \at[\bigg]{t=0} \\
        &=
        \frac{\D}{\D t}
        \pi
        \bigl(
            g \exp(t \Ad_g^{-1} \xi)
        \bigr)v
        \at[\bigg]{t=0} \\
        &=
        \bigl(
            \Lie(\Ad_g^{-1} \xi)
            \pi_v
        \bigr)(g) \\
        &=
        \pi(g)
        \circ
        T_\E
        \bigl(
            \Ad^{-1}_g \xi
        \bigr)
        \at[\Big]{v}.
    \end{align}
    This proves \eqref{eq:ActionsDoNotQuiteCommute} by inductive application.
\end{proof}

\begin{corollary}
    \label{cor:SmoothVectorsSeminorms}
    \index{Smooth vectors!Seminorms}
    \index{Seminorms!Smooth vectors}
    Let $G$ be a Lie group and $\pi \colon G \longrightarrow \GL(V)$ a representation on a sequentially complete Hausdorff locally convex space $V$. The seminorms
    \begin{equation}
        \seminorm{r}_{\xi,K,\seminorm{p}}
        \colon
        \Cinfty(\pi) \longrightarrow [0,\infty), \quad
        \seminorm{r}_{\xi,K,\seminorm{p}}(v)
        \coloneqq
        \max_{g \in K}
        \seminorm{p}
        \bigl(
            \pi(g)
            T_\E \pi \xi
            v
        \bigr)
    \end{equation}
    for $\xi \in \Universal^\bullet(\liealg{g}_\C)$, compact sets $K \subseteq G$ and
    $\seminorm{p} \in \cs(V)$ constitute a defining system of seminorms for the
    topology $\Cinfty(\pi)$ inherits from $\Cinfty(G,V)$.
\end{corollary}
\begin{proof}
    Insert \eqref{eq:Intertwining} into \eqref{eq:SmoothVectorsSeminorms}.
\end{proof}

Next, we prove variants of Lemma~\ref{lem:OrbitMapContinuityAtE} and
Lemma~\ref{lem:SemidirectProductContinuousVectors} for smoothness, which provide
effective methods for determining the space $\Cinfty(\pi)$ for a given strongly
continuous Lie group representation $\pi$.
\begin{lemma}
    \label{lem:OrbitMapSmoothnessAtE}
    \index{Smooth vectors!At group unit}
    Let $G$ be a Lie group and $\pi \colon G \longrightarrow \GL(V)$ be a
    representation on some sequentially complete locally convex space $V$ and
    fix $v \in V$. Then $\pi_v \colon G \longrightarrow V$ is smooth if and only if it
    is smooth at the group unit. In this case,
    \begin{equation}
        \label{eq:OrbitMapDerivatives}
        \Lie(\xi)
        \pi_v(g)
        =
        \pi(g)
        \circ
        \Lie(\xi)\pi_v
        \at[\Big]{\E}
    \end{equation}
    holds for all $g \in G$ and $\xi \in \Universal^\bullet(\liealg{g}_\C)$.
\end{lemma}
\begin{proof}
    We check the claim for a single derivative first.
    To this end, let $g \in G$ and $\xi \in \liealg{g}$.
    Consider the difference quotient
    \begin{equation}
        \frac{\pi_v(g \exp(t\xi)) - \pi_v(g)}{t}
        =
        \pi(g)
        \biggl(
            \frac{\pi_v(\exp(t\xi)) - \pi_v(\E)}{t}
        \biggr)
        =
        \pi(g)
        \biggl(
            \frac{\pi_v(\exp(t\xi)) - v}{t}
        \biggr),
    \end{equation}
    which is convergent in $V$ by differentiablity of $\pi_v$ at the group unit and
    continuity of
    \begin{equation}
        \pi(g)
        \colon
        V \longrightarrow V.
    \end{equation}
    Consequently, differentiability of $\pi_v$ at $\E$ implies $\pi_v \in
    \Fun[1](G,V)$. Combining \eqref{eq:LieAlgebraRep} with~\eqref{eq:Intertwining}
    -- both of which indeed require only differentiability and not smoothness -- then
    leads to
    \begin{equation}
        \pi(g)
        \circ
        \Lie(\xi)\pi_v
        \at[\Big]{\E}
        =
        \pi(g)
        \circ
        T_\E \pi \xi
        \at[\Big]{v}
        =
        \Lie(\xi)
        \pi_v(g),
    \end{equation}
    which entails \eqref{eq:OrbitMapDerivatives}. For higher derivatives, we may proceed inductively and replace $\pi_v$ with $\Lie(\xi)\pi_v$ throughout the discussion, which completes the proof.
\end{proof}

As typically not all continuous vectors are smooth by
Example~\ref{ex:TranslationRepContinuousSmooth}, we may not naively copy
and paste Lemma~\ref{lem:SemidirectProductContinuousVectors} into the
smooth category. The idea is that, if we have a sufficiently invariant candidate
$W \subseteq V$ for a space consisting of smooth vectors, it suffices to check
smoothness for each factor of a semidirect product.
\begin{lemma}
    \label{lem:SemidirectProductSmoothVectors}
    \index{Smooth vectors!Generators}
    \index{Smooth vectors!Factorization}
    Let $G$ and $N$ be Lie groups and $\varphi \colon G \times N
    \longrightarrow N$ be a smooth action of~$G$ by automorphisms of $N$.
    Moreover, assume that
    \begin{equation*}
        \pi
        \colon
        N \rtimes_\varphi G
        \longrightarrow
        \GL(V)
    \end{equation*}
    is a strongly continuous group representation on some sequentially complete
    locally convex space $V$ with a $\pi \at{G}$-invariant subspace $W \subseteq V$. If
    \begin{equation}
        W
        \subseteq
        \Cinfty
        \bigl(\pi \at{N}\bigr)
        \cap
        \Cinfty
        \bigl(\pi \at{G}\bigr),
    \end{equation}
    then $W \subseteq \Cinfty(\pi)$.
\end{lemma}
\begin{proof}
    By \eqref{eq:OrbitMapSemidirectFactorization}, we have
    \begin{equation}
        \pi_v(h,g)
        =
        \pi^N_{\pi_v^G(g)}(h)
        \qquad
        \textrm{for all }
        v \in V, g \in G, h \in N,
    \end{equation}
    where $\pi^N \coloneqq \pi \at{N}$ and $\pi^G \coloneqq \pi \at{G}$. If $v \in
    W$, then the~$\pi^G$-invariance of $W$ yields
    \begin{equation}
        \pi_v(g)
        =
        \pi(g)v
        \in
        W
        \subseteq
        \Cinfty(\pi^N)
    \end{equation}
    and thus $\pi_v$ is smooth at~$(h,g)$ by smoothness of $\pi^G_v$ at $g$
    and $\pi^N_{\pi_v^G(g)}$ at $h$.
\end{proof}

Note that $W = \Cinfty(\pi)$ is always $\pi$-invariant by
Proposition~\ref{prop:SmoothVectorsRestriction},
\ref{item:SmoothVectorsRestrictionPiInvariance}. Thus a sufficiently educated guess
for $W$ does always recover all smooth vectors in
Lemma~\ref{lem:SemidirectProductSmoothVectors}. In concrete examples, one typically
has good candidates for $W$ and may then use the finiteness of the seminorms from
Corollary~\ref{cor:SmoothVectorsSeminorms} to prove the converse inclusion. We
conclude our discussion of smooth vectors with a somewhat lengthy example, which is
well-known and discussed using machinery from Fourier analysis in
\cite[Sec.~6]{goodman:1969a}. We shall instead treat it by showcasing the abstract
machinery we have developed throughout this section combined with some
rather elementary means.
\begin{example}[Schrödinger representation and Schwartz functions]
    \index{Heisenberg group}
    \index{Schrödinger representation}
    \index{Schwartz!functions}
    \index{Strong continuity!Schrödinger representation}
    \index{Smooth vectors!Schrödinger representation}
    \label{ex:Schroedinger}
    \; \\
    Let $G \coloneqq \gls{Heisenberg}$ be the Heisenberg group in one spacial
    dimension, see e.g. \cite[Ex.~9.5.20]{hilgert.neeb:2012a} for a concrete
    representation. As a set, $G$ is given by $\R^3$ with basis $(q,p,1)$. We
    shall use the multiplication
    \begin{equation}
        (q,p,\alpha)
        \cdot
        (q',p',\alpha')
        \coloneqq
        (q + q', p + p', qp' + \alpha + \alpha')
    \end{equation}
    with unit $(0,0,0)$ and inverses
    \begin{equation}
        (q,p,\alpha)^{-1}
        =
        (-q,-p,-\alpha+qp').
    \end{equation}
    Alternatively, one may view $H_1$ as the semidirect product $\R \ltimes_\gamma
    \R^2$ twisted by the group morphism
    \begin{equation}
        \gamma
        \colon
        \R \longrightarrow \Aut(\R^2), \quad
        \gamma(q)(p,\alpha)
        \coloneqq
        (p, qp+\alpha).
    \end{equation}
    Indeed,
    \begin{align}
        \gamma(q)
        \circ
        \gamma(q')
        \at[\Big]{(p,\alpha)}
        &=
        \gamma(q)
        (p,pq'+\alpha)
        =
        (p,pq+pq'+\alpha)
%        =
%        (p, (q+q')p + \alpha)
        =
        \gamma(q+q')
        \at[\Big]{(p,\alpha)}, \\
        \gamma(q)
        (p+p',\alpha+\alpha')
        &=
        (p+p',q(p+p') + \alpha + \alpha')
        =
        \gamma(q)(p,\alpha)
        +
        \gamma(q)(p',\alpha')
    \end{align}
    and
    \begin{equation}
        \bigl(
            q, (p,\alpha)
        \bigr)
        \cdot
        \bigl(
            q',(p',\alpha)
        \bigr)
        =
        \bigl(
            q + q',
            (\gamma(q)(p',\alpha))
            +
            \alpha'
        \bigr)
        =
        \bigl(
            q + q',
            qp' + \alpha + \alpha'
        \bigr)
    \end{equation}
    hold for all $q,q',p,p',\alpha,\alpha' \in \R$. This observation is convenient, as
    it allows us to utilize both Lemma~\ref{lem:SemidirectProductContinuousVectors} and
    Lemma~\ref{lem:SemidirectProductSmoothVectors}. We suppress the superfluous
    brackets in the sequel. Note also that $\R^2$ is itself a semidirect product over
    the identity -- a direct product -- and so we may decompose even
    further.

    The Lie algebra of $H_1$ is again $\R^3$, say with basis $(\mathfrak{q},\mathfrak{p},\mathfrak{1})$ with only non-trivial Lie bracket given by
    \begin{equation}
        \label{eq:HeisenbergLie}
        [\mathfrak{q},\mathfrak{p}]
        =
        \mathfrak{1}.
    \end{equation}
    In particular, the Heisenberg group is nilpotent of step $2$. We consider its
    Schrödinger representations
    \begin{equation}
        \label{eq:SchroedingerRepresentation}
        \pi
        \colon
        H_1 \longrightarrow \GL\bigl(\Ltwo(\R)\bigr), \quad
        \pi(q,p,\alpha)
        \psi
        \at[\Big]{x}
        \coloneqq
        e^{\I (\hbar \alpha + px)}
        \cdot
        \psi
        (x + q)
    \end{equation}
    on the Hilbert space $\Ltwo(\R)$ of equivance classes of square integrable
    functions on the real line. Note that the representations are indexed by a
    parameter $\hbar \neq 0$. Clearly, each of the linear mappings $\pi(q,p,\alpha)$
    is a unitary transformation of $\Ltwo(\R)$. Notably, restricting~$\pi$ to the
    subgroup $\R
    \times \{0\} \times \{0\}$ yields the translation representation $r^*$ of the
    group $\R$ from Example~\ref{ex:TranslationRepContinuousSmooth}, except that it
    now acts on $\Ltwo(\R)$ instead of $\Continuous(\R)$.

    Our first goal is to prove that $\pi$ is strongly continuous. To this end, we consider the subrepresentations $\pi_{\pos}, \pi_{\mom}$ and $\pi_\Center$ of $\R$ defined by
    \begin{equation}
        \label{eq:SchroedingerSubrepresentations}
        \pi_{\pos}(q)
        \coloneqq
        \pi(q,0,0), \quad
        \pi_{\mom}(p)
        \coloneqq
        \pi(0,p,0)
        \quad \textrm{and} \quad
        \pi_\Center(\alpha)
        \coloneqq
        \pi(0,0,\alpha),
    \end{equation}
    which directly correspond to the decomposition of $H_1$ into a semidirect product
    of three copies of $\R$. By Lemma~\ref{lem:SemidirectProductContinuousVectors} it
    suffices to check strong continuity of the three subrepresentations. As $\pi$ is
    unitary, the sets $\pi(K) \subseteq \Linear(\Ltwo(\R))$ are even bounded in
    operator norm. Namely, contained in the closed unit ball. By
    Corollary~\ref{cor:StrongContinuityOnDense} it suffices to check that $\pi_\phi$
    is continuous for $\phi$ from the dense subspace $\gls{TestFunctionsSmooth}(\R)$
    of test
    functions
    on $\R$.\footnote{The space $\Cinfty_c(\R)$ is defined as the space of smooth
    functions $\phi \in \Cinfty(\R)$ such that the \emph{support}
    \begin{equation}
        \supp(\phi)
        \coloneqq
        \bigl\{
            x \in \R
            \colon
            f(x) \neq 0
        \bigr\}^\cl
    \end{equation}
    is compact. Here, the support $\supp(\phi)$ is the smallest closed subset of
    $\R$ such that $\phi$ vanishes identically on its complement.}

    We fix $\phi \in \Cinfty_c(\R)$ and study its restricted orbit maps
    $\pi_{\pos,\phi}$, $\pi_{\mom,\phi}$ and $\pi_{\Center,\phi}$.

    We begin with $\pi_{\pos}$. By compactness of $\supp(\phi)$, the restriction $\phi \at{K}$ is uniformly continuous. As $\phi$ vanishes identically outside of $K$, this implies uniform continuity of $\phi$ on all of $\R$. Thus, given $\epsilon > 0$ there exists a $\delta > 0$ such that $\abs{\phi(x) - \phi(y)} \le \epsilon$ for all $x,y \in \R$ with $\abs{x-y} \le \delta$. If now $(q_n)_n \subseteq \R$ is a zero sequence, then $\abs{q_n} \le \delta$ for almost all $n \in \N$. This yields the estimate
    \begin{equation}
        \norm[\big]
        {
            \pi_{\pos}(q_n)\phi
            -
            \phi
        }^2_2
        =
        \int_\R
        \abs[\big]
        {
            \phi(x + q_n)
            -
            \phi(x)
        }^2
        \D x
        \le
        \epsilon^2
        \cdot
        \gls{measure}
        \bigl(
            \supp(\phi) + \Ball_\delta(0)
        \bigr)
    \end{equation}
    for almost all $n \in \N$, where $\operatorname{meas}(M)$ denotes the Lebesgue
    measure of a measurable subset~$M \subseteq \R$. Consequently,
    $\pi_{\pos}(q_n)\phi \rightarrow \phi$ in
    $\Ltwo(\R)$. That is, the orbit map $\pi_\phi$ is continuous at the group unit. By
    Lemma~\ref{lem:OrbitMapContinuityAtE}, this implies continuity of $\pi_{\pos,
    \phi}$ on all of~$\R$. Note that these arguments possess straightforward
    generalizions to the square integrable functions $\Ltwo(G)$ on any Lie group
    $G$ and with respect to the Haar measure without complications.

    Next, we consider $\pi_{\mom}$. Let $(p_n)_n \subseteq \R$ be a zero sequence and note
    \begin{equation}
        \label{eq:SchroedingerMomentumContinuity}
        \norm[\big]
        {
            \pi_{\mom}(p_n)\phi
            -
            \phi
        }^2_2
        =
        \int_{\R}
        \abs[\big]
        {\phi(x)}^2
        \cdot
        \abs[\big]
        {
            1 - \E^{\I p_n x}
        }^2
        \D x
        =
        \int_{\supp(\phi)}
        \abs[\big]
        {\phi(x)}^2
        \cdot
        \abs[\big]
        {
            1 - \E^{\I p_n x}
        }^2
        \D x
    \end{equation}
    for $n \in \N$. It thus remains to estimate
    \begin{equation}
        \abs[\big]
        {
            1 - \E^{\I p_n x}
        }^2
        =
        \abs[\big]
        {
            \E^{0} - \E^{\I p_n x}
        }^2
    \end{equation}
    for almost all $n \in \N$ and $x \in \supp(\phi)$. Consider the family $\{f_x\}_{x \in \supp(\phi)}$ with
    \begin{equation}
        f_x
        \colon
        \R \longrightarrow \C, \quad
        f_x(p)
        \coloneqq
        \E^{\I p x}.
    \end{equation}
    Each of the $f_x$ is smooth with derivative $f'_x(p) = \I x f_x(p)$ and thus we may estimate
    \begin{equation}
        \label{eq:SchroedingerMomentumAuxiliary}
        \abs[\big]
        {
            1 - \E^{\I p_n x}
        }
        =
        \abs[\big]
        {
            f_x(0) - f_x(p_n)
        }
        \le
        \int_0^{p_n}
        \abs[\big]
        {f_x'(t)}
        \D t
        \le
        \int_0^{p_n}
        \abs{x}
        \D t
        =
        \abs{p_n x}
    \end{equation}
    for almost all $n \in \N$. Returning to \eqref{eq:SchroedingerMomentumContinuity}, we arrive at
    \begin{equation}
        \norm[\big]
        {
            \pi_{\mom}(p_n)\phi
            -
            \phi
        }^2_2
        \le
        \abs{p_n}^2
        \int_{\supp(\phi)}
        \abs{x}^2
        \cdot
        \abs[\big]
        {\phi(x)}^2
        \D x
        =
        \abs{p_n}^2
        \cdot
        \norm[\big]
        {x^2 \phi}_2^2
    \end{equation}
    for almost all $n \in \N$. As $x^2 \phi \in \Cinfty_c(\R)$ again, this completes the proof of the strong continuity of $\pi_{\pos}$.

    The central subrepresentation $\pi_\Center$ is much simpler. Setting $x = 1$ and
    $p_n = \alpha \in \R$ in~\eqref{eq:SchroedingerMomentumAuxiliary} leads to
    \begin{equation}
        \abs[\big]
        {
            1 - \E^{\I \hbar \alpha}
        }
        \le
        \abs{\hbar \alpha}.
    \end{equation}
    Thus, if $(\alpha_n)_n \subseteq \R$ is a zero sequence, then
    \begin{equation}
        \norm[\big]
        {
            \pi_{\Center}(\alpha_n)\phi
            -
            \phi
        }^2
        =
        \abs[\big]
        {
            1 - \E^{\I \hbar \alpha_n}
        }
        \cdot
        \norm[\big]{\phi}_2
        \le
        \abs{\hbar}
        \cdot
        \abs{\alpha_n}
        \cdot
        \norm[\big]{\phi}_2.
    \end{equation}
    Consequently, $\pi_\Center$ is also strongly continuous. Putting everything
    together, Lemma~\ref{lem:SemidirectProductContinuousVectors} yields the strong
    continuity of the Schrödinger representation \eqref{eq:SchroedingerRepresentation}.

    Having established strong continuity, our next goal is to prove
    \begin{equation}
        \label{eq:SchroedingerSmoothIsSchwartz}
        \Cinfty(\pi)
        =
        \Schwartz(\R)
    \end{equation}
    as locally convex spaces, where we endow $\Cinfty(\pi)$ with the topology
    inherited from the space $\Cinfty(\R,\Ltwo(\R))$ and
    \begin{equation}
        \index{Schwartz!Laurent}
        \index{Schwartz!Space}
        \label{eq:Schwartz}
        \gls{Schwartz}(\R)
        \coloneqq
        \Bigl\{
            \phi \in \Cinfty(\R)
            \;\big|\;
            \forall_{n,m \in \N_0}
            \colon
            \seminorm{s}_{n,m}(\phi)
            \coloneqq
            \sup_{x \in \R}
            \abs[\big]
            {x^m f^{(n)}(x)}
            <
            \infty
        \Bigr\}
    \end{equation}
    is the usual Schwartz\footnote{Laurent Schwartz (1915-2002) was a French
    functional
    analyst. His theory of distributions -- in the sense of generalized functions --
    allows for a mathematically coherent and elegant treatment of theoretical
    electrodynamics and other classical field theories. They moreover serve as an
    indispensable tool in the explicit solution theory of partial differential
    equations.} space on $\R$ with defining
    system of seminorms
    $\{\seminorm{s}_{n,m}\}_{n,m \in \N_0}$. Strictly speaking, as $\Cinfty(\pi)
    \subseteq \Ltwo(\R)$, the equality \eqref{eq:SchroedingerSmoothIsSchwartz} does
    not reference $\Schwartz(\R)$ as defined in \eqref{eq:Schwartz}, but rather the
    corresponding \emph{equivalence classes} within $\Ltwo(\R)$, i.e. up to equality
    almost everywhere. In particular, one should think of the suprema in
    $\seminorm{s}_{n,m}$ as essential suprema instead and demand smoothness only
    almost everywhere. We shall suppress this subtlety in the sequel.

    As before, we split $\pi$ into the subrepresentations
    \eqref{eq:SchroedingerSubrepresentations} and apply
    Lemma~\ref{lem:SemidirectProductSmoothVectors} with the candidate $W \coloneqq
    \Cinfty_c(\R)$ first. Taking another look at
    \eqref{eq:SchroedingerRepresentation}, it is clear that $\pi$ and thus in
    particular all three subrepresentations leave $\Cinfty_c(\R)$ invariant. Note that
    we may suppress the Lie exponential when considering each subrepresentation, as it
    reduces to the identity mapping. Fix again $\phi \in \Cinfty_c(\R)$.

    Motivated by \eqref{eq:RightTranslationDifferenceQuotient} we prove that the difference quotient
    \begin{equation}
        D_q \phi
        \coloneqq
        \frac{\pi_{\pos}(q)\phi - \phi}{q}
    \end{equation}
    converges to the derivative $\phi'$ for $q \rightarrow 0$ within $\Ltwo(\R)$. To see this, let $\epsilon > 0$ and use the uniform continuity of $\phi' \in \Cinfty_c(\R)$ to choose $\delta > 0$ such that
    \begin{equation}
        \phi'(t)
        -
        \phi'(x)
        \le
        \abs[\big]
        {
            \phi'(t)
            -
            \phi'(x)
        }
        \le
        \epsilon
        \qquad
        \textrm{whenever}
        \quad
        \abs{t - x} \le \delta.
    \end{equation}
    Writing $K \coloneqq \supp(\phi) + B_\delta(0)^\cl$ and
    \begin{equation*}
        D_q(x)
        -
        \phi'(x)
        =
        \frac{1}{q}
        \int_x^{x+q}
        \phi'(t)
        \D t
        -
        \phi'(x)
        =
        \frac{1}{q}
        \int_x^{x+q}
        \phi'(t)
        -
        \phi'(x)
        \D t
    \end{equation*}
    for $x,q \in \R$ yields the estimate
    \begin{equation}
        \norm[\big]
        {
            D_q \phi
            -
            \phi'
        }_2^2
        \le
        \int_K
        \frac{1}{q^2}
        \biggl(
            \int_x^{x+q}
            \abs[\big]
            {
                \phi'(t)
                -
                \phi'(x)
            }
            \D t
        \biggr)^2
        \D x
        \le
        \operatorname{meas}(K)
        \cdot
        \epsilon^2
    \end{equation}
    for $\abs{q} \le \delta$. Hence, $D_q \phi \rightarrow \phi'$ for $q \rightarrow 0$ within $\Ltwo(\R)$. That is, $\pi_{\pos,\phi} \in \Fun[1](\R,\Ltwo(\R))$ by Lemma~\ref{lem:OrbitMapSmoothnessAtE} and we get the generator
    \begin{equation}
        \label{eq:SchroedingerPositionGenerator}
        \Lie(1)
        \pi_{\pos,\phi}
        (\E)
        =
        T_0 \pi_\pos 1 \phi
        =
        \phi'
        \in
        \Cinfty_c(\R).
    \end{equation}
    Consequently, we get $\phi \in \Cinfty(\pi_{\pos})$ by induction and we have
    proved the inclusion
    \begin{equation}
        \Cinfty_c(\R) \subseteq \Cinfty(\pi_{\pos}).
    \end{equation}
    Next, we check that the difference quotient
    \begin{equation}
        M_p(\phi)
        \coloneqq
        \frac{\pi_\mom(p) \phi - \phi}{p}
        \at[\bigg]{x}
        =
        \frac{e^{\I px} \phi(x) - \phi(x)}{p}
        =
        \frac{e^{\I px} - 1}{p}
        \cdot
        \phi(x)
    \end{equation}
    converges to $x \mapsto \I x \phi(x)$ within $\Ltwo(\R)$. To see this, we estimate
    \begin{align}
        \norm[\big]
        {M_p(\phi) - \I x \phi}^2_2
        &=
        \int_\R
        \abs[\big]
        {\phi(x)}^2
        \cdot
        \abs[\bigg]
        {
            \frac{e^{\I p x} - 1}{p}
            -
            \I x
        }^2
        \D x \\
        &\le
        \int_\R
        \abs[\big]
        {\phi(x)}^2
        \cdot
        \abs[\bigg]
        {
            \frac{1}{p}
            \sum_{k=2}^{\infty}
            \frac{(\I px)^k}{k!}
        }^2
        \D x \\
        &\le
        \int_\R
        \abs[\big]
        {\phi(x)}^2
        \cdot
        \abs[\bigg]
        {
            \sum_{k=2}^{\infty}
            \frac{(\I x)^{2k}}{k!}
        }
        \cdot
        \abs[\bigg]
        {
            \sum_{k=1}^\infty
            p^{2k}
        }
        \D x \\
        &\le
        \int_\R
        \abs[\big]
        {\phi(x)}^2
        \cdot
        \exp(2\abs{x})
        \cdot
        \frac{\abs{p}^2}{1 - \abs{p}^2} \\
        &=
        \norm[\big]
        {
            \phi \cdot \exp
        }_2^2
        \cdot
        \frac{\abs{p}^2}{1 - \abs{p}^2} \\
        &\overset{p \rightarrow 0}{\longrightarrow}
        0
    \end{align}
    for $0 < \abs{p} < 1$ by the Cauchy-Schwarz inequality. Our assumption $\phi \in
    \Cinfty_c(\R)$ ensures that $\phi \cdot \exp$ is still square integrable. By
    Lemma~\ref{lem:OrbitMapSmoothnessAtE}, we have shown
    \begin{equation}
        \pi_{\mom, \phi}
        \in
        \Fun[1]
        \bigl(
            \R, \Ltwo(\R)
        \bigr)
    \end{equation}
    and it follows that the corresponding infinitesimal generator is given by
    \begin{equation}
        \label{eq:SchroedinerMomentumGenerator}
        \Lie(1)
        \pi_{\mom, \phi}(\E)
        =
        T_0
        \pi_\mom
        1
        \phi
        =
        \I x \phi.
    \end{equation}
    As $\I x \phi \in \Cinfty_c(\R)$ again, induction yields $\Cinfty_c(\R) \subseteq \Cinfty(\pi_\mom)$.

    Repeating this type of consideration for the central subrepresentation
    moreover leads to
    \begin{equation}
        \norm[\bigg]
        {\frac{\pi_\Center(\alpha)\phi - \phi}{\alpha} - \I \hbar \phi}^2_2
        =
        \norm{\phi}_2^2
        \cdot
        \abs[\bigg]
        {
            \frac{1}{\alpha}
            \sum_{k=2}^{\infty}
            \frac{(\I \hbar \alpha)^k}{k!}
        }
        \overset{\alpha \rightarrow 0}{\longrightarrow}
        0
    \end{equation}
    by virtue of the continuity of the entire function defined by the power series. In
    particular, this works for any $\phi \in \Ltwo(\R)$ and thus $\Ltwo(\R) =
    \Cinfty(\pi_\Center)$ with generator given by $\phi \mapsto \I \hbar \phi$, i.e. a
    multiple of the identity. In passing, we note that the infinitesimal generators
    implement the canonical commutation relation
    \begin{equation}
        [T_0 \pi_\pos 1,
        T_0 \pi_\mom 1]
        =
        \frac{\D}{\D x}
        \circ
        M_{\I x}
        -
        M_{\I x}
        \circ
        \frac{\D}{\D x}
        =
        \I
        \cdot
        \id
        =
        \frac{1}{\hbar}
        T_0 \pi_\Center 1
    \end{equation}
    in accordance with \eqref{eq:HeisenbergLie}.

    Putting everything together, Lemma~\ref{lem:SemidirectProductSmoothVectors} proves
    the inclusion
    \begin{equation}
        \Cinfty_c(\R)
        \subseteq
        \Cinfty(\pi).
    \end{equation}
    We proceed by taking a closer look at the $\Cinfty(\pi)$-seminorms from
    \eqref{eq:SmoothVectorsSeminorms}. Indeed, as $\pi$ is a unitary representation
    and by virtue of \eqref{eq:SchroedingerPositionGenerator},
    \eqref{eq:SchroedinerMomentumGenerator} and
    Corollary~\ref{cor:SmoothVectorsSeminorms}, they are explicitly given by
    \begin{equation}
        \biggl(
       \int_\R
       \abs[\big]
       {
        \bigl(
            \hbar^m
            x^\ell
            \phi^{(k)}(x)
        \bigr)
       }^2
       \D x
       \biggr)^{1/2}
       =
       \biggl(
       \int_\R
       \abs[\big]
       {
           \bigl(
           \hbar^m
           x^\ell
           \phi^{(k)}(x)
           \bigr)
       }^2
       \D x
       \biggr)^{1/2}
    \end{equation}
    for $\phi \in \Cinfty(\pi)$ and $k,\ell,m \in \N_0$. We may reduce this to the defining system
    \begin{equation}
        \seminorm{r}_{n,m}(\phi)
        \coloneqq
        \biggl(
        \int_\R
        \abs[\big]
        {
            x^m
            \phi^{(n)}(x)
        }^2
        \D x
        \biggr)^{1/2}
    \end{equation}
    for $n,m \in \N_0$ and $\phi \in \Cinfty(\pi)$, which already resembles the seminorms $\seminorm{s}_{m,n}$ from \eqref{eq:Schwartz}. Indeed, we check that both systems generate the same topology on $\Cinfty_c(\R)$. To this end, we estimate
    \begin{align}
        \seminorm{r}_{n,m}(\phi)^2
        &=
        \int_\R
        \abs[\big]
        {
            x^m
            \phi^{(n)}(x)
        }^2
        \D x \\
        &\le
        \sup_{t \in \R}
        \abs[\big]
        {
            (1+t^2)
            t^m
            \phi^{(n)}(t)
        }^2
        \int_\R
        \frac{1}{(1 + x^2)^2}
        \D x \\
        &=
        \frac{\pi}{2}
        \bigl(
            \seminorm{s}_{n,m}(\phi)
            +
            \seminorm{s}_{n,m+2}(\phi)
        \bigr)^2.
    \end{align}
    For the converse estimate, we fix $n,m \in \N_0$ as well as $\phi \in \Cinfty_c(\R)$ and write
    \begin{equation}
        x^m \cdot \phi^{(n)}(x)
        =
        \int_{-\infty}^x
        \frac{\D}{\D t}
        (t^m \phi^{(n)}(t))
        \D t
        =
        m
        \int_{-\infty}^x
        t^{m-1}
        \phi^{(n)}(t)
        +
        t^m
        \phi^{(n+1)}(t)
        \D t
    \end{equation}
    for $x \in \R$. This yields
    \begin{align}
        \seminorm{s}_{n,m}(\phi)
        &\le
        m
        \int_\R
        \abs[\big]
        {
            {x}^{m-1}
            \phi^{(n)}(x)
        }
        \D x
        +
        \int_\R
        \abs[\big]
        {
            {x}^{m}
            \phi^{(n+1)}(x)
        }
        \D x \\
        &=
        m
        \int_\R
        \abs[\big]
        {
            (1+x)^2
            {x}^{m-1}
            \phi^{(n)}(x)
        }
        \cdot
        \frac{1}{1 + x^2}
        \D x
        +
        \int_\R
        \abs[\big]
        {
            (1+x^2)
            {x}^{m}
            \phi^{(n+1)}(x)
        }
        \cdot
        \frac{1}{1+x^2}
        \D x \\
        &\le
        \biggl(
            \frac{m \pi}{2}
            \int_\R
            \abs[\big]
            {
                (1 + x^2)
                x^{m-1}
                \phi^{(n)}(x)
            }^2
            \D x
        \biggr)^{1/2}
        +
        \biggl(
            \frac{\pi}{2}
            \int_\R
            \abs[\big]
            {
                (1 + x^2)
                x^{m}
                \phi^{(n+1)}(x)
            }^2
            \D x
        \biggr)^{1/2} \\
        &\le
        \frac{m\pi}{2}
        \cdot
        \bigl(
            \seminorm{r}_{m-1,n}(\phi)
            +
            \seminorm{r}_{m+1,n}(\phi)
        \bigr)
        +
        \frac{\pi}{2}
        \cdot
        \bigl(
            \seminorm{r}_{m,n+1}(\phi)
            +
            \seminorm{r}_{m+2,n+1}(\phi)
        \bigr),
    \end{align}
    where we have used the Cauchy-Schwarz and the Minkowski inequalities.
    Having established estimates in both directions, we have thus shown that the
    {$\Cinfty(\pi)$-topology} coincides with the $\Schwartz$-topology on
    $\Cinfty_c(\R)$. But $\Cinfty_c(\R)$ is dense in $\Schwartz(\R)$ with respect
    to the~$\Schwartz$-topology and $\Cinfty(\pi)$ is complete by
    Lemma~\ref{lem:SmoothVectorsTopology}, so the desired inclusion
    \begin{equation}
        \Schwartz(\R)
        \subseteq
        \Cinfty(\pi)
    \end{equation}
    follows. Conversely, every $\phi \in \Cinfty(\pi)$ necessarily
    fulfils $\seminorm{r}_{n,m}(\phi) < \infty$ and thus by our estimate also
    $\seminorm{s}_{n,m}(\phi) < \infty$.\footnote{And is in particular, necessarily,
    almost everywhere smooth as a
    function.} This yields the converse inclusion $\Cinfty(\pi) \subseteq
    \Schwartz(\R)$ and finally completes the proof of
    \eqref{eq:SchroedingerSmoothIsSchwartz}.

    This yields a conceptually quite pleasing interpretation of the Schwartz space as
    the natural habitat of naive quantum mechanics. Indeed, by the Stone-von Neumann
    Theorem \cite[Thm.~14.8]{hall:2013a}, the Schrödinger representations are, up to
    unitary equivalence, the only irreducible strongly continuous unitary
    representations of $H_1$, for which the center
    \begin{equation}
        \Center(H_1)
        =
        \{0\}
        \times
        \{0\}
        \times
        \R
    \end{equation}
    acts in a non-trivial manner. Moreover, representations corresponding to
    different $\hbar$ are unitarily inequivalent.

    \index{Gårding space!Schrödinger representation}
    We conclude the discussion of the Schrödinger representation with some remarks on
    the corresponding Gårding spaces. As in terms of formulas, $\pi_\pos$ is
    essentially the translation representation from
    Example~\ref{ex:TranslationRepContinuousSmooth}, the integrals \eqref{eq:Garding}
    for $\pi_\pos$ reproduce convolution with $\phi \in \Cinfty_c(\phi)$. The
    inclusion $\mathfrak{G}(\pi_\pos) \subseteq \Cinfty(\pi_\pos)$ then corresponds
    -- up to replacing the kernel function $\phi$ with its reflection
    $\phi(-\argument)$ -- to the fact that the convolution product
    \begin{equation}
        \gls{Convolution}
        \colon
        \Cinfty_c(\R)
        \tensor
        \Ltwo(\R)
        \longrightarrow
        \Cinfty(\pi_\pos)
        =
        W^{\infty,2}(\R)
    \end{equation}
    is well defined, where
    \begin{equation}
        W^{\infty,2}(\R)
        \coloneqq
        \bigl\{
            \phi \in \Cinfty(\R)
            \;\big|\;
            \forall_{k \in \N_0}
            \colon
            \phi^{(k)}
            \in
            \Ltwo(\R)
        \bigr\}
        \subseteq
        \Cinfty(\R)
    \end{equation}
    is the Fréchet space obtained from intersecting all Sobolev\footnote{Sergei
    Lvovich Sobolev (1908-1989) was a Soviet functional analyst working on partial
    differential equations. He introduced the notion of weak derivatives, laying the
    foundation for Laurent Schwartz's theory of distributions.}-Hilbert spaces
    $W^{k,2}(\R)$, a discussion of which can be found in
    \cite[Sec.~V.2]{werner:2002a}.\footnote{By virtue of the Sobolev embedding
    Theorem~\cite{hoermander:1990a} all functions in $W^{\infty,2}(\R)$
    always possess a smooth representative. This space was introduced and studied by
    \cite[Chapter~VI,
    §8]{schwartz:1966a} in the context of general $\Lp$-spaces, where it is denoted by
    $\mathcal{D}_{\Ltwo}$.} For the momentum representation, we get the integrals
    \begin{equation}
        \index{Paley-Wiener Theorem}
        \int_\R
        \phi(p)
        \cdot
        \bigl(
            \pi_\mom(p)
            \psi
        \bigr)
        \D p
        \at[\Big]{x}
        =
        \int_\R
        \phi(p)
        \cdot
        e^{i px}
        \cdot
        \psi(x)
        \D p
        =
        \hat{\phi}(x)
        \cdot
        \psi(x)
    \end{equation}
    for $\phi \in \Cinfty_c(\R)$, $\psi \in \Ltwo(\R)$ and $x \in \R$, where $\hat{\phi}$ is one of the many conventions for the Fourier transform of $\phi$. That is, the operator associated to $\phi$ acts by multiplication with~$\hat{\phi}$; in particular, $\psi$ only has a minor role.
    In passing, we note that $\hat{\phi}$ has an extension to an entire function of
    finite order by virtue of the Paley-Wiener Theorem \cite[Thm.~I \&
    IV]{paley.wiener:1987a} or
    \cite[Thm.~7.1.3]{hoermander:1990a}.\footnote{Which, in particular, states that
    any such
    function may be obtained as the Fourier transformation of a suitable test function
    and relates the order of the resulting entire function with the geometry of the
    support of the corresponding test function. We will return to this in
    Example~\ref{ex:PaleyWiener}.}
    The fact that $\hat{\phi} \cdot \psi$ is smooth with respect to $\pi_\mom$ implies that the corresponding seminorms are finite, i.e.
    \begin{equation}
        \int_\R
        \abs[\big]
        {
            x^m
            \cdot
            \hat{\phi}(x)
            \cdot
            \psi(x)
            \D x
        }^2
        =
        \norm[\big]
        {
            x^m
            \cdot
            \hat{\phi}
            \cdot
            \psi
        }_2^2
        <
        \infty
    \end{equation}
    for $m \in \N_0$. One may prove this by elementary means by using that
    \begin{equation}
        \Cinfty_c(\R) \subseteq \Schwartz(\R)
        \quad \textrm{and} \quad
        \widehat{\Schwartz}(\R)
        =
        \Schwartz(\R)
        \subseteq
        \Linfty(\R).
    \end{equation}
    Finally and similarly, we get
    \begin{equation}
        \int_\R
        \phi(\alpha)
        \cdot
        \bigl(
        \pi_\Center(\alpha)
        \psi
        \bigr)
        \D \alpha
        =
        \psi
        \cdot
        \int_\R
        \phi(\alpha)
        e^{i\hbar \alpha}
        \D \alpha
        =
        \psi
        \cdot
        \hat{\phi}(\hbar)
    \end{equation}
    for $\phi \in \Cinfty_c(\R)$ and $\psi \in \Ltwo(\R)$. This reflects the fact that any element of $\Ltwo(\R)$ is already smooth with respect to $\pi_\Center$ and that $\Ltwo(\R)$ is indeed a vector space.
\end{example}

Throughout our considerations about the smooth vectors of the Schrödinger representation, we came across and utilized the following criterion.
\begin{corollary}
    \index{Smooth vectors!Induction}
    Let $G$ be a Lie group and $\pi \colon G \longrightarrow \GL(V)$ be a strongly
    continuous representation on some sequentially complete locally convex space $V$.
    If $W \subseteq V$ is subspace invariant under $\pi$ such that $\pi_w \in
    \Fun[1](G,V)$ for all $w \in W$, then $W \subseteq \Cinfty(\pi)$.
\end{corollary}

\section{Entire Functions as Strongly Entire Vectors}
\label{sec:EntireVectors}
\epigraph{Once you had a good excuse, you opened the door to bad
excuses.}{\emph{Thud!} -- Terry Pratchett}
% !TeX root = ../Dissertation.tex

In Theorem~\ref{thm:ActionDifferentiation}, we have seen that every strongly
continuous Lie group representation induces a Lie algebra representation on its
space of smooth vectors, which is automatically continuous by finite
dimensionality. Our next goal is to investigate when we
can reverse this process. That is, we seek conditions, under which a Lie algebra
representation comes from a strongly continuous group representation.

Here, one principal problem is that $\Cinfty(\pi)$ is typically strictly smaller than
the original space $V$. This means that one should try to extend the integrated action
to a larger space as a second step. These questions naturally lead to the concept of
\emph{analytic} vectors, as it was conceived in \cite{nelson:1959a}. We begin with the
well-known equivalence of several feasible definitions, which is an amalgamation of
\cite[Sec.~2]{goodman:1969a}, \cite[Ch.~4.4]{warner:1972a} and
\cite[Prop.~6.2.21]{niestijl:2023a}.
\begin{proposition}[Analytic vectors]
    \index{Vector!Analytic}
    \index{Analytic vectors}
    \label{prop:AnalyticVectors}
    Let $\varrho \colon \liealg{g} \longrightarrow \Linear(V)$ be a Lie algebra
    representation on a sequentially complete Hausdorff locally
    convex space $V$. Then the following are equivalent for a vector $v \in V$:
    \begin{propositionlist}
%        \item \label{item:OrbitMapRealAnalytic}
%        The orbit map $\varrho_v \colon \liealg{g} \longrightarrow \Linear(V)$ from
%        \eqref{eq:OrbitMapAlgebra} is real analytic.
        \item \label{item:OneVariableTaylorConvergence}
        There exists a zero neighbourhood $U$ in $\liealg{g}$ such that the
        exponential series
        \begin{equation}
            \label{eq:OneVariableTaylorConvergence}
            \exp
            \bigl(
                \varrho(\xi)
            \bigr)v
            \coloneqq
            \sum_{k=0}^{\infty}
            \frac{1}{k!}
            \varrho(\xi)^k v
        \end{equation}
        converges in $V$ for every $\xi \in U$.
        \item \label{item:OneVariableTaylorMajorant}
        There exists a zero neighbourhood $U \subseteq \liealg{g}$ such that the one
        variable Taylor majorants
        \begin{equation}
            \label{eq:OneVariableTaylorMajorant}
            \sum_{k=0}^{\infty}
            \frac{1}{k!}
            \seminorm{p}
            \bigl(
            \varrho(\xi)^k v
            \bigr)
        \end{equation}
        are finite for all $\xi \in U$ and $\seminorm{p} \in \cs(V)$.
        \item \label{item:AnalyticToHolomorphicLocal}
        There exists a zero neighbourhood $U_\C \subseteq \liealg{g}_\C$ such that
        \eqref{eq:OneVariableTaylorConvergence} extends to a holomorphic function with
        respect to $\xi \in U_{\C}$.\footnote{By
        Proposition~\ref{prop:HartogRedux} we do not have to distinguish between
        Gâteaux and Fréchet holomorphy here.}
    \end{propositionlist}
    If $\varrho = T_\E \pi$ and $V \subseteq \Cinfty(\pi)$ for some strongly
    continuous group representation
    \begin{equation}
        \pi \colon G \longrightarrow \GL(W)
    \end{equation}
    of a connected Lie group $G$ on a locally convex space $W$, then
    \ref{item:OneVariableTaylorConvergence}--\ref{item:AnalyticToHolomorphicLocal}
     are equivalent to:
    \begin{propositionlist}
        \setcounter{enumi}{3}
        \item \label{item:OrbitMapAnalytic}
        The orbit map $\pi_v \colon G \longrightarrow V$ from
        \eqref{eq:OrbitMapGroup} is real analytic.
        \item \label{item:OneVariableTaylorConvergenceGroup}
        The one variable Taylor formula
        \begin{equation}
            \label{eq:OneVariableTaylorConvergenceGroup}
            \sum_{k=0}^{\infty}
            \frac{1}{k!}
            \varrho(\xi)^k v
            =
            \pi(\xi)
            v
            =
            \pi_v (\exp \xi)
        \end{equation}
        holds for every $\xi$ in some zero neighbourhood $U \subseteq \liealg{g}$.
    \end{propositionlist}
\end{proposition}
\begin{proof}
    By sequential completeness of $V$, \ref{item:OneVariableTaylorMajorant} implies
    \ref{item:OneVariableTaylorConvergence}. Assume now
    \ref{item:OneVariableTaylorConvergence} and note
    that~\eqref{eq:OneVariableTaylorConvergence} is a series of
    $k$-homogeneous polynomials by linearity of $\varrho \colon \liealg{g}
    \longrightarrow \Linear(V)$ and thus may be viewed as a power series with values
    in $V$ in the sense of Corollary~\ref{cor:GateauxPowerSeries}, proving
    \ref{item:AnalyticToHolomorphicLocal}. Conversely, assuming
    \ref{item:AnalyticToHolomorphicLocal} and denoting the holomorphic extension of
    \eqref{eq:OneVariableTaylorConvergence} by $f \colon U_\C \longrightarrow V$, the
    identity principle implies that
    \begin{equation}
        f(\xi)
        =
        \sum_{k=0}^{\infty}
        \frac{1}{k!}
        \varrho(\xi)^k v
        =
        \sum_{k=0}^{\infty}
        P_k(\xi)
    \end{equation}
    with Taylor polynomials $P_k(\xi) = \tfrac{1}{k!} \varrho(\xi)^k v$ for all $k \in
    \N_0$. Invoking the vector valued version of Hartogs' Theorem from
    Proposition~\ref{prop:HartogRedux}, we get Fréchet holomorphy and thus continuity
    of the mapping $f$. We apply the locally convex Cauchy estimates
    from Proposition~\ref{prop:CauchyEstimatesLocallyConvex} to infer the convergence
    of \eqref{eq:OneVariableTaylorMajorant}. Indeed, let $B \subseteq U_\C$ be
    a closed ball with respect to some auxiliary norm. Then $B$ is balanced as well as
    bounded. Thus,~\eqref{eq:CauchyEstimatesLocallyConvex} yields the estimate
    \begin{equation}
        \sup_{\xi \in B/2}
        \seminorm{p}
        \bigl(
            P_k(\xi)
        \bigr)
        \le
        2^{-k}
        \cdot
        \sup_{\xi \in B}
        \seminorm{p}
        \bigl(
            f(\xi)
        \bigr)
        \qquad
        \textrm{for all $k \in \N_0$.}
    \end{equation}
    This proves the convergence of \eqref{eq:OneVariableTaylorMajorant} for all $\xi
    \in U \coloneqq (B/2)^\interior \cap \liealg{g}$ at once, completing the first
    part of the proof.

    Assume now that $\varrho = T_\E \pi$ as described above. As the
    Lie exponential constitutes a chart of $G$,
    \ref{item:OneVariableTaylorConvergenceGroup} implies the analyticity of $\pi_v$ at
    the group unit. By the classical version of Taylor's Theorem, we get smoothness of
    $v$ at the group unit and thus $v \in \Cinfty(\pi)$ by virtue of
    Lemma~\ref{lem:OrbitMapSmoothnessAtE}.
    Continuity and linearity of $\pi(g)$ for any $g \in G$ and
    \eqref{eq:OneVariableTaylorConvergenceGroup} yield
    \begin{align}
        \pi_v
        \bigl(
            g\exp \xi
        \bigr)
        &=
        \pi(g)
        \pi_v
        (
            \exp \xi
        ) \\
        &=
        \pi(g)
        \sum_{k=0}^{\infty}
        \frac{1}{k!}
        \varrho(\xi)^k v \\
        &=
        \sum_{k=0}^{\infty}
        \frac{1}{k!}
        \bigl(
            \pi(g)
            \circ
            T_\E \pi (\xi)^k
        \bigr)
        v \\
        &=
        \sum_{k=0}^{\infty}
        \frac{1}{k!}
        \bigl(
        \Lie(\xi)^k
        \pi_v
        \bigr)(g)
    \end{align}
    for all $\xi \in U$, where we have also invested \eqref{eq:Intertwining}. Thus
    \ref{item:OrbitMapAnalytic} and
    \ref{item:OneVariableTaylorConvergenceGroup} are equivalent. Finally, the
    equivalence of \ref{item:OneVariableTaylorConvergence} and
    \ref{item:OneVariableTaylorConvergenceGroup} is obvious in view of $\varrho = T_\E
    \pi$.
\end{proof}

In the sequel, we call a vector analytic for $\varrho$ resp.~$\pi$ if one -- and thus
all -- of the conditions in Proposition~\ref{prop:AnalyticVectors} are fulfilled. The
spaces of all such vectors are denoted by \gls{AnalyticVectorAlgebra}
and~\gls{AnalyticVectorGroup}, respectively. Our proof reveals that
\begin{equation}
    \label{eq:AnalyticVectorsAreSmooth}
    \Comega(\pi)
    \subseteq
    \Cinfty(\pi).
\end{equation}
By definition, we moreover have $\Comega(\pi) = \Comega(T_\E \pi)$. Along the way,
we have proved the following version of the one variable Lie-Taylor formula for
orbit maps.
\begin{corollary}[Lie-Taylor, orbit maps]
    \index{Lie-Taylor!Orbit maps}
    \index{Orbit maps!Lie Taylor}
    \index{Analytic vectors!Lie-Taylor}
    Let $\pi \colon G \longrightarrow \GL(V)$ be a strongly continuous representation
    of a Lie group $G$ on a sequentially complete Hausdorff locally convex space $V$.
    Then
    \begin{equation}
        \label{eq:OrbitMapLieTaylor}
        \pi_v
        \bigl(
            g\exp \xi
        \bigr)
        =
        \sum_{k=0}^{\infty}
        \frac{1}{k!}
        \bigl(
            \Lie(\xi)^k
            \pi_v
        \bigr)
        (g)
        =
        \exp
        \bigl(
            \Lie(\xi)
        \bigr)
        \pi_v(g)
    \end{equation}
    holds for all $v \in \Comega(\pi)$, $g \in G$ and $\xi$ in some sufficiently small
    zero neighbourhood $U \subseteq \liealg{g}$.
\end{corollary}

We cannot resist to at least state Nelson's Theorem, which may
be viewed as a real analytic version of Gårding's Theorem for unitary representations,
see again Remark~\ref{rem:Garding}.
\begin{theorem}[Nelson's Theorem, {\cite[Theorem~3]{nelson:1959a}}]
    \index{Nelson, Edward}
    \index{Analytic vectors!Unitary representation}
    \index{Analytic vectors!Density}
    \label{thm:Nelson}
    Let $\pi \colon G \longrightarrow \GL(H)$ be a unitary representation of a Lie
    group $G$ on a Hilbert space $H$. Then the set of analytic vectors $\Comega(\pi)$
    is dense in $H$.
\end{theorem}

Consequently, one may often reduce considerations concerning unitary representations
to analytic vectors. We proceed by studying the restrictions of the representations $\pi$
and $\varrho$ to their spaces of analytic vectors.
\begin{corollary}[Invariance of $\Comega(\pi)$]
    \label{cor:AnalyticInvariance}
    \index{Invariance!Analytic vectors}
    \index{Analytic vectors!Invariance}
    Let $\pi \colon G \longrightarrow \GL(V)$ be a strongly continuous Lie group
    representation on a sequentially complete locally convex Hausdorff space $V$. Then
    $\pi$ and $T_\E \pi$ leave $\Comega(\pi)$ invariant, i.e. restrict to
    representations on $\Comega(\pi)$.
\end{corollary}
\begin{proof}
    Let $v \in \Comega(\pi)$, $g \in G$ and write $w \coloneqq \pi(g)v$. Then, taking
    another look at~\eqref{eq:SmoothVectorInvarianceProof} and using
    Proposition~\ref{prop:AnalyticVectors}, \ref{item:OrbitMapAnalytic}, we see that
    $\pi_{w}$ is real analytic as the composition of real analytic maps. Let now $v
    \in \Comega(\varrho)$ and $\xi \in \liealg{g}$. Using
    \eqref{eq:ActionsDoNotQuiteCommute}, we compute
    \begin{equation}
        \pi_{\varrho(\xi)v}(g)
        =
        \pi(g)
        T_\E \pi(\eta)
        v
        =
        T_\E \pi
        \bigl(
            \Ad_g \xi
        \bigr)
        \pi(g)v
        =
        \bigl(
            T_\E \pi
            \bigl(
            \Ad_g \xi
            \bigr)
            \circ
            \pi_v
        \bigr)(g)
    \end{equation}
    for $g \in G$. Thus, $\pi_{\varrho(\xi)v}$ is real analytic as a composition of
    such maps; the additional dependence on $g$ is indeed harmless as the
    adjoint
    representation is always real analytic.\footnote{The correspondence of Lie groups
    and Lie algebras even provides the explicit power series expansion in terms of an
    exponential series.}
\end{proof}

In principle, one could now proceed as for smooth vectors and identify $v \in
\Comega(\pi)$ with~$\pi_v \in \Comega(G,V)$ to obtain a locally convex topology, see
again Lemma~\ref{lem:SmoothVectorsTopology}. However, as the natural topology
$\Comega(G,V)$ is rather unpleasant from a locally convex point of view\footnote{It
involves an inductive limit, essentially over spaces of the type \eqref{eq:ComegaR},
which is typically not strict.}, we refrain from doing so.
\begin{example}[Translation Representation, II: Analytic vectors]
    \index{Translation representation!Analytic vectors}
    \index{Analytic vectors!Translation representation}
    \label{ex:TranslationRepAnalytic}
    We return to the translation representation, see again
    Example~\ref{ex:TranslationRepContinuousSmooth}.
    In analogy with~\eqref{eq:TranslationSmoothVectors} we claim that
    \begin{equation}
        \Comega(r^*)
        =
        \Comega(G).
    \end{equation}
    Indeed, it is the content of \eqref{eq:LieTaylorOneVariable} that $\Comega(G)
    \subseteq \Comega(r^*)$. If, conversely, $\phi \in \Comega(r^*)$, then~$\phi \in
    \Cinfty(G)$ by a combination of \eqref{eq:AnalyticVectorsAreSmooth} and
    \eqref{eq:TranslationSmoothVectors}. Moreover,
    \eqref{eq:OneVariableTaylorMajorant} takes the form
    \begin{equation}
        \phi
        \bigl(
        g
        \exp(\xi)
        \bigr)
        =
        r^*_\phi
        \bigl(
        \exp(\xi)
        \bigr)
        \at[\Big]{g}
        =
        \sum_{k=0}^{\infty}
        \frac{1}{k!}
        T_\E r^*(\xi)^k
        \phi
        \at[\Big]{g}
        =
        \sum_{k=0}^{\infty}
        \frac{1}{k!}
        \bigl(
        \Lie(\xi)^k
        \phi
        \bigr)(g),
    \end{equation}
    which we recognize as \eqref{eq:LieTaylorOneVariable}. This yields $\phi \in
    \Comega(G)$ and provides an alternative proof of the Lie-Taylor formula.
\end{example}

Taking another look at Proposition~\ref{prop:AnalyticVectors}, we see that the
conditions are still too weak to integrate a general Lie algebra representation
$\varrho$ to a Lie group representation on all of~$G$. The problem is already visible
for the translation representation and exactly what we have described in
Remark~\ref{rem:ExponentialOperators}. A more sophisticated example can be found in
\cite[Sec.~8, Example]{goodman:1969a}. The resolution of this particular obstruction
also comes in the same form: Demanding that we may choose $U = \liealg{g}$ throughout
Proposition~\ref{prop:AnalyticVectors} yields the notion of \emph{entire} vectors as
it was first introduced in \cite[Sec.~2]{goodman:1969a} for Banach spaces $V$.
\begin{proposition}[Entire vectors, {\cite[Prop.~2.2.7]{niestijl:2023b}}]
    \index{Entire vectors}
    \index{Vector!Entire}
    \label{prop:EntireVectors}
    Let $\varrho \colon \liealg{g} \longrightarrow \Linear(V)$ be a Lie algebra
    representation on a sequentially complete Hausdorff locally
    convex space $V$. Then the following are equivalent for a vector $v \in V$:
    \begin{propositionlist}
        \item \label{item:OneVariableTaylorConvergenceGlobal}
        The exponential series
        \begin{equation}
            \label{eq:OneVariableTaylorConvergenceGlobal}
            \exp
            \bigl(
            \varrho(\xi)
            \bigr)v
            =
            \sum_{k=0}^{\infty}
            \frac{1}{k!}
            \varrho(\xi)^k v
        \end{equation}
        converges in $V$ for every $\xi \in \liealg{g}$.
        \item \label{item:OneVariableTaylorMajorantGlobal}
        The one variable Taylor majorants
        \begin{equation}
            \label{eq:OneVariableTaylorMajorantGlobal}
            \sum_{k=0}^{\infty}
            \frac{1}{k!}
            \seminorm{p}
            \bigl(
            \varrho(\xi)^k v
            \bigr)
        \end{equation}
        are finite for all $\xi \in \liealg{g}$ and $\seminorm{p} \in \cs(V)$.
        \item \label{item:AnalyticToHolomorphicGlobal}
        The function \eqref{eq:OneVariableTaylorConvergence} extends to an entire
        function with respect to $\xi$ on $\liealg{g}_\C$.
    \end{propositionlist}
    If $\varrho = T_\E \pi$ and $V \subseteq \Cinfty(\pi)$ for some strongly
    continuous group representation
    \begin{equation}
        \pi \colon G \longrightarrow \GL(W)
    \end{equation}
    of a connected Lie group $G$ on a locally convex space $W$, then
    \ref{item:OneVariableTaylorConvergenceGlobal} through
    \ref{item:AnalyticToHolomorphicGlobal} are moreover equivalent to:
    \begin{propositionlist}
        \setcounter{enumi}{3}
        \item \label{item:OrbitMapEntireEverywhere}
        For every $g \in G$, the infinitesimal orbit map
        \begin{equation}
            \label{eq:OrbitInfinitesimal}
            \liealg{g}_\C
            \ni
            \xi
            \longrightarrow
            \pi
            \bigl(
                g \exp \xi
            \bigr)v
            \in
            V
        \end{equation}
        is an entire function.\footnote{Again, either in the sense of Gâteaux or Fréchet by
        Proposition~\ref{prop:HartogRedux}, see also Corollary~\ref{cor:GateauxEntire}.}
        \item \label{item:OrbitMapEntire}
        The infinitesimal orbit map \eqref{eq:OrbitInfinitesimal} for $g = \E$ is an
        entire function.
        \item \label{item:OneVariableTaylorConvergenceGroupGlobal}
        The one variable Taylor formula
        \begin{equation}
            \label{eq:OneVariableTaylorConvergenceGroupGlobal}
            \sum_{k=0}^{\infty}
            \frac{1}{k!}
            \varrho(\xi)^k v
            =
            \pi(\exp \xi)
            v
            =
            \pi_v
            (
            \exp \xi
            )
        \end{equation}
        holds for every $\xi \in \liealg{g}$.
    \end{propositionlist}
\end{proposition}
\begin{proof}
    Identical to Proposition~\ref{prop:AnalyticVectors} with $U = \liealg{g}$ throughout.
\end{proof}

In the sequel, we call a vector $v \in V$ \emph{entire} for $\varrho$ resp.~$\pi$  if
it fulfils one -- and thus all -- of the conditions in
Proposition~\ref{prop:EntireVectors}. Moreover, we denote the spaces of all such
vectors by~\gls{EntireVectorAlgebra} and \gls{EntireVectorGroup}, respectively.
Notably, this somewhat baroque notation differs from the one we use for the
entire functions $\Entire(G)$ on $G$. As we shall see, this is not an oversight.

The following example arises from a combination of the discussion in
\cite[Prop.~5.1]{goodman:1969a} and
\cite[Ch.~VI,~§7]{katznelson:1968a}, which also contains a detailed proof.
\begin{example}[Translations on $\Ltwo(\R)$]
    \index{Paley, Raymond}
    \index{Wiener, Norbert}
    \index{Paley-Wiener Theorem}
    \index{Entire vectors!Translations on $\Ltwo(\R)$}
    \label{ex:PaleyWiener}
    Consider the space of square integrable functions $V \coloneqq \Ltwo(\R)$ on $\R$
    endowed with the translation action, that is $\pi_{\pos}$ from
    Example~\ref{ex:Schroedinger}. Then $f \in \underline{\Entire}(\pi_{\pos})$
    if and only if is the restriction of a holomorphic function $F$ defined on a strip $\R
    \times (-r,r)$ fulfilling
    \begin{equation}
        \sup_{\abs{y} \le s}
        \int_{-\infty}^\infty
        \abs[\big]
        {
            F(x+\I y)
        }^2
        \D x
        <
        \infty
        \qquad
        \textrm{for all }
        s < t.
    \end{equation}
    Indeed, this boils down to a version of the classical
    Paley\footnote{Raymond Paley (1907-1933) was a British
    mathematician, who tragically died in a skiing accident on Deception Pass. He
    nevertheless managed to publish a total of 26 papers, making contributions to many
    areas of analysis.}-Wiener\footnote{Norbert Wiener (1894-1964) was an American
    mathematician working on stochastic processes with applications in electronic
    engineering and communications. He coined the term
    \emph{cybernetics} in an attempt to bring together the different sciences.}
    Theorem \cite[Thm.~I~\&~IV]{paley.wiener:1987a}, a modern treatment
    of which is \cite[Thm.~7.3.1]{hoermander:1990a}.
\end{example}

Returning to \eqref{eq:OneVariableTaylorMajorantGlobal}, we have also
already found a natural system of seminorms for $\underline{\Entire}(\varrho)$, namely
\begin{equation}
    \index{Seminorms!Entire vectors}
    \label{eq:SeminormsEntire}
    \gls{SeminormsEntireVectors}
    (v)
    \coloneqq
    \sum_{k=0}^{\infty}
    \frac{c^k}{k!}
    \sup_{\ell=1,\ldots,n}
    \seminorm{p}
    \bigl(
        \varrho(\basis{e}_\ell)^k
        v
    \bigr),
\end{equation}
where $c \ge 0$, $\seminorm{p} \in \cs(V)$ and $(\basis{e}_1, \ldots, \basis{e}_n)$ is
a fixed basis of $\liealg{g}$. These seminorms also appear as
\cite[(2.1)]{goodman:1969a} and one readily verifies that the resulting topology is
independent of the chosen basis by retracing Lemma~\ref{lem:MajorantVsBasis} in this
more general situation. In particular, this also endows
$\underline{\Entire}(\pi)$ with a locally convex topology by differentiation. That is,
by considering~$\varrho = T_\E \pi$, see again Theorem~\ref{thm:ActionDifferentiation}.
The system may be chosen countable if $V$ admits a countable defining system of
seminorms. Mimicking the proof of Proposition~\ref{prop:AnalyticVectors} allows us
to identify the locally convex topology induced by the seminorms
\eqref{eq:SeminormsEntire} as the topology of locally uniform convergence of the
associated mappings \eqref{eq:OneVariableTaylorConvergenceGlobal} on the complexified
Lie algebra $\liealg{g}_\C$. This is essentially a version of
Proposition~\ref{lem:SmoothVectorsTopology} for entire vectors. In analogy to
\eqref{eq:BoundedSeminorms}, we write
\begin{equation}
    \gls{SeminormsEntireVectorsExtensions}
    (f)
    \coloneqq
    \sup_{\xi \in K}
    \seminorm{p}
    \bigl(
        f(\xi)
    \bigr)
\end{equation}
for compact subsets $K \subseteq \liealg{g}_\C$, $f \in
\Continuous(\liealg{g}_\C,V)$
and $\seminorm{p} \in \cs(V)$.
\begin{proposition}
    \index{Entire vectors!Topology}
    \label{prop:EntireVectorsTopology}
    Let $\varrho \colon \liealg{g} \longrightarrow \Linear(V)$ be a Lie algebra
    representation on a sequentially complete Hausdorff locally
    convex space $V$. Then the mapping
    \begin{equation}
        \label{eq:EntireVectorsAsFunctions}
        \iota
        \colon
        \underline{\Entire}(\varrho)
        \longrightarrow
        \Holomorphic(\liealg{g}_\C,V), \quad
        \iota(v)
        \at{\xi}
        \coloneqq
        \exp
        \bigl(
            \varrho(\xi)
        \bigr)v
    \end{equation}
    is a linear embedding, that is a continuous injective linear mapping such that the
    subspace topology inherited from the inclusion
    \begin{equation}
        \iota
        \bigl(
            \underline{\Entire}(\varrho)
        \bigr)
        \subseteq
        \Holomorphic
        (\liealg{g}_\C,V)
    \end{equation}
    coincides with the $\underline{\Entire}(\varrho)$-topology.
\end{proposition}
\begin{proof}
    The linearity of $\iota$ is clear and injectivity follows at once from $\iota(v)
    \at{0} = v$. Let now $v \in \underline{\Entire}(\varrho)$, $K \coloneqq
    \Ball_r(0)^\cl \subseteq \liealg{g}_\C$ the closed $\ell^\infty$-ball of radius $r
    > 0$ and $\seminorm{p} \in \cs(V)$. Then
    \begin{align}
        \seminorm{q}_{K,\seminorm{p}}
        \bigl(
            \iota(v)
        \bigr)
        &=
        \sup_{\xi \in K}
        \seminorm{p}
        \biggl(
            \sum_{k=0}^{\infty}
            \frac{1}{k!}
            \varrho(\xi)^k
            v
        \biggr) \\
        &\le
        \sum_{k=0}^{\infty}
        \frac{1}{k!}
        \sup_{\xi \in K}
        \seminorm{p}
        \bigl(
            \varrho(\xi)^k
            v
        \bigr) \\
        &=
        \sum_{k=0}^{\infty}
        \frac{1}{k!}
        \sup_{\xi \in K}
        \seminorm{p}
        \biggl(
        \varrho
        \Bigl(
            \sum_{j=1}^{n}
            \varrho(\basis{e}_j \xi^j)
        \Bigr)^k
        v
        \biggr) \\
        &\le
        \sum_{k=0}^{\infty}
        \frac{(nr)^k}{k!}
        \sup_{j=1,\ldots,n}
        \seminorm{p}
        \bigl(
        \varrho(\basis{e}_j)^k
        v
        \bigr) \\
        &=
        \seminorm{q}_{nr,\seminorm{p}}
        (v).
    \end{align}
    For the converse inequality, recall that $\iota(v)$ is holomorphic with Taylor
    polynomials
    \begin{equation}
        P_k(\xi)
        =
        \frac{1}{k!}
        \varrho(\xi)^k
        v
        \qquad
        \textrm{for all }
        k \in \N_0.
    \end{equation}
    Thus, the locally convex Cauchy estimates \eqref{eq:CauchyEstimatesLocallyConvex}
    with $B \coloneqq K$ from above yield
    \begin{equation}
        \sup_{\xi \in K}
        \seminorm{p}
        \bigl(
            P_k(\xi)
        \bigr)
        \le
        r^{-k}
        \sup_{\xi \in rK}
        \seminorm{p}
        \bigl(
            \iota(v)
        \bigr)
        =
        r^{-k}
        \cdot
        \seminorm{q}_{rK,\seminorm{p}}
        \bigl(
            \iota(v)
        \bigr)
    \end{equation}
    for all $k \in \N_0$ and $r > 0$. As $\basis{e}_1, \ldots, \basis{e}_n \in K$, we
    therefore have
    \begin{equation}
        \seminorm{q}_{c,\seminorm{p}}
        (v)
        =
        \sum_{k=0}^{\infty}
        \frac{c^k}{k!}
        \sup_{\ell=1,\ldots,n}
        \seminorm{p}
        \bigl(
        \varrho(\basis{e}_\ell)^k
        v
        \bigr)
        \le
        \seminorm{q}_{(2c)K,\seminorm{p}}
        \bigl(
            \iota(v)
        \bigr)
        \sum_{k=0}^{\infty}
        \frac{c^k}{(2c)^k}
        =
        2
        \cdot
        \seminorm{q}_{(2c)K,\seminorm{p}}
        \bigl(
            \iota(v)
        \bigr).
    \end{equation}
    This completes the proof.
\end{proof}

Consequently, $\underline{\Entire}(\varrho)$ inherits most of the pleasant functional
analytic properties of $\Holomorphic(\liealg{g}_\C, V)$ such as nuclearity by
\cite[(50.3)]{treves:2006a} and completeness for complete $V$.
\begin{corollary}
    Let $\varrho \colon \liealg{g} \longrightarrow \Linear(V)$ be a Lie algebra
    representation on a complete Hausdorff locally convex space $V$. Then the
    space of entire vectors $\underline{\Entire}(\varrho)$ is complete with respect
    to the locally convex topology induced by the
    seminorms~\eqref{eq:SeminormsEntire}. If~$V$ is merely sequentially
    complete, the same is true for $\underline{\Entire}(\varrho)$.
\end{corollary}
\begin{proof}
    By Proposition~\ref{prop:EntireVectorsTopology}, it suffices to prove the
    closedness of $\iota(\underline{\Entire}(\varrho))$ within the ambient complete
    space $\Holomorphic(\liealg{g}_\C,V)$. To this end, let $(v_\alpha)_{\alpha \in J}
    \subseteq \underline{\Entire}(\varrho)$ be such that
    \begin{equation}
        \iota(v_\alpha)
        \longrightarrow
        \phi
        \in
        \Holomorphic(\liealg{g}_\C, V)
    \end{equation}
    locally uniformly. As this entails convergence at the origin, we get $v_\alpha
    \rightarrow \phi(0)$. By continuity of \eqref{eq:EntireVectorsAsFunctions}, this
    in turn means that $\iota(v_\alpha) \rightarrow \iota(\phi(0))$ and thus
    \begin{equation}
        \phi(\xi)
        =
        \exp
        \bigl(
            \varrho
            (\xi)
        \bigr)
        \phi(0)
        \qquad
        \textrm{for all }
        \xi \in \liealg{g}_\C,
    \end{equation}
    proving $\phi \in \iota(\underline{\Entire}(\varrho))$.
\end{proof}

\begin{remark}[Existence of entire vectors]
    \index{Entire vectors!Existence}
    \index{Entire vectors!Obstructions}
    Unlike for smooth and analytic vectors, see again Remark~\ref{rem:Garding} and
    Theorem~\ref{thm:Nelson}, there are quite general situations prohibiting the
    existence of \emph{any} entire vectors besides the zero vector. Indeed, by
    \cite[Thm.~8.1]{goodman:1969a} any unitary representation $\pi$ of a
    non compact simple Lie group\footnote{Such as the projective special unitary
    groups.},
    which does not contain the identity representation as a
    subrepresentation\footnote{This essentially states that all
    elements of the group act in a non trivial fashion.}, fulfils
    $\underline{\Entire}(\pi) = \{0\}$. This pathology in particular includes any
    group containing such a group as a subgroup. For another concrete example without
    non trivial entire vectors, we refer to \cite[Sec.~7]{goodman:1969a}. In this sense,
    Nelson's Theorem~\ref{thm:Nelson} is best possible.
\end{remark}

\begin{remark}[Entire vectors on nilpotent Lie groups]
    \index{Entire vectors!Nilpotent group}
    \index{Nilpotent!Entire vectors}
    If $G$ is nilpotent and connected, see again
    Proposition~\ref{prop:ComplexificationNilpotent},
    this point of view moreover provides a positive answer to the integration problem
    $\varrho \leadsto \pi$. Indeed, in this case $\liealg{g}$ and $G$ are
    diffeomorphic by means of the Lie exponential and
    \eqref{eq:OneVariableTaylorConvergenceGlobal} thus defines a function on all of
    $G$. As $\varrho$ is a Lie algebra morphism by assumption, this indeed yields a
    group morphism $\pi \colon G \longrightarrow \GL(V)$. This may be viewed as the
    infinite dimensional Lie correspondence for the group $\GL(V)$ and Lie algebra
    $\Linear(V)$. As an aside, we note that this also matches nicely with complexification:
    We have already seen that with $G$, also $G_\C$ is nilpotent and thus $\pi$ may be
    defined directly on $G_\C$ instead by Proposition~\ref{prop:EntireVectors},
    \ref{item:AnalyticToHolomorphicGlobal}. Formally, one may derive this by the naive
    application of the the universal property \eqref{eq:UniversalComplexificationDiagram}
    of the universal $G_\C$, now for the infinite dimensional target $\GL(V)$ of invertible
    continuous linear self-maps of $V$.
\end{remark}

Outside of this very particular situation, one may proceed as follows for connected
groups~$G$. Choose an open neighbourhood $U$ of the origin in $\liealg{g}_\C$ such
that $\exp \at{U}$ is a diffeomorphism onto its image. By openness of $\exp(U)
\subseteq G$, every $g \in G$ may then be written as
\begin{equation}
    g
    =
    \exp(\xi_1)
    \cdots
    \exp(\xi_k)
\end{equation}
with -- not at all unique -- $\xi_1, \ldots, \xi_k \in U$. Then one would
be lead to define
\begin{equation}
    \label{eq:IntegrationNaive}
    \pi(g)v
    \coloneqq
    \exp
    \bigl(
        \varrho(\xi_1)
    \bigr)
    \cdots
    \exp
    \bigl(
        \varrho(\xi_k)
    \bigr)
    v
    \qquad
    \textrm{for }
    v \in \underline{\Entire}(\varrho).
\end{equation}
This approach has two problems. Firstly, it is
not at all clear whether this depends on the choice of $\xi_1, \ldots, \xi_k$.
This is ultimately a purely topological obstruction, which may be avoided by
assuming simple connectedness of $G$, i.e. the absence of monodromy.
Concretely, we will resolve this by the Lie theoretic incarnation of the Monodromy
Theorem, as it can be found in \cite[Prop.~9.5.8]{hilgert.neeb:2012a}; this particular
formulation has the advantage that no differentiable structure on the codomain is
required.\footnote{The
group $\GL(V)$ surprisingly ill-behaved for general locally convex spaces $V$. Indeed, if
$V$ is not normable, then $\GL(V) \subseteq \Linear(V)$ is not open and there exists no
topological vector space structure such that composition becomes continuous, see
\cite[Satz~2]{maissen:1962a} for the latter result and \cite{neeb:2006a} for a
comprehensive discussion. That being said, if $V$ is Banach, then $\GL(V)$ indeed
carries the structure of a \emph{linear} infinite dimensional Lie group.}

We return to the problem at hand: Secondly, one needs the invariance of
the space of entire vectors $\underline{\Entire}(\varrho)$ under the exponentiated
action
$\exp(\varrho(\argument))$ to formulate the composition. Indeed, our proof of
Corollary~\ref{cor:AnalyticInvariance} relied heavily on the presence of $\pi$.

This brings us to our final notion of regularity: \emph{Strongly} entire vectors.
To facilitate comprehensive formulas, we extend our Lie algebra representation
$\varrho$ to an algebra morphism
\begin{equation}
    \varrho
    \colon
    \Universal^\bullet(\liealg{g}_\C) \longrightarrow
    \Linear(V)
\end{equation}
by means of the universal property of the universal enveloping algebra
$\Universal^\bullet(\liealg{g}_\C)$, as we did for the Lie
derivative in \eqref{eq:LieDerivativeOnEnveloping} and
\eqref{eq:LieDerivativeOfMulti}. Moreover, we write
\begin{equation}
    \varrho(\alpha)
    \coloneqq
    \varrho
    \bigl(
        \basis{e}_{\alpha_k}
        \tensor \cdots \tensor
        \basis{e}_{\alpha_1}
    \bigr)
    \qquad
    \textrm{for }
    \alpha \in \N_n^k.
\end{equation}
\begin{definition}[Strongly entire vector]
    \index{Vector!Strongly entire}
    \index{Strongly entire vectors}
    Let $\varrho \colon \liealg{g} \longrightarrow \Linear(V)$ be a Lie algebra
    representation on a sequentially complete Hausdorff locally
    convex space~$V$.
    \begin{definitionlist}
        \item A vector $v \in V$ is called strongly entire if
        \begin{equation}
            \label{eq:StronglyEntireSeminorms}
            \index{Seminorms!Strongly entire vectors}
            \index{Strongly entire vectors!Seminorms}
            \seminorm{q}_{c,\seminorm{p}}
            (v)
            \coloneqq
            \sum_{k=0}^{\infty}
            \frac{c^k}{k!}
            \sum_{\alpha \in \N_n^k}
            \seminorm{p}
            \bigl(
            \varrho(\alpha)
            v
            \bigr)
        \end{equation}
        is finite for all $c \ge 0$, $\seminorm{p} \in \cs(V)$ and some basis
        $(\basis{e}_1, \ldots,
        \basis{e}_n)$ of $\liealg{g}$.
        \item We denote the space of strongly entire vectors by
        \gls{StronglyEntireVectorAlgebra} and endow it with the locally convex
        topology induced by the seminorms \eqref{eq:StronglyEntireSeminorms}.
        \item If $\pi \colon G \longrightarrow \GL(V)$ is a strongly continuous Lie
        group representation, then we define
        \begin{equation}
            \gls{StronglyEntireVectorGroup}
            \coloneqq
            \Entire(T_\E \pi).
        \end{equation}
    \end{definitionlist}
\end{definition}

Retracing our argument in Lemma~\ref{lem:MajorantVsBasis}, we see that
the resulting notion is independent on the choice of basis. By what we have shown and
the definitions, we arrive at the -- disregarding the lack of topology for the analytic
vectors -- continuous inclusions
\begin{equation}
    V
    =
    \Continuous(\varrho)
    \supseteq
    \Cinfty(\varrho)
    \supseteq
    \Comega(\varrho)
    \supseteq
    \underline{\Entire}(\varrho)
    \supseteq
    \Entire(\varrho)
\end{equation}
and
\begin{equation}
    \index{Vectors!Inclusions}
    \label{eq:VectorInclusions}
    V
    =
    \Continuous(\pi)
    \supseteq
    \Cinfty(\pi)
    \supseteq
    \Comega(\pi)
    \supseteq
    \underline{\Entire}(\pi)
    \supseteq
    \Entire(\pi).
\end{equation}
We chose this particular version of the definition to stress the connection with
$\Entire(G)$. Thus, our first goal is to connect our definition with alternatives used
in the literature such as \cite{goodman:1969a, goodman:1970a, goodman:1971a},
\cite{penney:1974a} or
\cite{niestijl:2023a}. To this end, we note the following elementary lemma, which
states that we may replace the inner
summation in \eqref{eq:StronglyEntireSeminorms} with a supremum over any bounded
subset within $\liealg{g}$.
\begin{lemma}
    Let $G$ be a connected Lie group and
    \begin{equation}
        \index{Seminorms!Strongly entire vectors}
        \label{eq:StronglyEntireAlternative}
        \tilde{\seminorm{q}}_{B,\seminorm{p}}
        (v)
        \coloneqq
        \sum_{k=0}^{\infty}
        \frac{1}{k!}
        \sup_{\xi_1, \ldots, \xi_k \in B}
        \seminorm{p}
        \bigl(
            \varrho(\xi_k \tensor \cdots \tensor \xi_1)v
        \bigr)
    \end{equation}
    for bounded subsets $B \subseteq \liealg{g}$ and $\seminorm{p} \in \cs(V)$.
    Then the collection of all $\tilde{\seminorm{q}}_{B,\seminorm{p}}$ generates the
    topology of $\Entire(\varrho)$.
\end{lemma}
\begin{proof}
    Let $v \in \Entire(\varrho)$, $B \subseteq \liealg{g}$ be a bounded set and
    $\seminorm{p} \in \cs(V)$. Without loss of generality, we may assume that $B$ is
    contained in the unit ball $\Ball_1(0)$ for the $\ell^\infty$-norm corresponding
    to some basis $(\basis{e}_1, \ldots, \basis{e}_n)$ of $\liealg{g}$ by virtue of
    its finite dimensionality. We estimate
    \begin{align}
        \tilde{\seminorm{q}}_{B,\seminorm{p}}
        (v)
        &\le
        \sum_{k=0}^{\infty}
        \frac{1}{k!}
        \sup_{\xi_1, \ldots, \xi_k \in \Ball_1(0)}
        \seminorm{p}
        \bigl(
            \varrho(\xi_k \tensor \cdots \tensor \xi_1)v
        \bigr) \\
        &=
        \sum_{k=0}^{\infty}
        \frac{1}{k!}
        \sup_{\xi_1, \ldots, \xi_k \in \Ball_1(0)}
        \sum_{j_1, \ldots, j_k = 1}^{n}
        \abs[\big]
        {\xi^{j_1} \cdots \xi^{j_k}}
        \cdot
        \seminorm{p}
        \bigl(
            \varrho(\basis{e}_{j_k} \tensor \cdots \tensor \basis{e}_{j_1})v
        \bigr) \\
        &\le
        \sum_{k=0}^{\infty}
        \frac{1}{k!}
        \sum_{j_1, \ldots, j_k = 1}^{n}
        \sup_{\xi_1, \ldots, \xi_k \in \Ball_1(0)}
        \abs[\big]
        {\xi^{j_1} \cdots \xi^{j_k}}
        \cdot
        \seminorm{p}
        \bigl(
            \varrho(\basis{e}_{j_k} \tensor \cdots \tensor \basis{e}_{j_1})v
        \bigr) \\
        &=
        \sum_{k=0}^{\infty}
        \frac{1}{k!}
        \sum_{j_1, \ldots, j_k = 1}^{n}
        \seminorm{p}
        \bigl(
        \varrho(\basis{e}_{j_k} \tensor \cdots \tensor \basis{e}_{j_1})v
        \bigr) \\
        &=
        \seminorm{q}_{1,\seminorm{p}}(v).
    \end{align}
    For the converse inequality, we note $\basis{e}_1, \ldots, \basis{e}_n \in
    \Ball_1(0)$, which is bounded. Hence, if $c \ge 0$, then
    \begin{equation}
        \seminorm{q}_{c,\seminorm{p}}(v)
        \le
        \sum_{k=0}^{\infty}
        \frac{1}{k!}
        \sum_{j_1, \ldots, j_k = 1}^{n}
        \sup_{\xi_1, \ldots, \xi_k \in \Ball_c(0)}
        \seminorm{p}
        \bigl(
        \varrho(\xi_k \tensor \cdots \tensor \xi_1)v
        \bigr)
        =
        \tilde{\seminorm{q}}_{\Ball_{nc}(0),\seminorm{p}}
        (v).
        \tag*{\qed}
    \end{equation}
\end{proof}

In the setting of infinite dimensional groups, this leads to two versions of entire
vectors, both of which turn out to be useful: One with respect to bounded sets $B$ in
\eqref{eq:StronglyEntireAlternative} and another for compact ones. In
\cite[Def.~6.2.19]{niestijl:2023a} the resulting spaces are denoted by
the symbols $\mathcal{H}_\varrho^{\mathcal{O}_b}$
and~$\mathcal{H}_\varrho^{\mathcal{O}}$,
respectively. We refer to \cite[Sec.~6.3]{niestijl:2023a} for a comprehensive
discussion of these spaces in the setting of \emph{unitary representations} of regular
infinite dimensional Lie groups. It
turns out that quite a few of the results hold also in our slightly different situation, where
we can however no longer guarantee the existence of strongly entire vectors. In
particular, our definition is consistent with the original one given in
\cite[p.~61]{goodman:1969a}. The following is a generalization of
\cite[Thm~6.3.10~(3)]{niestijl:2023a} and crucially does not assume the presence of a
integrated representation $\pi$.
\begin{proposition}[Invariance of $\Entire(\varrho)$]
    \index{Invariance!Strongly entire vectors}
    \index{Strongly entire vectors!Invariance}
    \label{prop:StronglyEntireInvariance}
    Let $\varrho \colon \liealg{g} \longrightarrow \Linear(V)$ be a Lie algebra
    representation on a sequentially complete locally convex
    Hausdorff space $V$. Then $\varrho$ leaves $\Entire(\varrho)$ invariant and
    restricts to a Lie algebra representation
    \begin{equation}
        \label{eq:StronglyEntireInvariance}
        \varrho
        \colon
        \liealg{g}
        \longrightarrow
        \Linear
        \bigl(
            \Entire(\varrho)
        \bigr),
    \end{equation}
    whose space of strongly entire vectors is precisely $\Entire(\varrho)$.
\end{proposition}
\begin{proof}
    Let $v \in \Entire(\varrho)$, $\xi \in \liealg{g}$, $B \subseteq \liealg{g}$
    be a bounded subset and $\seminorm{p} \in \cs(V)$. Then
    \begin{equation}
        B_\xi
        \coloneqq
        B \cup \{\xi\}
        \subseteq
        \liealg{g}
    \end{equation}
    is still bounded and we may estimate
    \begin{align}
        \tilde{\seminorm{q}}_{B, \seminorm{p}}
        \bigl(
            \varrho(\xi)v
        \bigr)
        &=
        \sum_{k=0}^{\infty}
        \frac{1}{k!}
        \sup_{\xi_1, \ldots, \xi_k \in B}
        \seminorm{p}
        \bigl(
            \varrho(\xi_k \tensor \cdots \tensor \xi_1) \varrho(\xi)v
        \bigr) \\
        &\le
        \sum_{k=0}^{\infty}
        \frac{k+1}{k+1}
        \frac{1}{k!}
        \sup_{\xi_1, \ldots, \xi_{k+1} \in B_\xi}
        \seminorm{p}
        \bigl(
        \varrho(\xi_{k+1} \tensor \cdots \tensor \xi_1)v
        \bigr) \\
        &\le
        \sum_{k=0}^{\infty}
        \frac{2^{k+1}}{(k+1)!}
        \sup_{\xi_1, \ldots, \xi_{k+1} \in B_\xi}
        \seminorm{p}
        \bigl(
        \varrho(\xi_{k+1} \tensor \cdots \tensor \xi_1)v
        \bigr) \\
        &=
        \sum_{k=1}^{\infty}
        \frac{2^{k}}{k!}
        \sup_{\xi_1, \ldots, \xi_{k} \in B_\xi}
        \seminorm{p}
        \bigl(
        \varrho(\xi_{k} \tensor \cdots \tensor \xi_1)v
        \bigr) \\
        &\le
        \tilde{\seminorm{q}}_{2B_\xi, \seminorm{p}}(v),
    \end{align}
    which proves $\varrho(\xi)v \in \Entire(\varrho)$ and the continuity of the linear
    mapping
    \begin{equation}
        \varrho(\xi)
        \colon
        \Entire(\varrho)
        \longrightarrow
        \Entire(\varrho).
    \end{equation}
    Taking $c = 0$ in \eqref{eq:StronglyEntireSeminorms} shows the inclusion
    \begin{equation}
        \cs(V)
        \at[\Big]{\Entire(\varrho)}
        \subseteq
        \cs
        \bigl(
            \Entire(\varrho)
        \bigr),
    \end{equation}
    which implies that every strongly entire vector for
    \eqref{eq:StronglyEntireInvariance} is strongly entire for $\varrho$ as a
    representation on $V$. It remains to establish that every $v \in \Entire(\varrho)$ is
    strongly entire also with respect to the $\Entire(\varrho)$-topology.
    We once again make use of the alternative
    seminorms~\eqref{eq:StronglyEntireAlternative}. To this end, let $B',B''
    \subseteq V$ be bounded and $\seminorm{p} \in \cs(V)$. Then $B \coloneqq B' \cup
    B''$ is bounded, as well, and by assumption
    $\tilde{\seminorm{q}}_{B,\seminorm{p}}(v) < \infty$. This implies that the
    function $f \colon \C \longrightarrow \C$ defined by the power series
    \begin{equation}
        f(z)
        =
        \sum_{k=0}^\infty
        \frac{z^k}{k!}
        \sup_{\xi_1, \ldots, \xi_k \in B}
        \seminorm{p}
        \Bigl(
            \varrho
            \bigl(
                \xi_1 \tensor \cdots \tensor \xi_k
            \bigr)v
        \Bigr)
    \end{equation}
    is entire in the classical sense. Taylor expanding around $z_0 = 1$ leads to the expression
    \begin{equation}
        f(z)
        =
        \sum_{k=0}^\infty
        \frac{(z-1)^k}{k!}
        f^{(k)}(1)
        =
        \sum_{k=0}^\infty
        (z-1)^k
        \sum_{\ell=k}^\infty
        \frac{1}{(\ell-k)!}
        \sup_{\xi_1, \ldots, \xi_\ell \in B}
        \seminorm{p}
        \Bigl(
        \varrho
        \bigl(
            \xi_1 \tensor \cdots \tensor \xi_\ell
        \bigr)v
        \Bigr).
    \end{equation}
    Consequently, we obtain the estimate
    \begin{align}
        \tilde{\seminorm{q}}_{B',\tilde{\seminorm{q}}_{B'',\seminorm{p}}}(v)
        &=
        \sum_{k=0}^{\infty}
        \frac{1}{k!}
        \sup_{\xi_1, \ldots, \xi_k \in B'}
        \tilde{\seminorm{q}}_{B'',\seminorm{p}}
        \bigl(
            \varrho(\xi_k \tensor \cdots \tensor \xi_1)v
        \bigr) \\
        &=
        \sum_{k=0}^{\infty}
        \frac{1}{k!}
        \sup_{\xi_1, \ldots, \xi_{k} \in B'}
        \sum_{\ell=0}^{\infty}
        \frac{1}{\ell!}
        \sup_{\eta_1, \ldots, \eta_{\ell} \in B''}
        \seminorm{p}
        \bigl(
            \varrho(\eta_\ell \tensor \cdots \tensor \eta_1)
            \varrho(\xi_k \tensor \cdots \tensor \xi_1)v
        \bigr) \\
        &\le
        \sum_{k=0}^{\infty}
        \frac{1}{k!}
        \sum_{\ell=0}^{\infty}
        \frac{1}{\ell!}
        \sup_{\xi_1, \ldots, \xi_{k+\ell} \in B}
        \seminorm{p}
        \bigl(
            \varrho(\xi_{k+\ell} \tensor \cdots \tensor \xi_1)v
        \bigr) \\
        &=
        \sum_{k=0}^{\infty}
        \frac{1}{k!}
        \sum_{m=k}^{\infty}
        \frac{1}{(m-k)!}
        \sup_{\xi_1, \ldots, \xi_{m} \in B}
        \seminorm{p}
        \bigl(
        \varrho(\xi_{m} \tensor \cdots \tensor \xi_1)v
        \bigr) \\
        &\le
        f(2) \\
        &<
        \infty.
    \end{align}
    This completes the proof.
\end{proof}

\begin{corollary}[Integration on $\Entire(\varrho)$]
    \index{Strongly entire vectors!Integration}
    \label{cor:StronglyEntireIntegration}
    Let $\varrho \colon \liealg{g} \longrightarrow \Linear(V)$ be a Lie algebra
    representation on a sequentially complete locally convex
    Hausdorff space $V$. Then there
    exists a unique strongly continuous group representation
    \begin{equation}
        \label{eq:StronglyEntireIntegration}
        \pi
        \colon
        \widetilde{G}
        \longrightarrow
        \GL
        \bigl(
            \Entire(\varrho)
        \bigr)
    \end{equation}
    such that $T_\E \pi = \varrho \at{\Entire(\varrho)}$, where $\widetilde{G}$
    denotes the universal covering group of $G$.
\end{corollary}
\begin{proof}
    We have discussed the existence of $\widetilde{G}$ in Theorem~\ref{thm:Lie3}. Let
    $U \subseteq \liealg{g}_\C$ be an absolutely convex neighbourhood of zero such that
    the Baker-Campbell-Hausdorff series converges on $U \times U$ and define
    $\pi(\exp \xi)$ by \eqref{eq:OneVariableTaylorConvergenceGlobal} for all $\xi \in
    U$. By \cite[Cor.~6.3.11]{niestijl:2023a}, we get
    \begin{equation}
        \pi
        \bigl(
            \exp \xi \exp \eta
        \bigr)
        =
        \pi
        \bigl(
            \exp
            \operatorname{BCH}(\xi,\eta)
        \bigr)
        =
        \pi(\xi)
        \circ
        \pi(\eta)
        \qquad
        \textrm{for all }
        \xi, \eta \in U.
    \end{equation}
    By simple connectedness of $\widetilde{G}$ and the
    Monodromy Theorem as formulated in \cite[Prop.~9.5.8]{hilgert.neeb:2012a}, our
    mapping $\pi$ thus extends to a unique group morphism
    \begin{equation}
        \pi
        \colon
        \tilde{G} \longrightarrow
        \GL
        \bigl(
            \Entire(\varrho)
        \bigr).
    \end{equation}
    By construction, $\pi$ is real analytic with Taylor series around the group unit
    given by~\eqref{eq:OneVariableTaylorConvergenceGlobal}. This implies that $\pi$ is
    smooth with derivative $\varrho$ by a combination of \eqref{eq:LieAlgebraRep}
    with~\eqref{eq:OneVariableTaylorConvergenceGlobal}.
\end{proof}

Thus the notion of strong entirety turns out to be powerful enough to solve the
integration problem that motivated our study of analytic vectors. As $\pi$ is a group
morphism,~\eqref{eq:IntegrationNaive} turns out to be correct a posteriori. In the
sequel, we will work with the group representation directly. Incorporating our
considerations on universal complexifications from
Section~\ref{sec:UniversalComplexification}, we may now complexify by utilizing that
\eqref{eq:StronglyEntireIntegration} may be represented by the power series
\eqref{eq:OneVariableTaylorConvergenceGlobal}. This generalizes the universal property
\eqref{eq:UniversalComplexificationDiagram} to the infinite dimensional situation.
\begin{corollary}[Complexification on $\Entire(\varrho)$]
    \label{cor:StronglyEntireComplexification}
    \index{Strongly entire vectors!Complexification}
    Let $\pi \colon G \longrightarrow \GL(V)$ be a strongly continuous
    representation of a connected and simply connected Lie group on a sequentially
    complete locally convex space $V$. Write $\eta \colon G \longrightarrow G_\C$ for
    the universal complexification of~$G$. Then there exists a unique strongly
    continuous and holomorphic group representation
    \begin{equation}
        \label{eq:StronglyEntireComplexification}
        \pi_\C
        \colon
        G_\C
        \longrightarrow
        \GL
        \bigl(
            \Entire(\varrho)
        \bigr)
    \end{equation}
    such that $\pi_\C \circ \eta = \pi$.
\end{corollary}
\begin{proof}
    By Corollary~\ref{cor:UniversalComplexificationProperties},
    \ref{item:UniversalComplexificationSimplyConnected2} the universal
    complexification $G_\C$ is again connected and simply connected and
    \begin{equation}
        T_\E \eta
        \colon
        \liealg{g} \longrightarrow \liealg{g}_\C
    \end{equation}
    is the canonical embedding. Retracing the proof of
    Corollary~\ref{cor:StronglyEntireIntegration}
    and using that \eqref{eq:OneVariableTaylorConvergenceGlobal} converges on all of
    $\liealg{g}_\C$, the Monodromy Theorem \cite[Prop.~9.5.8]{hilgert.neeb:2012a}
    yields a unique group morphism as in \eqref{eq:StronglyEntireComplexification}. It
    is holomorphic, as the Lie exponential constitutes a holomorphic chart and
    \eqref{eq:OneVariableTaylorConvergenceGlobal} is a power series. Finally,
    \begin{equation}
        T_\E
        \bigl(
            \pi_\C \circ \eta
        \bigr)
        =
        (T_\E \pi)_\C \circ T_\E \eta
        =
        T_\E \pi
    \end{equation}
    by the universal property of vector space complexification. By connectedness of
    $G$, this implies $\pi_\C \circ \eta = \pi$ as desired.
\end{proof}

The following considerations warrant our notation and generalize
Theorem~\ref{thm:RepresentativeFunctions} from representations on finite dimensional
spaces to representations on locally convex spaces, but comes at the cost of demanding
strong entirety of the involved vectors. Indeed, this is automatic if $V$ is finite
dimensional.
\begin{lemma}
    \label{lem:StronglyEntireFiniteDimensional}
    \index{Strongly entire vectors!Finite dimensions}
    Let $\pi \colon G \longrightarrow \GL_k(\C)$ be a representation.\footnote{One
    should keep in mind that representations on finite dimensional spaces are always
    continuous group morphisms by their very definition. As finite dimensional spaces
    are always barrelled, Lemma~\ref{lem:StrongContinuityContinuityOfAction} moreover
    implies that strong continuity and continuity coincide in this case.}
    Then $\Entire(\pi) = \C^k$. In particular, the inclusions
    \eqref{eq:VectorInclusions} are equalities of locally convex spaces.\footnote{Once
    again, disregarding the absence of a topology on $\Comega(\C^k)$.}
\end{lemma}
\begin{proof}
    We first note that for every $g \in G$, the corresponding mapping $\pi(g) \colon
    \C^k \longrightarrow \C^k$ is automatically continuous by its linearity. The
    crucial point is that the Lie algebra morphism
    \begin{equation}
        \varrho
        \coloneqq
        T_\E \pi
        \colon
        \liealg{g} \longrightarrow \Linear(\C^k)
    \end{equation}
    is, in particular, a linear mapping between finite dimensional vector spaces.
    Choosing a suitable auxiliary norm on $\liealg{g}$ and operator norms otherwise,
    we may achieve that the operator norm $\norm{\varrho} \le 1$ and $(\basis{e}_1,
    \ldots, \basis{e}_n) \subseteq \liealg{g}$ is a basis with $\norm{\basis{e}_k} =
    1$ for $k=1,\ldots,n$. Consequently, the submultiplicativity of both operator
    norms yields
    \begin{equation}
        \norm[\big]
        {
            \varrho(\alpha)v
        }
        =
        \norm[\big]
        {
            \bigl(
                \varrho(\basis{e}_{\alpha_k})
                \circ \cdots \circ
                \varrho(\basis{e}_{\alpha_1})
            \bigr)
            v
        }
        \le
        \norm{\varrho}
        \cdot
        \norm[\big]{\basis{e}_{\alpha_k}}
        \cdots
        \norm{\varrho}
        \cdot
        \norm[\big]{\basis{e}_{\alpha_1}}
        \cdot
        \norm{v}
        \le
        \norm{v}
    \end{equation}
    for any $\alpha \in \N_n^k$ and $v \in \C^k$. Plugging this into the seminorms
    \eqref{eq:StronglyEntireSeminorms}, we get
    \begin{equation}
        \seminorm{q}_{c,\norm{\argument}}(v)
        =
        \sum_{k=0}^{\infty}
        \frac{c^k}{k!}
        \sum_{\alpha \in \N_n^k}
        \norm[\big]
        {
            \varrho(\alpha)v
        }
        \le
        \norm{v}
        \sum_{k=0}^{\infty}
        \frac{(n \cdot c)^k}{k!}
        =
        \norm{v}
        \cdot
        \exp(nc)
        <
        \infty
    \end{equation}
    for $c \ge 0$ and $v \in V$. Thus, $\Entire(\pi) = \Entire(\varrho) = V$.
    For the final statement, we note that~$V = \C^k$ is certainly Hausdorff and thus
    the same is true for all the spaces in \eqref{eq:VectorInclusions}. Hence, all of
    them carry the unique Hausdorff locally convex topology of finite dimensional
    vector spaces, see \cite[I.~3.2]{schaefer:1999a}.
\end{proof}

Lemma~\ref{lem:StronglyEntireFiniteDimensional} may be seen as a generalization of
the fact that every continuous group morphism between finite dimensional Lie groups is
already smooth by \cite[Thm.~9.2.16]{hilgert.neeb:2012a} and, indeed, analytic
by the Baker-Campbell-Hausdorff formula.
\begin{example}
    \label{ex:FiniteDimensionsStronglyEntire}
    The continuity assumption in Lemma~\ref{lem:StronglyEntireFiniteDimensional} is
    necessary. Indeed, consider the additive reals $G \coloneqq (\R, +)$ and let $\Phi
    \colon \R \longrightarrow \R$ be a discontinuous $\field{Q}$-linear mapping. Then
    $\Phi$ is a discontinuous group morphism, say at the point $x_0 \in \R$. Moreover,
    identifying $x \in \R$ with the multiplication operator
    \begin{equation}
        m_x
        \colon
        \R \longrightarrow \R, \quad
        m_x(y)
        \coloneqq
        x \cdot y
    \end{equation}
    provides a group isomorphism $m \colon \R \longrightarrow \GL_1(\R)$. Combining
    both morphisms thus yields a representation
    \begin{equation}
        \pi
        \colon
        \R \longrightarrow \GL_1(\R), \quad
        \pi(x)y
        \coloneqq
        m_{\Phi(x)}(y)
    \end{equation}
    of $\R$ on itself, as $\pi(x) \colon \R \longrightarrow \R$ is continuous by its
    linearity. The orbit map corresponding to $y \in \R$ is however given by
    \begin{equation}
        \pi_y
        (x)
        =
        m_{\Phi(x)} y,
    \end{equation}
    which is discontinuous at $x = x_0$ for $y \coloneqq 1$. Complexifying all
    structures and extending~$\Phi$ by the identity on the imaginary parts provides a
    complex version of this phenomenon.
\end{example}

We turn towards the aforementioned generalization of
Theorem~\ref{thm:RepresentativeFunctions}.
\begin{theorem}[Universality of $\Entire_0(G)$]
    \label{thm:Universality}
    \index{Universality of $\Entire_0(G)$}
    \index{Entire functions!Universality}
    \index{Strongly entire vectors!Entire functions}
    Let $\pi \colon G \longrightarrow \GL(V)$ be a strongly continuous Lie group
    representation on a sequentially complete locally convex Hausdorff space $V$ and $v \in
    \Entire(\pi)$ a strongly entire vector. Then the corresponding representative
    functions
    \begin{equation}
        \pi_{v,\varphi}
        \colon
        G
        \longrightarrow
        \C, \quad
        \pi_{v,\varphi}(g)
        \coloneqq
        \varphi
        \bigl(
            \pi(g)v
        \bigr)
        =
        \bigl(
            \varphi
            \circ
            \pi_v
        \bigr)(g)
    \end{equation}
    are elements of $\Entire_0(G)$ for all $\varphi \in \Entire(T_\E \pi)'$.
\end{theorem}
\begin{proof}
    As usual, we write $\varrho \coloneqq T_\E \pi$ for the differentiated
    representation, see again Theorem~\ref{thm:ActionDifferentiation}.
    Let $\varphi \in \Entire(\varrho)'$.\footnote{This includes, in particular, all of $V'$. As
    $\Entire(\varrho)$ typically constitutes a proper subspace of $V$, there may however
    be
    additional continuous linear functionals available.} Using the continuity of $\varphi$
    and
    \eqref{eq:Intertwining}, we compute
    \begin{equation}
        \Lie(\xi)
        \pi_{v,\varphi}
        \at[\Big]{\E}
        =
        \varphi
        \bigl(
            (
                \Lie(\xi)
                \pi_{v}
            )(\E)
        \bigr)
        =
        \varphi
        \bigl(
            (
                \pi(\E)
                \circ
                \varrho(\xi)
            )v
        \bigr)
        =
        \varphi
        \bigl(
            \varrho(\xi)v
        \bigr)
    \end{equation}
    for $\xi \in \Universal^\bullet(\liealg{g}_\C)$. If now $v \in \Entire(\pi)$, then $\seminorm{p}
    \coloneqq \abs{\varphi} \in \cs(V)$. Consequently, we get for $c \ge 0$
    \begin{equation}
        \seminorm{q}_{0,c}
        \bigl(
            \pi_{v,\varphi}
        \bigr)
        =
        \sum_{k=0}^{\infty}
        \frac{c^k}{k!}
        \sum_{\alpha \in \N_n^k}
        \abs[\Big]
        {
            \bigl(
                \Lie(\alpha)\pi_{v,\varphi}
            \bigr)(\E)
        }
        \le
        \sum_{k=0}^{\infty}
        \frac{c^k}{k!}
        \sum_{\alpha \in \N_n^k}
        \abs[\Big]
        {
            \varphi
            \bigl(
                \varrho(\alpha)v
            \bigr)
        }
        =
        \seminorm{q}_{c,\seminorm{p}}
        (v)
        <
        \infty
    \end{equation}
    by virtue of $v \in \Entire(\pi)$. Hence, $\pi_{v,\varphi} \in \Entire_0(G)$.
\end{proof}

Combining Theorem~\ref{thm:Universality} with
Lemma~\ref{lem:StronglyEntireFiniteDimensional} provides another proof of
Theorem~\ref{thm:RepresentativeFunctions}.

Conversely, one may wonder whether $\pi_{v,\varphi} \in \Entire_0(G)$ for all $\varphi
\in V'$ -- or $\varphi \in \Entire(\varrho)'$ -- already implies $v \in \Entire(\pi)$. One
may think of this condition as a \emph{weak} version of strong entirety. The principal
problem one is faced with is the fact that a general seminorm $\seminorm{p} \in
\cs(V)$ cannot be majorized by the modulus of any $\varphi \in V'$. We leave the
problem as a question.
\begin{question}
    \index{Questions!Weak strong entirety}
    Does the converse of Theorem~\ref{thm:Universality} hold?
\end{question}

For strongly continuous unitary representations of matrix Lie groups on Hilbert
spaces, there is a positive answer by a combination of \cite[Cor.~I.6]{penney:1974a}
and Theorem~\ref{thm:Restriction}. The particular case $\pi = r^*$ yields now the
following.
\begin{example}[Translation representation, III: Strongly entire vectors]
    \index{Translation representation!Strongly entire vectors}
    \index{Strongly entire vectors!Translation representation}
     \; \\
    Let $G$ be a connected Lie group and consider $V \coloneqq
    \Continuous(G)$, for now
    endowed with the topology of pointwise convergence. Then
    \eqref{eq:StronglyEntireSeminorms} becomes
    \begin{equation}
        \label{eq:StronglyEntireVsEntireFunction}
        \seminorm{q}_{c,g}(\phi)
        =
        \sum_{k=0}^{\infty}
        \frac{c^k}{k!}
        \sum_{\alpha \in \N_0^k}
        \abs[\Big]
        {
            \bigl(
                \Lie(\alpha)\phi
            \bigr)(g)
        }
    \end{equation}
    for any $c \ge 0$ and $g \in G$ by virtue of
    \eqref{eq:RightTranslationDifferenceQuotient}. In particular, we recognize
    $\seminorm{q}_{c,\E} = \seminorm{q}_{0,c}$ as the seminorms
    inducing the $0$-topology from \eqref{eq:EntireSeminorms}. This implies
    $\Entire(r^*) \subseteq \Entire(G)$. Combining Theorem~\ref{thm:Symmetries},
    \ref{item:TranslationInvariance} with \eqref{eq:TranslationVsDifferentiation} we
    actually get equality
    \begin{equation}
        \label{eq:TranslationRepresentationStronglyEntire}
        \Entire(r^*)
        =
        \Entire(G).
    \end{equation}
    This statement is essentially \cite[Theorem~4.17, \textit{v.)}]{heins.roth.waldmann:2023a}.
    Remarkably, \eqref{eq:TranslationRepresentationStronglyEntire} also holds if one endows
    $\Continuous(G)$ with its natural topology of uniform convergence on compact subsets by
    \cite[Theorem~4.17, \textit{vi.)}]{heins.roth.waldmann:2023a}. Using again the Dirac
    functionals\footnote{That is, $\delta_g \colon \Continuous(G) \longrightarrow \C$ defined
    by $\delta_g(\phi) \coloneqq \phi(g)$.}, the representative
    functions corresponding to $\phi \in \Entire(r^*)$ are given by
    \begin{equation}
        r^*_{\phi,g_0}
        (g)
        =
        \phi
        \bigl(
            r_{g_0}(g)
        \bigr)
        =
        \phi
        \bigl(
            g g_0
        \bigr)
        \qquad
        \textrm{for all }
        g,g_0 \in G.
    \end{equation}
    Taking $g_0 = \E$ thus reproduces $\phi$, now as an element of
    $\Entire(G)$, via the universality statement in Theorem~\ref{thm:Universality}.
    This also provides another proof of the inclusion~$\Entire(r^*) \subseteq
    \Entire(G)$.
\end{example}

The attentive reader will have noticed that we have not discussed the entire vectors
of the translation representation. The reason for this is that it is generally not at
all clear whether
\begin{equation}
    \label{eq:EntireVsStronglyEntire}
    \Entire(\pi)
    \subsetneq
    \underline{\Entire}(\pi).
\end{equation}
Indeed, under certain assumptions, equality holds in \eqref{eq:EntireVsStronglyEntire}
by Penney's\footnote{Richard Penney is a harmonic analyst and
    Professor Emeritus at Purdue University. His doctoral thesis studied entire vectors
    under the supervision of Roe Goodman, who initiated their
    study in \cite{goodman:1969a,goodman:1970a,goodman:1971a}.} Theorem.
\begin{theorem}[Penney's Theorem, {\cite[Thm.~I.3]{penney:1974a}}]
    \index{Penney, Richard}
    \index{Strongly entire vectors!Penney's Theorem}
    \index{Entire vectors!Penney's Theorem}
    \index{Solvable Lie group!Penney's Theorem}
    Let $G$ be a connected, simply connected and solvable\footnote{There are subgroups
    \begin{equation}
            G
            =
            G_k
            \acts
            G_{k-1}
            \acts \cdots \acts
            G_1
            =
            \{\E\}
    \end{equation}
    such that the quotients $G_{j} / G_{j-1}$ are abelian for $j=2,\ldots,k$.} Lie
    group. Let moreover $H$ be a Hilbert space and
    \begin{equation}
        \pi
        \colon
        G
        \longrightarrow
        U(H)
    \end{equation}
    be a strongly continuous unitary representation of $G$. Then
    \begin{equation}
        \underline{\Entire}(\pi)
        =
        \Entire(\pi).
    \end{equation}
\end{theorem}
\begin{proof}[Idea]
    First of all one may use Corollary~\ref{cor:StronglyEntireComplexification} to
    extend $\pi$ to a holomorphic representation of $G_\C$. The principal idea is then
    that one may reorder
    \begin{equation}
        \varrho
        (\alpha)
        =
        \varrho
        \bigl(
            \basis{e}_k
            \tensor \cdots \tensor
            \basis{e}_k
            +
            \basis{e}_{k-1}
            \tensor \cdots \tensor
            \basis{e}_{k-1}
            +
            \cdots
            +
            \basis{e}_1
            \tensor \cdots \tensor
            \basis{e}_1
        \bigr)
        +
        \varrho(\textrm{lower order})
    \end{equation}
    for any $\alpha \in \N_n^N$ by utilizing that $\varrho$ is a Lie algebra morphism.
    In general, this results in an extraordinary amount of terms and the bookkeeping
    becomes unfeasible quickly.\footnote{This may be wonderfully illustrated by working
    out the combinatorics for the $ax+b$ group in one dimension, a pleasure the author
    will not indulge in again within this lifetime.} The additional
    assumption of solvability makes this approach viable. The concrete statement of
    \cite[Thm.~I.3]{penney:1974a} then subsumes the Lie theoretic Cauchy
    estimates, as we have used them in the proof of Theorem~\ref{thm:Restriction} for
    $\pi = r^*$.
\end{proof}

For more on this problem and its analogue for representations of infinite dimensional
Fréchet Lie groups, we refer to \cite[Rem.~3.7]{neeb:2011a}. Assuming the uniform form
of the
Cauchy estimates as concluded in \cite[Thm.~I.3]{penney:1974a} notably suffices to
construct the holomorphic extensions by \cite[Rem.~I.7]{penney:1974a}, at least if the
universal complexification morphism $\eta$ is injective. For this reason, some authors
use this then as the definition for strong entirety.
Consequently, the extension and restriction statements we have derived within
Theorem~\ref{thm:Extension}
and
Theorem~\ref{thm:Restriction} may be viewed as variants of Penney's results on
holomorphic extensions \cite[Cor.~I.2]{penney:1974a} to not necessarily unitary
representations of matrix Lie groups.

\chapter{Backmatter}
\epigraph{Whatever in creation exists without my knowledge exists without my
consent.}{\emph{Blood Meridian} -- Cormac McCarthy}
\newpage
\addcontentsline{toc}{section}{Acknowledgements}
\phantomsection
% !TeX root = ../Dissertation.tex

\Large
\textbf{Acknowledgements}
\normalsize
\bigskip

There is a multitude of important people and flavours of otherwordly entities, to
which I would like to extend heartfelt instances of ``Thank you!''. The ostensibly
straightforward and provably possible
task of finding a suitable, at least vaguely linear ordering within this collection
has been astonishingly difficult.

That being said, I am deeply grateful for the infinite patience and wisdom of my
doctoral advisors, Oliver Roth and Stefan Waldmann. Your unwavering faith and
unrivalled ability to put up with and discuss my nonsense until it became partial
sense has been instrumental throughout the past years. Your efforts have afforded me
invaluable freedom both with regards to teaching and to research. Thank you for always
having my back and taking my invariably fractal thoughts seriously. This
applies also to the time of my studies. Moreover, I would like to thank you for
affording me the time to revisit many of my leftover problems while writing my
dissertation; turning a three months project into a thirteen month one. The
alternative would have been vastly easier, but also considerably more dull, and would
have deprived us of most of the interesting results within the text.

Much gratitude is owed to my family, or what remains thereof: My mother Sabine Heins,
my girlfriend Johanna Fladung and our cat, Kugelblitz. Thank you for tolerating my
repeated absence, greatly varying degrees of distress and occasional nihilistic
megalomania. Your support has been a necessary part of this whole enterprise. I might
be horrible at showing it sometimes, but I am oh-so-very-glad to have all of you.

Moreover, I would like to thank my colleagues, Annika Moucha, Jürgen Grahl and
Matthias Schötz for always being willing to discuss about mathematics, the world,
life, the universe itself and many books, which invariably leads to my favourite
author Sir Terry Pratchett -- in case anybody has still not made this connection by
some cosmic coincidence. Stay
as wonderful as you are. Or somehow break the laws of nature to become even better.
Yes, that would be acceptable, as well.

Finally, there are my friends, whose faith in my abilities ranges between disturbing
and reassuring. Particular thanks go to Ruth Knüppel, for always being excited about
everything, including overworking herself, to Priska Dieterle for always entertaining
random mathematical ideas during Ramen and boardgame nights, to Marc
Technau for taking care of me during my repeated visits to Paderborn, and to my
raiding group for tolerating my many absences in the name of science.

\flushright{Michael}

\newpage
\bibliographystyle{nchairx}
\phantomsection
\addcontentsline{toc}{section}{Bibliography}
\bibliography{Dissertation,dqarticle,dqbook,dqthesis,dqbook,dqarticle}

\newpage
\phantomsection
\renewcommand*\pagelistname{Page} % renaming third column
\printnoidxglossary[type=symbols, title={List of
Symbols},style=super3colheaderborder]

\newpage
\phantomsection
\label{sec:LocallyConvexVocabulary}
% !TeX root = ../Dissertation.tex

\begin{center}
    \index{Balanced set}
    \index{Circled set}
    \index{Convex set}
    \index{Absolutely convex set}
    \index{Absorbing set}
    \index{Bornivorous set}
    \index{Barrel}
    \index{Baire space}
    \index{Reflexivity}
    \index{Fréchet!space}
    \index{Gâteaux!Holomorphic}
    \index{Fréchet!Holomorphic}
    \index{Boundedness!Set}
    \index{Boundedness!Mapping}
    \index{Bounded mapping}
    \index{Set!Balanced}
    \index{Set!Circled}
    \index{Set!Absolutely convex}
    \index{Set!Absorbing}
    \index{Set!Bornivorous}
    \index{Set!Barrel}
    \index{Set!Bounded}
    \Large \textbf{Terminology from locally convex analysis} \\[1cm]
    \normalsize
    \addcontentsline{toc}{section}{Locally Convex Vocabulary}
    \begin{tabular}{p{.25\textwidth} p{.67\textwidth}}
        \toprule
        $S \subseteq V$ is called...       &         ... if \\
        \midrule
         balanced or circled & for all $v \in
        S$ and $\lambda \in \C$ with $\abs{\lambda} =
        1$ also $\lambda \cdot v \in S$. \\
         convex & for all $v,w \in S$ also the line segment $[v,w]
        \subseteq S$. \\
        absolutely convex & it is balanced and convex. \\
        bounded & $\sup_{v \in S} \seminorm{p}(v) < \infty$ for all $\seminorm{p} \in
        \cs(V)$. \\
        absorbing & for every $v \in V$ there is $r > 0$ such
        that $v \in rS$. \\
        bornivorous & for every bounded $B \subseteq V$ there
        is $r > 0$ such that $B
        \subseteq rS$. \\
         barrel & it is absolutely convex, closed and absorbing. \\
        \bottomrule
    \end{tabular}

    \vspace{1cm}

    \begin{tabular}{p{.25\textwidth} p{.67\textwidth}}
        \toprule
        $V$ is called...       &         ... if \\
        \midrule
        Hausdorff & for every $v \in V$ there exists $\seminorm{p} \in \cs(V)$ with
        $\seminorm{p}(v) > 0$. \\
         Baire& every countable intersection of dense open subsets
        is dense. \\
        barrelled & every barrel is a zero neighbourhood. \\
        bornological & every absolutely convex bornivorous subset is a zero
        neighbourhood. \\
         reflexive & the canonical mapping $\iota_1 \colon V
        \longrightarrow
        (V'_\beta)'_\beta$ is an isomorphism of locally convex spaces. \\
        Montel& every bounded and closed subset is compact. \\
        nuclear & $V \tensor_\pi W \cong V \tensor_\epsilon W$ for all locally convex
        spaces $W$. \\
         Fréchet& $V$ is complete and its topology may be induced
        by a countable system
        of seminorms. \\
        LF (Limit of Fréchet)& $V$ is a strict countable direct limit of Fréchet
        spaces $F_n$ with
        $F_{n+1} \setminus F_n \neq \emptyset$. \\
        normable & the topology of $V$ may be induced by a norm. \\
        \bottomrule
    \end{tabular}

    \vspace{1cm}

    \begin{tabular}{p{.25\textwidth} p{.67\textwidth}}
        \toprule
        $f \colon U \to W$ is called...       &         ... if \\
        \midrule
         Gâteaux holomorphic & $f \at{F \cap U}$ is
        holomorphic for all finite dimensional
        $F \subseteq V$. \\
         Fréchet holomorphic & it is Gâteaux holomorphic
        and continuous. \\
         bounded & it maps bounded sets to bounded sets. \\
        locally bounded & every $v \in U$ has an open neighbourhood $U_0$ such that
        $f(U_0) \subseteq W$ is bounded. \\
        \bottomrule
    \end{tabular}
\end{center}

\newpage
\phantomsection
\addcontentsline{toc}{section}{Index}
\printindex

\end{document}